\documentclass[english,leqno]{article}

\usepackage{latexsym}
\usepackage{linearb}
\usepackage[LGR, OT1]{fontenc}
\usepackage{amsmath,amsthm}
\usepackage{amssymb}
\usepackage{MnSymbol}
\usepackage{graphicx}
\usepackage{marvosym}
\usepackage{eufrak}
\usepackage[margin=1.25in,dvips]{geometry}
\usepackage{color}
\usepackage{comment}
\usepackage{mathrsfs}\usepackage[LGR,T1]{fontenc}
\usepackage[latin9]{inputenc}
\usepackage{geometry}
\usepackage{babel}
\usepackage{enumitem}
\usepackage{amsthm}
\usepackage{amstext}
\usepackage{amssymb}
\usepackage{esint}
\usepackage[unicode=true,pdfusetitle,
 bookmarks=true,bookmarksnumbered=false,bookmarksopen=false,
 breaklinks=false,pdfborder={0 0 1},backref=false,colorlinks=false]
 {hyperref}
\usepackage{breakurl}
\makeatletter

\DeclareRobustCommand{\greektext}{%
  \fontencoding{LGR}\selectfont\def\encodingdefault{LGR}}
\DeclareRobustCommand{\textgreek}[1]{\leavevmode{\greektext #1}}
\DeclareFontEncoding{LGR}{}{}
\DeclareTextSymbol{\~}{LGR}{126}

\numberwithin{equation}{section}
\numberwithin{figure}{section}
\usepackage{enumitem}		% customizable list environments
\newcommand{\lyxaddress}[1]{
\par {\raggedright #1
\vspace{1.4em}
\noindent\par}
}
  \theoremstyle{remark}
  \newtheorem*{rem*}{\protect\remarkname}
  \theoremstyle{plain}
  \newtheorem*{thm*}{\protect\theoremname}
\theoremstyle{plain}
\newtheorem{thm}{\protect\theoremname}[section]
  \theoremstyle{plain}
  \newtheorem*{cor*}{\protect\corollaryname}
  \theoremstyle{definition}
  \newtheorem*{example*}{\protect\examplename}
  \theoremstyle{definition}
  \newtheorem{defn}[thm]{\protect\definitionname}
  \theoremstyle{plain}
  \newtheorem{lem}[thm]{\protect\lemmaname}
  \theoremstyle{plain}
  \newtheorem{cor}[thm]{\protect\corollaryname}
 \theoremstyle{definition}
 \newtheorem*{defn*}{\protect\definitionname}
  \theoremstyle{plain}
  \newtheorem*{lem*}{\protect\lemmaname}
  \theoremstyle{plain}
  \newtheorem{prop}[thm]{\protect\propositionname}

\newtheorem{innercustomthm}{Theorem}
\newenvironment{customthm}[1]
  {\renewcommand\theinnercustomthm{#1}\innercustomthm}
  {\endinnercustomthm}

\let\OldItem\item% remember the previous definition
\newcommand{\MyItem}[2][]{}%
\newenvironment{MyDescription}[1][]{%
    \renewcommand{\item}[2][]{%
        \begin{enumerate}[#1,label={##1},ref={##1}]%
            \OldItem {##2}%
        \end{enumerate}%
    }%
}{%
}%

\makeatother

\usepackage{babel}
  \providecommand{\corollaryname}{Corollary}
  \providecommand{\definitionname}{Definition}
  \providecommand{\examplename}{Example}
  \providecommand{\lemmaname}{Lemma}
  \providecommand{\propositionname}{Proposition}
  \providecommand{\remarkname}{Remark}
  \providecommand{\theoremname}{Theorem}
\providecommand{\theoremname}{Theorem}

\begin{document}
\title{The $r^{p}$-weighted energy method of Dafermos and Rodnianski\\ in general asymptotically flat spacetimes and applications}

\author{Georgios Moschidis}

\maketitle

\lyxaddress{Princeton University, Department of Mathematics, Fine Hall, Washington
Road, Princeton, NJ 08544, United States, \tt gm6@math.princeton.edu}
\begin{abstract}
In \cite{DafRod7}, Dafermos and Rodnianski presented a novel approach
to establish uniform decay rates for solutions $\text{\textgreek{f}}$
to the scalar wave equation $\square_{g}\text{\textgreek{f}}=0$ on
Minkowski, Schwarzschild and other asymptotically flat backgrounds.
This paper generalises the methods and results of \cite{DafRod7}
to a broad class of asymptotically flat spacetimes $(\mathcal{M},g)$,
including Kerr spacetimes in the full subextremal range $|a|<M$,
but also radiating spacetimes with no exact symmetries in general
dimension $d+1$, $d\ge3$. As a soft corollary, it is shown that
the Friedlander radiation field for $\text{\textgreek{f}}$ is well
defined on future null infinity. Moreover, polynomial decay rates
are established for $\text{\textgreek{f}}$, provided that an integrated
local energy decay statement (possibly with a finite loss of derivatives)
holds and the near region of $(\mathcal{M},g)$ satisfies some mild
geometric conditions. The latter conditions allow for $(\mathcal{M},g)$
to be the exterior of a black hole spacetime with a non-degenerate
event horizon (having possibly complicated topology) or the exterior
of a compact moving obstacle in an ambient globally hyperbolic spacetime
satisfying suitable geometric conditions. 
\end{abstract}
\tableofcontents{}

\section{Introduction}

The covariant wave equation 
\begin{equation}
\square_{g}\text{\textgreek{f}}=\frac{1}{\sqrt{-g}}\partial_{\text{\textgreek{m}}}\big(g^{\text{\textgreek{m}\textgreek{n}}}\sqrt{-g}\partial_{\text{\textgreek{n}}}\text{\textgreek{f}}\big)=0,\label{eq:WaveEquation}
\end{equation}
where $g$ is the Lorentzian metric of a background manifold $\mathcal{M}$,
arises in various areas of mathematical physics, including fluid mechanics,
where $g$ is the so called acoustical metric of a fluid in motion,
as well as general relativity, in which case $g$ corresponds to the
spacetime metric of a $3+1$ dimensional model of our universe.

Of fundamental importance in most settings where equation (\ref{eq:WaveEquation})
appears is the case where the background $(\mathcal{M},g)$ is flat
or almost flat, that is, when $\mathcal{M}=\mathbb{R}^{d+1}$ and
$g$ is the Minkowski metric $\text{\textgreek{h}}$ 
\begin{equation}
\text{\textgreek{h}}=-dt^{2}+(dx^{1})^{2}+\ldots+(dx^{d})^{2}\label{eq:MinkowskiMetric}
\end{equation}
(in the usual $(t,x^{1},\ldots,x^{d})$ coordinates of $\mathbb{R}^{d+1}$)
or small perturbations of it, respectively. These are the simplest
settings for which the stability properties (i.\,e.~uniform boundedness
and decay properties) of solutions to (\ref{eq:WaveEquation}) have
been studied extensively. 

One of the most successful approaches for obtaining decay estimates
for solutions $\text{\textgreek{f}}$ to (\ref{eq:WaveEquation})
on flat or almost flat backgounds has been the so called \emph{vector
field method} (see e.\,g.~\cite{Sogge2008}), which utilises the
vector fields generating the conformal isometries of Minkowski spacetime
in two ways:
\begin{enumerate}
\item As \emph{multipliers}: For any conformally Killing vector field $X$
of $(\mathbb{R}^{d+1},\text{\textgreek{h}})$, one can multiply equation
(\ref{eq:WaveEquation}) with $X(\text{\textgreek{f}})+w\text{\textgreek{f}}$
(where $w$ is a smooth function on $\mathbb{R}^{d+1}$ depending
on the choice of $X$) and then integrate the resulting expression
over a domain $\text{\textgreek{W}}$ of $\mathbb{R}^{d+1}$ bounded
by two achronal hypersurfaces $\mathcal{S}_{1},\mathcal{S}_{2}$,
with $\mathcal{S}_{2}$ being in the future of $\mathcal{S}_{1}$
(e.\,g.~$\text{\textgreek{W}}$ can be of the form $\{0\le t\le T\}$).
For the right choices of $X,w$, performing an integration by parts
yields an identity of the form 
\begin{equation}
\int_{\mathcal{S}_{2}}\mathcal{E}_{X,w}[\text{\textgreek{f}}]=\int_{\mathcal{S}_{1}}\mathcal{E}_{X,w}[\text{\textgreek{f}}],\label{eq:EnergyIdentity}
\end{equation}
where $\mathcal{E}_{X,w}[\text{\textgreek{f}}]$ is a positive definite
weighted quadratic expression in $\text{\textgreek{f}}$ and its first
derivatives. Notice that the identity (\ref{eq:EnergyIdentity}) only
contains terms on the boundary of $\text{\textgreek{W}}$ and can
be interpreted as an estimate of the ``final'' energy norm $\int_{\mathcal{S}_{2}}\mathcal{E}_{X,w}[\text{\textgreek{f}}]$
in terms of the ``initial'' energy norm $\int_{\mathcal{S}_{1}}\mathcal{E}_{X,w}[\text{\textgreek{f}}]$.
This approach can be traced back to Morawetz (see \cite{Morawetz1962}).
\item As \emph{commutation vector fields}: For certain elements $X$ of
the algebra of conformally Killing vector fields of $\mathbb{R}^{d+1}$,
the commutator $[\square_{g},X]$ is is either $0$ or a multiple
of $\square_{g}$. Thus, equation (\ref{eq:WaveEquation}) is also
satisfied by $X\text{\textgreek{f}}$ (or even higher derivatives
of $\text{\textgreek{f}}$), and this fact allows the establishment
of $L^{2}$ estimates for higher order derivatives of $\text{\textgreek{f}}$,
which in turn yield pointwise decay estimates for $\text{\textgreek{f}}$
itself through suitable global Sobolev inequalities. This approach
was initiated and developed by Klainerman (see e.\,g.~\cite{Klainerman1985},
\cite{Klainerman1987}).
\end{enumerate}
The vector field method has turned out to be especially fruitful in
the study of non linear variants of (\ref{eq:WaveEquation}), culminating
in the proof of the non linear stability of Minkowski spacetime in
\cite{Christodoulou1993}.

Preceding the use of conformally Killing vector fields $X$ as multipliers
for equation (\ref{eq:WaveEquation}), Morawetz \cite{Morawetz1961}
utilised more general first order operators generating ``positive
bulk terms'' in $\text{\textgreek{f}}$, i.\,e.~estimates for the
$L^{2}$ norm of $\text{\textgreek{f}}$ integrated over spacetime.
In particular, studying the decay properties of solutions $\text{\textgreek{f}}$
to equation (\ref{eq:WaveEquation}) on the exterior of a compact
star-shaped obstacle $\mathcal{O}$ in $\mathbb{R}^{d}$ with reflecting
boundary conditions on $\partial\mathcal{O}$, Morawetz derived an
\emph{integrated local energy decay }statement for $\text{\textgreek{f}},$
that is an estimate of the form 
\begin{equation}
\int_{0}^{\infty}\int_{\{t=\text{\textgreek{t}}\}\cap\{r\le R\}}\Big(|\partial\text{\textgreek{f}}|^{2}+|\text{\textgreek{f}}|^{2}\Big)\, dxd\text{\textgreek{t}}+\int_{0}^{\infty}\int_{\{t=\text{\textgreek{t}}\}\cap\partial\mathcal{O}}|\partial\text{\textgreek{f}}|^{2}\, d\text{\textgreek{sv}}d\text{\textgreek{t}}\lesssim\int_{\{t=0\}}|\partial\text{\textgreek{f}}|^{2}\, dx.\label{eq:ILED Morawetz}
\end{equation}
 This estimate was obtained in \cite{Morawetz1961} by using the (not
conformally Killing) radial vector field $\partial_{r}$ as a multiplier
for (\ref{eq:WaveEquation}). 

The exterior of a compact obstacle $\mathcal{O}$ in $\mathbb{R}^{d}$
(where suitable boundary conditions for solutions $\text{\textgreek{f}}$
to (\ref{eq:WaveEquation}) are imposed on the boundary $\partial\mathcal{O}$
of $\mathcal{O}$) is already an example of a background for equation
(\ref{eq:WaveEquation}) which is not a globally small perturbation
of Minkowski spacetime. More complicated examples far from Minkowski
include spacetimes $(\mathcal{M}^{d+1},g$), $d\ge3$, which contain
black hole regions, like Schwarzschild or Kerr (see \cite{DafRod6}).
Such backgrounds are of particular interest to general relativity.
One common feature that the exterior of a compact obstacle $\mathcal{O}$
in flat space and the exterior of a black hole spacetime share is
the fact that they are naturally separated into two regions where
different geometric mechanisms contribute to the long time behaviour
of solutions to (\ref{eq:WaveEquation}) on them:
\begin{itemize}
\item In the ``near'' region of these backgrounds, the long time behaviour
of solutions to (\ref{eq:WaveEquation}) is strongly affected by the
characteristics of the null geodesic flow, such as the existence of
trapped null geodesics which are reflected on the obstacle or orbit
around the black hole. In the black hole case, the existence of such
geodesics is unavoidable. A further geometric aspect of a black hole
spacetime $(\mathcal{M},g)$ which is absent in the obstacle case
is the so called \emph{event horizon} $\mathcal{H}$. In most interesting
examples, the geometric structure of $\mathcal{H}$ leads to the celebrated
\emph{red-shift} effect, which forces ``wave packets'' travelling
along the null generators of $\mathcal{H}$ to decay fast. For this
reason, the null geodesics spanning $\mathcal{H}$ are not considered
trapped in this case.%
\footnote{They are considered trapped, however, in the case when $\mathcal{H}$
is degenerate and the red-shift effect is absent, which happens in
extremal black hole spacetimes. See \cite{Aretakis2011a}.%
}
\item In the ``far away'' region of these backgrounds, there exists a
coordinate chart $(t,x^{1},\ldots,x^{d})$ in which the metric $g$
is pointwise close to the Minkowski metric $\text{\textgreek{h}}$
(\ref{eq:MinkowskiMetric}) and tends to it along all outgoing null
directions (of course, in the exterior of a compact obstacle in flat
space, $g$ is identically equal to the Minkowski metric $\text{\textgreek{h}}$
in this region). Thus, setting $r=\sqrt{(x^{1})^{2}+\ldots+(x^{d})^{2}}$,
the area of the $\{r,t=const\}$ surfaces increases to infinity along
the outgoing null directions, and this fact serves as a decay mechanism
for solutions to (\ref{eq:WaveEquation}). In particular, the quantity
$r^{\frac{d-1}{2}}\text{\textgreek{f}}$ is expected to have a finite
limit on \emph{future null infinity }$\mathcal{I}^{+}$, provided
that $\text{\textgreek{f}}$ arises from suitably decaying initial
data (see \cite{Friedlander2001} and Section \ref{sec:FriedlanderRadiation}).
Notice that on a general asymptotically flat spacetime $(\mathcal{M},g)$,
with its asymptotically flat region foliated by a set of outgoing
null hypersurfaces $\{\mathcal{S}_{\text{\textgreek{t}}}\}_{\text{\textgreek{t}}\in\mathbb{R}}$,
$\mathcal{I}^{+}$ can be abstractly defined and is parametrised by
the ``points at infinity'' of the null geodesics generating $\{\mathcal{S}_{\text{\textgreek{t}}}\}_{\text{\textgreek{t}}\in\mathbb{R}}$.
\end{itemize}
The issue of matching the estimates obtained for solutions to (\ref{eq:WaveEquation})
in different regions of a black hole spacetime implicitly appeared
in \cite{DafRod2,DafRod4,DafRod5,DafRod6,BlueSterb,Tat,TatToh1,AndBlue1,BlueSof1},
where definitive boundedness and decay estimates were established
for solutions to (\ref{eq:WaveEquation}) on Schwarzschild and very
slowly rotating Kerr exterior spacetimes (i.\,e.~for Kerr spacetimes
with angular momentum $a$ and mass $M$ satisfying the relation $|a|\ll M$).
This was achieved by the use of a Morawetz-type integrated local energy
decay statement, in conjunction with an adaptation of techniques previously
applied on flat spacetime.

In \cite{DafRod7}, Dafermos and Rodnianski suggested a more flexible
strategy for proving polynomial decay estimates for solutions to (\ref{eq:WaveEquation}),
which is explicitly tied to the aforementioned partition of a general
asymptotically flat spacetime. This approach makes use of first order
multipliers producing \emph{both} positive boundary terms (like in
(\ref{eq:EnergyIdentity})) and positive bulk terms (like in (\ref{eq:ILED Morawetz})),
and each term contains weights which grow towards $\mathcal{I}^{+}$
\emph{but are time-translation invariant}. For the sake of simplicity
of our exposition, we will discuss here the approach of \cite{DafRod7}
restricted to the case of Schwarzschild spacetime.

On the exterior of Schwarzschild spacetime $(\mathcal{M}_{Sch},g_{M})$
of mass $M$, fix the $(u,v,\text{\textgreek{sv}})$ double null coordinate
system (where $u=\frac{1}{2}(t-r^{*})$ and $v=\frac{1}{2}(t+r^{*})$,
see \cite{DafRod6}) and let $\{\mathcal{S}_{\text{\textgreek{t}}}\}_{\text{\textgreek{t}}\in\mathbb{R}}$
be a foliation of $\mathcal{M}_{Sch}$ by spacelike hypersurfaces
terminating at $\mathcal{I}^{+}$ (see Section \ref{sec:GeometrySpacetimes}
for the relevant definition), such that $\mathcal{S}_{\text{\textgreek{t}}_{2}}$
is in the future domain of dependence of $\mathcal{S}_{\text{\textgreek{t}}_{1}}$
when $\text{\textgreek{t}}_{2}>\text{\textgreek{t}}_{1}$. Let also
$T$ denote the stationary Killing field of $(\mathcal{M}_{Sch},g_{M})$,
and let $N$ be a globally timelike vector field on $\mathcal{M}_{Sch}$
such that $[T,N]=0$ and $T\equiv N$ in the far away region $\{r\gg1\}$.
Then the following estimates hold for solutions $\text{\textgreek{f}}$
to (\ref{eq:WaveEquation}) (See Section \ref{sub:Notational-conventions}
for the notations on vector field currents):
\begin{description}
\item [{Non~degenerate~energy~boundedness:}] For any $\text{\textgreek{t}}_{1}<\text{\textgreek{t}}_{2}$:
\begin{equation}
\int_{\mathcal{S}_{\text{\textgreek{t}}_{2}}}J_{\text{\textgreek{m}}}^{N}(\text{\textgreek{f}})n_{\mathcal{S}_{\text{\textgreek{t}}_{2}}}^{\text{\textgreek{m}}}\le C\cdot\int_{\mathcal{S}_{\text{\textgreek{t}}_{1}}}J_{\text{\textgreek{m}}}^{N}(\text{\textgreek{f}})n_{\mathcal{S}_{\text{\textgreek{t}}_{1}}}^{\text{\textgreek{m}}},\label{eq:ModelBoundedness}
\end{equation}
 where $n_{\mathcal{S}_{\text{\textgreek{t}}}}$ is the future directed
unit normal on the leaves of the foliation $\{\mathcal{S}_{\text{\textgreek{t}}}\}$,
and the constant $C$ in (\ref{eq:ModelBoundedness}) depends only
on the precise choice of the foliation $\{\mathcal{S}_{\text{\textgreek{t}}}\}_{\text{\textgreek{t}}\in\mathbb{R}}$
and the vectro field $N$. See \cite{DafRod2} for a proof of (\ref{eq:ModelBoundedness}).
\item [{Integrated~local~energy~decay~in~the~near~region:}] There
exists an $m>0$, such that for any $R>0$ and $\text{\textgreek{t}}\in\mathbb{R}$:
\begin{equation}
\int_{\mathcal{D}^{+}(\mathcal{S}_{\text{\textgreek{t}}})\cap\{r\le R\}}\big(|\partial\text{\textgreek{f}}|^{2}+|\text{\textgreek{f}}|^{2}\big)\, dg_{\mathcal{M}}\le C(R)\cdot\sum_{j=0}^{m}\int_{\mathcal{S}_{\text{\textgreek{t}}}}J_{\text{\textgreek{m}}}^{N}(T^{j}\text{\textgreek{f}})n_{\mathcal{S}_{\text{\textgreek{t}}}}^{\text{\textgreek{m}}},\label{eq:ModelILED}
\end{equation}
where $dg_{\mathcal{M}}$ is the spacetime volume form, $n_{\mathcal{S}_{\text{\textgreek{t}}}}$
is the future directed unit normal on $\mathcal{S}_{\text{\textgreek{t}}}$
and the constant $C(R)$ depends only on $R$ and the precise choice
of the foliation $\{\mathcal{S}_{\text{\textgreek{t}}}\}_{\text{\textgreek{t}}\in\mathbb{R}}$.
This was established in \cite{BlueSof1,DafRod2,DafRod4}.

\begin{rem*}
Notice that (\ref{eq:ModelILED}) is actually valid for $m=1$. However,
due to the existence of trapped null geodesics on $(\mathcal{M}_{Sch},g_{M})$,
the requirement that $m>0$ is necessary in this case. Notice also
that it is the red shift effect that allows the integrand in the left
hand side of (\ref{eq:ModelILED}) to be non-degenerate up to the
event horizon $\mathcal{H}$ of $(\mathcal{M}_{Sch},g_{M})$ (see
\cite{DafRod6}).
\end{rem*}
\end{description}
Using as ingredients the estimates (\ref{eq:ModelBoundedness}) and
(\ref{eq:ModelILED}), the novel approach of \cite{DafRod7} for establishing
polynomial decay rates for solutions $\text{\textgreek{f}}$ to (\ref{eq:WaveEquation})
lies in the proof of a hierarchy of $r^{p}$-weighted energy estimates
for $\text{\textgreek{f}}$ in a neighborhood of $\mathcal{I}^{+}$
and the repeated use of the pigeonhole principle on the resulting
set of estimates in order to obtain polynomial decay rates for various
weighted energies of $\text{\textgreek{f}}$. In particular, the following
result was established in \cite{DafRod7}:
\begin{thm*}
(Dafermos-Rodnianski \cite{DafRod7}, specialised here to Schwarzschild)
On Schwarzschild exterior spacetime $(\mathcal{M}_{Sch},g_{M})$,
the following statements hold for any solution $\text{\textgreek{f}}$
to the wave equation (\ref{eq:WaveEquation}):

1. An $r^{p}$-weighted energy hierarchy of the form 
\begin{multline}
\int_{\{r\ge R\}\cap\{u=\text{\textgreek{t}}_{2}\}}r^{p}|\partial_{v}(r\text{\textgreek{f}})|^{2}\, dvd\text{\textgreek{sv}}+\int_{\{r\ge R\}\cap\mathcal{D}_{\text{\textgreek{t}}_{1}}^{\text{\textgreek{t}}_{2}}}r^{p-1}\big(p|\partial_{v}(r\text{\textgreek{f}})|^{2}+(2-p)|r^{-1}\partial_{\text{\textgreek{sv}}}(r\text{\textgreek{f}})|^{2}\big)\, dudvd\text{\textgreek{sv}}\lesssim\\
\lesssim\int_{\{r\ge R\}\cap\{u=\text{\textgreek{t}}_{1}\}}r^{p}|\partial_{v}(r\text{\textgreek{f}})|^{2}\, dvd\text{\textgreek{sv}}+\int_{\{r\sim R\}\cap\mathcal{D}_{\text{\textgreek{t}}_{1}}^{\text{\textgreek{t}}_{2}}}\big(|\partial\text{\textgreek{f}}|^{2}+|\text{\textgreek{f}}|^{2}\big),\label{eq:ModelR^pHierarchy}
\end{multline}
 for $p\in[0,2]$ holds, where $\mathcal{D}_{\text{\textgreek{t}}_{1}}^{\text{\textgreek{t}}_{2}}=\{\text{\textgreek{t}}_{1}\le u\le\text{\textgreek{t}}_{2}\}$
for $\text{\textgreek{t}}_{1}<\text{\textgreek{t}}_{2}$, and the
derivatives are considered with respect to the double null coordinate
system $(u,v,\text{\textgreek{sv}})$ on $\mathcal{M}_{Sch}$. The
hierarchy (\ref{eq:ModelR^pHierarchy}) is stable under suitable perturbations
of the background metric.

2. Let $\bar{t}$ be a time function on $\mathcal{M}_{Sch}$ with
spacelike level sets intersecting the future event horizon $\mathcal{H}^{+}$
and terminating at null infinity $\mathcal{I}^{+}$, such that $T(\bar{t})=1$.
In view of (\ref{eq:ModelILED}), (\ref{eq:ModelBoundedness}) and
(\ref{eq:ModelR^pHierarchy}), $\bar{t}^{-1}$ polynomial decay estimates
hold for $\text{\textgreek{f}}$, provided its initial data on $\mathcal{S}_{0}$
(or on the hypersurface $\{t=0\}$, where $t$ is the usual Schwarzschild
exterior time coordinate) are sufficiently smooth and decaying.

3. (Schlue \cite{Schlue2013}) In the near region of $(\mathcal{M}_{Sch},g_{M})$,
$\bar{t}^{-\frac{3}{2}+\text{\textgreek{d}}}$ polynomial decay rates
for $\text{\textgreek{f}}$ hold, provided its initial data on $\mathcal{S}_{0}$
(or on the hypersurface $\{t=0\}$) are sufficiently smooth and decaying.
\end{thm*}
See \cite{DafRod7} for a more detailed description of the above result
and an explanation of how the proof immediately carries over to a
certain wider class of spacetimes. 

The goal of the present paper is to introduce a broad class of asymptotically
flat Lorentzian manifolds $(\mathcal{M}^{d+1},g)$, $d\ge3$, on which
the methods of \cite{DafRod7,Schlue2013} (suitably adapted) can be
generalised. In particular, this class (described in Section \ref{sec:GeometrySpacetimes})
is broad enough to include spacetimes which radiate Bondi mass through
future null infinity $\mathcal{I}^{+}$ and are allowed to have a
timelike boundary $\partial_{tim}\mathcal{M}$ with compact spacelike
cross-sections (modeling the boundary of a compact, possibly moving,
obstacle in an ambient globally hyperbolic spacetime). An increasing
hierarchy of geometric conditions will be imposed on this class of
spacetimes, with each additional set of conditions leading to additional
decay estimates for solutions $\text{\textgreek{f}}$ to the wave
equation (\ref{eq:WaveEquation}) on $(\mathcal{M},g)$. These conditions
are partly motivated by the geometric structure of Kerr spacetime
(and perturbations of it). 

In particular, we will establish the following three results, each
following from the previous under additional assumptions on the structure
of $(\mathcal{M},g)$:
\begin{thm*}
Let $(\mathcal{M}^{d+1},g)$, $d\ge3$, be a Lorentzian manifold with
the asymptotics (\ref{eq:RoughAsymptotics}), possibly with non-empty
timelike boundary $\partial_{tim}\mathcal{M}$ with compact spacelike
cross-sections. Then the following statements hold for any solution
$\text{\textgreek{f}}$ to the wave equation (\ref{eq:WaveEquation})
on $(\mathcal{M},g)$:

\begin{description}

\item[\hspace{0.3cm} 1.  Weighted energy hierarchy.] An $r^{p}$-weighted
energy hierarchy holds, similar to (\ref{eq:ModelR^pHierarchy}).
See Theorems \ref{thm:NewMethodFinalStatementHyperboloids} and \ref{thm:NewMethodDu+DvPhi}.

\item[\hspace{0.3cm} 2. Slow polynomial decay.] \underline{Assume}
that an integrated local energy decay statement of the form (\ref{eq:ModelILED})
holds for solutions $\text{\textgreek{f}}$ to (\ref{eq:WaveEquation})
on $(\mathcal{M},g)$ (satisfying suitable boundary conditions on
$\partial_{tim}\mathcal{M}$, if non empty). Then $\bar{t}^{-1}$
polynomial decay estimates hold for \textgreek{f}, provided its initial
data are sufficiently smooth and decaying, where $\bar{t}$ is a suitably
defined time function on $\mathcal{M}$. See Theorem \ref{thm:FirstPointwiseDecayNewMethod}.

\item[\hspace{0.3cm} 3. Improved polynomial decay.] \underline{Assume},
in addition to the previous integrated local energy decay assumption,
that $(\mathcal{M},g)$ possesses two vector fields $\{T,K\}$ (not
necessarily distinct) with timelike span and with slowly decaying
in time deformation tensor. Then provided the initial data for $\text{\textgreek{f}}$
are sufficiently smooth and decaying (and that suitable boundary conditions
have been imposed on $\partial_{tim}\mathcal{M}$):\hfill \\ \hspace{1.5cm}
- In case $d$ is odd, a $\bar{t}^{-\frac{d}{2}}$ decay rate for
\textgreek{f} and a $\bar{t}^{-\frac{d+1}{2}}$ decay rate for the
derivatives of $\text{\textgreek{f}}$ hold.\hfill \\ \hspace{1.5cm}
- In case $d$ is even, a $\bar{t}^{-\frac{d}{2}+\text{\textgreek{d}}}$
decay rate for \textgreek{f} and its derivatives holds.\hfill \\
See Theorem \ref{thm:ImprovedDecayEnergy}.

\end{description}
\end{thm*}
See also Sections \ref{sub:IntroductionHierarchy}, \ref{sub:IntroductionFirstPolynomialDecay}
and \ref{sub:IntroductionimprovedPolynomialDecay} for a more detailed
statement of Parts 1, 2 and 3 of the above theorem. 
\begin{rem*}
We should note that in fact, the integrated local energy estimate
assumed in Parts 2 and 3 of the above theorem is weaker than (\ref{eq:ModelILED}),
as we allow for an additional $\int_{\mathcal{D}^{+}(\mathcal{S}_{\text{\textgreek{t}}})}r^{-1}J_{\text{\textgreek{m}}}^{N}(T^{j}\text{\textgreek{f}})n_{\mathcal{S}}^{\text{\textgreek{m}}}$
summand on the right hand side. On general spacetimes $(\mathcal{M},g)$
with $g$ having radiating asymptotics (without satisfying any special
monotonicity condition), this additional ``error'' term appears
necessary for (\ref{eq:ModelILED}) to hold (see Sections \ref{sec:Morawetz}
and \ref{sec:Firstdecay}). Furthermore, in Part 3 above we can relax
the condition that the deformation tensors of $T,K$ decay in time,
replacing this with the statement that they are merely \emph{uniformly
$\text{\textgreek{e}}$-small}, provided there is no loss of derivatives
in the assumed integrated local energy decay estimate. In this case,
however, there is an extra $O(\text{\textgreek{e}})$ loss in the
exponents of $\bar{t}$ in the related decay estimates. See also the
remark in Section \ref{sub:IntroductionimprovedPolynomialDecay}.
\end{rem*}
As an application of Part 1 of the above theorem, we will establish
that solutions to (\ref{eq:WaveEquation}) on general asymptotically
flat spacetimes (without any assumptions posed on the structure of
their near region) have a well defined radiation field on future null
infinity $\mathcal{I}^{+}$:

\begin{customthm}{(Existence of radiation field at $\mathcal{I}^+$)}
Let $(\mathcal{M}^{d+1},g)$, $d\ge3$, be a Lorentzian manifold with
the asymptotics (\ref{eq:RoughAsymptotics}). Then for any smooth
solution $\text{\textgreek{f}}$ to (\ref{eq:WaveEquation}) with
suitably decaying intial data on a spacelike hypersurface $\text{\textgreek{S}}$
of $\mathcal{M}$ which is asymptotically of the form $\{t=const\}$,
the Friedlander radiation field $\text{\textgreek{F}}_{\mathcal{I}^{+}}$
of $\text{\textgreek{f}}$ on future null infinity: 
\begin{equation}
\text{\textgreek{F}}_{\mathcal{I}^{+}}(u,\text{\textgreek{sv}})=\lim_{r\rightarrow+\infty}\big(\text{\textgreek{W}}\cdot\text{\textgreek{f}}(u,r,\text{\textgreek{sv}})\big),
\end{equation}
where $\text{\textgreek{W}}=r^{\frac{d-1}{2}}\big(1+O(r^{-1})\big)$,
exists and is a smooth function of $(u,\text{\textgreek{sv}})$. See
Theorem \ref{thm:FriedlanderRadiation}. \end{customthm}

The assumption of an integrated local energy decay estimate for solutions
$\text{\textgreek{f}}$ to (\ref{eq:WaveEquation}), stated in Part
2 of the above theorem, does not hold on general spacetimes $(\mathcal{M},g)$
without restricting the structure of their trapped set. In particular,
in the case when $(\mathcal{M},g)$ contains a stably trapped null
geodesic, the local energy of $\text{\textgreek{f}}$ will not decay
faster than logarithmically, see e.\,g.~\cite{Ralston1971}. Hence,
in that case, no ILED statement with finite loss of derivatives (i.\,e.~of
the form (\ref{eq:ModelILED})) can hold on $(\mathcal{M},g)$. 

Even in the case where no ILED statement holds, however, the $r^{p}$-weighted
energy hierarchy (\ref{eq:ModelR^pHierarchy}) can still yield decay
estimates for $\text{\textgreek{f}}$ provided some decay estimate
for the local energy of $\text{\textgreek{f}}$ can be established.
In \cite{Moschidisb}, it is shown that on a general class of stationary
and asymptotically flat spacetimes $(\mathcal{M},g)$, the local energy
of solutions $\text{\textgreek{f}}$ to (\ref{eq:WaveEquation}) decays
logarithmically in time. Combining Part 1 of the above theorem with
the logarithmic local energy decay estimate established in \cite{Moschidisb},
we will thus be able to infer that the energy of $\text{\textgreek{f}}$
through a hyperboloidal foliation of $\mathcal{M}$ decays logarithmically
in time:

\begin{customthm}{(Logarithmic decay of the energy flux through a hyperboloidal foliation, \cite{Moschidisb})}
Let $(\mathcal{M}^{d+1},g)$, $d\ge3$, be a globally hyperbolic spacetime
with a Cauchy hypersurface $\text{\textgreek{S}}$. 

Assume that $(\mathcal{M},g)$ is stationary, with stationary Killing
field $T$, and asymptotically flat. If $\mathcal{M}$ contains a
black hole region bounded by an event horizon $\mathcal{H}$, assume
that $\mathcal{H}$ has positive surface gravity and that the ergoregion
(i.\,e.~the set where $g(T,T)>0$) is ``small'' (see \cite{Moschidisb}
for the precise statement of these assumptions). Finally, assume that
an energy boundedness statement of the form (\ref{eq:ModelBoundedness})
holds for solutions to $\square\text{\textgreek{f}}=0$ on the domain
of outer communications $\mathcal{D}$ of $\mathcal{M}$. 

It then follows that the energy flux through a $T$-translated hyperboloidal
foliation of $\mathcal{M}$ terminating at $\mathcal{I}^{+}$ of any
smooth solution $\text{\textgreek{f}}$ to (\ref{eq:WaveEquation})
on $(\mathcal{M},g)$ with suitably decaying initial data on a Cauchy
hypersurface $\text{\textgreek{S}}$ of $\mathcal{M}$ decays at least
logarithmically in time. See \cite{Moschidisb}.\end{customthm}

We will now give some examples of spacetimes $(\mathcal{M},g)$ satisfying
the assumptions of the above theorem. On these spacetimes, polynomial
decay rates for solutions to (\ref{eq:WaveEquation}) will be inferred
as a result of Parts 2 and 3 of the above theorem. 

Our first example will be the exterior region of a subextremal Kerr
spacetime (with parameters $a$, $M$ in the fulll subextremal range
$|a|<M$). This satisfies all the geometric assumptions of Parts 1,
2 and 3 of the above theorem. We should remark that, in fact, our
assumption on the properties of the vector fields $T,K$ of Part 3
of the above theorem was motivated by the geometric properties of
the subextremal Kerr family. In view of the integrated local energy
decay statement and the energy boundedness estimate established in
\cite{DafRodSchlap}, we will be able to infer Corollary 3.1 of \cite{DafRodSchlap}:

\begin{customthm}{(Polynomial decay on subextremal Kerr exterior for $|a|<M$, \cite{DafRodSchlap})}
Corollary 3.1 of \cite{DafRodSchlap} holds, that is to say, a $\bar{t}^{-\frac{3}{2}}$
pointwise decay rate for $\text{\textgreek{f}}$ and $\bar{t}^{-2}$
decay rate for the derivatives of $\text{\textgreek{f}}$ hold for
solutions $\text{\textgreek{f}}$ to the wave equation (\ref{eq:WaveEquation})
on subextremal Kerr spacetimes in the full parameter range $|a|<M$.\end{customthm}

See Section \ref{sub:KerrDecay} for a precise statement of this result.

Notice also that, in view of the integrated local energy decay estimate
established in \cite{Laul2015}, the results of the present paper
also imply a $\bar{t}^{-2+\text{\textgreek{d}}}$ decay estimate for
solutions $\text{\textgreek{f}}$ to (\ref{eq:WaveEquation}) on very
slowly rotating $4+1$ dimensional Myers--Perry spacetimes.

For our second example, we will first need to introduce a definition:
A metric $g$ on $\mathbb{R}^{d+1}$ will be called a \emph{radiating
uniformly small perturbation} of Minkowski spacetime $(\mathbb{R}^{d+1},\text{\textgreek{h}})$
if it has the asymptotics (\ref{eq:RoughAsymptotics}), and moreover
there exists a small $\text{\textgreek{e}}_{0}>0$ such that $r\cdot(g-\text{\textgreek{h}})$
and all its derivatives are $\text{\textgreek{e}}_{0}$-globally small,
with each differentiation of this tensor with respect to $\partial_{t}$
\underline{except} for the first one yielding additional decay in
terms of $|u|$ (see (\ref{eq:Uniform boundedness}) and (\ref{eq:DeformationTensorTAwayDecay-1-1-1})
for a more precise definition). For such spacetimes, the geometric
assumptions of Parts 1, 2 and 3 are satisfied and an integrated loacal
energy decay estimate of the form (\ref{eq:ModelILED}) without loss
of derivatives holds (in view of the stability to small perturbations
of the estimates provided by the $\partial_{r}$-Morawetz current,
combined with the estimates of Section \ref{sec:Morawetz} of the
present paper). Examples of such spacetimes include the vacuum dynamical
perturbations of Minkowski spacetime considered in \cite{Christodoulou1993}.

We will infer the following result:

\begin{customthm}{(Improved polynomial decay on radiating uniformly small perturbations of Minkowski)}
If $(\mathbb{R}^{d+1},g)$ is a radiating uniformly small perturbation
of Minkowski spacetime and $\text{\textgreek{e}}_{0}$ is small enough,
then any solution $\text{\textgreek{f}}$ to $\square_{g}\text{\textgreek{f}}=0$
on $(\mathbb{R}^{d+1},g)$ with suitably decaying initial data on
$\{t=0\}$ will satisfy a $\bar{t}^{-\frac{d}{2}+O(\text{\textgreek{e}}_{0})}$
decay estimate. If, in addition, the deformation tensor of the vector
field $\partial_{t}$ is $O(\bar{t}^{-\text{\textgreek{d}}_{0}})$
decaying for some $\text{\textgreek{d}}_{0}$, then $\text{\textgreek{f}}$
will satisfy a $\bar{t}^{-\frac{d}{2}}$ decay rate. \end{customthm}

See Section \ref{sub:PerturbationsMinkowski} for a precise statement
of this result. Let us remark that this theorem extends a recent result
of Oliver \cite{Oliver2014}. 

Our final example will concern the class of radiating black hole exterior
spacetimes $(\mathcal{M},g)$ dynamically settling down to the exterior
region of a subextremal Kerr spacetime. In order to present our example
in the most simple form that can be deduced without computation from
previous results, we will retrict ourselves to spacetimes $(\mathcal{M},g)$
settling down to Schwarzschild exterior at a sufficiently fast polynomial
rate. This class includes the dynamical vacuum spacetimes constructed
in \cite{Dafermos2013} (which actually approach Schwarzschild at
an exponential rate). 

The energy current yielding the integrated local energy decay statement
for Schwarzschild exterior constructed in \cite{DafRod4}, combined
with the estimates of Section \ref{sec:Morawetz} of the present paper
and the fast rate at which $g$ approaches the Schwarzschild metric
$g_{M}$, immediately imply that an integrated local energy decay
statement of the form (\ref{eq:ModelILED}) also holds on $(\mathcal{M},g)$.
Furthermore, it is straightforward to check that $(\mathcal{M},g)$
satisfies the assumptions of Parts 1, 2 and 3 of the above Theorem
(in view of the fast approach to the Schwarzschild exterior metric,
which satisfies these assumptions). Thus, on these spacetimes we will
be able to infer the following result:

\begin{customthm}{(Improved polynomial decay on dynamical, radiating black hole spacetimes)}
If $(\mathcal{M}^{3+1},g)$ is a radiating black hole spacetime settling
down to a Schwarzschild exterior at a sufficiently fast polynomial
decay rate (such us the ones constructed in \cite{Dafermos2013}),
then any solution $\text{\textgreek{f}}$ to $\square_{g}\text{\textgreek{f}}=0$
on $(\mathcal{M},g)$ with suitably decaying initial data on a Cauchy
hypersurface will satisfy a $\bar{t}^{-\frac{3}{2}}$ decay estimate.\end{customthm}

We will discuss in more detail the results of this paper and their
applications in the next sections of the introduction. But first,
we will briefly review the ``old'' approach of using the conformal
isometries of Minkowski spacetime for establishing decay rates for
solutions to (\ref{eq:WaveEquation}) on asymptotically flat spacetimes,
and compare it to the method of \cite{DafRod7}.

\subsection{Comparison of the two approaches}

\subsubsection{The ``old'' approach}

The use of first order operators as multipliers and commutators for
(\ref{eq:WaveEquation}) has been implemented extensively during the
last 50 years to deal with linear and non linear wave equations on
small perturbations of Minkowski spacetime $(\mathbb{R}^{d+1},\text{\textgreek{h}})$. 

Following Morawetz (see e.\,g.~\cite{Morawetz1962}), one way to
obtain decay for the local energy of solutions $\text{\textgreek{f}}$
to (\ref{eq:WaveEquation}) is to apply the conformal Killing field
$Z$ of ($\mathbb{R}^{d+1},\text{\textgreek{h}})$ 
\begin{equation}
Z=(t^{2}+r^{2})\partial_{t}+2tr\partial_{r}\label{eq:MorawetzKVectorField}
\end{equation}
as multiplier for (\ref{eq:WaveEquation}). On Minkowski spacetime
itself for $d\ge3$, the vector field $Z$ gives rise to a conserved
positive definite energy norm $E_{Z}[\text{\textgreek{f}}](t)$ with
weights growing in $t$. In particular, one can bound: 
\begin{equation}
E_{Z}[\text{\textgreek{f}}](\text{\textgreek{t}})\gtrsim\int_{\{t=\text{\textgreek{t}}\}}\Big(v^{2}|\partial_{v}\text{\textgreek{f}}|^{2}+u^{2}|\partial_{u}\text{\textgreek{f}}|^{2}+|\partial_{\text{\textgreek{sv}}}\text{\textgreek{f}}|^{2}+|\text{\textgreek{f}}|^{2}\Big)\, dx,\label{eq:BoundKEnergy}
\end{equation}
where $(t,r,\text{\textgreek{sv}})$ is the usual polar coordinate
system on Minkowski space $\mathbb{R}^{3+1}$, $v=t+r$, $u=t-r$
and $dx$ denotes the usual integration measure on $\{t=const\}$
slices of $\mathbb{R}^{d+1}$ (see also Section \ref{sub:Notational-conventions}
for the $\text{\textgreek{sv}}$ notation). Thus, the preservation
of $E_{Z}[\text{\textgreek{f}}](t)$ and the growth in time of the
weights in the expression (\ref{eq:BoundKEnergy}) can be used to
establish polynomial decay in time estimates for the $L^{2}$ norm
of certain derivatives of $\text{\textgreek{f}}$.

The above approach has been also implemented in the treatment of the
wave equation (\ref{eq:WaveEquation}) on the complement of a compact
obstacle $\mathcal{O}$ in flat space, with suitable boundary conditions
imposed on the boundary of $\mathcal{O}$. In \cite{Morawetz1961,Morawetz1962},
for instance, pointwise polynomial decay rates were established for
solutions $\text{\textgreek{f}}$ to (\ref{eq:WaveEquation}) on the
complement of a star shaped obstacle with Dirichlet boundary conditions,
and this was achieved with the use of the conformally Killing vector
field $Z$ and the radial vector field $\partial_{r}$ as multipliers
for equation (\ref{eq:WaveEquation}). Moreover, the use of $\partial_{r}$
as a multiplier for (\ref{eq:WaveEquation}) yielded the integrated
local energy decay statement (\ref{eq:ILED Morawetz}).

Another method for obtaining refined pointwise decay rates for solutions
$\text{\textgreek{f}}$ to (\ref{eq:WaveEquation}) on flat spacetime
is the commutation vector field method, introduced by Klainerman:
By commuting equation (\ref{eq:WaveEquation}) with the generators
of the isometries of $(\mathbb{R}^{d+1},\text{\textgreek{h}})$ plus
the dilation vector fieldand the dilation vector field $S$: 
\begin{equation}
S=t\partial_{t}+r\partial_{r}\label{eq:SVectorFieldMinkowski}
\end{equation}
 ($(t,r)$ being the usual time and radius coordinates on Minkowski
space), and using the conservation of the $E_{Z}$ energy norm (\ref{eq:BoundKEnergy})
on $(\mathbb{R}^{d+1},\text{\textgreek{h}})$ together with a modified
version of the Sobolev embedding theorem due to Klainerman (see \cite{Klainerman1985,Klainerman1987}),
one can attain a pointwise decay estimate for $\text{\textgreek{f}}$
($d\ge3)$: 
\begin{equation}
|\partial_{u}^{l_{1}}\partial_{v}^{l_{2}}(r^{-1}\partial_{\text{\textgreek{sv}}})^{l_{3}}\text{\textgreek{f}}|\lesssim(1+|t-r|)^{-\frac{1}{2}-l_{1}}(1+t+r)^{-\frac{d-1}{2}-l_{2}-l_{3}}\sqrt{\mathcal{E}_{l_{1},l_{2},l_{3}}},\label{eq:BestDecaySobolev}
\end{equation}
 where $l_{1},l_{2},l_{3}\ge0$ are integers and $\mathcal{E}_{l_{1},l_{2},l_{3}}$
is a weighted higher order energy norm of the initial data for $\text{\textgreek{f}}$
on $\{t=0\}$. See \cite{Sogge2008} for more details on the commutation
vector field approach. 

Notice that the $t^{-\frac{d}{2}}$ decay rate for $|\text{\textgreek{f}}|$
in the region $\{r\lesssim1\}$ provided by (\ref{eq:BestDecaySobolev})
guarantees that $|\text{\textgreek{f}}(t,x)|$ is integrable in $t$
in dimensions $d\ge3$, and this fact is of fundamental importance
in the treatment of non linear variants of the wave equation (\ref{eq:WaveEquation}). 

The aforementioned techniques have been also extended to the exterior
of black hole spacetimes, such as the Schwarzschild and very slowly
rotating (i.\,e~with $|a|\ll M$) Kerr exterior spacetimes, see
\cite{DafRod2,DafRod4,AndBlue1,BlueSterb}. In these works, a variant
of the conformally Killing vector field $Z$ of Minkowski spacetime
was constructed and used, but this construction came at a cost: Since
$Z$ is not a conformally Killing vector field on these black hole
spacetimes, decay estimates obtained in this way for solutions to
(\ref{eq:WaveEquation}) were coupled with error terms in the near
region of the spacetimes under consideration, and these error terms
carried weights growing in time. 

In view also of the unavoidable presence of trapping in the near region
of a black hole spacetime, the error terms associated to the use of
the modified $Z$ vector field as a multiplier for (\ref{eq:WaveEquation})
required additional effort in order to be controlled. An essential
step towards controlling these error terms was the establishment of
an integrated local energy decay statement of the form (\ref{eq:ModelILED}),
with the use of carefully chosen first order multipliers for (\ref{eq:WaveEquation})
capturing the red-shift effect near the horizon $\mathcal{H}$ and
the structure of the trapped set in the near region $\{r\lesssim1\}$
(these multipliers being equal to $\partial_{r}$ plus a lower order
correction in the far away region $\{r\gg1\}$). See \cite{BlueSof1,DafRod2,DafRod4,DafRod6,DafRod9,AndBlue1,TatToh1}. 

The above approach of using an adaptation of the Morawetz $Z$ vector
field and an integrated local energy decay statement yielded $t^{-1}$
decay estimates for solutions $\text{\textgreek{f}}$ to (\ref{eq:WaveEquation})
on Schwarzschild exterior spacetimes and $t^{-1+\text{\textgreek{d}}(a)}$
decay estimates on slowly rotating Kerr exterior spacetimes, with
$\text{\textgreek{d}}(a)\rightarrow0$ as $a\rightarrow0$ (see \cite{DafRod6}).
In \cite{Luk2010b,Luk2}, Luk was able to obtain improved $t^{-\frac{3}{2}+\text{\textgreek{d}}}$
decay estimates for $\text{\textgreek{f}}$ in the near region of
these backgrounds by commuting the wave equation (\ref{eq:WaveEquation})
with an analogue of the dilation vector field $S$ (\ref{eq:SVectorFieldMinkowski})
of Minkowski spacetime.

Let us note at this point that the vector field approach has been
effectively applied in the case of non linear wave equations on a
radiating spacetime which is globally close to $(\mathbb{R}^{3+1},\text{\textgreek{h}})$:
This can be viewed as a corollary of the monumental proof of the non
linear stability of Minkowski spacetime in the context of the Einstein
equations, by Christodoulou and Klainerman (see \cite{Christodoulou1993}).
These techniques have also been applied in the study of non linear
wave equations on black hole spacetimes (see the work of Luk \cite{Luk2010a}).
See also \cite{Oliver2014} for the treatment of the linear wave equation
(\ref{eq:WaveEquation}) on radiating spacetimes which are globally
close to $(\mathbb{R}^{3+1},\text{\textgreek{h}})$ (where, among
other decay results, a $t^{-\frac{3}{2}}$ decay rate in the near
region is established).

The difficulties in extending the ``old'' approach of establishing
decay estimates for solutions to (\ref{eq:WaveEquation}) on more
general black hole spacetimes led the authors of \cite{DafRod7} to
suggest a more flexible approach that does not involve multipliers
and commutators with weights growing in time. This is the approach
that we will now discuss.

\subsubsection{The $r^{p}$-weighted energy method}

The crux of the new method of obtaining decay estimates for solutions
to (\ref{eq:WaveEquation}) introduced in \cite{DafRod7} lies in
the establishment of a hierarchy of estimates for $r^{p}$-weighted
energies, $0\le p\le2$, using as a multiplier an $r^{p}$-weighted
outgoing null vector field. On Minkowski spacetime, this hierarchy
of estimates takes the following form for any solution $\text{\textgreek{f}}$
to the wave equation (\ref{eq:WaveEquation}) and any $\text{\textgreek{t}}_{1}\le\text{\textgreek{t}}_{2}$,
$R>0$: 
\begin{multline}
\int_{\{u=\text{\textgreek{t}}_{2}\}\cap\{r\ge R\}}r^{p}\cdot|\partial_{v}(r\text{\textgreek{f}})|^{2}\, dvd\text{\textgreek{sv}}+\int_{\{\text{\textgreek{t}}_{1}\le u\le\text{\textgreek{t}}_{2}\}\cap\{r\ge R\}}r^{p-1}\Big(p|\partial_{v}(r\text{\textgreek{f}})|^{2}+(2-p)|r^{-1}\partial_{\text{\textgreek{sv}}}(r\text{\textgreek{f}})|^{2}\Big)\, dudvd\text{\textgreek{sv}}\le\\
\le\int_{\{u=\text{\textgreek{t}}_{2}\}\cap\{r\ge R\}}r^{p}\cdot|\partial_{v}(r\text{\textgreek{f}})|^{2}\, dvd\text{\textgreek{sv}}+\int_{\{\text{\textgreek{t}}_{1}\le u\le\text{\textgreek{t}}_{2}\}\cap\{r=R\}}r^{p}\cdot\Big(|r^{-1}\partial_{\text{\textgreek{sv}}}(r\text{\textgreek{f}})|^{2}-|\partial_{v}(r\text{\textgreek{f}})|^{2}\Big)\, dud\text{\textgreek{sv}}.\label{eq:PHierarchy}
\end{multline}
In the above, $u=t-r$, $v=t+r$. For the $\text{\textgreek{sv}}$
notation on the angular variables, see Section \ref{sub:Notational-conventions}.
Moreover, the right hand side of (\ref{eq:PHierarchy}) also controls
the angular derivatives of the radiation field of $\text{\textgreek{f}}$
on future null infinity $\mathcal{I}^{+}$, but we have dropped these
terms for simplicity. A similar expression is also valid on Schwarzschild
spacetimes and suitable perturbations, see \cite{DafRod7}.

The importance of the hierarchy (\ref{eq:PHierarchy}) lies in the
fact that the left hand side of (\ref{eq:PHierarchy}) contains a
positive definite bulk term, while the ``error'' term in the near
region (namely the last term of the right hand side) does not carry
weights growing in $t$. Thus, combining (\ref{eq:PHierarchy}) with
the integrated local energy decay statement (\ref{eq:ModelILED})
and the energy boundedness estimate (\ref{eq:ModelBoundedness}),
the authors of \cite{DafRod7} were able to obtain uniform polynomial
decay rates for $r^{\frac{1}{2}}\text{\textgreek{f}}$ and $r\text{\textgreek{f}}$
in terms of $u$.

A noteable aspect of this novel approach of \cite{DafRod7} is that
decay rates for $\text{\textgreek{f}}$ are obtained by repeatedly
applying the pigeonhole principle on the positive definite bulk term
(i.\,e.~the second term of the left hand side) controlled in (\ref{eq:PHierarchy}).
This is in contrast to the older approach (described in the previous
section), which yielded decay rates for $\text{\textgreek{f}}$ by
establishing uniform bounds for $t$-weighted energy norms of $\text{\textgreek{f}}$
on suitable hypersurfaces.

Moreover, the new method of \cite{DafRod7} allows one to obtain the
result of Luk (\cite{Luk2}), namely to establish improved (i.\,e.
$t^{-\frac{3}{2}+\text{\textgreek{d}}}$) polynomial decay estimates
for $\text{\textgreek{f}}$ in the near region of Schwarzschild exterior:
This was achieved by Schlue in \cite{Schlue2013}, where it was established
that commuting (\ref{eq:WaveEquation}) with the outgoing null vector
field $\partial_{v}$ (as well as the generators of the isometries
of Schwarzschild) leads to an improvement of the $p$-hierarchy (\ref{eq:PHierarchy}).%
\footnote{Note also that \cite{Schlue2013} also deals with the case of higher
dimensional Schwarzschild spacetimes.%
} In particular, it was established that higher order $\partial_{v}$
and $r^{-1}\partial_{\text{\textgreek{sv}}}$ derivatives of $\text{\textgreek{f}}$
satisfy (\ref{eq:PHierarchy}) for larger values of $p$. This better
decay rate in $r$ was then translated into a better decay rate in
$u$ by using the expression of the wave equation (\ref{eq:WaveEquation})
as well as the pigeonhole principle argument of \cite{DafRod7}. 

A variant of this novel approach has also been effectively used for
the case of the Einstein equations themselves: In \cite{Holzegel2010,Dafermos2013},
the authors establish an $r^{p}$-weighted energy hierarchy, a proper
tensorial analogue of (\ref{eq:PHierarchy}), for radiating solutions
to the Einstein equations $Ric(g)=0$ that approach the Schwarzschild
exterior in the future.

\subsection{Statement of the main results}

The present paper introduces a broad class of asymptotically flat
spacetimes $(\mathcal{M}^{d+1},g)$, $d\ge3$, on which the techniques
of \cite{DafRod7} and \cite{Schlue2013} can be generalised. See
Section \ref{sec:GeometrySpacetimes} for a more detailed discussion
of the class of spacetimes under consideration. Notice that this class
of metrics includes spacetimes with non-constant Bondi mass at null
infinity, see e.\,g.~\cite{Bondi1962,Sachs1962}, such as the dynamical
vacuum perturbations of Minkowski spacetime (see \cite{Christodoulou1993}).
Moreover, spacetimes in this class are allowed to have a timelike
boundary $\partial_{tim}\mathcal{M}$ with compact spacelike cross-sections. 

We will now proceed to briefly review the results established in the
following sections of the paper.

\subsubsection{\label{sub:IntroductionHierarchy}The $r^{p}$-weighted energy hierarchy
in dimensions $d\ge3$}

In this Section, all results will be stated on an asymptotic region
$\mathcal{N}_{af}\subset\mathcal{M}$ of a general radiating asymptotically
flat%
\footnote{Let us also remark that the results of this paper also hold on spacetimes
which are asymptotically conic instead of asymptotically flat. However,
we will not pursue this issue further in this paper.%
} spacetime $(\mathcal{M},g)$. In particular, let $(\mathcal{N}_{af}^{d+1},g)$,
$d\ge3$, be a Lorentzian manifold diffeomorphic to $\mathbb{R}\times[R_{0},+\infty)\times\mathbb{S}^{d-1}$
for some $R_{0}>0$, on which a single $(u,r,\text{\textgreek{sv}})$
coordinate chart has been fixed. Assume that in this chart $g$ takes
the form 
\begin{align}
g= & -4\Big(1-\frac{2M(u,\text{\textgreek{sv}})}{r}+O(r^{-1-a})\Big)du^{2}-\Big(4+O(r^{-1-a})\Big)dudr+r^{2}\cdot\Big(g_{\mathbb{S}^{d-1}}+O(r^{-1})\Big)+O(1)dud\text{\textgreek{sv}}+\label{eq:RoughAsymptotics}\\
 & +O(r^{-a})drd\text{\textgreek{sv}}+O(r^{-2-a})dr^{2},\nonumber 
\end{align}
where $M(u,\text{\textgreek{sv}})$ is a bounded and sufficiently
regular function of $u,\text{\textgreek{sv}}$. Notice that $(\mathcal{N}_{af},g)$
is not in general globally hyperbolic.

We will extend the hierarchy (\ref{eq:PHierarchy}) to $(\mathcal{N}_{af},g)$
as follows (see the remark below for explanation of the notation):
\begin{thm}
\label{thm:NewMethodModelIntro}Let $\mathcal{S}_{1},\mathcal{S}_{2}$
be two spacelike hyperboloidal hypersurfaces of $(\mathcal{N}_{af}^{d+1},g)$,
$d\ge3$, terminating at $\mathcal{I}^{+}$, such that $\mathcal{S}_{2}\subset J^{+}(\mathcal{S}_{1})$.
Then for any $0<p\le2$, any given $0<\text{\textgreek{h}}<a$, $0<\text{\textgreek{d}}<1$
and $R>0$ large, the following inequality is true for any smooth
function $\text{\textgreek{f}}:\mathcal{N}_{af}\rightarrow\mathbb{C}$
(setting also $\text{\textgreek{F}}\doteq\text{\textgreek{W}}\cdot\text{\textgreek{f}}$,
where $\text{\textgreek{W}}=(-\det(g))^{\frac{1}{4}}$ ):
\begin{equation}
\begin{split}\int_{\mathcal{S}_{2}\cap\{r\gtrsim R\}}r^{p}|\partial_{r}\text{\textgreek{F}}|^{2}\, dv & d\text{\textgreek{sv}}+\int_{\mathcal{S}_{2}\cap\{r\gtrsim R\}}\Big(r^{p}|r^{-1}\partial_{\text{\textgreek{sv}}}\text{\textgreek{F}}|^{2}+\max\big\{(d-3),r^{-\text{\textgreek{d}}}\big\}\cdot r^{p-2}|\text{\textgreek{F}}|^{2}\Big)\, dud\text{\textgreek{sv}}+\\
+\int_{J^{+}(\mathcal{S}_{1})\cap J^{-}(\mathcal{S}_{2})\cap\{r\gtrsim R\}} & \text{\textgreek{q}}_{R}\cdot\Big(pr^{p-1}|\partial_{r}\text{\textgreek{F}}|^{2}+(2-p)r^{p-1}|r^{-1}\partial_{\text{\textgreek{sv}}}\text{\textgreek{F}}|^{2}+\max\big\{(2-p)(d-3),r^{-\text{\textgreek{d}}}\big\}\cdot r^{p-3}|\text{\textgreek{F}}|^{2}\Big)\, dudvd\text{\textgreek{sv}}\lesssim_{p,\text{\textgreek{h},\textgreek{d}}}\\
\lesssim_{p,\text{\textgreek{h},\textgreek{d}}} & \int_{\mathcal{S}_{1}\cap\{r\gtrsim R\}}r^{p}|\partial_{r}\text{\textgreek{F}}|^{2}\, dvd\text{\textgreek{sv}}+\int_{\mathcal{S}_{2}\cap\{r\gtrsim R\}}\Big(r^{p}|r^{-1}\partial_{\text{\textgreek{sv}}}\text{\textgreek{F}}|^{2}+\max\big\{(d-3),r^{-\text{\textgreek{d}}}\big\}\cdot r^{p-2}|\text{\textgreek{F}}|^{2}\Big)\, dud\text{\textgreek{sv}}+\\
 & +\int_{J^{+}(\mathcal{S}_{1})\cap J^{-}(\mathcal{S}_{2})\cap\{r\sim R\}}\big(r^{p}|\partial\text{\textgreek{F}}|^{2}+r^{p-2}|\text{\textgreek{F}}|^{2}\big)\, dudvd\text{\textgreek{sv}}+\int_{\mathcal{S}_{1}\cap\{r\gtrsim R\}}J_{\text{\textgreek{m}}}^{\partial_{u}}(\text{\textgreek{f}})n_{\mathcal{S}}^{\text{\textgreek{m}}}+\\
 & +\int_{J^{+}(\mathcal{S}_{1})\cap J^{-}(\mathcal{S}_{2})\cap\{r\gtrsim R\}}\text{\textgreek{q}}_{R}\cdot(r^{p+1}+r^{1+\text{\textgreek{h}}})\cdot|\text{\textgreek{W}}\cdot\square_{g}\text{\textgreek{f}}|^{2}\, dudvd\text{\textgreek{sv}}.
\end{split}
\label{eq:newMethodInitialStatementHyperboloids}
\end{equation}
 In the above, the constants implicit in the $\lesssim_{p,\text{\textgreek{h},\textgreek{d}}}$
notation depend only on $p,\text{\textgreek{h}},\text{\textgreek{d}}$
and on the geometry of $(\mathcal{N}_{af},g)$. The partial derivatives
$\partial_{r},\partial_{\text{\textgreek{sv}}}$ are considered with
respect to the cooordinate chart $(u,r,\text{\textgreek{sv}})$ and
the notation $\partial_{\text{\textgreek{sv}}}$ is explained in Section
\ref{sub:Notational-conventions}.\end{thm}
\begin{rem*}
Notice that the $dud\text{\textgreek{sv}}$ volume form on the hyperboloidal
hypersurfaces $\mathcal{S}_{i}$ degenerates as $r\rightarrow+\infty$
when compared to $dvd\text{\textgreek{sv}}$. The notion of a spacelike
hyperboloidal hypersurface terminating at $\mathcal{I}^{+}$ is given
in Section \ref{sub:Spacelike-hyperboloidal-hypersurfaces}. For the
notations on vector field currents see Section \ref{sub:Notational-conventions}. 
\end{rem*}
For a more detailed statement of the above result, see Theorems \ref{thm:NewMethodFinalStatementHyperboloids}
and \ref{thm:NewMethodFinalStatement} in Section \ref{sec:The-new-method}. 

We will also establish the following improved $r^{p}$-weighted hierarchy
for higher derivatives of $\text{\textgreek{f}}$:
\begin{thm}
\label{thm:ImprovedNewMethodIntro}With the notations as in Theorem
\ref{thm:NewMethodModelIntro}, for any $k\in\mathbb{N}$, any $2k-2<p\le2k$,
any given $0<\text{\textgreek{h}}<a$ and $R>0$ large, the following
inequality is true for any smooth function $\text{\textgreek{f}}:\mathcal{N}_{af}\rightarrow\mathbb{C}$
(setting also $\text{\textgreek{F}}\doteq\text{\textgreek{W}}\cdot\text{\textgreek{f}}$
): 
\begin{equation}
\begin{split}\int_{\mathcal{S}_{2}\cap\{r\gtrsim R\}}r^{p}|\partial_{r}\nabla_{\mathcal{S}}^{k-1}\text{\textgreek{F}}|^{2}\, dvd\text{\textgreek{sv}}+\int_{\mathcal{S}_{2}\cap\{r\gtrsim R\}} & r^{p}|r^{-1}\partial_{\text{\textgreek{sv}}}\nabla_{\mathcal{S}}^{k-1}\text{\textgreek{F}}|^{2}\, dud\text{\textgreek{sv}}+\\
+\int_{J^{+}(\mathcal{S}_{1})\cap J^{-}(\mathcal{S}_{2})\cap\{r\gtrsim R\}}\text{\textgreek{q}}_{R}\cdot\Big(r^{p-1}|\partial_{r} & \nabla_{\mathcal{S}}^{k-1}\text{\textgreek{F}}|^{2}+(2k-p)r^{p-1}|r^{-1}\partial_{\text{\textgreek{sv}}}\nabla_{\mathcal{S}}^{k-1}\text{\textgreek{F}}|^{2}\Big)\, dudvd\text{\textgreek{sv}}\lesssim_{p,\text{\textgreek{h}}}\\
\lesssim_{p,\text{\textgreek{h}}}\sum_{j=1}^{k}\int_{\mathcal{S}_{1}\cap\{r\gtrsim R\}}r^{p-2(k-j)}|\partial_{r} & \nabla_{\mathcal{S}}^{j-1}\text{\textgreek{F}}|^{2}\, dvd\text{\textgreek{sv}}+\sum_{j=1}^{k}\int_{\mathcal{S}_{1}\cap\{r\gtrsim R\}}\Big(r^{p-2(k-j)}|r^{-1}\partial_{\text{\textgreek{sv}}}\nabla_{\mathcal{S}}^{j-1}\text{\textgreek{F}}|^{2}+r^{p-2k}|\text{\textgreek{F}}|^{2}\Big)\, dud\text{\textgreek{sv}}+\\
+\sum_{j=0}^{k}\int_{J^{+}(\mathcal{S}_{1})\cap J^{-}(\mathcal{S}_{2})\cap\{r\sim R\}} & r^{p-2(k-j)}|\partial^{j}\text{\textgreek{F}}|^{2}\, dudvd\text{\textgreek{sv}}+\sum_{j=1}^{k}\int_{\mathcal{S}_{1}\cap\{r\gtrsim R\}}J_{\text{\textgreek{m}}}^{\partial_{u}}(\nabla_{\mathcal{S}}^{j-1}\text{\textgreek{f}})n_{\mathcal{S}}^{\text{\textgreek{m}}}+\\
+\sum_{j=1}^{k}\int_{J^{+}(\mathcal{S}_{1})\cap J^{-}(\mathcal{S}_{2})\cap\{r\gtrsim R\}} & \text{\textgreek{q}}_{R}\cdot(r^{p+1-2(k-j)}+r^{1+\text{\textgreek{h}}})\cdot|\nabla_{\mathcal{S}}^{j-1}(\text{\textgreek{W}}\square_{g}\text{\textgreek{f}})|^{2}\, dudvd\text{\textgreek{sv}}.
\end{split}
\label{eq:newMethodInitialStatementHyperboloids-1}
\end{equation}
 In the above, 
\begin{equation}
|\nabla_{S}^{j}\text{\textgreek{y}}|^{2}\doteq\sum_{j_{1}+j_{2}+j_{3}=j}|r^{-j_{2}-j_{3}}\partial_{r}^{j_{1}}\partial_{\text{\textgreek{sv}}}^{j_{2}}\partial_{u}^{j_{3}}\text{\textgreek{y}}|^{2}
\end{equation}
 and the constants implicit in the $\lesssim_{p,\text{\textgreek{h}}}$
notation depend only on $p,\text{\textgreek{h}}$ and on the geometry
of $(\mathcal{N}_{af},g)$. 
\end{thm}
See Section \ref{sec:The-improved--hierarchy} for a more detailed
statement of the above result.

\subsubsection{\label{sub:IntroductionFirstPolynomialDecay}A $\bar{t}^{-1}$ polynomial
decay estimate for solutions to the wave equation $\square\text{\textgreek{f}}=0$ }

In this section, we will be concerned with obtaining results for (\ref{eq:WaveEquation})
on the whole spacetime $(\mathcal{M},g)$ (and not merely the asymptotic
region $\mathcal{N}_{af}$). Provided that an integrated local energy
decay statement (possibly with loss of derivatives) holds for solutions
to $\square_{g}\text{\textgreek{f}}=F$ on a spacetime $(\mathcal{M}^{d+1},g)$,
$d\ge3$, with $g$ asymptotically of the form (\ref{eq:RoughAsymptotics}),
we will establish polynomial decay rates for $\text{\textgreek{f}}$
with respect to a hyperboloidal foliation of $\mathcal{M}$.

In particular, let $(\mathcal{M}^{d+1},g)$, $d\ge3$, be a Lorentzian
manifold with possibly non-empty boundary $\partial\mathcal{M}$,
which can be split as 
\begin{equation}
\partial\mathcal{M}=\partial_{tim}\mathcal{M}\oplus\partial_{hor}\mathcal{M},
\end{equation}
where $\partial_{tim}\mathcal{M}$ is smooth and timelike and $\partial_{hor}\mathcal{M}$
is piecewise smooth and null. Assume also that $(\mathcal{M},g)$
is globally hyperbolic as a manifold with timelike boundary, which
means that the double $(\tilde{\mathcal{M}}_{tim},g)$ of $(\mathcal{M},g)$
across $\partial_{tim}\mathcal{M}$ is globally hyperbolic. 

Suppose that $(\mathcal{M},g)$ satisfies the following geometric
assumptions:

\begin{MyDescription}[leftmargin=3.0em]

\item [(GM1)]{ \label{enu:Asymptotic-flatness}\emph{Asymptotic flatness:}
$(\mathcal{M},g)$ is asymptotically flat in the sense that there
exists an open subset $\mathcal{N}_{af}\subset\mathcal{M}$ such that
each connected component of $\mathcal{N}_{af}$ is mapped diffeomorphically
on $\mathbb{R}\times(R_{0},+\infty)\times\mathbb{S}^{d-1}$ through
a coordinate chart $(u,r,\text{\textgreek{sv}})$, and in this coordinate
chart $g$ has the form (\ref{eq:RoughAsymptotics}).}

\item [(GM2)]{\label{enu:NearRegion}\emph{Existence of a well behaved
time function: }There exists a function $\bar{t}:\mathcal{M}\rightarrow\mathbb{R}$
with level sets which are spacelike hyperboloids terminating at future
null infinity, such that on each component of $\mathcal{N}_{af}$
the difference $|\bar{t}-u|$ is bounded, and the foliation $\{\bar{t}=const\}$
in the region $\mathcal{M}\backslash\mathcal{N}_{af}$ is sufficiently
``regular'' (see Section \ref{sec:Firstdecay} for the precise relevant
assumptions on $\bar{t}$). }

\end{MyDescription}

\noindent In fact, the precise description of Assumptions \ref{enu:Asymptotic-flatness}--\ref{enu:NearRegion}
is more complicated, and requires the splitting of these assumptions
into a larger number of statements: see Assumptions \ref{enu:AsymptoticFlatness}--\ref{enu:SobolevInequality}
of Sections \ref{sub:GeometricAssumptionsFriedlander} and \ref{sub:GeometricAssumptionsAndConstructions}. 

For convenience, we define the globally timelike vector field $N$
so that $N\equiv grad(\bar{t})$ on $\mathcal{M}\backslash\mathcal{N}_{af}$
and $N=\partial_{u}$ on each connected component of the region $\mathcal{N}_{af}\cap\{r\gg1\}$.

Suppose also that on $(\mathcal{M},g)$ the following integrated local
energy decay estimate holds: 

\begin{MyDescription}[leftmargin=3.0em]

\item [(ILED1)]{\label{enu:Integrated-local-energy}\emph{Integrated
local energy decay with polynomial loss of derivatives:} There exists
an integer $k\ge0$, such that for any solution $\text{\textgreek{f}}$
to $\square\text{\textgreek{f}}=F$ with suitable boundary conditions
on $\partial_{tim}\mathcal{M}$ , any $m\in\mathbb{N}$, $0\le\text{\textgreek{t}}_{1}\le\text{\textgreek{t}}_{2}$,
$\text{\textgreek{h}}>0$ and $R>0$:
\begin{equation}
\begin{split}\sum_{j=0}^{m}\int_{\{\text{\textgreek{t}}_{1}\le\bar{t}\le\text{\textgreek{t}}_{2}\}\cap\{r\le R\}}|\nabla^{j}\text{\textgreek{f}}|^{2} & +\sum_{j=1}^{m}\int_{\{\text{\textgreek{t}}_{1}\le\bar{t}\le\text{\textgreek{t}}_{2}\}\cap\partial_{tim}\mathcal{M}}|\nabla^{j}\text{\textgreek{f}}|^{2}\le\\
\le & C_{m,\text{\textgreek{h}}}(R)\sum_{j=0}^{m+k-1}\int_{\{\bar{t}=\text{\textgreek{t}}_{1}\}}J_{\text{\textgreek{m}}}^{N}(N^{j}\text{\textgreek{f}})\bar{n}^{\text{\textgreek{m}}}+C_{m,\text{\textgreek{h}}}\sum_{j=0}^{m-1}\int_{\{\text{\textgreek{t}}_{1}\le\bar{t}\le\text{\textgreek{t}}_{2}\}\cap\{r\ge R\}}r^{-1}J_{\text{\textgreek{m}}}^{N}(N^{j}\text{\textgreek{f}})\bar{n}^{\text{\textgreek{m}}}+\\
 & +C_{m,\text{\textgreek{h}}}(R)\sum_{j=0}^{m+k-1}\int_{\{\text{\textgreek{t}}_{1}\le\bar{t}\le\text{\textgreek{t}}_{2}\}}r^{1+\text{\textgreek{h}}}|\nabla^{j}F|^{2},
\end{split}
\label{eq:ILEDwithloss-1}
\end{equation}
 where $C_{m,\text{\textgreek{h}}}(R)$ depends only on $m$,$\text{\textgreek{h}}$,$R$
and the geometry of $(\mathcal{M},g)$.}

\end{MyDescription}

See Section \ref{sub:IntegratedLocalEnergyDecay} for a more detailed
description of Assumption \ref{enu:Integrated-local-energy}. For
an alternative to the Assumption \ref{enu:Integrated-local-energy},
see the remarks in Section \ref{sub:IntegratedLocalEnergyDecay} (and
in particular (\ref{eq:IntegratedLocalEnergyDecayImprovedDecayAltLoss})
and (\ref{eq:IntegratedLocalEnergyDecayImprovedDecayAlt})).

On any spacetime $(\mathcal{M},g)$ satisfying the above assumptions
we will establish the following decay statement for solutions to (\ref{eq:WaveEquation}):
\begin{thm}
\label{thm:FirstDecayIntroduction}Let $(\mathcal{M}^{d+1},g)$ satisfy
the geometric assumptions \ref{enu:Asymptotic-flatness} and \ref{enu:NearRegion},
and the integrated local energy decay assumption \ref{enu:Integrated-local-energy},
and let $\bar{t}$ and $N$ be as above. Then the following decay
estimates hold for any $0<\text{\textgreek{d}}<1$, any $\text{\textgreek{t}}\ge0$
and any solution $\text{\textgreek{f}}$ to the inhomogeneous wave
equation $\square_{g}\text{\textgreek{f}}=F$ with suitably decaying
inital data on $\{\bar{t}=0\}$ (and satisfying suitable boundary
conditions on $\partial_{tim}\mathcal{M}$):

\begin{equation}
\int_{\{\bar{t}=\text{\textgreek{t}}\}}J_{\text{\textgreek{m}}}^{N}(\text{\textgreek{f}})\bar{n}^{\text{\textgreek{m}}}\le\frac{C_{\text{\textgreek{d}}}}{\text{\textgreek{t}}^{2-\text{\textgreek{d}}}}\mathcal{E}^{2,k}[\text{\textgreek{f}}](0)+\mathcal{F}^{2,k,d,\text{\textgreek{d}}}[F](\text{\textgreek{t}}),\label{eq:Final decay for first energy-1-1-1-1}
\end{equation}
\begin{equation}
\sup_{\{\bar{t}=\text{\textgreek{t}}\}}r^{d-2}\cdot|\text{\textgreek{f}}|^{2}\le\frac{C_{\text{\textgreek{d}}}}{\text{\textgreek{t}}^{2-\text{\textgreek{d}}}}\mathcal{E}^{2,k,d}[\text{\textgreek{f}}](0)+\mathcal{F}^{2,k,d,\text{\textgreek{d}}}[F](\text{\textgreek{t}})\label{eq:DecayNotFriedlanderRadiation-1-1-1}
\end{equation}
and 
\begin{equation}
\sup_{\{\bar{t}=\text{\textgreek{t}}\}}r^{d-1}\cdot|\text{\textgreek{f}}|^{2}\le\frac{C}{\text{\textgreek{t}}}\mathcal{E}^{2,k,d}[\text{\textgreek{f}}](0)+\mathcal{F}^{1,k,d}[F](\text{\textgreek{t}}).\label{eq:DecayFriedlanderRadiation-1-1-1}
\end{equation}

Furthermore, in case the vector field $T=\partial_{u}$ in the coordinate
chart $(u,r,\text{\textgreek{sv}})$ in the asymptotically flat region
$\{r\gg1\}$ of $\mathcal{M}$ satisfies for some (small) $\text{\textgreek{d}}_{0}>0$
and any $k\in\mathbb{N}$:
\begin{align}
\mathcal{L}_{T}^{k}g=O(|u|{}^{-\text{\textgreek{d}}_{0}})\Big\{ & O(r^{-1-a})dvdu+O(r)d\text{\textgreek{sv}}d\text{\textgreek{sv}}+O(1)dud\text{\textgreek{sv}}+\label{eq:DeformationTensorTAwaySlow-1}\\
 & +O(r^{-a})dvd\text{\textgreek{sv}}+O(r^{-1})du^{2}+O(r^{-2-a})dv^{2}\Big\}\nonumber 
\end{align}
and the second term of the right hand side of the integrated local
energy decay estimate (\ref{eq:ILEDwithloss-1}) is replaced by
\begin{equation}
C_{m,\text{\textgreek{h}}}\sum_{j=0}^{m-1}\int_{\{\text{\textgreek{t}}_{1}\le\bar{t}\le\text{\textgreek{t}}_{2}\}\cap\{r\ge R\}}|u|^{-\text{\textgreek{d}}_{0}}r^{-1}J_{\text{\textgreek{m}}}^{N}(N^{j}\text{\textgreek{f}})\bar{n}^{\text{\textgreek{m}}},
\end{equation}
 then the $\text{\textgreek{d}}$-loss in the decay estimates (\ref{eq:Final decay for first energy-1-1-1-1})
and (\ref{eq:DecayNotFriedlanderRadiation-1-1-1}) can be removed:
\begin{equation}
\int_{\{\bar{t}=\text{\textgreek{t}}\}}J_{\text{\textgreek{m}}}^{N}(\text{\textgreek{f}})\bar{n}^{\text{\textgreek{m}}}\le\frac{C}{\text{\textgreek{t}}^{2}}\mathcal{E}^{2,k}[\text{\textgreek{f}}](0)+\mathcal{F}^{2,k,d}[F](\text{\textgreek{t}})\label{eq:Final decay for first energy-1-1-1-1-1}
\end{equation}
and
\begin{equation}
\sup_{\{\bar{t}=\text{\textgreek{t}}\}}r^{d-2}\cdot|\text{\textgreek{f}}|^{2}\le\frac{C}{\text{\textgreek{t}}^{2}}\mathcal{E}^{2,k,d}[\text{\textgreek{f}}](0)+\mathcal{F}^{2,k,d}[F](\text{\textgreek{t}}).\label{eq:DecayNotFriedlanderRadiation-1-1-1-1}
\end{equation}

\end{thm}
See Theorem \ref{thm:FirstPointwiseDecayNewMethod} (and the remark
below it) in Section~\ref{sec:Firstdecay} for a more detailed statement
of the above result and the definition of the weighted energy norms
of $\text{\textgreek{f}}$ and $F$ in the right hand side of (\ref{eq:Final decay for first energy-1-1-1-1}),
(\ref{eq:DecayNotFriedlanderRadiation-1-1-1}) and (\ref{eq:DecayFriedlanderRadiation-1-1-1}).
\begin{rem*}
Let us remark that the initial weighted energy norm on the hyperboloid
$\{\bar{t}=0\}$ in the right hand sides of (\ref{eq:Final decay for first energy-1-1-1-1})--(\ref{eq:DecayNotFriedlanderRadiation-1-1-1-1})
can be readily replaced by a similar weighted norm on a hypersurface
$\text{\textgreek{S}}$ terminating at spacelike infinity (e.\,g.~a
hypersurface which in the asymptotically flat region is of the form
$\{t=const\}$). In that case, the source spacetime energy norms $\mathcal{F}[F]$
in (\ref{eq:Final decay for first energy-1-1-1-1})--(\ref{eq:DecayNotFriedlanderRadiation-1-1-1-1})
are replaced by similar weighted spacetime norms of $F$ over the
region $J^{+}(\text{\textgreek{S}})\cap\{\bar{t}\le\text{\textgreek{t}}\}$. 
\end{rem*}

\subsubsection{\label{sub:IntroductionimprovedPolynomialDecay}An improved $\bar{t}^{-\frac{d}{2}}$
polynomial decay estimate for solutions to the wave equation $\square\text{\textgreek{f}}=0$
in dimensions $d\ge3$}

Finally, we will also be able to establish improved polynomial decay
rates for $\text{\textgreek{f}}$, under some additional restrictions
on the spacetimes $(\mathcal{M},g)$. Let $(\mathcal{M},g)$, satisfy
the geometric assumptions \ref{enu:Asymptotic-flatness} and \ref{enu:NearRegion}
of the previous section, as well as the integrated local energy decay
assumption \ref{enu:Integrated-local-energy}. Let also $\bar{t}$
and $N$ be as in the statement of Theorem \ref{thm:FirstDecayIntroduction}. 

Assume furthermore that $(\mathcal{M},g)$ satisfies the following
two geometric conditions:

\begin{MyDescription}[leftmargin=3.0em]

\item[(GM3)]{\label{enu:KillingFieldsTimelikeSpan} There exist two
smooth vector fields $T,K$ (not necessarily distinct) on $(\mathcal{M},g)$
such that:

\begin{MyDescription}[leftmargin=0.5em]

\item[1.]{$d\bar{t}(T)=d\bar{t}(K)=1$}

\item[2.]{The span of $\{T,K\}$ is everywhere timelike on $\mathcal{M}\backslash\mathcal{H}^{+}$
(where $\mathcal{H}^{+}$ is the future event horizon of $(\mathcal{M},g)$,
which is required to be a subset of $\partial\mathcal{M}_{hor}$).
}

\item[3.]{In the coordinate chart $(u,r,\text{\textgreek{sv}})$
on each connected component of the region $r\gg1$, $T=\partial_{u}$
and $K=T+\text{\textgreek{F}}$ (where $\text{\textgreek{F}}$ is
the generator of a rotation of $\mathbb{S}^{d-1}$, allowed to be
identically $0$).}

\item[4.]{The vector fields $T$ and $K$ are almost Killing in the
sense that there exists a small $\text{\textgreek{d}}_{0}>0$ such
that their deformation tensor satisfies the $O(\bar{t}{}^{-\text{\textgreek{d}}_{0}})$
decay estimates (\ref{eq:DeformationTensorTF}) and (\ref{eq:DeformationTensorTAway}).}

\end{MyDescription}

}

\item [(GM4)]{\label{enu:NonDegenerateHorizon-1} The span of $\{T,K\}$
is tangential to the future event horizon $\mathcal{H}^{+}$ of $(\mathcal{M},g)$
(if non-empty). Moreover, $\mathcal{H}^{+}$ is non-degenerate with
respect to $K$, in the sense that $K$ satisfies $g(K,K)=0$ and
$d\big(g(K,K)\big)\neq0$ on $\mathcal{H}^{+}$.%
\footnote{Hence, $K$ should be viewed as the analogue of the Hawking vector
field of the Kerr spacetime.%
}}

\item [(GM5)]{\label{enu:UniformityOfEstimates}The constants in the
elliptic, Sobolev and Gagliardo--Nirenberg type estimates on the leaves
of the foliation $\{\bar{t}=\text{\textgreek{t}}\}$ stated in Section
\ref{sub:GeometricAssumtionsImprovedDecay} can be chosen to be independent
of $\text{\textgreek{t}}\ge0$.}

\end{MyDescription}
\begin{rem*}
Assumption \ref{enu:UniformityOfEstimates} holds automatically on
spacetimes $(\mathcal{M},g)$ which are near stationary or time periodic.
\end{rem*}
Again, the precise description of Assumptions \ref{enu:KillingFieldsTimelikeSpan}--\ref{enu:UniformityOfEstimates}
is actually more complicated, and will require the splitting of these
assumptions into a larger number of statements: see Assumptions \ref{enu:Transversality}--\ref{enu:BoundednessVolume}
in Section \ref{sub:GeometricAssumtionsImprovedDecay}. 

We will also assume that the following stronger form of Assumption
\ref{enu:Integrated-local-energy} holds:

\begin{MyDescription}[leftmargin=3.0em]

\item [(ILED2)]{\label{enu:Integrated-local-energy-Imp}\emph{Integrated
local energy decay with polynomial loss of derivatives:} There exists
an integer $k\ge0$, such that for any solution $\text{\textgreek{f}}$
to $\square\text{\textgreek{f}}=F$ with suitable boundary conditions
on $\partial_{tim}\mathcal{M}$, any $m\in\mathbb{N}$, $0\le\text{\textgreek{t}}_{1}\le\text{\textgreek{t}}_{2}$,
$\text{\textgreek{h}}>0$, $R>0$ and any integers $i_{1},i_{2}\ge0$
we can bound:
\begin{equation}
\begin{split}\sum_{j=0}^{m}\int_{\{\text{\textgreek{t}}_{1}\le\bar{t}\le\text{\textgreek{t}}_{2}\}\cap\{r\le R\}}|\nabla^{j}(T^{i_{1}}K^{i_{2}}\text{\textgreek{f}})|^{2} & +\sum_{j=1}^{m}\int_{\{\text{\textgreek{t}}_{1}\le\bar{t}\le\text{\textgreek{t}}_{2}\}\cap\partial_{tim}\mathcal{M}}|\nabla^{j}(T^{i_{1}}K^{i_{2}}\text{\textgreek{f}})|^{2}\le\\
\le & C_{m,\text{\textgreek{h}},i_{1},i_{2}}(R)\sum_{j=0}^{m+k-1}\int_{\{\bar{t}=\text{\textgreek{t}}_{1}\}}J_{\text{\textgreek{m}}}^{N}(N^{j}T^{i_{1}}K^{i_{2}}\text{\textgreek{f}})\bar{n}^{\text{\textgreek{m}}}+\\
 & \hphantom{C}+C_{m,\text{\textgreek{h}},i_{1},i_{2}}\sum_{j=0}^{m-1}\int_{\{\text{\textgreek{t}}_{1}\le\bar{t}\le\text{\textgreek{t}}_{2}\}\cap\{r\ge R\}}|\bar{t}|^{-\text{\textgreek{d}}_{0}}r^{-1}J_{\text{\textgreek{m}}}^{N}(N^{j}T^{i_{1}}K^{i_{2}}\text{\textgreek{f}})\bar{n}^{\text{\textgreek{m}}}+\\
 & \hphantom{C}+C_{m,\text{\textgreek{h}},i_{1},i_{2}}(R)\sum_{j=0}^{m+k-1}\int_{\{\text{\textgreek{t}}_{1}\le\bar{t}\le\text{\textgreek{t}}_{2}\}}r^{1+\text{\textgreek{h}}}|\nabla^{j}\square(T^{i_{1}}K^{i_{2}}\text{\textgreek{f}})|^{2}.
\end{split}
\label{eq:ILEDwithloss-1-1}
\end{equation}
}

\end{MyDescription}

See Section \ref{sub:GeometricAssumtionsImprovedDecay} for a more
detailed description of Assumption \ref{enu:Integrated-local-energy-Imp}

For spacetimes $(\mathcal{M},g)$ as above, we will infer the following
improved decay result:
\begin{thm}
\label{thm:ImprovedDecayIntroduction}Let $(\mathcal{M}^{d+1},g)$
satisfy the geometric assumptions \ref{enu:Asymptotic-flatness},
\ref{enu:NearRegion}, \ref{enu:KillingFieldsTimelikeSpan}, \ref{enu:NonDegenerateHorizon-1}
and \ref{enu:UniformityOfEstimates} and the integrated local energy
decay assumption \ref{enu:Integrated-local-energy-Imp}, and let $\bar{t}$
and $N$ be as above. Then for any integer $1\le q\le\lfloor\frac{d+1}{2}\rfloor$,
any $0<\text{\textgreek{e}}\ll\text{\textgreek{d}}_{0}$, any $\text{\textgreek{t}}\ge0$
and any solution $\text{\textgreek{f}}$ to the inhomogeneous wave
equation $\square_{g}\text{\textgreek{f}}=F$ with suitably decaying
intial data on $\{\bar{t}=0\}$ (and satisfying suitable boundary
conditions on $\partial_{tim}\mathcal{M}$) the following estimates
hold: 
\begin{equation}
\mathcal{E}_{en}^{(0,q)}[\text{\textgreek{f}}](\text{\textgreek{t}})+\int_{\text{\textgreek{t}}}^{+\infty}\mathcal{E}_{en}^{(-1+\text{\textgreek{e}},q)}[\text{\textgreek{f}}](s)\, ds\lesssim_{m,\text{\textgreek{e}}}\text{\textgreek{t}}^{-2q+C\text{\textgreek{e}}}\mathcal{E}_{in}^{(2q,k,\text{\textgreek{d}}_{0})}[\text{\textgreek{f}}](0)+\mathcal{F}_{\text{\textgreek{e}}}^{(q,k,m,\text{\textgreek{d}}_{0})}[F](\text{\textgreek{t}})\label{eq:ImprovedNonDegenerateEnergyDecay-1}
\end{equation}
and 
\begin{equation}
\mathcal{E}_{en,deg}^{(0,q)}[\text{\textgreek{f}}](\text{\textgreek{t}})\lesssim_{m,\text{\textgreek{e}}}\text{\textgreek{t}}^{-2q}\mathcal{E}_{in}^{(2q,k,\text{\textgreek{d}}_{0})}[\text{\textgreek{f}}](0)+\mathcal{F}_{deg,\text{\textgreek{e}}}^{(q,k,m,\text{\textgreek{d}}_{0})}[F](\text{\textgreek{t}}).\label{eq:ImprovedEnergyDecay-1}
\end{equation}
In the above, $\mathcal{E}_{en}^{(a,q)}[\text{\textgreek{f}}](\text{\textgreek{t}})$
is the non degenerate $L^{2}$ norm on $\{\bar{t}=\text{\textgreek{t}}\}$
of all derivatives of $\text{\textgreek{f}}$ of order $q$, with
$r^{a}$ weights near infinity, and $\mathcal{E}_{en,deg}^{(a,q)}[\text{\textgreek{f}}](\text{\textgreek{t}})$
is similar to $\mathcal{E}_{en}^{(a,q)}[\text{\textgreek{f}}](\text{\textgreek{t}})$
but with a degeneracy on $\mathcal{H}^{+}$. See Section \ref{sec:Improved-polynomial-decay}
for a more precise definition of these norms.

Moreover, the following pointwise decay rates for $\text{\textgreek{f}}$
are established: 

1. In case the dimension $d$ is odd, we can bound: 
\begin{equation}
\sup_{\{\bar{t}=\text{\textgreek{t}}\}}\big|\text{\textgreek{f}}\big|\lesssim_{m,\text{\textgreek{e}}}\text{\textgreek{t}}^{-\frac{d}{2}}\sqrt{\mathcal{E}_{0,d}[\text{\textgreek{f}}](0)}+\mathcal{F}_{pw,\text{\textgreek{e}}}^{(q,k,0,\text{\textgreek{d}}_{0})}[F](\text{\textgreek{t}}),\label{eq:ImprovedPointwiseDecayAway-1}
\end{equation}

and for any integer $m\ge1$: 
\begin{equation}
\sup_{\{\bar{t}=\text{\textgreek{t}}\}}\big|\nabla^{m}\text{\textgreek{f}}\big|_{h}\lesssim_{m,\text{\textgreek{e}}}\text{\textgreek{t}}^{-\frac{d+1}{2}}\sqrt{\mathcal{E}_{m+2,d}[\text{\textgreek{f}}](0)}+\mathcal{F}_{pw,\text{\textgreek{e}}}^{(q,k,m+2,\text{\textgreek{d}}_{0})}[F](\text{\textgreek{t}}).\label{eq:ImprovedPointwiseDecayAwayDerivative}
\end{equation}

2. In case the dimension $d$ is even, for any integer $m\ge0$ we
can bound:
\begin{equation}
\sup_{\{\bar{t}=\text{\textgreek{t}}\}}\big|\nabla^{m}\text{\textgreek{f}}\big|_{h}\lesssim_{m,\text{\textgreek{e}}}\text{\textgreek{t}}^{-\frac{d}{2}+C\text{\textgreek{e}}}\sqrt{\mathcal{E}_{m,d}[\text{\textgreek{f}}](0)}+\mathcal{F}_{pw,\text{\textgreek{e}}}^{(q,k,\text{\textgreek{m}},\text{\textgreek{d}}_{0})}[F](\text{\textgreek{t}}).\label{eq:ImprovedPointwiseDecay-1}
\end{equation}
For the definition of the weighted energy norms of $\text{\textgreek{f}}$
and $F$ appearing in the right hand sides of the inequalities above,
see Section~\ref{sec:Improved-polynomial-decay}.
\end{thm}
See Theorem~\ref{thm:ImprovedDecayEnergy} and Corollary~\ref{cor:ImprovedPointwiseDecay}
in Section~\ref{sec:Improved-polynomial-decay} for more details.
\begin{rem*}
Let us remark at this point that in the case the integrated local
energy decay statement in Assumption \ref{enu:Integrated-local-energy}
does not lose derivatives (i.\,e.~$k=0$), we can relax the assumption
that the deformation tensors of $T$ and $K$ decay like $\bar{t}^{-\text{\textgreek{d}}_{0}}$
(i.\,e.~(\ref{eq:DeformationTensorTF}) and (\ref{eq:DeformationTensorTAway}))
by replacing it with a uniform $\text{\textgreek{e}}_{0}$-smallness
assumption (i.\,e.~(\ref{eq:DeformationTensorUniformDecay-1-1})--(\ref{eq:DeformationTensorKAwayDecay-1-1}))
for some $\text{\textgreek{e}}_{0}>0$. In this case, we can still
obtain (\ref{eq:ImprovedNonDegenerateEnergyDecay-1}), (\ref{eq:ImprovedPointwiseDecay-1}),
(\ref{eq:ImprovedEnergyDecay-1}) and (\ref{eq:ImprovedPointwiseDecayAway-1}),
at a cost of an $O(\text{\textgreek{e}}_{0})$ loss in the exponent
of $\text{\textgreek{t}}$ in all these inequalities. Thus, in the
absence of trapping, the $r^{p}$-weighted energy method of \cite{DafRod7}
is robust enough to yield the full ``improved'' polynomial hierarchy
on spacetimes that do not settle down to a stationary background.
See also the remark below Theorem \ref{thm:ImprovedDecayEnergy}.

We should also notice that in the case when the vector fields $T$
and $K$ are exactly Killing, the proof of Theorem \ref{thm:ImprovedDecayEnergy}
yields that for any solution $\text{\textgreek{f}}$ to $\square\text{\textgreek{f}}=0$
with compactly supported initial data and any integer $k\ge0$: 
\begin{equation}
|T^{k}\text{\textgreek{f}}|\lesssim_{k}\bar{t}{}^{-1-k}.\label{eq:FastDecay}
\end{equation}
 Therefore, using the frequency cut-off techniques of \cite{DafRod5}
or \cite{Moschidisb}, from (\ref{eq:FastDecay}) (and the corresponding
statement for decay of the energy of $T^{k}\text{\textgreek{f}}$
on the foliation $\{\bar{t}=\text{\textgreek{t}}\}$) we can deduce
that for any $\text{\textgreek{w}}_{0}>0$, $\text{\textgreek{f}}_{\ge\text{\textgreek{w}}_{0}}$
decays superpolynomially in $\bar{t}$ (where $\text{\textgreek{f}}_{\ge\text{\textgreek{w}}_{0}}$
is the part of $\text{\textgreek{f}}$ supported in the frequency
range $|\text{\textgreek{w}}|>\text{\textgreek{w}}_{0}$ with respect
to the $\bar{t}$ variable in a coordinate chart where $T=\partial_{\bar{t}}$).

Finally, as before, we should note that the initial weighted energy
norm on the hyperboloid $\{\bar{t}=0\}$ in the right hand sides of
(\ref{eq:ImprovedNonDegenerateEnergyDecay-1})--(\ref{eq:ImprovedPointwiseDecay-1})
can be readily replaced by a similar weighted norm on a hypersurface
$\text{\textgreek{S}}$ terminating at spacelike infinity. In that
case, the source spacetime energy norms $\mathcal{F}[F]$ in (\ref{eq:ImprovedNonDegenerateEnergyDecay-1})--(\ref{eq:ImprovedPointwiseDecay-1})
are replaced by similar weighted spacetime norms of $F$ over the
region $J^{+}(\text{\textgreek{S}})\cap\{\bar{t}\le\text{\textgreek{t}}\}$. 
\end{rem*}

\subsection{Applications of the $r^{p}$-weighted energy method}

We will now discuss some applications of Theorems \ref{thm:NewMethodModelIntro}--\ref{thm:ImprovedDecayIntroduction}.

\subsubsection{The Friedlander radiation field for solutions to the wave equation
$\square\text{\textgreek{f}}=0$}

On any product Lorentzian manifold of the form $(\mathcal{M}^{d+1},g)=(\mathbb{R}\times\mathcal{S}^{d},-dt^{2}+\bar{g})$,%
\footnote{As usual for product Lorentzian manifolds, $t$ will denote the projection
onto the first factor of $\mathbb{R}\times\mathcal{S}$.%
} where $(\mathcal{S}^{d},\bar{g})$ is an asymptotically Euclidean
Riemannian manifold for $d\ge2$, Friedlander \cite{Friedlander2001}
has established that for any smooth solution $\text{\textgreek{f}}$
to the wave equation $\square_{g}\text{\textgreek{f}}=0$ on $(\mathcal{M},g)$
with compactly supported initial data on $\{t=0\}$, $r^{\frac{d-1}{2}}\cdot\text{\textgreek{f}}$
has a well defined and smooth limit on future null infinity. This
limit is called the future radiation field of $\text{\textgreek{f}}$.
In order to deduce this result, Friedlander utilised the Penrose compactification
method. 

As a soft corollary of the hierarchy of $r^{p}$-weighted estimates
(\ref{eq:newMethodInitialStatementHyperboloids}), we will extend
the result of Friedlander to more general asymptotically flat spacetimes
$(\mathcal{M}^{d+1},g)$, not necessarily of product type, with $d\ge3$:
\begin{thm}
\label{thm:FriedlanderRadiationIntroduction}Let $(\mathcal{M}^{d+1},g)$,
$d\ge3$, be a Lorentzian manifold with with the asymptotics (\ref{eq:RoughAsymptotics}),
in the sense that each connected component of an open subset $\mathcal{N}_{af,\mathcal{M}}$
of $\mathcal{M}$ is mapped diffeomorphically on $\mathbb{R}\times(R_{0},+\infty)\times\mathbb{S}^{d-1}$
through a $(u,r,\text{\textgreek{sv}})$ coordinate chart, in which
$g$ has the form (\ref{eq:RoughAsymptotics}). Then for any smooth
solution $\text{\textgreek{f}}$ to the inhomogeneous wave equation
$\square_{g}\text{\textgreek{f}}=F$ on $(\mathcal{M},g)$ with $(\text{\textgreek{f}},\partial\text{\textgreek{f}})|_{\{t=0\}}$
and $F$ suitably decaying in $r$, the limit 
\begin{equation}
\Big(\lim_{r\rightarrow+\infty}\text{\textgreek{W}}\cdot\text{\textgreek{f}}\Big)(u,r,\text{\textgreek{sv}})\doteq\text{\textgreek{F}}_{\mathcal{I}^{+}}(u,\text{\textgreek{sv}}),
\end{equation}
where $\text{\textgreek{W}}=r^{\frac{d-1}{2}}\big(1+O(r^{-1})\big)$,
exists on all connected components of $\mathcal{N}_{af,\mathcal{M}}$
and defines a smooth function on $\mathbb{R}\times\mathbb{S}^{d-1}$.
Moreover, the following limit exists and is finite for all integers
$j_{1},j_{2},j_{3}\ge0$ 
\begin{equation}
\lim_{r\rightarrow+\infty}\Big(r^{j_{1}}\partial_{r}^{j_{1}}\partial_{\text{\textgreek{sv}}}^{j_{2}}\partial_{u}^{j_{3}}(\text{\textgreek{W}}\text{\textgreek{f}})\Big)<+\infty,\label{eq:LimitTransversalDerivativeI+}
\end{equation}
 where the coordinate derivatives $\partial_{r}$, $\partial_{\text{\textgreek{sv}}}$
and $\partial_{u}$ are considered with respect to the $(u,r,\text{\textgreek{sv}})$
coordinate system in the region $\{r\gg1\}$.
\end{thm}
This result will be established in Section \ref{sec:FriedlanderRadiation}.
For the required decay rates for the initial data of $\text{\textgreek{f}}$
and the source term $F$, see the statement of Theorem \ref{thm:FriedlanderRadiation}. 
\begin{rem*}
Notice that Theorem \ref{thm:FriedlanderRadiationIntroduction} applies
also on spacetimes $(\mathcal{M},g)$ where the decay rate of $g$
near the future null infinity does not allow for a smooth conformal
compactification of the spacetime. Let us also notice that we actually
expect the limit (\ref{eq:LimitTransversalDerivativeI+}) to be identically
$0$ when $j_{1}\ge1$ and $\text{\textgreek{f}}$ solves $\square\text{\textgreek{f}}=0$
with compactly supported initial data, but we do not establish this
fact here.%
\footnote{In case $(\mathcal{M},g)$ admits a conformal compactification near
$\mathcal{I}^{+}$, the stronger statement $\lim_{r\rightarrow+\infty}\Big(r^{2j_{1}}\partial_{r}^{j_{1}}\partial_{\text{\textgreek{sv}}}^{j_{2}}\partial_{u}^{j_{3}}(\text{\textgreek{W}}\text{\textgreek{f}})\Big)<+\infty$
is known to hold for solutions $\text{\textgreek{f}}$ to $\square\text{\textgreek{f}}=0$
with compactly supported initial data.%
}
\end{rem*}
The above result will be established in Section \ref{sec:FriedlanderRadiation}.
For the required decay rates for the initial data of $\text{\textgreek{f}}$
and the source terms $F$, see the statement of Theorem \ref{thm:FriedlanderRadiation}.

\subsubsection{Logarithmic hyperboloidal energy decay for solutions to the wave
equation $\square\text{\textgreek{f}}=0$ on a general class of stationary
asymptotically flat spacetimes}

As described before, the results of the present paper have been used
in our \cite{Moschidisb} to establish the following result:
\begin{thm*}
(Corollary 2.2 of \cite{Moschidisb}). Let $(\mathcal{M}^{d+1},g)$,
$d\ge3$, be a globally hyperbolic spacetime with a Cauchy hypersurface
$\text{\textgreek{S}}$. 

Assume that $(\mathcal{M},g)$ is stationary, with stationary Killing
field $T$, and asymptotically flat. If $\mathcal{M}$ contains a
black hole region bounded by an event horizon $\mathcal{H}$, assume
that $\mathcal{H}$ has positive surface gravity and that the ergoregion
(i.\,e. the set where $g(T,T)>0$) is ``small'' (see \cite{Moschidisb}
for the precise statement of these assumptions). Finally, assume that
an energy boundedness statement is true for solutions to $\square\text{\textgreek{f}}=0$
on the domain of outer communications $\mathcal{D}$ of $\mathcal{M}$. 

It then follows that any smooth solution $\text{\textgreek{f}}$ to
$\square_{g}\text{\textgreek{f}}=0$ on $\mathcal{M}$ with suitably
decaying initial data on a Cauchy hypersurface $\text{\textgreek{S}}$
of $\mathcal{M}$ satisfies on $\mathcal{D}$ for any integer $m>0$:

\begin{equation}
E_{hyp}(\bar{t})\le\frac{C_{m}}{\big\{\log(2+\bar{t})\big\}^{2m}}\big(E_{hyp}^{(m)}(0)+E_{w,hyp}(0)\big).\label{eq:CariqatureInequality}
\end{equation}
 In the above, $\bar{t}\ge0$ is a suitable time function on $J^{+}(\text{\textgreek{S}})\cap\mathcal{D}$
with hyperboloidal level sets, satisfying $T(\bar{t})=1$, and $E_{hyp}(\bar{t})$
is the energy flux of $\text{\textgreek{f}}$ with respect to the
level sets of the time function$\bar{t}$. $E_{hyp}^{(m)}(0)$ is
the energy of the first $m$ derivatives of $\text{\textgreek{f}}$
at $\{\bar{t}=0\}$, while $E_{w}(0)$ is a suitable weighted energy
of $\text{\textgreek{f}}$ at $\{\bar{t}=0\}$. The constant $C$
on the right hand side depends on the geometry of $(\mathcal{D},g)$
and the precise choice of the function $\bar{t}$, while in addition
to that, $C_{m}$ also depends on the number $m$ of derivatives of
$\text{\textgreek{f}}$ in $E^{(m)}(0)$.
\end{thm*}
For a more detailed statement of the above result, see \cite{Moschidisb}.

\subsubsection{\label{sub:KerrDecay}Polynomial decay for solutions to the wave
equation $\square\text{\textgreek{f}}=0$ on the exterior of subextremal
Kerr spacetimes for $|a|<M$}

In the next three sections, we will introduce some examples of spacetimes
$(\mathcal{M},g)$ which satisfy the assumptions of Theorems \ref{thm:NewMethodModelIntro}--\ref{thm:ImprovedDecayIntroduction}. 

Our first such example will be the exterior of a subextremal Kerr
spacetime $(\mathcal{M}_{a,M},g_{a,M})$ with parameters lying in
the full subextremal range $|a|<M$. Notice that this spacetime satisfies
all the geometric assumptions of Theorems \ref{thm:NewMethodModelIntro}--\ref{thm:ImprovedDecayIntroduction}.
In fact, the form of the assumptions of Theorem \ref{thm:ImprovedDecayIntroduction}
was motivated by the geometry of the subextremal Kerr family. In \cite{DafRod9},
the authors have established an energy boundedness and integrated
local energy decay statement for solutions to (\ref{eq:WaveEquation})
on $(\mathcal{M}_{a,M},g_{a,M})$. As already noted in \cite{DafRod9},
by applying Theorems \ref{thm:FirstPointwiseDecayNewMethod} and \ref{thm:ImprovedDecayEnergy}
one can thus readily upgrade these results to polynomial decay estimates
for solutions to (\ref{eq:WaveEquation}), and therefore establish
Corollary 3.1 of \cite{DafRod9}, which we state here with the notation
of \cite{DafRod9}:
\begin{cor*}
(Corollary 3.1 of \cite{DafRodSchlap}) Let $(\mathcal{M}_{a,M},g_{a,M})$
be the exterior of a Kerr black hole spacetime of mass $M$ and angular
momentum $a$, such that $|a|<M$. Let $\tilde{\text{\textgreek{S}}}_{0}$
be a smooth spacelike hypersurface of $(\mathcal{M}_{a,M},g_{a,M})$
intersecting transversally the future event horizon $\mathcal{H}^{+}$
and terminating at fututre null infinity $\mathcal{I}^{+}$. Let also
$\tilde{\text{\textgreek{S}}}_{\text{\textgreek{t}}}$ denote the
image of $\tilde{\text{\textgreek{S}}}_{0}$ under the flow of the
stationary Killing field $T$ of $(\mathcal{M}_{a,M},g_{a,M})$ (see
\cite{DafRodSchlap}), and let $N$ be a globally timelike, future
directed and $T$-invariant vector field on $\mathcal{M}_{a,M}$ coinciding
with $T$ in the region $\{r\gg1\}$. 

Then, for any $\text{\textgreek{d}}>0$, there exists a constant $C=C(a,M,\tilde{\text{\textgreek{S}}}_{0},\text{\textgreek{d}})>0$
such that for any smooth solution $\text{\textgreek{f}}$ to the wave
equation (\ref{eq:WaveEquation}) on $J^{+}(\tilde{\text{\textgreek{S}}}_{0})\subset(\mathcal{M}_{a,M},g_{a,M})$
with suitably decaying initial data on $\tilde{\text{\textgreek{S}}}_{0}$
the following energy decay estimates hold:

\begin{equation}
\int_{\tilde{\text{\textgreek{S}}}_{\text{\textgreek{t}}}}J_{\text{\textgreek{m}}}^{N}(\text{\textgreek{f}})n_{\tilde{\text{\textgreek{S}}}}^{\text{\textgreek{m}}}\le C\cdot E\text{\textgreek{t}}^{-2},
\end{equation}
\begin{equation}
\int_{\tilde{\text{\textgreek{S}}}_{\text{\textgreek{t}}}}J_{\text{\textgreek{m}}}^{N}(N\text{\textgreek{f}})n_{\tilde{\text{\textgreek{S}}}}^{\text{\textgreek{m}}}\le C\cdot E\text{\textgreek{t}}^{-4+\text{\textgreek{d}}}.
\end{equation}
 Moreover, the following pointwise decay estimates hold: 
\begin{equation}
\sup_{\tilde{\text{\textgreek{S}}}_{\text{\textgreek{t}}}}r\cdot\big|\text{\textgreek{f}}|\le C\sqrt{E}\cdot\text{\textgreek{t}}^{-\frac{1}{2}},
\end{equation}
\begin{equation}
\sup_{\tilde{\text{\textgreek{S}}}_{\text{\textgreek{t}}}}\big|\text{\textgreek{f}}|\le C\sqrt{E}\cdot\text{\textgreek{t}}^{-\frac{3}{2}}
\end{equation}
 and 
\begin{equation}
\sup_{\tilde{\text{\textgreek{S}}}_{\text{\textgreek{t}}}}\big(\big|N\text{\textgreek{f}}|+|\nabla_{\tilde{\text{\textgreek{S}}}}\text{\textgreek{f}}|\big)\le C\sqrt{\text{\textgreek{E}}}\cdot\text{\textgreek{t}}^{-2}.
\end{equation}
 In the above, $E$ denotes a suitable higher order weighted energy
norm of the intial data of $\text{\textgreek{f}}$ on $\tilde{\text{\textgreek{S}}}_{0}$,
and is not necessarily the same quantity in all of the above estimates. 
\end{cor*}
See \cite{DafRodSchlap} for more details.
\begin{rem*}
Notice that the slowly rotating $4+1$ dimensional Myers--Perry spacetimes
satisfy all the geometric assumptions of Theorem \ref{thm:ImprovedDecayIntroduction}.
Therefore, in view of the integrated local energy decay estimate established
in \cite{Laul2015}, Theorem \ref{thm:ImprovedDecayIntroduction}
implies that any solution $\text{\textgreek{f}}$ to (\ref{eq:WaveEquation})
on a slowly rotating $4+1$ dimensional Myers--Perry spacetime with
suitably decaying initial data satisfies a $\bar{t}^{-2+\text{\textgreek{d}}}$
pointwise decay estimate.
\end{rem*}

\subsubsection{\label{sub:PerturbationsMinkowski}Improved polynomial decay on radiating
uniformly small perturbations of Minkowski spacetime}

For our second example of a spacetime satisfying the assumptions of
Theorems \ref{thm:NewMethodModelIntro}--\ref{thm:ImprovedDecayIntroduction},
we will need to introduce a definition: We will define a metric $g$
on $\mathbb{R}^{d+1}$, $d\ge3$, to be a \emph{radiating uniformly
small perturbation of Minkowski spacetime} if there exists a (small)
$\text{\textgreek{e}}_{0}>0$ and an $R>0$ such that, in the $(u,r,\text{\textgreek{sv}})$
coordinate system on $\mathbb{R}^{d+1}$ in the region $\{r\ge R\}$,
$g$ is of the form (\ref{eq:MetricUR}) for some $0<a\le1$, and
moreover:
\begin{itemize}
\item For any integers $m_{1},m_{2}\ge0$ we have the global bound: 
\begin{equation}
\sup_{\mathbb{R}^{d+1}}|\mathcal{L}_{T}^{m_{1}}\nabla_{e}^{m_{2}}(g-\text{\textgreek{h}})|_{e}\lesssim_{m_{1},m_{2}}\text{\textgreek{e}}_{0}(1+r)^{-1}\min\{1,|u|^{1-m_{1}}\}\label{eq:Uniform boundedness}
\end{equation}

\item In the region $\{r\ge R\}$ we can estimate for any $m\ge1$ in the
$(u,r,\text{\textgreek{sv}})$ coordinate system: 
\begin{align}
\mathcal{L}_{T}^{m}g=O_{m}(\text{\textgreek{e}}_{0}\min\{1,|u|^{1-m}\})\Big\{ & O(r^{-1-a})drdu+O(r)d\text{\textgreek{sv}}d\text{\textgreek{sv}}+O(1)dud\text{\textgreek{sv}}+\label{eq:DeformationTensorTAwayDecay-1-1-1}\\
 & +O(r^{-a})drd\text{\textgreek{sv}}+O(r^{-1})du^{2}+O(r^{-2-a})dr^{2}\Big\}.\nonumber 
\end{align}

\end{itemize}
In the above, $T$ is the vector field $\partial_{t}$ in the Cartesian
coordinate system $(t,x^{1},\ldots,x^{d})$ on $\mathbb{R}^{d+1}$,
$e$ is the usual Euclidean metric on $\mathbb{R}^{d+1}$ and $\nabla_{e}$
is the flat connection on $\mathbb{R}^{d+1}$. Notice that if $\text{\textgreek{e}}_{0}$
is smaller than an absolute constant, $T$ is everywhere timelike
and furthermore $(\mathbb{R}^{d+1},g)$ can not contain any trapped
geodesics. In fact, if $\text{\textgreek{e}}_{0}$ is small enough,
the $\partial_{r}$-Morawetz current of Minkowski spacetime (combined
with the estimates of Section \ref{sec:Morawetz} of the present paper)
yields an integrated local energy decay estimate of the form (\ref{eq:ILEDwithloss-1})
without loss of derivatives. Furthermore, the rest of the geometric
Assumptions of Theorems \ref{thm:NewMethodModelIntro}--\ref{thm:ImprovedDecayIntroduction}
are satisfied, except for the assumption on the $\bar{t}^{-\text{\textgreek{d}}_{0}}$
decay of deformation tensors of $T$ and $K$ which is replaced by
a uniform $\text{\textgreek{e}}_{0}$-smallness assumption (see the
remark below Theorem \ref{thm:ImprovedDecayIntroduction}).

For such a spacetime $(\mathbb{R}^{d+1},g)$, we will fix $\mathcal{S}\subset\mathbb{R}^{d+1}$
to be a smooth spacelike hypersurface of $(\mathbb{R}^{d+1},g)$ which
terminates at $\mathcal{I}^{+}$, and let $\bar{t}:\mathbb{R}^{d+1}\rightarrow\mathbb{R}$
be defined by the condition $T(\bar{t})=1$ and $\bar{t}|_{\mathcal{S}}=0$. 

One can deduce from \cite{Christodoulou1993} that dynamical solutions
of the vacuum Einstein equations arising from initial data which are
close to the ones for Minkowski spacetime are included in this class. 

The following pointwise decay estimate for solutions to the wave equation
on radiating uniformly small perturbations of Minkowski spacetime
is a straightforward application of Theorems \ref{thm:FirstDecayIntroduction}
and \ref{thm:ImprovedDecayIntroduction}:
\begin{cor*}
Let $(\mathbb{R}^{d+1},g)$, $d\ge3$, be a uniformly small perturbation
of Minkowski spacetime, in the sense that for some $\text{\textgreek{e}}_{0}>0$
and $R>0$, $g$ is of the form (\ref{eq:MetricUR}) in the region
$\{r\ge R\}$ and (\ref{eq:Uniform boundedness}) and (\ref{eq:DeformationTensorTAwayDecay-1-1-1})
hold. Let also $\bar{t}:(\mathbb{R}^{d+1},g)\rightarrow\mathbb{R}$
be constructed as above. Then, provided $\text{\textgreek{e}}_{0}$
is smaller than an absolute constant, for any solution $\text{\textgreek{f}}$
to the wave equation $\square_{g}\text{\textgreek{f}}=0$ on $(\mathbb{R}^{d+1},g)$
and any $\text{\textgreek{t}}\ge0$ we can bound 
\begin{equation}
\sup_{\{\bar{t}=\text{\textgreek{t}}\}}|\text{\textgreek{f}}|\lesssim\text{\textgreek{t}}^{-\frac{d}{2}+O(\text{\textgreek{e}}_{0})}\sqrt{\mathcal{E}_{0,d}[\text{\textgreek{f}}](0)}\label{eq:ImprovedPointwiseMinkowski}
\end{equation}
and, for any integer $m\ge1$: 
\begin{equation}
\sup_{\{\bar{t}=\text{\textgreek{t}}\}}|\nabla^{m}\text{\textgreek{f}}|_{e}\lesssim\text{\textgreek{t}}^{-\frac{d+1}{2}+O(\text{\textgreek{e}}_{0})}\sqrt{\mathcal{E}_{m+2,d}[\text{\textgreek{f}}](0)}.\label{eq:ImprovedDerivativesMinkowski}
\end{equation}

In case the following stronger assumptions on the deformation tensor
of $T$ hold for some $\text{\textgreek{d}}_{0}>0$ and any $m_{1}\ge1$,
$m_{2}\ge0$ in place of (\ref{eq:Uniform boundedness}) and (\ref{eq:DeformationTensorTAwayDecay-1-1-1}):
\begin{equation}
\sup_{\mathbb{R}^{d+1}}|r\mathcal{L}_{T}^{m_{1}}\nabla_{e}^{m_{2}}(g-\text{\textgreek{h}})|_{e}\lesssim_{m_{1},m_{2}}|u|^{1-m_{1}-\text{\textgreek{d}}_{0}}\label{eq:Uniform boundedness-1}
\end{equation}
and in the region $\{r\ge R\}$ for any $m\ge1$: 
\begin{align}
\mathcal{L}_{T}^{m}g=O_{m}(|u|^{1-m-\text{\textgreek{d}}_{0}})\Big\{ & O(r^{-1-a})drdu+O(r)d\text{\textgreek{sv}}d\text{\textgreek{sv}}+O(1)dud\text{\textgreek{sv}}+\label{eq:DeformationTensorTAwayDecay-1-1-1-1}\\
 & +O(r^{-a})drd\text{\textgreek{sv}}+O(r^{-1})du^{2}+O(r^{-2-a})dr^{2}\Big\},\nonumber 
\end{align}
then (\ref{eq:ImprovedPointwiseMinkowski}) and (\ref{eq:ImprovedDerivativesMinkowski})
can be upgraded to 
\begin{equation}
\sup_{\{\bar{t}=\text{\textgreek{t}}\}}|\text{\textgreek{f}}|\lesssim\text{\textgreek{t}}^{-\frac{d}{2}}\sqrt{\mathcal{E}_{0,d}[\text{\textgreek{f}}](0)}\label{eq:ImprovedPointwiseMinkowski-1}
\end{equation}
and, for any integer $m\ge1$: 
\begin{equation}
\sup_{\{\bar{t}=\text{\textgreek{t}}\}}|\nabla^{m}\text{\textgreek{f}}|_{e}\lesssim\text{\textgreek{t}}^{-\frac{d+1}{2}}\sqrt{\mathcal{E}_{m+2,d}[\text{\textgreek{f}}](0)}.\label{eq:ImprovedDerivativesMinkowski-1}
\end{equation}
 For the definition of the initial energy norms $\mathcal{E}_{0,d}[\text{\textgreek{f}}](0)$
and $\mathcal{E}_{m+2,d}[\text{\textgreek{f}}](0)$ on the hypersurfaces
$\{\bar{t}=0\}$ (which can also be replaced by norms on $\{t=0\}$),
see Section \ref{sub:ShorthandEnergyNorms}.\end{cor*}
\begin{rem*}
Notice that the above corollary extends a recent result of Oliver
\cite{Oliver2014}.
\end{rem*}

\subsubsection{Improved polynomial decay on dynamical, radiating black hole spacetimes}

A final example of a class of spacetimes satisfying the assumptions
of Theorems \ref{thm:NewMethodModelIntro}--\ref{thm:ImprovedDecayIntroduction}
will concern the exterior region of black hole spacetimes dynamically
settling down to a subextremal Kerr spacetime. Here, we will restrict
ourselves only to spacetimes $(\mathcal{M}_{Sch},g)$ (where $\mathcal{M}_{Sch}$
has the differentiable structure of the Schwarzschild exterior) settling
down to the Schwarzschild exterior spacetime $(\mathcal{M}_{Sch},g_{M})$
for some $M>0$ at a sufficiently fast polynomial rate. In particular,
we will assume that we can fix a double null foliation on $(\mathcal{M}_{Sch},g)$
such that the components of $g$ with respect to this foliation approach
the components of the Schwarzschild metric $g_{M}$ at a sufficiently
fast polynomial rate towards ``timelike infinity''. This class of
spacetimes includes, in particular, the radiating spacetimes constructed
in \cite{Dafermos2013}, which approach the Schwarschild metric at
an exponential rate. We will not provide more details of this setup
here, but instead we will refer the reader to \cite{Dafermos2013}.
The reason for this restriction is that it is straightforward to check
(essentially without calculation) that these spacetimes $(\mathcal{M}_{Sch},g)$
satisfy the assumptions of Theorems \ref{thm:NewMethodModelIntro}--\ref{thm:ImprovedDecayIntroduction}
(we will omit the details). 

On spacetimes $(\mathcal{M}_{Sch},g)$ as above, the energy current
yielding the integrated local energy decay statement for the Schwarzschild
exterior $(\mathcal{M}_{Sch},g_{M})$ constructed in \cite{DafRod4},
combined with the estimates of Section \ref{sec:Morawetz} of the
present paper and the fast rate at which $g$ approaches the Schwarzschild
metric $g_{M}$, imply that an integrated local energy decay statement
of the form (\ref{eq:ILEDwithloss-1-1}) also holds on $(\mathcal{M}_{Sch},g)$.
Furthermore, $(\mathcal{M}_{Sch},g)$ also satisfies the rest of the
geometric assumptions of Theorem \ref{thm:ImprovedDecayIntroduction}.
Therefore, as an application of Theorem \ref{thm:ImprovedDecayIntroduction},
we obtain the following result:
\begin{cor*}
Let $(\mathcal{M}_{Sch},g)$ be a radiating spacetime approaching
$(\mathcal{M}_{Sch},g_{M})$ in the future at a sufficiently fast
polynomial rate (in the sense described above). Let also $\bar{t}:\mathcal{M}_{Sch}\rightarrow\mathbb{R}$
be a function with spacelike level sets intersecting $\mathcal{H}^{+}$
transversally and terminating at $\mathcal{I}^{+}$, such that $T(\bar{t})=1$
(where $T$ is the Schwarzschild stationary Killing field). Then for
any solution $\text{\textgreek{f}}$ to $\square_{g}\text{\textgreek{f}}=0$
on $(\mathcal{M}_{Sch},g)$ with suitably decaying initial data on
a Cauchy hypersurface, the following pointwise decay estimates hold:
\begin{equation}
\sup_{\{\bar{t}=\text{\textgreek{t}}\}}\big|\text{\textgreek{f}}\big|\lesssim_{m,\text{\textgreek{e}}}\text{\textgreek{t}}^{-\frac{3}{2}}\sqrt{\mathcal{E}_{0}[\text{\textgreek{f}}](0)}\label{eq:ImprovedPointwiseDecayAway-1-1}
\end{equation}
and for any integer $m\ge1$: 
\begin{equation}
\sup_{\{\bar{t}=\text{\textgreek{t}}\}}\big|\nabla^{m}\text{\textgreek{f}}\big|\lesssim_{m}\text{\textgreek{t}}^{-2}\sqrt{\mathcal{E}_{m+2}[\text{\textgreek{f}}](0)}.\label{eq:ImprovedPointwiseDecayAwayDerivative-1}
\end{equation}
 For the definition of the initial energy norms $\mathcal{E}_{0}[\text{\textgreek{f}}](0)$
and $\mathcal{E}_{m+2}[\text{\textgreek{f}}](0)$ on the hypersurfaces
$\{\bar{t}=0\}$ (which can also be replaced by norms on a Cauchy
hypersurface), see Section \ref{sub:ShorthandEnergyNorms}\end{cor*}
\begin{rem*}
We should notice that the spacetimes constructed in \cite{Dafermos2013}
are only $C^{l}$ on the future event horizon $\mathcal{H}^{+}$ for
some sufficiently large $l$, but not $C^{\infty}$. However, Theorem
\ref{thm:ImprovedDecayIntroduction} still applies in this case, and
the above Corollary holds provided the integer $m$ in (\ref{eq:ImprovedPointwiseDecayAwayDerivative-1})
is restricted to lie below some constant $C(l)$ depending on the
order of differentiability of $g$ on $\mathcal{H}^{+}$.
\end{rem*}

\subsection{Outline of the paper and technical comments}

In this section, we will describe briefly how the current paper is
organised and we will sketch the difficulties arising in the proof
of the main statements. In particular, we will point out the new difficulties
that appear in comparison to \cite{DafRod7,Schlue2013}. The reader
might find it helpful to return to this Section after viewing the
detailed setup of the Theorems in Sections \ref{sec:Morawetz}--\ref{sec:Improved-polynomial-decay}.

The geometry of the asymptotically fat region $\mathcal{N}_{af}$
of the spacetimes $(\mathcal{M},g)$ under consideration is introduced
in Section \ref{sec:GeometrySpacetimes}. In this region, a function
$\bar{t}$ with hyperboloidal level sets is constructed. It is also
shown that in $\mathcal{N}_{af}$, the wave operator takes the form:
\begin{align}
\text{\textgreek{W}}\cdot\square\text{\textgreek{f}}=-\big(1 & +O(r^{-1-a})\big)\cdot\partial_{u}\partial_{v}(\text{\textgreek{W}}\text{\textgreek{f}})+r^{-2}\text{\textgreek{D}}_{g_{\mathbb{S}^{d-1}}+O(r^{-1})}(\text{\textgreek{W}}\text{\textgreek{f}})-\label{eq:ConformalWaveOperator-3}\\
 & -\frac{(d-1)(d-3)}{4}r^{-2}\cdot(\text{\textgreek{W}}\text{\textgreek{f}})+Err(\text{\textgreek{W}}\text{\textgreek{f}}),\nonumber 
\end{align}
where $\text{\textgreek{W}}=r^{\frac{d-1}{2}}\big(1+O(r^{-1})\big)$
and the ``error'' terms $Err(\text{\textgreek{W}}\text{\textgreek{f}})$
have the form (\ref{eq:ErrTerms}). Notice that the particular choice
of the factor $\text{\textgreek{W}}$ serves to eliminate some terms
in the expression for $Err(\text{\textgreek{W}}\text{\textgreek{f}})$
which would be ``problematic'' in the derivation of the $r^{p}$-weighted
energy estimates (\ref{eq:newMethodInitialStatementHyperboloids})
and (\ref{eq:newMethodInitialStatementHyperboloids-1}) (such terms
would appear, for instance, if one substituted $\text{\textgreek{W}}$
by $r^{\frac{d-1}{2}}$ in the case when $\partial_{u}M\neq0$ in
(\ref{eq:RoughAsymptotics}))

In Section \ref{sec:Morawetz}, we establish $\partial_{r}$-Morawetz
and $J^{T}$-energy boundedness estimates of the form
\begin{equation}
\begin{split}\int_{\{\text{\textgreek{t}}_{1}\le\bar{t}\le\text{\textgreek{t}}_{2}\}\cap\{r\ge R\}}r^{-1-\text{\textgreek{h}}}\big(|\partial & \text{\textgreek{f}}|^{2}+r^{-2}|\text{\textgreek{f}}|^{2}\big)\lesssim_{\text{\textgreek{h}}}\int_{\{\bar{t}=\text{\textgreek{t}}_{1}\}\cap\{r\ge R\}}J_{\text{\textgreek{m}}}^{T}(\text{\textgreek{f}})\bar{n}^{\text{\textgreek{m}}}+\int_{\{\text{\textgreek{t}}_{1}\le\bar{t}\le\text{\textgreek{t}}_{2}\}\cap\{r\sim R\}}\big(|\partial\text{\textgreek{f}}|^{2}+r^{-2}|\text{\textgreek{f}}|^{2}\big)+\\
 & +\int_{\{\text{\textgreek{t}}_{1}\le\bar{t}\le\text{\textgreek{t}}_{2}\}\cap\{r\ge R\}}r^{1+\text{\textgreek{h}}}|\square\text{\textgreek{f}}|^{2}+\int_{\{\text{\textgreek{t}}_{1}\le\bar{t}\le\text{\textgreek{t}}_{2}\}\cap\{r\ge R\}}r^{-1}\big(|\partial_{v}\text{\textgreek{f}}|^{2}+|r^{-1}\partial_{\text{\textgreek{sv}}}\text{\textgreek{f}}|^{2}+r^{-2}|\text{\textgreek{f}}|^{2}\big)
\end{split}
\label{eq:SketchMorawetz}
\end{equation}
 and 
\begin{equation}
\begin{split}\int_{\{\bar{t}=\text{\textgreek{t}}_{1}\}\cap\{r\ge R\}}J_{\text{\textgreek{m}}}^{T}(\text{\textgreek{f}} & )\bar{n}^{\text{\textgreek{m}}}\lesssim_{\text{\textgreek{h}}}\int_{\{\bar{t}=\text{\textgreek{t}}_{1}\}\cap\{r\ge R\}}J_{\text{\textgreek{m}}}^{T}(\text{\textgreek{f}})\bar{n}^{\text{\textgreek{m}}}+\int_{\{\text{\textgreek{t}}_{1}\le\bar{t}\le\text{\textgreek{t}}_{2}\}\cap\{r\sim R\}}\big(|\partial\text{\textgreek{f}}|^{2}+r^{-2}|\text{\textgreek{f}}|^{2}\big)+\\
 & +\int_{\{\text{\textgreek{t}}_{1}\le\bar{t}\le\text{\textgreek{t}}_{2}\}\cap\{r\ge R\}}r^{1+\text{\textgreek{h}}}|\square\text{\textgreek{f}}|^{2}+\int_{\{\text{\textgreek{t}}_{1}\le\bar{t}\le\text{\textgreek{t}}_{2}\}\cap\{r\ge R\}}r^{-1}\big(|\partial_{v}\text{\textgreek{f}}|^{2}+|r^{-1}\partial_{\text{\textgreek{sv}}}\text{\textgreek{f}}|^{2}+r^{-2}|\text{\textgreek{f}}|^{2}\big),
\end{split}
\label{eq:SketchBoundedness}
\end{equation}
respectively. We notice that the last terms of the right hand sides
of (\ref{eq:SketchMorawetz}) and (\ref{eq:SketchBoundedness}) appear
due to the radiating asymptotics of (\ref{eq:RoughAsymptotics}),
and can be completely dropped in the case the spacetime is non radiating
or when the radiating components of (\ref{eq:RoughAsymptotics}) satisfy
some special monotonicity conditions (which are satisfied in the case
when $\partial_{u}M\le0$ in (\ref{eq:RoughAsymptotics}) and the
spacetime is spherically symmetric). 

In Section \ref{sec:The-new-method}, the $r^{p}$-weighted energy
hierarchy (\ref{eq:newMethodInitialStatementHyperboloids}) is established.
This is achieved by multiplying the expression (\ref{eq:ConformalWaveOperator-3})
by $r^{p}\partial_{v}(\text{\textgreek{W}}\text{\textgreek{f}})$
and then integrating by parts (in the top order terms) over a region
of the form $\{\text{\textgreek{t}}_{1}\le\bar{t}\le\text{\textgreek{t}}_{2}\}$
(athough regions of different ``shape'' are also treated). In this
integration by parts procedure, the error terms occuring from the
$Err(\text{\textgreek{W}}\text{\textgreek{f}})$ summands are controlled
with the help of the already positive definite terms in the resulting
expression, after adding to it the estimates (\ref{eq:SketchMorawetz})
and (\ref{eq:SketchBoundedness}), using also a Hardy-type inequality
for the zeroth order terms. It is in this procedure that the elimination
of the ``worst'' terms in $Err(\text{\textgreek{W}}\text{\textgreek{f}})$,
resulting from the precise choice of the factor $\text{\textgreek{W}}$
in (\ref{eq:ConformalWaveOperator-3}), is important.

In Section \ref{sec:The-improved--hierarchy}, the higher order $r^{p}$-weighted
energy hierarchy (\ref{eq:newMethodInitialStatementHyperboloids-1})
is established. This is achieved by commuting equation (\ref{eq:ConformalWaveOperator-3})
with $\partial_{v}$, $r^{-1}\partial_{\text{\textgreek{sv}}}$ and
$r^{-1}\partial_{u}$, and repeating the proof leading to (\ref{eq:newMethodInitialStatementHyperboloids}),
after noticing that the first two commutation vector fields lead to
the appearence of some new bulk terms with favorable sign. Notice
that in this procedure, multiple integrations by parts are performed
also on lower order terms, in order to guarantee that (\ref{eq:newMethodInitialStatementHyperboloids-1})
is valid even at the (upper) endpoint $p=2k$.%Mhpws na sygkrinw me Volker? 

In Section \ref{sec:FriedlanderRadiation}, Theorem \ref{thm:FriedlanderRadiationIntroduction}
concerning the existence of the Friedlander radiation field is formulated
and established with the use of the boundary terms controlled by the
$r^{p}$-weighted energy hierarchy (\ref{eq:newMethodInitialStatementHyperboloids-1})
(for $2k-1<p\le2$), combined simply with the fundamental theorem
of calculus. As a corollary of Theorem \ref{thm:FriedlanderRadiationIntroduction},
it is shown the $L^{2}$ norm of certain derivatives of the radiation
field of $\text{\textgreek{f}}$ on $\mathcal{I}^{+}$ are also controlled
by the right hand side of (\ref{eq:newMethodInitialStatementHyperboloids})
and (\ref{eq:newMethodInitialStatementHyperboloids-1}).

In Section \ref{sec:Firstdecay}, the geometric conditions and the
integrated local energy decay assumption (consistent with the right
hand side error term in (\ref{eq:SketchMorawetz})) on the spacetimes
$(\mathcal{M},g)$, on which Theorem \ref{thm:FirstDecayIntroduction}
(concerning the $\bar{t}^{-1}$ decay estimates for $\text{\textgreek{f}}$)
applies, are introduced. The proof of Theorem \ref{thm:FirstDecayIntroduction}
then follows by applying the pigeonhole principle on the positive
bulk terms of the hierarchy (\ref{eq:newMethodInitialStatementHyperboloids}),
as was first done in \cite{DafRod7}. Notice that, since no energy
boundedness statement is a priori assumed, in order to obtain the
final decay estimate in this procedure, an energy boundedness estimate
with loss of derivatives (and right hand side error terms consistent
with (\ref{eq:SketchBoundedness})) is established.%
\footnote{Notice, however, that this energy estimate can yield a ``true''
energy boundedness estimate only in the case where the spacetime $(\mathcal{M},g)$
is non-radiating or when the radiating components of (\ref{eq:RoughAsymptotics})
satisfy some special monotonicity conditions.%
}

Finally, in Section \ref{sec:Improved-polynomial-decay}, the extra
geometric conditions on $(\mathcal{M},g)$, required for Theorem \ref{thm:ImprovedDecayIntroduction}
(concerning the $\bar{t}^{-\frac{d}{2}}$ decay estimates for $\text{\textgreek{f}}$)
to hold, are formulated. The proof of Theorem \ref{thm:ImprovedDecayIntroduction}
then follows by repeated applying the pigeonhole principle argument
on the higher order hierarchy (\ref{eq:newMethodInitialStatementHyperboloids-1})
for higher $T$ and $\text{\textgreek{q}}\cdot K$ derivatives of
$\text{\textgreek{f}}$ (where $\text{\textgreek{q}}$ is a suitable
compactly supported cut-off function). In each step in this procedure,
the wave equation is used to substitute derivatives of $\text{\textgreek{f}}$
tangential to the hyperboloids $\{\bar{t}=const\}$ with $T$ and
$\text{\textgreek{q}}\cdot K$ derivatives of $\text{\textgreek{f}}$,
in a fashion similar to \cite{Schlue2013}. In the end, however, the
wave equation is used again to transform decay estimates of $T$ and
$\text{\textgreek{q}}\cdot K$ derivatives of $\text{\textgreek{f}}$
into decay estimates of certain elliptic operators on $\{\bar{t}=const\}$
applied on $\text{\textgreek{f}}$, and then the elliptic estimates
of Section \ref{sec:Elliptic-estimates} of the Appendix yield decay
estimates for the energy of all higher order derivatives of $\text{\textgreek{f}}$.
Pointwise decay estimates for $\text{\textgreek{f}}$ then follow
by applying the Gagliardo--Nirenberg type estimates of Section \ref{sub:GagliardoNirenberg}.
Thus, our method for extracting $\bar{t}^{-\frac{d}{2}}$ decay estimates
differs substantially from the method implemented in \cite{Schlue2013}
(which yielded $\bar{t}^{-\frac{3}{2}+\text{\textgreek{d}}}$ decay
estimates). For a more detailed sketch of the proof of Theorem \ref{thm:ImprovedDecayIntroduction},
see Section \ref{sub:SketchOfProof}. Let us remark that the elliptic
estimates of Section \ref{sec:Elliptic-estimates} are also used to
control error terms arising from the commutations with the truncated
vector field $\text{\textgreek{q}}\cdot K$.

\subsection{Acknowledgements}

I would like to express my gratitude to my advisor Mihalis Dafermos
for suggesting this problem to me, offering his advice and support
and providing comments, ideas and assistance while this paper was
being written. I would also like to thank Igor Rodnianski for his
invaluable ideas and suggestions, as well as for our discussions on
many aspects of the problem. Finally, I would like to thank Stefanos
Aretakis, Yakov Shlapentokh-Rothman and Yannis Angelopoulos for many
insightful conversations and important comments on preliminary versions
of this paper.

\section{\label{sub:Notational-conventions}Notational conventions}

\subsection{Conventions on constants and inequality symbols}

We will use capital letters (e.\,g.~$C$) to denote ``large''
constants, namely constants that appear on the right hand side of
inequalities, and hence can be replaced by larger ones without affecting
the validity of the inequality. Lower case letters (e.\,g.~$c$)
will be used to denote small constants (which can similarly freely
be replaced by smaller ones). Moreover, the same character might be
used to denote different constants even in adjacent lines or formulas. 

We will not keep track of the dependence of constants on the specific
geometric aspects of our spacetime, except for some very specific
cases. However, we will always keep track of the dependence of all
constants on each parameter that has not been fixed. Once a parameter
is fixed (which will be clearly stated in the text), we will feel
free to drop the dependence of constants on it. 

The notation $f_{1}\lesssim f_{2}$ for two real functions $f_{1},f_{2}$
will as usual mean that there exists some $C>0$, such that $f_{1}\le C\cdot f_{2}$.
Of course, it should be stated clearly in each case whether this constant
$C$ depends on any free parameters. If nothing is stated regarding
the dependence of this constant on parameters, it should be assumed
that it only depends on the geometry of the background under consideration.

We will also write $f_{1}\sim f_{2}$ if $f_{1}\lesssim f_{2}$ and
$f_{2}\lesssim f_{1}$. Moreover, $f_{1}\ll f_{2}$ will mean that
the quotient $\frac{|f_{1}|}{|f_{2}|}$ can be bounded from above
by some sufficiently small positive constant, the magnitude and the
dependence of which on variable parameters will be clear in each case
from the context. Furthermore, for any function $f:\mathcal{M}\rightarrow[0,+\infty)$
defined on some set $\mathcal{M}$, we will denote with $\{f\gg1\}$
the subset $\{f\ge R\}$ of $\mathcal{M}$ for some constant $R\gg1$.

For functions $f_{1},f_{2}$ of some variable $x$ taking values in
a semi-infinite interval $[a,+\infty)$, writing $f_{1}=o(f_{2})$
will imply that $\frac{f_{1}}{f_{2}}$ can be bounded by some continuous
function $h:[a,+\infty)\rightarrow\mathbb{R}_{+}$ such that $h(x)\rightarrow0$
as $x\rightarrow+\infty$. Again, the dependence of this bound $h$
on any free parameter will be clear from the context.

\subsection{Convention on connections and volume form notations}

We will frequently denote the natural connection of a pseudo-Riemannian
manifold $(\mathcal{N},h_{\mathcal{N}})$ as $\nabla^{h_{\mathcal{N}}}$
or $\nabla_{h_{\mathcal{N}}}$, and the associated volume form as
$dh_{\mathcal{N}}$. If $h_{\mathcal{N}}$ is Riemannian, we will
denote the associated norm on $\oplus_{n,m\in\mathbb{N}}\Big(\otimes^{n}T\mathcal{N}\otimes^{m}T^{*}\mathcal{N}\Big)$
with $\big|\cdot\big|_{h_{\mathcal{N}}}$. 

For any integer $j\ge0$, $\big(\nabla^{h_{\mathcal{N}}}\big)^{j}$
or $\nabla_{h_{\mathcal{N}}}^{j}$ will as usual denote the higher
order operator 
\begin{equation}
\underbrace{\nabla_{h_{\mathcal{N}}}\cdots\nabla_{h_{\mathcal{N}}}}_{j\mbox{ times}}.
\end{equation}
 Notice that the above product is not symmetrised. We will always
use Latin characters to denote such powers of covariant derivative
operators, while Greek characters will be used for the indices of
a tensor in an abstract index notation. 
\begin{example*}
Under these conventions, for a $(n,m)$-tensor $k$ and a function
$u$ on a pseudo-Riemannian manifold $(\mathcal{N},h_{\mathcal{N}})$,
the quantity 
\[
k_{\text{\textgreek{b}}_{1}\ldots\text{\textgreek{b}}_{m}}^{\text{\textgreek{a}}_{1}\ldots\text{\textgreek{a}}_{n}}\cdot\big(\nabla_{h_{\mathcal{N}}}^{n+m}\big)_{\hphantom{\text{\textgreek{b}}_{1}\ldots\text{\textgreek{b}}_{m}}\text{\textgreek{a}}_{1}\ldots\text{\textgreek{a}}_{n}}^{\text{\textgreek{b}}_{1}\ldots\text{\textgreek{b}}_{m}}u
\]
 denotes a contraction of $k$ with the higher order derivative $\nabla_{h_{\mathcal{N}}}^{n+m}u$
of $u$, where the metric $h_{\mathcal{N}}$ was used to raise the
first $m$ indices of $\nabla_{h_{\mathcal{N}}}^{n+m}u$. Notice that
the abstract index notation used above is independent of the choice
of an underlying coordinate chart for the indices.
\end{example*}

\subsection{Conventions on notations for derivatives on $\mathbb{S}^{d-1}$}

In this paper we will frequently work in polar coordinates, and hence
it will prove convenient to introduce some shorthand notation regarding
iterated derivatives on the unit sphere $\mathbb{S}^{d-1}$, $d\ge2$.

We will denote with $g_{\mathbb{S}^{d-1}}$ the usual round metric
on the sphere $\mathbb{S}^{d-1}$, which is the induced metric on
the unit sphere of $\mathbb{R}^{d}$. The metric $g_{\mathbb{S}^{d-1}}$
extends naturally to an inner product on the tensor bundle $\oplus_{n,m\in\mathbb{N}}\Big(\otimes^{n}T\mathbb{S}^{d-1}\otimes^{m}T^{*}\mathbb{S}^{d-1}\Big)$
over $\mathbb{S}^{d-1}$; we will still denote this inner product
as $g_{\mathbb{S}^{d-1}}$. For any tensor field $\mathfrak{m}$ on
$\mathbb{S}^{d-1}$, $|\mathfrak{m}|_{g_{\mathbb{S}^{d-1}}}$ will
as usual denote the norm of $\mathfrak{m}$ with respect to $g_{\mathbb{S}^{d-1}}$.
We will also denote with $\nabla^{\mathbb{S}^{d-1}}$ or $\nabla_{\mathbb{S}^{d-1}}$
the covariant derivative associated with $g_{\mathbb{S}^{d-1}}$.
We will denote with $\text{\textgreek{D}}_{g_{\mathbb{S}^{d-1}}}$
the Laplace--Beltrami operator on $(\mathbb{S}^{d-1},g_{\mathbb{S}^{d-1}})$.
For any smooth $(k_{1},k_{2})$-tensor field $\mathfrak{m}$ on $\mathbb{S}^{d-1}$,
$\big(\nabla^{\mathbb{S}^{d-1}}\big)^{k}\mathfrak{m}$ (or $\nabla_{\mathbb{S}^{d-1}}^{k}\mathfrak{m}$)
will denote the $(k_{1},k_{2}+k)$-tensor field on $\mathbb{S}^{d-1}$
obtained after applying the operator $\nabla^{\mathbb{S}^{d-1}}$
on $\mathfrak{m}$ $k$ times.

We will frequently work on regions $\mathcal{U}$ of a spacetime $\mathcal{M}^{d+1}$
where there exists a coordinate ``chart'' mapping $\mathcal{U}$
diffeomorphically onto $\mathbb{R}_{+}\times\mathbb{R}_{+}\times\mathbb{S}^{d-1}$.
In any such a coordinate ``chart'', $\text{\textgreek{sv}}$ will
denote the projection $\text{\textgreek{sv}}:\mathcal{U}\rightarrow\mathbb{S}^{d-1}$.
Notice that with this notation, for any $x\in\mathcal{M}$, $\text{\textgreek{sv}}(x)$
is a point on $\mathbb{S}^{d-1}$ and not just the coordinates of
this point in a coordinate chart on $\mathbb{S}^{d-1}$. We will \emph{not
}need to fix a coordinate atlas on $\mathbb{S}^{d-1}$. We will use
the same $\text{\textgreek{sv}}$ notation also for the spherical
variable of a polar coordinate ``chart'' on codimension $1$ submanifolds
of $\mathcal{M}$ (the range of such a ``chart'' will be $\mathbb{R}_{+}\times\mathbb{S}^{d-1}$).
For instance, $(r,\text{\textgreek{sv}}):\{t=0\}\rightarrow\mathbb{R}_{+}\times\mathbb{S}^{d-1}$
will denote the usual polar coordinate chart on the hyperplane $\{t=0\}$
of $\mathbb{R}^{d+1}$. 

For any function $h$ on a subset $\mathcal{U}$ of a spacetime $\mathcal{M}$
covered by a polar coordinate chart $(u_{1},u_{2},\text{\textgreek{sv}}):\mathcal{U}\rightarrow\mathbb{R}_{+}\times\mathbb{R}_{+}\times\mathbb{S}^{d-1}$
and any $\text{\textgreek{a}}_{1},\text{\textgreek{a}}_{2}\in\mathbb{R}_{+}$,
$h(\text{\textgreek{a}}_{1},\text{\textgreek{a}}_{2},\cdot)$ defines
a function on $\mathbb{S}^{d-1}$. In this way, the $\nabla^{\mathbb{S}^{d-1}}$
differential operator on $\mathbb{S}^{d-1}$ is extended to a tangential
differential operator on the hypersurfaces $\{u_{1},u_{2}=const\}$
of $\mathcal{U}$. Notice, of course, that this operator is tied to
the specific choice of the polar coordinate chart $(u_{1},u_{2},\text{\textgreek{sv}})$.

We will now introduce some schematic notation for derivatives on $\mathbb{S}^{d-1}$
(and the associated tangential operators on the hypersurfaces $\{u_{1},u_{2}=const\}$
in a $(u_{1},u_{2},\text{\textgreek{sv}})$ coordinate chart on a
spacetime $\mathcal{M}$). For any function $h$ on $\mathbb{S}^{d-1}$
and any $k\in\mathbb{N}$, we will frequently denote the $k$-th order
derivative $\nabla_{\mathbb{S}^{d-1}}^{k}h$ as $\partial_{\text{\textgreek{sv}}}^{k}h$,
and we will also use the following notation for the norm of this tensor:
\begin{equation}
|\partial_{\text{\textgreek{sv}}}^{k}h|\doteq\big|\nabla_{\mathbb{S}^{d-1}}^{k}h\big|_{\mathbb{S}^{d-1}}.
\end{equation}
 Moreover, for any \underline{symmetric} $(k,0)$-tensor $a$ on
$\mathbb{S}^{d-1}$, we will use the following schematic notation
for the contraction of $\big(\nabla^{\mathbb{S}^{d-1}}\big)^{k}h$
with $a$:
\begin{equation}
a\cdot\partial_{\text{\textgreek{sv}}}^{k}h\doteq a^{\text{\textgreek{i}}_{1}\ldots\text{\textgreek{i}}_{k}}(\nabla_{\mathbb{S}^{d-1}}^{k})_{\text{\textgreek{i}}_{1}\ldots\text{\textgreek{i}}_{k}}h\label{eq:ContractionOneFunction}
\end{equation}
(see the previous section for the notations on powers of covariant
derivatives and the abstract index notation). We will use the same
notation for the contraction of the product of derivatives of two
or more functions: For any set of $n$ functions $h_{1},\ldots,h_{n}$
on $\mathbb{S}^{d-1}$ and any set $(j_{1},\ldots j_{n})$ of non
negative integers, for any $(\sum_{k=1}^{n}j_{k},0$)-tensor $a$
on $\mathbb{S}^{d-1}$ which is symmetric in any pair of indices lying
in the same one of the intervals $I_{m}=\Big(\sum_{k=1}^{m-1}j_{k}+1,\sum_{k=1}^{m}j_{k}\Big)$
for each $m\in\{1,\ldots n\}$ (but not necessarily symmetric in pairs
of indices lying in different $I_{m}$'s), we will denote 
\begin{equation}
a\cdot\partial_{\text{\textgreek{sv}}}^{j_{1}}h_{1}\cdots\partial_{\text{\textgreek{sv}}}^{j_{n}}h_{n}\doteq a^{\text{\textgreek{i}}_{1}\ldots\text{\textgreek{i}}_{\sum_{k=1}^{n}j_{k}}}\cdot\big(\nabla_{\mathbb{S}^{d-1}}^{j_{1}}\big)_{\text{\textgreek{i}}_{1}\ldots\text{\textgreek{i}}_{j_{1}}}h_{1}\cdots\big(\nabla_{\mathbb{S}^{d-1}}^{j_{n}}\big)_{\text{\textgreek{i}}_{\sum_{k=1}^{n-1}j_{k}+1}\ldots\text{\textgreek{i}}_{\sum_{k=1}^{n}j_{k}}}h_{n}.\label{eq:ContractionProductFunctions}
\end{equation}

The same notation (\ref{eq:ContractionOneFunction}) and (\ref{eq:ContractionProductFunctions})
will also apply when $h$, $h_{1},\ldots,h_{n}$ are tensor fields
on $\mathbb{S}^{d-1}$.

Depending on the context, $d\text{\textgreek{sv}}$ will be used to
denote either the usual volume form on $(\mathbb{S}^{d-1},g_{\mathbb{S}^{d-1}}$)
or a $1$-form on $\mathbb{S}^{d-1}$ satisfying for any $k\in\mathbb{N}$
the bound $\big|\big(\nabla^{\mathbb{S}^{d-1}}\big)^{k}d\text{\textgreek{sv}}\big|_{g_{\mathbb{S}^{d-1}}}\le10^{k}$.
Similarly, $d\text{\textgreek{sv}}d\text{\textgreek{sv}}$ will denote
a symmetric $(2,0)$-tensor on $\mathbb{S}^{d-1}$ satisfying for
any $k\in\mathbb{N}$ the bound $\big|\big(\nabla^{\mathbb{S}^{d-1}}\big)^{k}(d\text{\textgreek{sv}}d\text{\textgreek{sv}})\big|_{g_{\mathbb{S}^{d-1}}}\le10^{k}$.

As an example, the above notation will allow us to perform the following
integration by parts procedure on $\mathbb{S}^{d-1}$ for any function
$f$ and any tensor $a$ with the aforementioned symmetries: 
\begin{equation}
\int_{\mathbb{S}^{d-1}}a\cdot\partial_{\text{\textgreek{sv}}}f\cdot\partial_{\text{\textgreek{sv}}}\partial_{\text{\textgreek{sv}}}f\, d\text{\textgreek{sv}}=-\frac{1}{2}\int_{\mathbb{S}^{d-1}}(\tilde{e}_{1}\partial_{\text{\textgreek{sv}}}a+\tilde{e}_{2}a)\cdot\partial_{\text{\textgreek{sv}}}f\cdot\partial_{\text{\textgreek{sv}}}f\, d\text{\textgreek{sv}},\label{eq:Integrationbypartsspheregeneralrule-1}
\end{equation}
 for some smooth contracting tensors $\tilde{e}_{1}$, $\tilde{e}_{2}$
which are bounded (in any $C^{k}$ norm) with bounds depending only
on the tensor type of $a$. 

We will frequently use the notation (\ref{eq:ContractionOneFunction})
and (\ref{eq:ContractionProductFunctions}) in cases where we do not
have an explicit form for the contracting tensor $a$, but we merely
have bounds for the norm of $a$ and its derivatives. It is for this
reason that we choose to use a notation which apparently loses information
regarding the structure of the underlying expression.

Notice that in a polar coordinate chart $(u_{1},u_{2},\text{\textgreek{sv}}):\mathcal{U}\rightarrow\mathbb{R}_{+}\times\mathbb{R}_{+}\times\mathbb{S}^{d-1}$,
the following commutation relation holds for any function $h$: 
\begin{equation}
\mathcal{L}_{\partial_{u_{i}}}\nabla^{\mathbb{S}^{d-1}}h=\nabla^{\mathbb{S}^{d-1}}\partial_{u_{i}}h,
\end{equation}
where $\partial_{u_{i}}$ is the coordinate vector field associated
to the coordinate function $u_{i}$, $i=1,2$. Therefore, we will
frequently denote 
\begin{equation}
\mathcal{L}_{\partial_{u_{i}}}\nabla^{\mathbb{S}^{d-1}}h\doteq\partial_{u_{i}}\partial_{\text{\textgreek{sv}}}h,
\end{equation}
and this will allow us to commute $\partial_{u_{i}}$ with $\partial_{\text{\textgreek{sv}}}$,
as if $\partial_{\text{\textgreek{sv}}}$ was a regular coordinate
vector field.

\subsection{Convention on the $O(\cdot)$ notation}

For the conventions regarding the use the $O(\cdot)$ notation on
asymptotically flat spacetimes in this paper, see the beginning Section
\ref{sec:GeometrySpacetimes}.

\subsection{Conventions on integration}

When we integrate over open subsets of a Lorentzian manifold $(\mathcal{M},g)$
using the natural volume form $\text{\textgreek{w}}$ associated to
$g$, we will often drop the volume form in the expression for the
integral. Recall that $\text{\textgreek{w}}$ is expressed as 
\[
\text{\textgreek{w}}=\sqrt{-det(g)}dx^{0}\cdots dx^{d}
\]
 in any local coordinate chart $(x^{0},x^{1},x^{2},\ldots,x^{d})$.
The same rule will apply when integrating over any spacelike hypersurface
$\mathcal{S}$ of $(\mathcal{M},g)$ using the natural volume form
of its induced Riemannian metric.

\subsection{\noindent Notations for vector field multipliers and currents:}

Since we will use the language of currents and vector field multipliers
in order to establish the desired estimates, let us briefly recall
first the required notation: On any Lorentzian manifold $(\mathcal{M},g)$,
associated to the wave operator $\square_{g}=\frac{1}{\sqrt{-det(g)}}\partial_{\text{\textgreek{m}}}\Big(\sqrt{-det(g)}\cdot g^{\text{\textgreek{m}\textgreek{n}}}\partial_{\text{\textgreek{n}}}\Big)$
is the energy momentum tensor $T$, which for any smooth function
$\text{\textgreek{y}}:\mathcal{M}\rightarrow\mathbb{C}$ takes the
form

\begin{equation}
T_{\text{\textgreek{m}\textgreek{n}}}(\text{\textgreek{y}})=\frac{1}{2}\Big(\partial_{\text{\textgreek{m}}}\text{\textgreek{y}}\cdot\partial_{\text{\textgreek{n}}}\bar{\text{\textgreek{y}}}+\partial_{\text{\textgreek{m}}}\bar{\text{\textgreek{y}}}\cdot\partial_{\text{\textgreek{n}}}\text{\textgreek{y}}\Big)-\frac{1}{2}\big(\partial^{\text{\textgreek{l}}}\text{\textgreek{y}}\cdot\partial_{\text{\textgreek{l}}}\bar{\text{\textgreek{y}}}\big)g_{\text{\textgreek{m}\textgreek{n}}}.
\end{equation}

Given any continuous and piecewise $C^{1}$ vector field $X$ on $\mathcal{M}$,
we can define almost everywhere the associated currents

\begin{equation}
J_{\text{\textgreek{m}}}^{X}(\text{\textgreek{y}})=T_{\text{\textgreek{m}\textgreek{n}}}(\text{\textgreek{y}})X^{\text{\textgreek{m}}},
\end{equation}
\begin{equation}
K^{X}(\text{\textgreek{y}})=T_{\text{\textgreek{m}\textgreek{n}}}(\text{\textgreek{y}})\nabla^{\text{\textgreek{m}}}X^{\text{\textgreek{n}}}.
\end{equation}
 The following divergence identity then holds almost everywhere:

\begin{equation}
\nabla^{\text{\textgreek{m}}}J_{\text{\textgreek{m}}}^{X}(\text{\textgreek{y}})=K^{X}(\text{\textgreek{y}})+Re\Big\{(\square_{g}\text{\textgreek{y}})\cdot X\bar{\text{\textgreek{y}}}\Big\}.
\end{equation}

\section{\label{sec:GeometrySpacetimes}Geometry of the asymptotically flat
regions $(\mathcal{N}_{af},g)$}

In Sections \ref{sec:Morawetz}, \ref{sec:The-new-method} and \ref{sec:The-improved--hierarchy},
we will work in $d+1$ dimensional, smooth and time oriented Lorentzian
manifolds $(\mathcal{N}_{af}^{d+1},g)$ diffeomorphic to $\mathbb{R}\times(\mathbb{R}^{d}\backslash B_{R})$
for $d\ge3$ (where $B_{R}$ is the closed Euclidian ball of radius
$R$). The manifolds $(\mathcal{N}_{af},g)$ will serve as models
of the asymptotically flat region of more general asymptotically flat
spacetimes, and will appear as open subsets of the Lorentzian manifolds
studied in Sections \ref{sec:FriedlanderRadiation}, \ref{sec:Firstdecay}
and \ref{sec:Improved-polynomial-decay}. Let us notice that $(\mathcal{N}_{af},g)$
will not in general be globally hyperbolic.

We will fix a global coordinate chart $(t,r,\text{\textgreek{sv}}):\mathcal{N}_{af}\rightarrow\mathbb{R}\times(\mathbb{R}^{d}\backslash B_{R})$,
where $t\in\mathbb{R}$ is the projection to the first factor of $\mathbb{R}\times(\mathbb{R}^{d}\backslash B_{R})$,
and $(r,\text{\textgreek{sv}})\in(\mathbb{R}^{+},\mathbb{S}^{d-1})$
are the usual polar coordinates on $\mathbb{R}^{d}\backslash B_{R}$.
Moreover, it will be useful for us to define the coordinate function
$u=t-r$, and introduce the coordinate system $(u,r,\text{\textgreek{sv}})$.

In the $(u,r,\text{\textgreek{sv}})$ coordinate system, we will adopt
the following notation for the derivatives of functions $h$ on $\mathcal{N}_{af}$:
We will write 
\begin{equation}
h=O(r^{b})
\end{equation}
for some $b\in\mathbb{R}$ if for any integer $k\ge0$ we can bound:
\begin{equation}
\sum_{k_{1}+k_{2}+k_{3}=k}\big|r^{k_{1}}\partial_{r}^{k_{1}}\partial_{\text{\textgreek{sv}}}^{k_{2}}\partial_{u}^{k_{3}}h\big|\le C_{k}\cdot r^{b}\label{eq:BigONotation}
\end{equation}
 for some constant $C_{k}$ depending only on $k$ and the function
$h$ itself. By replacing the coordinate derivatives in (\ref{eq:BigONotation})
with the connection derivatives in the associated directions, the
same $O(\cdot)$ notation will apply when $h$ takes values in some
vector bundle $\mathcal{E}$ over $\mathcal{N}_{af}$ with a fixed
Hermitian metric to measure the size of the norms in the left hand
side of (\ref{eq:BigONotation}) and a fixed compatible connection.
For instance, the case when $\mathcal{E}=\otimes^{n_{1}}T\mathbb{S}^{d-1}\otimes^{n_{2}}T^{*}\mathbb{S}^{d-1}$
(equipped with the natural metric and connection arising from the
standard round $g_{\mathbb{S}^{d-1}}$) will appear in the text.

In the coordinate chart $(u,r,\text{\textgreek{sv}})$ on $\mathcal{N}_{af}$,
the Lorentzian metric $g$ will be assumed to take the following radiative
form for some $0\le a\le1$: 
\begin{align}
g=-4\Big(1- & \frac{2M(u,\text{\textgreek{sv}})}{r}+O(r^{-1-a})\Big)du^{2}-\Big(4+O(r^{-1-a})\Big)dudr+r^{2}\cdot\Big(g_{\mathbb{S}^{d-1}}+h_{\mathbb{S}^{d-1}}\Big)+\label{eq:MetricUR}\\
 & +\Big(h_{3}^{as}(u,\text{\textgreek{sv}})+O(r^{-a})\Big)dud\text{\textgreek{sv}}+O(r^{-a})drd\text{\textgreek{sv}}+O(r^{-2-a})dr^{2},\nonumber 
\end{align}
where $M(u,\text{\textgreek{sv}})$ and $h_{3}^{as}(u,\text{\textgreek{sv}})$
are real functions on $\mathbb{R}\times\mathbb{S}^{d-1}$ with all
their derivatives uniformly bounded and $h_{\mathbb{S}^{d-1}}$ is
a symmetric $(2,0)$-tensor on $\mathbb{S}^{d-1}$ satisfying the
bound $h_{\mathbb{S}^{d-1}}=O(r^{-1})$.

Notice that this class of metrics includes the Bondi radiating spacetimes
(see i.\,e. \cite{Bondi1962,Sachs1962}). Due to the form (\ref{eq:MetricUR})
of the metric, the vector field $\partial_{r}$ in the $(u,r,\text{\textgreek{sv}})$
coordinate system is almost null, but not necessarily null, and $u$
is not necessarily an optical function. 

It will also be convenient to express the metric in an \emph{almost}
double null coordinate system $(u,v,\text{\textgreek{sv}})$, where
$v=u+r$. We easily calculate from (\ref{eq:MetricUR}) that in the
$(u,v,\text{\textgreek{sv}})$ coordinate system, the metric has the
following form: 
\begin{align}
g=-\Big(4+ & O(r^{-1-a})\Big)dvdu+r^{2}\cdot\Big(g_{\mathbb{S}^{d-1}}+h_{\mathbb{S}^{d-1}}\Big)+\Big(h^{as}(u,\text{\textgreek{sv}})+O(r^{-a})\Big)dud\text{\textgreek{sv}}+\label{eq:MetricUV}\\
 & +O(r^{-a})dvd\text{\textgreek{sv}}+4\Big(-\frac{2M(u,\text{\textgreek{sv}})}{r}+O(r^{-1-a})\Big)du^{2}+O(r^{-2-a})dv^{2},\nonumber 
\end{align}
 where $h_{\mathbb{S}^{d-1}}=O(r^{-1})$.

Notice that differentiation with respect to $\partial_{v}$ in the
$(u,v,\text{\textgreek{sv}})$ coordinate system is the same as differentiation
with respect to $\partial_{r}$ in the $(u,r,\text{\textgreek{sv}})$
chart. Notice also that differentiation with respect to $\partial_{u}$
is not the same in the two coordinate systems.

Let us also define the function $t\doteq2u+r=u+v$, which has spacelike
level sets, at least for $r\gg1$. In the $(t,r,\text{\textgreek{sv}})$
coordinate system, we easily calculate from (\ref{eq:MetricUR}) that
the metric $g$ has the following expression: 
\begin{align}
g=-\Big(1- & \frac{2M(u,\text{\textgreek{sv}})}{r}+O(r^{-1-a})\Big)\cdot dt^{2}-\Big(\frac{4M(u,\text{\textgreek{sv}})}{r}+O(r^{-1-a})\Big)dtdr+\Big(1+\frac{2M(u,\text{\textgreek{sv}})}{r}+O(r^{-1-a})\Big)dr{}^{2}+\label{eq:metricTR}\\
 & +r^{2}\cdot\Big(g_{\mathbb{S}^{d-1}}+O(r^{-1})\Big)+\Big(-\frac{1}{2}h^{as}(u,\text{\textgreek{sv}})+O(r^{-a})\Big)(dr-dt)d\text{\textgreek{sv}}+O(r^{-a})dtd\text{\textgreek{sv}}.\nonumber 
\end{align}

Notice that the vector field $T=\partial_{t}$ in the $(t,r,\text{\textgreek{sv}})$
coordinate chart will not necessarily be a Killing vector field for
$\mathcal{N}_{af}$, but it will certainly be timelike for $r\gg1$.
Notice also that in this coordinate system, differentiation with respect
to either $\partial_{t}$ or $\partial_{r}$ does not improve the
decay rate in $r$ of $O(r^{-c})$ functions ($\partial_{r}$ in this
system is different than $\partial_{r}$ in the $(u,r,\text{\textgreek{sv}})$
chart).

It is important to remark that if $R\gg1$, in the region $\{r\ge R\}$
of $\mathcal{N}_{af}$ we can estimate $dudvd\text{\textgreek{sv}}\sim r^{-(d-1)}\cdot dvol_{g}$.
In this region, we will also use the notation 
\begin{equation}
|\partial\text{\textgreek{f}}|^{2}\doteq|\partial_{v}\text{\textgreek{f}}|^{2}+|\partial_{u}\text{\textgreek{f}}|^{2}+|r^{-1}\partial_{\text{\textgreek{sv}}}\text{\textgreek{f}}|^{2}
\end{equation}
 for the ``coordinate Euclidean norm'' of the gradient of any differentiable
function $\text{\textgreek{f}}:\mathcal{N}_{af}\rightarrow\mathbb{C}$.

\subsection{\label{sub:Spacelike-hyperboloidal-hypersurfaces}Spacelike hyperboloidal
hypersurfaces terminating at $\mathcal{I}^{+}$}

It will be convenient to introduce a family of spacelike hypersurfaces
``terminating at future null infinity'' (in a sense that will be
made precise shortly). These hypersurfaces will make easier the extraction
of information regarding the radiating properties of solutions $\text{\textgreek{f}}$
to the wave equation $\square_{g}\text{\textgreek{f}}=0$ on $(\mathcal{N}_{af},g)$.

We first introduce the following definition:
\begin{defn}
We define the future null infinity $\mathcal{I}^{+}$ of $(\mathcal{N}_{af},g)$
as the abstract limit of the hypersurfaces $\{v=v_{n}\}$ as $v_{n}\rightarrow+\infty$
in the $(u,v,\text{\textgreek{sv}})$ coordinate system. In particular,
a function $\text{\textgreek{Y}}$ on $\mathcal{I}^{+}$ will always
be defined as the limit $\lim_{v\rightarrow+\infty}\tilde{\text{\textgreek{Y}}}$
in the $(u,v,\text{\textgreek{sv}})$ coordinate system for some function
$\tilde{\text{\textgreek{Y}}}$ on $\mathcal{N}_{af}$. We will also
set $\mathcal{I}^{+}(\text{\textgreek{t}}_{1},\text{\textgreek{t}}_{2})=\mathcal{I}^{+}\cap\{\text{\textgreek{t}}_{1}\le u\le t_{2}\}$,
which is to be understood as the abstract limit of the hypersurfaces
$\{v=v_{n}\}\cap\{\text{\textgreek{t}}_{1}\le u\le\text{\textgreek{t}}_{2}\}$
as $v_{n}\rightarrow+\infty$ in the $(u,v,\text{\textgreek{sv}})$
coordinate system. 
\end{defn}
We will now define the notion of a hypersurface terminating at $\mathcal{I}^{+}$:
\begin{defn}
\label{DefinitionTerminatingAtI+}Let $\mathcal{S}$ be an achronal
inextendible hypersurface of $(\mathcal{N}_{af},g)$. We will say
that $\mathcal{S}$ terminates at future null infinity $\mathcal{I}^{+}$
if the coordinate function $u$ restricted on $\mathcal{S}\cap\{r\gg1\}$
is bounded.
\end{defn}
Notice that due to the form (\ref{eq:MetricUR}) of the metric in
the $(u,r,\text{\textgreek{sv}})$, and in particular due to the fact
that in this coordinate system $g(\partial_{r},\partial_{r})=O(r^{-2-a})$
and $g(\partial_{r},\partial_{\text{\textgreek{sv}}})=O(r^{-a})$,
we can easily calculate that for any $R\gg1$ and $\text{\textgreek{t}}\in\mathbb{R}$
the boundary hypersurface $\partial J^{+}\big(\{r\le R\}\cap\{t=\text{\textgreek{t}}\}\big)$
is an achronal hypersurface terminating at future null infinity. 

We can now construct a function $\bar{t}:\mathcal{N}_{af}\rightarrow\mathbb{R}$
with level sets that are spacelike hypersurfaces%
\footnote{at least for $r\gg1$%
} terminating at future null infinity. More precisely, we define for
any fixed $0<\text{\textgreek{h}}'<1+a$

\noindent 
\begin{equation}
\bar{t}_{\text{\textgreek{h}}'}\doteq u-\frac{1}{1+r^{\text{\textgreek{h}'}}}.\label{eq:Hyperboloids}
\end{equation}

\noindent The level sets of $\bar{t}_{\text{\textgreek{h}}'}$ are
spacelike hypersurfaces for $r$ large enough depending on \textgreek{h}'.
This follows from the computation:
\begin{align*}
g^{\text{\textgreek{m}\textgreek{n}}}\partial_{\text{\textgreek{m}}}\bar{t}_{\text{\textgreek{h}}'}\cdot\partial_{\text{\textgreek{n}}}\bar{t}_{\text{\textgreek{h}}'} & =g^{\text{\textgreek{m}\textgreek{n}}}\partial_{\text{\textgreek{m}}}u\cdot\partial_{\text{\textgreek{n}}}u+2g^{\text{\textgreek{m}\textgreek{n}}}\partial_{\text{\textgreek{m}}}u\cdot\partial_{\text{\textgreek{n}}}(\frac{-1}{1+r^{\text{\textgreek{h}'}}})+g^{\text{\textgreek{m}\textgreek{n}}}\partial_{\text{\textgreek{m}}}(\frac{-1}{1+r^{\text{\textgreek{h}'}}})\cdot\partial_{\text{\textgreek{n}}}(\frac{-1}{1+r^{\text{\textgreek{h}'}}})=\\
 & =-2\text{\textgreek{h}}'r^{-1-\text{\textgreek{h}}'}+O_{\text{\textgreek{h}}'}(r^{-2-a}+r^{-2-\text{\textgreek{h}'}})<0.
\end{align*}
 for $r$ large enough in terms of $\text{\textgreek{h}}'$. 

Moreover, $|\bar{t}_{\text{\textgreek{h}}'}-u|\le1$, and hence the
level sets of $\bar{t}_{\text{\textgreek{h}}'}$ are spacelike hypersurfaces
terminating at future null infinity, according to the definition \ref{DefinitionTerminatingAtI+}. 

We will frequently simply write $\bar{t}$ in place of $\bar{t}_{\text{\textgreek{h}}'}$,
since $\text{\textgreek{h}}'$ will be considered fixed throughout
the next Sections.

For any $\text{\textgreek{t}}\in\mathbb{R}$ and $\text{\textgreek{t}}_{1}\le\text{\textgreek{t}}_{2}\in\mathbb{R}$,
we will denote 
\begin{equation}
\mathcal{S}_{\text{\textgreek{t}}}\doteq\{\bar{t}=\text{\textgreek{t}}\}\cap\{r\ge R_{3}\}\label{eq:HyperboloidNotation}
\end{equation}
 and 
\begin{equation}
\mathcal{R}(\text{\textgreek{t}}_{1},\text{\textgreek{t}}_{2})\doteq\{\text{\textgreek{t}}_{1}\le\bar{t}\le\text{\textgreek{t}}_{2}\}\cap\{r\ge R_{3}\},\label{eq:HyperboloidalRegion}
\end{equation}
where $R_{3}=R_{3}(\text{\textgreek{h}}')$ is fixed large enough
so that in the region $\{r\ge R_{3}\}$ the level sets of $\bar{t}_{\text{\textgreek{h}}'}$
are spacelike. 

Due to the definition of $\bar{t}_{\text{\textgreek{h}}'}$, the image
of $\mathcal{S}_{\text{\textgreek{t}}_{1}}$ under the flow of the
vector field $T=\partial_{t}$ in the $(t,r,\text{\textgreek{sv}})$
coordinate system for time $\text{\textgreek{t}}$ is precisely $\mathcal{S}_{\text{\textgreek{t}}_{1}+\text{\textgreek{t}}}$. 

Note that the $\mathcal{S}_{\text{\textgreek{t}}}$ can be regularly
parametrized both by $(u,\text{\textgreek{sv}})$ and by $(v,\text{\textgreek{sv}})$,
and the corresponding coordinate volume forms $dud\text{\textgreek{sv}}$
and $dvd\text{\textgreek{sv}}$ satisfy the relation 
\begin{equation}
dud\text{\textgreek{sv}}\sim_{\text{\textgreek{h}}'}r^{-1-\text{\textgreek{h}}'}dvd\text{\textgreek{sv}}.
\end{equation}

\subsection{Expression of the wave equation in the coordinate chart $(u,v,\text{\textgreek{sv}})$}

\textgreek{I}n the coordinate chart $(u,v,\text{\textgreek{sv}})$
on $\mathcal{N}_{af}$, the metric $g$ takes the form (\ref{eq:MetricUV}),
and we can easily compute that in these coordinates: 
\begin{equation}
det(g)=-4r^{2(d-1)}\cdot\big(1+O(r^{-1})\big)
\end{equation}
and the inverse of the metric has the following form: 
\begin{align}
g^{-1}=-\Big(1 & +O(r^{-1-a})\Big)\partial_{u}\partial_{v}+r^{-2}\Big(g_{\mathbb{S}^{d-1}}^{-1}-h_{\mathbb{S}^{d-1}}^{inv}\Big)+r^{-2}\Big(\frac{1}{2}h^{as}(u,\text{\textgreek{sv}})+O(r^{-a})\Big)\partial_{v}\partial_{\text{\textgreek{sv}}}+\label{eq:InverseMetricUV}\\
 & +O(r^{-2-a})\partial_{u}\partial_{\text{\textgreek{sv}}}+\Big(-\frac{2M(u,\text{\textgreek{sv}})}{r}+O(r^{-1-a})\Big)\partial_{v}^{2}+O(r^{-2-a})\partial_{u}^{2},\nonumber 
\end{align}
where $h_{\mathbb{S}^{d-1}}^{inv}=O(r^{-1})$.

Setting 
\begin{equation}
\text{\textgreek{W}}\doteq\big(-\frac{1}{4}det(g)\big)^{\frac{1}{4}}=r^{\frac{d-1}{2}}\big(1+O(r^{-1})\big),\label{eq:ConformalFactorOmega}
\end{equation}
the wave operator then takes the following form: 

\begin{align}
\text{\textgreek{W}}\cdot\square\text{\textgreek{f}}=-\big(1 & +O(r^{-1-a})\big)\cdot\partial_{u}\partial_{v}(\text{\textgreek{W}}\text{\textgreek{f}})+r^{-2}\text{\textgreek{D}}_{g_{\mathbb{S}^{d-1}}+h_{\mathbb{S}^{d-1}}}(\text{\textgreek{W}}\text{\textgreek{f}})-\label{eq:ConformalWaveOperator}\\
 & -\frac{(d-1)(d-3)}{4}r^{-2}\cdot(\text{\textgreek{W}}\text{\textgreek{f}})+Err(\text{\textgreek{W}}\text{\textgreek{f}}),\nonumber 
\end{align}
 where the $Err(\text{\textgreek{F}})$ term is of the form (in our
schematic notation of Section \ref{sub:Notational-conventions}):
\begin{align}
Err(\text{\textgreek{F}})= & O(r^{-2-a})\cdot\partial_{u}^{2}\text{\textgreek{F}}+O(r^{-1})\cdot\partial_{v}^{2}\text{\textgreek{F}}+O(r^{-2-a})\partial_{u}\partial_{\text{\textgreek{sv}}}\text{\textgreek{F}}+O(r^{-2})\partial_{v}\partial_{\text{\textgreek{sv}}}\text{\textgreek{F}}+\label{eq:ErrTerms}\\
 & +O(r^{-3-a})\partial_{\text{\textgreek{sv}}}\partial_{\text{\textgreek{sv}}}\text{\textgreek{F}}+O(r^{-2-a})\partial_{u}\text{\textgreek{F}}+O(r^{-1-a})\partial_{v}\text{\textgreek{F}}+O(r^{-2-a})\cdot\partial_{\text{\textgreek{sv}}}\text{\textgreek{F}}+O(r^{-3})\text{\textgreek{F}}.\nonumber 
\end{align}

Notice the similarity of the expression (\ref{eq:ConformalWaveOperator})
with the expression of the wave operator in the double null coordinates
of Minkowski spacetime.

\section{\label{sec:Morawetz}Some $\partial_{r}$- Morawetz type and energy
boundedness estimates}

In this Section, we will establish some estimates of Morawetz type,
controlling the behaviour of solutions to $\square\text{\textgreek{f}}=F$
on our asymptotically flat model $(\mathcal{N}_{af},g)$. 

Recall that in the $(u,v,\text{\textgreek{sv}})$ coordinate system,
the wave equation takes the form (\ref{eq:ConformalWaveOperator}):
\begin{align}
\text{\textgreek{W}}\cdot\square\text{\textgreek{f}}=-\big(1 & +O(r^{-1-a})\big)\cdot\partial_{u}\partial_{v}(\text{\textgreek{W}}\text{\textgreek{f}})+r^{-2}\text{\textgreek{D}}_{g_{\mathbb{S}^{d-1}}+h_{\mathbb{S}^{d-1}}}(\text{\textgreek{W}}\text{\textgreek{f}})-\label{eq:ConformalWaveOperator-2}\\
 & -\frac{(d-1)(d-3)}{4}r^{-2}\cdot(\text{\textgreek{W}}\text{\textgreek{f}})+Err(\text{\textgreek{W}}\text{\textgreek{f}}),\nonumber 
\end{align}
 where 
\begin{equation}
\text{\textgreek{W}}=r^{\frac{d-1}{2}}\big(1+O(r^{-1})\big)
\end{equation}
and the $Err(\text{\textgreek{F}})$ term is of the form: 
\begin{align}
Err(\text{\textgreek{F}})= & O(r^{-2-a})\cdot\partial_{u}^{2}\text{\textgreek{F}}+O(r^{-1})\cdot\partial_{v}^{2}\text{\textgreek{F}}+O(r^{-2-a})\partial_{u}\partial_{\text{\textgreek{sv}}}\text{\textgreek{F}}+O(r^{-2})\partial_{v}\partial_{\text{\textgreek{sv}}}\text{\textgreek{F}}+\label{eq:ErrorTermsMorawetz}\\
 & +O(r^{-3-a})\partial_{\text{\textgreek{sv}}}\partial_{\text{\textgreek{sv}}}\text{\textgreek{F}}+O(r^{-2-a})\partial_{u}\text{\textgreek{F}}+O(r^{-1-a})\partial_{v}\text{\textgreek{F}}+O(r^{-2-a})\cdot\partial_{\text{\textgreek{sv}}}\text{\textgreek{F}}+O(r^{-3})\text{\textgreek{F}}.\nonumber 
\end{align}

We will also adopt the following convention: We will say that a function
$\text{\textgreek{f}}:N_{af}\rightarrow\mathbb{C}$ has compact support
in space if there exists a continuous function $h:\mathbb{R}\rightarrow\mathbb{R}_{+}$
such that $supp(\text{\textgreek{f}})\subseteq\{r\le h(t)\}$.

\subsection{A first $\partial_{r}$- Morawetz type estimate}

We will establish the following Morawetz-type lemma in the region
$\{r\gg1\}$:
\begin{lem}
\label{lem:MorawetzDrLemmaHyperboloids}For any given $0<\text{\textgreek{h}}<a$,
there exists an $R=R(\text{\textgreek{h}})>0$ and constants $C(\text{\textgreek{h}}),c(\text{\textgreek{h}})>0$
such that for any smooth function $\text{\textgreek{f}}:\mathcal{N}_{af}\rightarrow\mathbb{C}$,
any $\text{\textgreek{t}}_{1}\le\text{\textgreek{t}}_{2}$, any $T^{*}>0$
and any smooth cut-off function $\text{\textgreek{q}}:\mathcal{N}_{af}\rightarrow[0,1]$
supported in $\{r\ge R\}$ we can bound:
\begin{equation}
\begin{split}\int_{\mathcal{R}(\text{\textgreek{t}}_{1},\text{\textgreek{t}}_{2})\cap\{t\le T^{*}\}}\text{\textgreek{q}}\cdot\Big(r^{-1-\text{\textgreek{h}}}\Big(|\partial_{u}\text{\textgreek{f}}|^{2} & +|\partial_{v}\text{\textgreek{f}}|^{2}\Big)+r^{-1}|r^{-1}\partial_{\text{\textgreek{sv}}}\text{\textgreek{f}}|^{2}+r^{-3-\text{\textgreek{h}}}|\text{\textgreek{f}}|^{2}\Big)\le\\
\le & C(\text{\textgreek{h}})\cdot\int_{\{supp(\partial\text{\textgreek{q}})\}\cap\mathcal{R}(\text{\textgreek{t}}_{1},\text{\textgreek{t}}_{2})\cap\{t\le T^{*}\}}|\partial\text{\textgreek{q}}|\cdot\big(|\partial\text{\textgreek{f}}|^{2}+r^{-2}|\text{\textgreek{f}}|^{2}\big)+\\
 & +C(\text{\textgreek{h}})\cdot\sum_{i=1}^{2}\int_{\mathcal{S}_{\text{\textgreek{t}}_{i}}\cap\{t\le T^{*}\}}\text{\textgreek{q}}\Big(|\partial_{v}\text{\textgreek{f}}|^{2}+|r^{-1}\partial_{\text{\textgreek{sv}}}\text{\textgreek{f}}|^{2}+r^{-1-\text{\textgreek{h}}'}|\partial_{u}\text{\textgreek{f}}|^{2}+r^{-2}|\text{\textgreek{f}}|^{2}\Big)\, r^{d-1}dvd\text{\textgreek{sv}}+\\
 & +C(\text{\textgreek{h}})\cdot\int_{\{t=T^{*}\}\cap\mathcal{R}(\text{\textgreek{t}}_{1},\text{\textgreek{t}}_{2})}\text{\textgreek{q}}\big(|\partial\text{\textgreek{f}}|^{2}+r^{-2}|\text{\textgreek{f}}|^{2}\big)\, r^{d-1}dvd\text{\textgreek{sv}}+\\
 & +\int_{\mathcal{R}(\text{\textgreek{t}}_{1},\text{\textgreek{t}}_{2})\cap\{t\le T^{*}\}}\text{\textgreek{q}}\cdot Re\big\{\big(O(1)(\partial_{v}-\partial_{u})\bar{\text{\textgreek{f}}}+O(r^{-1})\bar{\text{\textgreek{f}}}\big)\cdot\square_{g}\text{\textgreek{f}}\big\}+\\
 & +C(\text{\textgreek{h}})\int_{\mathcal{R}(\text{\textgreek{t}}_{1},\text{\textgreek{t}}_{2})\cap\{t\le T^{*}\}}\text{\textgreek{q}}\cdot r^{-1}\Big(|\partial_{v}\text{\textgreek{f}}|^{2}+|r^{-1}\partial_{\text{\textgreek{sv}}}\text{\textgreek{f}}|^{2}+r^{-2}|\text{\textgreek{f}}|^{2}\Big).
\end{split}
\label{eq:MorawetzGeneralCaseRadiativeHyperboloids}
\end{equation}
\end{lem}
\begin{rem*}
In case the radiative components of the metric satisfy the bounds
$\partial_{u}M\le0$ and $|\partial_{u}h_{as}|+|r\partial_{u}h_{\mathbb{S}^{d-1}}|\ll-(\partial_{u}M)+O(r^{-a})$
(which includes the non-radiating case $\partial_{u}M=0$, $h_{as}=0$
and $h_{\mathbb{S}^{d-1}}=O(r^{-1-a})$), the last term of the right
hand side of (\ref{eq:MorawetzGeneralCaseRadiativeHyperboloids})
can be omitted. Furthermore, in case the $T$ vector field  satisfies
(\ref{eq:DeformationTensorTAway}) for $m=1$, then the last term
of the right hand side of (\ref{eq:MorawetzGeneralCaseRadiative})
is replaced by 
\begin{equation}
\int_{\mathcal{R}(\text{\textgreek{t}}_{1},\text{\textgreek{t}}_{2})\cap\{t\le T^{*}\}}\text{\textgreek{q}}\cdot\bar{t}^{-\text{\textgreek{d}}_{0}}r^{-1}\Big(|\partial_{v}\text{\textgreek{f}}|^{2}+|r^{-1}\partial_{\text{\textgreek{sv}}}\text{\textgreek{f}}|^{2}+r^{-2}|\text{\textgreek{f}}|^{2}\Big).
\end{equation}
\end{rem*}
\begin{proof}
Without loss of generality, we can assume that $\text{\textgreek{f}}$
is real valued. 

Let us consider the function $f:[0,+\infty)\rightarrow(0,+\infty)$
defined as: 
\begin{equation}
f(r)\doteq\frac{r^{\text{\textgreek{h}}}}{1+r^{\text{\textgreek{h}}}}.
\end{equation}
Setting $\text{\textgreek{F}}=\text{\textgreek{W}}\text{\textgreek{f}}$
and multiplying equation (\ref{eq:ConformalWaveOperator-2}) with
\[
\text{\textgreek{q}}\cdot f(r)\cdot\big(\partial_{v}-\partial_{u}\big)\text{\textgreek{F}},
\]
 we obtain after integrating over$\mathcal{R}(\text{\textgreek{t}}_{1},\text{\textgreek{t}}_{2})\cap\{t\le T^{*}\}$
(with $dudvd\text{\textgreek{sv}}$ used as a volume form): 
\begin{equation}
\begin{split}\int_{\mathcal{R}(\text{\textgreek{t}}_{1},\text{\textgreek{t}}_{2})\cap\{t\le T^{*}\}}\text{\textgreek{q}}f(r)\cdot\big(\partial_{v}-\partial_{u}\big)\text{\textgreek{F}}\cdot\text{\textgreek{W}}\square\text{\textgreek{f}}\, dudvd\text{\textgreek{sv}}=\\
=\int_{\mathcal{R}(\text{\textgreek{t}}_{1},\text{\textgreek{t}}_{2})\cap\{t\le T^{*}\}}\text{\textgreek{q}}f(r)\cdot\big(\partial_{v}-\partial_{u}\big)\text{\textgreek{F}}\cdot\Big\{ & -\big(1+O_{1}(r^{-1-a})\big)\cdot\partial_{u}\partial_{v}\text{\textgreek{F}}+r^{-2}\text{\textgreek{D}}_{g_{\mathbb{S}^{d-1}}+h_{\mathbb{S}^{d-1}}}\text{\textgreek{F}}-\\
 & -\frac{(d-1)(d-3)}{4}r^{-2}\cdot\text{\textgreek{F}}+Err(\text{\textgreek{F}})\Big\}\, dudvd\text{\textgreek{sv}}.
\end{split}
\label{eq:BeforeIntegrationByPartsMorawetz}
\end{equation}
 Using the expression 
\[
\text{\textgreek{D}}_{g_{\mathbb{S}^{d-1}}+h_{\mathbb{S}^{d-1}}}=\text{\textgreek{D}}_{g_{\mathbb{S}^{d-1}}}+O(r^{-1})\partial_{\text{\textgreek{sv}}}\partial_{\text{\textgreek{sv}}}+O(r^{-1})\partial_{\text{\textgreek{sv}}},
\]
we obtain after integrating by parts in $\partial_{u}$, $\partial_{v}$
and $\partial_{\text{\textgreek{sv}}}$ and absorbing the error terms
in the $Err$ summand (and recalling that $\partial_{v}r=-\partial_{u}r=1$):
\begin{equation}
\begin{split}-\int_{\mathcal{R}(\text{\textgreek{t}}_{1},\text{\textgreek{t}}_{2})\cap\{t\le T^{*}\}} & \text{\textgreek{q}}f(r)\cdot\big(\partial_{v}-\partial_{u}\big)\text{\textgreek{F}}\cdot\text{\textgreek{W}}\square\text{\textgreek{f}}\, dudvd\text{\textgreek{sv}}=\\
= & \int_{\mathcal{R}(\text{\textgreek{t}}_{1},\text{\textgreek{t}}_{2})\cap\{t\le T^{*}\}}\text{\textgreek{q}}\cdot\frac{1}{2}\Big\{\big(f^{\prime}+O(r^{-1-a})f\big)|\partial_{v}\text{\textgreek{F}}|^{2}+\big(f^{\prime}+O(r^{-1-a})f\big)|\partial_{u}\text{\textgreek{F}}|^{2}+\\
 & \hphantom{\int_{\mathcal{R}(\text{\textgreek{t}}_{1},\text{\textgreek{t}}_{2})\cap\{t\le T^{*}\}}}+2\big(2r^{-1}(1+O(r^{-1}))f-(1+O(r^{-a}))f^{\prime}\big)\Big(|r^{-1}\partial_{\text{\textgreek{sv}}}\text{\textgreek{F}}|^{2}+\frac{(d-1)(d-3)}{4}|r^{-1}\text{\textgreek{F}}|^{2}\Big)\Big\}\, dudvd\text{\textgreek{sv}}+\\
 & +\int_{\mathcal{R}(\text{\textgreek{t}}_{1},\text{\textgreek{t}}_{2})\cap\{t\le T^{*}\}}\text{\textgreek{q}}f\cdot\big(\partial_{v}-\partial_{u}\big)\text{\textgreek{F}}\cdot\big(Err(\text{\textgreek{F}})+O(r^{-3})\partial_{\text{\textgreek{sv}}}\partial_{\text{\textgreek{sv}}}\text{\textgreek{F}}\big)\, dudvd\text{\textgreek{sv}}+\\
 & +\int_{\mathcal{R}(\text{\textgreek{t}}_{1},\text{\textgreek{t}}_{2})\cap\{t\le T^{*}\}}O(|\partial\text{\textgreek{q}}|)\cdot\big(|\partial\text{\textgreek{F}}|^{2}+r^{-2}|\text{\textgreek{F}}|^{2}\big)\, dudvd\text{\textgreek{sv}}+\\
 & +\sum_{i=1}^{2}\int_{\mathcal{S}_{\text{\textgreek{t}}_{1}}\cap\{t\le T^{*}\}}\text{\textgreek{q}}\cdot O(1)\big(|\partial_{v}\text{\textgreek{F}}|^{2}+|r^{-1}\partial_{\text{\textgreek{sv}}}\text{\textgreek{F}}|^{2}+r^{-1-\text{\textgreek{h}}'}|\partial_{u}\text{\textgreek{F}}|^{2}+r^{-2}|\text{\textgreek{F}}|^{2}\big)\, dvd\text{\textgreek{sv}}+\\
 & +\int_{\mathcal{R}(\text{\textgreek{t}}_{1},\text{\textgreek{t}}_{2})\cap\{t=T^{*}\}\cap\{r\ge R\}}O(1)\big(|\partial\text{\textgreek{F}}|^{2}+r^{-2}|\text{\textgreek{F}}|^{2}\big)\, dvd\text{\textgreek{sv}}.
\end{split}
\label{eq:BeforeIntegrationByPartsMorawetz-1}
\end{equation}

Since $\text{\textgreek{F}}=\text{\textgreek{W}}\text{\textgreek{f}}$,
we calculate by applying the product rule and expanding the square:
\begin{equation}
\begin{split}\int_{\mathcal{R}(\text{\textgreek{t}}_{1},\text{\textgreek{t}}_{2})\cap\{t\le T^{*}\}}\text{\textgreek{q}}\cdot\frac{1}{2}\big(f^{\prime}+O(r^{-1-a})f\big)|\partial_{v}\text{\textgreek{F}}|^{2}\, dudv & d\text{\textgreek{sv}}=\\
=\int_{\mathcal{R}(\text{\textgreek{t}}_{1},\text{\textgreek{t}}_{2})\cap\{t\le T^{*}\}}\text{\textgreek{q}}\cdot\frac{1}{2}\big(f^{\prime}+O(r^{-1-a})f\big)\big( & (\partial_{v}\text{\textgreek{f}})^{2}+r^{-2}(\frac{(d-1)^{2}}{4}+O(r^{-1}))\text{\textgreek{f}}^{2}+\\
 & +r^{-1}(\frac{d-1}{2}+O(r^{-1}))\partial_{v}(\text{\textgreek{f}}^{2})\big)\,\text{\textgreek{W}}^{2}dudvd\text{\textgreek{sv}}.
\end{split}
\end{equation}
 By integrating by parts in the $\partial_{v}(\text{\textgreek{f}}^{2})$
term, since for all $x>0$: 
\begin{equation}
-\frac{d}{dx}\big(x^{d-2}\frac{df}{dx}\big)+\frac{d-1}{2}x^{d-3}\frac{df}{dx}\ge c_{\text{\textgreek{h}}}x^{d-4+\text{\textgreek{h}}}(1+x^{\text{\textgreek{h}}})^{-2}-C_{\text{\textgreek{h}}}(d-3)x^{d-4}f,
\end{equation}
we obtain for any $\text{\textgreek{d}}_{0}\le1$ due to the form
of $f$ if $R$ is large enough in terms of $\text{\textgreek{h}}$:
\begin{align}
\int_{\mathcal{R}(\text{\textgreek{t}}_{1},\text{\textgreek{t}}_{2})\cap\{t\le T^{*}\}}\text{\textgreek{q}}\cdot\frac{1}{2}\big(f^{\prime}+O(r^{-1-a})f\big)|\partial_{v}\text{\textgreek{F}}|^{2}\, dudvd\text{\textgreek{sv}} & \ge\text{\textgreek{d}}_{0}\Big(c_{\text{\textgreek{h}}}\int_{\mathcal{R}(\text{\textgreek{t}}_{1},\text{\textgreek{t}}_{2})\cap\{t\le T^{*}\}}\text{\textgreek{q}}\cdot r^{-1-\text{\textgreek{h}}}\big((\partial_{v}\text{\textgreek{f}})^{2}+r^{-2}\text{\textgreek{f}}^{2}\big)\,\text{\textgreek{W}}^{2}dudvd\text{\textgreek{sv}}-\label{eq:DisguisedHardy}\\
 & -C_{\text{\textgreek{h}}}(d-3)\int_{\mathcal{R}(\text{\textgreek{t}}_{1},\text{\textgreek{t}}_{2})\cap\{t\le T^{*}\}}\text{\textgreek{q}}\cdot r^{-3}|\text{\textgreek{F}}|^{2}\, dudvd\text{\textgreek{sv}}-\nonumber \\
 & -C\int_{\mathcal{R}(\text{\textgreek{t}}_{1},\text{\textgreek{t}}_{2})\cap\{t\le T^{*}\}}|\partial\text{\textgreek{q}}|\cdot r^{-2}|\text{\textgreek{F}}|^{2}\, dudvd\text{\textgreek{sv}}-\nonumber \\
 & -C\sum_{i=1}^{2}\int_{\mathcal{S}_{\text{\textgreek{t}}_{i}}\cap\{t\le T^{*}\}}\text{\textgreek{q}}\cdot r^{-2}|\text{\textgreek{F}}|^{2}\, dvd\text{\textgreek{sv}}-C\int_{\mathcal{R}(\text{\textgreek{t}}_{1},\text{\textgreek{t}}_{2})\cap\{t=T^{*}\}}\text{\textgreek{q}}\cdot r^{-2}|\text{\textgreek{F}}|^{2}\, dvd\text{\textgreek{sv}}\Big).\nonumber 
\end{align}
 Notice that the $(d-3)\int\text{\textgreek{q}}\cdot r^{-3}|\text{\textgreek{F}}|^{2}$
error term in the right hand side (\ref{eq:DisguisedHardy}) can be
controlled by the corresponding term in the right hand side of (\ref{eq:BeforeIntegrationByPartsMorawetz-1}),
provided that $\text{\textgreek{d}}_{0}$ is small enough in terms
of $\text{\textgreek{h}}$.

Using, now the expression (\ref{eq:ErrorTermsMorawetz}) for $Err(\text{\textgreek{F}})$,
we can readily bound after integrating by parts in the highest order
terms (and in the $\text{\textgreek{F}}\partial_{u}\text{\textgreek{F}}=\frac{1}{2}\partial_{u}(\text{\textgreek{F}}^{2})$
term) and using a Cauchy--Schwarz inequality: 
\begin{equation}
\begin{split}\int_{\mathcal{R}(\text{\textgreek{t}}_{1},\text{\textgreek{t}}_{2})\cap\{t\le T^{*}\}}\text{\textgreek{q}}f(r)\cdot\big(\partial_{v}- & \partial_{u}\big)\text{\textgreek{F}}\cdot\big(Err(\text{\textgreek{F}})+O(r^{-3})\partial_{\text{\textgreek{sv}}}^{2}\text{\textgreek{F}}\big)\, dudvd\text{\textgreek{sv}}\le\\
\le & C_{\text{\textgreek{h}}}\int_{\mathcal{R}(\text{\textgreek{t}}_{1},\text{\textgreek{t}}_{2})\cap\{t\le T^{*}\}}\text{\textgreek{q}}\Big\{ O(r^{-1-a})|\partial_{u}\text{\textgreek{F}}|^{2}+O(r^{-1})\Big(|\partial_{v}\text{\textgreek{F}}|^{2}+|r^{-1}\partial_{\text{\textgreek{sv}}}\text{\textgreek{F}}|^{2}+r^{-2}|\text{\textgreek{F}}|^{2}\Big)\Big\}\, dudvd\text{\textgreek{sv}}+\\
 & +C_{\text{\textgreek{h}}}\int_{\mathcal{R}(\text{\textgreek{t}}_{1},\text{\textgreek{t}}_{2})\cap\{t\le T^{*}\}}|\partial\text{\textgreek{q}}|\big(|\partial\text{\textgreek{F}}|^{2}+r^{-2}|\text{\textgreek{F}}|^{2}\big)\, dudvd\text{\textgreek{sv}}+\\
 & +C_{\text{\textgreek{h}}}\sum_{i=1}^{2}\int_{\mathcal{S}_{\text{\textgreek{t}}_{1}}\cap\{t\le T^{*}\}}\text{\textgreek{q}}\cdot O(1)\big(|\partial_{v}\text{\textgreek{F}}|^{2}+|r^{-1}\partial_{\text{\textgreek{sv}}}\text{\textgreek{F}}|^{2}+r^{-1-\text{\textgreek{h}}'}|\partial_{u}\text{\textgreek{F}}|^{2}+r^{-2}|\text{\textgreek{F}}|^{2}\big)\, dvd\text{\textgreek{sv}}+\\
 & +C_{\text{\textgreek{h}}}\int_{\mathcal{R}(\text{\textgreek{t}}_{1},\text{\textgreek{t}}_{2})\cap\{t=T^{*}\}}\text{\textgreek{q}}\cdot\big(|\partial\text{\textgreek{F}}|^{2}+r^{-2}|\text{\textgreek{F}}|^{2}\big)\, dvd\text{\textgreek{sv}}.
\end{split}
\label{eq:BoundsErrorMorawetz}
\end{equation}

Therefore, from (\ref{eq:BeforeIntegrationByPartsMorawetz-1}), (\ref{eq:DisguisedHardy})
and (\ref{eq:BoundsErrorMorawetz}) for $\text{\textgreek{d}}_{0}$
small enough in terms of $\text{\textgreek{h}}$ we obtain the desired
result (\ref{eq:MorawetzGeneralCaseImproved}) if $R$ is large enough
in terms of $\text{\textgreek{h}}$.
\end{proof}
We can also establish the following variant of Lemma \ref{lem:MorawetzDrLemmaHyperboloids}
in the region bounded by a pair of $\{t=const\}$ hypersurfaces: 
\begin{lem}
\label{lem:MorawetzDrLemmaGeneralCase}For any given $0<\text{\textgreek{h}}<a$,
there exists an $R=R(\text{\textgreek{h}})>0$ and constants $C(\text{\textgreek{h}}),c(\text{\textgreek{h}})>0$
such that for any smooth function $\text{\textgreek{f}}:\mathcal{N}_{af}\rightarrow\mathbb{C}$
with compact support in space, any $t_{1}\le t_{2}$ and any smooth
cut-off function $\text{\textgreek{q}}:\mathcal{N}_{af}\rightarrow[0,1]$
supported in $\{r\ge R\}$ we can bound: 
\begin{equation}
\begin{split}\int_{\{t_{1}\le t\le t_{2}\}}\text{\textgreek{q}}\cdot\Big(r^{-1-\text{\textgreek{h}}}\Big(|\partial_{u}\text{\textgreek{f}}|^{2}+ & |\partial_{v}\text{\textgreek{f}}|^{2}\Big)+r^{-1}|r^{-1}\partial_{\text{\textgreek{sv}}}\text{\textgreek{f}}|^{2}+r^{-3-\text{\textgreek{h}}}|\text{\textgreek{f}}|^{2}\Big)\le\\
\le & C_{\text{\textgreek{h}}}\int_{\{t_{1}\le t\le t_{2}\}}|\partial\text{\textgreek{q}}|\cdot\big(|\partial\text{\textgreek{f}}|^{2}+r^{-2}|\text{\textgreek{f}}|^{2}\big)+C_{\text{\textgreek{h}}}\sum_{i=1}^{2}\int_{\{t=t_{i}\}\cap\{r\ge R\}}J_{\text{\textgreek{m}}}^{T}(\text{\textgreek{f}})n^{\text{\textgreek{m}}}+\\
 & +\int_{\{t_{1}\le t\le t_{2}\}}\text{\textgreek{q}}\cdot Re\big\{\big(O_{\text{\textgreek{h}}}(1)(\partial_{v}-\partial_{u})\bar{\text{\textgreek{f}}}+O_{\text{\textgreek{h}}}(r^{-1})\bar{\text{\textgreek{f}}}\big)\cdot\square_{g}\text{\textgreek{f}}\big\}+\\
 & +C_{\text{\textgreek{h}}}\int_{\{t_{1}\le t\le t_{2}\}}\text{\textgreek{q}}\cdot r^{-1}\Big(|\partial_{v}\text{\textgreek{f}}|^{2}+|r^{-1}\partial_{\text{\textgreek{sv}}}\text{\textgreek{f}}|^{2}+r^{-2}|\text{\textgreek{f}}|^{2}\Big).
\end{split}
\label{eq:MorawetzGeneralCaseRadiative}
\end{equation}
\end{lem}
\begin{rem*}
In case the radiative components of the metric satisfy the bounds
$\partial_{u}M\le0$ and $|\partial_{u}h_{as}|+|r\partial_{u}h_{\mathbb{S}^{d-1}}|\ll-(\partial_{u}M)+O(r^{-a})$
(which includes the non-radiative case $\partial_{u}M=0$, $h_{as}=0$
and $h_{\mathbb{S}^{d-1}}=O(r^{-1-a})$), the last term of the right
hand side of (\ref{eq:MorawetzGeneralCaseRadiative}) can be omitted.

The proof of this lemma is identical to the one for Lemma \ref{lem:MorawetzDrLemmaHyperboloids}
(the only difference being the domain of $\mathcal{N}_{af}$ over
which integrations by parts take place, and an application of a Hardy
type inequality for the boundary terms at $t=t_{1},t_{2}$). Hence,
the proof will be omitted.
\end{rem*}

\subsection{\label{sub:ImprovedMorawetz}An improved $\partial_{r}$- Morawetz
type estimate}

By a more careful choice of the function $f$ used in the proof of
Lemma \ref{lem:MorawetzDrLemmaHyperboloids}, we can obtain improved
control of the spacetime integral of $|\nabla\text{\textgreek{f}}|^{2},|\text{\textgreek{f}}|^{2}$
over any given compact subset of $\mathcal{R}(\text{\textgreek{t}}_{1},\text{\textgreek{t}}_{2})\cap\{r\ge R\}$,
at the expense of having to introduce a larger constant in the dependence
on the initial energy of $\text{\textgreek{f}}$, but without such
a loss in the $r\sim R$ boundary terms. A related microlocal construction
in the case of the subextremal Kerr family can be seen in \cite{DafRodSchlap}. 
\begin{lem}
\label{lem:ImprovedMorawetzGeneralCaseHyperboloids}For any given
$0<\text{\textgreek{h}}<a$, and any $R>0$ sufficiently large in
terms of $\text{\textgreek{h}}$, any $R_{c}\ge R$, any $\text{\textgreek{t}}_{1}\le\text{\textgreek{t}}_{2}$,
any $T^{*}>0$, any function $\text{\textgreek{q}}:\mathcal{N}_{af}\rightarrow[0,1]$
supported in $\{r\ge R\}$ with $\partial\text{\textgreek{q}}$ supported
in $\{R\le r\le R_{c}\}$ and any smooth function $\text{\textgreek{f}}:\mathcal{N}_{af}\rightarrow\mathbb{C}$
we can bound: 
\begin{equation}
\begin{split}\int_{\mathcal{R}(\text{\textgreek{t}}_{1},\text{\textgreek{t}}_{2})\cap\{t\le T^{*}\}\cap\{r\le R_{c}\}}\text{\textgreek{q}}\cdot\Big(\big(|\partial_{u} & \text{\textgreek{f}}|^{2}+|\partial_{v}\text{\textgreek{f}}|^{2}\big)+|r^{-1}\partial_{\text{\textgreek{sv}}}\text{\textgreek{f}}|^{2}+r^{-2}|\text{\textgreek{f}}|^{2}\Big)+\\
+R_{c}\int_{\mathcal{R}(\text{\textgreek{t}}_{1},\text{\textgreek{t}}_{2})\cap\{t\le T^{*}\}\cap\{r\ge R_{c}\}} & \text{\textgreek{q}}\cdot\Big(r^{-1-\text{\textgreek{h}}}\big(|\partial_{u}\text{\textgreek{f}}|^{2}+|\partial_{v}\text{\textgreek{f}}|^{2}\big)+r^{-1}|r^{-1}\partial_{\text{\textgreek{sv}}}\text{\textgreek{f}}|^{2}+r^{-3-\text{\textgreek{h}}}|\text{\textgreek{f}}|^{2}\Big)\le\\
\le & C(\text{\textgreek{h}})\cdot\int_{\mathcal{R}(\text{\textgreek{t}}_{1},\text{\textgreek{t}}_{2})\cap\{t\le T^{*}\}}|\partial\text{\textgreek{q}}|\cdot r\cdot\big(|\partial\text{\textgreek{f}}|^{2}+r^{-2}|\text{\textgreek{f}}|^{2}\big)+\\
 & +C(\text{\textgreek{h}},R_{c})\cdot\sum_{i=1}^{2}\int_{\mathcal{S}_{\text{\textgreek{t}}_{i}}\cap\{t\le T^{*}\}}\text{\textgreek{q}}\Big(|\partial_{v}\text{\textgreek{f}}|^{2}+|r^{-1}\partial_{\text{\textgreek{sv}}}\text{\textgreek{f}}|^{2}+r^{-1-\text{\textgreek{h}}'}|\partial_{u}\text{\textgreek{f}}|^{2}+r^{-2}|\text{\textgreek{f}}|^{2}\Big)\, r^{d-1}dvd\text{\textgreek{sv}}+\\
 & +C(\text{\textgreek{h}},R_{c})\cdot\int_{\{t=T^{*}\}\cap\mathcal{R}(\text{\textgreek{t}}_{1},\text{\textgreek{t}}_{2})}\text{\textgreek{q}}\big(|\partial\text{\textgreek{f}}|^{2}+r^{-2}|\text{\textgreek{f}}|^{2}\big)\, r^{d-1}dvd\text{\textgreek{sv}}+\\
 & +\int_{\mathcal{R}(\text{\textgreek{t}}_{1},\text{\textgreek{t}}_{2})\cap\{t\le T^{*}\}}\text{\textgreek{q}}\cdot Re\big\{\big(O_{R_{c},\text{\textgreek{h}}}(1)(\partial_{v}-\partial_{u})\bar{\text{\textgreek{f}}}+O_{R_{c},\text{\textgreek{h}}}(r^{-1})\bar{\text{\textgreek{f}}}\big)\cdot\square_{g}\text{\textgreek{f}}\big\}+\\
 & +C(\text{\textgreek{h}})R_{c}\cdot\int_{\mathcal{R}(\text{\textgreek{t}}_{1},\text{\textgreek{t}}_{2})\cap\{t\le T^{*}\}}\text{\textgreek{q}}\cdot r^{-1}\Big(|\partial_{v}\text{\textgreek{f}}|^{2}+|r^{-1}\partial_{\text{\textgreek{sv}}}\text{\textgreek{f}}|^{2}+r^{-2}|\text{\textgreek{f}}|^{2}\Big).
\end{split}
\label{eq:MorawetzGeneralCaseImprovedHyperboloids}
\end{equation}
\end{lem}
\begin{rem*}
In case the radiative components of the metric satisfy the bounds
$\partial_{u}M\le0$ and $|\partial_{u}h_{as}|+|r\partial_{u}h_{\mathbb{S}^{d-1}}|\ll-(\partial_{u}M)+O(r^{-a})$
(which includes the non-radiative case $\partial_{u}M=0$, $h_{as}=0$
and $h_{\mathbb{S}^{d-1}}=O(r^{-1-a})$), the last term of the right
hand side of (\ref{lem:ImprovedMorawetzGeneralCaseHyperboloids})
can be omitted. Furthermore, in case the $T$ vector field  satisfies
(\ref{eq:DeformationTensorTAway}) for $m=1$, then the last term
of the right hand side of (\ref{eq:MorawetzGeneralCaseImprovedHyperboloids})
is replaced by 
\begin{equation}
\int_{\mathcal{R}(\text{\textgreek{t}}_{1},\text{\textgreek{t}}_{2})\cap\{t\le T^{*}\}}\text{\textgreek{q}}\cdot\bar{t}^{-\text{\textgreek{d}}_{0}}r^{-1}\Big(|\partial_{v}\text{\textgreek{f}}|^{2}+|r^{-1}\partial_{\text{\textgreek{sv}}}\text{\textgreek{f}}|^{2}+r^{-2}|\text{\textgreek{f}}|^{2}\Big).
\end{equation}

Note that the constant in front of the boundary term in the $r\sim R$
region does not depend on $R_{c}$. This is where the importance of
this lemma lies, and this is where Lemma (\ref{lem:MorawetzDrLemmaHyperboloids})
would fail to give a similar statement.\end{rem*}
\begin{proof}
Without loss of generality, we will assume that $\text{\textgreek{f}}$
is real valued. 

Let us consider the smooth function $f_{R_{c}}:\mathcal{N}_{af}\rightarrow(0,+\infty)$
defined as $f_{R_{c}}=R_{c}g\circ(\frac{r}{R_{c}})$, where $g:(0,+\infty)\rightarrow(0,+\infty)$
is a smooth, increasing and concave function satisfying: 
\begin{equation}
g(x)=\begin{cases}
x, & x\le1,\\
2, & x\ge2.
\end{cases}\label{eq:PropertiesImprovedcurrent}
\end{equation}
 We will follow the proof of Lemma \ref{lem:MorawetzDrLemmaHyperboloids},
but we will use as a multiplier for equation (\ref{eq:ConformalWaveOperator-2})
the function 
\begin{equation}
f_{imp}=f\cdot f_{R_{c}},
\end{equation}
 where 
\begin{equation}
f=\frac{r^{\text{\textgreek{h}}}}{1+r^{\text{\textgreek{h}}}}
\end{equation}
 is the seed function used in the proof of Lemma \ref{lem:MorawetzDrLemmaHyperboloids}.

Notice that for $f_{imp}$ we can calculate (since $\text{\textgreek{W}}^{2}=r^{d-1}(1+O(r^{-1}))$):
\begin{equation}
\partial_{v}f_{imp}=-\partial_{u}f_{imp}=\big(1+O(r^{-\text{\textgreek{h}}})\big)\cdot g^{\prime}(\frac{r}{R_{c}})+\text{\textgreek{h}}r^{-1-\text{\textgreek{h}}}R_{c}\big(1+O(r^{-1})\big)\cdot g(\frac{r}{R_{c}})
\end{equation}
 and:
\begin{align*}
\partial_{v}\big(r^{-1}\text{\textgreek{W}}^{2}\partial_{u}f_{imp}\big)= & -r^{d-2}R_{c}^{-1}\big(1+O(r^{-1})\big)\cdot g^{\prime\prime}(\frac{r}{R_{c}})-\\
 & -\big((d-2)r^{d-3}+2\text{\textgreek{h}}r^{d-3-\text{\textgreek{h}}}\big)\cdot\big(1+O(r^{-1})\big)\cdot g^{\prime}(\frac{r}{R_{c}})+\\
 & +\text{\textgreek{h}}\big(\text{\textgreek{h}}-(d-3)\big)r^{d-4-\text{\textgreek{h}}}R_{c}\big(1+O(r^{-1})\big)\cdot g(\frac{r}{R_{c}}).
\end{align*}
Therefore, since $g$ was assumed to be increasing and concave and
satisfies (\ref{eq:PropertiesImprovedcurrent}), we can bound: 
\begin{equation}
\partial_{v}\big(r^{-1}\text{\textgreek{W}}^{2}\partial_{u}f_{imp}\big)-r^{-2}\text{\textgreek{W}}^{2}\partial_{u}f_{imp}\ge\begin{cases}
\big(-(d-3)r^{-2}+\text{\textgreek{h}}r^{-2-\text{\textgreek{h}}}+O(r^{-3})\big)\text{\textgreek{W}}^{2}, & r\le R_{c}\\
R_{c}\big(-(d-3)r^{-3}+\text{\textgreek{h}}^{2}r^{-3-\text{\textgreek{h}}}+O(r^{-4})\big)\text{\textgreek{W}}^{2}, & r\ge R_{c}.
\end{cases}\label{eq:BoundZerothOrderCoefficient}
\end{equation}

Repeating the proof of Lemma \ref{lem:MorawetzDrLemmaHyperboloids},
setting $\text{\textgreek{F}}=\text{\textgreek{W}}\text{\textgreek{f}}$
and multiplying equation (\ref{eq:ConformalWaveOperator-2}) with
$\text{\textgreek{q}}\cdot f_{imp}\cdot\big(\partial_{v}-\partial_{u}\big)\text{\textgreek{F}}$,
we obtain after integrating over $\mathcal{R}(\text{\textgreek{t}}_{1},\text{\textgreek{t}}_{2})\cap\{t\le T^{*}\}$
(with $dudvd\text{\textgreek{sv}}$ used as a volume form) and integrating
by parts in $\partial_{u}$, $\partial_{v}$ and $\partial_{\text{\textgreek{sv}}}$:
\begin{equation}
\begin{split}-\int_{\mathcal{R}(\text{\textgreek{t}}_{1},\text{\textgreek{t}}_{2})\cap\{t\le T^{*}\}}\text{\textgreek{q}}f_{imp} & \cdot\big(\partial_{v}-\partial_{u}\big)\text{\textgreek{F}}\cdot\text{\textgreek{W}}\square\text{\textgreek{f}}\, dudvd\text{\textgreek{sv}}=\\
= & \int_{\mathcal{R}(\text{\textgreek{t}}_{1},\text{\textgreek{t}}_{2})\cap\{t\le T^{*}\}}\text{\textgreek{q}}\cdot\frac{1}{2}\Big\{\big(-\partial_{u}f_{imp}+O(r^{-1-a})f_{imp}\big)|\partial_{v}\text{\textgreek{F}}|^{2}+\big(\partial_{v}f_{imp}+O(r^{-1-a})f_{imp}\big)|\partial_{u}\text{\textgreek{F}}|^{2}+\\
 & \hphantom{\int_{\mathcal{R}(\text{\textgreek{t}}_{1},\text{\textgreek{t}}_{2})\cap\{t\le T^{*}\}}\text{\textgreek{q}}\cdot\frac{1}{2}\Big\{}+\mathcal{A}_{imp}\Big(|r^{-1}\partial_{\text{\textgreek{sv}}}\text{\textgreek{F}}|^{2}+\frac{(d-1)(d-3)}{4}|r^{-1}\text{\textgreek{F}}|^{2}\Big)\Big\}\, dudvd\text{\textgreek{sv}}+\\
 & +\int_{\mathcal{R}(\text{\textgreek{t}}_{1},\text{\textgreek{t}}_{2})\cap\{t\le T^{*}\}}\text{\textgreek{q}}f_{imp}\cdot\big(\partial_{v}-\partial_{u}\big)\text{\textgreek{F}}\cdot\big(Err(\text{\textgreek{F}})+O(r^{-3})\partial_{\text{\textgreek{sv}}}^{2}\text{\textgreek{F}}\big)\, dudvd\text{\textgreek{sv}}+\\
 & +\int_{\mathcal{R}(\text{\textgreek{t}}_{1},\text{\textgreek{t}}_{2})\cap\{t\le T^{*}\}}O(f_{imp}|\partial\text{\textgreek{q}}|)\cdot\big(|\partial\text{\textgreek{F}}|^{2}+r^{-2}|\text{\textgreek{F}}|^{2}\big)\, dudvd\text{\textgreek{sv}}+\\
 & +\sum_{i=1}^{2}\int_{\mathcal{S}_{\text{\textgreek{t}}_{i}}\cap\{t\le T^{*}\}}\text{\textgreek{q}}O(f_{imp})\big(|\partial_{v}\text{\textgreek{F}}|^{2}+|r^{-1}\partial_{\text{\textgreek{sv}}}\text{\textgreek{F}}|^{2}+r^{-1-\text{\textgreek{h}}'}|\partial_{u}\text{\textgreek{F}}|^{2}+r^{-2}|\text{\textgreek{F}}|^{2}\big)\, dvd\text{\textgreek{sv}}+\\
 & +\int_{\mathcal{R}(\text{\textgreek{t}}_{1},\text{\textgreek{t}}_{2})\cap\{t=T^{*}\}}O(f_{imp})\big(|\partial\text{\textgreek{F}}|^{2}+r^{-2}|\text{\textgreek{F}}|^{2}\big)\, dvd\text{\textgreek{sv}},
\end{split}
\label{eq:BeforeIntegrationByPartsMorawetz-1-1}
\end{equation}
where 
\begin{equation}
\mathcal{A}_{imp}\doteq2\big(2r^{-1}(1+O(r^{-1}))f_{imp}-\frac{1}{2}(1+O(r^{-a}))\big(\partial_{v}-\partial_{u}\big)f_{imp}\big).
\end{equation}

Since $\text{\textgreek{F}}=\text{\textgreek{W}}\text{\textgreek{f}}$,
we calculate by applying the product rule and expanding the square:
\begin{equation}
\begin{split}\int_{\mathcal{R}(\text{\textgreek{t}}_{1},\text{\textgreek{t}}_{2})\cap\{t\le T^{*}\}} & \text{\textgreek{q}}\cdot\frac{1}{2}\big(-\partial_{u}f_{imp}+O(r^{-1-a})f_{imp}\big)|\partial_{v}\text{\textgreek{F}}|^{2}\, dudvd\text{\textgreek{sv}}=\\
= & \int_{\mathcal{R}(\text{\textgreek{t}}_{1},\text{\textgreek{t}}_{2})\cap\{t\le T^{*}\}}\text{\textgreek{q}}\cdot\big(-\partial_{u}f_{imp}+O(r^{-1-a})f_{imp}\big)\Big\{(\partial_{v}\text{\textgreek{f}})^{2}+r^{-2}(1+O(r^{-1}))\text{\textgreek{f}}^{2}+\\
 & \hphantom{\int_{\mathcal{R}(\text{\textgreek{t}}_{1},\text{\textgreek{t}}_{2})\cap\{t\le T^{*}\}}\text{\textgreek{q}}\cdot\big(-\partial_{u}f_{imp}+O(r^{-1-a})f_{imp}\big)\Big\{}+r^{-1}(1+O(r^{-1}))\partial_{v}(\text{\textgreek{f}}^{2})\Big\}\,\text{\textgreek{W}}^{2}dudvd\text{\textgreek{sv}}.
\end{split}
\label{eq:ExpandingTheSquare}
\end{equation}
 By integrating by parts in the $\partial_{v}(\text{\textgreek{f}}^{2})$
term, we thus obtain (due to (\ref{eq:BoundZerothOrderCoefficient}))
for any $\text{\textgreek{d}}_{0}\le1$: 
\begin{equation}
\begin{split}\int_{\mathcal{R}(\text{\textgreek{t}}_{1},\text{\textgreek{t}}_{2})\cap\{t\le T^{*}\}}\text{\textgreek{q}}\cdot\frac{1}{2}\big(-\partial_{u}f_{imp}+O( & r^{-1-a})f_{imp}\big)|\partial_{v}\text{\textgreek{F}}|^{2}\, dudvd\text{\textgreek{sv}}\ge\\
\ge\text{\textgreek{d}}_{0}\Big\{ & c_{\text{\textgreek{h}}}\int_{\mathcal{R}(\text{\textgreek{t}}_{1},\text{\textgreek{t}}_{2})\cap\{t\le T^{*}\}\cap\{r\le R_{c}\}}\text{\textgreek{q}}\cdot\big((\partial_{v}\text{\textgreek{f}})^{2}+r^{-2-\text{\textgreek{h}}}\text{\textgreek{f}}^{2}\big)\,\text{\textgreek{W}}^{2}dudvd\text{\textgreek{sv}}+\\
 & +c_{\text{\textgreek{h}}}R_{c}\int_{\mathcal{R}(\text{\textgreek{t}}_{1},\text{\textgreek{t}}_{2})\cap\{t\le T^{*}\}\cap\{r\ge R_{c}\}}\text{\textgreek{q}}\cdot r^{-1-\text{\textgreek{h}}}\big((\partial_{v}\text{\textgreek{f}})^{2}+r^{-2}\text{\textgreek{f}}^{2}\big)\,\text{\textgreek{W}}^{2}dudvd\text{\textgreek{sv}}-\\
 & -C_{\text{\textgreek{h}}}(d-3)\int_{\mathcal{R}(\text{\textgreek{t}}_{1},\text{\textgreek{t}}_{2})\cap\{t\le T^{*}\}\cap\{r\le R_{c}\}}\text{\textgreek{q}}\cdot r^{-2}|\text{\textgreek{F}}|^{2}\, dudvd\text{\textgreek{sv}}-\\
 & -C_{\text{\textgreek{h}}}(d-3)R_{c}\int_{\mathcal{R}(\text{\textgreek{t}}_{1},\text{\textgreek{t}}_{2})\cap\{t\le T^{*}\}\cap\{r\ge R_{c}\}}\text{\textgreek{q}}\cdot r^{-2}|\text{\textgreek{F}}|^{2}\, dudvd\text{\textgreek{sv}}-\\
 & -\int_{\mathcal{R}(\text{\textgreek{t}}_{1},\text{\textgreek{t}}_{2})\cap\{t\le T^{*}\}}O(f_{imp}|\partial\text{\textgreek{q}}|)\cdot r^{-2}|\text{\textgreek{F}}|^{2}\, dudvd\text{\textgreek{sv}}-\\
 & -C(R_{c})\sum_{i=1}^{2}\int_{\mathcal{S}_{\text{\textgreek{t}}_{i}}\cap\{t\le T^{*}\}}\text{\textgreek{q}}r^{-2}|\text{\textgreek{F}}|^{2}\, dvd\text{\textgreek{sv}}-C(R_{c})\int_{\mathcal{R}(\text{\textgreek{t}}_{1},\text{\textgreek{t}}_{2})\cap\{t=T^{*}\}}\text{\textgreek{q}}r^{-2}|\text{\textgreek{F}}|^{2}\, dvd\text{\textgreek{sv}}\Big\}.
\end{split}
\label{eq:DisguisedHardy-1}
\end{equation}

Moreover, using a trivial Cauchy--Schwarz inequality for the $\text{\textgreek{f}}\partial_{v}\text{\textgreek{f}}$
term we can bound from (\ref{eq:ExpandingTheSquare}): 
\begin{multline}
\int_{\mathcal{R}(\text{\textgreek{t}}_{1},\text{\textgreek{t}}_{2})\cap\{t\le T^{*}\}\cap\{r\le R_{c}\}}\text{\textgreek{q}}\cdot\frac{1}{2}\big(-\partial_{u}f_{imp}+O(r^{-1-a})f_{imp}\big)|\partial_{v}\text{\textgreek{F}}|^{2}\, dudvd\text{\textgreek{sv}}\ge\\
\ge c_{\text{\textgreek{h}}}\int_{\mathcal{R}(\text{\textgreek{t}}_{1},\text{\textgreek{t}}_{2})\cap\{t\le T^{*}\}\cap\{r\le R_{c}\}}\text{\textgreek{q}}\cdot\big(c_{\text{\textgreek{h}}}r^{-2}\text{\textgreek{f}}^{2}-C_{\text{\textgreek{h}}}(\partial_{v})^{2}\big)\,\text{\textgreek{W}}^{2}dudvd\text{\textgreek{sv}},\label{eq:TrivialCauchySchwarz}
\end{multline}
 and thus by adding to (\ref{eq:DisguisedHardy-1}) a small multiple
of (\ref{eq:TrivialCauchySchwarz}) we obtain: 
\begin{equation}
\begin{split}\int_{\mathcal{R}(\text{\textgreek{t}}_{1},\text{\textgreek{t}}_{2})\cap\{t\le T^{*}\}}\text{\textgreek{q}}\cdot\frac{1}{2}\big(-\partial_{u}f_{imp}+O( & r^{-1-a})f_{imp}\big)|\partial_{v}\text{\textgreek{F}}|^{2}\, dudvd\text{\textgreek{sv}}\ge\\
\ge\text{\textgreek{d}}_{0}\Big\{ & c_{\text{\textgreek{h}}}\int_{\mathcal{R}(\text{\textgreek{t}}_{1},\text{\textgreek{t}}_{2})\cap\{t\le T^{*}\}\cap\{r\le R_{c}\}}\text{\textgreek{q}}\cdot\big((\partial_{v}\text{\textgreek{f}})^{2}+r^{-2}\text{\textgreek{f}}^{2}\big)\,\text{\textgreek{W}}^{2}dudvd\text{\textgreek{sv}}+\\
 & +c_{\text{\textgreek{h}}}R_{c}\int_{\mathcal{R}(\text{\textgreek{t}}_{1},\text{\textgreek{t}}_{2})\cap\{t\le T^{*}\}\cap\{r\ge R_{c}\}}\text{\textgreek{q}}\cdot r^{-1-\text{\textgreek{h}}}\big((\partial_{v}\text{\textgreek{f}})^{2}+r^{-2}\text{\textgreek{f}}^{2}\big)\,\text{\textgreek{W}}^{2}dudvd\text{\textgreek{sv}}-\\
 & -C_{\text{\textgreek{h}}}(d-3)\int_{\mathcal{R}(\text{\textgreek{t}}_{1},\text{\textgreek{t}}_{2})\cap\{t\le T^{*}\}\cap\{r\le R_{c}\}}\text{\textgreek{q}}\cdot r^{-2}|\text{\textgreek{F}}|^{2}\, dudvd\text{\textgreek{sv}}-\\
 & -C_{\text{\textgreek{h}}}(d-3)R_{c}\int_{\mathcal{R}(\text{\textgreek{t}}_{1},\text{\textgreek{t}}_{2})\cap\{t\le T^{*}\}\cap\{r\ge R_{c}\}}\text{\textgreek{q}}\cdot r^{-2}|\text{\textgreek{F}}|^{2}\, dudvd\text{\textgreek{sv}}-\\
 & -\int_{\mathcal{R}(\text{\textgreek{t}}_{1},\text{\textgreek{t}}_{2})\cap\{t\le T^{*}\}}O(f_{imp}|\partial\text{\textgreek{q}}|)\cdot r^{-2}|\text{\textgreek{F}}|^{2}\, dudvd\text{\textgreek{sv}}-\\
 & -C(R_{c})\sum_{i=1}^{2}\int_{\mathcal{S}_{\text{\textgreek{t}}_{i}}\cap\{t\le T^{*}\}}\text{\textgreek{q}}r^{-2}|\text{\textgreek{F}}|^{2}\, dvd\text{\textgreek{sv}}-C(R_{c})\int_{\mathcal{R}(\text{\textgreek{t}}_{1},\text{\textgreek{t}}_{2})\cap\{t=T^{*}\}}\text{\textgreek{q}}r^{-2}|\text{\textgreek{F}}|^{2}\, dvd\text{\textgreek{sv}}\Big\}.
\end{split}
\label{eq:DisguisedHardy-1-1}
\end{equation}

Using, now the expression (\ref{eq:ErrorTermsMorawetz}) for $Err(\text{\textgreek{F}})$,
we can readily bound after integrating by parts in the highest order
terms (and in the $\text{\textgreek{F}}\partial_{u}\text{\textgreek{F}}=\frac{1}{2}\partial_{u}(\text{\textgreek{F}}^{2})$
term) and using a Cauchy--Schwarz inequality: 
\begin{equation}
\begin{split}\int_{\mathcal{R}(\text{\textgreek{t}}_{1},\text{\textgreek{t}}_{2})\cap\{t\le T^{*}\}} & \text{\textgreek{q}}f_{imp}\cdot\big(\partial_{v}-\partial_{u}\big)\text{\textgreek{F}}\cdot\big(Err(\text{\textgreek{F}})+O(r^{-3})\partial_{\text{\textgreek{sv}}}^{2}\text{\textgreek{F}}\big)\, dudvd\text{\textgreek{sv}}\le\\
\le\, & C_{\text{\textgreek{h}}}\int_{\mathcal{R}(\text{\textgreek{t}}_{1},\text{\textgreek{t}}_{2})\cap\{t\le T^{*}\}\cap\{r\le R_{c}\}}\text{\textgreek{q}}\Big\{ O(r^{-a})|\partial_{u}\text{\textgreek{F}}|^{2}+\Big(|\partial_{v}\text{\textgreek{F}}|^{2}+|r^{-1}\partial_{\text{\textgreek{sv}}}\text{\textgreek{F}}|^{2}+r^{-2}|\text{\textgreek{F}}|^{2}\Big)\Big\}\, dudvd\text{\textgreek{sv}}+\\
 & +C_{\text{\textgreek{h}}}R_{c}\int_{\mathcal{R}(\text{\textgreek{t}}_{1},\text{\textgreek{t}}_{2})\cap\{t\le T^{*}\}\cap\{r\ge R_{c}\}}\text{\textgreek{q}}\cdot\Big\{ O(r^{-1-a})|\partial_{u}\text{\textgreek{F}}|^{2}+O(r^{-1})\Big(|\partial_{v}\text{\textgreek{F}}|^{2}+|r^{-1}\partial_{\text{\textgreek{sv}}}\text{\textgreek{F}}|^{2}+r^{-2}|\text{\textgreek{F}}|^{2}\Big)\Big\}\, dudvd\text{\textgreek{sv}}+\\
 & +C_{\text{\textgreek{h}}}\int_{\mathcal{R}(\text{\textgreek{t}}_{1},\text{\textgreek{t}}_{2})\cap\{t\le T^{*}\}}O(f_{imp}|\partial\text{\textgreek{q}}|)\big(|\partial\text{\textgreek{F}}|^{2}+r^{-2}|\text{\textgreek{F}}|^{2}\big)\, dudvd\text{\textgreek{sv}}+\\
 & +C_{\text{\textgreek{h}}}(R_{c})\sum_{i=1}^{2}\int_{\mathcal{S}_{\text{\textgreek{t}}_{i}}\cap\{t\le T^{*}\}}\text{\textgreek{q}}O(f_{imp})\big(|\partial_{v}\text{\textgreek{F}}|^{2}+|r^{-1}\partial_{\text{\textgreek{sv}}}\text{\textgreek{F}}|^{2}+r^{-1-\text{\textgreek{h}}'}|\partial_{u}\text{\textgreek{F}}|^{2}+r^{-2}|\text{\textgreek{F}}|^{2}\big)\, dvd\text{\textgreek{sv}}+\\
 & +C_{\text{\textgreek{h}}}(R_{c})\int_{\mathcal{R}(\text{\textgreek{t}}_{1},\text{\textgreek{t}}_{2})\cap\{t=T^{*}\}}O(f_{imp})\big(|\partial\text{\textgreek{F}}|^{2}+r^{-2}|\text{\textgreek{F}}|^{2}\big)\, dvd\text{\textgreek{sv}}.
\end{split}
\label{eq:BoundsErrorMorawetz-1}
\end{equation}

Therefore, from (\ref{eq:BeforeIntegrationByPartsMorawetz-1-1}),
(\ref{eq:DisguisedHardy-1}) and (\ref{eq:BoundsErrorMorawetz-1})
for $\text{\textgreek{d}}_{0}$ small enough in terms of $\text{\textgreek{h}}$,
we obtain the desired inequality if $R$ is large enough in terms
of $\text{\textgreek{h}}$: 
\begin{equation}
\begin{split}\int_{\mathcal{R}(\text{\textgreek{t}}_{1},\text{\textgreek{t}}_{2})\cap\{t\le T^{*}\}\cap\{r\le R_{c}\}}\text{\textgreek{q}}\cdot\Big(\big(|\partial_{u} & \text{\textgreek{f}}|^{2}+|\partial_{v}\text{\textgreek{f}}|^{2}\big)+|r^{-1}\partial_{\text{\textgreek{sv}}}\text{\textgreek{f}}|^{2}+r^{-2}|\text{\textgreek{f}}|^{2}\Big)+\\
+R_{c}\int_{\mathcal{R}(\text{\textgreek{t}}_{1},\text{\textgreek{t}}_{2})\cap\{t\le T^{*}\}\cap\{r\ge R_{c}\}} & \text{\textgreek{q}}\cdot\Big(r^{-1-\text{\textgreek{h}}}\big(|\partial_{u}\text{\textgreek{f}}|^{2}+|\partial_{v}\text{\textgreek{f}}|^{2}\big)+r^{-1}|r^{-1}\partial_{\text{\textgreek{sv}}}\text{\textgreek{f}}|^{2}+r^{-3-\text{\textgreek{h}}}|\text{\textgreek{f}}|^{2}\Big)\le\\
\le & C(\text{\textgreek{h}})\cdot\int_{\mathcal{R}(\text{\textgreek{t}}_{1},\text{\textgreek{t}}_{2})\cap\{t\le T^{*}\}}|\partial\text{\textgreek{q}}|\cdot r\cdot\big(|\partial\text{\textgreek{f}}|^{2}+r^{-2}|\text{\textgreek{f}}|^{2}\big)+\\
 & +C(\text{\textgreek{h}},R_{c})\cdot\sum_{i=1}^{2}\int_{\mathcal{S}_{\text{\textgreek{t}}_{i}}\cap\{t\le T^{*}\}}\text{\textgreek{q}}\Big(|\partial_{v}\text{\textgreek{f}}|^{2}+|r^{-1}\partial_{\text{\textgreek{sv}}}\text{\textgreek{f}}|^{2}+r^{-1-\text{\textgreek{h}}'}|\partial_{u}\text{\textgreek{f}}|^{2}+r^{-2}|\text{\textgreek{f}}|^{2}\Big)\, r^{d-1}dvd\text{\textgreek{sv}}+\\
 & +C(\text{\textgreek{h}},R_{c})\cdot\int_{\{t=T^{*}\}\cap\mathcal{R}(\text{\textgreek{t}}_{1},\text{\textgreek{t}}_{2})}\text{\textgreek{q}}\big(|\partial\text{\textgreek{f}}|^{2}+r^{-2}|\text{\textgreek{f}}|^{2}\big)\, r^{d-1}dvd\text{\textgreek{sv}}+\\
 & +\int_{\mathcal{R}(\text{\textgreek{t}}_{1},\text{\textgreek{t}}_{2})\cap\{t\le T^{*}\}}\text{\textgreek{q}}\cdot\big(O_{R_{c},\text{\textgreek{h}}}(1)(\partial_{v}-\partial_{u})\text{\textgreek{f}}+O_{R_{c},\text{\textgreek{h}}}(r^{-1})\text{\textgreek{f}}\big)\cdot\square_{g}\text{\textgreek{f}}+\\
 & +C(\text{\textgreek{h}})R_{c}\cdot\int_{\mathcal{R}(\text{\textgreek{t}}_{1},\text{\textgreek{t}}_{2})\cap\{t\le T^{*}\}}\text{\textgreek{q}}\cdot r^{-1}\Big(|\partial_{v}\text{\textgreek{f}}|^{2}+|r^{-1}\partial_{\text{\textgreek{sv}}}\text{\textgreek{f}}|^{2}+r^{-2}|\text{\textgreek{f}}|^{2}\Big).
\end{split}
\label{eq:MorawetzGeneralCaseImprovedHyperboloids-1}
\end{equation}

\end{proof}
In the same way, we can establish a similar lemma in the region spanned
by two hypersurfaces of the form $\{t=const\}$.
\begin{lem}
\label{lem:ImprovedMorawetzGeneralCase}For any given $0<\text{\textgreek{h}}<a$,
and any $R>0$ sufficiently large in terms of $\text{\textgreek{h}}$,
any $R_{c}\ge R$, any $t_{1}\le t_{2}$, any function $\text{\textgreek{q}}:\mathcal{N}_{af}\rightarrow[0,1]$
supported in $\{r\ge R\}$ with $\partial\text{\textgreek{q}}$ supported
in $\{R\le r\le R_{c}\}$ and any smooth function $\text{\textgreek{f}}:\mathcal{N}_{af}\rightarrow\mathbb{C}$
with compact support in space we can bound: 
\begin{equation}
\begin{split}\int_{\{t_{1}\le t\le t_{2}\}\cap\{r\le R_{c}\}}\text{\textgreek{q}}\cdot\Big(\big(|\partial_{u}\text{\textgreek{f}}|^{2}+|\partial_{v} & \text{\textgreek{f}}|^{2}\big)+|r^{-1}\partial_{\text{\textgreek{sv}}}\text{\textgreek{f}}|^{2}+r^{-2}|\text{\textgreek{f}}|^{2}\Big)+\\
+R_{c}\int_{\{t_{1}\le t\le t_{2}\}\cap\{r\ge R_{c}\}}\text{\textgreek{q}}\cdot\Big(r^{-1-\text{\textgreek{h}}}\big( & |\partial_{u}\text{\textgreek{f}}|^{2}+|\partial_{v}\text{\textgreek{f}}|^{2}\big)+r^{-1}|r^{-1}\partial_{\text{\textgreek{sv}}}\text{\textgreek{f}}|^{2}+r^{-3-\text{\textgreek{h}}}|\text{\textgreek{f}}|^{2}\Big)\le\\
\le & C(\text{\textgreek{h}})\cdot\int_{\{t_{1}\le t\le t_{2}\}}|\partial\text{\textgreek{q}}|\cdot r\cdot\big(|\partial\text{\textgreek{f}}|^{2}+r^{-2}|\text{\textgreek{f}}|^{2}\big)+\\
 & +C(\text{\textgreek{h}},R_{c})\cdot\int_{\{t=t_{1}\}\cap\{r\ge R\}}|\partial\text{\textgreek{f}}|^{2}+C(\text{\textgreek{h}},R_{c})\cdot\int_{\{t=t_{2}\}\cap\{r\ge R\}}|\partial\text{\textgreek{f}}|^{2}+\\
 & +\int_{\{t_{1}\le t\le t_{2}\}}\text{\textgreek{q}}\cdot Re\big\{\big(O_{R_{c},\text{\textgreek{h}}}(1)(\partial_{v}-\partial_{u})\bar{\text{\textgreek{f}}}+O_{R_{c},\text{\textgreek{h}}}(r^{-1})\bar{\text{\textgreek{f}}}\big)\cdot\square_{g}\text{\textgreek{f}}\big\}+\\
 & +C(\text{\textgreek{h}})R_{c}\cdot\int_{\{t_{1}\le t\le t_{2}\}}\text{\textgreek{q}}\cdot r^{-1}\Big(|\partial_{v}\text{\textgreek{f}}|^{2}+|r^{-1}\partial_{\text{\textgreek{sv}}}\text{\textgreek{f}}|^{2}+r^{-2}|\text{\textgreek{f}}|^{2}\Big).
\end{split}
\label{eq:MorawetzGeneralCaseImproved}
\end{equation}
\end{lem}
\begin{rem*}
In case the radiative components of the metric satisfy the bounds
$\partial_{u}M\le0$ and $|\partial_{u}h_{as}|+|r\partial_{u}h_{\mathbb{S}^{d-1}}|\ll-(\partial_{u}M)+O(r^{-a})$
(which includes the non-radiative case $\partial_{u}M=0$, $h_{as}=0$
and $h_{\mathbb{S}^{d-1}}=O_{1}(r^{-1-a})$), the last term of the
right hand side of (\ref{lem:ImprovedMorawetzGeneralCase}) can be
omitted.
\end{rem*}

\subsection{\label{sub:EnergyBoundedness}Estimates for the $J^{T}$-energy}

In this Section, we will establish some useful estimates for the energy
associated to the timelike vector field $T=\partial_{u}+\partial_{v}$
in the coordinate chart $(u,v,\text{\textgreek{sv}})$. Since we have
not assumed $(\mathcal{N}_{af},g)$ to be stationary, the $J^{T}$-energy
current will not be in general conserved. 
\begin{lem}
\label{lem:Boundedness Hyperboloids}For any smooth function $\text{\textgreek{f}}:\mathcal{N}_{af}\rightarrow\mathbb{C}$,
any $\text{\textgreek{q}}:\mathcal{N}_{af}\rightarrow[0,1]$ supported
on $\{r\ge R\}$ for some $R>0$ large in terms of the geometry of
$(\mathcal{N}_{af},g)$, any $\text{\textgreek{t}}_{1}\le\text{\textgreek{t}}_{2}$
and any $T^{*}>0$ the following estimate is true: 
\begin{equation}
\begin{split}\int_{\mathcal{S}_{\text{\textgreek{t}}_{2}}\cap\{t\le T^{*}\}} & \text{\textgreek{q}}\cdot J_{\text{\textgreek{m}}}^{T}(\text{\textgreek{f}})\bar{n}^{\text{\textgreek{m}}}+\int_{\{t=T^{*}\}\cap\mathcal{R}(\text{\textgreek{t}}_{1},\text{\textgreek{t}}_{2})}\Big(|\partial\text{\textgreek{f}}|^{2}+r^{-2}|\text{\textgreek{f}}|^{2}\Big)\le\\
\le & C\cdot\Bigg\{\int_{\mathcal{S}_{\text{\textgreek{t}}_{i}}\cap\{t\le T^{*}\}}\text{\textgreek{q}}\cdot J_{\text{\textgreek{m}}}^{T}(\text{\textgreek{f}})\bar{n}^{\text{\textgreek{m}}}+\int_{\mathcal{R}(\text{\textgreek{t}}_{1},\text{\textgreek{t}}_{2})\cap\{t\le T^{*}\}}\text{\textgreek{q}}\cdot r^{-1}\Big(|\partial_{v}\text{\textgreek{f}}|^{2}+|r^{-1}\partial_{\text{\textgreek{sv}}}\text{\textgreek{f}}|^{2}+|r^{-1}\text{\textgreek{f}}|^{2}\Big)+\\
 & +\int_{\mathcal{R}(\text{\textgreek{t}}_{1},\text{\textgreek{t}}_{2})\cap\{t\le T^{*}\}}O(|\partial\text{\textgreek{q}}|)\cdot\big(|\partial\text{\textgreek{f}}|^{2}+r^{-1}|\text{\textgreek{f}}|^{2}\big)+\int_{\mathcal{R}(\text{\textgreek{t}}_{1},\text{\textgreek{t}}_{2})\cap\{t\le T^{*}\}}\text{\textgreek{q}}O(1)\cdot Re\Big\{\text{\textgreek{W}}^{-1}\big(\partial_{v}+\partial_{u}\big)(\text{\textgreek{W}}\bar{\text{\textgreek{f}}})\cdot\square\text{\textgreek{f}}\Big\}\Bigg\}.
\end{split}
\label{eq:BoundednessGeneralRadiative}
\end{equation}
 for some constant $C>0$ depending only on the geometry of $(\mathcal{N}_{af},g)$.
In the above, $\bar{n}$ and $n$ denote the future directed unit
normals to the hypersurfaces $\mathcal{S}_{\text{\textgreek{t}}}$
and $\{t=const\}$ respectively.\end{lem}
\begin{rem*}
In case the radiative components of the metric satisfy the bounds
$\partial_{u}M\le0$ and $|\partial_{u}h_{as}|+|r\partial_{u}h_{\mathbb{S}^{d-1}}|\ll-(\partial_{u}M)+O(r^{-a})$
(which includes the non-radiating case $\partial_{u}M=0$, $h_{as}=0$
and $h_{\mathbb{S}^{d-1}}=O(r^{-1-a})$), the second term of the right
hand side of (\ref{eq:BoundednessGeneralRadiative}) can be omitted.
Furthermore, in case the $T$ vector field satisfies (\ref{eq:DeformationTensorTAway})
for $m=1$, then the second term of the right hand side of (\ref{eq:BoundednessGeneralRadiative})
is replaced by 
\begin{equation}
\int_{\mathcal{R}(\text{\textgreek{t}}_{1},\text{\textgreek{t}}_{2})\cap\{t\le T^{*}\}}\text{\textgreek{q}}\cdot\bar{t}^{-\text{\textgreek{d}}_{0}}r^{-1}\Big(|\partial_{v}\text{\textgreek{f}}|^{2}+|r^{-1}\partial_{\text{\textgreek{sv}}}\text{\textgreek{f}}|^{2}+r^{-2}|\text{\textgreek{f}}|^{2}\Big).
\end{equation}
\end{rem*}
\begin{proof}
Without loss of generality, we will assume that $\text{\textgreek{f}}$
is real valued. 

Setting $\text{\textgreek{F}}=\text{\textgreek{W}}\text{\textgreek{f}}$
and multiplying equation (\ref{eq:ConformalWaveOperator-2}) with
$\text{\textgreek{q}}\cdot\big(\partial_{v}+\partial_{u}\big)\text{\textgreek{F}}$,
we obtain after integrating over $\mathcal{R}(\text{\textgreek{t}}_{1},\text{\textgreek{t}}_{2})\cap\{t\le T^{*}\}$
(with $dudvd\text{\textgreek{sv}}$ used as a volume form): 
\begin{equation}
\begin{split}\int_{\mathcal{R}(\text{\textgreek{t}}_{1},\text{\textgreek{t}}_{2})\cap\{t\le T^{*}\}}\text{\textgreek{q}}\big(\partial_{v}+\partial_{u}\big)\text{\textgreek{F}}\cdot & \text{\textgreek{W}}\square\text{\textgreek{f}}\, dudvd\text{\textgreek{sv}}=\\
= & \int_{\mathcal{R}(\text{\textgreek{t}}_{1},\text{\textgreek{t}}_{2})\cap\{t\le T^{*}\}}\text{\textgreek{q}}\big(\partial_{v}+\partial_{u}\big)\text{\textgreek{F}}\cdot\Big\{-\big(1+O_{1}(r^{-1-a})\big)\cdot\partial_{u}\partial_{v}\text{\textgreek{F}}+r^{-2}\text{\textgreek{D}}_{g_{\mathbb{S}^{d-1}}+h_{\mathbb{S}^{d-1}}}\text{\textgreek{F}}-\\
 & \hphantom{\int_{\mathcal{R}(\text{\textgreek{t}}_{1},\text{\textgreek{t}}_{2})\cap\{t\le T^{*}\}}\text{\textgreek{q}}\big(\partial_{v}+\partial_{u}\big)\text{\textgreek{F}}\cdot\Big\{}-\frac{(d-1)(d-3)}{4}r^{-2}\cdot\text{\textgreek{F}}+Err(\text{\textgreek{F}})\Big\}\, dudvd\text{\textgreek{sv}}.
\end{split}
\label{eq:BeforeIntegrationByPartsBoundedness}
\end{equation}
Using the expression 
\[
\text{\textgreek{D}}_{g_{\mathbb{S}^{d-1}}+h_{\mathbb{S}^{d-1}}}=\text{\textgreek{D}}_{g_{\mathbb{S}^{d-1}}}+O_{1}(r^{-1})\partial_{\text{\textgreek{sv}}}\partial_{\text{\textgreek{sv}}}+O(r^{-1})\partial_{\text{\textgreek{sv}}},
\]
we obtain after integrating by parts in $\partial_{u}$, $\partial_{v}$
and $\partial_{\text{\textgreek{sv}}}$ and absorbing the error terms
in the $Err$ summand (and recalling that $\partial_{v}r=-\partial_{u}r=1$):
\begin{equation}
\begin{split}-\int_{\mathcal{R}(\text{\textgreek{t}}_{1},\text{\textgreek{t}}_{2})\cap\{t\le T^{*}\}} & \text{\textgreek{q}}\big(\partial_{v}+\partial_{u}\big)\text{\textgreek{F}}\cdot\text{\textgreek{W}}\square\text{\textgreek{f}}\, dudvd\text{\textgreek{sv}}=\\
= & \int_{\mathcal{S}_{\text{\textgreek{t}}_{2}}\cap\{t\le T^{*}\}}\frac{1}{2}\text{\textgreek{q}}\Big((1+O(r^{-1}))\cdot\big(\partial_{v}\text{\textgreek{F}}\big)^{2}+\big|r^{-1}\partial_{\text{\textgreek{sv}}}\text{\textgreek{F}}\big|^{2}+\frac{(d-1)(d-3)}{4}|r^{-1}\text{\textgreek{F}}|^{2}\Big)\, dvd\text{\textgreek{sv}}+\\
 & +\int_{\mathcal{S}_{\text{\textgreek{t}}_{2}}\cap\{t\le T^{*}\}}\frac{1}{2}\text{\textgreek{q}}\Big((1+O(r^{-1}))\cdot\big(\partial_{u}\text{\textgreek{F}}\big)^{2}+\big|r^{-1}\partial_{\text{\textgreek{sv}}}\text{\textgreek{F}}\big|^{2}+\frac{(d-1)(d-3)}{4}|r^{-1}\text{\textgreek{F}}|^{2}\Big)\, dud\text{\textgreek{sv}}+\\
 & +\int_{\{t=T^{*}\}\cap\mathcal{R}(\text{\textgreek{t}}_{1},\text{\textgreek{t}}_{2})}\frac{1}{2}\text{\textgreek{q}}\Big((1+O(r^{-1}))\cdot\big(\partial_{v}\text{\textgreek{F}}\big)^{2}+\big|r^{-1}\partial_{\text{\textgreek{sv}}}\text{\textgreek{F}}\big|^{2}+\frac{(d-1)(d-3)}{4}|r^{-1}\text{\textgreek{F}}|^{2}\Big)\, dvd\text{\textgreek{sv}}+\\
 & +\int_{\{t=T^{*}\}\cap\mathcal{R}(\text{\textgreek{t}}_{1},\text{\textgreek{t}}_{2})}\frac{1}{2}\text{\textgreek{q}}\Big((1+O(r^{-1}))\cdot\big(\partial_{u}\text{\textgreek{F}}\big)^{2}+\big|r^{-1}\partial_{\text{\textgreek{sv}}}\text{\textgreek{F}}\big|^{2}+\frac{(d-1)(d-3)}{4}|r^{-1}\text{\textgreek{F}}|^{2}\Big)\, dud\text{\textgreek{sv}}-\\
 & -\int_{\mathcal{S}_{\text{\textgreek{t}}_{1}}\cap\{t\le T^{*}\}}\frac{1}{2}\text{\textgreek{q}}\Big((1+O(r^{-1}))\cdot\big(\partial_{v}\text{\textgreek{F}}\big)^{2}+\big|r^{-1}\partial_{\text{\textgreek{sv}}}\text{\textgreek{F}}\big|^{2}+\frac{(d-1)(d-3)}{4}|r^{-1}\text{\textgreek{F}}|^{2}\Big)\, dvd\text{\textgreek{sv}}-\\
 & -\int_{\mathcal{S}_{\text{\textgreek{t}}_{1}}\cap\{t\le T^{*}\}}\frac{1}{2}\text{\textgreek{q}}\Big((1+O(r^{-1}))\cdot\big(\partial_{u}\text{\textgreek{F}}\big)^{2}+\big|r^{-1}\partial_{\text{\textgreek{sv}}}\text{\textgreek{F}}\big|^{2}+\frac{(d-1)(d-3)}{4}|r^{-1}\text{\textgreek{F}}|^{2}\Big)\, dud\text{\textgreek{sv}}+\\
 & +\int_{\mathcal{R}(\text{\textgreek{t}}_{1},\text{\textgreek{t}}_{2})\cap\{t\le T^{*}\}}\text{\textgreek{q}}\cdot\Big\{ O(r^{-1-a})|\partial_{v}\text{\textgreek{F}}|^{2}+O(r^{-1-a})|\partial_{u}\text{\textgreek{F}}|^{2}+\\
 & \hphantom{+\int_{\mathcal{R}(\text{\textgreek{t}}_{1},\text{\textgreek{t}}_{2})\cap\{t\le T^{*}\}}\text{\textgreek{q}}\cdot\Big\{}+O(r^{-1})\Big(|r^{-1}\partial_{\text{\textgreek{sv}}}\text{\textgreek{F}}|^{2}+\frac{(d-1)(d-3)}{2}|r^{-1}\text{\textgreek{F}}|^{2}\Big)\Big\}\, dudvd\text{\textgreek{sv}}+\\
 & +\int_{\mathcal{R}(\text{\textgreek{t}}_{1},\text{\textgreek{t}}_{2})\cap\{t\le T^{*}\}}\text{\textgreek{q}}\big(\partial_{v}+\partial_{u}\big)\text{\textgreek{F}}\cdot Err(\text{\textgreek{F}})\, dudvd\text{\textgreek{sv}}+\int_{\mathcal{R}(\text{\textgreek{t}}_{1},\text{\textgreek{t}}_{2})\cap\{t\le T^{*}\}}O(|\partial\text{\textgreek{q}}|)\cdot\big(|\partial\text{\textgreek{F}}|^{2}+r^{-2}|\text{\textgreek{F}}|^{2}\big)\, dudvd\text{\textgreek{sv}}.
\end{split}
\label{eq:BeforeIntegrationByPartsBoundedness-1}
\end{equation}

Since on $\mathcal{S}_{\text{\textgreek{t}}_{i}}$ we have $dvd\text{\textgreek{sv}}\sim r^{-1-\text{\textgreek{h}}'}dud\text{\textgreek{sv}}$,
we obtain from (\ref{eq:BeforeIntegrationByPartsBoundedness-1}):
\begin{equation}
\begin{split}\int_{\mathcal{S}_{\text{\textgreek{t}}_{2}}\cap\{t\le T^{*}\}}\text{\textgreek{q}}\Big(\big(\partial_{v}\text{\textgreek{F}}\big)^{2}+ & r^{-1-\text{\textgreek{h}}'}\big(\partial_{u}\text{\textgreek{F}}\big)^{2}+\big|r^{-1}\partial_{\text{\textgreek{sv}}}\text{\textgreek{F}}\big|^{2}+\frac{(d-1)(d-3)}{4}|r^{-1}\text{\textgreek{F}}|^{2}\Big)\, dvd\text{\textgreek{sv}}+\\
+\int_{\{t=T^{*}\}\cap\mathcal{R}(\text{\textgreek{t}}_{1},\text{\textgreek{t}}_{2})}\Big(\big(\partial_{v} & \text{\textgreek{F}}\big)^{2}+\big(\partial_{u}\text{\textgreek{F}}\big)^{2}+\big|r^{-1}\partial_{\text{\textgreek{sv}}}\text{\textgreek{F}}\big|^{2}+\frac{(d-1)(d-3)}{4}|r^{-1}\text{\textgreek{F}}|^{2}\Big)\, dvd\text{\textgreek{sv}}\le\\
\le C\cdot\Bigg\{ & \int_{\mathcal{S}_{\text{\textgreek{t}}_{1}}\cap\{t\le T^{*}\}}\text{\textgreek{q}}\Big(\big(\partial_{v}\text{\textgreek{F}}\big)^{2}+r^{-1-\text{\textgreek{h}}'}\big(\partial_{u}\text{\textgreek{F}}\big)^{2}+\big|r^{-1}\partial_{\text{\textgreek{sv}}}\text{\textgreek{F}}\big|^{2}+\frac{(d-1)(d-3)}{4}|r^{-1}\text{\textgreek{F}}|^{2}\Big)\, dvd\text{\textgreek{sv}}+\\
 & +\int_{\mathcal{R}(\text{\textgreek{t}}_{1},\text{\textgreek{t}}_{2})\cap\{t\le T^{*}\}}\text{\textgreek{q}}\cdot\Big\{ O(r^{-1-a})|\partial_{v}\text{\textgreek{F}}|^{2}+O(r^{-1-a})|\partial_{u}\text{\textgreek{F}}|^{2}+\\
 & \hphantom{+\int_{\mathcal{R}(\text{\textgreek{t}}_{1},\text{\textgreek{t}}_{2})\cap\{t\le T^{*}\}}\text{\textgreek{q}}\cdot\Big\{}+O(r^{-1})\Big(|r^{-1}\partial_{\text{\textgreek{sv}}}\text{\textgreek{F}}|^{2}+\frac{(d-1)(d-3)}{4}|r^{-1}\text{\textgreek{F}}|^{2}\Big)\Big\}\, dudvd\text{\textgreek{sv}}+\\
 & +\int_{\mathcal{R}(\text{\textgreek{t}}_{1},\text{\textgreek{t}}_{2})\cap\{t\le T^{*}\}}\text{\textgreek{q}}\big(\partial_{v}+\partial_{u}\big)\text{\textgreek{F}}\cdot Err(\text{\textgreek{F}})\, dudvd\text{\textgreek{sv}}+\\
 & +\int_{\mathcal{R}(\text{\textgreek{t}}_{1},\text{\textgreek{t}}_{2})\cap\{t\le T^{*}\}}O(|\partial\text{\textgreek{q}}|)\cdot\big(|\partial\text{\textgreek{F}}|^{2}+r^{-2}|\text{\textgreek{F}}|^{2}\big)\, dudvd\text{\textgreek{sv}}+\\
 & +\int_{\mathcal{R}(\text{\textgreek{t}}_{1},\text{\textgreek{t}}_{2})\cap\{t\le T^{*}\}}\text{\textgreek{q}}O(1)\cdot\big(\partial_{v}+\partial_{u}\big)\text{\textgreek{F}}\cdot\text{\textgreek{W}}\square\text{\textgreek{f}}\, dudvd\text{\textgreek{sv}}\Bigg\}.
\end{split}
\label{eq:BetterShapeForBoundedness}
\end{equation}
 Using a Hardy type inequality of the form established in Lemma \ref{lem:HardyGeneral}
(notice that it is a $1$-dimesnional Hardy inequality, since the
volume form is $dvd\text{\textgreek{sv}}$), we can bound: 
\begin{multline}
\int_{\mathcal{S}_{\text{\textgreek{t}}_{i}}\cap\{t\le T^{*}\}}\text{\textgreek{q}}r^{-2}|\text{\textgreek{F}}|^{2}\, dvd\text{\textgreek{sv}}+\int_{\mathcal{S}_{\text{\textgreek{t}}_{i}}\cap\{t=T^{*}\}}r^{-1}|\text{\textgreek{F}}|^{2}\, d\text{\textgreek{sv}}\le\\
\le C\cdot\bigg\{\int_{\mathcal{S}_{\text{\textgreek{t}}_{i}}\cap\{t\le T^{*}\}}\text{\textgreek{q}}\Big(\big(\partial_{v}\text{\textgreek{F}}\big)^{2}+r^{-1-\text{\textgreek{h}}'}\big(\partial_{u}\text{\textgreek{F}}\big)^{2}\Big)\, dvd\text{\textgreek{sv}}+\int_{\mathcal{S}_{\text{\textgreek{t}}_{i}}}|\partial\text{\textgreek{q}}|\cdot r^{-1}|\text{\textgreek{F}}|^{2}\, dvd\text{\textgreek{sv}}\Bigg\}\label{eq:HardyHyperboloid}
\end{multline}
 and 
\begin{equation}
\int_{\{t=T^{*}\}\cap\mathcal{R}(\text{\textgreek{t}}_{1},\text{\textgreek{t}}_{2})}r^{-2}|\text{\textgreek{F}}|^{2}\, dvd\text{\textgreek{sv}}\le C\cdot\Bigg\{\int_{\{t=T^{*}\}\cap\mathcal{R}(\text{\textgreek{t}}_{1},\text{\textgreek{t}}_{2})}\Big(\big(\partial_{v}\text{\textgreek{F}}\big)^{2}+\big(\partial_{u}\text{\textgreek{F}}\big)^{2}\Big)\, dvd\text{\textgreek{sv}}+\sum_{i=1}^{2}\int_{\mathcal{S}_{\text{\textgreek{t}}_{i}}\cap\{t=T^{*}\}}r^{-1}|\text{\textgreek{F}}|^{2}\, d\text{\textgreek{sv}}\Bigg\}.\label{eq:HardyPlate}
\end{equation}
 Therefore, in view of (\ref{eq:HardyHyperboloid}) and (\ref{eq:HardyPlate}),
(\ref{eq:BetterShapeForBoundedness}) can be improved into 
\begin{equation}
\begin{split}\int_{\mathcal{S}_{\text{\textgreek{t}}_{2}}\cap\{t\le T^{*}\}}\text{\textgreek{q}}\Big(\big(\partial_{v}\text{\textgreek{F}}\big)^{2}+ & r^{-1-\text{\textgreek{h}}'}\big(\partial_{u}\text{\textgreek{F}}\big)^{2}+\big|r^{-1}\partial_{\text{\textgreek{sv}}}\text{\textgreek{F}}\big|^{2}+|r^{-1}\text{\textgreek{F}}|^{2}\Big)\, dvd\text{\textgreek{sv}}+\\
+\int_{\{t=T^{*}\}\cap\mathcal{R}(\text{\textgreek{t}}_{1},\text{\textgreek{t}}_{2})}\Big(\big(\partial_{v} & \text{\textgreek{F}}\big)^{2}+\big(\partial_{u}\text{\textgreek{F}}\big)^{2}+\big|r^{-1}\partial_{\text{\textgreek{sv}}}\text{\textgreek{F}}\big|^{2}+|r^{-1}\text{\textgreek{F}}|^{2}\Big)\, dvd\text{\textgreek{sv}}\le\\
\le C\cdot\Bigg\{ & \int_{\mathcal{S}_{\text{\textgreek{t}}_{1}}\cap\{t\le T^{*}\}}\text{\textgreek{q}}\Big(\big(\partial_{v}\text{\textgreek{F}}\big)^{2}+r^{-1-\text{\textgreek{h}}'}\big(\partial_{u}\text{\textgreek{F}}\big)^{2}+\big|r^{-1}\partial_{\text{\textgreek{sv}}}\text{\textgreek{F}}\big|^{2}+|r^{-1}\text{\textgreek{F}}|^{2}\Big)\, dvd\text{\textgreek{sv}}+\\
 & +\int_{\mathcal{R}(\text{\textgreek{t}}_{1},\text{\textgreek{t}}_{2})\cap\{t\le T^{*}\}}\text{\textgreek{q}}\cdot\Big\{ O(r^{-1-a})|\partial_{v}\text{\textgreek{F}}|^{2}+O(r^{-1-a})|\partial_{u}\text{\textgreek{F}}|^{2}+\\
 & \hphantom{+\int_{\mathcal{R}(\text{\textgreek{t}}_{1},\text{\textgreek{t}}_{2})\cap\{t\le T^{*}\}}\text{\textgreek{q}}\cdot\Big\{}+O(r^{-1})\Big(|r^{-1}\partial_{\text{\textgreek{sv}}}\text{\textgreek{F}}|^{2}+\frac{(d-1)(d-3)}{4}|r^{-1}\text{\textgreek{F}}|^{2}\Big)\Big\}\, dudvd\text{\textgreek{sv}}+\\
 & +\int_{\mathcal{R}(\text{\textgreek{t}}_{1},\text{\textgreek{t}}_{2})\cap\{t\le T^{*}\}}\text{\textgreek{q}}\big(\partial_{v}+\partial_{u}\big)\text{\textgreek{F}}\cdot Err(\text{\textgreek{F}})\, dudvd\text{\textgreek{sv}}+\\
 & +\int_{\mathcal{R}(\text{\textgreek{t}}_{1},\text{\textgreek{t}}_{2})\cap\{t\le T^{*}\}}O(|\partial\text{\textgreek{q}}|)\cdot\big(|\partial\text{\textgreek{F}}|^{2}+r^{-2}|\text{\textgreek{F}}|^{2}\big)\, dudvd\text{\textgreek{sv}}+\sum_{i=1}^{2}\int_{\mathcal{S}_{\text{\textgreek{t}}_{i}}\cap\{t\le T^{*}\}}|\partial\text{\textgreek{q}}|\cdot|\text{\textgreek{F}}|^{2}\, dvd\text{\textgreek{sv}}+\\
 & +\int_{\mathcal{R}(\text{\textgreek{t}}_{1},\text{\textgreek{t}}_{2})\cap\{t\le T^{*}\}}\text{\textgreek{q}}O(1)\cdot\big(\partial_{v}+\partial_{u}\big)\text{\textgreek{F}}\cdot\text{\textgreek{W}}\square\text{\textgreek{f}}\, dudvd\text{\textgreek{sv}}\Bigg\}.
\end{split}
\label{eq:BetterShapeForBoundedness-1}
\end{equation}

Using the expression (\ref{eq:ErrorTermsMorawetz}) for $Err(\text{\textgreek{F}})$,
we can readily bound after integrating by parts in the highest order
terms (and in the $\text{\textgreek{F}}\partial_{u}\text{\textgreek{F}}=\frac{1}{2}\partial_{u}(\text{\textgreek{F}}^{2})$
term) and using a Cauchy--Schwarz inequality: 
\begin{equation}
\begin{split}\int_{\mathcal{R}(\text{\textgreek{t}}_{1},\text{\textgreek{t}}_{2})\cap\{t\le T^{*}\}}\text{\textgreek{q}}\big(\partial_{v}+\partial_{u}\big)\text{\textgreek{F}}\cdot & Err(\text{\textgreek{F}})\, dudvd\text{\textgreek{sv}}\le\\
\le & C_{\text{\textgreek{h}}}\int_{\mathcal{R}(\text{\textgreek{t}}_{1},\text{\textgreek{t}}_{2})\cap\{t\le T^{*}\}}\text{\textgreek{q}}\Big\{ O(r^{-1-a})|\partial_{u}\text{\textgreek{F}}|^{2}+O(r^{-1})\Big(|\partial_{v}\text{\textgreek{F}}|^{2}+\text{\textgreek{d}}_{1}|r^{-1}\partial_{\text{\textgreek{sv}}}\text{\textgreek{F}}|^{2}+r^{-2}|\text{\textgreek{F}}|^{2}\Big)\Big\}\, dudvd\text{\textgreek{sv}}+\\
 & +C_{\text{\textgreek{h}}}\Big(\int_{\mathcal{R}(\text{\textgreek{t}}_{1},\text{\textgreek{t}}_{2})\cap\{t\le T^{*}\}}|\partial\text{\textgreek{q}}|\big(|\partial\text{\textgreek{F}}|^{2}+r^{-2}|\text{\textgreek{F}}|^{2}\big)\, dudvd\text{\textgreek{sv}}+\\
 & +\int_{\mathcal{S}_{\text{\textgreek{t}}_{2}}\cap\{t\le T^{*}\}}\text{\textgreek{q}}O(r^{-a})\Big(\big(\partial_{v}\text{\textgreek{F}}\big)^{2}+r^{-1-\text{\textgreek{h}}'}\big(\partial_{u}\text{\textgreek{F}}\big)^{2}+\big|r^{-1}\partial_{\text{\textgreek{sv}}}\text{\textgreek{F}}\big|^{2}+|r^{-1}\text{\textgreek{F}}|^{2}\Big)\, dvd\text{\textgreek{sv}}+\\
 & +\int_{\{t=T^{*}\}\cap\mathcal{R}(\text{\textgreek{t}}_{1},\text{\textgreek{t}}_{2})}O(r^{-a})\Big(\big(\partial_{v}\text{\textgreek{F}}\big)^{2}+\big(\partial_{u}\text{\textgreek{F}}\big)^{2}+\big|r^{-1}\partial_{\text{\textgreek{sv}}}\text{\textgreek{F}}\big|^{2}+|r^{-1}\text{\textgreek{F}}|^{2}\Big)\, dvd\text{\textgreek{sv}}+\\
 & +\int_{\mathcal{S}_{\text{\textgreek{t}}_{1}}\cap\{t\le T^{*}\}}\text{\textgreek{q}}O(r^{-a})\Big(\big(\partial_{v}\text{\textgreek{F}}\big)^{2}+r^{-1-\text{\textgreek{h}}'}\big(\partial_{u}\text{\textgreek{F}}\big)^{2}+\big|r^{-1}\partial_{\text{\textgreek{sv}}}\text{\textgreek{F}}\big|^{2}+|r^{-1}\text{\textgreek{F}}|^{2}\Big)\, dvd\text{\textgreek{sv}}.
\end{split}
\label{eq:BoundsErrorBoundedness}
\end{equation}

Therefore, from (\ref{eq:BetterShapeForBoundedness-1}), (\ref{eq:BoundsErrorBoundedness}),
(\ref{eq:HardyHyperboloid}) and the trace inequality 
\begin{equation}
\sum_{i=1}^{2}\int_{\mathcal{S}_{\text{\textgreek{t}}_{i}}\cap\{t\le T^{*}\}}|\partial\text{\textgreek{q}}|\cdot r^{-1}|\text{\textgreek{F}}|^{2}\, dvd\text{\textgreek{sv}}\le\int_{\mathcal{R}(\text{\textgreek{t}}_{1},\text{\textgreek{t}}_{2})\cap\{t\le T^{*}\}}O(|\partial\text{\textgreek{q}}|)\cdot\big(r^{-1}|\partial\text{\textgreek{F}}|^{2}+r^{-1}|\text{\textgreek{F}}|^{2}\big)\, dudvd\text{\textgreek{sv}},
\end{equation}
we obtain if $R$ is large enough in terms of the geometry of $(\mathcal{N}_{af},g)$:
\begin{equation}
\begin{split}\int_{\mathcal{S}_{\text{\textgreek{t}}_{2}}\cap\{t\le T^{*}\}} & \text{\textgreek{q}}\cdot J_{\text{\textgreek{m}}}^{T}(\text{\textgreek{f}})\bar{n}^{\text{\textgreek{m}}}+\int_{\{t=T^{*}\}\cap\mathcal{R}(\text{\textgreek{t}}_{1},\text{\textgreek{t}}_{2})}\Big(J_{\text{\textgreek{m}}}^{T}(\text{\textgreek{f}})n^{\text{\textgreek{m}}}+r^{-2}|\text{\textgreek{f}}|^{2}\Big)\le\\
\le C\cdot\Bigg\{ & \int_{\mathcal{S}_{\text{\textgreek{t}}_{1}}\cap\{t\le T^{*}\}}\text{\textgreek{q}}\cdot J_{\text{\textgreek{m}}}^{T}(\text{\textgreek{f}})\bar{n}^{\text{\textgreek{m}}}+\int_{\mathcal{R}(\text{\textgreek{t}}_{1},\text{\textgreek{t}}_{2})\cap\{t\le T^{*}\}}\text{\textgreek{q}}\cdot\Big\{ r^{-1-a}J_{\text{\textgreek{m}}}^{T}(\text{\textgreek{f}})n^{\text{\textgreek{m}}}+r^{-1}\Big(|\partial_{v}\text{\textgreek{f}}|^{2}+|r^{-1}\partial_{\text{\textgreek{sv}}}\text{\textgreek{f}}|^{2}+|r^{-1}\text{\textgreek{f}}|^{2}\Big)\Big\}+\\
 & +\int_{\mathcal{R}(\text{\textgreek{t}}_{1},\text{\textgreek{t}}_{2})\cap\{t\le T^{*}\}}O(|\partial\text{\textgreek{q}}|)\cdot\big(|\partial\text{\textgreek{f}}|^{2}+r^{-1}|\text{\textgreek{f}}|^{2}\big)+\int_{\mathcal{R}(\text{\textgreek{t}}_{1},\text{\textgreek{t}}_{2})\cap\{t\le T^{*}\}}\text{\textgreek{q}}O(1)\cdot Re\Big\{\text{\textgreek{W}}^{-1}\big(\partial_{v}+\partial_{u}\big)(\text{\textgreek{W}}\bar{\text{\textgreek{f}}})\cdot\square\text{\textgreek{f}}\Big\}\Bigg\}.
\end{split}
\label{eq:BoundednessGeneralRadiative-2}
\end{equation}

Using (\ref{eq:MorawetzGeneralCaseRadiativeHyperboloids}) for some
fixed $\text{\textgreek{h}}<a$, we thus obtain (\ref{eq:BoundednessGeneralRadiative})
from (\ref{eq:BoundednessGeneralRadiative-2}) provided $R$ is large
in terms of the geometry of $(\mathcal{N}_{af},g)$.
\end{proof}
We can also establish the following generalisation of Lemma \ref{lem:Boundedness Hyperboloids}:
\begin{lem}
\label{lem:BoundednessGeometric}For any smooth function $\text{\textgreek{f}}:\mathcal{N}_{af}\rightarrow\mathbb{C}$,
any $\text{\textgreek{q}}:\mathcal{N}_{af}\rightarrow[0,1]$ supported
on $\{r\ge R\}$ for some $R>0$ large in terms of the geometry of
$(\mathcal{N}_{af},g)$, and any two smooth, spacelike hypersurfaces
$\mathcal{S}_{1},\mathcal{S}_{2}$ of $\mathcal{N}_{af}$ intersecting
the region $\{r\le R\}$ such that $\mathcal{S}_{2}$ lies in the
future domain of dependence of $\mathcal{S}_{1}\cup\{r\le R\}$, the
following estimate is true: 
\begin{align}
\int_{\mathcal{S}_{2}}\text{\textgreek{q}}\cdot J_{\text{\textgreek{m}}}^{T}(\text{\textgreek{f}})n_{\mathcal{S}_{2}}^{\text{\textgreek{m}}}\le C\cdot\Bigg\{ & \int_{\mathcal{S}_{1}}\text{\textgreek{q}}\cdot J_{\text{\textgreek{m}}}^{T}(\text{\textgreek{f}})n_{\mathcal{S}_{1}}^{\text{\textgreek{m}}}+\int_{J^{+}(\mathcal{S}_{1})\cap J^{-}(\mathcal{S}_{2})}\text{\textgreek{q}}\cdot r^{-1}\Big(|\partial_{v}\text{\textgreek{f}}|^{2}+|r^{-1}\partial_{\text{\textgreek{sv}}}\text{\textgreek{f}}|^{2}+|r^{-1}\text{\textgreek{f}}|^{2}\Big)+\label{eq:BoundednessGeneralRadiative-1}\\
 & +\int_{J^{+}(\mathcal{S}_{1})\cap J^{-}(\mathcal{S}_{2})}O(|\partial\text{\textgreek{q}}|)\cdot\big(|\partial\text{\textgreek{f}}|^{2}+r^{-1}|\text{\textgreek{f}}|^{2}\big)+\nonumber \\
 & +\int_{J^{+}(\mathcal{S}_{1})\cap J^{-}(\mathcal{S}_{2})}\text{\textgreek{q}}O(1)\cdot Re\Big\{\text{\textgreek{W}}^{-1}\big(\partial_{v}+\partial_{u}\big)(\text{\textgreek{W}}\bar{\text{\textgreek{f}}})\cdot\square\text{\textgreek{f}}\Big\}\Bigg\}.\nonumber 
\end{align}
 for some constant $C>0$ depending only on the geometry of $(\mathcal{N}_{af},g)$
and the precise form of $\mathcal{S}_{1},\mathcal{S}_{2}$ but independent
of translation of these hypersurfaces by the flow of $T=\partial_{u}+\partial_{v}$.
In the above, $n_{\mathcal{S}_{i}}$ is the future directed unit normal
of $\mathcal{S}_{i}$. Integration over $\mathcal{S}_{i}$ is performed
using the volume form of the induced metric, while integration over
$J^{+}(\mathcal{S}_{1})\cap J^{-}(\mathcal{S}_{2})$ is performed
using the natural volume form of $g$.\end{lem}
\begin{rem*}
Again, in case the radiative components of the metric satisfy the
bounds $\partial_{u}M\le0$ and $|\partial_{u}h_{as}|+|r\partial_{u}h_{\mathbb{S}^{d-1}}|\ll-(\partial_{u}M)+O(r^{-a})$
(which includes the non-radiating case $\partial_{u}M=0$, $h_{as}=0$
and $h_{\mathbb{S}^{d-1}}=O(r^{-1-a})$), the second term of the right
hand side of (\ref{eq:BoundednessGeneralRadiative-1}) can be omitted.
Furthermore, in case the $T$ vector field  satisfies (\ref{eq:DeformationTensorTAway})
for $m=1$, then the second term of the right hand side of (\ref{eq:BoundednessGeneralRadiative-1})
is replaced by 
\begin{equation}
\int_{J^{+}(\mathcal{S}_{1})\cap J^{-}(\mathcal{S}_{2})}\text{\textgreek{q}}\cdot\bar{t}^{-\text{\textgreek{d}}_{0}}r^{-1}\Big(|\partial_{v}\text{\textgreek{f}}|^{2}+|r^{-1}\partial_{\text{\textgreek{sv}}}\text{\textgreek{f}}|^{2}+r^{-2}|\text{\textgreek{f}}|^{2}\Big).
\end{equation}
\end{rem*}
\begin{proof}
The proof of Lemma \ref{lem:BoundednessGeometric} follows in exactly
the same way as that of Lemma \ref{lem:Boundedness Hyperboloids},
by integrating $\text{\textgreek{q}}\big(\partial_{v}+\partial_{u}\big)\text{\textgreek{F}}\cdot\text{\textgreek{W}}\square\text{\textgreek{f}}$
over $J^{+}(\mathcal{S}_{1})\cap J^{-}(\mathcal{S}_{2})$ in place
of $\mathcal{R}(\text{\textgreek{t}}_{1},\text{\textgreek{t}}_{2})$.
Hence, the details will be omitted.
\end{proof}
Finally, we will also establish the following improvement of Lemma
\ref{lem:Boundedness Hyperboloids} for higher order derivatives of
$\text{\textgreek{f}}$, which will be used in Section \ref{sec:Firstdecay}: 
\begin{lem}
\label{lem:Boundedness Hyperboloids Higher Order}For any smooth function
$\text{\textgreek{f}}:\mathcal{N}_{af}\rightarrow\mathbb{C}$, any
$\text{\textgreek{q}}:\mathcal{N}_{af}\rightarrow[0,1]$ supported
on $\{r\ge R\}$ for some $R>0$ large in terms of the geometry of
$(\mathcal{N}_{af},g)$, any $\text{\textgreek{t}}_{1}\le\text{\textgreek{t}}_{2}$
and any $T^{*}>0$ the following estimate is true: 
\begin{equation}
\begin{split}\sum_{j=0}^{k}\sum_{j_{1}+j_{2}+j_{3}=j}\Bigg\{\int_{\mathcal{S}_{\text{\textgreek{t}}_{2}}\cap\{t\le T^{*}\}} & \text{\textgreek{q}}\cdot r^{2j_{1}}J_{\text{\textgreek{m}}}^{T}\big(\text{\textgreek{W}}^{-1}\partial_{v}^{j_{1}}\partial_{\text{\textgreek{sv}}}^{j_{2}}\partial_{u}^{j_{3}}(\text{\textgreek{W}}\text{\textgreek{f}})\big)\bar{n}^{\text{\textgreek{m}}}+\int_{\{t=T^{*}\}\cap\mathcal{R}(\text{\textgreek{t}}_{1},\text{\textgreek{t}}_{2})}r^{2j_{1}}\Big(J_{\text{\textgreek{m}}}^{T}\big(\text{\textgreek{W}}^{-1}\partial_{v}^{j_{1}}\partial_{\text{\textgreek{sv}}}^{j_{2}}\partial_{u}^{j_{3}}(\text{\textgreek{W}}\text{\textgreek{f}})\big)n^{\text{\textgreek{m}}}\Big)\Bigg\}\le\\
\le C_{k}\cdot\sum_{j=0}^{k}\sum_{j_{1}+j_{2}+j_{3}=j}\Bigg\{ & \int_{\mathcal{S}_{\text{\textgreek{t}}_{1}}\cap\{t\le T^{*}\}}\text{\textgreek{q}}\cdot r^{2j_{1}}J_{\text{\textgreek{m}}}^{T}\big(\text{\textgreek{W}}^{-1}\partial_{v}^{j_{1}}\partial_{\text{\textgreek{sv}}}^{j_{2}}\partial_{u}^{j_{3}}(\text{\textgreek{W}}\text{\textgreek{f}})\big)\bar{n}^{\text{\textgreek{m}}}+\\
 & +\int_{\mathcal{R}(\text{\textgreek{t}}_{1},\text{\textgreek{t}}_{2})\cap\{t\le T^{*}\}}\text{\textgreek{q}}\cdot r^{2j_{1}-1}\text{\textgreek{W}}^{-2}\Big(|\partial_{v}^{j_{1}+1}\partial_{\text{\textgreek{sv}}}^{j_{2}}\partial_{u}^{j_{3}}(\text{\textgreek{W}}\text{\textgreek{f}})|^{2}+r^{-2}|\partial_{v}^{j_{1}}\partial_{\text{\textgreek{sv}}}^{j_{2}+1}\partial_{u}^{j_{3}}(\text{\textgreek{W}}\text{\textgreek{f}})|^{2}+r^{-2}|\partial_{v}^{j_{1}}\partial_{\text{\textgreek{sv}}}^{j_{2}}\partial_{u}^{j_{3}}(\text{\textgreek{W}}\text{\textgreek{f}})|^{2}\Big)+\\
 & +\int_{\mathcal{R}(\text{\textgreek{t}}_{1},\text{\textgreek{t}}_{2})\cap\{t\le T^{*}\}}O(|\partial\text{\textgreek{q}}|)\cdot r^{2j_{1}}\text{\textgreek{W}}^{-2}\big(|\partial\partial_{v}^{j_{1}}\partial_{\text{\textgreek{sv}}}^{j_{2}}\partial_{u}^{j_{3}}(\text{\textgreek{W}}\text{\textgreek{f}})|^{2}+r^{-1}|\partial_{v}^{j_{1}}\partial_{\text{\textgreek{sv}}}^{j_{2}}\partial_{u}^{j_{3}}(\text{\textgreek{W}}\text{\textgreek{f}})|^{2}\big)+\\
 & +\int_{\mathcal{R}(\text{\textgreek{t}}_{1},\text{\textgreek{t}}_{2})\cap\{t\le T^{*}\}}\text{\textgreek{q}}O(1)r^{2j_{1}}\text{\textgreek{W}}^{-2}\cdot Re\Big\{\big(\partial_{v}+\partial_{u}\big)(\partial_{v}^{j_{1}}\partial_{\text{\textgreek{sv}}}^{j_{2}}\partial_{u}^{j_{3}}(\text{\textgreek{W}}\bar{\text{\textgreek{f}}}))\cdot\partial_{v}^{j_{1}}\partial_{\text{\textgreek{sv}}}^{j_{2}}\partial_{u}^{j_{3}}(\text{\textgreek{W}}\square\text{\textgreek{f}})\Big\}\Bigg\}.
\end{split}
\label{eq:BoundednessGeneralRadiativeHigherOrderWeighted}
\end{equation}
 for some constant $C_{k}>0$ depending only on the geometry of $(\mathcal{N}_{af},g)$
and $k$. \end{lem}
\begin{rem*}
Again, in case the radiative components of the metric satisfy the
bounds $\partial_{u}M\le0$ and $|\partial_{u}h_{as}|+|r\partial_{u}h_{\mathbb{S}^{d-1}}|\ll-(\partial_{u}M)+O(r^{-a})$
(which includes the non-radiating case $\partial_{u}M=0$, $h_{as}=0$
and $h_{\mathbb{S}^{d-1}}=O(r^{-1-a})$), the second term of the right
hand side of (\ref{eq:BoundednessGeneralRadiativeHigherOrderWeighted})
can be omitted. Furthermore, in case the $T$ vector field  satisfies
(\ref{eq:DeformationTensorTAway}) for $m=1$, then the integrand
of the second term of the right hand side of (\ref{eq:BoundednessGeneralRadiativeHigherOrderWeighted})
has an extra $\bar{t}^{-\text{\textgreek{d}}_{0}}$ factor.\end{rem*}
\begin{proof}
The proof of (\ref{eq:BoundednessGeneralRadiativeHigherOrderWeighted})
follows in exactly the same way as the proof of (\ref{eq:BoundednessGeneralRadiative}),
by integrating (for real $\text{\textgreek{f}}$) for any $0\le j\le k$
and any partition $j_{1},j_{2},j_{3}$ of $j$: 
\begin{equation}
\text{\textgreek{q}}\cdot\big(\partial_{v}+\partial_{u}\big)(r^{j_{1}}\partial_{v}^{j_{1}}\partial_{\text{\textgreek{sv}}}^{j_{2}}\partial_{u}^{j_{3}}(\text{\textgreek{W}}\text{\textgreek{f}}))\cdot\text{\textgreek{W}}\square\big(\text{\textgreek{W}}^{-1}r^{j_{1}}\partial_{v}^{j_{1}}\partial_{\text{\textgreek{sv}}}^{j_{2}}\partial_{u}^{j_{3}}(\text{\textgreek{W}}\text{\textgreek{f}})\big)
\end{equation}
 over $\mathcal{R}(\text{\textgreek{t}}_{1},\text{\textgreek{t}}_{2})\cap\{t\le T^{*}\}$
with $dudvd\text{\textgreek{sv}}$ as volume form, and then summing
the resulting estimates over all possible $j_{1},j_{2},j_{3}$, using
Lemma \ref{lem:Commutator expressions} to obtain an expression the
commutator of $\square$ with $\partial_{v},\partial_{\text{\textgreek{sv}}},\partial_{u}$.
\end{proof}

\section{\label{sec:The-new-method}The extension of the $r^{p}$-weighted
energy hierarchy}

In this section, we will generalise the $r^{p}$-weighted energy hierarchy
introduced in \cite{DafRod7}, so as to apply to the asymptotically
flat region of a broad class of stationary and asymptotically flat
spacetimes, modeled on the spacetime $(\mathcal{N}_{af},g)$. This
will be achieved by repeating the main steps of the method of \cite{DafRod7}
in the the case of Minkowski spacetime and controlling the error terms
with the use of Lemmas \ref{lem:MorawetzDrLemmaHyperboloids} and
\ref{lem:Boundedness Hyperboloids}.

\subsection{Statement of the $r^{p}$-weighted energy hierarchy}

We will establish the following decay estimate for solutions to $\square\text{\textgreek{f}}=F$
on $(\mathcal{N}_{af},g)$: 
\begin{thm}
\label{thm:NewMethodFinalStatementHyperboloids}For any $0<p\le2$,
any given $0<\text{\textgreek{h}}<a$ and $0<\text{\textgreek{d}}<1$,
any $R>0$ large enough in terms of $p,\text{\textgreek{h}},\text{\textgreek{d}}$
and any $\text{\textgreek{t}}_{1}\le\text{\textgreek{t}}_{2}$, the
following inequality is true for any smooth function $\text{\textgreek{f}}:\mathcal{N}_{af}\rightarrow\mathbb{C}$
and any smooth cut-off $\text{\textgreek{q}}_{R}:\mathcal{N}_{af}\rightarrow[0,1]$
supported in $\{r\ge R\}$:
\begin{equation}
\begin{split}\mathcal{E}_{bound,R;\text{\textgreek{d}}}^{(p)}[\text{\textgreek{f}}](\text{\textgreek{t}}_{2}) & +\int_{\text{\textgreek{t}}_{1}}^{\text{\textgreek{t}}_{2}}\mathcal{E}_{bulk,R,\text{\textgreek{h}};\text{\textgreek{d}}}^{(p-1)}[\text{\textgreek{f}}](\text{\textgreek{t}})\, d\text{\textgreek{t}}+\limsup_{T^{*}\rightarrow+\infty}\mathcal{E}_{\mathcal{I}^{+},T^{*};\text{\textgreek{d}}}^{(p)}[\text{\textgreek{f}}](\text{\textgreek{t}}_{1},\text{\textgreek{t}}_{2})\lesssim_{p,\text{\textgreek{h},\textgreek{d}}}\\
\lesssim_{p,\text{\textgreek{h},\textgreek{d}}} & \mathcal{E}_{bound,R;\text{\textgreek{d}}}^{(p)}[\text{\textgreek{f}}](\text{\textgreek{t}}_{1})+\int_{\mathcal{R}(t_{1},t_{2})}|\partial\text{\textgreek{q}}_{R}|\cdot\big(r^{p}|\partial\text{\textgreek{f}}|^{2}+r^{p-2}\cdot|\text{\textgreek{f}}|^{2}\big)+\\
 & +\int_{\mathcal{R}(t_{1},t_{2})}\text{\textgreek{q}}_{R}\cdot(r^{p+1}+r^{1+\text{\textgreek{h}}})\cdot|\square_{g}\text{\textgreek{f}}|^{2}\,\text{\textgreek{W}}^{2}dudvd\text{\textgreek{sv}}.
\end{split}
\label{eq:newMethodFinalStatementHyperboloids}
\end{equation}
 In the above, the constants implicit in the $\lesssim_{p,\text{\textgreek{h},\textgreek{d}}}$
notation depend only on $p,\text{\textgreek{h}},\text{\textgreek{d}}$
and on the geometry of $(\mathcal{N}_{af},g)$, and the $p$-weighted
energies are defined as follows: 
\begin{align}
\mathcal{E}_{bound,R;\text{\textgreek{d}}}^{(p)}[\text{\textgreek{f}}](\text{\textgreek{t}})= & \int_{\mathcal{S}_{\text{\textgreek{t}}}}\text{\textgreek{q}}_{R}\Big(r^{p}\big|\partial_{v}(\text{\textgreek{W}}\text{\textgreek{f}})\big|^{2}+r^{-1-\text{\textgreek{h}}'}\big(r^{p}\big|r^{-1}\partial_{\text{\textgreek{sv}}}(\text{\textgreek{W}}\text{\textgreek{f}})\big|^{2}+\big((d-3)r^{p-2}+\min\{r^{p-2},r^{-\text{\textgreek{d}}}\}\big)\big|\text{\textgreek{W}}\text{\textgreek{f}}\big|^{2}\Big)\, dvd\text{\textgreek{sv}}+\\
 & +\int_{\mathcal{S}_{\text{\textgreek{t}}}}\text{\textgreek{q}}_{R}J_{\text{\textgreek{m}}}^{T}(\text{\textgreek{f}})\bar{n}^{\text{\textgreek{m}}},\nonumber 
\end{align}
\begin{align}
\mathcal{E}_{bulk,R,\text{\textgreek{h}};\text{\textgreek{d}}}^{(p-1)}[\text{\textgreek{f}}](\text{\textgreek{t}})= & \int_{\mathcal{S}_{\text{\textgreek{t}}}}\text{\textgreek{q}}_{R}\Big(pr^{p-1}\big|\partial_{v}(\text{\textgreek{W}}\text{\textgreek{f}})\big|^{2}+\big\{\big((2-p)r^{p-1}+r^{p-1-\text{\textgreek{d}}}\big)\big|r^{-1}\partial_{\text{\textgreek{sv}}}(\text{\textgreek{W}}\text{\textgreek{f}})\big|^{2}+\\
 & \hphantom{\int_{\mathcal{S}_{\text{\textgreek{t}}}}\text{\textgreek{q}}_{R}\Big(}+\big((2-p)(d-3)r^{p-3}+\min\{r^{p-3},r^{-1-\text{\textgreek{d}}}\}\big)\big|\text{\textgreek{W}}\text{\textgreek{f}}\big|^{2}\big\}\Big)\, dvd\text{\textgreek{sv}}+\nonumber \\
 & +\int_{\mathcal{S}_{\text{\textgreek{t}}}}\text{\textgreek{q}}_{R}r^{-1-\text{\textgreek{h}}}|\partial_{u}(\text{\textgreek{W}}\text{\textgreek{f}})|^{2}\, dvd\text{\textgreek{sv}}\nonumber 
\end{align}
and \textup{
\begin{align}
\mathcal{E}_{\mathcal{I}^{+},T^{*},\text{\textgreek{d}}}^{(p)}[\text{\textgreek{f}}](\text{\textgreek{t}}_{1},\text{\textgreek{t}}_{2})= & \int_{\mathcal{R}(\text{\textgreek{t}}_{1},\text{\textgreek{t}}_{2})\cap\{t=T^{*}\}}\Big(r^{p}\big|r^{-1}\partial_{\text{\textgreek{sv}}}(\text{\textgreek{W}}\text{\textgreek{f}})\big|^{2}+\big((d-3)r^{p-2}+\min\{r^{p-2},r^{-\text{\textgreek{d}}}\}\big)\big|\text{\textgreek{W}}\text{\textgreek{f}}\big|^{2}\Big)\, dvd\text{\textgreek{sv}}+\\
 & +\int_{\mathcal{R}(\text{\textgreek{t}}_{1},\text{\textgreek{t}}_{2})\cap\{t=T^{*}\}}J_{\text{\textgreek{m}}}^{T}(\text{\textgreek{f}})n^{\text{\textgreek{m}}}.\nonumber 
\end{align}
}\end{thm}
\begin{rem*}
In fact, it follows from the proof of Theorem \ref{thm:NewMethodFinalStatementHyperboloids}
that (\ref{eq:newMethodFinalStatementHyperboloids}) can be improved
into: 
\begin{equation}
\begin{split}\Bigg|\Big(\mathcal{E}_{bound,R,\text{\textgreek{e}};\text{\textgreek{d}}}^{(p)}[\text{\textgreek{f}}](\text{\textgreek{t}}_{2}) & +\int_{\text{\textgreek{t}}_{1}}^{\text{\textgreek{t}}_{2}}\mathcal{E}_{bulk,R,\text{\textgreek{h}},\text{\textgreek{e}};\text{\textgreek{d}}}^{(p-1)}[\text{\textgreek{f}}](\text{\textgreek{t}})\, d\text{\textgreek{t}}+\limsup_{T^{*}\rightarrow+\infty}\mathcal{E}_{\mathcal{I}^{+},T^{*},\text{\textgreek{e}};\text{\textgreek{d}}}^{(p)}[\text{\textgreek{f}}](\text{\textgreek{t}}_{1},\text{\textgreek{t}}_{2})\Big)-\mathcal{E}_{bound,R,\text{\textgreek{e}};\text{\textgreek{d}}}^{(p)}[\text{\textgreek{f}}](\text{\textgreek{t}}_{1})\Bigg|\le\\
\le & C_{p,\text{\textgreek{h}},\text{\textgreek{d}}}\cdot\text{\textgreek{e}}\cdot\mathcal{E}_{bound,R,\text{\textgreek{e}};\text{\textgreek{d}}}^{(p)}[\text{\textgreek{f}}](\text{\textgreek{t}}_{1})+C_{p,\text{\textgreek{h}},\text{\textgreek{d}},\text{\textgreek{e}}}\int_{\mathcal{R}(t_{1},t_{2})}|\partial\text{\textgreek{q}}_{R}|\cdot\big(r^{p}|\partial\text{\textgreek{f}}|^{2}+r^{p-2}\cdot|\text{\textgreek{f}}|^{2}\big)+\\
 & +C_{p,\text{\textgreek{h}},\text{\textgreek{d}},\text{\textgreek{e}}}\int_{\mathcal{R}(t_{1},t_{2})}\text{\textgreek{q}}_{R}\cdot(r^{p+1}+r^{1+\text{\textgreek{h}}})\cdot|\square_{g}\text{\textgreek{f}}|^{2}\,\text{\textgreek{W}}^{2}dudvd\text{\textgreek{sv}}
\end{split}
\label{eq:newMethodFinalStatementHyperboloidsImp}
\end{equation}
 for any $0<\text{\textgreek{d}}<a$ and $0<\text{\textgreek{e}}<1$,
where 
\begin{align}
\mathcal{E}_{bound,R,\text{\textgreek{e}};\text{\textgreek{d}}}^{(p)}[\text{\textgreek{f}}](\text{\textgreek{t}})= & \int_{\mathcal{S}_{\text{\textgreek{t}}}}\text{\textgreek{q}}_{R}\Big(r^{p}\big|\partial_{v}(\text{\textgreek{W}}\text{\textgreek{f}})\big|^{2}+r^{-1-\text{\textgreek{h}}'}\big(r^{p}\big|r^{-1}\partial_{\text{\textgreek{sv}}}(\text{\textgreek{W}}\text{\textgreek{f}})\big|^{2}+\big((d-3)r^{p-2}+\text{\textgreek{e}}\cdot\min\{r^{p-2},r^{-\text{\textgreek{d}}}\}\big)\big|\text{\textgreek{W}}\text{\textgreek{f}}\big|^{2}\Big)\, dvd\text{\textgreek{sv}}+\\
 & +\int_{\mathcal{S}_{\text{\textgreek{t}}}}\text{\textgreek{q}}_{R}J_{\text{\textgreek{m}}}^{T}(\text{\textgreek{f}})\bar{n}^{\text{\textgreek{m}}},\nonumber 
\end{align}
\begin{align}
\mathcal{E}_{bulk,R,\text{\textgreek{h}},\text{\textgreek{e}};\text{\textgreek{d}}}^{(p-1)}[\text{\textgreek{f}}](\text{\textgreek{t}})= & \int_{\mathcal{S}_{\text{\textgreek{t}}}}\text{\textgreek{q}}_{R}\Big(pr^{p-1}\big|\partial_{v}(\text{\textgreek{W}}\text{\textgreek{f}})\big|^{2}+\big\{\big((2-p)r^{p-1}+r^{p-1-\text{\textgreek{d}}}\big)\big|r^{-1}\partial_{\text{\textgreek{sv}}}(\text{\textgreek{W}}\text{\textgreek{f}})\big|^{2}+\\
 & \hphantom{\int_{\mathcal{S}_{\text{\textgreek{t}}}}\text{\textgreek{q}}_{R}\Big(}+\big((2-p)(d-3)r^{p-3}+\text{\textgreek{e}}\cdot\min\{r^{p-3},r^{-1-\text{\textgreek{d}}}\}\big)\big|\text{\textgreek{W}}\text{\textgreek{f}}\big|^{2}\big\}\Big)\, dvd\text{\textgreek{sv}}+\nonumber \\
 & +\int_{\mathcal{S}_{\text{\textgreek{t}}}}\text{\textgreek{q}}_{R}r^{-1-\text{\textgreek{h}}}|\partial_{u}(\text{\textgreek{W}}\text{\textgreek{f}})|^{2}\, dvd\text{\textgreek{sv}}\nonumber 
\end{align}
and 
\begin{align}
\mathcal{E}_{\mathcal{I}^{+},T^{*},\text{\textgreek{e}};\text{\textgreek{d}}}^{(p)}[\text{\textgreek{f}}](\text{\textgreek{t}}_{1},\text{\textgreek{t}}_{2})= & \int_{\mathcal{R}(\text{\textgreek{t}}_{1},\text{\textgreek{t}}_{2})\cap\{t=T^{*}\}}\Big(r^{p}\big|r^{-1}\partial_{\text{\textgreek{sv}}}(\text{\textgreek{W}}\text{\textgreek{f}})\big|^{2}+\big((d-3)r^{p-2}+\text{\textgreek{e}}\cdot\min\{r^{p-2},r^{-\text{\textgreek{d}}}\}\big)\big|\text{\textgreek{W}}\text{\textgreek{f}}\big|^{2}\Big)\, dvd\text{\textgreek{sv}}+\\
 & +\int_{\mathcal{R}(\text{\textgreek{t}}_{1},\text{\textgreek{t}}_{2})\cap\{t=T^{*}\}}J_{\text{\textgreek{m}}}^{T}(\text{\textgreek{f}})n^{\text{\textgreek{m}}}.\nonumber 
\end{align}

\end{rem*}
Using Lemma \ref{lem:Commutator expressions} on the commutator of
$\square$ with $\partial_{v}$, $\partial_{\text{\textgreek{sv}}}$
and $\partial_{u}$, we can also obtain the following corollary of
Theorem \ref{thm:NewMethodFinalStatementHyperboloids}: 
\begin{cor}
\label{cor:NewMethodHigherDerivativesNotImproved}Using the notations
of Theorem \ref{thm:NewMethodFinalStatementHyperboloids}, for any
$0<p\le2$, any integer $m\ge0$, any given $0<\text{\textgreek{h}}<a$
and $0<\text{\textgreek{d}}<1$, any $R>0$ large enough in terms
of $p,\text{\textgreek{h}},\text{\textgreek{d}}$ and any $\text{\textgreek{t}}_{1}\le\text{\textgreek{t}}_{2}$,
the following inequality is true for any smooth function $\text{\textgreek{f}}:\mathcal{N}_{af}\rightarrow\mathbb{C}$
and any smooth cut-off $\text{\textgreek{q}}_{R}:\mathcal{N}_{af}\rightarrow[0,1]$
supported in $\{r\ge R\}$: 
\begin{equation}
\begin{split}\sum_{j=0}^{m-1}\sum_{j_{1}+j_{2}+j_{3}=j}\Big\{\mathcal{E}_{bound,R;\text{\textgreek{d}}}^{(p)}[r^{-j_{2}}\partial_{v}^{j_{1}}\partial_{\text{\textgreek{sv}}}^{j_{2}}\partial_{u}^{j_{3}}\text{\textgreek{f}}](\text{\textgreek{t}}_{2})+\int_{\text{\textgreek{t}}_{1}}^{\text{\textgreek{t}}_{2}}\mathcal{E}_{bulk,R,\text{\textgreek{h}};\text{\textgreek{d}}}^{(p-1)} & [r^{-j_{2}}\partial_{v}^{j_{1}}\partial_{\text{\textgreek{sv}}}^{j_{2}}\partial_{u}^{j_{3}}\text{\textgreek{f}}](\text{\textgreek{t}})\, d\text{\textgreek{t}}\\
+\limsup_{T^{*}\rightarrow+\infty}\mathcal{E}_{\mathcal{I}^{+},T^{*};\text{\textgreek{d}}}^{(p)}[r^{-j_{2}}\partial_{v}^{j_{1}}\partial_{\text{\textgreek{sv}}}^{j_{2}}\partial_{u}^{j_{3}}\text{\textgreek{f}}](\text{\textgreek{t}}_{1},\text{\textgreek{t}}_{2})\Big\}\lesssim_{p,m,\text{\textgreek{h},\textgreek{d}}} & \sum_{j=0}^{m-1}\sum_{j_{1}+j_{2}+j_{3}=j}\mathcal{E}_{bound,R;\text{\textgreek{d}}}^{(p)}[r^{-j_{2}}\partial_{v}^{j_{1}}\partial_{\text{\textgreek{sv}}}^{j_{2}}\partial_{u}^{j_{3}}\text{\textgreek{f}}](\text{\textgreek{t}}_{1})+\\
 & +\sum_{j=0}^{m-1}\int_{\mathcal{R}(t_{1},t_{2})}|\partial\text{\textgreek{q}}_{R}|\cdot\big(r^{p}|\partial^{j+1}\text{\textgreek{f}}|^{2}+r^{p-2}\cdot|\partial^{j}\text{\textgreek{f}}|^{2}\big)+\\
 & +\sum_{j=0}^{m-1}\int_{\mathcal{R}(t_{1},t_{2})}\text{\textgreek{q}}_{R}\cdot(r^{p+1}+r^{1+\text{\textgreek{h}}})\cdot|\partial^{j}(\square_{g}\text{\textgreek{f}})|^{2}\,\text{\textgreek{W}}^{2}dudvd\text{\textgreek{sv}}.
\end{split}
\label{eq:newMethodHigherDerivativesNotImproved}
\end{equation}

\end{cor}
The proof of Corollary \ref{cor:NewMethodHigherDerivativesNotImproved}
will be omitted. Notice that in comparison to the improved hierarchy
statement (\ref{eq:newMethodDu+DvPhi}), the estimate (\ref{eq:newMethodHigherDerivativesNotImproved})
holds for smaller values of $p$ and the terms of all orders in $\text{\textgreek{f}}$
appear with the same weight.

We can also establish the following variant of Theorem \ref{thm:NewMethodFinalStatementHyperboloids}
in the region bounded by two hypersurfaces of the form $\{t=const\}$:
\begin{thm}
\label{thm:NewMethodFinalStatement}For any $0<p\le2$, any given
$0<\text{\textgreek{h}}<a$ and $0<\text{\textgreek{d}}<1$, any $R>0$
large enough in terms of $p,\text{\textgreek{h}},\text{\textgreek{d}}$
and any $t_{1}\le t_{2}$, the following inequality is true for any
smooth function $\text{\textgreek{f}}:\mathcal{N}_{af}\rightarrow\mathbb{C}$
with compact support in space and any smooth cut-off $\text{\textgreek{q}}_{R}:\mathcal{N}_{af}\rightarrow[0,1]$
supported in $\{r\ge R\}$: 
\begin{equation}
\begin{split}\int_{\{t=t_{2}\}}\text{\textgreek{q}}_{R}\cdot\Big(r^{p}|\partial_{v}(\text{\textgreek{W}}\text{\textgreek{f}})|^{2} & +r^{p}|\frac{1}{r}\partial_{\text{\textgreek{sv}}}(\text{\textgreek{W}}\text{\textgreek{f}})|^{2}+\big((d-3)r^{p-2}+\min\{r^{p-2},r^{-\text{\textgreek{d}}}\}\big)|\text{\textgreek{W}}\text{\textgreek{f}}|^{2}\Big)\, dvd\text{\textgreek{sv}}+\int_{\{t=t_{2}\}}\text{\textgreek{q}}_{R}J_{\text{\textgreek{m}}}^{T}(\text{\textgreek{f}})n^{\text{\textgreek{m}}}+\\
+\int_{\{t_{1}\le t\le t_{2}\}}\text{\textgreek{q}}_{R}\cdot\Big(pr^{p-1}\big|\partial_{v} & (\text{\textgreek{W}}\text{\textgreek{f}})\big|^{2}+\big\{\big((2-p)r^{p-1}+r^{p-1-\text{\textgreek{d}}}\big)\big|r^{-1}\partial_{\text{\textgreek{sv}}}(\text{\textgreek{W}}\text{\textgreek{f}})\big|^{2}+\\
\hphantom{+\int_{\{t_{1}\le t\le t_{2}\}}\text{\textgreek{q}}_{R}\cdot\Big(}+\big((2- & p)(d-3)r^{p-3}+\min\{r^{p-3},r^{-1-\text{\textgreek{d}}}\}\big)\big|\text{\textgreek{W}}\text{\textgreek{f}}\big|^{2}\big\}+r^{-1-\text{\textgreek{h}}}|\partial_{u}(\text{\textgreek{W}}\text{\textgreek{f}})|^{2}\Big)\, dudvd\text{\textgreek{sv}}\lesssim_{p,\text{\textgreek{h}},\text{\textgreek{d}}}\\
\lesssim_{p,\text{\textgreek{h},\textgreek{d}}} & \int_{\{t=t_{1}\}}\text{\textgreek{q}}_{R}\cdot\Big(r^{p}|\partial_{v}(\text{\textgreek{W}}\text{\textgreek{f}})|^{2}+r^{p}|\frac{1}{r}\partial_{\text{\textgreek{sv}}}(\text{\textgreek{W}}\text{\textgreek{f}})|^{2}+\big((d-3)r^{p-2}+\min\{r^{p-2},r^{-\text{\textgreek{d}}}\}\big)|\text{\textgreek{W}}\text{\textgreek{f}}|^{2}\Big)\, dvd\text{\textgreek{sv}}+\\
 & +\int_{\{t=t_{1}\}}\text{\textgreek{q}}_{R}J_{\text{\textgreek{m}}}^{T}(\text{\textgreek{f}})n^{\text{\textgreek{m}}}+\int_{\{t_{1}\le t\le t_{2}\}}|\partial\text{\textgreek{q}}_{R}|\cdot\big(r^{p}|\partial\text{\textgreek{f}}|^{2}+r^{p-2}\cdot|\text{\textgreek{f}}|^{2}\big)+\\
 & +\int_{\{t_{1}\le t\le t_{2}\}}\text{\textgreek{q}}_{R}\cdot(r^{p+1}+r^{1+\text{\textgreek{h}}})\cdot|\square_{g}\text{\textgreek{f}}|^{2}\,\text{\textgreek{W}}^{2}dudvd\text{\textgreek{sv}}.
\end{split}
\label{eq:newMethodFinalStatement}
\end{equation}

In the above, the constants implicit in the $\lesssim_{p,\text{\textgreek{h},\textgreek{d}}}$
notation depend only on $p,\text{\textgreek{h}},\text{\textgreek{d}}$
and on the geometry of $(\mathcal{N}_{af},g)$. 
\end{thm}
The proof of Theorem \ref{thm:NewMethodFinalStatement} is identical
to the proof of Theorem \ref{thm:NewMethodFinalStatementHyperboloids}
(the only difference being the domain over which integrations are
performed), and hence it will be omitted.

\subsection{\label{sub:ProofNewMethod}Proof of Theorem \ref{thm:NewMethodFinalStatementHyperboloids}}

For the proof of Theorem \ref{thm:NewMethodFinalStatementHyperboloids},
we will need to introduce the following energy norms on the hyperboloids
$\{\bar{t}=\text{\textgreek{t}}\}$:

\begin{equation}
\mathcal{E}_{bound,R,T^{*}}^{(p)}[\text{\textgreek{f}}](\text{\textgreek{t}})=\int_{\mathcal{S}_{\text{\textgreek{t}}}\cap\{t\le T^{*}\}}\text{\textgreek{q}}_{R}\Big(r^{p}\big|\partial_{v}(\text{\textgreek{W}}\text{\textgreek{f}})\big|^{2}+r^{-1-\text{\textgreek{h}}'}\big(r^{p}\big|r^{-1}\partial_{\text{\textgreek{sv}}}(\text{\textgreek{W}}\text{\textgreek{f}})\big|^{2}+\frac{(d-1)(d-3)}{4}r^{p-2}\big|\text{\textgreek{W}}\text{\textgreek{f}}\big|^{2}\big)\Big)\, dvd\text{\textgreek{sv}},
\end{equation}
 
\begin{equation}
\mathcal{E}_{bulk,R,T^{*}}^{(p-1)}[\text{\textgreek{f}}](\text{\textgreek{t}})=\int_{\mathcal{S}_{\text{\textgreek{t}}}\cap\{t\le T^{*}\}}\text{\textgreek{q}}_{R}\Big(pr^{p-1}\big|\partial_{v}(\text{\textgreek{W}}\text{\textgreek{f}})\big|^{2}+(2-p)\big(r^{p-1}\big|r^{-1}\partial_{\text{\textgreek{sv}}}(\text{\textgreek{W}}\text{\textgreek{f}})\big|^{2}+\frac{(d-1)(d-3)}{4}r^{p-3}\big|\text{\textgreek{W}}\text{\textgreek{f}}\big|^{2}\big)\Big)\, dvd\text{\textgreek{sv}},
\end{equation}
\begin{align}
\mathcal{E}_{bulk,R,\text{\textgreek{h}},T^{*}}^{(p-1)}[\text{\textgreek{f}}](\text{\textgreek{t}})= & \int_{\mathcal{S}_{\text{\textgreek{t}}}\cap\{t\le T^{*}\}}\text{\textgreek{q}}_{R}\Big(pr^{p-1}\big|\partial_{v}(\text{\textgreek{W}}\text{\textgreek{f}})\big|^{2}+(2-p)\big(r^{p-1}\big|r^{-1}\partial_{\text{\textgreek{sv}}}(\text{\textgreek{W}}\text{\textgreek{f}})\big|^{2}+\frac{(d-1)(d-3)}{4}r^{p-3}\big|\text{\textgreek{W}}\text{\textgreek{f}}\big|^{2}\big)\Big)\, dvd\text{\textgreek{sv}}+\\
 & +\int_{\mathcal{S}_{\text{\textgreek{t}}}\cap\{t\le T^{*}\}}\text{\textgreek{q}}_{R}r^{-1-\text{\textgreek{h}}}|\partial_{u}(\text{\textgreek{W}}\text{\textgreek{f}})|^{2}\, dvd\text{\textgreek{sv}}\nonumber 
\end{align}
and 
\begin{equation}
\mathcal{E}_{\mathcal{I}^{+},R,T^{*}}^{(p)}[\text{\textgreek{f}}](\text{\textgreek{t}}_{1},\text{\textgreek{t}}_{2})=\int_{\mathcal{R}(\text{\textgreek{t}}_{1},\text{\textgreek{t}}_{2})\cap\{t=T^{*}\}}\text{\textgreek{q}}_{R}\Big(r^{p}\big|r^{-1}\partial_{\text{\textgreek{sv}}}(\text{\textgreek{W}}\text{\textgreek{f}})\big|^{2}+\frac{(d-1)(d-3)}{4}r^{p-2}\big|\text{\textgreek{W}}\text{\textgreek{f}}\big|^{2}\Big)\, dvd\text{\textgreek{sv}}.
\end{equation}

The main step in establishing Theorem \ref{thm:NewMethodFinalStatementHyperboloids}
is contained in the following Lemma:
\begin{lem}
\label{lem:NewMethodGeneralCase}For any $0<p\le2$, any given $0<\text{\textgreek{h}}<a$,
any $0<\text{\textgreek{e}}<1$, any $R>0$ large enough in terms
of $p$, $\text{\textgreek{h}}$ and $\text{\textgreek{e}}$, any
$\text{\textgreek{t}}_{1}\le\text{\textgreek{t}}_{2}$ and any $T^{*}\ge0$,
the following inequality holds for any smooth function $\text{\textgreek{f}}:\mathcal{N}_{af}\rightarrow\mathbb{C}$
and any smooth cut-off $\text{\textgreek{q}}_{R}:\mathcal{N}_{af}\rightarrow[0,1]$
supported in $\{r\ge R\}$:
\begin{equation}
\begin{split}\mathcal{E}_{bound,R,T^{*}}^{(p)}[\text{\textgreek{f}}](\text{\textgreek{t}}_{2})+\int_{\text{\textgreek{t}}_{1}}^{\text{\textgreek{t}}_{2}} & \mathcal{E}_{bulk,R,\text{\textgreek{h}},T^{*}}^{(p-1)}[\text{\textgreek{f}}](\text{\textgreek{t}})\, d\text{\textgreek{t}}+\mathcal{E}_{\mathcal{I}^{+},R,T^{*}}^{(p)}[\text{\textgreek{f}}](\text{\textgreek{t}}_{1},\text{\textgreek{t}}_{2})\le\\
\le\big(1+ & O_{p,\text{\textgreek{h}}}(\text{\textgreek{e}})\big)\mathcal{E}_{bound,R,T^{*}}^{(p)}[\text{\textgreek{f}}](\text{\textgreek{t}}_{1})+\\
 & +C_{p,\text{\textgreek{h}},\text{\textgreek{e}}}\int_{\mathcal{R}(\text{\textgreek{t}}_{1},\text{\textgreek{t}}_{2})\cap\{t\le T^{*}\}}|\partial\text{\textgreek{q}}_{R}|\cdot\big(r^{p}|\partial(\text{\textgreek{W}}\text{\textgreek{f}})|^{2}+r^{p-2}|\text{\textgreek{W}}\text{\textgreek{f}}|^{2}\big)\, dudvd\text{\textgreek{sv}}+\\
 & +C_{p,\text{\textgreek{h}},\text{\textgreek{e}}}\int_{\mathcal{S}_{\text{\textgreek{t}}_{1}}\cap\{t\le T^{*}\}}\text{\textgreek{q}}_{R}J_{\text{\textgreek{m}}}^{T}(\text{\textgreek{f}})\bar{n}^{\text{\textgreek{m}}}+C_{p,\text{\textgreek{h}},\text{\textgreek{e}}}\sum_{i=1}^{2}\int_{\mathcal{S}_{\text{\textgreek{t}}_{i}}\cap\{t\le T^{*}\}}\text{\textgreek{q}}_{R}\cdot r^{p-2-a}|\text{\textgreek{W}}\text{\textgreek{f}}|^{2}\, dvd\text{\textgreek{sv}}+\\
 & +C_{p,\text{\textgreek{h}},\text{\textgreek{e}}}\int_{\mathcal{R}(\text{\textgreek{t}}_{1},\text{\textgreek{t}}_{2})\cap\{t=T^{*}\}}\text{\textgreek{q}}_{R}\cdot r^{p-2-a}|\text{\textgreek{W}}\text{\textgreek{f}}|^{2}\, dvd\text{\textgreek{sv}}+\\
 & +C_{p,\text{\textgreek{h}},\text{\textgreek{e}}}\int_{\mathcal{R}(\text{\textgreek{t}}_{1},\text{\textgreek{t}}_{2})\cap\{t\le T^{*}\}}\text{\textgreek{q}}_{R}\cdot\max\{r^{p-3-a},r^{-3}\}|\text{\textgreek{W}}\text{\textgreek{f}}|^{2}\, dudvd\text{\textgreek{sv}}+\\
 & +C_{p,\text{\textgreek{h}},\text{\textgreek{e}}}\int_{\mathcal{R}(\text{\textgreek{t}}_{1},\text{\textgreek{t}}_{2})\cap\{t\le T^{*}\}}\text{\textgreek{q}}_{R}\cdot(r^{p+1}+r^{1+\text{\textgreek{h}}})|\text{\textgreek{W}}\square\text{\textgreek{f}}|^{2}\, dudvd\text{\textgreek{sv}},
\end{split}
\label{eq:NewMethodSimpleHyperboloids}
\end{equation}
 where $\bar{n}$ is the future directed unit normal to the hyperboloids
$\mathcal{S}_{\text{\textgreek{t}}}$.\end{lem}
\begin{proof}
Without loss of generality we will assume that $\text{\textgreek{f}}$
is real valued. We will set 
\begin{equation}
\text{\textgreek{F}}\doteq\text{\textgreek{W}}\cdot\text{\textgreek{f}}.
\end{equation}

Following the approach of \cite{DafRod7}, we start by multiplying
both sides of (\ref{eq:ConformalWaveOperator}) with $\text{\textgreek{q}}_{R}\cdot r^{p}\partial_{v}\text{\textgreek{F}}$.
We therefore compute, after applying the product rule for derivatives:
\begin{align}
\text{\textgreek{q}}_{R}\cdot r^{p}\partial_{v}\text{\textgreek{F}}\cdot\text{\textgreek{W}}\square\text{\textgreek{f}}= & -\frac{1}{2}\partial_{u}\big\{\text{\textgreek{q}}_{R}\cdot r^{p}(1+O_{1}(r^{-1-a}))|\partial_{v}\text{\textgreek{F}}|^{2}\big\}+\frac{1}{2}(\partial_{u}\{\text{\textgreek{q}}_{R}\cdot r^{p}(1+O_{1}(r^{-1-a}))\}|\partial_{v}\text{\textgreek{F}}|^{2}-\label{eq:IntegrationByParts}\\
 & -\frac{1}{2}\partial_{v}\big\{\text{\textgreek{q}}_{R}\cdot\Big(\frac{(d-1)(d-3)}{4}r^{p-2}\Big)\cdot\text{\textgreek{F}}^{2}\big\}+\nonumber \\
 & +\frac{1}{2}\partial_{v}\big\{\text{\textgreek{q}}_{R}\cdot\Big(\frac{(d-1)(d-3)}{4}r^{p-2}\Big)\big\}\text{\textgreek{F}}^{2}+\nonumber \\
 & +\text{\textgreek{q}}_{R}\cdot r^{p-2}\cdot\partial_{v}\text{\textgreek{F}}\cdot\text{\textgreek{D}}_{g_{\mathbb{S}^{d-1}}+h_{\mathbb{S}^{d-1}}}\text{\textgreek{F}}+\text{\textgreek{q}}_{R}\cdot r^{p}\partial_{v}\text{\textgreek{F}}\cdot Err(\text{\textgreek{F}})\nonumber 
\end{align}

Integration of (\ref{eq:IntegrationByParts}) over $\{\text{\textgreek{t}}_{1}\le\bar{t}\le\text{\textgreek{t}}_{2}\}\cap\{t\le T^{*}\}$
using the coordinate volume form $dudvd\text{\textgreek{sv}}$ (in
place of the geometric volume form $\text{\textgreek{W}}^{2}dudvd\text{\textgreek{sv}}$)
readily yields: 
\begin{equation}
\begin{split}\int_{\mathcal{S}_{\text{\textgreek{t}}_{2}}\cap\{t\le T^{*}\}}\frac{1}{2}\text{\textgreek{q}}_{R}\cdot\Big(r^{p}(1+O(r^{-1-a}) & )\cdot|\partial_{v}\text{\textgreek{F}}|^{2}\Big)\, dvd\text{\textgreek{sv}}+\int_{\mathcal{R}(\text{\textgreek{t}}_{1},\text{\textgreek{t}}_{2})\cap\{t=T^{*}\}}\frac{1}{2}\text{\textgreek{q}}_{R}\cdot\Big(r^{p}(1+O(r^{-1-a}))\cdot|\partial_{v}\text{\textgreek{F}}|^{2}\Big)\, dvd\text{\textgreek{sv}}+\\
+\int_{\mathcal{S}_{\text{\textgreek{t}}_{2}}\cap\{t\le T^{*}\}}\frac{1}{2}\text{\textgreek{q}}_{R}\cdot\Big(\frac{(d-1)(d-3)}{4} & r^{p-2}+O(r^{p-3})\Big)\text{\textgreek{F}}^{2}\, dud\text{\textgreek{sv}}+\int_{\mathcal{R}(\text{\textgreek{t}}_{1},\text{\textgreek{t}}_{2})\cap\{t=T^{*}\}}\frac{1}{2}\text{\textgreek{q}}_{R}\cdot\Big(\frac{(d-1)(d-3)}{4}r^{p-2}+O(r^{p-3})\Big)\text{\textgreek{F}}^{2}\, dud\text{\textgreek{sv}}+\\
+\int_{\mathcal{R}(t_{1},t_{2})\cap\{t\le T^{*}\}}\text{\textgreek{q}}_{R}\cdot\Big\{\frac{1}{2}pr^{p-1}\{1 & +O(r^{-1-a})\}|\partial_{v}\text{\textgreek{F}}|^{2}+\frac{1}{2}\text{\textgreek{q}}_{R}\cdot\Big((2-p)\frac{(d-1)(d-3)}{4}r^{p-3}+O(r^{p-4})\Big)\text{\textgreek{F}}^{2}\Big\}\, dudvd\text{\textgreek{sv}}+\\
+\int_{\mathcal{R}(t_{1},t_{2})\cap\{t\le T^{*}\}}\text{\textgreek{q}}_{R}\cdot r^{p-2}\cdot\partial_{v}\text{\textgreek{F}}\cdot & \text{\textgreek{D}}_{g_{\mathbb{S}^{d-1}}+h_{\mathbb{S}^{d-1}}}\text{\textgreek{F}}\, dudvd\text{\textgreek{sv}}+\int_{\mathcal{R}(t_{1},t_{2})\cap\{t\le T^{*}\}}\text{\textgreek{q}}_{R}\cdot r^{p}\partial_{v}\text{\textgreek{F}}\cdot Err(\text{\textgreek{F}})\, dudvd\text{\textgreek{sv}}=\\
= & \int_{\mathcal{S}_{\text{\textgreek{t}}_{1}}\cap\{t\le T^{*}\}}\frac{1}{2}\text{\textgreek{q}}_{R}\cdot\Big(r^{p}(1+O(r^{-1-a}))\cdot|\partial_{v}\text{\textgreek{F}}|^{2}\Big)\, dvd\text{\textgreek{sv}}+\\
 & +\int_{\mathcal{S}_{\text{\textgreek{t}}_{1}}\cap\{t\le T^{*}\}}\frac{1}{2}\text{\textgreek{q}}_{R}\cdot\Big(\frac{(d-1)(d-3)}{4}r^{p-2}+O(r^{p-3})\Big)\text{\textgreek{F}}^{2}\, dud\text{\textgreek{sv}}+\\
 & +\int_{\mathcal{R}(t_{1},t_{2})\cap\{t\le T^{*}\}}O(|\partial\text{\textgreek{q}}_{R}|)\cdot\big(O(r^{p})\cdot|\partial_{v}\text{\textgreek{F}}|^{2}+O(r^{p-2})\cdot\text{\textgreek{F}}^{2}\big)\, dudvd\text{\textgreek{sv}}-\\
 & -\int_{\mathcal{R}(t_{1},t_{2})\cap\{t\le T^{*}\}}\text{\textgreek{q}}_{R}\cdot r^{p}\partial_{v}\text{\textgreek{F}}\cdot\text{\textgreek{W}}\square\text{\textgreek{f}}\, dudvd\text{\textgreek{sv}}.
\end{split}
\label{eq:BeforeSphericalIntegration}
\end{equation}

Moreover, we obtain after integrating by parts in the spherical directions
(recalling that $h_{\mathbb{S}^{d-1}}=O(r^{-1})$):

\begin{multline}
\int_{\mathcal{R}(t_{1},t_{2})\cap\{t\le T^{*}\}}\text{\textgreek{q}}_{R}\cdot r^{p-2}\cdot\partial_{v}\text{\textgreek{F}}\cdot\text{\textgreek{D}}_{g_{\mathbb{S}^{d-1}}+h_{\mathbb{S}^{d-1}}}\text{\textgreek{F}}\, dudvd\text{\textgreek{sv}}=\int_{\mathcal{R}(t_{1},t_{2})\cap\{t\le T^{*}\}}\text{\textgreek{q}}_{R}\cdot r^{p-2}\cdot\frac{1}{2}\partial_{v}\big\{|\partial_{\text{\textgreek{sv}}}\text{\textgreek{F}}|^{2}\big\}\, dudvd\text{\textgreek{sv}}+\\
+\int_{\mathcal{R}(t_{1},t_{2})\cap\{t\le T^{*}\}}\text{\textgreek{q}}_{R}\cdot O(r^{p-3})\cdot\partial_{v}\text{\textgreek{F}}\cdot\partial_{\text{\textgreek{sv}}}\text{\textgreek{F}}\, dudvd\text{\textgreek{sv}}+\int_{\mathcal{R}(t_{1},t_{2})\cap\{t\le T^{*}\}}O(|\partial_{\text{\textgreek{sv}}}\text{\textgreek{q}}_{R}|)\cdot r^{p-2}\cdot\partial_{v}\text{\textgreek{F}}\cdot\partial_{\text{\textgreek{sv}}}\text{\textgreek{F}}\, dudvd\text{\textgreek{sv}}.\label{eq:SphericalIntegration}
\end{multline}
 Hence, integrating by parts in the first term of the right hand side
of (\ref{eq:SphericalIntegration}), substituting in (\ref{eq:BeforeSphericalIntegration})
and absorbing the $\int\text{\textgreek{q}}_{R}\cdot O(r^{p-3})\cdot\partial_{v}\text{\textgreek{F}}\cdot\partial_{\text{\textgreek{sv}}}\text{\textgreek{F}}$
summand of the right hand side of (\ref{eq:SphericalIntegration})
into the $Err$ term, we infer (provided $R$ is large enough in terms
of $\text{\textgreek{e}}$): 
\begin{equation}
\begin{split}\mathcal{E}_{bound,R,T^{*}}^{(p)}[\text{\textgreek{f}}](\text{\textgreek{t}}_{2})+ & \int_{\text{\textgreek{t}}_{1}}^{\text{\textgreek{t}}_{2}}\mathcal{E}_{bulk,R,T^{*}}^{(p-1)}[\text{\textgreek{f}}](\text{\textgreek{t}})\, d\text{\textgreek{t}}+\mathcal{E}_{\mathcal{I}^{+},R,T^{*}}^{(p)}[\text{\textgreek{f}}](\text{\textgreek{t}}_{1},\text{\textgreek{t}}_{2})\le\\
\le & (1+\text{\textgreek{e}})\mathcal{E}_{bound,R,T^{*}}^{(p)}[\text{\textgreek{f}}](\text{\textgreek{t}}_{1})+\int_{\mathcal{R}(\text{\textgreek{t}}_{1},\text{\textgreek{t}}_{2})\cap\{t\le T^{*}\}}|\partial\text{\textgreek{q}}_{R}|\cdot\big(O(r^{p})|\partial(\text{\textgreek{W}}\text{\textgreek{f}})|^{2}+O(r^{p-2})|\text{\textgreek{W}}\text{\textgreek{f}}|^{2}\big)\, dudvd\text{\textgreek{sv}}+\\
 & +C\cdot\sum_{i=1}^{2}\int_{\mathcal{S}_{\text{\textgreek{t}}_{i}}\cap\{t\le T^{*}\}}\text{\textgreek{q}}_{R}J_{\text{\textgreek{m}}}^{T}(\text{\textgreek{f}})\bar{n}^{\text{\textgreek{m}}}+C\cdot\int_{\mathcal{R}(\text{\textgreek{t}}_{1},\text{\textgreek{t}}_{2})\cap\{t=T^{*}\}}\text{\textgreek{q}}_{R}J_{\text{\textgreek{m}}}^{T}(\text{\textgreek{f}})\bar{n}^{\text{\textgreek{m}}}+\\
 & +\sum_{i=1}^{2}\int_{\mathcal{S}_{\text{\textgreek{t}}_{i}}\cap\{t\le T^{*}\}}\text{\textgreek{q}}_{R}\cdot O(r^{p-2-a})|\text{\textgreek{W}}\text{\textgreek{f}}|^{2}\, dvd\text{\textgreek{sv}}+\int_{\mathcal{R}(\text{\textgreek{t}}_{1},\text{\textgreek{t}}_{2})\cap\{t=T^{*}\}}\text{\textgreek{q}}_{R}\cdot O(r^{p-2-a})|\text{\textgreek{W}}\text{\textgreek{f}}|^{2}\, dvd\text{\textgreek{sv}}+\\
 & +\int_{\mathcal{R}(\text{\textgreek{t}}_{1},\text{\textgreek{t}}_{2})\cap\{t\le T^{*}\}}\text{\textgreek{q}}_{R}\cdot O(r^{p-3-a})|\text{\textgreek{W}}\text{\textgreek{f}}|^{2}\, dudvd\text{\textgreek{sv}}+\\
 & +\int_{\mathcal{R}(t_{1},t_{2})\cap\{t\le T^{*}\}}\text{\textgreek{q}}_{R}\cdot O(r^{p})\partial_{v}\text{\textgreek{F}}\cdot Err(\text{\textgreek{F}})\, dudvd\text{\textgreek{sv}}+\\
 & +\int_{\mathcal{R}(\text{\textgreek{t}}_{1},\text{\textgreek{t}}_{2})\cap\{t\le T^{*}\}}\text{\textgreek{q}}_{R}\cdot O(r^{p})\partial_{v}\text{\textgreek{F}}\cdot\text{\textgreek{W}}\square\text{\textgreek{f}}\, dudvd\text{\textgreek{sv}},
\end{split}
\label{eq:newMethod1-1}
\end{equation}
 where we denote for simplicity $|\partial\text{\textgreek{F}}|^{2}=|\partial_{v}\text{\textgreek{F}}|^{2}+|\partial_{u}\text{\textgreek{F}}|^{2}+|\frac{1}{r}\partial_{\text{\textgreek{sv}}}\text{\textgreek{F}}|^{2}$.
Notice that we have used the asymptotic relation $dud\text{\textgreek{sv}}\sim r^{-1-\text{\textgreek{h}}'}dvd\text{\textgreek{sv}}$
on the hyperboloids $\mathcal{S}_{\text{\textgreek{t}}}$.

The positive $\int_{\mathcal{R}(t_{1},t_{2})\cap\{t\le T^{*}\}}(2-p)\text{\textgreek{q}}_{R}r^{p-1}|\frac{1}{r}\partial_{\text{\textgreek{sv}}}\text{\textgreek{F}}|^{2}$
term of the left hand side of (\ref{eq:newMethod1-1}) vanishes for
$p=2$. Since it will prove useful to retain some control over angular
derivatives in the left hand side of (\ref{eq:newMethod1-1}) even
in the $p=2$ case, we will add to (\ref{eq:newMethod1-1}) the same
inequality but for $p'=p-a$ in place of $p$ in case $p\ge1$, obtaining:
\begin{equation}
\begin{split}\mathcal{E}_{bound,R,T^{*}}^{(p)}[\text{\textgreek{f}}](\text{\textgreek{t}}_{2}) & +\int_{\text{\textgreek{t}}_{1}}^{\text{\textgreek{t}}_{2}}\Big(\mathcal{E}_{bulk,R,T^{*}}^{(p-1)}[\text{\textgreek{f}}](\text{\textgreek{t}})+\mathcal{E}_{bulk,R,T^{*}}^{(p-1-a)}[\text{\textgreek{f}}](\text{\textgreek{t}})\Big)\, d\text{\textgreek{t}}+\mathcal{E}_{\mathcal{I}^{+},R,T^{*}}^{(p)}[\text{\textgreek{f}}](\text{\textgreek{t}}_{1},\text{\textgreek{t}}_{2})\le\\
\le & (1+\text{\textgreek{e}})\mathcal{E}_{bound,R,T^{*}}^{(p)}[\text{\textgreek{f}}](\text{\textgreek{t}}_{1})+\int_{\mathcal{R}(\text{\textgreek{t}}_{1},\text{\textgreek{t}}_{2})\cap\{t\le T^{*}\}}|\partial\text{\textgreek{q}}_{R}|\cdot\big(O(r^{p})|\partial(\text{\textgreek{W}}\text{\textgreek{f}})|^{2}+O(r^{p-2})|\text{\textgreek{W}}\text{\textgreek{f}}|^{2}\big)\, dudvd\text{\textgreek{sv}}+\\
 & +C\sum_{i=1}^{2}\int_{\mathcal{S}_{\text{\textgreek{t}}_{i}}\cap\{t\le T^{*}\}}\text{\textgreek{q}}_{R}J_{\text{\textgreek{m}}}^{T}(\text{\textgreek{f}})\bar{n}^{\text{\textgreek{m}}}+C\int_{\mathcal{R}(\text{\textgreek{t}}_{1},\text{\textgreek{t}}_{2})\cap\{t=T^{*}\}}\text{\textgreek{q}}_{R}J_{\text{\textgreek{m}}}^{T}(\text{\textgreek{f}})\bar{n}^{\text{\textgreek{m}}}+\\
 & +\sum_{i=1}^{2}\int_{\mathcal{S}_{\text{\textgreek{t}}_{i}}\cap\{t\le T^{*}\}}\text{\textgreek{q}}_{R}\cdot O(r^{p-2-a})|\text{\textgreek{W}}\text{\textgreek{f}}|^{2}\, dvd\text{\textgreek{sv}}+\int_{\mathcal{R}(\text{\textgreek{t}}_{1},\text{\textgreek{t}}_{2})\cap\{t=T^{*}\}}\text{\textgreek{q}}_{R}\cdot O(r^{p-2-a})|\text{\textgreek{W}}\text{\textgreek{f}}|^{2}\, dvd\text{\textgreek{sv}}+\\
 & +\int_{\mathcal{R}(\text{\textgreek{t}}_{1},\text{\textgreek{t}}_{2})\cap\{t\le T^{*}\}}\text{\textgreek{q}}_{R}\cdot O(r^{p-3-a})|\text{\textgreek{W}}\text{\textgreek{f}}|^{2}\, dudvd\text{\textgreek{sv}}+\\
 & +\int_{\mathcal{R}(t_{1},t_{2})\cap\{t\le T^{*}\}}\text{\textgreek{q}}_{R}\cdot O(r^{p})\partial_{v}\text{\textgreek{F}}\cdot Err(\text{\textgreek{F}})\, dudvd\text{\textgreek{sv}}+\\
 & +\int_{\mathcal{R}(\text{\textgreek{t}}_{1},\text{\textgreek{t}}_{2})\cap\{t\le T^{*}\}}\text{\textgreek{q}}_{R}\cdot O(r^{p})\partial_{v}\text{\textgreek{F}}\cdot\text{\textgreek{W}}\square\text{\textgreek{f}}\, dudvd\text{\textgreek{sv}}.
\end{split}
\label{eq:newMethod1}
\end{equation}

In order to reach (\ref{eq:NewMethodSimpleHyperboloids}) it only
remains to add to (\ref{eq:newMethod1}) the estimate of Lemma (\ref{lem:MorawetzDrLemmaGeneralCase})
and suitably absorb the $Err(\text{\textgreek{F}})$ term in (\ref{eq:newMethod1})
into the left hand side. Moreover, the energy boundary terms on $\mathcal{S}_{\text{\textgreek{t}}_{2}}$
and $\{t=T^{*}\}$ will be dealt with by using Lemma \ref{lem:Boundedness Hyperboloids}.
In particular, assuming without loss of generality that $\text{\textgreek{e}}\ll_{p,\text{\textgreek{h}}}1$,
if $R$ is large enough in terms of $\text{\textgreek{e}},p,\text{\textgreek{h}}$
the following estimate holds: 
\begin{equation}
\int_{\mathcal{R}(t_{1},t_{2})\cap\{t\le T^{*}\}}\text{\textgreek{q}}_{R}\cdot O(r^{p})\partial_{v}\text{\textgreek{F}}\cdot Err(\text{\textgreek{F}})\, dudvd\text{\textgreek{sv}}\lesssim_{p,\text{\textgreek{h}}}\mathcal{B}_{Err,\text{\textgreek{h}},\text{\textgreek{e}}}^{(p)}[\text{\textgreek{f}}](\text{\textgreek{t}}_{1},\text{\textgreek{t}}_{2}),\label{eq:BoundForTheErrTerms}
\end{equation}
 where 
\begin{align}
\mathcal{B}_{Err,\text{\textgreek{h}},\text{\textgreek{e}}}^{(p)}[\text{\textgreek{f}}](\text{\textgreek{t}}_{1},\text{\textgreek{t}}_{2})\doteq & \text{\textgreek{e}}\Big\{\int_{\text{\textgreek{t}}_{1}}^{\text{\textgreek{t}}_{2}}\Big(\mathcal{E}_{bulk,R,T^{*}}^{(p-1)}[\text{\textgreek{f}}](\text{\textgreek{t}})+\mathcal{E}_{bulk,R,T^{*}}^{(p-1-a)}[\text{\textgreek{f}}](\text{\textgreek{t}})\Big)\, d\text{\textgreek{t}}+\sum_{i=1}^{2}\mathcal{E}_{bound,R,T^{*}}^{(p)}[\text{\textgreek{f}}](\text{\textgreek{t}}_{i})+\mathcal{E}_{\mathcal{I}^{+},R,T^{*}}^{(p)}[\text{\textgreek{f}}](\text{\textgreek{t}}_{1},\text{\textgreek{t}}_{2})\Big\}+\label{eq:ShorthandBoundErrTerms}\\
 & +\int_{\mathcal{R}(\text{\textgreek{t}}_{1},\text{\textgreek{t}}_{2})\cap\{t\le T^{*}\}}|\partial\text{\textgreek{q}}_{R}|\cdot\big(r^{p}|\partial(\text{\textgreek{W}}\text{\textgreek{f}})|^{2}+r^{p-2}|\text{\textgreek{W}}\text{\textgreek{f}}|^{2}\big)\, dudvd\text{\textgreek{sv}}+\nonumber \\
 & +\sum_{i=1}^{2}\int_{\mathcal{S}_{\text{\textgreek{t}}_{i}}\cap\{t\le T^{*}\}}\text{\textgreek{q}}_{R}J_{\text{\textgreek{m}}}^{T}(\text{\textgreek{f}})\bar{n}^{\text{\textgreek{m}}}+\int_{\mathcal{R}(\text{\textgreek{t}}_{1},\text{\textgreek{t}}_{2})\cap\{t=T^{*}\}}\text{\textgreek{q}}_{R}J_{\text{\textgreek{m}}}^{T}(\text{\textgreek{f}})\bar{n}^{\text{\textgreek{m}}}+\nonumber \\
 & +\sum_{i=1}^{2}\int_{\mathcal{S}_{\text{\textgreek{t}}_{i}}\cap\{t\le T^{*}\}}\text{\textgreek{q}}_{R}\cdot r^{p-2-a}|\text{\textgreek{W}}\text{\textgreek{f}}|^{2}\, dvd\text{\textgreek{sv}}+\int_{\mathcal{R}(\text{\textgreek{t}}_{1},\text{\textgreek{t}}_{2})\cap\{t=T^{*}\}}\text{\textgreek{q}}_{R}\cdot r^{p-2-a}|\text{\textgreek{W}}\text{\textgreek{f}}|^{2}\, dvd\text{\textgreek{sv}}+\nonumber \\
 & +\int_{\mathcal{R}(\text{\textgreek{t}}_{1},\text{\textgreek{t}}_{2})\cap\{t\le T^{*}\}}\text{\textgreek{q}}_{R}\cdot\max\{r^{p-3-a},r^{-3}\}|\text{\textgreek{W}}\text{\textgreek{f}}|^{2}\, dudvd\text{\textgreek{sv}}+\nonumber \\
 & +\int_{\mathcal{R}(\text{\textgreek{t}}_{1},\text{\textgreek{t}}_{2})\cap\{t\le T^{*}\}}\text{\textgreek{q}}_{R}\cdot(r^{p+1}+r^{1+\text{\textgreek{h}}})|\text{\textgreek{W}}\square\text{\textgreek{f}}|^{2}\, dudvd\text{\textgreek{sv}}.\nonumber 
\end{align}
Notice that in view of Lemma \ref{lem:Boundedness Hyperboloids} and
the fact that $0<p\le2$, the following bound also holds: 
\begin{align}
\int_{\mathcal{S}_{\text{\textgreek{t}}_{2}}\cap\{t\le T^{*}\}}\text{\textgreek{q}}_{R}J_{\text{\textgreek{m}}}^{T}(\text{\textgreek{f}})\bar{n}^{\text{\textgreek{m}}}+ & \int_{\mathcal{R}(\text{\textgreek{t}}_{1},\text{\textgreek{t}}_{2})\cap\{t=T^{*}\}}\text{\textgreek{q}}_{R}J_{\text{\textgreek{m}}}^{T}(\text{\textgreek{f}})\bar{n}^{\text{\textgreek{m}}}\lesssim_{\text{\textgreek{h}}}\int_{\mathcal{S}_{\text{\textgreek{t}}_{1}}\cap\{t\le T^{*}\}}\text{\textgreek{q}}_{R}J_{\text{\textgreek{m}}}^{T}(\text{\textgreek{f}})\bar{n}^{\text{\textgreek{m}}}+\label{eq:FromBoundednessNewMethod}\\
 & +\text{\textgreek{e}}\int_{\text{\textgreek{t}}_{1}}^{\text{\textgreek{t}}_{2}}\Big(\mathcal{E}_{bulk,R,T^{*}}^{(p-1)}[\text{\textgreek{f}}](\text{\textgreek{t}})+\mathcal{E}_{bulk,R,T^{*}}^{(p-1-a)}[\text{\textgreek{f}}](\text{\textgreek{t}})\Big)\, d\text{\textgreek{t}}+\int_{\mathcal{R}(\text{\textgreek{t}}_{1},\text{\textgreek{t}}_{2})\cap\{t\le T^{*}\}}\text{\textgreek{q}}_{R}\cdot r^{-3}|\text{\textgreek{W}}\text{\textgreek{f}}|^{2}\, dudvd\text{\textgreek{sv}}+\nonumber \\
 & +\int_{\mathcal{R}(\text{\textgreek{t}}_{1},\text{\textgreek{t}}_{2})\cap\{t\le T^{*}\}}|\partial\text{\textgreek{q}}_{R}|\cdot\big(r^{p}|\partial(\text{\textgreek{W}}\text{\textgreek{f}})|^{2}+r^{p-2}|\text{\textgreek{W}}\text{\textgreek{f}}|^{2}\big)\, dudvd\text{\textgreek{sv}}+\nonumber \\
 & +\int_{\mathcal{R}(\text{\textgreek{t}}_{1},\text{\textgreek{t}}_{2})\cap\{t\le T^{*}\}}\text{\textgreek{q}}_{R}\cdot(r^{p+1}+r^{1+\text{\textgreek{h}}})|\text{\textgreek{W}}\square\text{\textgreek{f}}|^{2}\, dudvd\text{\textgreek{sv}}.\nonumber 
\end{align}
 Using (\ref{eq:newMethod1}), (\ref{eq:BoundFromMorawetzDuOnly}),
(\ref{eq:ShorthandBoundErrTerms}), (\ref{eq:BoundForTheErrTerms})
and (\ref{eq:FromBoundednessNewMethod}), inequality (\ref{eq:NewMethodSimpleHyperboloids})
follows readily, provided $\text{\textgreek{e}}$ has been chosen
small enough in terms of $p$ and $\text{\textgreek{h}}$ so that
the first line of the right hand side of (\ref{eq:ShorthandBoundErrTerms})
can be absorbed into the left hand side of (\ref{eq:newMethod1}).

Notice that 
\begin{equation}
\begin{split}\int_{\mathcal{R}(t_{1},t_{2})\cap\{t\le T^{*}\}}\text{\textgreek{q}}_{R}\cdot O(r^{p})\partial_{v}\text{\textgreek{F}}\cdot & Err(\text{\textgreek{F}})\, dudvd\text{\textgreek{sv}}=\\
= & \int_{\mathcal{R}(t_{1},t_{2})\cap\{t\le T^{*}\}}\text{\textgreek{q}}_{R}\partial_{v}\text{\textgreek{F}}\cdot\Big\{ O(r^{p-1})\partial_{v}^{2}\text{\textgreek{F}}+O(r^{p-2})\partial_{v}\partial_{\text{\textgreek{sv}}}\text{\textgreek{F}}+O(r^{p-2-a})\partial_{u}^{2}\text{\textgreek{F}}+\\
 & \hphantom{\int_{\mathcal{R}(t_{1},t_{2})\cap\{t\le T^{*}\}}\text{\textgreek{q}}_{R}\partial_{v}\text{\textgreek{F}}\cdot\Big\{}+O(r^{p-2-a})\partial_{u}\partial_{\text{\textgreek{sv}}}\text{\textgreek{F}}+O(r^{p-3-a})\partial_{\text{\textgreek{sv}}}\partial_{\text{\textgreek{sv}}}\text{\textgreek{F}}+O(r^{p-2-a})\partial_{\text{\textgreek{sv}}}\text{\textgreek{F}}+\\
 & \hphantom{\int_{\mathcal{R}(t_{1},t_{2})\cap\{t\le T^{*}\}}\text{\textgreek{q}}_{R}\partial_{v}\text{\textgreek{F}}\cdot\Big\{}+O(r^{p-1-a})\partial_{v}\text{\textgreek{F}}+O(r^{p-2-a})\partial_{u}\text{\textgreek{F}}+O(r^{p-3})\text{\textgreek{F}}\Big\}\, dudvd\text{\textgreek{sv}}.
\end{split}
\label{eq:ErrTerms-1}
\end{equation}
 In view of Lemma \ref{lem:MorawetzDrLemmaHyperboloids} and the fact
that $0<p\le2$, we can bound for any $\text{\textgreek{e}}>0$ provided
$R$ is large enough in terms of $\text{\textgreek{e}},\text{\textgreek{h}},p$:
\begin{align}
\int_{\mathcal{R}(\text{\textgreek{t}}_{1},\text{\textgreek{t}}_{2})\cap\{t\le T^{*}\}}\text{\textgreek{q}}_{R}\cdot r^{p-3-a}|\partial_{u}\text{\textgreek{F}}|^{2}\, dudvd\text{\textgreek{sv}} & \le\int_{\mathcal{R}(\text{\textgreek{t}}_{1},\text{\textgreek{t}}_{2})\cap\{t\le T^{*}\}}\text{\textgreek{q}}_{R}\cdot r^{-1-\text{\textgreek{h}}}|\partial_{u}\text{\textgreek{F}}|^{2}\, dudvd\text{\textgreek{sv}}\lesssim_{p,\text{\textgreek{h}}}\label{eq:BoundFromMorawetzDuOnly}\\
 & \lesssim_{p,\text{\textgreek{h}}}\mathcal{B}_{Err,\text{\textgreek{h}},\text{\textgreek{e}}}^{(p)}[\text{\textgreek{f}}](\text{\textgreek{t}}_{1},\text{\textgreek{t}}_{2}).\nonumber 
\end{align}
 We can also trivialy bound (if $R$ is large enough in terms of $\text{\textgreek{e}},\text{\textgreek{h}},p$):
\begin{equation}
\int_{\mathcal{R}(\text{\textgreek{t}}_{1},\text{\textgreek{t}}_{2})\cap\{t\le T^{*}\}}\text{\textgreek{q}}_{R}\cdot r^{p-1-a}\Big(|\partial_{v}\text{\textgreek{F}}|^{2}+|r^{-1}\partial_{\text{\textgreek{sv}}}\text{\textgreek{F}}|^{2}+r^{-2}|\text{\textgreek{F}}|^{2}\Big)\, dudvd\text{\textgreek{sv}}\lesssim_{p,\text{\textgreek{h}}}\mathcal{B}_{Err,\text{\textgreek{h}},\text{\textgreek{e}}}^{(p)}[\text{\textgreek{f}}](\text{\textgreek{t}}_{1},\text{\textgreek{t}}_{2}).\label{eq:BoundFromLeftHandSide}
\end{equation}
Therefore, applying a Cauchy--Schwarz ineqality, all the terms of
(\ref{eq:ErrTerms-1}) which do not contain second order derivatives
of $\text{\textgreek{F}}$ can be bounded by $\mathcal{B}_{Err,\text{\textgreek{h}},\text{\textgreek{e}}}^{(p)}[\text{\textgreek{f}}](\text{\textgreek{t}}_{1},\text{\textgreek{t}}_{2})$.

It thus remains to estimate the terms of $\int_{\mathcal{R}(t_{1},t_{2})}\text{\textgreek{q}}_{R}\cdot r^{p}\cdot\partial_{v}\text{\textgreek{F}}\cdot Err(\text{\textgreek{F}})\, dudvd\text{\textgreek{sv}}$
which contain second order derivatives of $\text{\textgreek{F}}$,
i.\,e.~we have to bound (provided again that $R$ is large enough
in terms of $p,\text{\textgreek{h}},\text{\textgreek{e}}$): 
\begin{align}
\int_{\mathcal{R}(\text{\textgreek{t}}_{1},\text{\textgreek{t}}_{2})\cap\{t\le T^{*}\}}\text{\textgreek{q}}_{R}\cdot\Big\{ & O(r^{p-1})\partial_{v}\text{\textgreek{F}}\cdot\partial_{v}^{2}\text{\textgreek{F}}+O(r^{p-2-a})\partial_{v}\text{\textgreek{F}}\cdot\partial_{u}^{2}\text{\textgreek{F}}+\label{eq:BoundSecondOrderErrTerms}\\
 & +O(r^{p-3-a})\partial_{v}\text{\textgreek{F}}\cdot\partial_{\text{\textgreek{sv}}}\partial_{\text{\textgreek{sv}}}\text{\textgreek{F}}+O(r^{p-2-a})\partial_{v}\text{\textgreek{F}}\cdot\partial_{u}\partial_{\text{\textgreek{sv}}}\text{\textgreek{F}}+O(r^{p-2})\partial_{v}\text{\textgreek{F}}\cdot\partial_{v}\partial_{\text{\textgreek{sv}}}\text{\textgreek{F}}\Big\}\, dudvd\text{\textgreek{sv}}\lesssim_{p,\text{\textgreek{h}}}\mathcal{B}_{Err,\text{\textgreek{h}},\text{\textgreek{e}}}^{(p)}[\text{\textgreek{f}}](\text{\textgreek{t}}_{1},\text{\textgreek{t}}_{2}).\nonumber 
\end{align}

\noindent In order to estimate these terms, we will perform some integrations
by parts over $\mathcal{R}(\text{\textgreek{t}}_{1},\text{\textgreek{t}}_{2})\cap\{t\le T^{*}\}$
which will reduce the maximum number of derivatives of $\text{\textgreek{F}}$
appearing in these expressions to just one. 

In particular, we proceed as follows:

\smallskip{}

\noindent 1. For the part of $\int_{\mathcal{R}(t_{1},t_{2})}\text{\textgreek{q}}_{R}\cdot r^{p}\cdot\partial_{v}\text{\textgreek{F}}\cdot Err(\text{\textgreek{F}})$
consisiting of terms of the form $\partial_{v}\text{\textgreek{F}}\cdot\partial_{i}\partial_{v}\text{\textgreek{F}}$,
namely: 
\begin{equation}
Err_{v,i,v}^{(p)}[\text{\textgreek{f}}](\text{\textgreek{t}}_{1},\text{\textgreek{t}}_{2})\doteq\int_{\mathcal{R}(\text{\textgreek{t}}_{1},\text{\textgreek{t}}_{2})\cap\{t\le T^{*}\}}\text{\textgreek{q}}_{R}\cdot\Big\{ O(r^{p-1})\partial_{v}\text{\textgreek{F}}\cdot\partial_{v}^{2}\text{\textgreek{F}}+O(r^{p-2})\cdot\partial_{v}\text{\textgreek{F}}\cdot\partial_{v}\partial_{\text{\textgreek{sv}}}\text{\textgreek{F}}\Big\}\, dudvd\text{\textgreek{sv}},
\end{equation}

\noindent after putting each summand in the form $\frac{1}{2}\partial_{i}\{(\partial_{v}\text{\textgreek{F}})^{2}\}$
and integrating by parts in $\partial_{i}$ we obtain 
\begin{align}
Err_{v,i,v}^{(p)}[\text{\textgreek{f}}](\text{\textgreek{t}}_{1},\text{\textgreek{t}}_{2})\lesssim & \int_{\mathcal{R}(\text{\textgreek{t}}_{1},\text{\textgreek{t}}_{2})\cap\{t\le T^{*}\}}\text{\textgreek{q}}_{R}\cdot r^{p-2}\cdot|\partial_{v}\text{\textgreek{F}}|^{2}\, dudvd\text{\textgreek{sv}}+\int_{\mathcal{R}(\text{\textgreek{t}}_{1},\text{\textgreek{t}}_{2})\cap\{t\le T^{*}\}}|\partial\text{\textgreek{q}}_{R}|\cdot r^{p-1}\cdot|\partial_{v}\text{\textgreek{F}}|^{2}\, dudvd\text{\textgreek{sv}}+\label{eq:BoundIVErrTerms}\\
 & +\sum_{i=1}^{2}\int_{\mathcal{S}_{\text{\textgreek{t}}_{i}}\cap\{t\le T^{*}\}}\text{\textgreek{q}}_{R}\cdot r^{p-2}|\partial_{v}\text{\textgreek{F}}|^{2}\, dvd\text{\textgreek{sv}}+\int_{\mathcal{R}(\text{\textgreek{t}}_{1},\text{\textgreek{t}}_{2})\cap\{t=T^{*}\}}\text{\textgreek{q}}_{R}\cdot r^{p-2}|\partial_{v}\text{\textgreek{F}}|^{2}\, dvd\text{\textgreek{sv}}.\nonumber 
\end{align}
 Thus, since $p\le2$, we immediately infer (provided again that $R$
is large enough in terms of $p,\text{\textgreek{h}},\text{\textgreek{e}}$)
that 
\begin{equation}
Err_{v,i,v}^{(p)}[\text{\textgreek{f}}](\text{\textgreek{t}}_{1},\text{\textgreek{t}}_{2})\lesssim_{p,\text{\textgreek{h}}}\mathcal{B}_{Err,\text{\textgreek{h}},\text{\textgreek{e}}}^{(p)}[\text{\textgreek{f}}](\text{\textgreek{t}}_{1},\text{\textgreek{t}}_{2}).\label{eq:FinalBoundIVErrTerms}
\end{equation}

\smallskip{}

\noindent 2. For the part of $\int_{\mathcal{R}(t_{1},t_{2})}\text{\textgreek{q}}_{R}\cdot r^{p}\cdot\partial_{v}\text{\textgreek{F}}\cdot Err(\text{\textgreek{F}})$
consisiting of terms of the form $\partial_{v}\text{\textgreek{F}}\partial_{i}^{2}\text{\textgreek{F}}$,
namely 
\begin{equation}
Err_{v,i,i}^{(p)}[\text{\textgreek{f}}](\text{\textgreek{t}}_{1},\text{\textgreek{t}}_{2})\doteq\int_{\mathcal{R}(\text{\textgreek{t}}_{1},\text{\textgreek{t}}_{2})\cap\{t\le T^{*}\}}\text{\textgreek{q}}_{R}\cdot\Big\{ O(r^{p-2-a})\cdot\partial_{v}\text{\textgreek{F}}\cdot\partial_{u}^{2}\text{\textgreek{F}}+O(r^{p-3-a})\cdot\partial_{v}\text{\textgreek{F}}\cdot\partial_{\text{\textgreek{sv}}}\partial_{\text{\textgreek{sv}}}\text{\textgreek{F}}\Big\}\, dudvd\text{\textgreek{sv}},
\end{equation}

\noindent we will perform two integrations by parts: By first integrating
by parts schematically as $\partial_{v}\text{\textgreek{F}}\partial_{i}^{2}\text{\textgreek{F}}\rightarrow\partial_{i}\partial_{v}\text{\textgreek{F}}\partial_{i}\text{\textgreek{F}}=\frac{1}{2}\partial_{v}\{\partial_{i}\text{\textgreek{F}}\}^{2}$
and then we integrating in $\partial_{v}$ as before, we can readily
bound recalling that $dud\text{\textgreek{sv}}\sim r^{-1-\text{\textgreek{h}}'}dvd\text{\textgreek{sv}}$
on $\mathcal{S}_{\text{\textgreek{t}}}$): 
\begin{align}
Err_{v,i,i}^{(p)}[\text{\textgreek{f}}](\text{\textgreek{t}}_{1},\text{\textgreek{t}}_{2})\lesssim & \int_{\mathcal{R}(\text{\textgreek{t}}_{1},\text{\textgreek{t}}_{2})\cap\{t\le T^{*}\}}\text{\textgreek{q}}_{R}\cdot\Big\{ r^{p-3-a}|\partial_{u}\text{\textgreek{F}}|^{2}+r^{p-2-a}|\partial_{u}\text{\textgreek{F}}|\cdot|\partial_{v}\text{\textgreek{F}}|^{2}+r^{p-4-a}|\partial_{\text{\textgreek{sv}}}\text{\textgreek{F}}|^{2}+r^{p-3-a}|\partial_{v}\text{\textgreek{F}}|\cdot|\partial_{\text{\textgreek{sv}}}\text{\textgreek{F}}|\Big\}\, dudvd\text{\textgreek{sv}}+\label{eq:BoundVIIErrTerms}\\
 & +\int_{\mathcal{R}(\text{\textgreek{t}}_{1},\text{\textgreek{t}}_{2})\cap\{t\le T^{*}\}}|\partial\text{\textgreek{q}}_{R}|\cdot r^{p-1}\cdot|\partial\text{\textgreek{F}}|^{2}\, dudvd\text{\textgreek{sv}}+\nonumber \\
 & +\sum_{i=1}^{2}\int_{\mathcal{S}_{\text{\textgreek{t}}_{i}}\cap\{t\le T^{*}\}}\text{\textgreek{q}}_{R}\Big\{ r^{p-2-a}|\partial_{v}\text{\textgreek{F}}|\cdot|\partial_{u}\text{\textgreek{F}}|+r^{p-4-a-\text{\textgreek{h}}'}|\partial_{\text{\textgreek{sv}}}\text{\textgreek{F}}|^{2}\Big\}\, dvd\text{\textgreek{sv}}+\nonumber \\
 & +\int_{\mathcal{R}(\text{\textgreek{t}}_{1},\text{\textgreek{t}}_{2})\cap\{t=T^{*}\}}\text{\textgreek{q}}_{R}\Big\{ r^{p-2-a}|\partial_{v}\text{\textgreek{F}}|\cdot|\partial_{u}\text{\textgreek{F}}|+r^{p-3-a}|\partial_{\text{\textgreek{sv}}}\text{\textgreek{F}}|^{2}\Big\}\, dvd\text{\textgreek{sv}}.\nonumber 
\end{align}
 Using a Cauchy--Schwarz inequality as well as inequality (\ref{eq:BoundFromMorawetzDuOnly}),
the first two lines of the right hand side of (\ref{eq:BoundVIIErrTerms})
can be bounded by $C_{p,\text{\textgreek{h}}}\cdot\mathcal{B}_{Err,\text{\textgreek{h}},\text{\textgreek{e}}}^{(p)}[\text{\textgreek{f}}](\text{\textgreek{t}}_{1},\text{\textgreek{t}}_{2})$.
Moreover, the last two lines of the right hand side of (\ref{eq:BoundVIIErrTerms})
can also be bounded by the boundry terms appearing in $\mathcal{B}_{Err,\text{\textgreek{h}},\text{\textgreek{e}}}^{(p)}[\text{\textgreek{f}}](\text{\textgreek{t}}_{1},\text{\textgreek{t}}_{2})$.
Thus, provided again that $R$ is large enough in terms of $p,\text{\textgreek{h}},\text{\textgreek{e}}$,
we have 
\begin{equation}
Err_{v,i,i}^{(p)}[\text{\textgreek{f}}](\text{\textgreek{t}}_{1},\text{\textgreek{t}}_{2})\lesssim_{p,\text{\textgreek{h}}}\mathcal{B}_{Err,\text{\textgreek{h}},\text{\textgreek{e}}}^{(p)}[\text{\textgreek{f}}](\text{\textgreek{t}}_{1},\text{\textgreek{t}}_{2}).\label{eq:FinalBoundVIIErrTerms}
\end{equation}

\smallskip{}

\noindent 3. Finally, it remains to bound the $\int O(r^{p-2-a})\cdot\partial_{v}\text{\textgreek{F}}\cdot\partial_{u}\partial_{\text{\textgreek{sv}}}\text{\textgreek{F}}$
summand in the expression (\ref{eq:ErrTerms-1}). We will perform
three integrations by parts schematically as 
\begin{equation}
\partial_{v}\text{\textgreek{F}}\cdot\partial_{\text{\textgreek{sv}}}\partial_{u}\text{\textgreek{F}}\rightarrow-\partial_{u}\partial_{v}\text{\textgreek{F}}\cdot\partial_{\text{\textgreek{sv}}}\text{\textgreek{F}}\rightarrow\partial_{u}\text{\textgreek{F}}\cdot\partial_{v}\partial_{\text{\textgreek{sv}}}\text{\textgreek{F}}\rightarrow-\partial_{\text{\textgreek{sv}}}\partial_{u}\text{\textgreek{F}}\cdot\partial_{v}\text{\textgreek{F}},
\end{equation}
and then move the resulting $\int O(r^{p-2-a})\cdot\partial_{v}\text{\textgreek{F}}\cdot\partial_{u}\partial_{\text{\textgreek{sv}}}\text{\textgreek{F}}$
bulk term (which is equal to $-1$ times the initial $O(r^{p-2-a})\cdot\partial_{v}\text{\textgreek{F}}\cdot\partial_{u}\partial_{\text{\textgreek{sv}}}\text{\textgreek{F}}$
term that we started with) from the right hand side to the left hand
side. This will provide an estimate of $\int O(r^{p-2-a})\cdot\partial_{v}\text{\textgreek{F}}\cdot\partial_{u}\partial_{\text{\textgreek{sv}}}\text{\textgreek{F}}$
by bulk and boundary terms which contain only first order derivatives
in $\text{\textgreek{F}}$. In particular, proceeding as described
we infer%
\footnote{Using again the fact that $dud\text{\textgreek{sv}}\sim r^{-1-\text{\textgreek{h}}'}dvd\text{\textgreek{sv}}$
on $\mathcal{S}_{\text{\textgreek{t}}}$%
} that 
\begin{align}
\int_{\mathcal{R}(\text{\textgreek{t}}_{1},\text{\textgreek{t}}_{2})\cap\{t\le T^{*}\}}\text{\textgreek{q}}_{R}O( & r^{p-2-a})\cdot\partial_{v}\text{\textgreek{F}}\cdot\partial_{u}\partial_{\text{\textgreek{sv}}}\text{\textgreek{F}}\, dudvd\text{\textgreek{sv}}\lesssim\int_{\mathcal{R}(\text{\textgreek{t}}_{1},\text{\textgreek{t}}_{2})\cap\{t\le T^{*}\}}\text{\textgreek{q}}_{R}r^{p-2-a}|\partial_{v}\text{\textgreek{F}}|\cdot|\partial_{\text{\textgreek{sv}}}\text{\textgreek{F}}|\, dudvd\text{\textgreek{sv}}+\\
 & +\int_{\mathcal{R}(\text{\textgreek{t}}_{1},\text{\textgreek{t}}_{2})\cap\{t\le T^{*}\}}\text{\textgreek{q}}_{R}r^{p-3-a}|\partial_{u}\text{\textgreek{F}}|\cdot|\partial_{\text{\textgreek{sv}}}\text{\textgreek{F}}|\, dudvd\text{\textgreek{sv}}+\int_{\mathcal{R}(\text{\textgreek{t}}_{1},\text{\textgreek{t}}_{2})\cap\{t\le T^{*}\}}\text{\textgreek{q}}_{R}r^{p-2-a}|\partial_{v}\text{\textgreek{F}}|\cdot|\partial_{u}\text{\textgreek{F}}|\, dudvd\text{\textgreek{sv}}+\nonumber \\
 & +\int_{\mathcal{R}(\text{\textgreek{t}}_{1},\text{\textgreek{t}}_{2})\cap\{t\le T^{*}\}}|\partial\text{\textgreek{q}}_{R}|r^{p}|\partial\text{\textgreek{F}}|^{2}\, dudvd\text{\textgreek{sv}}+\sum_{i=1}^{2}\int_{\mathcal{S}_{\text{\textgreek{t}}_{i}}\cap\{t\le T^{*}\}}\text{\textgreek{q}}_{R}r^{p-2-a}\big(|\partial_{v}\text{\textgreek{F}}|+r^{-1-\text{\textgreek{h}}'}|\partial_{u}\text{\textgreek{F}}|\big)\cdot|\partial_{\text{\textgreek{sv}}}\text{\textgreek{F}}|\, dvd\text{\textgreek{sv}}+\nonumber \\
 & +\int_{\mathcal{R}(\text{\textgreek{t}}_{1},\text{\textgreek{t}}_{2})\cap\{t=T^{*}\}}\text{\textgreek{q}}_{R}r^{p-2-a}\big(|\partial_{v}\text{\textgreek{F}}|+|\partial_{u}\text{\textgreek{F}}|\big)\cdot|\partial_{\text{\textgreek{sv}}}\text{\textgreek{F}}|\, dvd\text{\textgreek{sv}}.\nonumber 
\end{align}

\noindent Thus, after using a Cauchy--Schwarz inequality (as well
as the estimate (\ref{eq:BoundFromMorawetzDuOnly})), we obtain provided
that $R$ is large enough in terms of $p,\text{\textgreek{h}},\text{\textgreek{e}}$:
\begin{equation}
\int_{\mathcal{R}(\text{\textgreek{t}}_{1},\text{\textgreek{t}}_{2})\cap\{t\le T^{*}\}}\text{\textgreek{q}}_{R}O(r^{p-2-a})\cdot\partial_{v}\text{\textgreek{F}}\cdot\partial_{u}\partial_{\text{\textgreek{sv}}}\text{\textgreek{F}}\, dudvd\text{\textgreek{sv}}\lesssim_{p,\text{\textgreek{h}}}\mathcal{B}_{Err,\text{\textgreek{h}},\text{\textgreek{e}}}^{(p)}[\text{\textgreek{f}}](\text{\textgreek{t}}_{1},\text{\textgreek{t}}_{2}).\label{eq:BoundMixedErrTerms}
\end{equation}

\medskip{}

\noindent Inequalities (\ref{eq:FinalBoundIVErrTerms}), (\ref{eq:FinalBoundVIIErrTerms})
and (\ref{eq:BoundMixedErrTerms}) yield (\ref{eq:BoundForTheErrTerms}).
Therefore, the proof of the Lemma is complete.
\end{proof}
In order to control the zeroth order terms appearing in the right
hand side of (\ref{eq:NewMethodSimpleHyperboloids}), we will make
use of the following Hardy type inequality: 
\begin{lem}
\label{lem:NewMethodGeneralCaseHardy}For any $0<p<2$, any given
$0<\text{\textgreek{h}}<a$, any $R>0$ large enough in terms of $p,\text{\textgreek{h}}$,
any $\text{\textgreek{t}}_{1}\le\text{\textgreek{t}}_{2}$ and any
$T^{*}\ge0$, the following inequality is true for any smooth function
$\text{\textgreek{f}}:\mathcal{N}_{af}\rightarrow\mathbb{C}$ with
compact support in space and any smooth cut-off $\text{\textgreek{q}}_{R}:\mathcal{N}_{af}\rightarrow[0,1]$
supported in $\{r\ge R\}$:
\begin{equation}
\begin{split}\int_{\mathcal{S}_{\text{\textgreek{t}}_{2}}\cap\{t\le T^{*}\}}\text{\textgreek{q}}_{R}\cdot r^{p-2} & |\text{\textgreek{W}}\text{\textgreek{f}}|^{2}\, r^{-1-\text{\textgreek{h}}'}dvd\text{\textgreek{sv}}+\int_{\mathcal{R}(\text{\textgreek{t}}_{1},\text{\textgreek{t}}_{2})\cap\{t=T^{*}\}}\text{\textgreek{q}}_{R}\cdot r^{p-2}\cdot|\text{\textgreek{W}}\text{\textgreek{f}}|^{2}\, dvd\text{\textgreek{sv}}+\\
+\int_{\mathcal{R}(\text{\textgreek{t}}_{1},\text{\textgreek{t}}_{2})\cap\{t\le T^{*}\}} & \text{\textgreek{q}}_{R}\cdot r^{p-3}|\text{\textgreek{W}}\text{\textgreek{f}}|^{2}\, dudvd\text{\textgreek{sv}}\lesssim_{p,\text{\textgreek{h}}}\,\mathcal{E}_{bound,R,T^{*}}^{(p)}[\text{\textgreek{f}}](\text{\textgreek{t}}_{1})+\int_{\mathcal{S}_{\text{\textgreek{t}}_{1}}}\text{\textgreek{q}}_{R}\cdot r^{p-2}|\text{\textgreek{W}}\text{\textgreek{f}}|^{2}\, r^{-1-\text{\textgreek{h}}'}dvd\text{\textgreek{sv}}+\\
 & +\int_{\mathcal{R}(\text{\textgreek{t}}_{1},\text{\textgreek{t}}_{2})\cap\{t\le T^{*}\}}|\partial\text{\textgreek{q}}_{R}|\cdot\big(r^{p}|\partial(\text{\textgreek{W}}\text{\textgreek{f}})|^{2}+r^{p-2}|\text{\textgreek{W}}\text{\textgreek{f}}|^{2}\big)\, dudvd\text{\textgreek{sv}}+\int_{\mathcal{S}_{\text{\textgreek{t}}_{1}}\cap\{t\le T^{*}\}}\text{\textgreek{q}}_{R}J_{\text{\textgreek{m}}}^{T}(\text{\textgreek{f}})\bar{n}^{\text{\textgreek{m}}}+\\
 & +\int_{\mathcal{R}(\text{\textgreek{t}}_{1},\text{\textgreek{t}}_{2})\cap\{t\le T^{*}\}}\text{\textgreek{q}}_{R}\cdot(r^{p+1}+r^{1+\text{\textgreek{h}}})|\text{\textgreek{W}}\square\text{\textgreek{f}}|^{2}\, dudvd\text{\textgreek{sv}},
\end{split}
\label{eq:newMethodHardy}
\end{equation}

where $\bar{n}$ is the future directed unit normal on the hyperboloids
$\mathcal{S}_{\text{\textgreek{t}}}$. In the above, the constants
implicit in the $\lesssim_{p,\text{\textgreek{h}}}$ notation depend
only on $p,\text{\textgreek{h}}$ and on the geometry of $(\mathcal{N}_{af},g)$. \end{lem}
\begin{proof}
As in the proof of Lemma (\ref{lem:NewMethodGeneralCase}), we can
assume that $\text{\textgreek{f}}$ is real valued. We will also set
$\text{\textgreek{F}}\doteq\text{\textgreek{W}}\text{\textgreek{f}}$.

In dimensions $d\ge4$ (\ref{eq:newMethodHardy}) follows immediately
from Lemma (\ref{lem:NewMethodGeneralCase}), since for $d\ge4$ and
for $R$ large enough, the left hand side of (\ref{eq:NewMethodSimpleHyperboloids})
controls the left hand side of (\ref{eq:newMethodHardy}) and the
zeroth order terms in the right hand side of (\ref{eq:newMethodHardy})
can be absorbed into the left hand side. Therefore it only remains
to treat the case $d=3$.

The proof of this lemma follows the standard steps of proving a Hardy-type
inequality. Since $\partial_{v}r=1+o(1)$, we can write:

\begin{equation}
\int_{\mathcal{R}(\text{\textgreek{t}}_{1},\text{\textgreek{t}}_{2})\cap\{t\le T^{*}\}}\text{\textgreek{q}}_{R}\cdot r^{p-3}|\text{\textgreek{F}}|^{2}\, dudvd\text{\textgreek{sv}}=\int_{\mathcal{R}(\text{\textgreek{t}}_{1},\text{\textgreek{t}}_{2})\cap\{t\le T^{*}\}}\text{\textgreek{q}}_{R}\cdot\big\{\frac{1}{p-2}\partial_{v}(r^{p-2})+o(r^{p-3})\big\}|\text{\textgreek{F}}|^{2}\, dudvd\text{\textgreek{sv}}
\end{equation}
 and hence, if $R\gg1$: 

\begin{align}
\int_{\mathcal{R}(\text{\textgreek{t}}_{1},\text{\textgreek{t}}_{2})\cap\{t\le T^{*}\}}\text{\textgreek{q}}_{R}\cdot r^{p-3}(1+o(1))|\text{\textgreek{F}}|^{2}\, dudvd\text{\textgreek{sv}} & =\frac{1}{p-2}\int_{\mathcal{R}(\text{\textgreek{t}}_{1},\text{\textgreek{t}}_{2})\cap\{t\le T^{*}\}}\text{\textgreek{q}}_{R}\cdot\partial_{v}(r^{p-2})|\text{\textgreek{F}}|^{2}\, dudvd\text{\textgreek{sv}}=\label{eq:IntermediateHardyNewMethod}\\
 & =\frac{1}{p-2}\Big\{\int_{\mathcal{S}_{\text{\textgreek{t}}_{2}}\cap\{t\le T^{*}\}}\text{\textgreek{q}}_{R}\cdot r^{p-2}|\text{\textgreek{F}}|^{2}\, dud\text{\textgreek{sv}}+\int_{\mathcal{R}(\text{\textgreek{t}}_{1},\text{\textgreek{t}}_{2})\cap\{t=T^{*}\}}\text{\textgreek{q}}_{R}\cdot r^{p-2}|\text{\textgreek{F}}|^{2}\, dud\text{\textgreek{sv}}-\nonumber \\
 & \hphantom{=\frac{1}{p-2}\Big\{}-\int_{\mathcal{S}_{\text{\textgreek{t}}_{1}}\cap\{t\le T^{*}\}}\text{\textgreek{q}}_{R}\cdot r^{p-2}|\text{\textgreek{F}}|^{2}\, dud\text{\textgreek{sv}}-\nonumber \\
 & \hphantom{=\frac{1}{p-2}\Big\{}-\int_{\mathcal{R}(\text{\textgreek{t}}_{1},\text{\textgreek{t}}_{2})\cap\{t\le T^{*}\}}(\partial_{v}\text{\textgreek{q}}_{R})\cdot r^{p-2}|\text{\textgreek{F}}|^{2}\, dudvd\text{\textgreek{sv}}+\nonumber \\
 & \hphantom{=\frac{1}{p-2}\Big\{}+\int_{\mathcal{R}(\text{\textgreek{t}}_{1},\text{\textgreek{t}}_{2})\cap\{t\le T^{*}\}}\text{\textgreek{q}}_{R}\cdot r^{p-2}2\cdot\partial_{v}\text{\textgreek{F}}\cdot\text{\textgreek{F}}\, dudvd\text{\textgreek{sv}}\Big\},\nonumber 
\end{align}
the last equality following after performing an integration by parts
in $\partial_{v}$.

Notice that $\frac{1}{p-2}<0$ in the right hand side of (\ref{eq:IntermediateHardyNewMethod}).
Therefore, after applying a Cauchy-Schwarz inequality we obtain:
\begin{equation}
\begin{split}\int_{\mathcal{S}_{\text{\textgreek{t}}_{2}}\cap\{t\le T^{*}\}}\text{\textgreek{q}}_{R} & \cdot r^{p-2}|\text{\textgreek{F}}|^{2}\, dud\text{\textgreek{sv}}+\int_{\mathcal{R}(\text{\textgreek{t}}_{1},\text{\textgreek{t}}_{2})\cap\{t=T^{*}\}}\text{\textgreek{q}}_{R}\cdot r^{p-2}|\text{\textgreek{F}}|^{2}\, dud\text{\textgreek{sv}}+\int_{\mathcal{R}(\text{\textgreek{t}}_{1},\text{\textgreek{t}}_{2})\cap\{t\le T^{*}\}}\text{\textgreek{q}}_{R}\cdot r^{p-3}|\text{\textgreek{F}}|^{2}\, dudvd\text{\textgreek{sv}}\le\\
\le & C(p)\Big\{\int_{\mathcal{S}_{\text{\textgreek{t}}_{1}}\cap\{t\le T^{*}\}}\text{\textgreek{q}}_{R}\cdot r^{p-2}|\text{\textgreek{F}}|^{2}\, dud\text{\textgreek{sv}}+\int_{\mathcal{R}(\text{\textgreek{t}}_{1},\text{\textgreek{t}}_{2})\cap\{t\le T^{*}\}}|\partial\text{\textgreek{q}}_{R}|\cdot r^{p-2}|\text{\textgreek{F}}|^{2}\, dudvd\text{\textgreek{sv}}+\\
 & +\Big(\int_{\mathcal{R}(\text{\textgreek{t}}_{1},\text{\textgreek{t}}_{2})\cap\{t\le T^{*}\}}\text{\textgreek{q}}_{R}\cdot r^{p-1}|\partial_{v}\text{\textgreek{F}}|^{2}\, dudvd\text{\textgreek{sv}}\Big)^{1/2}\cdot\Big(\int_{\mathcal{R}(\text{\textgreek{t}}_{1},\text{\textgreek{t}}_{2})\cap\{t\le T^{*}\}}\text{\textgreek{q}}_{R}\cdot r^{p-3}|\text{\textgreek{F}}|^{2}\, dudvd\text{\textgreek{sv}}\Big)^{1/2}\Big\}.
\end{split}
\label{eq:AlmostThereHardyFinal}
\end{equation}

The second factor of the last term of the right hand side of (\ref{eq:AlmostThereHardyFinal})
can be absorbed into the left hand side, while the first factor of
the same term can be bounded by the left hand side of (\ref{eq:NewMethodSimpleHyperboloids}),
and thus:
\begin{equation}
\begin{split}\int_{\mathcal{S}_{\text{\textgreek{t}}_{2}}\cap\{t\le T^{*}\}} & \text{\textgreek{q}}_{R}\cdot r^{p-2}|\text{\textgreek{F}}|^{2}\, dud\text{\textgreek{sv}}+\int_{\mathcal{R}(\text{\textgreek{t}}_{1},\text{\textgreek{t}}_{2})\cap\{t=T^{*}\}}\text{\textgreek{q}}_{R}\cdot r^{p-2}|\text{\textgreek{F}}|^{2}\, dud\text{\textgreek{sv}}+\int_{\mathcal{R}(\text{\textgreek{t}}_{1},\text{\textgreek{t}}_{2})\cap\{t\le T^{*}\}}\text{\textgreek{q}}_{R}\cdot r^{p-3}|\text{\textgreek{F}}|^{2}\, dudvd\text{\textgreek{sv}}\lesssim_{p}\\
\lesssim_{p} & \int_{\mathcal{S}_{\text{\textgreek{t}}_{1}}\cap\{t\le T^{*}\}}\text{\textgreek{q}}_{R}\cdot r^{p-2}|\text{\textgreek{F}}|^{2}\, dud\text{\textgreek{sv}}+\mathcal{E}_{bound,R,T^{*}}^{(p)}[\text{\textgreek{f}}](\text{\textgreek{t}}_{1})+\int_{\mathcal{S}_{\text{\textgreek{t}}_{1}}\cap\{t\le T^{*}\}}\text{\textgreek{q}}_{R}J_{\text{\textgreek{m}}}^{T}(\text{\textgreek{f}})\bar{n}^{\text{\textgreek{m}}}+\\
 & +\int_{\mathcal{R}(\text{\textgreek{t}}_{1},\text{\textgreek{t}}_{2})\cap\{t\le T^{*}\}}|\partial\text{\textgreek{q}}_{R}|\cdot\big(r^{p}|\partial(\text{\textgreek{W}}\text{\textgreek{f}})|^{2}+r^{p-2}|\text{\textgreek{W}}\text{\textgreek{f}}|^{2}\big)\, dudvd\text{\textgreek{sv}}+\int_{\mathcal{R}(\text{\textgreek{t}}_{1},\text{\textgreek{t}}_{2})\cap\{t\le T^{*}\}}\text{\textgreek{q}}_{R}\cdot(r^{p+1}+r^{1+\text{\textgreek{h}}})|\text{\textgreek{W}}\square\text{\textgreek{f}}|^{2}\, dudvd\text{\textgreek{sv}}+\\
 & +\int_{\mathcal{S}_{\text{\textgreek{t}}_{2}}\cap\{t\le T^{*}\}}\text{\textgreek{q}}_{R}\cdot r^{p-2-a}|\text{\textgreek{F}}|^{2}\, dvd\text{\textgreek{sv}}+\int_{\mathcal{R}(\text{\textgreek{t}}_{1},\text{\textgreek{t}}_{2})\cap\{t=T^{*}\}}\text{\textgreek{q}}_{R}\cdot r^{p-2-a}|\text{\textgreek{F}}|^{2}\, dvd\text{\textgreek{sv}}+\\
 & +\int_{\mathcal{R}(\text{\textgreek{t}}_{1},\text{\textgreek{t}}_{2})\cap\{t\le T^{*}\}}\text{\textgreek{q}}_{R}\cdot\max\{r^{p-3-a},r^{-3}\}|\text{\textgreek{F}}|^{2}\, dudvd\text{\textgreek{sv}}.
\end{split}
\label{eq:AlmostThereHardyFinal-1}
\end{equation}
 If $R$ is large enough in terms of $p$, the last two lines of the
right hand side of (\ref{eq:AlmostThereHardyFinal-1}) can be absorbed
into the left hand side, yielding the desired inequality (\ref{eq:newMethodHardy}).
\end{proof}
\emph{Proof of Theorem \ref{thm:NewMethodFinalStatementHyperboloids}.
}The proof of Theorem \ref{thm:NewMethodFinalStatementHyperboloids}
follows readily by adding inequalities (\ref{eq:NewMethodSimpleHyperboloids})
for the given value of $p$ and (\ref{eq:newMethodHardy}) for $\min\{p,2-\text{\textgreek{d}}\}$
in place of $p$, using also (\ref{eq:BoundednessGeneralRadiative})
and (\ref{eq:MorawetzGeneralCaseRadiativeHyperboloids}), and letting
$T^{*}\rightarrow+\infty$. \qed

\section{\label{sec:The-improved--hierarchy}The improved $r^{p}$-weighted
energy hierarchy for higher order derivatives}

In \cite{Schlue2013}, Schlue established that on Schwarzschild exterior
spacetimes, commutation of the wave equation (\ref{eq:WaveEquation})
with the outgoing null vector field $\partial_{v}$ and the generators
of the spherical isometries leads to an improved $r^{p}$-weighted
hierarchy for $\partial_{v}(\text{\textgreek{W}}\text{\textgreek{f}})$
and $r^{-1}\partial_{\text{\textgreek{sv}}}(\text{\textgreek{W}}\text{\textgreek{f}})$.
This improvement of the $r$-weights in the hierarchy (\ref{eq:ModelR^pHierarchy})
for $\partial_{v}(\text{\textgreek{W}}\text{\textgreek{f}})$, $r^{-1}\partial_{\text{\textgreek{sv}}}(\text{\textgreek{W}}\text{\textgreek{f}})$
was fundamental in obtaining improved decay rates in $u$ for $T\text{\textgreek{f}}$,
and subsequently for $\text{\textgreek{f}}$ itself ($T$ being the
stationary Killing vector field of Schwarzschild spacetime). See \cite{Schlue2013}
for more details.

In this Section, we will extend and improve the results of \cite{Schlue2013}
to spacetime regions $(\mathcal{N}_{af},g)$ with $g$ of the form
(\ref{eq:MetricUR}). In particular, we will establish that higher
order derivatives of $\text{\textgreek{W}}\text{\textgreek{f}}$ in
directions tangential to the hyperboloids $\{\bar{t}=const\}$ satisfy
an $r^{p}$-weighted hierarchy similar to the one established in the
previous section for $\text{\textgreek{W}}\text{\textgreek{f}}$,
but for $p$ taking values larger than $2$. These improved estimates
will be crucial in the establishment of improved polynomial decay
rates for $\text{\textgreek{f}}$ in Section \ref{sec:Improved-polynomial-decay}. 
\begin{thm}
\label{thm:NewMethodDu+DvPhi}For any $k\in\mathbb{N}$, any $2k-2<p\le2k$
, any given $0<\text{\textgreek{h}}<a$ and $0<\text{\textgreek{d}}<1$,
any $R>0$ large enough in terms of $p,\text{\textgreek{h}},\text{\textgreek{d}},k$
and any $\text{\textgreek{t}}_{1}\le\text{\textgreek{t}}_{2}$, the
following inequality is true for any smooth function $\text{\textgreek{f}}:\mathcal{M}\rightarrow\mathbb{C}$
solving $\square_{g}\text{\textgreek{f}}=F$ and any smooth cut-off
$\text{\textgreek{q}}_{R}:\mathcal{M}\rightarrow[0,1]$ supported
in $\{r\ge R\}$: 
\begin{equation}
\begin{split}\mathcal{E}_{bound,R;\text{\textgreek{d}}}^{(p,k)}[\text{\textgreek{f}}](\text{\textgreek{t}}_{2})+ & \int_{\text{\textgreek{t}}_{1}}^{\text{\textgreek{t}}_{2}}\mathcal{E}_{bulk,R,\text{\textgreek{h}};\text{\textgreek{d}}}^{(p-1,k)}[\text{\textgreek{f}}](\text{\textgreek{t}})\, d\text{\textgreek{t}}+\limsup_{T^{*}\rightarrow+\infty}\mathcal{E}_{\mathcal{I}^{+},T^{*},\text{\textgreek{d}}}^{(p,k)}[\text{\textgreek{f}}](\text{\textgreek{t}}_{1},\text{\textgreek{t}}_{2})\lesssim_{p,\text{\textgreek{h}},\text{\textgreek{d}},k}\\
\lesssim_{p,\text{\textgreek{h}},\text{\textgreek{d}},k}\, & \mathcal{E}_{bound,R;\text{\textgreek{d}}}^{(p,k)}[\text{\textgreek{f}}](\text{\textgreek{t}}_{2})+\sum_{j=0}^{k}\int_{\mathcal{R}(\text{\textgreek{t}}_{1},\text{\textgreek{t}}_{2})}|\partial\text{\textgreek{q}}_{R}|\cdot r^{p-2(k-j)}|\partial^{j}\text{\textgreek{f}}|^{2}+\\
 & +\sum_{j=1}^{k}\sum_{k_{1}+k_{2}+k_{3}=j-1}\int_{\mathcal{R}(\text{\textgreek{t}}_{1},\text{\textgreek{t}}_{2})}\text{\textgreek{q}}_{R}\cdot(r^{p+1-2k_{3}-2(k-j)}+r^{1+\text{\textgreek{h}}}\big)\Big(|r^{-k_{2}}\partial_{v}^{k_{1}}\partial_{\text{\textgreek{sv}}}^{k_{2}}\partial_{u}^{k_{3}}(\text{\textgreek{W}}F)|^{2}\Big)\, dudvd\text{\textgreek{sv}}.
\end{split}
\label{eq:newMethodDu+DvPhi}
\end{equation}

In the above, the constants implicit in the $\lesssim_{p,\text{\textgreek{h}},k}$
notation depend only on $p,\text{\textgreek{h}},k$ and on the geometry
of $(\mathcal{M},g)$, and the higher order $p$-energy norms are
defined as 
\begin{align}
\mathcal{E}_{bound,R;\text{\textgreek{d}}}^{(p,k)}[\text{\textgreek{f}}](\text{\textgreek{t}})=\sum_{j=1}^{k}\sum_{k_{1}+k_{2}+k_{3}=j-1}\Big\{ & \int_{\mathcal{S}_{\text{\textgreek{t}}}}\text{\textgreek{q}}_{R}\Big(r^{p-2(k-j)}\big|r^{-k_{2}-k_{3}}\partial_{v}^{k_{1}+1}\partial_{\text{\textgreek{sv}}}^{k_{2}}\partial_{u}^{k_{3}}(\text{\textgreek{W}}\text{\textgreek{f}})\big|^{2}+r^{-1-\text{\textgreek{h}}'}\big(r^{p-2(k-j)}\big|r^{-k_{2}-k_{3}-1}\partial_{v}^{k_{1}}\partial_{\text{\textgreek{sv}}}^{k_{2}+1}\partial_{u}^{k_{3}}(\text{\textgreek{W}}\text{\textgreek{f}})\big|^{2}+\\
 & \hphantom{\int_{\mathcal{S}_{\text{\textgreek{t}}}}}+\big((d-3)r^{p-2-2(k-j)}+\min\{r^{p-2-2(k-j)},r^{-\text{\textgreek{d}}-2(k-j)}\}\big)\big|r^{-k_{2}-k_{3}}\partial_{v}^{k_{1}}\partial_{\text{\textgreek{sv}}}^{k_{2}}\partial_{u}^{k_{3}}(\text{\textgreek{W}}\text{\textgreek{f}})\big|^{2}\big)\Big)\, dvd\text{\textgreek{sv}}+\nonumber \\
 & +\int_{\mathcal{S}_{\text{\textgreek{t}}}}\text{\textgreek{q}}_{R}J_{\text{\textgreek{m}}}^{T}(r^{-k_{2}}\partial_{v}^{k_{1}}\partial_{\text{\textgreek{sv}}}^{k_{2}}\partial_{u}^{k_{3}}\text{\textgreek{f}})\bar{n}^{\text{\textgreek{m}}}\Big\}\nonumber 
\end{align}
 
\begin{align}
\mathcal{E}_{bulk,R,\text{\textgreek{h}};\text{\textgreek{d}}}^{(p-1,k)}[\text{\textgreek{f}}](\text{\textgreek{t}})=\sum_{j=1}^{k} & \sum_{k_{1}+k_{2}+k_{3}=j-1}\Big\{\int_{\mathcal{S}_{\text{\textgreek{t}}}}\text{\textgreek{q}}_{R}\Big(pr^{p-1-2(k-j)}\big|r^{-k_{2}-k_{3}}\partial_{v}^{k_{1}+1}\partial_{\text{\textgreek{sv}}}^{k_{2}}\partial_{u}^{k_{3}}(\text{\textgreek{W}}\text{\textgreek{f}})\big|^{2}+\\
 & \hphantom{\int_{\mathcal{S}_{\text{\textgreek{t}}}}}+\big\{\big((2k-p)r^{p-1-2(k-j)}+r^{p-1-\text{\textgreek{d}}-2(k-j)}\big)\big|r^{-k-1}\partial_{\text{\textgreek{sv}}}^{k_{1}+k_{2}+1}\partial_{u}^{k_{3}}(\text{\textgreek{W}}\text{\textgreek{f}})\big|^{2}+\nonumber \\
 & \hphantom{\int_{\mathcal{S}_{\text{\textgreek{t}}}}}+\big((2k-p)(d-3)r^{p-3-2(k-j)}+\min\{r^{p-3-2(k-j)},r^{-1-\text{\textgreek{d}}-2(k-j)}\}\big)\big|r^{-k_{2}-k_{3}}\partial_{v}^{k_{1}}\partial_{\text{\textgreek{sv}}}^{k_{2}}\partial_{u}^{k_{3}}(\text{\textgreek{W}}\text{\textgreek{f}})\big|^{2}\big\}\Big)\, dvd\text{\textgreek{sv}}+\nonumber \\
 & \hphantom{\int_{\mathcal{S}_{\text{\textgreek{t}}}}}+\int_{\mathcal{S}_{\text{\textgreek{t}}}}\text{\textgreek{q}}_{R}r^{-1-\text{\textgreek{h}}}|\partial_{u}^{j+1}(\text{\textgreek{W}}\text{\textgreek{f}})|^{2}\, dvd\text{\textgreek{sv}}\Big\}\nonumber 
\end{align}
and \textup{
\begin{align}
\mathcal{E}_{\mathcal{I}^{+},T^{*},\text{\textgreek{d}}}^{(p,k)}[\text{\textgreek{f}}](\text{\textgreek{t}}_{1},\text{\textgreek{t}}_{2})=\sum_{j=1}^{k}\sum_{k_{1}+k_{2}+k_{3}=j-1}\Big\{ & \int_{\mathcal{R}(\text{\textgreek{t}}_{1},\text{\textgreek{t}}_{2})\cap\{t=T^{*}\}}\Big(r^{p-2(k-j)}\big|r^{-k_{2}-k_{3}-1}\partial_{v}^{k_{1}}\partial_{\text{\textgreek{sv}}}^{k_{2}+1}\partial_{u}^{k_{3}}(\text{\textgreek{W}}\text{\textgreek{f}})\big|^{2}+\\
 & +\big((d-3)r^{p-2-2(k-j)}+\min\{r^{p-2-2(k-j)},r^{-\text{\textgreek{d}}-2(k-j)}\}\big)\big|r^{-k_{2}-k_{3}}\partial_{v}^{k_{1}}\partial_{\text{\textgreek{sv}}}^{k_{2}}\partial_{u}^{k_{3}}(\text{\textgreek{W}}\text{\textgreek{f}})\big|^{2}\Big)\, dvd\text{\textgreek{sv}}+\nonumber \\
 & +\int_{\mathcal{R}(\text{\textgreek{t}}_{1},\text{\textgreek{t}}_{2})\cap\{t=T^{*}\}}J_{\text{\textgreek{m}}}^{T}(r^{-k_{2}}\partial_{v}^{k_{1}}\partial_{\text{\textgreek{sv}}}^{k_{2}}\partial_{u}^{k_{3}}\text{\textgreek{f}})n^{\text{\textgreek{m}}}\Big\}.\nonumber 
\end{align}
}
\end{thm}
Applying Theorem \ref{thm:NewMethodDu+DvPhi} for $\partial^{j}\text{\textgreek{f}}$,
$j=0,\ldots m-1$ in place of $\text{\textgreek{f}}$ using Lemma
\ref{lem:Commutator expressions}, we will deduce the following estimate
(the proof of which is straightforward and will be omitted)
\begin{cor}
\label{cor:HigherOrderNewMethodAddition}Keeping the same notations
as in Theorem \ref{thm:NewMethodDu+DvPhi}, for any any integer $m\ge0$,
any $k\in\mathbb{N}$, any $2k-2<p\le2k$ , any given $0<\text{\textgreek{h}}<a$
and $0<\text{\textgreek{d}}<1$, any $R>0$ large enough in terms
of $p,\text{\textgreek{h}},\text{\textgreek{d}},k,m$ and any $\text{\textgreek{t}}_{1}\le\text{\textgreek{t}}_{2}$,
the following inequality is true for any smooth function $\text{\textgreek{f}}:\mathcal{M}\rightarrow\mathbb{C}$
solving $\square_{g}\text{\textgreek{f}}=F$ and any smooth cut-off
$\text{\textgreek{q}}_{R}:\mathcal{M}\rightarrow[0,1]$ supported
in $\{r\ge R\}$: 
\begin{equation}
\begin{split}\sum_{j=0}^{m-1}\sum_{j_{1}+j_{2}+j_{3}=j} & \Big\{\mathcal{E}_{bound,R;\text{\textgreek{d}}}^{(p,k)}[r^{-j_{2}}\partial_{v}^{j_{1}}\partial_{\text{\textgreek{sv}}}^{j_{2}}\partial_{u}^{j_{3}}\text{\textgreek{f}}](\text{\textgreek{t}}_{2})+\int_{\text{\textgreek{t}}_{1}}^{\text{\textgreek{t}}_{2}}\mathcal{E}_{bulk,R,\text{\textgreek{h}};\text{\textgreek{d}}}^{(p-1,k)}[r^{-j_{2}}\partial_{v}^{j_{1}}\partial_{\text{\textgreek{sv}}}^{j_{2}}\partial_{u}^{j_{3}}\text{\textgreek{f}}](\text{\textgreek{t}})\, d\text{\textgreek{t}}+\limsup_{T^{*}\rightarrow+\infty}\mathcal{E}_{\mathcal{I}^{+},T^{*},\text{\textgreek{d}}}^{(p,k)}[r^{-j_{2}}\partial_{v}^{j_{1}}\partial_{\text{\textgreek{sv}}}^{j_{2}}\partial_{u}^{j_{3}}\text{\textgreek{f}}](\text{\textgreek{t}}_{1},\text{\textgreek{t}}_{2})\Big\}\\
\lesssim_{p,\text{\textgreek{h}},\text{\textgreek{d}},k,m} & \sum_{j=0}^{m-1}\sum_{j_{1}+j_{2}+j_{3}=j}\mathcal{E}_{bound,R;\text{\textgreek{d}}}^{(p,k)}[r^{-j_{2}}\partial_{v}^{j_{1}}\partial_{\text{\textgreek{sv}}}^{j_{2}}\partial_{u}^{j_{3}}\text{\textgreek{f}}](\text{\textgreek{t}}_{2})+\sum_{j=0}^{k}\sum_{i=0}^{m-1}\int_{\mathcal{R}(\text{\textgreek{t}}_{1},\text{\textgreek{t}}_{2})}|\partial\text{\textgreek{q}}_{R}|\cdot r^{p-2(k-j)}|\partial^{j+i}\text{\textgreek{f}}|^{2}+\\
 & +\sum_{i_{1}+i_{2}+i_{3}\le m-1}\sum_{j=1}^{k}\sum_{k_{1}+k_{2}+k_{3}=j-1}\int_{\mathcal{R}(\text{\textgreek{t}}_{1},\text{\textgreek{t}}_{2})}\text{\textgreek{q}}_{R}\cdot(r^{p+1-2k_{3}-2(k-j)}+r^{1+\text{\textgreek{h}}}\big)\Big(|r^{-k_{2}-i_{2}}\partial_{v}^{k_{1}+i_{1}}\partial_{\text{\textgreek{sv}}}^{k_{2}+i_{2}}\partial_{u}^{k_{3}+i_{3}}(\text{\textgreek{W}}F)|^{2}\Big)\, dudvd\text{\textgreek{sv}}.
\end{split}
\label{eq:newMethodDu+DvPhiHigherOrderAddition}
\end{equation}

\end{cor}
The proof of Theorem \ref{thm:NewMethodDu+DvPhi} will be presented
in Section \ref{sub:Proof-of-Improved Hierarchy}. Before that, we
will obtain an expression for the equation satisfied by the derivatives
of $\text{\textgreek{W}}\text{\textgreek{f}}$ when $\text{\textgreek{f}}$
solves $\square\text{\textgreek{f}}=F$. This will be accomplished
in the following section.

\subsection{\label{sub:Commutator-expressions}Commutator expressions}

In this section, we will commute the wave operator $\square_{g}$
with the coordinate vector fields $\partial_{u}$ and $\partial_{v}$
in the $(u,v,\text{\textgreek{sv}})$ coordinate system, as well as
the first order operator $r^{-1}\nabla^{\mathbb{S}^{d-1}}$ ($\nabla^{\mathbb{S}^{d-1}}$
denoting here the gradient of a function on $(\mathbb{S}^{d-1},g_{\mathbb{S}^{d-1}})$)
in the asymptotically flat region $\{r\gg1\}$.

We will establish the following lemma: 
\begin{lem}
\label{lem:Commutator expressions}For any smooth function $\text{\textgreek{f}}:\mathcal{M}\rightarrow\mathbb{C}$
the following expressions are true in the region $\{r\gg1\}$ for
any $l\in\mathbb{N}$: 
\begin{align}
\text{\textgreek{W}}\square\big(\text{\textgreek{W}}^{-1}\partial_{v}^{l}(\text{\textgreek{W}}\text{\textgreek{f}})\big)= & \partial_{v}^{l}\big(\text{\textgreek{W}}\square\text{\textgreek{f}}\big)+\sum_{j=1}^{l}(-1)^{j+1}\Big\{\binom{l}{j}(1+j)!\cdot r^{-2-j}\big(\text{\textgreek{D}}_{g_{\mathbb{S}^{d-1}}+h_{\mathbb{S}^{d-1}}}\partial_{v}^{l-j}(\text{\textgreek{W}}\text{\textgreek{f}})-\frac{(d-1)(d-3)}{4}\partial_{v}^{l-j}(\text{\textgreek{W}}\text{\textgreek{f}})\big)\Big\}+\label{eq:CommutationDv}\\
 & +\sum_{j=0}^{l}r^{-(l-j)}Err(\partial_{v}^{j}(\text{\textgreek{W}}\text{\textgreek{f}})),\nonumber 
\end{align}
 
\begin{align}
\text{\textgreek{W}}\square\big(\text{\textgreek{W}}^{-1}(r^{-1}\nabla^{\mathbb{S}^{d-1}})^{l}(\text{\textgreek{W}}\text{\textgreek{f}})\big)= & \big(r^{-1}\nabla^{\mathbb{S}^{d-1}}\big)^{l}\big(\text{\textgreek{W}}\square\text{\textgreek{f}}\big)+l\cdot r^{-l-1}\partial_{u}\big((\nabla^{\mathbb{S}^{d-1}})^{l}(\text{\textgreek{W}}\text{\textgreek{f}})\big)-l\cdot r^{-l-1}\partial_{v}\big((\nabla^{\mathbb{S}^{d-1}})^{l}(\text{\textgreek{W}}\text{\textgreek{f}})\big)+\label{eq:CommutationDsigma}\\
 & +l(l+1)r^{-l-2}(\nabla^{\mathbb{S}^{d-1}})^{l}(\text{\textgreek{W}}\text{\textgreek{f}})+\sum_{j=0}^{l}r^{-(l-j)}Err\big((r^{-1}\nabla^{\mathbb{S}^{d-1}})^{j}(\text{\textgreek{W}}\text{\textgreek{f}})\big)\nonumber 
\end{align}
 and

\begin{align}
\text{\textgreek{W}}\square\big(\text{\textgreek{W}}^{-1}\partial_{u}^{l}(\text{\textgreek{W}}\text{\textgreek{f}})\big)= & (1+O(r^{-1-a}))\cdot\partial_{u}^{l}\big((1+O(r^{-1-a}))\cdot\text{\textgreek{W}}\square\text{\textgreek{f}}\big)+\label{eq:CommutationDu+Dv}\\
 & +\sum_{j=0}^{l-1}O(r^{-3})\partial_{\text{\textgreek{sv}}}\partial_{\text{\textgreek{sv}}}\big(\partial_{u}^{j}(\text{\textgreek{W}}\text{\textgreek{f}})\big)+\sum_{j=0}^{l}Err\big(\partial_{u}^{j}(\text{\textgreek{W}}\text{\textgreek{f}})\big),\nonumber 
\end{align}
 where $\nabla^{\mathbb{S}^{d-1}}$ denotes the gradient on $(\mathbb{S}^{d-1},g_{\mathbb{S}^{d-1}})$
and the $Err$ terms are of the form: 
\begin{align}
Err(\text{\textgreek{F}})= & O(r^{-2-a})\cdot\partial_{u}^{2}\text{\textgreek{F}}+O(r^{-1})\cdot\partial_{v}^{2}\text{\textgreek{F}}+O(r^{-2-a})\partial_{u}\partial_{\text{\textgreek{sv}}}\text{\textgreek{F}}+O(r^{-2})\partial_{v}\partial_{\text{\textgreek{sv}}}\text{\textgreek{F}}+\\
 & +O(r^{-3-a})\partial_{\text{\textgreek{sv}}}\partial_{\text{\textgreek{sv}}}\text{\textgreek{F}}+O(r^{-2-a})\partial_{u}\text{\textgreek{F}}+O(r^{-1-a})\partial_{v}\text{\textgreek{F}}+O(r^{-2-a})\cdot\partial_{\text{\textgreek{sv}}}\text{\textgreek{F}}+O(r^{-3})\text{\textgreek{F}}.\nonumber 
\end{align}
\end{lem}
\begin{proof}
Let us set 
\begin{equation}
\text{\textgreek{F}}\doteq\text{\textgreek{W}}\text{\textgreek{f}}.
\end{equation}
 In $(u,v,\text{\textgreek{sv}})$ coordinates in the region $\{r\gg1\}$,
the wave operator takes the following form according to (\ref{eq:ConformalWaveOperator}): 

\begin{align}
\text{\textgreek{W}}\cdot\square\text{\textgreek{f}}= & -\big(1+O(r^{-1-a})\big)\cdot\partial_{u}\partial_{v}\text{\textgreek{F}}+r^{-2}\text{\textgreek{D}}_{g_{\mathbb{S}^{d-1}}+h_{\mathbb{S}^{d-1}}}\text{\textgreek{F}}-\label{eq:ConformalWaveOperator-1}\\
 & -\frac{(d-1)(d-3)}{4}r^{-2}\cdot\text{\textgreek{F}}+Err(\text{\textgreek{F}}).\nonumber 
\end{align}

Differentiating (\ref{eq:ConformalWaveOperator-1}) $l$ times with
respect to $\partial_{v}$, we readily obtain: 
\begin{align}
\partial_{v}^{l}\big(\text{\textgreek{W}}\square\text{\textgreek{f}}\big)= & -\big(1+O(r^{-1-a})\big)\cdot\partial_{u}\partial_{v}(\partial_{v}^{l}\text{\textgreek{F}})+\\
 & +\sum_{j=0}^{l}(-1)^{j}\Big\{\binom{l}{j}(1+j)!\cdot r^{-2-j}\big(\text{\textgreek{D}}_{g_{\mathbb{S}^{d-1}}+h_{\mathbb{S}^{d-1}}}\partial_{v}^{l-j}\text{\textgreek{F}}-\frac{(d-1)(d-3)}{4}\partial_{v}^{l-j}\text{\textgreek{F}}\big)\Big\}+\nonumber \\
 & +\sum_{j=0}^{l}r^{-(l-j)}Err(\partial_{v}^{j}\text{\textgreek{F}}),\nonumber 
\end{align}
 which (due to the expression (\ref{eq:ConformalWaveOperator-1})
for the wave operator) can be rewritten in the desired form (\ref{eq:CommutationDv})
as: 
\begin{align}
\text{\textgreek{W}}\square\big(\text{\textgreek{W}}^{-1}\partial_{v}\text{\textgreek{F}}\big)= & \partial_{v}\big(\text{\textgreek{W}}\square\text{\textgreek{f}}\big)+\sum_{j=1}^{l}(-1)^{j+1}\Big\{\binom{l}{j}(1+j)!\cdot r^{-2-j}\big(\text{\textgreek{D}}_{g_{\mathbb{S}^{d-1}}+h_{\mathbb{S}^{d-1}}}\partial_{v}^{l-j}\text{\textgreek{F}}-\frac{(d-1)(d-3)}{4}\partial_{v}^{l-j}\text{\textgreek{F}}\big)\Big\}+\label{eq:CommutationDv-1}\\
 & +\sum_{j=0}^{l}r^{-(l-j)}Err(\partial_{v}^{j}\text{\textgreek{F}}).\nonumber 
\end{align}

Similarly, applying to (\ref{eq:ConformalWaveOperator-1}) $l$ times
the rescaled angular gradient $r^{-1}\nabla^{\mathbb{S}^{d-1}}$ to
(\ref{eq:ConformalWaveOperator-1}) (notice that $\nabla^{\mathbb{S}^{d-1}}$
commutes with $\text{\textgreek{D}}_{g_{\mathbb{S}^{d-1}}},\partial_{u},\partial_{v}$
and $r$), we calculate: 
\begin{align}
(r^{-1}\nabla^{\mathbb{S}^{d-1}})^{l}\big(\text{\textgreek{W}}\square\text{\textgreek{f}}\big)= & -\big(1+O(r^{-1-a})\big)\cdot\partial_{u}\partial_{v}(r^{-l}(\nabla^{\mathbb{S}^{d-1}})^{l}\text{\textgreek{F}})+r^{-2}\text{\textgreek{D}}_{g_{\mathbb{S}^{d-1}}+h_{\mathbb{S}^{d-1}}}(r^{-l}(\nabla^{\mathbb{S}^{d-1}})^{l}\text{\textgreek{F}})-\\
 & -\frac{(d-1)(d-3)}{2}r^{-2}\cdot(r^{-l}(\nabla^{\mathbb{S}^{d-1}})^{l}\text{\textgreek{F}})-l\cdot r^{-l-1}\partial_{u}\big((\nabla^{\mathbb{S}^{d-1}})^{l}\text{\textgreek{F}}\big)+l\cdot r^{-l-1}\partial_{v}\big((\nabla^{\mathbb{S}^{d-1}})^{l}\text{\textgreek{F}}\big)-\nonumber \\
 & -l(l+1)r^{-l-2}(\nabla^{\mathbb{S}^{d-1}})^{l}\text{\textgreek{F}}+\sum_{j=0}^{l}r^{-(l-j)}Err\big((r^{-1}\nabla^{\mathbb{S}^{d-1}}))^{j}\text{\textgreek{F}}\big),\nonumber 
\end{align}
 (where the $O(\cdot)$ terms in the $r^{-1}Err(\text{\textgreek{F}})$
summand should be considered to denote vector fields on $\mathbb{S}^{d-1}$
rather than functions). Thus, rearranging the terms we deduce that:
\begin{align}
\text{\textgreek{W}}\square\big(\text{\textgreek{W}}^{-1}(r^{-1}\nabla^{\mathbb{S}^{d-1}})^{l}\text{\textgreek{F}}\big)= & \big(r^{-1}\nabla^{\mathbb{S}^{d-1}}\big)^{l}\big(\text{\textgreek{W}}\square\text{\textgreek{f}}\big)+l\cdot r^{-l-1}\partial_{u}\big((\nabla^{\mathbb{S}^{d-1}})^{l}\text{\textgreek{F}}\big)-l\cdot r^{-l-1}\partial_{v}\big((\nabla^{\mathbb{S}^{d-1}})^{l}\text{\textgreek{F}}\big)+\label{eq:CommutationDsigma-1}\\
 & +l(l+1)r^{-l-2}(\nabla^{\mathbb{S}^{d-1}})^{l}\text{\textgreek{F}}+\sum_{j=0}^{l}r^{-(l-j)}Err\big((r^{-1}\nabla^{\mathbb{S}^{d-1}})^{j}\text{\textgreek{F}}\big)\nonumber 
\end{align}

After multiplying (\ref{eq:ConformalWaveOperator-1}) with $1+O(r^{-1-a})$
so as to make the coefficient of $\partial_{u}\partial_{v}\text{\textgreek{F}}$
equal to $1$, differentiating $l$ times with respect to $\partial_{u}$
we eventually obtain (since $\partial_{u}r=-1$) that 
\begin{align}
\text{\textgreek{W}}\square\big(\text{\textgreek{W}}^{-1}\partial_{u}^{l}\text{\textgreek{F}}\big)= & (1+O(r^{-1-a}))\cdot\partial_{u}^{l}\big((1+O(r^{-1-a}))\cdot\text{\textgreek{W}}\square\text{\textgreek{f}}\big)+\label{eq:CommutationDu+Dv-1}\\
 & +\sum_{j=0}^{l-1}O(r^{-3})\partial_{\text{\textgreek{sv}}}\partial_{\text{\textgreek{sv}}}\big(\partial_{u}^{j}\text{\textgreek{F}}\big)+\sum_{j=0}^{l}Err\big(\partial_{u}^{j}\text{\textgreek{F}}\big),\nonumber 
\end{align}

\end{proof}

\subsection{\label{sub:Proof-of-Improved Hierarchy}Proof of Theorem \ref{thm:NewMethodDu+DvPhi}}

Without loss of generality, we will assume that $\text{\textgreek{f}}$
is real valued. We will set 
\begin{equation}
\text{\textgreek{F}}\doteq\text{\textgreek{W}}\text{\textgreek{f}}.
\end{equation}

In the case $k=1$ the statement of Theorem \ref{thm:NewMethodDu+DvPhi}
reduces to the statement of Theorem \ref{thm:NewMethodFinalStatementHyperboloids}.
Thus, it suffices to assume that $k\ge2$. In order to avoid unnecessarily
complicated notations, we will only establish the case $k=2$, since
the case $k>2$ can be treated in exactly the same way (through induction
on $k$). 

Fix an $\text{\textgreek{e}}>0$ small enough in terms of $p,\text{\textgreek{h}},\text{\textgreek{d}}$.
We will assume without loss of generality that $R$ in terms of $\text{\textgreek{e}}$.
Repeating the proof of Lemma \ref{lem:NewMethodGeneralCase} for $\text{\textgreek{W}}^{-1}\partial_{v}\text{\textgreek{F}}$
in place of $\text{\textgreek{f}}$, but without absorbing terms of
the form $\int O(r^{p-1-a})|r^{-1}\partial_{\text{\textgreek{sv}}}\partial_{v}\text{\textgreek{F}}|^{2}$
by the left hand side of the resulting inequality, and without using
a Cauchy--Schwarz inequality to bound terms of the form $\partial^{2}\text{\textgreek{F}}\cdot\square\big(\text{\textgreek{W}}\partial_{v}\text{\textgreek{F}}\big)$,
we readily obtain the following inequality for any $T^{*}\ge0$ (using
the energy norm notation of the proof of Theorem \ref{thm:NewMethodFinalStatementHyperboloids}):
\begin{equation}
\begin{split}\mathcal{E}_{bound,R,T^{*}}^{(p)}[ & \text{\textgreek{W}}^{-1}\partial_{v}\text{\textgreek{F}}](\text{\textgreek{t}}_{2})+\int_{\text{\textgreek{t}}_{1}}^{\text{\textgreek{t}}_{2}}\mathcal{E}_{bulk,R,\text{\textgreek{h}},T^{*}}^{(p-1)}[\text{\textgreek{W}}^{-1}\partial_{v}\text{\textgreek{F}}](\text{\textgreek{t}})\, d\text{\textgreek{t}}+\mathcal{E}_{\mathcal{I}^{+},R,T^{*}}^{(p)}[\text{\textgreek{W}}^{-1}\partial_{v}\text{\textgreek{F}}](\text{\textgreek{t}}_{1},\text{\textgreek{t}}_{2})\le\\
\le & \big(1+O_{p,\text{\textgreek{h}}}(\text{\textgreek{e}})\big)\cdot\mathcal{E}_{bound,R,T^{*}}^{(p)}[\text{\textgreek{W}}^{-1}\partial_{v}\text{\textgreek{F}}](\text{\textgreek{t}}_{1})+C_{p,\text{\textgreek{h}},\text{\textgreek{e}}}\int_{\mathcal{R}(\text{\textgreek{t}}_{1},\text{\textgreek{t}}_{2})\cap\{t\le T^{*}\}}|\partial\text{\textgreek{q}}_{R}|\cdot\big(r^{p}|\partial\partial_{v}\text{\textgreek{F}}|^{2}+r^{p-2}|\partial_{v}\text{\textgreek{F}}|^{2}\big)\, dudvd\text{\textgreek{sv}}+\\
 & +C_{p,\text{\textgreek{h}},\text{\textgreek{e}}}\int_{\mathcal{R}(\text{\textgreek{t}}_{1},\text{\textgreek{t}}_{2})\cap\{t\le T^{*}\}}\text{\textgreek{q}}_{R}\cdot\Big(r^{p-1-a}|r^{-1}\partial_{\text{\textgreek{sv}}}\partial_{v}\text{\textgreek{F}}|^{2}+r^{p-3-a}|\partial_{u}\partial_{v}\text{\textgreek{F}}|^{2}+\max\{r^{p-3-a},r^{-3}\}|\partial_{v}\text{\textgreek{F}}|^{2}\Big)\, dudvd\text{\textgreek{sv}}+\\
 & +C_{p,\text{\textgreek{h}},\text{\textgreek{e}}}\int_{\mathcal{R}(\text{\textgreek{t}}_{1},\text{\textgreek{t}}_{2})\cap\{t\le T^{*}\}}\text{\textgreek{q}}_{R}\cdot\big(\frac{1}{2}f(r)(\partial_{v}-\partial_{u})(\text{\textgreek{W}}^{-1}\partial_{v}\text{\textgreek{F}})+\frac{(d-1)f(r)}{r}(\text{\textgreek{W}}^{-1}\partial_{v}\text{\textgreek{F}})\big)\cdot\square(\text{\textgreek{W}}^{-1}\partial_{v}\text{\textgreek{F}})\,\text{\textgreek{W}}^{2}dudvd\text{\textgreek{sv}}+\\
 & +C_{p,\text{\textgreek{h}},\text{\textgreek{e}}}\int_{\mathcal{R}(\text{\textgreek{t}}_{1},\text{\textgreek{t}}_{2})\cap\{t\le T^{*}\}}\text{\textgreek{q}}_{R}\cdot(\partial_{v}+\partial_{u})(\partial_{v}\text{\textgreek{F}})\cdot\text{\textgreek{W}}\square(\text{\textgreek{W}}^{-1}\partial_{v}\text{\textgreek{F}})\, dudvd\text{\textgreek{sv}}+\\
 & +\int_{\mathcal{R}(\text{\textgreek{t}}_{1},\text{\textgreek{t}}_{2})\cap\{t\le T^{*}\}}\text{\textgreek{q}}_{R}\cdot O_{p,\text{\textgreek{h}},\text{\textgreek{e}}}(r^{p-2-a})\cdot\partial_{\text{\textgreek{sv}}}\partial_{v}\text{\textgreek{F}}\cdot\text{\textgreek{W}}\square\big(\text{\textgreek{W}}^{-1}\partial_{v}\text{\textgreek{F}}\big)\, dudvd\text{\textgreek{sv}}-\\
 & -\int_{\mathcal{R}(\text{\textgreek{t}}_{1},\text{\textgreek{t}}_{2})\cap\{t\le T^{*}\}}\text{\textgreek{q}}_{R}\cdot r^{p}\partial_{v}^{2}\text{\textgreek{F}}\cdot\text{\textgreek{W}}\square\big(\text{\textgreek{W}}^{-1}\partial_{v}\text{\textgreek{F}}\big)\, dudvd\text{\textgreek{sv}}+C_{p,\text{\textgreek{h}},\text{\textgreek{e}}}Bound_{p,T^{*}}^{\partial_{v}}[\text{\textgreek{f}}](\text{\textgreek{t}}_{1},\text{\textgreek{t}}_{2}),
\end{split}
\label{eq:newMethodFinalFormZeroMassHyperboloids-1-1}
\end{equation}
 where $f(r)=\frac{r^{\text{\textgreek{h}}}}{1+r^{\text{\textgreek{h}}}}$
is the function used in the proof of Lemma \ref{lem:MorawetzDrLemmaHyperboloids}
and 
\begin{align}
Bound_{p,T^{*}}^{\partial_{v}}[\text{\textgreek{f}}](\text{\textgreek{t}}_{1},\text{\textgreek{t}}_{2})= & \sum_{i=1}^{2}\int_{\mathcal{S}_{\text{\textgreek{t}}_{i}}\cap\{t\le T^{*}\}}\text{\textgreek{q}}_{R}\cdot r^{p-2-a}\Big(|\partial_{u}\partial_{v}\text{\textgreek{F}}|^{2}+|r^{-1}\partial_{\text{\textgreek{sv}}}\partial_{v}\text{\textgreek{F}}|^{2}+|\partial_{v}\text{\textgreek{F}}|^{2}+r^{-2}\text{\textgreek{F}}^{2}\Big)\, dvd\text{\textgreek{sv}}+\\
 & +\int_{\mathcal{R}(\text{\textgreek{t}}_{1},\text{\textgreek{t}}_{2})\cap\{t=T^{*}\}}\text{\textgreek{q}}_{R}\cdot r^{p-2-a}\Big(|\partial_{u}\partial_{v}\text{\textgreek{F}}|^{2}+|r^{-1}\partial_{\text{\textgreek{sv}}}\partial_{v}\text{\textgreek{F}}|^{2}+|\partial_{v}\text{\textgreek{F}}|^{2}+r^{-2}\text{\textgreek{F}}^{2}\Big)\, dvd\text{\textgreek{sv}}+\nonumber \\
 & +\int_{\mathcal{S}_{\text{\textgreek{t}}_{1}}\cap\{t\le T^{*}\}}\text{\textgreek{q}}_{R}J_{\text{\textgreek{m}}}^{T}(\text{\textgreek{W}}^{-1}\partial_{v}\text{\textgreek{F}})\bar{n}^{\text{\textgreek{m}}}.\nonumber 
\end{align}
 Notice that since $2<p\le4$, the left hand side of (\ref{eq:newMethodFinalFormZeroMassHyperboloids-1-1})
is not positive definite, since the $\mathcal{E}_{bulk,R,T^{*}}^{(p-1)}$
term contains a summand of the form 
\begin{equation}
\int_{\mathcal{R}(\text{\textgreek{t}}_{1},\text{\textgreek{t}}_{2})\cap\{t\le T^{*}\}}(2-p)\Big(r^{p-1}|r^{-1}\partial_{\text{\textgreek{sv}}}\partial_{v}\text{\textgreek{F}}|^{2}+\frac{(d-1)(d-3)}{4}r^{p-3}|\partial_{v}\text{\textgreek{F}}|^{2}\Big)\, dudvd\text{\textgreek{sv}}
\end{equation}
 which has a negative sign.

We will show that the last term of the right hand side of (\ref{eq:newMethodFinalFormZeroMassHyperboloids-1-1})
can provide us with extra control over bulk terms of the form $\int r^{p-1}|r^{-1}\partial_{\text{\textgreek{sv}}}\partial_{v}\text{\textgreek{F}}|^{2}$.
These terms will then be moved to the left hand side, rendering it
positive definite for $2<p\le4$. 

According to Lemma \ref{lem:Commutator expressions}, $\text{\textgreek{W}}^{-1}\partial_{v}\text{\textgreek{F}}$
satisfies the following equation: 
\begin{align}
\text{\textgreek{W}}\square\big(\text{\textgreek{W}}^{-1}\partial_{v}\text{\textgreek{F}}\big)= & \partial_{v}(\text{\textgreek{W}}F)+2r^{-3}\text{\textgreek{D}}_{g_{\mathbb{S}^{d-1}}+h_{\mathbb{S}^{d-1}}}\text{\textgreek{F}}-\frac{(d-1)(d-3)}{2}r^{-3}\text{\textgreek{F}}+\label{eq:CommutationDv-2}\\
 & +Err(\partial_{v}\text{\textgreek{F}})+r^{-1}Err(\text{\textgreek{F}})\nonumber 
\end{align}
 Therefore (omitting the $dudvd\text{\textgreek{sv}}$ volume form
for the next few lines):
\begin{equation}
\begin{split}\int_{\mathcal{R}(\text{\textgreek{t}}_{1},\text{\textgreek{t}}_{2})\cap\{t\le T^{*}\}} & \text{\textgreek{q}}_{R}\cdot r^{p}\partial_{v}^{2}\text{\textgreek{F}}\cdot\text{\textgreek{W}}\square\big(\text{\textgreek{W}}^{-1}\partial_{v}\text{\textgreek{F}}\big)=\int_{\mathcal{R}(\text{\textgreek{t}}_{1},\text{\textgreek{t}}_{2})\cap\{t\le T^{*}\}}\text{\textgreek{q}}_{R}\cdot r^{p}\partial_{v}^{2}\text{\textgreek{F}}\cdot\partial_{v}(\text{\textgreek{W}}F)+\\
 & +2\int_{\mathcal{R}(\text{\textgreek{t}}_{1},\text{\textgreek{t}}_{2})\cap\{t\le T^{*}\}}\text{\textgreek{q}}_{R}\cdot r^{p-3}\partial_{v}^{2}\text{\textgreek{F}}\cdot\text{\textgreek{D}}_{g_{\mathbb{S}^{d-1}}+h_{\mathbb{S}^{d-1}}}\text{\textgreek{F}}-\frac{(d-1)(d-3)}{2}\int_{\mathcal{R}(\text{\textgreek{t}}_{1},\text{\textgreek{t}}_{2})\cap\{t\le T^{*}\}}\text{\textgreek{q}}_{R}\cdot r^{p-3}\partial_{v}^{2}\text{\textgreek{F}}\cdot\text{\textgreek{F}}+\\
 & +\int_{\mathcal{R}(\text{\textgreek{t}}_{1},\text{\textgreek{t}}_{2})\cap\{t\le T^{*}\}}\text{\textgreek{q}}_{R}\cdot r^{p}\partial_{v}^{2}\text{\textgreek{F}}\cdot\big(Err(\partial_{v}\text{\textgreek{F}})+r^{-1}Err(\text{\textgreek{F}})\big).
\end{split}
\label{eq:BeforeProducingTheGoodAngularTerms}
\end{equation}

Integrating by parts in $\partial_{\text{\textgreek{sv}}},\partial_{v}$,
we have: (note, again, that the volume form used here is $dudvd\text{\textgreek{sv}}$)
\begin{equation}
\begin{split}\int_{\mathcal{R}(\text{\textgreek{t}}_{1},\text{\textgreek{t}}_{2})\cap\{t\le T^{*}\}}\text{\textgreek{q}}_{R}\cdot r^{p-3} & \partial_{v}^{2}\text{\textgreek{F}}\cdot\text{\textgreek{D}}_{g_{\mathbb{S}^{d-1}}+h_{\mathbb{S}^{d-1}}}\text{\textgreek{F}}=\int_{\mathcal{R}(\text{\textgreek{t}}_{1},\text{\textgreek{t}}_{2})\cap\{t\le T^{*}\}}\text{\textgreek{q}}_{R}\cdot r^{p-1}|r^{-1}\partial_{\text{\textgreek{sv}}}\partial_{v}\text{\textgreek{F}}|^{2}-\\
 & -\frac{(p-3)(p-4)}{2}\int_{\mathcal{R}(\text{\textgreek{t}}_{1},\text{\textgreek{t}}_{2})\cap\{t\le T^{*}\}}\text{\textgreek{q}}_{R}\cdot r^{p-3}|r^{-1}\partial_{\text{\textgreek{sv}}}\text{\textgreek{F}}|^{2}+\\
 & +\int_{\mathcal{R}(\text{\textgreek{t}}_{1},\text{\textgreek{t}}_{2})\cap\{t\le T^{*}\}}|\partial\text{\textgreek{q}}_{R}|\cdot\big(O(r^{p})|\partial^{2}\text{\textgreek{F}}|^{2}+O(r^{p-2})|\partial\text{\textgreek{F}}|^{2}\big)+\\
 & +\int_{\mathcal{R}(\text{\textgreek{t}}_{1},\text{\textgreek{t}}_{2})\cap\{t\le T^{*}\}}\text{\textgreek{q}}_{R}\cdot O(r^{p})\partial_{v}\text{\textgreek{F}}\cdot Err(\text{\textgreek{F}})+O_{p,\text{\textgreek{h}}}(1)\cdot Bound_{p,T^{*},\text{\textgreek{e}}}^{(n)}[\text{\textgreek{f}}](\text{\textgreek{t}}_{1},\text{\textgreek{t}}_{2})
\end{split}
\label{eq:PositivityAngularTerms}
\end{equation}
 and 
\begin{equation}
\begin{split}\int_{\mathcal{R}(\text{\textgreek{t}}_{1},\text{\textgreek{t}}_{2})\cap\{t\le T^{*}\}}\text{\textgreek{q}}_{R}\cdot r^{p-3} & \partial_{v}^{2}\text{\textgreek{F}}\cdot\text{\textgreek{F}}=-\int_{\mathcal{R}(\text{\textgreek{t}}_{1},\text{\textgreek{t}}_{2})\cap\{t\le T^{*}\}}\text{\textgreek{q}}_{R}\cdot r^{p-3}|\partial_{v}\text{\textgreek{F}}|^{2}-\\
 & -\frac{(p-3)(p-4)}{2}\int_{\mathcal{R}(\text{\textgreek{t}}_{1},\text{\textgreek{t}}_{2})\cap\{t\le T^{*}\}}\text{\textgreek{q}}_{R}\cdot r^{p-5}|\text{\textgreek{F}}|^{2}+\\
 & +\int_{\mathcal{R}(\text{\textgreek{t}}_{1},\text{\textgreek{t}}_{2})\cap\{t\le T^{*}\}}|\partial\text{\textgreek{q}}_{R}|\cdot\big(O(r^{p-2})|\partial\text{\textgreek{F}}|^{2}+O(r^{p-4})|\text{\textgreek{F}}|^{2}\big)+\\
 & +\int_{\mathcal{R}(\text{\textgreek{t}}_{1},\text{\textgreek{t}}_{2})\cap\{t\le T^{*}\}}\text{\textgreek{q}}_{R}\cdot O(r^{p})\partial_{v}\text{\textgreek{F}}\cdot Err(\text{\textgreek{F}})+O_{p,\text{\textgreek{h}}}(1)\cdot Bound_{p,T^{*},\text{\textgreek{e}}}^{(n)}[\text{\textgreek{f}}](\text{\textgreek{t}}_{1},\text{\textgreek{t}}_{2}),
\end{split}
\label{eq:PositivityAngularTerms-1}
\end{equation}
 where 
\begin{align}
Bound_{p,T^{*},\text{\textgreek{e}}}^{(n)}[\text{\textgreek{f}}](\text{\textgreek{t}}_{1},\text{\textgreek{t}}_{2})= & \text{\textgreek{e}}\cdot\sum_{i=1}^{2}\Big(\mathcal{E}_{bound,R,T^{*}}^{(p)}[\text{\textgreek{W}}^{-1}\partial_{v}\text{\textgreek{F}}](\text{\textgreek{t}}_{i})+\mathcal{E}_{bound,R,T^{*}}^{(p)}[\text{\textgreek{W}}^{-1}r^{-1}\partial_{\text{\textgreek{sv}}}\text{\textgreek{F}}](\text{\textgreek{t}}_{i})+\mathcal{E}_{bound,R,T^{*}}^{(p-2)}[\text{\textgreek{W}}^{-1}\partial_{u}\text{\textgreek{F}}](\text{\textgreek{t}}_{i})\Big)+\label{eq:BoundaryTermsIntegrationByParts}\\
 & +\text{\textgreek{e}}^{-1}\sum_{i=1}^{2}\mathcal{E}_{bound,R,T^{*}}^{(p-2)}[\text{\textgreek{f}}](\text{\textgreek{t}}_{i}).\nonumber 
\end{align}

Therefore, using (\ref{eq:CommutationDv-2}), (\ref{eq:BeforeProducingTheGoodAngularTerms}),
(\ref{eq:PositivityAngularTerms}) and (\ref{eq:PositivityAngularTerms-1}),
as well as an integration by parts scheme similar to the one implemented
in the proof of Lemma \ref{lem:NewMethodGeneralCase}, we can estimate:
\begin{equation}
\begin{split}C_{p,\text{\textgreek{h}},\text{\textgreek{e}}}\int_{\mathcal{R}(\text{\textgreek{t}}_{1},\text{\textgreek{t}}_{2})\cap\{t\le T^{*}\}}\text{\textgreek{q}}_{R} & \cdot\big(\frac{1}{2}f(r)(\partial_{v}-\partial_{u})(\text{\textgreek{W}}^{-1}\partial_{v}\text{\textgreek{F}})+\frac{(d-1)f(r)}{r}(\text{\textgreek{W}}^{-1}\partial_{v}\text{\textgreek{F}})\big)\cdot\square(\text{\textgreek{W}}^{-1}\partial_{v}\text{\textgreek{F}})\,\text{\textgreek{W}}^{2}dudvd\text{\textgreek{sv}}+\\
+C_{p,\text{\textgreek{h}},\text{\textgreek{e}}}\int_{\mathcal{R}(\text{\textgreek{t}}_{1},\text{\textgreek{t}}_{2})\cap\{t\le T^{*}\}} & \text{\textgreek{q}}_{R}\cdot(\partial_{v}+\partial_{u})(\partial_{v}\text{\textgreek{F}})\cdot\text{\textgreek{W}}\square(\text{\textgreek{W}}^{-1}\partial_{v}\text{\textgreek{F}})\, dudvd\text{\textgreek{sv}}+\\
+\int_{\mathcal{R}(\text{\textgreek{t}}_{1},\text{\textgreek{t}}_{2})\cap\{t\le T^{*}\}} & \text{\textgreek{q}}_{R}\cdot O(r^{p-2-a})\cdot\partial_{\text{\textgreek{sv}}}\partial_{v}\text{\textgreek{F}}\cdot\text{\textgreek{W}}\square\big(\text{\textgreek{W}}^{-1}\partial_{v}\text{\textgreek{F}}\big)\, dudvd\text{\textgreek{sv}}-\\
-\int_{\mathcal{R}(\text{\textgreek{t}}_{1},\text{\textgreek{t}}_{2})\cap\{t\le T^{*}\}} & \text{\textgreek{q}}_{R}\cdot r^{p}\partial_{v}^{2}\text{\textgreek{F}}\cdot\text{\textgreek{W}}\square\big(\text{\textgreek{W}}^{-1}\partial_{v}\text{\textgreek{F}}\big)\, dudvd\text{\textgreek{sv}}\le\\
\le & -2\int_{\mathcal{R}(\text{\textgreek{t}}_{1},\text{\textgreek{t}}_{2})\cap\{t\le T^{*}\}}\text{\textgreek{q}}_{R}\cdot r^{p-1}|r^{-1}\partial_{\text{\textgreek{sv}}}\partial_{v}\text{\textgreek{F}}|^{2}\, dudvd\text{\textgreek{sv}}-\frac{(d-1)(d-3)}{2}\int_{\mathcal{R}(\text{\textgreek{t}}_{1},\text{\textgreek{t}}_{2})\cap\{t\le T^{*}\}}\text{\textgreek{q}}_{R}\cdot r^{p-3}|\partial_{v}\text{\textgreek{F}}|^{2}+\\
 & +(p-3)(p-4)\int_{\mathcal{R}(\text{\textgreek{t}}_{1},\text{\textgreek{t}}_{2})\cap\{t\le T^{*}\}}\text{\textgreek{q}}_{R}\cdot\Big(r^{p-3}|r^{-1}\partial_{\text{\textgreek{sv}}}\text{\textgreek{F}}|^{2}+\frac{(d-1)(d-3)}{4}r^{p-5}|\text{\textgreek{F}}|^{2}\Big)\, dudvd\text{\textgreek{sv}}+\\
 & +C_{p,\text{\textgreek{h}},\text{\textgreek{e}}}\int_{\mathcal{R}(\text{\textgreek{t}}_{1},\text{\textgreek{t}}_{2})\cap\{t\le T^{*}\}}|\partial\text{\textgreek{q}}_{R}|\cdot\big(r^{p}|\partial^{2}\text{\textgreek{F}}|^{2}+r^{p-2}|\partial\text{\textgreek{F}}|^{2}\big)\, dudvd\text{\textgreek{sv}}+\\
 & +C_{p,\text{\textgreek{h}},\text{\textgreek{e}}}\int_{\mathcal{R}(\text{\textgreek{t}}_{1},\text{\textgreek{t}}_{2})\cap\{t\le T^{*}\}}\text{\textgreek{q}}_{R}\cdot\big(r^{1+\text{\textgreek{h}}}+r^{p+1}\big)\cdot|\partial_{v}(\text{\textgreek{W}}F)|^{2}\, dudvd\text{\textgreek{sv}}+\\
 & +C_{p,\text{\textgreek{h}},\text{\textgreek{e}}}Bound_{p,T^{*}}^{\partial_{v}}[\text{\textgreek{f}}](\text{\textgreek{t}}_{1},\text{\textgreek{t}}_{2})+C_{p,\text{\textgreek{h}}}Bound_{p,T^{*},\text{\textgreek{e}}}^{(n)}[\text{\textgreek{f}}](\text{\textgreek{t}}_{1},\text{\textgreek{t}}_{2})+\mathcal{B}_{p,T^{*},\text{\textgreek{h}},\text{\textgreek{e}}}^{(\partial_{v})}[\text{\textgreek{f}}](\text{\textgreek{t}}_{1},\text{\textgreek{t}}_{2}),
\end{split}
\label{eq:BoundNonHomogeneousTerms}
\end{equation}
 where we have set: 
\begin{align}
\mathcal{B}_{p,T^{*},\text{\textgreek{h}},\text{\textgreek{e}}}^{(\partial_{v})}[\text{\textgreek{f}}](\text{\textgreek{t}}_{1},\text{\textgreek{t}}_{2})= & C_{p,\text{\textgreek{h}},\text{\textgreek{e}}}\int_{\mathcal{R}(\text{\textgreek{t}}_{1},\text{\textgreek{t}}_{2})\cap\{t\le T^{*}\}}\text{\textgreek{q}}_{R}\cdot\big(r^{1+\text{\textgreek{h}}}+r^{p-1-a}\big)\cdot r^{-2}|r^{-2}\partial_{\text{\textgreek{sv}}}\partial_{\text{\textgreek{sv}}}\text{\textgreek{F}}|^{2}\, dudvd\text{\textgreek{sv}}+\label{eq:BoundRightHandSideDv}\\
 & +C_{p,\text{\textgreek{h}},\text{\textgreek{e}}}\int_{\mathcal{R}(\text{\textgreek{t}}_{1},\text{\textgreek{t}}_{2})\cap\{t\le T^{*}\}}\text{\textgreek{q}}_{R}\cdot r^{p-1-a}|r^{-1}\partial_{\text{\textgreek{sv}}}\partial_{v}\text{\textgreek{F}}|^{2}\, dudvd\text{\textgreek{sv}}+\nonumber \\
 & +C_{p,\text{\textgreek{h}},\text{\textgreek{e}}}\int_{\mathcal{R}(\text{\textgreek{t}}_{1},\text{\textgreek{t}}_{2})\cap\{t\le T^{*}\}}\text{\textgreek{q}}_{R}\cdot r^{p-3-a}|\partial_{u}\partial_{v}\text{\textgreek{F}}|^{2}\, dudvd\text{\textgreek{sv}}+\nonumber \\
 & +C_{p,\text{\textgreek{h}},\text{\textgreek{e}}}\int_{\mathcal{R}(\text{\textgreek{t}}_{1},\text{\textgreek{t}}_{2})\cap\{t\le T^{*}\}}\text{\textgreek{q}}_{R}\cdot\big(r^{1+\text{\textgreek{h}}}+r^{p-1-a}\big)\cdot r^{-2}|Err(\text{\textgreek{F}})|^{2}\, dudvd\text{\textgreek{sv}}+\nonumber \\
 & +\int_{\mathcal{R}(\text{\textgreek{t}}_{1},\text{\textgreek{t}}_{2})\cap\{t\le T^{*}\}}\text{\textgreek{q}}_{R}\cdot O_{p,\text{\textgreek{h}},\text{\textgreek{e}}}(r^{p})\partial_{v}^{2}\text{\textgreek{F}}\cdot Err(\partial_{v}\text{\textgreek{F}})\, dudvd\text{\textgreek{sv}}+\nonumber \\
 & +\int_{\mathcal{R}(\text{\textgreek{t}}_{1},\text{\textgreek{t}}_{2})\cap\{t\le T^{*}\}}\text{\textgreek{q}}_{R}\cdot O_{p,\text{\textgreek{h}},\text{\textgreek{e}}}(r^{p-2-a})\partial_{v}\partial_{\text{\textgreek{sv}}}\text{\textgreek{F}}\cdot Err(\text{\textgreek{F}})\, dudvd\text{\textgreek{sv}}+\nonumber \\
 & +\int_{\mathcal{R}(\text{\textgreek{t}}_{1},\text{\textgreek{t}}_{2})\cap\{t\le T^{*}\}}\text{\textgreek{q}}_{R}\cdot O_{p,\text{\textgreek{h}},\text{\textgreek{e}}}(f(r))\partial_{u}\partial_{v}\text{\textgreek{F}}\cdot Err(\partial_{v}\text{\textgreek{F}})\, dudvd\text{\textgreek{sv}}+\nonumber \\
 & +O_{p,\text{\textgreek{h}}}(\text{\textgreek{e}})\int_{\mathcal{R}(\text{\textgreek{t}}_{1},\text{\textgreek{t}}_{2})\cap\{t\le T^{*}\}}\text{\textgreek{q}}_{R}\cdot r^{p-1}\big(|\partial_{v}^{2}\text{\textgreek{F}}|^{2}+r^{-2}|\partial_{v}\text{\textgreek{F}}|^{2}\big)\, dudvd\text{\textgreek{sv}}.\nonumber 
\end{align}

Returning to (\ref{eq:newMethodFinalFormZeroMassHyperboloids-1-1}),
(\ref{eq:BoundNonHomogeneousTerms}) implies that 
\begin{equation}
\begin{split}\mathcal{E}_{bound,R,T^{*}}^{(p)}[ & \text{\textgreek{W}}^{-1}\partial_{v}\text{\textgreek{F}}](\text{\textgreek{t}}_{2})+\int_{\text{\textgreek{t}}_{1}}^{\text{\textgreek{t}}_{2}}\mathcal{E}_{bulk,R,\text{\textgreek{h}},T^{*}}^{(p-1)}[\text{\textgreek{W}}^{-1}\partial_{v}\text{\textgreek{F}}](\text{\textgreek{t}})\, d\text{\textgreek{t}}+\mathcal{E}_{\mathcal{I}^{+},R,T^{*}}^{(p)}[\text{\textgreek{W}}^{-1}\partial_{v}\text{\textgreek{F}}](\text{\textgreek{t}}_{1},\text{\textgreek{t}}_{2})\le\\
\le & \big(1+O_{p,\text{\textgreek{h}}}(\text{\textgreek{e}})\big)\mathcal{E}_{bound,R,T^{*}}^{(p)}[\text{\textgreek{W}}^{-1}\partial_{v}\text{\textgreek{F}}](\text{\textgreek{t}}_{1})+C_{p,\text{\textgreek{h}},\text{\textgreek{e}}}\int_{\mathcal{R}(\text{\textgreek{t}}_{1},\text{\textgreek{t}}_{2})\cap\{t\le T^{*}\}}|\partial\text{\textgreek{q}}_{R}|\cdot\big(r^{p}|\partial\partial_{v}\text{\textgreek{F}}|^{2}+r^{p-2}|\partial_{v}\text{\textgreek{F}}|^{2}\big)\, dudvd\text{\textgreek{sv}}+\\
 & +C_{p,\text{\textgreek{h}},\text{\textgreek{e}}}\int_{\mathcal{R}(\text{\textgreek{t}}_{1},\text{\textgreek{t}}_{2})\cap\{t\le T^{*}\}}\text{\textgreek{q}}_{R}\cdot\max\{r^{p-3-a},r^{-3}\}|\partial_{v}\text{\textgreek{F}}|^{2}\, dudvd\text{\textgreek{sv}}+\\
 & +C_{p,\text{\textgreek{h}},\text{\textgreek{e}}}\int_{\mathcal{R}(\text{\textgreek{t}}_{1},\text{\textgreek{t}}_{2})\cap\{t\le T^{*}\}}\text{\textgreek{q}}_{R}\cdot\big(r^{1+\text{\textgreek{h}}}+r^{p+1}\big)\cdot|\partial_{v}(\text{\textgreek{W}}F)|^{2}\, dudvd\text{\textgreek{sv}}+\\
 & +(p-3)(p-4)\int_{\mathcal{R}(\text{\textgreek{t}}_{1},\text{\textgreek{t}}_{2})\cap\{t\le T^{*}\}}\text{\textgreek{q}}_{R}\cdot\Big(r^{p-3}|r^{-1}\partial_{\text{\textgreek{sv}}}\text{\textgreek{F}}|^{2}+\frac{(d-1)(d-3)}{4}r^{p-5}|\text{\textgreek{F}}|^{2}\Big)\, dudvd\text{\textgreek{sv}}+\\
 & +C_{p,\text{\textgreek{h}},\text{\textgreek{e}}}Bound_{p,T^{*}}^{\partial_{v}}[\text{\textgreek{f}}](\text{\textgreek{t}}_{1},\text{\textgreek{t}}_{2})+C_{p,\text{\textgreek{h}}}Bound_{p,T^{*},\text{\textgreek{e}}}^{(n)}[\text{\textgreek{f}}](\text{\textgreek{t}}_{1},\text{\textgreek{t}}_{2})+\mathcal{B}_{p,T^{*},\text{\textgreek{h}},\text{\textgreek{e}}}^{(\partial_{v})}[\text{\textgreek{f}}](\text{\textgreek{t}}_{1},\text{\textgreek{t}}_{2})-\\
 & -2\int_{\mathcal{R}(\text{\textgreek{t}}_{1},\text{\textgreek{t}}_{2})\cap\{t\le T^{*}\}}\text{\textgreek{q}}_{R}\cdot r^{p-1}|r^{-1}\partial_{\text{\textgreek{sv}}}\partial_{v}\text{\textgreek{F}}|^{2}\, dudvd\text{\textgreek{sv}}-\frac{(d-1)(d-3)}{2}\int_{\mathcal{R}(\text{\textgreek{t}}_{1},\text{\textgreek{t}}_{2})\cap\{t\le T^{*}\}}\text{\textgreek{q}}_{R}\cdot r^{p-1}|\partial_{v}\text{\textgreek{F}}|^{2}.
\end{split}
\label{eq:newMethodHyperboloidsAlmostDv}
\end{equation}
 Notice that the last two terms of the right hand side of (\ref{eq:newMethodHyperboloidsAlmostDv})
can be moved to the left hand side, thus providing us with extra control
over bulk terms of the form $\int\text{\textgreek{q}}_{R}\cdot r^{p-1}|r^{-1}\partial_{\text{\textgreek{sv}}}\partial_{v}\text{\textgreek{F}}|^{2}$
and (in dimensions $d>3$) $\int\text{\textgreek{q}}_{R}\cdot r^{p-3}|\partial_{v}\text{\textgreek{F}}|^{2}$:
\begin{equation}
\begin{split}\mathcal{E}_{bound,R,T^{*}}^{(p)}[ & \text{\textgreek{W}}^{-1}\partial_{v}\text{\textgreek{F}}](\text{\textgreek{t}}_{2})+\int_{\text{\textgreek{t}}_{1}}^{\text{\textgreek{t}}_{2}}\mathcal{E}{}_{bulk,R,\text{\textgreek{h}},T^{*}}^{(p-1,\partial_{v})}[\text{\textgreek{W}}^{-1}\partial_{v}\text{\textgreek{F}}](\text{\textgreek{t}})\, d\text{\textgreek{t}}+\mathcal{E}_{\mathcal{I}^{+},R,T^{*}}^{(p)}[\text{\textgreek{W}}^{-1}\partial_{v}\text{\textgreek{F}}](\text{\textgreek{t}}_{1},\text{\textgreek{t}}_{2})\le\\
\le & \big(1+O_{p,\text{\textgreek{h}}}(\text{\textgreek{e}})\big)\mathcal{E}_{bound,R,T^{*}}^{(p)}[\text{\textgreek{W}}^{-1}\partial_{v}\text{\textgreek{F}}](\text{\textgreek{t}}_{1})+C_{p,\text{\textgreek{h}},\text{\textgreek{e}}}\int_{\mathcal{R}(\text{\textgreek{t}}_{1},\text{\textgreek{t}}_{2})\cap\{t\le T^{*}\}}|\partial\text{\textgreek{q}}_{R}|\cdot\big(r^{p}|\partial\partial_{v}\text{\textgreek{F}}|^{2}+r^{p-2}|\partial_{v}\text{\textgreek{F}}|^{2}\big)\, dudvd\text{\textgreek{sv}}+\\
 & +C_{p,\text{\textgreek{h}},\text{\textgreek{e}}}\int_{\mathcal{R}(\text{\textgreek{t}}_{1},\text{\textgreek{t}}_{2})\cap\{t\le T^{*}\}}\text{\textgreek{q}}_{R}\cdot\max\{r^{p-3-a},r^{-3}\}|\partial_{v}\text{\textgreek{F}}|^{2}\, dudvd\text{\textgreek{sv}}+\\
 & +C_{p,\text{\textgreek{h}},\text{\textgreek{e}}}\int_{\{\text{\textgreek{t}}_{1}\le\bar{t}\le\text{\textgreek{t}}_{2}\}\cap\{t\le T^{*}\}}\text{\textgreek{q}}_{R}\cdot\big(r^{1+\text{\textgreek{h}}}+r^{p+1}\big)\cdot|\partial_{v}(\text{\textgreek{W}}F)|^{2}\, dudvd\text{\textgreek{sv}}+\\
 & +(p-3)(p-4)\int_{\mathcal{R}(\text{\textgreek{t}}_{1},\text{\textgreek{t}}_{2})\cap\{t\le T^{*}\}}\text{\textgreek{q}}_{R}\cdot\Big(r^{p-3}|r^{-1}\partial_{\text{\textgreek{sv}}}\text{\textgreek{F}}|^{2}+\frac{(d-1)(d-3)}{4}r^{p-5}|\text{\textgreek{F}}|^{2}\Big)\, dudvd\text{\textgreek{sv}}+\\
 & +C_{p,\text{\textgreek{h}},\text{\textgreek{e}}}Bound_{p,T^{*}}^{\partial_{v}}[\text{\textgreek{f}}](\text{\textgreek{t}}_{1},\text{\textgreek{t}}_{2})+C_{p,\text{\textgreek{h}}}Bound_{p,T^{*},\text{\textgreek{e}}}^{(n)}[\text{\textgreek{f}}](\text{\textgreek{t}}_{1},\text{\textgreek{t}}_{2})+\mathcal{B}_{p,T^{*},\text{\textgreek{h}},\text{\textgreek{e}}}^{(\partial_{v})}[\text{\textgreek{f}}](\text{\textgreek{t}}_{1},\text{\textgreek{t}}_{2}),
\end{split}
\label{eq:newMethodHyperboloidsDvImprovedFinal}
\end{equation}
 where we have set 
\begin{align}
\mathcal{E}{}_{bulk,R,\text{\textgreek{h}},T^{*}}^{(p-1,\partial_{v})}[\text{\textgreek{y}}](\text{\textgreek{t}})\doteq & \int_{\mathcal{S}_{\text{\textgreek{t}}}\cap\{t\le T^{*}\}}\text{\textgreek{q}}_{R}\Big(pr^{p-1}\big|\partial_{v}(\text{\textgreek{W}}\text{\textgreek{y}})\big|^{2}+(6-p)\big(r^{p-1}\big|r^{-1}\partial_{\text{\textgreek{sv}}}(\text{\textgreek{W}}\text{\textgreek{y}})\big|^{2}+\frac{(d-1)(d-3)}{4}r^{p-3}\big|\text{\textgreek{W}}\text{\textgreek{y}}\big|^{2}\big)\Big)\, dvd\text{\textgreek{sv}}+\label{eq:ImprovedAngularBulk}\\
 & +\int_{\mathcal{S}_{\text{\textgreek{t}}}\cap\{t\le T^{*}\}}\text{\textgreek{q}}_{R}r^{-1-\text{\textgreek{h}}}|\partial_{u}(\text{\textgreek{W}}\text{\textgreek{y}})|^{2}\, dvd\text{\textgreek{sv}}.\nonumber 
\end{align}
 Notice the $(6-p)$ factor in the $\mathcal{E}{}_{bulk,R,\text{\textgreek{h}},T^{*}}^{(p-1,\partial_{v})}$
bulk energy norm which makes it positive definite for $2<p\le4$.

We can extract a similar inequality for $r^{-1}\partial_{\text{\textgreek{sv}}}\text{\textgreek{F}}$
in place of $\partial_{v}\text{\textgreek{F}}$. Repeating as before
the proof of Lemma \ref{lem:NewMethodGeneralCase} for $\text{\textgreek{W}}^{-1}\big(r^{-1}\nabla^{\mathbb{S}^{d-1}}\text{\textgreek{F}}\big)$
in place of $\text{\textgreek{f}}$, we obtain:
\begin{equation}
\begin{split}\mathcal{E}_{bound,R,T^{*}}^{(p)}[ & \text{\textgreek{W}}^{-1}r^{-1}\partial_{\text{\textgreek{sv}}}\text{\textgreek{F}}](\text{\textgreek{t}}_{2})+\int_{\text{\textgreek{t}}_{1}}^{\text{\textgreek{t}}_{2}}\mathcal{E}_{bulk,R,\text{\textgreek{h}},T^{*}}^{(p-1)}[\text{\textgreek{W}}^{-1}r^{-1}\partial_{\text{\textgreek{sv}}}\text{\textgreek{F}}](\text{\textgreek{t}})\, d\text{\textgreek{t}}+\mathcal{E}_{\mathcal{I}^{+},R,T^{*}}^{(p)}[\text{\textgreek{W}}^{-1}r^{-1}\partial_{\text{\textgreek{sv}}}\text{\textgreek{F}}](\text{\textgreek{t}}_{1},\text{\textgreek{t}}_{2})\le\\
\le & \big(1+O_{p,\text{\textgreek{h}}}(\text{\textgreek{e}})\big)\mathcal{E}_{bound,R,T^{*}}^{(p)}[\text{\textgreek{W}}^{-1}r^{-1}\partial_{\text{\textgreek{sv}}}\text{\textgreek{F}}](\text{\textgreek{t}}_{1})+C_{p,\text{\textgreek{h}},\text{\textgreek{e}}}\int_{\mathcal{R}(\text{\textgreek{t}}_{1},\text{\textgreek{t}}_{2})\cap\{t\le T^{*}\}}|\partial\text{\textgreek{q}}_{R}|\cdot\big(r^{p-2}|\partial\partial_{\text{\textgreek{sv}}}\text{\textgreek{F}}|^{2}+r^{p-4}|\partial_{\text{\textgreek{sv}}}\text{\textgreek{F}}|^{2}\big)\, dudvd\text{\textgreek{sv}}+\\
 & +C_{p,\text{\textgreek{h}},\text{\textgreek{e}}}\int_{\mathcal{R}(\text{\textgreek{t}}_{1},\text{\textgreek{t}}_{2})\cap\{t\le T^{*}\}}\text{\textgreek{q}}_{R}\cdot\Big(r^{p-1-a}|r^{-2}\partial_{\text{\textgreek{sv}}}\partial_{\text{\textgreek{sv}}}\text{\textgreek{F}}|^{2}+r^{p-3-a}|r^{-1}\partial_{u}\partial_{\text{\textgreek{sv}}}\text{\textgreek{F}}|^{2}+\max\{r^{p-3-a},r^{-3}\}|r^{-1}\partial_{\text{\textgreek{sv}}}\text{\textgreek{F}}|^{2}\Big)\, dudvd\text{\textgreek{sv}}+\\
 & +C_{p,\text{\textgreek{h}},\text{\textgreek{e}}}\int_{\mathcal{R}(\text{\textgreek{t}}_{1},\text{\textgreek{t}}_{2})\cap\{t\le T^{*}\}}\text{\textgreek{q}}_{R}\cdot\frac{1}{2}f\cdot(\partial_{v}-\partial_{u})(\text{\textgreek{W}}^{-1}r^{-1}\nabla^{\mathbb{S}^{d-1}}\text{\textgreek{F}})\cdot\square(\text{\textgreek{W}}^{-1}r^{-1}\nabla^{\mathbb{S}^{d-1}}\text{\textgreek{F}})\,\text{\textgreek{W}}^{2}dudvd\text{\textgreek{sv}}+\\
 & +C_{p,\text{\textgreek{h}},\text{\textgreek{e}}}\int_{\mathcal{R}(\text{\textgreek{t}}_{1},\text{\textgreek{t}}_{2})\cap\{t\le T^{*}\}}\text{\textgreek{q}}_{R}\cdot\frac{(d-1)f}{r}(\text{\textgreek{W}}^{-1}r^{-1}\nabla^{\mathbb{S}^{d-1}}\text{\textgreek{F}})\cdot\square(\text{\textgreek{W}}^{-1}r^{-1}\nabla^{\mathbb{S}^{d-1}}\text{\textgreek{F}})\,\text{\textgreek{W}}^{2}dudvd\text{\textgreek{sv}}+\\
 & +C_{p,\text{\textgreek{h}},\text{\textgreek{e}}}\int_{\mathcal{R}(\text{\textgreek{t}}_{1},\text{\textgreek{t}}_{2})\cap\{t\le T^{*}\}}\text{\textgreek{q}}_{R}\cdot(\partial_{v}+\partial_{u})(r^{-1}\nabla^{\mathbb{S}^{d-1}}\text{\textgreek{F}})\cdot\text{\textgreek{W}}\square(\text{\textgreek{W}}^{-1}r^{-1}\nabla^{\mathbb{S}^{d-1}}\text{\textgreek{F}})\, dudvd\text{\textgreek{sv}}+\\
 & +\int_{\mathcal{R}(\text{\textgreek{t}}_{1},\text{\textgreek{t}}_{2})\cap\{t\le T^{*}\}}\text{\textgreek{q}}_{R}\cdot O(r^{p-3-a})\cdot\partial_{\text{\textgreek{sv}}}^{2}\text{\textgreek{F}}\cdot\text{\textgreek{W}}\square\big(\text{\textgreek{W}}^{-1}r^{-1}\partial_{\text{\textgreek{sv}}}\text{\textgreek{F}}\big)\, dudvd\text{\textgreek{sv}}-\\
 & -\int_{\mathcal{R}(\text{\textgreek{t}}_{1},\text{\textgreek{t}}_{2})\cap\{t\le T^{*}\}}\text{\textgreek{q}}_{R}\cdot r^{p}\partial_{v}(r^{-1}\nabla^{\mathbb{S}^{d-1}}\text{\textgreek{F}})\cdot\text{\textgreek{W}}\square\big(\text{\textgreek{W}}^{-1}r^{-1}\nabla^{\mathbb{S}^{d-1}}\text{\textgreek{F}}\big)\, dudvd\text{\textgreek{sv}}+C_{p,\text{\textgreek{h}},\text{\textgreek{e}}}Bound_{p,T^{*}}^{r^{-1}\partial_{\text{\textgreek{sv}}}}[\text{\textgreek{f}}](\text{\textgreek{t}}_{1},\text{\textgreek{t}}_{2}),
\end{split}
\label{eq:newMethodFinalFormZeroMassHyperboloidsDs-1-1}
\end{equation}
 where the multiplication between derivatives of the $\mathbb{S}^{d-1}$
gradient $\nabla^{\mathbb{S}^{d-1}}\text{\textgreek{F}}$ of $\text{\textgreek{F}}$
in the last lines of (\ref{eq:newMethodFinalFormZeroMassHyperboloidsDs-1-1})
is performed with respect to the usual metric $g_{\mathbb{S}^{d-1}}$
of $\mathbb{S}^{d-1}$, and 
\begin{align}
Bound_{p,T^{*}}^{r^{-1}\partial_{\text{\textgreek{sv}}}}[\text{\textgreek{f}}](\text{\textgreek{t}}_{1},\text{\textgreek{t}}_{2})= & \sum_{i=1}^{2}\int_{\mathcal{S}_{\text{\textgreek{t}}_{i}}\cap\{t\le T^{*}\}}\text{\textgreek{q}}_{R}\cdot r^{p-2-a}\Big(|r^{-1}\partial_{u}\partial_{\text{\textgreek{sv}}}\text{\textgreek{F}}|^{2}+|r^{-2}\partial_{\text{\textgreek{sv}}}\partial_{\text{\textgreek{sv}}}\text{\textgreek{F}}|^{2}+|r^{-1}\partial_{\text{\textgreek{sv}}}\text{\textgreek{F}}|^{2}+r^{-2}\text{\textgreek{F}}^{2}\Big)\, dvd\text{\textgreek{sv}}+\\
 & +\int_{\mathcal{R}(\text{\textgreek{t}}_{1},\text{\textgreek{t}}_{2})\cap\{t=T^{*}\}}\text{\textgreek{q}}_{R}\cdot r^{p-2-a}\Big(|r^{-1}\partial_{u}\partial_{\text{\textgreek{sv}}}\text{\textgreek{F}}|^{2}+|r^{-2}\partial_{\text{\textgreek{sv}}}\partial_{\text{\textgreek{sv}}}\text{\textgreek{F}}|^{2}+|r^{-1}\partial_{\text{\textgreek{sv}}}\text{\textgreek{F}}|^{2}+r^{-2}\text{\textgreek{F}}^{2}\Big)\, dvd\text{\textgreek{sv}}+\nonumber \\
 & +\int_{\mathcal{S}_{\text{\textgreek{t}}_{1}}\cap\{t\le T^{*}\}}\text{\textgreek{q}}_{R}J_{\text{\textgreek{m}}}^{T}(\text{\textgreek{W}}^{-1}r^{-1}\partial_{\text{\textgreek{sv}}}\text{\textgreek{F}})\bar{n}^{\text{\textgreek{m}}}.\nonumber 
\end{align}

Using Lemma \ref{lem:Commutator expressions}, we have: 
\begin{align}
\text{\textgreek{W}}\square\big(\text{\textgreek{W}}^{-1}(r^{-1}\partial_{\text{\textgreek{sv}}}\text{\textgreek{F}})\big)= & r^{-1}\partial_{\text{\textgreek{sv}}}\big(\text{\textgreek{W}}\square\text{\textgreek{f}}\big)+r^{-2}\partial_{u}\partial_{\text{\textgreek{sv}}}\text{\textgreek{F}}-r^{-2}\partial_{v}\partial_{\text{\textgreek{sv}}}\text{\textgreek{F}}+\label{eq:CommutationDsigma-2}\\
 & +Err(r^{-1}\partial_{\text{\textgreek{sv}}}\text{\textgreek{F}})+r^{-1}Err(\text{\textgreek{F}}).\nonumber 
\end{align}

\begin{rem*}
The term $r^{-2}\partial_{u}\partial_{\text{\textgreek{sv}}}\text{\textgreek{F}}$
\footnote{Recall that this term is actually $r^{-2}\mathcal{L}_{\partial_{u}}(\nabla^{\mathbb{S}^{d-1}}\text{\textgreek{F}})$,
see Section \ref{sub:Notational-conventions}%
} in (\ref{eq:CommutationDsigma-2}) will provide us with improved
control over $|\partial_{\text{\textgreek{sv}}}\partial_{\text{\textgreek{sv}}}\text{\textgreek{F}}|^{2}$
bulk terms, in the same way that the term $r^{-3}\text{\textgreek{D}}_{g_{\mathbb{S}^{d-1}}+h_{\mathbb{S}^{d-1}}}\text{\textgreek{F}}$
in (\ref{eq:CommutationDv-2}) eventually provided us with improved
control over $|\partial_{\text{\textgreek{sv}}}\partial_{v}\text{\textgreek{F}}|^{2}$
bulk terms through (\ref{eq:BeforeProducingTheGoodAngularTerms})
and (\ref{eq:PositivityAngularTerms}). We should also notice that
the term $-r^{-2}\partial_{v}\partial_{\text{\textgreek{sv}}}\text{\textgreek{F}}$
has a bad sign, and will result in the appearence of bulk terms of
the form $\int\text{\textgreek{q}}_{R}r^{p-1}|\partial_{\text{\textgreek{sv}}}\partial_{v}\text{\textgreek{F}}|^{2}$
with a non-convenient sign. All the terms with a bad sign, however,
will be controlled by the corresponding terms in the left hand side
of (\ref{eq:newMethodHyperboloidsDvImprovedFinal}) plus a multiple
of (\ref{eq:newMethodFinalStatementHyperboloids}) for $p-2$ in place
of $p$, provided $p\le4$.
\end{rem*}
Integrating by parts in $\partial_{v}$ and the spherical directions
(omitting the $dudvd\text{\textgreek{sv}}$ volume form in the next
few lines), and then once again in $\partial_{u}$ for the error term,
we compute: 
\begin{equation}
\begin{split}\int_{\mathcal{R}(\text{\textgreek{t}}_{1},\text{\textgreek{t}}_{2})\cap\{t\le T^{*}\}}\text{\textgreek{q}}_{R}\cdot(1+ & C\cdot r^{-a})r^{p-2}\partial_{v}(r^{-1}\nabla^{\mathbb{S}^{d-1}}\text{\textgreek{F}}\big)\cdot\partial_{u}\nabla^{\mathbb{S}^{d-1}}\text{\textgreek{F}}=\\
= & \int_{\mathcal{R}(\text{\textgreek{t}}_{1},\text{\textgreek{t}}_{2})\cap\{t\le T^{*}\}}\text{\textgreek{q}}_{R}\cdot(1+O(r^{-a}))r^{p-1}\big(r^{-2}\text{\textgreek{D}}_{g_{\mathbb{S}^{d-1}}}\text{\textgreek{F}}\big)\cdot\partial_{u}\partial_{v}\text{\textgreek{F}}+\\
 & +\int_{\mathcal{R}(\text{\textgreek{t}}_{1},\text{\textgreek{t}}_{2})\cap\{t\le T^{*}\}}\text{\textgreek{q}}_{R}\cdot\Big(\frac{(p-2)(4-p)}{2}r^{p-3}+O(r^{p-3-a})\Big)\big|r^{-1}\partial_{\text{\textgreek{sv}}}\text{\textgreek{F}}\big|^{2}+\\
 & +\int_{\mathcal{R}(\text{\textgreek{t}}_{1},\text{\textgreek{t}}_{2})\cap\{t\le T^{*}\}}|\partial\text{\textgreek{q}}_{R}|\cdot\big(O(r^{p})|\partial^{2}\text{\textgreek{F}}|^{2}+O(r^{p-2})|\partial\text{\textgreek{F}}|^{2}\big)+\\
 & +O_{p,\text{\textgreek{h}}}(1)\cdot Bound_{p,T^{*},\text{\textgreek{e}}}^{(n)}[\text{\textgreek{f}}](\text{\textgreek{t}}_{1},\text{\textgreek{t}}_{2}).
\end{split}
\label{eq:PositivityDsDs}
\end{equation}

However, since $\square_{g}\text{\textgreek{f}}=F$, the following
equality holds: 
\begin{align}
\big(1+O(r^{-1-a})\big)\partial_{u}\partial_{v}\text{\textgreek{F}}= & -\text{\textgreek{W}}F+r^{-2}\text{\textgreek{D}}_{g_{\mathbb{S}^{d-1}}+h_{\mathbb{S}^{d-1}}}\text{\textgreek{F}}-\label{eq:DuDvPhi}\\
 & -\Big(\frac{(d-1)(d-3)}{4}\cdot r^{-2}\Big)\text{\textgreek{F}}+Err(\text{\textgreek{F}}).\nonumber 
\end{align}

Using (\ref{eq:DuDvPhi}) to substitute $\partial_{u}\partial_{v}\text{\textgreek{F}}$
in (\ref{eq:PositivityDsDs}), as well as the following ellliptic-type
estimate on $\mathbb{S}^{d-1}$ 
\begin{equation}
\int_{\mathbb{S}^{d-1}}|\text{\textgreek{D}}_{g_{\mathbb{S}^{d-1}}}\text{\textgreek{F}}|^{2}\, d\text{\textgreek{sv}}\ge\int_{\mathbb{S}^{d-1}}|\partial_{\text{\textgreek{sv}}}\partial_{\text{\textgreek{sv}}}\text{\textgreek{F}}|^{2}\, d\text{\textgreek{sv}}+(d-2)\int_{\mathbb{S}^{d-1}}|\partial_{\text{\textgreek{sv}}}\text{\textgreek{F}}|^{2},
\end{equation}
 we infer after integrating by parts in $\partial_{\text{\textgreek{sv}}}$
in the term $\int\text{\textgreek{q}}_{R}\cdot\frac{(d-1)(d-3)}{2}\cdot r^{p-5}\text{\textgreek{D}}_{g_{\mathbb{S}^{d-1}}}\text{\textgreek{F}}\cdot\text{\textgreek{F}}$:
\begin{equation}
\begin{split}\int_{\mathcal{R}(\text{\textgreek{t}}_{1},\text{\textgreek{t}}_{2})\cap\{t\le T^{*}\}} & \text{\textgreek{q}}_{R}\cdot r^{p-2}\partial_{v}(r^{-1}\nabla^{\mathbb{S}^{d-1}}\text{\textgreek{F}}\big)\cdot\partial_{u}\nabla^{\mathbb{S}^{d-1}}\text{\textgreek{F}}\ge\\
\ge & \int_{\mathcal{R}(\text{\textgreek{t}}_{1},\text{\textgreek{t}}_{2})\cap\{t\le T^{*}\}}\text{\textgreek{q}}_{R}\cdot r^{p-1}|r^{-2}\partial_{\text{\textgreek{sv}}}\partial_{\text{\textgreek{sv}}}\text{\textgreek{F}}|^{2}+\\
 & +\int_{\mathcal{R}(\text{\textgreek{t}}_{1},\text{\textgreek{t}}_{2})\cap\{t\le T^{*}\}}\text{\textgreek{q}}_{R}\cdot r^{p-1-a}\cdot O\big(|r^{-2}\partial_{\text{\textgreek{sv}}}\partial_{\text{\textgreek{sv}}}\text{\textgreek{F}}|^{2}+|r^{-1}\partial_{\text{\textgreek{sv}}}\text{\textgreek{F}}|^{2}\big)+\\
 & +\int_{\mathcal{R}(\text{\textgreek{t}}_{1},\text{\textgreek{t}}_{2})\cap\{t\le T^{*}\}}\text{\textgreek{q}}_{R}\cdot\Big(\frac{(d-1)(d-3)+2(p-2)(4-p)}{4}+(d-2)+O(r^{-a})\Big)\cdot r^{p-3}|r^{-1}\partial_{\text{\textgreek{sv}}}\text{\textgreek{F}}|^{2}+\\
 & +\int_{\mathcal{R}(\text{\textgreek{t}}_{1},\text{\textgreek{t}}_{2})\cap\{t\le T^{*}\}}O(|\partial\text{\textgreek{q}}_{R}|)\cdot\big(O(r^{p})|\partial^{2}\text{\textgreek{F}}|^{2}+O(r^{p-2})|\partial\text{\textgreek{F}}|^{2}\big)+\\
 & +\int_{\mathcal{R}(\text{\textgreek{t}}_{1},\text{\textgreek{t}}_{2})\cap\{t\le T^{*}\}}\text{\textgreek{q}}_{R}\cdot O(r^{p-1})\big(r^{-2}\partial_{\text{\textgreek{sv}}}\partial_{\text{\textgreek{sv}}}\text{\textgreek{F}}\big)\cdot Err(\text{\textgreek{F}})+\\
 & +\int_{\mathcal{R}(\text{\textgreek{t}}_{1},\text{\textgreek{t}}_{2})\cap\{t\le T^{*}\}}\text{\textgreek{q}}_{R}\cdot O(r^{p-1})\cdot\partial_{v}\big(r^{-1}\partial_{\text{\textgreek{sv}}}\text{\textgreek{F}}\big)\cdot\text{\textgreek{W}}F+\\
 & +C_{p,\text{\textgreek{h}},\text{\textgreek{e}}}Bound_{p,T^{*}}^{r^{-1}\partial_{\text{\textgreek{sv}}}}[\text{\textgreek{f}}](\text{\textgreek{t}}_{1},\text{\textgreek{t}}_{2})+C_{p,\text{\textgreek{h}}}Bound_{p,T^{*},\text{\textgreek{e}}}^{(n)}[\text{\textgreek{f}}](\text{\textgreek{t}}_{1},\text{\textgreek{t}}_{2}).
\end{split}
\label{eq:PositivityDsDs-1}
\end{equation}
 This should be considered as the analogue of (\ref{eq:PositivityAngularTerms}).
Proceeding therefore as before, from (\ref{eq:newMethodFinalFormZeroMassHyperboloidsDs-1-1}),
(\ref{eq:CommutationDsigma-2}) and (\ref{eq:PositivityDsDs-1}) we
can extract the following analogue of (\ref{eq:newMethodHyperboloidsDvImprovedFinal}):
\begin{equation}
\begin{split}\mathcal{E}_{bound,R,T^{*}}^{(p)}[ & \text{\textgreek{W}}^{-1}r^{-1}\partial_{\text{\textgreek{sv}}}\text{\textgreek{F}}](\text{\textgreek{t}}_{2})+\int_{\text{\textgreek{t}}_{1}}^{\text{\textgreek{t}}_{2}}\mathcal{E}{}_{bulk,R,\text{\textgreek{h}},T^{*}}^{(p-1,\partial_{\text{\textgreek{sv}}})}[\text{\textgreek{W}}^{-1}r^{-1}\partial_{\text{\textgreek{sv}}}\text{\textgreek{F}}](\text{\textgreek{t}})\, d\text{\textgreek{t}}+\mathcal{E}_{\mathcal{I}^{+},R,T^{*}}^{(p)}[\text{\textgreek{W}}^{-1}r^{-1}\partial_{\text{\textgreek{sv}}}\text{\textgreek{F}}](\text{\textgreek{t}}_{1},\text{\textgreek{t}}_{2})\le\\
\le & \big(1+O_{p,\text{\textgreek{h}}}(\text{\textgreek{e}})\big)\mathcal{E}_{bound,R,T^{*}}^{(p)}[\text{\textgreek{W}}^{-1}r^{-1}\partial_{\text{\textgreek{sv}}}\text{\textgreek{F}}](\text{\textgreek{t}}_{1})+C_{p,\text{\textgreek{h}},\text{\textgreek{e}}}\int_{\mathcal{R}(\text{\textgreek{t}}_{1},\text{\textgreek{t}}_{2})\cap\{t\le T^{*}\}}|\partial\text{\textgreek{q}}_{R}|\cdot\big(r^{p-2}|\partial\partial_{\text{\textgreek{sv}}}\text{\textgreek{F}}|^{2}+r^{p-4}|\partial_{\text{\textgreek{sv}}}\text{\textgreek{F}}|^{2}\big)\, dudvd\text{\textgreek{sv}}+\\
 & +C_{p,\text{\textgreek{h}},\text{\textgreek{e}}}\int_{\mathcal{R}(\text{\textgreek{t}}_{1},\text{\textgreek{t}}_{2})\cap\{t\le T^{*}\}}\text{\textgreek{q}}_{R}\cdot\max\{r^{p-3-a},r^{-3}\}|r^{-1}\partial_{\text{\textgreek{sv}}}\text{\textgreek{F}}|^{2}\, dudvd\text{\textgreek{sv}}+\\
 & +C_{p,\text{\textgreek{h}},\text{\textgreek{e}}}\int_{\mathcal{R}(\text{\textgreek{t}}_{1},\text{\textgreek{t}}_{2})\cap\{t\le T^{*}\}}\text{\textgreek{q}}_{R}\cdot\big(r^{1+\text{\textgreek{h}}}+r^{p+1}\big)\cdot\big(|r^{-1}\partial_{\text{\textgreek{sv}}}(\text{\textgreek{W}}F)|^{2}+r^{-2}|\text{\textgreek{W}}F|^{2}\big)\, dudvd\text{\textgreek{sv}}+\\
 & +C_{p,\text{\textgreek{h}},\text{\textgreek{e}}}\cdot\int_{\mathcal{R}(\text{\textgreek{t}}_{1},\text{\textgreek{t}}_{2})\cap\{t\le T^{*}\}}\text{\textgreek{q}}_{R}\cdot r^{p-1-a}|r^{-1}\partial_{\text{\textgreek{sv}}}\partial_{v}\text{\textgreek{F}}|^{2}\, dudvd\text{\textgreek{sv}}+\\
 & +\int_{\mathcal{R}(\text{\textgreek{t}}_{1},\text{\textgreek{t}}_{2})\cap\{t\le T^{*}\}}\text{\textgreek{q}}_{R}\cdot(1+O(r^{-a}))r^{p-1}\big(r^{-1}\partial_{v}\nabla^{\mathbb{S}^{d-1}}\text{\textgreek{F}}\big)\cdot\big(\partial_{v}(r^{-1}\nabla^{\mathbb{S}^{d-1}}\text{\textgreek{F}})\big)\, dudvd\text{\textgreek{sv}}+\\
 & +C_{p,\text{\textgreek{h}},\text{\textgreek{e}}}Bound_{p,T^{*}}^{r^{-1}\partial_{\text{\textgreek{sv}}}}[\text{\textgreek{f}}](\text{\textgreek{t}}_{1},\text{\textgreek{t}}_{2})+C_{p,\text{\textgreek{h}}}Bound_{p,T^{*},\text{\textgreek{e}}}^{(n)}[\text{\textgreek{f}}](\text{\textgreek{t}}_{1},\text{\textgreek{t}}_{2})+\mathcal{B}_{p,T^{*},\text{\textgreek{h}},\text{\textgreek{e}}}^{(r^{-1}\partial_{\text{\textgreek{sv}}})}[\text{\textgreek{f}}](\text{\textgreek{t}}_{1},\text{\textgreek{t}}_{2}),
\end{split}
\label{eq:newMethodHyperboloidsDsImprovedFinal}
\end{equation}
 where we have set 
\begin{align}
\mathcal{E}{}_{bulk,R,\text{\textgreek{h},}T^{*}}^{(p-1,\partial_{\text{\textgreek{sv}}})}[\text{\textgreek{y}}](\text{\textgreek{t}})\doteq & \int_{\mathcal{S}_{\text{\textgreek{t}}}\cap\{t\le T^{*}\}}\text{\textgreek{q}}_{R}\Big(pr^{p-1}\big|\partial_{v}(\text{\textgreek{W}}\text{\textgreek{y}})\big|^{2}+\max\big\{(4-p),r^{-\frac{a}{2}}\big\}\cdot\big(r^{p-1}\big|r^{-1}\partial_{\text{\textgreek{sv}}}(\text{\textgreek{W}}\text{\textgreek{y}})\big|^{2}+\frac{(d-1)(d-3)}{4}r^{p-3}\big|\text{\textgreek{W}}\text{\textgreek{y}}\big|^{2}\big)\Big)\, dvd\text{\textgreek{sv}}+\\
 & +\int_{\mathcal{S}_{\text{\textgreek{t}}}\cap\{t\le T^{*}\}}\text{\textgreek{q}}_{R}r^{-1-\text{\textgreek{h}}}|\partial_{u}(\text{\textgreek{W}}\text{\textgreek{y}})|^{2}\, dvd\text{\textgreek{sv}}\nonumber 
\end{align}
and $\mathcal{B}_{p,T^{*},\text{\textgreek{h}},\text{\textgreek{e}}}^{(r^{-1}\partial_{\text{\textgreek{sv}}})}[\text{\textgreek{f}}](\text{\textgreek{t}}_{1},\text{\textgreek{t}}_{2})$
is of the form: 
\begin{align}
\mathcal{B}_{p,T^{*},\text{\textgreek{h}},\text{\textgreek{e}}}^{(r^{-1}\partial_{\text{\textgreek{sv}}})}[\text{\textgreek{f}}](\text{\textgreek{t}}_{1},\text{\textgreek{t}}_{2})= & C_{p,\text{\textgreek{h}},\text{\textgreek{e}}}\int_{\mathcal{R}(\text{\textgreek{t}}_{1},\text{\textgreek{t}}_{2})\cap\{t\le T^{*}\}}\text{\textgreek{q}}_{R}\cdot r^{-1+\text{\textgreek{h}}}\cdot\Big(|r^{-1}\partial_{\text{\textgreek{sv}}}\partial_{u}\text{\textgreek{F}}|^{2}+|r^{-1}\partial_{\text{\textgreek{sv}}}\partial_{v}\text{\textgreek{F}}|^{2}\Big)\, dudvd\text{\textgreek{sv}}+\label{eq:BoundRightHandSideDs}\\
 & +C_{p,\text{\textgreek{h}},\text{\textgreek{e}}}\int_{\mathcal{R}(\text{\textgreek{t}}_{1},\text{\textgreek{t}}_{2})\cap\{t\le T^{*}\}}\text{\textgreek{q}}_{R}\cdot r^{p-1-a}\big(|r^{-2}\partial_{\text{\textgreek{sv}}}\partial_{\text{\textgreek{sv}}}\text{\textgreek{F}}|^{2}+|r^{-1}\partial_{\text{\textgreek{sv}}}\text{\textgreek{F}}|^{2}\big)\, dudvd\text{\textgreek{sv}}+\nonumber \\
 & +C_{p,\text{\textgreek{h}},\text{\textgreek{e}}}\int_{\mathcal{R}(\text{\textgreek{t}}_{1},\text{\textgreek{t}}_{2})\cap\{t\le T^{*}\}}\text{\textgreek{q}}_{R}\cdot r^{p-3-a}\big(|r^{-1}\partial_{u}\partial_{\text{\textgreek{sv}}}\text{\textgreek{F}}|^{2}\, dudvd\text{\textgreek{sv}}+\nonumber \\
 & +C_{p,\text{\textgreek{h}},\text{\textgreek{e}}}\int_{\mathcal{R}(\text{\textgreek{t}}_{1},\text{\textgreek{t}}_{2})\cap\{t\le T^{*}\}}\text{\textgreek{q}}_{R}\cdot\big(r^{1+\text{\textgreek{h}}}+r^{p-1-a}\big)\cdot r^{-2}|Err(\text{\textgreek{F}})|^{2}\, dudvd\text{\textgreek{sv}}+\nonumber \\
 & +\int_{\mathcal{R}(\text{\textgreek{t}}_{1},\text{\textgreek{t}}_{2})\cap\{t\le T^{*}\}}\text{\textgreek{q}}_{R}\cdot O_{p,\text{\textgreek{h}},\text{\textgreek{e}}}(r^{p-1})\partial_{v}\partial_{\text{\textgreek{sv}}}\text{\textgreek{F}}\cdot\big(Err(r^{-1}\partial_{\text{\textgreek{sv}}}\text{\textgreek{F}})+r^{-1}Err(\text{\textgreek{F}})\big)\, dudvd\text{\textgreek{sv}}+\nonumber \\
 & +\int_{\mathcal{R}(\text{\textgreek{t}}_{1},\text{\textgreek{t}}_{2})\cap\{t\le T^{*}\}}\text{\textgreek{q}}_{R}\cdot O_{p,\text{\textgreek{h}},\text{\textgreek{e}}}(r^{p-3})\partial_{\text{\textgreek{sv}}}\partial_{\text{\textgreek{sv}}}\text{\textgreek{F}}\cdot Err(\text{\textgreek{F}})\, dudvd\text{\textgreek{sv}}+\nonumber \\
 & +C_{p,\text{\textgreek{h}},\text{\textgreek{e}}}\int_{\mathcal{R}(\text{\textgreek{t}}_{1},\text{\textgreek{t}}_{2})\cap\{t\le T^{*}\}}\text{\textgreek{q}}_{R}\cdot f(r)\partial_{u}(r^{-1}\partial_{\text{\textgreek{sv}}}\text{\textgreek{F}})\cdot Err(r^{-1}\partial_{\text{\textgreek{sv}}}\text{\textgreek{F}})\, dudvd\text{\textgreek{sv}}.\nonumber 
\end{align}

We can now add (\ref{eq:newMethodHyperboloidsDvImprovedFinal}) and
(\ref{eq:newMethodHyperboloidsDsImprovedFinal}) so as to obtain 
\begin{equation}
\begin{split}\mathcal{E}_{bound,R,T^{*}}^{(p)}[\text{\textgreek{W}}^{-1} & \partial_{v}\text{\textgreek{F}}](\text{\textgreek{t}}_{2})+\mathcal{E}_{bound,R,T^{*}}^{(p)}[\text{\textgreek{W}}^{-1}r^{-1}\partial_{\text{\textgreek{sv}}}\text{\textgreek{F}}](\text{\textgreek{t}}_{2})+\int_{\text{\textgreek{t}}_{1}}^{\text{\textgreek{t}}_{2}}\Big(\mathcal{E}{}_{bulk,R,\text{\textgreek{h}},T^{*}}^{(p-1,\partial_{v})}[\text{\textgreek{W}}^{-1}\partial_{v}\text{\textgreek{F}}](\text{\textgreek{t}})+\mathcal{E}{}_{bulk,R,\text{\textgreek{h}},T^{*}}^{(p-1,\partial_{\text{\textgreek{sv}}})}[\text{\textgreek{W}}^{-1}r^{-1}\partial_{\text{\textgreek{sv}}}\text{\textgreek{F}}](\text{\textgreek{t}})\Big)\, d\text{\textgreek{t}}+\\
+\mathcal{E}_{\mathcal{I}^{+},R,T^{*}}^{(p)}[\text{\textgreek{W}}^{-1} & \partial_{v}\text{\textgreek{F}}](\text{\textgreek{t}}_{1},\text{\textgreek{t}}_{2})+\mathcal{E}_{\mathcal{I}^{+},R,T^{*}}^{(p)}[\text{\textgreek{W}}^{-1}r^{-1}\partial_{\text{\textgreek{sv}}}\text{\textgreek{F}}](\text{\textgreek{t}}_{1},\text{\textgreek{t}}_{2})\le\\
\le & \big(1+O_{p,\text{\textgreek{h}}}(\text{\textgreek{e}})\big)\Big(\mathcal{E}_{bound,R,T^{*}}^{(p)}[\text{\textgreek{W}}^{-1}\partial_{v}\text{\textgreek{F}}](\text{\textgreek{t}}_{1})+\mathcal{E}_{bound,R,T^{*}}^{(p)}[\text{\textgreek{W}}^{-1}r^{-1}\partial_{\text{\textgreek{sv}}}\text{\textgreek{F}}](\text{\textgreek{t}}_{1})\Big)+\\
 & +C_{p,\text{\textgreek{h}},\text{\textgreek{e}}}\int_{\mathcal{R}(\text{\textgreek{t}}_{1},\text{\textgreek{t}}_{2})\cap\{t\le T^{*}\}}|\partial\text{\textgreek{q}}_{R}|\cdot\big(r^{p}|\partial^{2}\text{\textgreek{F}}|^{2}+r^{p-2}|\partial\text{\textgreek{F}}|^{2}\big)\, dudvd\text{\textgreek{sv}}+\\
 & +C_{p,\text{\textgreek{h}},\text{\textgreek{e}}}\int_{\mathcal{R}(\text{\textgreek{t}}_{1},\text{\textgreek{t}}_{2})\cap\{t\le T^{*}\}}\text{\textgreek{q}}_{R}\cdot\max\{r^{p-3-a},r^{-3}\}\Big(|\partial_{v}\text{\textgreek{F}}|^{2}+|r^{-1}\partial_{\text{\textgreek{sv}}}\text{\textgreek{F}}|^{2}\Big)\, dudvd\text{\textgreek{sv}}+\\
 & +C_{p,\text{\textgreek{h}},\text{\textgreek{e}}}\int_{\mathcal{R}(\text{\textgreek{t}}_{1},\text{\textgreek{t}}_{2})\cap\{t\le T^{*}\}}\text{\textgreek{q}}_{R}\cdot\big(r^{1+\text{\textgreek{h}}}+r^{p+1}\big)\cdot\Big(|\partial_{v}(\text{\textgreek{W}}F)|^{2}+|r^{-1}\partial_{\text{\textgreek{sv}}}(\text{\textgreek{W}}F)|^{2}+r^{-2}|\text{\textgreek{W}}F|^{2}\Big)\, dudvd\text{\textgreek{sv}}+\\
 & +(p-3)(p-4)\int_{\mathcal{R}(\text{\textgreek{t}}_{1},\text{\textgreek{t}}_{2})\cap\{t\le T^{*}\}}\text{\textgreek{q}}_{R}\cdot\Big(r^{p-3}|r^{-1}\partial_{\text{\textgreek{sv}}}\text{\textgreek{F}}|^{2}+\frac{(d-1)(d-3)}{4}r^{p-5}|\text{\textgreek{F}}|^{2}\Big)\, dudvd\text{\textgreek{sv}}+\\
 & +\int_{\mathcal{R}(\text{\textgreek{t}}_{1},\text{\textgreek{t}}_{2})\cap\{t\le T^{*}\}}\text{\textgreek{q}}_{R}\cdot(1+O(r^{-a}))r^{p-1}\big(r^{-1}\partial_{v}\nabla^{\mathbb{S}^{d-1}}\text{\textgreek{F}}\big)\cdot\big(\partial_{v}(r^{-1}\nabla^{\mathbb{S}^{d-1}}\text{\textgreek{F}})\big)\, dudvd\text{\textgreek{sv}}+\\
 & +C_{p,\text{\textgreek{h}},\text{\textgreek{e}}}Bound_{p,T^{*}}^{\partial_{v},\partial_{\text{\textgreek{sv}}}}[\text{\textgreek{f}}](\text{\textgreek{t}}_{1},\text{\textgreek{t}}_{2})+C_{p,\text{\textgreek{h}}}Bound_{p,T^{*},\text{\textgreek{e}}}^{(n)}[\text{\textgreek{f}}](\text{\textgreek{t}}_{1},\text{\textgreek{t}}_{2})+\mathcal{B}_{p,T^{*},\text{\textgreek{h}},\text{\textgreek{e}}}^{(\partial_{v})}[\text{\textgreek{f}}](\text{\textgreek{t}}_{1},\text{\textgreek{t}}_{2})+\mathcal{B}_{p,T^{*},\text{\textgreek{h}},\text{\textgreek{e}}}^{(r^{-1}\partial_{\text{\textgreek{sv}}})}[\text{\textgreek{f}}](\text{\textgreek{t}}_{1},\text{\textgreek{t}}_{2}),
\end{split}
\label{eq:newMethodHyperboloidsDsDvImprovedAlmostFinal}
\end{equation}
 where 
\begin{align}
Bound_{p,T^{*}}^{\partial_{v},\partial_{\text{\textgreek{sv}}}}[\text{\textgreek{f}}](\text{\textgreek{t}}_{1},\text{\textgreek{t}}_{2})= & C_{p,\text{\textgreek{h}},\text{\textgreek{e}}}\sum_{i=1}^{2}\int_{\mathcal{S}_{\text{\textgreek{t}}_{i}}\cap\{t\le T^{*}\}}\text{\textgreek{q}}_{R}\cdot r^{p-2-a}\Big(|\partial_{u}\partial_{v}\text{\textgreek{F}}|^{2}+|r^{-1}\partial_{\text{\textgreek{sv}}}\partial_{v}\text{\textgreek{F}}|^{2}+|\partial_{v}\text{\textgreek{F}}|^{2}+r^{-2}|\text{\textgreek{F}}|^{2}\Big)\, dvd\text{\textgreek{sv}}+\\
 & +C_{p,\text{\textgreek{h}},\text{\textgreek{e}}}\sum_{i=1}^{2}\int_{\mathcal{S}_{\text{\textgreek{t}}_{i}}\cap\{t\le T^{*}\}}\text{\textgreek{q}}_{R}\cdot r^{p-2-a}\Big(|r^{-1}\partial_{u}\partial_{\text{\textgreek{sv}}}\text{\textgreek{F}}|^{2}+|r^{-2}\partial_{\text{\textgreek{sv}}}\partial_{\text{\textgreek{sv}}}\text{\textgreek{F}}|^{2}+|r^{-1}\partial_{\text{\textgreek{sv}}}\text{\textgreek{F}}|^{2}\Big)\, dvd\text{\textgreek{sv}}+\nonumber \\
 & +C_{p,\text{\textgreek{h}},\text{\textgreek{e}}}\int_{\mathcal{R}(\text{\textgreek{t}}_{1},\text{\textgreek{t}}_{2})\cap\{t=T^{*}\}}\text{\textgreek{q}}_{R}\cdot r^{p-2-a}\Big(|\partial_{u}\partial_{v}\text{\textgreek{F}}|^{2}+|r^{-1}\partial_{\text{\textgreek{sv}}}\partial_{v}\text{\textgreek{F}}|^{2}+|\partial_{v}\text{\textgreek{F}}|^{2}+r^{-2}|\text{\textgreek{F}}|^{2}\Big)\, dvd\text{\textgreek{sv}}+\nonumber \\
 & +C_{p,\text{\textgreek{h}},\text{\textgreek{e}}}\int_{\mathcal{R}(\text{\textgreek{t}}_{1},\text{\textgreek{t}}_{2})\cap\{t=T^{*}\}}\text{\textgreek{q}}_{R}\cdot r^{p-2-a}\Big(|r^{-1}\partial_{u}\partial_{\text{\textgreek{sv}}}\text{\textgreek{F}}|^{2}+|r^{-2}\partial_{\text{\textgreek{sv}}}\partial_{\text{\textgreek{sv}}}\text{\textgreek{F}}|^{2}+|r^{-1}\partial_{\text{\textgreek{sv}}}\text{\textgreek{F}}|^{2}\Big)\, dvd\text{\textgreek{sv}}+\nonumber \\
 & +C_{p,\text{\textgreek{h}},\text{\textgreek{e}}}\int_{\mathcal{S}_{\text{\textgreek{t}}_{1}}\cap\{t\le T^{*}\}}\text{\textgreek{q}}_{R}J_{\text{\textgreek{m}}}^{T}(\text{\textgreek{W}}^{-1}\partial_{v}\text{\textgreek{F}})\bar{n}^{\text{\textgreek{m}}}+C_{p,\text{\textgreek{h}},\text{\textgreek{e}}}\int_{\mathcal{S}_{\text{\textgreek{t}}_{1}}\cap\{t\le T^{*}\}}\text{\textgreek{q}}_{R}J_{\text{\textgreek{m}}}^{T}(\text{\textgreek{W}}^{-1}r^{-1}\partial_{\text{\textgreek{sv}}}\text{\textgreek{F}})\bar{n}^{\text{\textgreek{m}}}.\nonumber 
\end{align}

Since $2<p\le4$, the second to the end line of the right hand side
of (\ref{eq:newMethodHyperboloidsDsDvDuImprovedAlmostFinal}) can
be absorbed by the left hand side after using a Cauchy--Schwarz inequality,
since 
\begin{multline}
\int_{\mathcal{R}(\text{\textgreek{t}}_{1},\text{\textgreek{t}}_{2})\cap\{t\le T^{*}\}}\text{\textgreek{q}}_{R}\cdot(1+O(r^{-a}))r^{p-1}\big(r^{-1}\partial_{v}\nabla^{\mathbb{S}^{d-1}}\text{\textgreek{F}}\big)\cdot\big(\partial_{v}(r^{-1}\nabla^{\mathbb{S}^{d-1}}\text{\textgreek{F}})\big)\, dudvd\text{\textgreek{sv}}\le\\
\le\frac{1}{2}\int_{\mathcal{R}(\text{\textgreek{t}}_{1},\text{\textgreek{t}}_{2})\cap\{t\le T^{*}\}}\text{\textgreek{q}}_{R}\cdot\Big((6-p)r^{p-1}|r^{-1}\partial_{\text{\textgreek{sv}}}\partial_{v}\text{\textgreek{F}}|^{2}+pr^{p-1}|\partial_{v}(r^{-1}\partial_{\text{\textgreek{sv}}}\text{\textgreek{F}})|^{2}\Big)\, dvd\text{\textgreek{sv}}.
\end{multline}
Moreover, the third to the end line of the right hand side of (\ref{eq:newMethodHyperboloidsDsDvDuImprovedAlmostFinal})
can be bounded by the left hand side of (\ref{eq:newMethodFinalFormZeroMassHyperboloidsDu})
for $p-2$ in place of $p$ (notice the importance of the $(p-4)$
factor appearing in front of this term). Therefore, after adding to
(\ref{eq:newMethodHyperboloidsDsDvDuImprovedAlmostFinal}) a large
multiple of (\ref{eq:NewMethodSimpleHyperboloids}) for $p-2$ in
place of $p$ and a large multiple of (\ref{eq:newMethodHardy}) for
$p-2-\text{\textgreek{d}}$ in place of $p$, we obtain: 
\begin{equation}
\begin{split}\mathcal{E}_{bound,R,T^{*}}^{(p)}[\text{\textgreek{W}}^{-1}\partial_{v} & \text{\textgreek{F}}](\text{\textgreek{t}}_{2})+\mathcal{E}_{bound,R,T^{*}}^{(p)}[\text{\textgreek{W}}^{-1}r^{-1}\partial_{\text{\textgreek{sv}}}\text{\textgreek{F}}](\text{\textgreek{t}}_{2})+\mathcal{E}_{bound,R,T^{*}}^{(p-2)}[\text{\textgreek{f}}](\text{\textgreek{t}}_{2})+\\
+\int_{\text{\textgreek{t}}_{1}}^{\text{\textgreek{t}}_{2}}\Big(\mathcal{E}{}_{bulk,R,\text{\textgreek{h}},T^{*}}^{(p-1,\partial_{v})} & [\text{\textgreek{W}}^{-1}\partial_{v}\text{\textgreek{F}}](\text{\textgreek{t}})+\mathcal{E}{}_{bulk,R,\text{\textgreek{h}},T^{*}}^{(p-1,\partial_{\text{\textgreek{sv}}})}[\text{\textgreek{W}}^{-1}r^{-1}\partial_{\text{\textgreek{sv}}}\text{\textgreek{F}}](\text{\textgreek{t}})+\mathcal{E}{}_{bulk,R,\text{\textgreek{h}},T^{*}}^{(p-3,\partial_{\text{\textgreek{sv}}})}[\text{\textgreek{f}}](\text{\textgreek{t}})\Big)\, d\text{\textgreek{t}}+\\
+\mathcal{E}_{\mathcal{I}^{+},R,T^{*}}^{(p)}[\text{\textgreek{W}}^{-1}\partial_{v}\text{\textgreek{F}}] & (\text{\textgreek{t}}_{1},\text{\textgreek{t}}_{2})+\mathcal{E}_{\mathcal{I}^{+},R,T^{*}}^{(p)}[\text{\textgreek{W}}^{-1}r^{-1}\partial_{\text{\textgreek{sv}}}\text{\textgreek{F}}](\text{\textgreek{t}}_{1},\text{\textgreek{t}}_{2})+\mathcal{E}_{\mathcal{I}^{+},R,T^{*}}^{(p-2)}[\text{\textgreek{f}}](\text{\textgreek{t}}_{1},\text{\textgreek{t}}_{2})\le\\
\le & \big(1+O_{p,\text{\textgreek{h}}}(\text{\textgreek{e}})\big)\Big(\mathcal{E}_{bound,R,T^{*}}^{(p)}[\text{\textgreek{W}}^{-1}\partial_{v}\text{\textgreek{F}}](\text{\textgreek{t}}_{1})+\mathcal{E}_{bound,R,T^{*}}^{(p)}[\text{\textgreek{W}}^{-1}r^{-1}\partial_{\text{\textgreek{sv}}}\text{\textgreek{F}}](\text{\textgreek{t}}_{1})+C\cdot\mathcal{E}_{bound,R,T^{*}}^{(p-2)}[\text{\textgreek{f}}](\text{\textgreek{t}}_{1})\Big)+\\
 & +C_{p,\text{\textgreek{h}},\text{\textgreek{e}}}\int_{\mathcal{R}(\text{\textgreek{t}}_{1},\text{\textgreek{t}}_{2})\cap\{t\le T^{*}\}}|\partial\text{\textgreek{q}}_{R}|\cdot\big(r^{p}|\partial^{2}\text{\textgreek{F}}|^{2}+r^{p-2}|\partial\text{\textgreek{F}}|^{2}+r^{p-4}|\text{\textgreek{F}}|^{2}\big)\, dudvd\text{\textgreek{sv}}+\\
 & +C_{p,\text{\textgreek{h}},\text{\textgreek{e}}}\int_{\mathcal{R}(\text{\textgreek{t}}_{1},\text{\textgreek{t}}_{2})\cap\{t\le T^{*}\}}\text{\textgreek{q}}_{R}\cdot\big(r^{1+\text{\textgreek{h}}}+r^{p+1}\big)\cdot\Big(|\partial_{v}(\text{\textgreek{W}}F)|^{2}+|r^{-1}\partial_{\text{\textgreek{sv}}}(\text{\textgreek{W}}F)|^{2}+r^{-2}|\text{\textgreek{W}}F|^{2}\Big)\, dudvd\text{\textgreek{sv}}+\\
 & +C_{p,\text{\textgreek{h}},\text{\textgreek{e}}}Bound_{p,T^{*}}^{\partial_{v},\partial_{\text{\textgreek{sv}}}}[\text{\textgreek{f}}](\text{\textgreek{t}}_{1},\text{\textgreek{t}}_{2})+C_{p,\text{\textgreek{h}}}Bound_{p,T^{*},\text{\textgreek{e}}}^{(n)}[\text{\textgreek{f}}](\text{\textgreek{t}}_{1},\text{\textgreek{t}}_{2})+\mathcal{B}_{p,T^{*},\text{\textgreek{h}},\text{\textgreek{e}}}^{(\partial_{v})}[\text{\textgreek{f}}](\text{\textgreek{t}}_{1},\text{\textgreek{t}}_{2})+\mathcal{B}_{p,T^{*},\text{\textgreek{h}},\text{\textgreek{e}}}^{(r^{-1}\partial_{\text{\textgreek{sv}}})}[\text{\textgreek{f}}](\text{\textgreek{t}}_{1},\text{\textgreek{t}}_{2}).
\end{split}
\label{eq:newMethodHyperboloidsDsDvImprovedAlmostFinal-1}
\end{equation}

Finally, in order to absorb the error terms of the last line the right
hand side of (\ref{eq:newMethodHyperboloidsDsDvImprovedAlmostFinal-1})
into the left hand side, we will need to add to (\ref{eq:newMethodHyperboloidsDsDvImprovedAlmostFinal-1})
a constant multiple of the estimate (\ref{eq:newMethodFinalStatementHyperboloids})
for $\text{\textgreek{W}}^{-1}\partial_{u}\text{\textgreek{F}}$ in
place of $\text{\textgreek{F}}$. \textgreek{B}y following again the
proof of Lemma \ref{lem:NewMethodGeneralCase} for $\text{\textgreek{W}}^{-1}\partial_{u}\text{\textgreek{F}}$
in place of $\text{\textgreek{f}}$ for $p-2$ in place of $p$, we
obtain: 
\begin{equation}
\begin{split}\mathcal{E}_{bound,R,T^{*}}^{(p-2)}[ & \text{\textgreek{W}}^{-1}\partial_{u}\text{\textgreek{F}}](\text{\textgreek{t}}_{2})+\int_{\text{\textgreek{t}}_{1}}^{\text{\textgreek{t}}_{2}}\mathcal{E}_{bulk,R,\text{\textgreek{h}},T^{*}}^{(p-3)}[\text{\textgreek{W}}^{-1}\partial_{u}\text{\textgreek{F}}](\text{\textgreek{t}})\, d\text{\textgreek{t}}+\mathcal{E}_{\mathcal{I}^{+},R,T^{*}}^{(p-2)}[\text{\textgreek{W}}^{-1}\partial_{u}\text{\textgreek{F}}](\text{\textgreek{t}}_{1},\text{\textgreek{t}}_{2})\le\\
\le & \big(1+O_{p,\text{\textgreek{h}}}(\text{\textgreek{e}})\big)\mathcal{E}_{bound,R,T^{*}}^{(p-2)}[\text{\textgreek{W}}^{-1}\partial_{u}\text{\textgreek{F}}](\text{\textgreek{t}}_{1})+C_{p,\text{\textgreek{h}},\text{\textgreek{e}}}\int_{\mathcal{R}(\text{\textgreek{t}}_{1},\text{\textgreek{t}}_{2})\cap\{t\le T^{*}\}}|\partial\text{\textgreek{q}}_{R}|\cdot\big(r^{p-2}|\partial\partial_{u}\text{\textgreek{F}}|^{2}+r^{p-4}|\partial_{u}\text{\textgreek{F}}|^{2}\big)\, dudvd\text{\textgreek{sv}}+\\
 & +C_{p,\text{\textgreek{h}},\text{\textgreek{e}}}\int_{\mathcal{R}(\text{\textgreek{t}}_{1},\text{\textgreek{t}}_{2})\cap\{t\le T^{*}\}}\text{\textgreek{q}}_{R}\cdot\Big(r^{p-3-a}|r^{-1}\partial_{\text{\textgreek{sv}}}\partial_{u}\text{\textgreek{F}}|^{2}+\max\{r^{p-5-a},r^{-3}\}|\partial_{u}\text{\textgreek{F}}|^{2}\Big)\, dudvd\text{\textgreek{sv}}+\\
 & +C_{p,\text{\textgreek{h}},\text{\textgreek{e}}}\int_{\mathcal{R}(\text{\textgreek{t}}_{1},\text{\textgreek{t}}_{2})\cap\{t\le T^{*}\}}\text{\textgreek{q}}_{R}\cdot\big(\frac{1}{2}f(r)(\partial_{v}-\partial_{u})(\text{\textgreek{W}}^{-1}\partial_{u}\text{\textgreek{F}})+\frac{(d-1)f(r)}{r}(\text{\textgreek{W}}^{-1}\partial_{u}\text{\textgreek{F}})\big)\cdot\square(\text{\textgreek{W}}^{-1}\partial_{u}\text{\textgreek{F}})\,\text{\textgreek{W}}^{2}dudvd\text{\textgreek{sv}}+\\
 & +C_{p,\text{\textgreek{h}},\text{\textgreek{e}}}\int_{\mathcal{R}(\text{\textgreek{t}}_{1},\text{\textgreek{t}}_{2})\cap\{t\le T^{*}\}}\text{\textgreek{q}}_{R}\cdot(\partial_{v}+\partial_{u})(\partial_{u}\text{\textgreek{F}})\cdot\text{\textgreek{W}}\square(\text{\textgreek{W}}^{-1}\partial_{u}\text{\textgreek{F}})\, dudvd\text{\textgreek{sv}}+\\
 & +\int_{\mathcal{R}(\text{\textgreek{t}}_{1},\text{\textgreek{t}}_{2})\cap\{t\le T^{*}\}}\text{\textgreek{q}}_{R}\cdot O(r^{p-4-a})\cdot\partial_{\text{\textgreek{sv}}}\partial_{u}\text{\textgreek{F}}\cdot\text{\textgreek{W}}\square\big(\text{\textgreek{W}}^{-1}\partial_{u}\text{\textgreek{F}}\big)\, dudvd\text{\textgreek{sv}}-\\
 & -\int_{\mathcal{R}(\text{\textgreek{t}}_{1},\text{\textgreek{t}}_{2})\cap\{t\le T^{*}\}}\text{\textgreek{q}}_{R}\cdot r^{p-2}\partial_{v}\partial_{u}\text{\textgreek{F}}\cdot\text{\textgreek{W}}\square\big(\text{\textgreek{W}}^{-1}\partial_{u}\text{\textgreek{F}}\big)\, dudvd\text{\textgreek{sv}}+C_{p,\text{\textgreek{h}},\text{\textgreek{e}}}Bound_{p-2,T^{*}}^{\partial_{u}}[\text{\textgreek{f}}](\text{\textgreek{t}}_{1},\text{\textgreek{t}}_{2}),
\end{split}
\label{eq:newMethodFinalFormZeroMassHyperboloidsDu}
\end{equation}
 where 
\begin{align}
Bound_{p-2,T^{*}}^{\partial_{u}}[\text{\textgreek{f}}](\text{\textgreek{t}}_{1},\text{\textgreek{t}}_{2})= & \sum_{i=1}^{2}\int_{\mathcal{S}_{\text{\textgreek{t}}_{i}}\cap\{t\le T^{*}\}}\text{\textgreek{q}}_{R}\cdot r^{p-4-a}|\partial_{u}\text{\textgreek{F}}|^{2}\, dvd\text{\textgreek{sv}}+\int_{\mathcal{R}(\text{\textgreek{t}}_{1},\text{\textgreek{t}}_{2})\cap\{t=T^{*}\}}\text{\textgreek{q}}_{R}\cdot r^{p-4-a}|\partial_{u}\text{\textgreek{F}}|^{2}\, dvd\text{\textgreek{sv}}+\\
 & +\int_{\mathcal{S}_{\text{\textgreek{t}}_{1}}\cap\{t\le T^{*}\}}\text{\textgreek{q}}_{R}J_{\text{\textgreek{m}}}^{T}(\text{\textgreek{W}}^{-1}\partial_{u}\text{\textgreek{F}})\bar{n}^{\text{\textgreek{m}}}.\nonumber 
\end{align}

Adding to (\ref{eq:newMethodHyperboloidsDsDvImprovedAlmostFinal-1})
a large multiple of the estimate (\ref{eq:newMethodFinalFormZeroMassHyperboloidsDu})
(implementing also a Hardy type inequality for the $\partial_{v}\partial_{u}\text{\textgreek{F}}$
term), and using the fact that according to Lemma \ref{lem:Commutator expressions}
\begin{align}
\text{\textgreek{W}}\square\big(\text{\textgreek{W}}^{-1}\partial_{u}\text{\textgreek{F}}\big)= & (1+O(r^{-1-a}))\cdot\big(\partial_{u}(\text{\textgreek{W}}F)+O(r^{-1-a})\text{\textgreek{W}}F\big)+\label{eq:CommutationDu-1-1}\\
 & +O(r^{-3})\partial_{\text{\textgreek{sv}}}\partial_{\text{\textgreek{sv}}}\text{\textgreek{F}}+\sum_{j=0}^{1}Err\big(\partial_{u}^{j}\text{\textgreek{F}}\big),\nonumber 
\end{align}
we thus obtain: 
\begin{equation}
\begin{split}\mathcal{E}_{bound,R,T^{*}}^{(p)}[\text{\textgreek{W}}^{-1} & \partial_{v}\text{\textgreek{F}}](\text{\textgreek{t}}_{2})+\mathcal{E}_{bound,R,T^{*}}^{(p)}[\text{\textgreek{W}}^{-1}r^{-1}\partial_{\text{\textgreek{sv}}}\text{\textgreek{F}}](\text{\textgreek{t}}_{2})++\mathcal{E}_{bound,R,T^{*}}^{(p-2)}[\text{\textgreek{W}}^{-1}\partial_{u}\text{\textgreek{F}}](\text{\textgreek{t}}_{2})+\mathcal{E}_{bound,R,T^{*}}^{(p-2)}[\text{\textgreek{f}}](\text{\textgreek{t}}_{2})+\\
+\int_{\text{\textgreek{t}}_{1}}^{\text{\textgreek{t}}_{2}}\Big(\mathcal{E}{}_{bulk,R,\text{\textgreek{h}},T^{*}}^{(p-1,\partial_{v})} & [\text{\textgreek{W}}^{-1}\partial_{v}\text{\textgreek{F}}](\text{\textgreek{t}})+\mathcal{E}{}_{bulk,R,\text{\textgreek{h}},T^{*}}^{(p-1,\partial_{\text{\textgreek{sv}}})}[\text{\textgreek{W}}^{-1}r^{-1}\partial_{\text{\textgreek{sv}}}\text{\textgreek{F}}](\text{\textgreek{t}})+\mathcal{E}_{bulk,R,\text{\textgreek{h}},T^{*}}^{(p-3)}[\text{\textgreek{W}}^{-1}\partial_{u}\text{\textgreek{F}}](\text{\textgreek{t}})+\mathcal{E}{}_{bulk,R,\text{\textgreek{h}},T^{*}}^{(p-3,\partial_{\text{\textgreek{sv}}})}[\text{\textgreek{f}}](\text{\textgreek{t}})\Big)\, d\text{\textgreek{t}}+\\
+\mathcal{E}_{\mathcal{I}^{+},R,T^{*}}^{(p)}[\text{\textgreek{W}}^{-1} & \partial_{v}\text{\textgreek{F}}](\text{\textgreek{t}}_{1},\text{\textgreek{t}}_{2})+\mathcal{E}_{\mathcal{I}^{+},R,T^{*}}^{(p)}[\text{\textgreek{W}}^{-1}r^{-1}\partial_{\text{\textgreek{sv}}}\text{\textgreek{F}}](\text{\textgreek{t}}_{1},\text{\textgreek{t}}_{2})+\mathcal{E}_{\mathcal{I}^{+},R,T^{*}}^{(p-2)}[\text{\textgreek{W}}^{-1}\partial_{u}\text{\textgreek{F}}](\text{\textgreek{t}}_{1},\text{\textgreek{t}}_{2})+\mathcal{E}_{\mathcal{I}^{+},R,T^{*}}^{(p-2)}[\text{\textgreek{f}}](\text{\textgreek{t}}_{1},\text{\textgreek{t}}_{2})\le\\
\le & \big(1+O_{p,\text{\textgreek{h}}}(\text{\textgreek{e}})\big)\Big(\mathcal{E}_{bound,R,T^{*}}^{(p)}[\text{\textgreek{W}}^{-1}\partial_{v}\text{\textgreek{F}}](\text{\textgreek{t}}_{1})+\mathcal{E}_{bound,R,T^{*}}^{(p)}[\text{\textgreek{W}}^{-1}r^{-1}\partial_{\text{\textgreek{sv}}}\text{\textgreek{F}}](\text{\textgreek{t}}_{1})\Big)\\
 & +C_{p,\text{\textgreek{h}}}\cdot\Big(\mathcal{E}_{bound,R,T^{*}}^{(p-2)}[\text{\textgreek{f}}](\text{\textgreek{t}}_{1})+\mathcal{E}_{bound,R,T^{*}}^{(p-2)}[\text{\textgreek{W}}^{-1}\partial_{u}\text{\textgreek{F}}](\text{\textgreek{t}}_{1})\Big)+\\
 & +C_{p,\text{\textgreek{h}},\text{\textgreek{e}}}\int_{\mathcal{R}(\text{\textgreek{t}}_{1},\text{\textgreek{t}}_{2})\cap\{t\le T^{*}\}}|\partial\text{\textgreek{q}}_{R}|\cdot\big(r^{p}|\partial^{2}\text{\textgreek{F}}|^{2}+r^{p-2}|\partial\text{\textgreek{F}}|^{2}+r^{p-4}|\text{\textgreek{F}}|^{2}\big)\, dudvd\text{\textgreek{sv}}+\\
 & +C_{p,\text{\textgreek{h}},\text{\textgreek{e}}}\int_{\mathcal{R}(\text{\textgreek{t}}_{1},\text{\textgreek{t}}_{2})\cap\{t\le T^{*}\}}\text{\textgreek{q}}_{R}\cdot r^{p+1}\Big(|\partial_{v}(\text{\textgreek{W}}F)|^{2}+|r^{-1}\partial_{\text{\textgreek{sv}}}(\text{\textgreek{W}}F)|^{2}+r^{-2}|\partial_{u}(\text{\textgreek{W}}F)|^{2}+r^{-2}|\text{\textgreek{W}}F|^{2}\Big)\, dudvd\text{\textgreek{sv}}+\\
 & +C_{p,\text{\textgreek{h}},\text{\textgreek{e}}}\int_{\mathcal{R}(\text{\textgreek{t}}_{1},\text{\textgreek{t}}_{2})\cap\{t\le T^{*}\}}\text{\textgreek{q}}_{R}\cdot r^{1+\text{\textgreek{h}}}\Big(|\partial_{v}(\text{\textgreek{W}}F)|^{2}+|r^{-1}\partial_{\text{\textgreek{sv}}}(\text{\textgreek{W}}F)|^{2}+|\partial_{u}(\text{\textgreek{W}}F)|^{2}+|\text{\textgreek{W}}F|^{2}\Big)\, dudvd\text{\textgreek{sv}}+\\
 & +C_{p,\text{\textgreek{h}},\text{\textgreek{e}}}Bound_{p-2,T^{*}}^{\partial_{u}}[\text{\textgreek{f}}](\text{\textgreek{t}}_{1},\text{\textgreek{t}}_{2})+C_{p,\text{\textgreek{h}},\text{\textgreek{e}}}Bound_{p,T^{*}}^{\partial_{v},\partial_{\text{\textgreek{sv}}}}[\text{\textgreek{f}}](\text{\textgreek{t}}_{1},\text{\textgreek{t}}_{2})+C_{p,\text{\textgreek{h}}}Bound_{p,T^{*},\text{\textgreek{e}}}^{(n)}[\text{\textgreek{f}}](\text{\textgreek{t}}_{1},\text{\textgreek{t}}_{2})+\\
 & +\mathcal{B}_{p,T^{*},\text{\textgreek{h}},\text{\textgreek{e}}}^{(\partial_{v})}[\text{\textgreek{f}}](\text{\textgreek{t}}_{1},\text{\textgreek{t}}_{2})+\mathcal{B}_{p,T^{*},\text{\textgreek{h}},\text{\textgreek{e}}}^{(r^{-1}\partial_{\text{\textgreek{sv}}})}[\text{\textgreek{f}}](\text{\textgreek{t}}_{1},\text{\textgreek{t}}_{2})+\mathcal{B}_{p-2,T^{*},\text{\textgreek{h}},\text{\textgreek{e}}}^{(\partial_{u})}[\text{\textgreek{f}}](\text{\textgreek{t}}_{1},\text{\textgreek{t}}_{2}).
\end{split}
\label{eq:newMethodHyperboloidsDsDvDuImprovedAlmostFinal}
\end{equation}

where 
\begin{align}
\mathcal{B}_{p-2,T^{*},\text{\textgreek{h}},\text{\textgreek{e}}}^{(\partial_{u})}[\text{\textgreek{f}}](\text{\textgreek{t}}_{1},\text{\textgreek{t}}_{2})= & C_{p,\text{\textgreek{h}},\text{\textgreek{e}}}\int_{\mathcal{R}(\text{\textgreek{t}}_{1},\text{\textgreek{t}}_{2})\cap\{t\le T^{*}\}}\text{\textgreek{q}}_{R}\cdot f(r)\Big\{(\partial_{v}-\partial_{u})\big(\text{\textgreek{W}}^{-1}\partial_{u}\text{\textgreek{F}}\big)+O(r^{-1})\big(\text{\textgreek{W}}^{-1}\partial_{u}\text{\textgreek{F}}\big)\Big\}\cdot\label{eq:BoundRightHandSideDu}\\
 & \hphantom{C_{p,\text{\textgreek{h}},\text{\textgreek{e}}}\int_{\mathcal{R}(\text{\textgreek{t}}_{1},\text{\textgreek{t}}_{2})\cap\{t\le T^{*}\}}}\cdot\Big\{ O(r^{-3})\partial_{\text{\textgreek{sv}}}\partial_{\text{\textgreek{sv}}}\text{\textgreek{F}}+Err(\partial_{u}\text{\textgreek{F}})+Err(\text{\textgreek{F}})\Big\}\,\text{\textgreek{W}}dudvd\text{\textgreek{sv}}+\nonumber \\
 & +\int_{\mathcal{R}(\text{\textgreek{t}}_{1},\text{\textgreek{t}}_{2})\cap\{t\le T^{*}\}}\text{\textgreek{q}}_{R}\cdot O(r^{p-3-a})\partial_{\text{\textgreek{sv}}}\partial_{u}\text{\textgreek{F}}\cdot\Big\{ O(r^{-3})\partial_{\text{\textgreek{sv}}}\partial_{\text{\textgreek{sv}}}\text{\textgreek{F}}+Err(\partial_{u}\text{\textgreek{F}})+Err(\text{\textgreek{F}})\Big\}\, dudvd\text{\textgreek{sv}}-\nonumber \\
 & -\int_{\mathcal{R}(\text{\textgreek{t}}_{1},\text{\textgreek{t}}_{2})\cap\{t\le T^{*}\}}\text{\textgreek{q}}_{R}\cdot O(r^{p-2})\partial_{v}\partial_{u}\text{\textgreek{F}}\cdot\Big\{ O(r^{-3})\partial_{\text{\textgreek{sv}}}\partial_{\text{\textgreek{sv}}}\text{\textgreek{F}}+Err(\partial_{u}\text{\textgreek{F}})+Err(\text{\textgreek{F}})\Big\}\, dudvd\text{\textgreek{sv}}.\nonumber 
\end{align}

After adding to (\ref{eq:newMethodHyperboloidsDsDvDuImprovedAlmostFinal})
a large (in terms of $p,\text{\textgreek{h}},\text{\textgreek{e}}$)
multiple of (\ref{eq:NewMethodSimpleHyperboloids}) for $p-2$ in
place of $p$ and of (\ref{eq:newMethodHardy}) for $\min\{p-2,2-\text{\textgreek{d}}\}$
in place of $p$, the 
\[
C_{p,\text{\textgreek{h}}}Bound_{p,T^{*},\text{\textgreek{e}}}^{(n)}[\text{\textgreek{f}}](\text{\textgreek{t}}_{1},\text{\textgreek{t}}_{2})+\mathcal{B}_{p,T^{*},\text{\textgreek{h}},\text{\textgreek{e}}}^{(\partial_{v})}[\text{\textgreek{f}}](\text{\textgreek{t}}_{1},\text{\textgreek{t}}_{2})+\mathcal{B}_{p,T^{*},\text{\textgreek{h}},\text{\textgreek{e}}}^{(r^{-1}\partial_{\text{\textgreek{sv}}})}[\text{\textgreek{f}}](\text{\textgreek{t}}_{1},\text{\textgreek{t}}_{2})+\mathcal{B}_{p-2,T^{*},\text{\textgreek{h}},\text{\textgreek{e}}}^{(\partial_{u})}[\text{\textgreek{f}}](\text{\textgreek{t}}_{1},\text{\textgreek{t}}_{2})
\]
summand of the right hand side of (\ref{eq:newMethodHyperboloidsDsDvDuImprovedAlmostFinal})
can be absorbed by the positive terms of the left hand side using
an integration by parts scheme similar to the one used in the proof
of Theorem \ref{thm:NewMethodFinalStatementHyperboloids} for the
top order terms and a simple Cauchy--Schwarz inequality for the lower
order terms. Moreover (in view also of Lemma \ref{lem:BoundednessGeometric},
as well as the fact that $\text{\textgreek{h}}'<1+a$), the summand
\[
C_{p,\text{\textgreek{h}},\text{\textgreek{e}}}Bound_{p-2,T^{*}}^{\partial_{u}}[\text{\textgreek{f}}](\text{\textgreek{t}}_{1},\text{\textgreek{t}}_{2})+C_{p,\text{\textgreek{h}},\text{\textgreek{e}}}Bound_{p,T^{*}}^{\partial_{v},\partial_{\text{\textgreek{sv}}}}[\text{\textgreek{f}}](\text{\textgreek{t}}_{1},\text{\textgreek{t}}_{2})
\]
can be bounded by 
\[
\int_{\mathcal{S}_{\text{\textgreek{t}}_{1}}}\text{\textgreek{q}}_{R}J_{\text{\textgreek{m}}}^{T}(r^{-k_{2}}\partial_{v}^{k_{1}}\partial_{\text{\textgreek{sv}}}^{k_{2}}\partial_{u}^{k_{3}}\text{\textgreek{f}})\bar{n}^{\text{\textgreek{m}}}
\]
plus $O_{p,\text{\textgreek{h}}}(\text{\textgreek{e}})$ times the
existing terms in the left and right hand side of (\ref{eq:newMethodHyperboloidsDsDvDuImprovedAlmostFinal})
and (\ref{eq:newMethodHardy}) (for $\min\{p-2,2-\text{\textgreek{d}}\}$
in place of $p$). Therefore, provided $\text{\textgreek{e}}$ has
been fixed small in terms of $p,\text{\textgreek{h}},\text{\textgreek{d}}$
and $R$ is large enough in terms of $p,\text{\textgreek{h}},\text{\textgreek{d}}$,
the desired inequality (\ref{eq:newMethodDu+DvPhi}) for $k=1$ readily
follows from (\ref{eq:newMethodHyperboloidsDsDvDuImprovedAlmostFinal}),
(\ref{eq:NewMethodSimpleHyperboloids}) (for $p-2$ in place of $p$)
and (\ref{eq:newMethodHardy}) (for $\min\{p-2,2-\text{\textgreek{d}}\}$
in place of $p$) after letting $T^{*}\rightarrow+\infty$.\qed

\section{\label{sec:FriedlanderRadiation}Friedlander radiation field on $\mathcal{I}^{+}$}

In this section, we will establish the existence of the Friedlander
radiation field on future null infinity $\mathcal{I}^{+}$ for solutions
$\text{\textgreek{f}}$ to $\square\text{\textgreek{f}}=F$ on general
asymptotically flat spacetimes $(\mathcal{M}^{d+1},g)$, $d\ge3$,
with the asymptotics (\ref{eq:MetricUR}), provided the source term
$F$ decays suitably fast in terms of $r$. This result is essentially
a ``soft'' corollary of the results of the previous sections.

\subsection{\label{sub:GeometricAssumptionsFriedlander}Assumptions on the spacetimes
under consideration}

Let $(\mathcal{M}^{d+1},g)$, $d\ge3$, be a time oriented smooth
Lorentzian manifold, possibly with non empty piecewise smooth boundary
$\partial\mathcal{M}$. We will assume the following condition on
the asymptotics of $(\mathcal{M},g)$:

\begin{MyDescription}[leftmargin=3.0em]

\item [(G1)]{ \label{enu:AsymptoticFlatness} \emph{Asymptotic flatness.}
There exists an open subset $\mathcal{N}_{af,\mathcal{M}}\subset\mathcal{M}$
such that each connected component of $\mathcal{N}_{af,\mathcal{M}}$
is diffeomorphic to $\mathbb{R}\times(1,+\infty)\times\mathbb{S}^{d-1}$.
Fixing such a diffeomorphism (i.\,e.~a coordinate chart) on each
component, we will denote with $(u,r,\text{\textgreek{sv}})$ the
associated coordinate functions. Furthermore, we assume that on each
component of $\mathcal{N}_{af,\mathcal{M}}$, the metric $g$ takes
the form (\ref{eq:MetricUR}):
\begin{align}
g= & -4\Big(1-\frac{2M(u,\text{\textgreek{sv}})}{r}+O(r^{-1-a})\Big)du^{2}-\Big(4+O(r^{-1-a})\Big)dudr+r^{2}\cdot\Big(g_{\mathbb{S}^{d-1}}+O(r^{-1})\Big)+\label{eq:MetricUR-2}\\
 & +\big(h_{3}^{as}(u,\text{\textgreek{sv}})+O(r^{-a})\big)dud\text{\textgreek{sv}}+O(r^{-a})drd\text{\textgreek{sv}}+O(r^{-2-a})dr^{2}.\nonumber 
\end{align}
 Thus, $(\mathcal{N}_{af,\mathcal{M}},g)$ should be thought of as
a disjoint union of copies of the manifold $(\mathcal{N}_{af},g)$
of Section \ref{sec:GeometrySpacetimes}.

}\end{MyDescription}
\begin{rem*}
We do not assume that $(\mathcal{M},g)$ is globally hyperbolic. 
\end{rem*}
We will also set $t=u+r$ and $v=t-r$ on $\mathcal{N}_{af,\mathcal{M}}$.
Let also $T$ denote the vector field $\partial_{t}$ in the $(t,r,\text{\textgreek{sv}})$
coordinate system on each component of $\mathcal{N}_{af,\mathcal{M}}$.
We extend $r$ as a smooth function on the whole of $\mathcal{M}$
by the requirement that $0\le r\le2$ on $\mathcal{M}\backslash\mathcal{N}_{af,\mathcal{M}}$.
Finally, we also construct the function $\bar{t}_{\text{\textgreek{h}}'}:\{r\gg1\}\subset\mathcal{N}_{af,\mathcal{M}}\rightarrow\mathbb{R}$
for some fixed $0<\text{\textgreek{h}}'<1+a$ as in Section \ref{sec:GeometrySpacetimes}.

\subsection{Existence of the Friedlander radiation field on future null infinity}

We will prove the following result:
\begin{thm}
\label{thm:FriedlanderRadiation}Let $(\mathcal{M}^{d+1},g),\, d\ge3$
satisfy the assumption \ref{enu:AsymptoticFlatness}, and let $R_{0}>0$
be large in terms of the geometry of $(\mathcal{M},g)$ (this implies
that $\{r\ge R_{0}\}\subset\mathcal{N}_{af,\mathcal{M}}$). Let $\text{\textgreek{f}}:\mathcal{M}\rightarrow\mathbb{C}$
be a smooth function solving $\square_{g}\text{\textgreek{f}}=F$,
such that for some $0<\text{\textgreek{d}}<1$ and for any integer$0\le j\le\lceil\frac{d}{2}\rceil$
and $\text{\textgreek{t}}\in\mathbb{R}$ the following quantity is
finite on each connected component of $\mathcal{N}_{af,\mathcal{M}}$:
\begin{equation}
\mathcal{E}_{in}^{(1+\text{\textgreek{d}}+2j,j+1)}[\text{\textgreek{f}}](0)+\sum_{i=1}^{j+1}\sum_{k_{1}+k_{2}+k_{3}=i-1}\int_{J^{+}(\{t=0\})\cap J^{-}(\{\bar{t}=\text{\textgreek{t}}\})\cap\{r\ge R_{0}\}}r^{\text{\textgreek{d}}+2i}\Big(|r^{-k_{2}-k_{3}}\partial_{v}^{k_{1}}\partial_{\text{\textgreek{sv}}}^{k_{2}}\partial_{u}^{k_{3}}(\text{\textgreek{W}}F)|^{2}\Big)\, dudvd\text{\textgreek{sv}}<+\infty\label{eq:NormFinitenessFriedlanderRadiation}
\end{equation}
where 
\begin{align}
\mathcal{E}_{in}^{(p,k)}[\text{\textgreek{f}}](0)=\sum_{j=1}^{k}\sum_{k_{1}+k_{2}+k_{3}=j-1}\Big\{ & \int_{\{t=0\}\cap\{r\ge R_{0}\}}\Big(r^{p-2(k-j)}\big|r^{-k_{2}-k_{3}}\partial_{v}^{k_{1}+1}\partial_{\text{\textgreek{sv}}}^{k_{2}}\partial_{u}^{k_{3}}(\text{\textgreek{W}}\text{\textgreek{f}})\big|^{2}+\label{eq:NewMethodBoundaryNorm-1}\\
 & \hphantom{\int_{\{t=0\}\cap\{r}}+r^{-1-\text{\textgreek{h}}'}\big(r^{p-2(k-j)}\big|r^{-k_{2}-k_{3}-1}\partial_{v}^{k_{1}}\partial_{\text{\textgreek{sv}}}^{k_{2}+1}\partial_{u}^{k_{3}}(\text{\textgreek{W}}\text{\textgreek{f}})\big|^{2}+\nonumber \\
 & \hphantom{\int_{\{t=0\}\cap\{r}}+\big((d-3)r^{p-2-2(k-j)}+\min\{r^{p-2-2(k-j)},r^{-\text{\textgreek{d}}-2(k-j)}\}\big)\big|r^{-k_{2}-k_{3}}\partial_{v}^{k_{1}}\partial_{\text{\textgreek{sv}}}^{k_{2}}\partial_{u}^{k_{3}}(\text{\textgreek{W}}\text{\textgreek{f}})\big|^{2}\big)\Big)\, dvd\text{\textgreek{sv}}+\nonumber \\
 & +\int_{\{t=0\}\cap\{r\ge R_{0}\}}J_{\text{\textgreek{m}}}^{T}(r^{-k_{2}}\partial_{v}^{k_{1}}\partial_{\text{\textgreek{sv}}}^{k_{2}}\partial_{u}^{k_{3}}\text{\textgreek{f}})\bar{n}^{\text{\textgreek{m}}}\Big\}\nonumber 
\end{align}
and the derivatives are considered with respect to the $(u,v,\text{\textgreek{sv}})$
coordinate system on each connected component of $\{r\ge R_{0}\}$.
Then for each connected component of $\mathcal{N}_{af,\mathcal{M}}$,
the Friedlander radiation field 
\begin{equation}
\text{\textgreek{F}}_{\mathcal{I}^{+}}(u,\text{\textgreek{sv}})\doteq\lim_{v\rightarrow+\infty}\text{\textgreek{W}}\cdot\text{\textgreek{f}}(u,v,\text{\textgreek{sv}})\label{eq:RadiationField}
\end{equation}
 exists on $\mathbb{R}\times\mathbb{S}^{d-1}$. 

Moreover, if (\ref{eq:NormFinitenessFriedlanderRadiation}) holds
for all integers $0\le j\le l+\lceil\frac{d}{2}\rceil$ for some integer
$l\ge0$, then $\text{\textgreek{F}}_{\mathcal{I}^{+}}$ is $C^{l-1}$
in $(u,\text{\textgreek{sv}})$ and for all integers $j_{1},j_{2},j_{3}\ge0$
with $j_{1}+j_{2}+j_{3}\le l$ and any $(u,\text{\textgreek{sv}})\in\mathbb{R}\times\mathbb{S}^{d-1}$
the limit 
\begin{equation}
\lim_{v\rightarrow+\infty}r^{j_{1}}\partial_{v}^{j_{1}}\partial_{\text{\textgreek{sv}}}^{j_{3}}\partial_{u}^{j_{3}}(\text{\textgreek{W}}\text{\textgreek{f}})(u,v,\text{\textgreek{sv}})\label{eq:LimitHigherOrderFriedlander}
\end{equation}
 is finite. 

In particular, if (\ref{eq:NormFinitenessFriedlanderRadiation}) holds
for all integers $j\ge0$ (for instance when $\text{\textgreek{f}}$
solves $\square\text{\textgreek{f}}=0$ with smooth and compactly
supported initial data on $\{t=0\}$), then $\text{\textgreek{F}}_{\mathcal{I}^{+}}$
is smooth and (\ref{eq:LimitHigherOrderFriedlander}) exists for all
integers $j_{1},j_{2},j_{3}\ge0$.\end{thm}
\begin{rem*}
In the case $j_{1}\ge1$ we actually expect that the limit (\ref{eq:LimitHigherOrderFriedlander})
is identically $0$ when $\text{\textgreek{f}}$ solves $\square\text{\textgreek{f}}=0$
with compactly supported initial data. This expectation is justified
by the fact follows from the fact that on spacetimes admitting a smooth
conformal compactification of future null infinity, the following
stronger statement (in comparison to (\ref{eq:LimitHigherOrderFriedlander}))
holds for all integers $j_{1},j_{2},j_{3}\ge0$:
\begin{equation}
\lim_{v\rightarrow+\infty}r^{2j_{1}}\partial_{v}^{j_{1}}\partial_{\text{\textgreek{sv}}}^{j_{3}}\partial_{u}^{j_{3}}(\text{\textgreek{W}}\text{\textgreek{f}})(u,v,\text{\textgreek{sv}})<+\infty.\label{eq:LimitHigherOrderFriedlander-1}
\end{equation}
 Notice, however, that in our case we are also including spacetimes
which do not necessarily admit a smooth conformal compactification
at $\mathcal{I}^{+}$.\end{rem*}
\begin{proof}
We will assume without loss of generality that $\text{\textgreek{f}}$
is real valued. We will also work on a fixed connected component of
$\mathcal{N}_{af,\mathcal{M}}$.

Let $\{v_{n}\}_{n\in\mathbb{N}}$ be an increasing sequence of positive
real numbers tending to $+\infty$. Let us also fix a smooth function
$\text{\textgreek{q}}_{R}:\mathcal{M}\rightarrow[0,1]$ such that
$\text{\textgreek{q}}_{R}\equiv0$ on $\{r\le R\}$ and $\text{\textgreek{q}}_{R}\equiv1$
on $\{r\ge R+1\}$ for some large fixed constant $R\ge R_{0}\gg1$. 

By repeating the proof of Theorem \ref{thm:NewMethodDu+DvPhi} for
$p=1+\text{\textgreek{d}}$ in the spacetime region 
\begin{equation}
\mathcal{D}_{0,\text{\textgreek{t}}}\doteq J^{+}(\{t=0\})\cap J^{-}(\{\bar{t}=\text{\textgreek{t}}\})\cap\{r\ge R_{0}\}
\end{equation}
 (instead of the region bounded by two hyperboloidal hypersurfaces)
for any $\text{\textgreek{t}}\in\mathbb{R}$, we can readily bound
for any $k\in\mathbb{N}$ and $\text{\textgreek{t}}\in\mathbb{R}$:
\begin{align}
\mathcal{E}_{bound,R;\text{\textgreek{d}}}^{(1+\text{\textgreek{d}}+2(k-1),k)}[\text{\textgreek{f}}](\text{\textgreek{t}})\lesssim_{k,\text{\textgreek{d}}} & \mathcal{E}_{in}^{(1+\text{\textgreek{d}}+2(k-1),k)}[\text{\textgreek{f}}](0)+\sum_{j=0}^{k}\int_{\mathcal{D}_{0,\text{\textgreek{t}}}}|\partial\text{\textgreek{q}}_{R}|\cdot r^{\text{\textgreek{d}}+2j-1}|\partial^{j}\text{\textgreek{f}}|^{2}+\label{eq:FromNewMethod-2}\\
 & +\sum_{j=1}^{k}\sum_{k_{1}+k_{2}+k_{3}=j-1}\int_{\mathcal{D}_{0,\text{\textgreek{t}}}}\text{\textgreek{q}}_{R}\cdot r^{\text{\textgreek{d}}+2j-2k_{3}}\Big(|r^{-k_{2}}\partial_{v}^{k_{1}}\partial_{\text{\textgreek{sv}}}^{k_{2}}\partial_{v}^{k_{1}}\partial_{u}^{k_{3}}(\text{\textgreek{W}}F)|^{2}\Big)\, dudvd\text{\textgreek{sv}},\nonumber 
\end{align}
 where 
\begin{align}
\mathcal{E}_{bound,R;\text{\textgreek{d}}}^{(p,k)}[\text{\textgreek{f}}](\text{\textgreek{t}})=\sum_{j=1}^{k}\sum_{k_{1}+k_{2}+k_{3}=j-1}\Big\{ & \int_{\{\bar{t}=\text{\textgreek{t}}\}}\text{\textgreek{q}}_{R}\Big(r^{p-2(k-j)}\big|r^{-k_{2}-k_{3}}\partial_{v}^{k_{1}+1}\partial_{\text{\textgreek{sv}}}^{k_{2}}\partial_{u}^{k_{3}}(\text{\textgreek{W}}\text{\textgreek{f}})\big|^{2}+\label{eq:NewMethodBoundaryNorm}\\
 & \hphantom{\int_{\{\bar{t}=\text{\textgreek{t}}\}}}+r^{-1-\text{\textgreek{h}}'}\big(r^{p-2(k-j)}\big|r^{-k_{2}-k_{3}-1}\partial_{v}^{k_{1}}\partial_{\text{\textgreek{sv}}}^{k_{2}+1}\partial_{u}^{k_{3}}(\text{\textgreek{W}}\text{\textgreek{f}})\big|^{2}+\nonumber \\
 & \hphantom{\int_{\{\bar{t}=\text{\textgreek{t}}\}}}+\big((d-3)r^{p-2-2(k-j)}+\min\{r^{p-2-2(k-j)},r^{-\text{\textgreek{d}}-2(k-j)}\}\big)\big|r^{-k_{2}-k_{3}}\partial_{v}^{k_{1}}\partial_{\text{\textgreek{sv}}}^{k_{2}}\partial_{u}^{k_{3}}(\text{\textgreek{W}}\text{\textgreek{f}})\big|^{2}\big)\Big)\, dvd\text{\textgreek{sv}}+\nonumber \\
 & +\int_{\{\bar{t}=\text{\textgreek{t}}\}}\text{\textgreek{q}}_{R}J_{\text{\textgreek{m}}}^{T}(r^{-k_{2}}\partial_{v}^{k_{1}}\partial_{\text{\textgreek{sv}}}^{k_{2}}\partial_{u}^{k_{3}}\text{\textgreek{f}})\bar{n}^{\text{\textgreek{m}}}\Big\}\nonumber 
\end{align}
 (the $\partial_{v},\partial_{\text{\textgreek{sv}}},\partial_{u}$
coordinate vector fields are a priori only defined in the region $\{r\ge R\}$,
but since the integrand is multiplied with the cut-off $\text{\textgreek{q}}_{R}$
the expression in (\ref{eq:NewMethodBoundaryNorm}) is well defined).

Using the fundamental theorem of calculus and the expression (\ref{eq:Hyperboloids})
for $\bar{t}$, we can bound for any $C^{1}$ function $\text{\textgreek{y}}:\mathcal{M}\rightarrow\mathbb{R}$
and any $\text{\textgreek{t}}\in\mathbb{R}$, $n_{0}\in\mathbb{N}$:
\begin{equation}
\begin{split}\sum_{n=1}^{n_{0}}\int_{\mathbb{S}^{d-1}}\Big|\text{\textgreek{q}}_{R}\text{\textgreek{W}}\text{\textgreek{y}}|_{\{\bar{t}=\text{\textgreek{t}}\}}(v_{n+1},\text{\textgreek{sv}}) & -\text{\textgreek{q}}_{R}\text{\textgreek{W}}\text{\textgreek{y}}|_{\{\bar{t}=\text{\textgreek{t}}\}}(v_{n},\text{\textgreek{sv}})\Big|\, d\text{\textgreek{sv}}\lesssim\int_{\{\bar{t}=\text{\textgreek{t}}\}\cap\{v\le v_{n_{0}}\}}\big(|\partial_{v}(\text{\textgreek{q}}_{R}\text{\textgreek{W}}\text{\textgreek{y}})|+r^{-1}|\partial_{u}(\text{\textgreek{q}}_{R}\text{\textgreek{W}}\text{\textgreek{y}})|\big)\, dvd\text{\textgreek{sv}}\lesssim_{\text{\textgreek{d}}}\\
 & \lesssim_{\text{\textgreek{d}}}\int_{\{\bar{t}=\text{\textgreek{t}}\}}|\partial\text{\textgreek{q}}_{R}|\cdot|\text{\textgreek{W}}\text{\textgreek{y}}|\, d\text{\textgreek{sv}}+\Big(\int_{\{\bar{t}=\text{\textgreek{t}}\}\cap\{v\le v_{n_{0}}\}}r^{1+\text{\textgreek{d}}}\text{\textgreek{q}}_{R}\big(|\partial_{v}(\text{\textgreek{W}}\text{\textgreek{y}})|^{2}+r^{-2}|\partial_{u}(\text{\textgreek{W}}\text{\textgreek{y}})|^{2}\big)\, dvd\text{\textgreek{sv}}\Big)^{1/2}.
\end{split}
\label{eq:FundamentalTheoremOfCalculus}
\end{equation}
In the above, $\text{\textgreek{q}}_{R}\text{\textgreek{W}}\text{\textgreek{y}}|_{\{\bar{t}=\text{\textgreek{t}}\}}$
is considered as a function of $(v,\text{\textgreek{sv}})$, since
$(v,\text{\textgreek{sv}})$ is a valid parametrisation of $\{\bar{t}=\text{\textgreek{t}}\}\cap\{r\ge R\}$
where $\text{\textgreek{q}}_{R}\text{\textgreek{W}}\text{\textgreek{y}}$
is supported. Moreover, the coordinate derivatives $\partial_{u},\partial_{v}$
in the right hand side are defined in the $(u,v,\text{\textgreek{sv}})$
coordinate system in the region $\{r\ge R\}$.

Therefore, from (\ref{eq:FromNewMethod-2}), (\ref{eq:NewMethodBoundaryNorm})
and (\ref{eq:FundamentalTheoremOfCalculus}) for $\text{\textgreek{W}}^{-1}r^{j_{1}}\partial_{v}^{j_{1}}\partial_{\text{\textgreek{sv}}}^{j_{3}}\partial_{u}^{j_{3}}(\text{\textgreek{W}}\text{\textgreek{f}})$
in place of $\text{\textgreek{y}}$ we readily deduce that for any
$k\in\mathbb{N}$ and any $n_{0}\in\mathbb{N}$: 
\begin{equation}
\begin{split}\Big(\sum_{n=1}^{n_{0}}\sum_{j=0}^{k-1}\sum_{j_{1}+j_{2}+j_{3}=j}\int_{\mathbb{S}^{d-1}} & \Big|\text{\textgreek{q}}_{R}r^{j_{1}}\partial_{v}^{j_{1}}\partial_{\text{\textgreek{sv}}}^{j_{3}}\partial_{u}^{j_{3}}(\text{\textgreek{W}}\text{\textgreek{f}})|_{\{\bar{t}=\text{\textgreek{t}}\}}(v_{n+1},\text{\textgreek{sv}})-\text{\textgreek{q}}_{R}r^{j_{1}}\partial_{v}^{j_{1}}\partial_{\text{\textgreek{sv}}}^{j_{3}}\partial_{u}^{j_{3}}(\text{\textgreek{W}}\text{\textgreek{f}})|_{\{\bar{t}=\text{\textgreek{t}}\}}(v_{n},\text{\textgreek{sv}})\Big|\, d\text{\textgreek{sv}}\Big)^{2}\lesssim_{k,\text{\textgreek{d}}}\\
\lesssim_{k,\text{\textgreek{d}}} & \mathcal{E}_{in}^{(1+\text{\textgreek{d}}+2(k-1),k)}[\text{\textgreek{f}}](0)+\sum_{j=0}^{k}\int_{\mathcal{D}_{0,\text{\textgreek{t}}}}|\partial\text{\textgreek{q}}_{R}|\cdot r^{\text{\textgreek{d}}+2j-1}|\partial^{j}\text{\textgreek{f}}|^{2}+\Big(\sum_{j=0}^{k-1}\int_{\{\bar{t}=\text{\textgreek{t}}\}}|\partial\text{\textgreek{q}}_{R}|r^{j}|\partial^{j}\text{\textgreek{f}}|\Big)^{2}+\\
 & +\sum_{j=1}^{k}\sum_{k_{1}+k_{2}+k_{3}=j-1}\int_{\mathcal{D}_{0,\text{\textgreek{t}}}}\text{\textgreek{q}}_{R}\cdot r^{\text{\textgreek{d}}+2j-2k_{3}}\Big(|r^{-k_{2}}\partial_{v}^{k_{1}}\partial_{\text{\textgreek{sv}}}^{k_{2}}\partial_{v}^{k_{1}}\partial_{u}^{k_{3}}(\text{\textgreek{W}}F)|^{2}\Big)\, dudvd\text{\textgreek{sv}}.
\end{split}
\label{eq:FromNewMethod&Fundmental}
\end{equation}

Using the Sobolev inequality on $\mathbb{S}^{d-1}$ 
\begin{equation}
||\text{\textgreek{y}}||_{L^{\infty}(\mathbb{S}^{d-1})}^{2}\lesssim\sum_{j=0}^{\lceil\frac{d}{2}\rceil}\int_{\mathbb{S}^{d-1}}|\partial_{\text{\textgreek{sv}}}^{j}\text{\textgreek{y}}|^{2}\, d\text{\textgreek{sv}},
\end{equation}
 from (\ref{eq:FromNewMethod&Fundmental}) (and the fact that $||\text{\textgreek{y}}||_{L^{1}(\mathbb{S}^{d-1})}\lesssim_{d}||\text{\textgreek{y}}||_{L^{2}(\mathbb{S}^{d-1})}$)
we infer that for any $k_{0},n_{0}\in\mathbb{N}$: 
\begin{equation}
\begin{split}\sup_{\text{\textgreek{sv}}\in\mathbb{S}^{d-1}}\Big\{\sum_{n=1}^{n_{0}}\sum_{j=0}^{k_{0}}\sum_{j_{1}+j_{2}+j_{3}=j}\Big| & \text{\textgreek{q}}_{R}r^{j_{1}}\partial_{v}^{j_{1}}\partial_{\text{\textgreek{sv}}}^{j_{3}}\partial_{u}^{j_{3}}(\text{\textgreek{W}}\text{\textgreek{f}})|_{\{\bar{t}=\text{\textgreek{t}}\}}(v_{n+1},\text{\textgreek{sv}})-\text{\textgreek{q}}_{R}r^{j_{1}}\partial_{v}^{j_{1}}\partial_{\text{\textgreek{sv}}}^{j_{3}}\partial_{u}^{j_{3}}(\text{\textgreek{W}}\text{\textgreek{f}})|_{\{\bar{t}=\text{\textgreek{t}}\}}(v_{n},\text{\textgreek{sv}})\Big|^{2}\Big\}\lesssim_{k,\text{\textgreek{d}}}\\
\lesssim_{k,\text{\textgreek{d}}} & \mathcal{E}_{in}^{(1+\text{\textgreek{d}}+2(k_{0}+\lceil\frac{d}{2}\rceil),k_{0}+\lceil\frac{d+2}{2}\rceil)}[\text{\textgreek{f}}](0)+\sum_{j=0}^{k_{0}+\lceil\frac{d+2}{2}\rceil}\int_{\mathcal{D}_{0,\text{\textgreek{t}}}}|\partial\text{\textgreek{q}}_{R}|\cdot r^{\text{\textgreek{d}}+2j-1}|\partial^{j}\text{\textgreek{f}}|^{2}+\Big(\sum_{j=0}^{k_{0}+\lceil\frac{d}{2}\rceil}\int_{\{\bar{t}=\text{\textgreek{t}}\}}|\partial\text{\textgreek{q}}_{R}|r^{j}|\partial^{j}\text{\textgreek{f}}|\Big)^{2}+\\
 & +\sum_{j=1}^{k_{0}+\lceil\frac{d+2}{2}\rceil}\sum_{k_{1}+k_{2}+k_{3}=j-1}\int_{\mathcal{D}_{0,\text{\textgreek{t}}}}\text{\textgreek{q}}_{R}\cdot r^{\text{\textgreek{d}}+2j-2k_{3}}\Big(|r^{-k_{2}}\partial_{v}^{k_{1}}\partial_{\text{\textgreek{sv}}}^{k_{2}}\partial_{v}^{k_{1}}\partial_{u}^{k_{3}}(\text{\textgreek{W}}F)|^{2}\Big)\, dudvd\text{\textgreek{sv}}.
\end{split}
\label{eq:AfterSobolev}
\end{equation}

In view of (\ref{eq:NormFinitenessFriedlanderRadiation}) and the
fact that $\mathcal{D}_{0,\text{\textgreek{t}}}\cap supp(\partial\text{\textgreek{q}}_{R})$
and $\{\bar{t}=\text{\textgreek{t}}\}\cap supp(\partial\text{\textgreek{q}}_{R})$
are compact subsets of $\mathcal{M}$, from (\ref{eq:AfterSobolev})
we deduce that for all $\text{\textgreek{t}}_{0}\in\mathbb{R}$ and
all $k_{0}\in\mathbb{N}$ such that (\ref{eq:NormFinitenessFriedlanderRadiation})
holds for all $k\le k_{0}+\lceil\frac{d}{2}\rceil$, the following
quantity is finite independently of $n_{0}$: 
\begin{equation}
\sup_{\text{\textgreek{sv}}\in\mathbb{S}^{d-1},\text{\textgreek{t}}\le\text{\textgreek{t}}_{0}}\Big\{\sum_{n=1}^{n_{0}}\sum_{j=0}^{k_{0}}\sum_{j_{1}+j_{2}+j_{3}=j}\Big|\text{\textgreek{q}}_{R}r^{j_{1}}\partial_{v}^{j_{1}}\partial_{\text{\textgreek{sv}}}^{j_{3}}\partial_{u}^{j_{3}}(\text{\textgreek{W}}\text{\textgreek{f}})|_{\{\bar{t}=\text{\textgreek{t}}\}}(v_{n+1},\text{\textgreek{sv}})-\text{\textgreek{q}}_{R}r^{j_{1}}\partial_{v}^{j_{1}}\partial_{\text{\textgreek{sv}}}^{j_{3}}\partial_{u}^{j_{3}}(\text{\textgreek{W}}\text{\textgreek{f}})|_{\{\bar{t}=\text{\textgreek{t}}\}}(v_{n},\text{\textgreek{sv}})\Big|^{2}\Big\}\le C_{k_{0},\text{\textgreek{d}}}[\text{\textgreek{f}}](\text{\textgreek{t}}_{0})<+\infty.\label{eq:CauchySequence}
\end{equation}
 Therefore, by letting $n_{0}\rightarrow+\infty$ (\ref{eq:CauchySequence})
yields that for any $\text{\textgreek{t}}\in\mathbb{R}$, $\text{\textgreek{sv}}\in\mathbb{S}^{d-1}$,
any $k_{0}\in\mathbb{N}$ such that (\ref{eq:NormFinitenessFriedlanderRadiation})
holds for all $k\le k_{0}+\lceil\frac{d}{2}\rceil$ and any $j_{1}+j_{2}+j_{3}=k_{0}$,
the sequence 
\begin{equation}
\Big\{ r^{j_{1}}\partial_{v}^{j_{1}}\partial_{\text{\textgreek{sv}}}^{j_{3}}\partial_{u}^{j_{3}}(\text{\textgreek{W}}\text{\textgreek{f}})|_{\{\bar{t}=\text{\textgreek{t}}\}}(v_{n},\text{\textgreek{sv}})\Big\}_{n\in\mathbb{N}}\label{eq:SequenceToBeCauchy}
\end{equation}
 is a Cauchy sequence (and hence (\ref{eq:LimitHigherOrderFriedlander})
follows). 

Moreover, since $\bar{t}-u=O(r^{-\text{\textgreek{h}}'})$, from (\ref{eq:CauchySequence})
we infer the limit 
\begin{equation}
\text{\textgreek{F}}_{\mathcal{I}^{+}}(u,\text{\textgreek{sv}})=\lim_{v\rightarrow+\infty}\text{\textgreek{W}}\text{\textgreek{f}}(v,u,\text{\textgreek{sv}})
\end{equation}
 exists and is a $C^{k_{0}-1}$ function of $(u,\text{\textgreek{sv}})$
(if $k_{0}\neq0$). 
\end{proof}

\subsection{Estimates for $\text{\textgreek{F}}_{\mathcal{I}^{+}}$ provided
by Lemma \ref{thm:NewMethodFinalStatementHyperboloids}}

The following corollary is a straightforward consequence of Theorems
\ref{thm:NewMethodDu+DvPhi} and \ref{thm:FriedlanderRadiation}:
\begin{cor}
\label{lem:NewMethodFinalStatementHyperboloidsFriedlander}Let $\mathcal{N}$
be any connected component of $\mathcal{N}_{af,\mathcal{M}}$. Then
for any $k\in\mathbb{N}$, any $2k-2<p\le2k$ , any given $0<\text{\textgreek{h}}<a$
and $0<\text{\textgreek{d}}<1$, any $R>0$ large enough in terms
of $p,\text{\textgreek{h}},\text{\textgreek{d}},k$, any $\text{\textgreek{t}}_{1}\le\text{\textgreek{t}}_{2}$
and any smooth cut-off $\text{\textgreek{q}}_{R}:\mathcal{M}\rightarrow[0,1]$
supported in $\{r\ge R\}\cap\mathcal{N}$, the following inequality
holds for any smooth function $\text{\textgreek{f}}:\mathcal{M}\rightarrow\mathbb{C}$
solving $\square_{g}\text{\textgreek{f}}=F$ with suitably decaying
initial data on $\{t=0\}$: 
\begin{equation}
\begin{split}\sum_{k_{1}+k_{2}=k}\int_{\mathcal{I}^{+}\cap\{\text{\textgreek{t}}_{1}\le u\le\text{\textgreek{t}}_{2}\}} & r^{p-2k}\big|\partial_{\text{\textgreek{sv}}}^{k_{1}}\partial_{u}^{k_{2}}\text{\textgreek{F}}_{\mathcal{I}^{+}}\big|^{2}\, dud\text{\textgreek{sv}}\lesssim_{p,\text{\textgreek{h}},\text{\textgreek{d}},k}\\
\lesssim_{p,\text{\textgreek{h}},\text{\textgreek{d}},k} & \mathcal{E}_{bound,R;\text{\textgreek{d}}}^{(p,k)}[\text{\textgreek{f}}](\text{\textgreek{t}}_{2})+\sum_{j=0}^{k}\int_{\mathcal{R}(\text{\textgreek{t}}_{1},\text{\textgreek{t}}_{2})}|\partial\text{\textgreek{q}}_{R}|\cdot r^{p-2(k-j)}|\partial^{j}\text{\textgreek{f}}|^{2}+\\
 & +\sum_{j=1}^{k}\sum_{k_{1}+k_{2}+k_{3}=j-1}\int_{\mathcal{R}(\text{\textgreek{t}}_{1},\text{\textgreek{t}}_{2})}\text{\textgreek{q}}_{R}\cdot(r^{p+1-2k_{3}-2(k-j)}+r^{1+\text{\textgreek{h}}}\big)\Big(|r^{-k_{2}}\partial_{v}^{k_{1}}\partial_{\text{\textgreek{sv}}}^{k_{2}}\partial_{v}^{k_{1}}\partial_{u}^{k_{3}}(\text{\textgreek{W}}F)|^{2}\Big)\, dudvd\text{\textgreek{sv}},
\end{split}
\label{eq:BoundNewMethodFriedlander}
\end{equation}

where we have adopted the convention 
\[
r^{p-2k}|_{\mathcal{I}^{+}}=\begin{cases}
1, & p=2k\\
0, & p<2k.
\end{cases}
\]

\end{cor}

\section{\label{sec:Firstdecay}Polynomial decay $\bar{t}^{-1}$ for solutions
to $\square_{g}\text{\textgreek{f}}=0$}

In this Section, we will generalise the results of \cite{DafRod7}
by showing that on asymptotically flat spacetimes $\mathcal{M}$ with
possibly non-empty timelike boundary $\partial_{tim}\mathcal{M}$,
a $\bar{t}{}^{-1}$ polynomial decay rate hols for solutions $\text{\textgreek{f}}$
to $\square\text{\textgreek{f}}=0$ with suitable boundary conditions
on $\partial_{tim}\mathcal{M}$, provided some specific geometric
conditions on the interior region of $(\mathcal{M},g)$ hold and assuming
that an integrated local energy decay statement (possibly with loss
of derivatives) holds for $\text{\textgreek{f}}$.

\subsection{\label{sub:AssumptionsFirstDecay}Assumptions on the class of spacetimes
$(\mathcal{M},g)$ under consideration}

\subsubsection{\label{sub:GeometricAssumptionsAndConstructions}Geometric Assumptions
on $(\mathcal{M},g)$ and related geometric constructions}

Let $(\mathcal{M}^{d+1},g)$, $d\ge3$, be a smooth Lorentzian manifold
with possibly non empty piecewise smooth boundary $\partial\mathcal{M}$.
We assume that $(\mathcal{M},g)$ satisfies the Assumption \ref{enu:AsymptoticFlatness}
on asymptotic flatness. We will now proceed to state a few more assumptions
on the geometric structure of $(\mathcal{M},g)$, and present some
geometric constructions that will be used later.

\paragraph{Assumptions on the causal structure of $(\mathcal{M},g)$ and $(\partial\mathcal{M},g|_{\partial\mathcal{M}})$}

Since we will need to establish some global estimates for solutions
to the wave equation (\ref{eq:WaveEquation}) on $(\mathcal{M},g)$,
we will need to impose some conditions on the causal structure of
$(\mathcal{M},g)$ and its boundary.

\begin{MyDescription}[leftmargin=3.0em]

\item [(G2)]{ \label{enu:PartitionBoundary} \emph{Partition of the
boundary. }We assume that the boundary $\partial\mathcal{M}$ (if
non-empty) can be split into two components (not necessarily connected)
\begin{equation}
\partial\mathcal{M}=\partial_{tim}\mathcal{M}\cup\partial_{hor}\mathcal{M},
\end{equation}
 where $(\partial_{tim}\mathcal{M},g|_{\partial_{tim}\mathcal{M}})$
is a smooth Lorentzian manifold (i.\,e.~$\partial_{tim}\mathcal{M}$
is a smooth timelike hypersurface with respect to $g$) and $(\partial_{hor}\mathcal{M},g|_{\partial_{hor}\mathcal{M}})$
is piecewise smooth and degenerate pseudo-Riemannian manifold (i.\,e.~$\partial_{hor}\mathcal{M}$
is a null hypersurface with respect to $g$).}

\item [(G3)]{ \label{enu:GlobalHyperbolicity} \emph{Global Hyperbolicity.
}Let $\tilde{\mathcal{M}}_{tim}$ denote the double of $\mathcal{M}$
along $\partial_{tim}\mathcal{M}$. We will denote as $i_{or}:\mathcal{M}\rightarrow\tilde{\mathcal{M}}_{tim}$
the natural inclusion of $\mathcal{M}$ into $\tilde{\mathcal{M}}_{tim}$,
while $i_{ref}:\mathcal{M}\rightarrow\tilde{\mathcal{M}}_{tim}$ will
denote the reflection map along $\partial_{tim}\mathcal{M}$. We assume
that $\tilde{\mathcal{M}}_{tim}$ is globally hyperbolic. Let $\tilde{\text{\textgreek{S}}}$
be a Cauchy hypersurface of $\tilde{\mathcal{M}}_{tim}$. We will
denote with $\text{\textgreek{S}}$ the restriction of $\tilde{\text{\textgreek{S}}}$
on $\mathcal{M}$. We will also fix a time function $t$ associated
with $\tilde{\text{\textgreek{S}}}$ on $\tilde{\mathcal{M}}_{tim}$,
i.\,e.~$g(\nabla t,\nabla t)<0$ on $\tilde{\mathcal{M}}_{tim}$
and $\tilde{\text{\textgreek{S}}}\equiv\{t=0\}$. Notice that with
the help of $t$ we can identify $\tilde{\mathcal{M}}_{tim}$ with
$\mathbb{R}\times\tilde{\text{\textgreek{S}}}$ and $\mathcal{M}$
with $\mathbb{R}\times\text{\textgreek{S}}$.}

\item [(G4)] {\label{enu:DomainOfOuterCommunications}\emph{Domain
of outer communications.} Let $\mathcal{N}_{af,\mathcal{M}}$ be the
open subset of $\mathcal{M}$ defined in Assumption \ref{enu:AsymptoticFlatness},
where $g$ has the asymptotically flat form (\ref{eq:MetricUR-2}),
and let $\tilde{\mathcal{N}}_{af,\mathcal{M}}=i_{or}(\mathcal{N}_{af,\mathcal{M}})\cup i_{ref}(\mathcal{N}_{af,\mathcal{M}})$.
Having identified $\tilde{\mathcal{M}}_{tim}$ with $\mathbb{R}\times\tilde{\text{\textgreek{S}}}$,
we assume that $\tilde{\mathcal{N}}_{af,\mathcal{M}}=\mathbb{R}\times(\tilde{\text{\textgreek{S}}}\backslash K)$
for some compact $K\subset\tilde{\text{\textgreek{S}}}$, and that
$\tilde{\mathcal{N}}_{af,\mathcal{M}}$ has a finite number of connected
components.%
\footnote{If $\partial\mathcal{M}\neq\emptyset$, then it is necessary that
$K\cap\partial\mathcal{M}\neq\emptyset$.%
} Moreover, we assume that the domain of dependence of $\mathbb{R}\times(\tilde{\text{\textgreek{S}}}\backslash K)$
is the whole of $\tilde{\mathcal{M}}_{tim}$. }

\end{MyDescription}

In view of Assumption \ref{enu:DomainOfOuterCommunications}, $i_{or}(\partial_{hor}\mathcal{M})\cup i_{ref}(\partial_{hor}\mathcal{M})$
(if non-empty) constitutes the \emph{event horizon} of $\tilde{\mathcal{M}}_{tim}$.
From now on we will use the notation $\mathcal{H}$ for $\partial_{hor}\mathcal{M}$
and we will call $\mathcal{H}$ the event horizon of $\mathcal{M}$.
Using the fact that $i_{or}(\mathcal{H})\cup i_{ref}(\mathcal{H})$
constitutes the event horizon of a globally hyperbolic spacetime,
we can define the future event horizon $\mathcal{H}^{+}$ and the
past event horizon $\mathcal{H}^{-}$ of $\mathcal{M}$ by the requirement
that 
\begin{itemize}
\item $\mathcal{H}^{+}$,$\mathcal{H}^{-}$ are piecewise smooth achronal
hypersurfaces (possibly with boundary and not necessarily connected)
such that $\mathcal{H}^{+}\cup\mathcal{H}^{-}=\mathcal{H}$
\item $\mathcal{H}^{-}\subset J^{-}(\mathcal{H}^{+})$
\end{itemize}
We will also assume that $\mathcal{H}^{+}$ is smooth (if non-empty).

\paragraph*{Assumptions on the existence of a well behaved foliation by hyperboloidal
hypersurfaces}

Fixing $\text{\textgreek{h}}'=1$, we will assume that the function
$\bar{t}_{\text{\textgreek{h}}'}$ originally defined on the subset
$\{r\gg1\}$ of $\mathcal{N}_{af,\mathcal{M}}$ can be extended as
a smooth function on $\mathcal{M}$ satisfying the following conditions: 

\begin{MyDescription}[leftmargin=3.0em]

\item [(G5)]{\label{enu:Hyperboloids}$\bar{t}$ is given by the relation
(\ref{eq:Hyperboloids}) in the region $\{r\ge R\}$ of each connected
component of $\mathcal{N}_{af,\mathcal{M}}$ for some $R\gg1$.}

\item [(G6)]{ \label{enu:DomainDependence} For any $0\le\text{\textgreek{t}}_{1}\le\text{\textgreek{t}}_{2}$,
$\{\bar{t}=\text{\textgreek{t}}_{2}\}$ is contained in the future
domain of dependence of $\{\bar{t}=\text{\textgreek{t}}_{1}\}$.%
\footnote{The future domain of dependence $\mathcal{D}^{+}(\mathcal{B})$ of
a set $\mathcal{B}\subset\mathcal{M}$ is defined as the set of all
points $p\in\mathcal{M}$ such that all past inextendible causal curves
$\text{\textgreek{g}}$ emanating from $p$ intersect $\mathcal{B}$,
where now $\text{\textgreek{g}}$ is \underline{not} considered past
inextendible if it has a past endpoint $q$ on $\partial_{tim}\mathcal{M}$,
since from $q$ one can further extend $\text{\textgreek{g}}$ by
a causal path inside $J^{-}(q)\backslash\partial_{tim}\mathcal{M}$.%
}}

\item [(G7)]{\label{enu:Spacelike}$g(\nabla\bar{t},\nabla\bar{t})<0$
everywhere on $\mathcal{M}\cap\{\bar{t}\ge0\}$, where $\nabla\bar{t}$
denotes the gradient of $\bar{t}$ with respect to $g$. Moreover,
$-C\le g(\nabla\bar{t},\nabla\bar{t})\le-c<0$ in the region $\{\bar{t}\ge0\}\backslash\mathcal{N}_{af,\mathcal{M}}$,
for some $C,c>0$.}

\end{MyDescription}

It will be convenient to have a globally defined future directed timelike
vector field $N$ adjusted to the choice of our foliation $\{\bar{t}=const\}$.
Therefore, we will fix $N$ to be a timelike future directed vector
field on $\mathcal{M}$ such that $N\equiv\big(-g(\nabla\bar{t},\nabla\bar{t})\big)^{-1}\cdot\nabla\bar{t}$
on $\{\bar{t}\ge0\}\backslash\mathcal{N}_{af,\mathcal{M}}$, $N\equiv T$
in the region $\{\bar{t}\ge0\}\cap\{r\ge2R\}$ of  $\mathcal{N}_{af,\mathcal{M}}$,
and the relations $-C\le g(N,N)\le-C^{-1}<0$ and $d\bar{t}(N)=1$
hold everywhere on $\mathcal{M}$ for some $C>0$. The existence of
such a vector field follows from time orientability of $\mathcal{M}$
and the convexity of the set 
\[
\mathfrak{F}_{p}=\{X\in T_{p}\mathcal{M}|g(X,X)<0\mbox{ and }d\bar{t}(X)=1\}
\]
 for each $p\in\mathcal{M}$.

We also extend the function $r$ (defined originally in the region
$\mathcal{N}_{af,\mathcal{M}}$) as a Morse function (but not necessarily
as a coordinate function) on the whole of $\mathcal{M}$, under the
assumption that $r\ge0$ everywhere on $\mathcal{M}$, $r\equiv0$
on $\partial\mathcal{M}$ and $N(r)=0$. In this way, the asymptotically
flat region $\mathcal{N}_{af,\mathcal{M}}$ of $\mathcal{M}$ will
correspond to the region $\{r\gg1\}$.

Assumption \ref{enu:Spacelike} implies that, if $dg_{\bar{t}}$ denotes
the volume form of the induced Riemannian metric $g_{\bar{t}}$ on
the $\{\bar{t}=const\}$ hypersurfaces, then there exists a $C>0$
such that for any measurable function $f:\mathcal{M}\rightarrow[0,+\infty)$
and any $0\le\text{\textgreek{t}}_{1}\le\text{\textgreek{t}}_{2}$
we have the equivalence 
\begin{equation}
\int_{\{\text{\textgreek{t}}_{1}\le\bar{t}\le\text{\textgreek{t}}_{2}\}\cap\{r\le R\}}f\, dvol_{g}\sim_{R}\int_{\text{\textgreek{t}}1}^{\text{\textgreek{t}}_{2}}\Big(\int_{\{\bar{t}=\text{\textgreek{sv}}\}\cap\{r\le R\}}f\, dg_{\bar{t}}\Big)d\text{\textgreek{sv}}.
\end{equation}

Notice that Assumptions \ref{enu:DomainDependence} and \ref{enu:Spacelike}
also imply that for $\bar{t}\ge0$ the level sets of the extended
$\bar{t}$ intersect transversely $\mathcal{H}^{+}$ (if $\mathcal{H}^{+}\neq\emptyset$).%
\footnote{As an example, on Schwarzschild exterior the function $\bar{t}$ could
not have been chosen to coincide with the coordinate function $t$
in a neighborhood horizon, but it can coincide with $t^{*}$ (see
i.\,e. \cite{DafRod6}).%
} It will be useful to denote 
\begin{equation}
\mathcal{H}_{\text{\textgreek{t}}}\doteq\mathcal{H}^{+}\cap\{\bar{t}=\text{\textgreek{t}}\}.
\end{equation}
 We will also denote 
\begin{equation}
\partial_{tim}\mathcal{M}^{\text{\textgreek{t}}}\doteq\partial_{tim}\mathcal{M}\cap\{\bar{t}=\text{\textgreek{t}}\}
\end{equation}
and 
\begin{equation}
\partial\mathcal{M}_{\text{\textgreek{t}}}\doteq\mathcal{H}_{\text{\textgreek{t}}}\cup\partial_{tim}\mathcal{M}^{\text{\textgreek{t}}}.
\end{equation}

Without loss of generality, we also assume that the function $r$
has been extended in such a way in the region $\{r\lesssim1\}$ so
that $dr\neq0$ on $\mathcal{H}\cap\{\bar{t}\ge0\}$. We will also
use the shorthand notation 
\begin{equation}
r_{+}\doteq\big(1+r^{2}\big)^{1/2}.
\end{equation}

Finally, we will also need to assume that the deformation tensor of
$N$ and its derivatives are bounded on $\{\bar{t}\ge0\}$ when measured
with the reference Riemannian metric (\ref{eq:ReferenceRiemannian metric})
(that will be constructed in a moment):

\begin{MyDescription}[leftmargin=3.0em]

\item [(G8)]{\label{enu:DeformationTensor}For any $l\in\mathbb{N}$,
there exists a $C_{l}>0$ such that 
\begin{equation}
\sup_{\{\bar{t}\ge0\}}\sum_{j=0}^{l}\big|\nabla_{g}^{l-j}\mathcal{L}_{N}^{j}g\big|_{h}\le C_{l}.\label{eq:BoundednessDeformationTensor}
\end{equation}
}

\end{MyDescription}
\begin{rem*}
Assumption \ref{enu:DeformationTensor} holds in the case when the
spacetime $\mathcal{M}$ is near stationary or time periodic. Moreover,
this is an assumption regarding the structure of the foliation in
the region $\{r\lesssim1\}$.
\end{rem*}
It will be convenient to fix a vector field $Y$ in a neighborhood
of $\mathcal{H}^{+}\cup\partial_{tim}\mathcal{M}$ so that for any
$\text{\textgreek{t}}\ge0$, $Y$ is tangent to $\{\bar{t}=\text{\textgreek{t}}\}$,
orthogonal to $\mathcal{H}_{\text{\textgreek{t}}}$ and $\partial_{tim}\mathcal{M}^{\text{\textgreek{t}}}$
and satisfies $g(Y,Y)|_{\mathcal{H}_{\text{\textgreek{t}}}}=1$ and
$g(Y,Y)|_{\partial_{tim}\mathcal{M}^{\text{\textgreek{t}}}}=1$.

\paragraph*{Boundary conditions on $\partial_{tim}\mathcal{M}$ for $\square\text{\textgreek{f}}=F$}

Assumptions \ref{enu:DomainDependence} and \ref{enu:Spacelike} guarantee
that for any $\text{\textgreek{t}}\ge0$, we can solve the inhomogeneous
wave equation 
\begin{equation}
\square_{g}\text{\textgreek{f}}=F\label{eq:InhomogeneousWaveEquation}
\end{equation}
 on $J^{+}(\{\bar{t}=\text{\textgreek{t}}\})$ with Cauchy initial
data on $\{\bar{t}=\text{\textgreek{t}}\}$, provided suitable boundary
conditions (e.\,g.~Dirichlet conditions) have been imposed on $\partial_{tim}\mathcal{M}$.
In particular, we will introduce the following definition:
\begin{defn*}
We define the class of \emph{admissible boundary conditions on }$\partial_{tim}\mathcal{M}$
to be the set $\mathcal{C}_{adm}$ of all families of linear functions
\begin{equation}
\mathcal{F}_{\text{\textgreek{t}}}:C^{\infty}(\partial_{tim}\mathcal{M}^{\text{\textgreek{t}}})\times C^{\infty}(\partial_{tim}\mathcal{M}^{\text{\textgreek{t}}})\rightarrow C^{\infty}(\partial_{tim}\mathcal{M}^{\text{\textgreek{t}}})\label{eq:BoundaryFunctional}
\end{equation}
 depending smoothly on $\text{\textgreek{t}}\ge0$ such that for any
$\text{\textgreek{t}}_{0}\ge0$, any $F\in C^{\infty}(\{\bar{t}\ge\text{\textgreek{t}}_{0}\})$
and any $\text{\textgreek{f}}_{0},\text{\textgreek{f}}_{1}\in C^{\infty}(\{\bar{t}=\text{\textgreek{t}}_{0}\})$,
the inititial-boundary value problem 
\begin{equation}
\begin{cases}
\square\text{\textgreek{f}}=F & \mbox{ on }\{\bar{t}\ge\text{\textgreek{t}}_{0}\}\\
(\text{\textgreek{f}}|_{\bar{t}=\text{\textgreek{t}}_{0}},N\text{\textgreek{f}}|_{\bar{t}=\text{\textgreek{t}}_{0}})=(\text{\textgreek{f}}_{0},\text{\textgreek{f}}_{1})\\
\mathcal{F}_{\text{\textgreek{t}}}(\text{\textgreek{f}}|_{\partial_{tim}\mathcal{M}},Y\text{\textgreek{f}}|_{\partial_{tim}\mathcal{M}})=0 & \mbox{ for }\text{\textgreek{t}}\ge\text{\textgreek{t}}_{0}
\end{cases}
\end{equation}
is well posed.
\end{defn*}
Notice that the usual Dirichlet and Neumann boundary condition belong
to the class $\mathcal{C}_{adm}$, corresponding to $\mathcal{F}_{\text{\textgreek{t}}}=Id\oplus0$
and $\mathcal{F}_{\text{\textgreek{t}}}=0\oplus Id$ respectively.

\paragraph*{Construction of the reference Riemannian metric $h$ on $\mathcal{M}$}

In Section \ref{sec:RiemannianMetric} of the Appendix we establish
the existence of a natural Riemannian metric $h_{\text{\textgreek{t}},N}$
defined on the hypersurfaces $\{\bar{t}=\text{\textgreek{t}}\}$ (wich
we will sometimes denote with $h_{N}$ for simplicity), associated
to $g$ and $N$ (and distinct from the induced metric $g_{\bar{t}}$,
which degenerates as on approaches $\mathcal{I}^{+}$). 
\begin{rem*}
Notice that $h_{N}$ is non singular up to $\mathcal{H}_{\text{\textgreek{t}}}$,
since $N$ is timelike everywhere up to $\mathcal{H}^{+}$. Thus,
in the language of Section \ref{sec:Elliptic-estimates} of the Appendix,
$h_{N}$ corresponds to the $\tilde{h}$ metric of that section.
\end{rem*}
We will extend $h_{N}$ to a Riemannian metric $h$ on $\mathcal{M}$
by setting 
\begin{equation}
h\doteq(d\bar{t})^{2}+h_{N}.\label{eq:ReferenceRiemannian metric}
\end{equation}
This Riemannian metric will be used to measure the norms of tensors
on $\mathcal{M}$. Moreover, we will denote with $h_{\mathcal{H}_{\text{\textgreek{t}}}}$
the Riemannian metric induced by $h$ on $\mathcal{H}_{\text{\textgreek{t}}}$.

Due to the expression (\ref{eq:MetricUV}) for $g$ in the region
$\{r\gg1\}$ and (\ref{eq:SplittingWaveOperator}) for $\square_{g}$,
we can bound for any smooth $\text{\textgreek{f}}:\mathcal{M}\rightarrow\mathbb{C}$
and any $l\in\mathbb{N}$: 
\begin{align}
\big|\nabla_{h_{\text{\textgreek{t}},N}}^{l}(\text{\textgreek{D}}_{h_{\text{\textgreek{t}},N},N}\text{\textgreek{f}})\big|_{h_{\text{\textgreek{t}},N}}^{2}\le C(\text{\textgreek{t}})\Big\{ & \sum_{j=1}^{2}r_{+}^{-2}|\nabla_{h_{\text{\textgreek{t}},N}}^{l}(N^{j}\text{\textgreek{f}})|_{h_{\text{\textgreek{t}},N}}^{2}+\sum_{i=0}^{l}\big(|\nabla_{h_{\text{\textgreek{t}},N}}^{i+1}(N\text{\textgreek{f}})|_{h_{\text{\textgreek{t}},N}}^{2}+|\nabla_{h_{\text{\textgreek{t}},N}}^{i+1}\text{\textgreek{f}}|_{h_{\text{\textgreek{t}},N}}^{2}\big)+\label{eq:UsingTheWaveEquation}\\
 & +\big|\nabla_{h_{\text{\textgreek{t}},N}}^{l}(\square\text{\textgreek{f}})\big|_{h_{\text{\textgreek{t}},N}}^{2}\Big\},\nonumber 
\end{align}
where the operator $\text{\textgreek{D}}_{h_{\text{\textgreek{t}},N},N}$
on $\{\bar{t}=\text{\textgreek{t}}\}$ is defined as: 
\begin{equation}
\text{\textgreek{D}}_{h_{\text{\textgreek{t}},N},N}=\frac{1}{\sqrt{-g(N,N)}}div_{h_{\text{\textgreek{t}},N}}\big(\sqrt{-g(N,N)}d\big).\label{eq:PerturbedLaplacian-2}
\end{equation}
 Notice that in the right hand side of (\ref{eq:UsingTheWaveEquation})
there is no term of order $l+2$ in the spatial derivatives (i.\,e.~$\nabla_{h_{\text{\textgreek{t}},\text{\textgreek{N}}}}$).

\paragraph*{Assumptions on the uniformity of elliptic, Poincare and Sobolev type
estimates on the leaves of the foliation $\{\bar{t}=\text{\textgreek{t}}\}$}

We will also need to ensure that we can establish elliptic, Poincare
and Sobolev type estimates on the leaves of the foliation $\{\bar{t}=\text{\textgreek{t}}\}$
with constants that do not depend on $\text{\textgreek{t}}$. We will
assume the following uniformity condition on $Y$:

\begin{MyDescription}[leftmargin=3.0em]

\item [(G9)]{\label{enu:UniformityY}For any $l\in\mathbb{N}$ the
following uniform bound holds: 
\begin{equation}
\sup_{\{\bar{t}\ge0\}}\big|\nabla_{g}^{l}Y\big|_{h}\le C_{l}.
\end{equation}
}

\end{MyDescription}

According to Proposition \ref{Prop:NonDegenerateEllipticEstimates}
of the Appendix and the estimate (\ref{eq:UsingTheWaveEquation}),
without imposing any extra assumptions the following statement holds: 
\begin{lem*}
\label{lem:NonDegenerateEllipticEstimates-1}For any integer $l\ge2$
and any $\text{\textgreek{b}}\in[0,1)$ we can bound for any $\text{\textgreek{t}}\ge0$
and any $\text{\textgreek{f}}\in C^{\infty}(\mathcal{M})$ satisfying
for any $j_{1}+j_{2}\le l$ the finite radiation field condition $\limsup_{r\rightarrow+\infty}\big|r^{\frac{d-1}{2}+j_{1}}\nabla_{h_{\text{\textgreek{t}},N}}^{j_{1}}(N^{j_{2}}\text{\textgreek{f}})\big|_{h_{\text{\textgreek{t}},N}}<+\infty$:
\begin{equation}
\begin{split}\int_{\{\bar{t}=\text{\textgreek{t}}\}}r_{+}^{-\text{\textgreek{b}}}|\nabla_{g}^{l}\text{\textgreek{f}}|_{h}^{2}\, & dh_{N}\le C_{\text{\textgreek{b}},n_{0}}(\text{\textgreek{t}})\Bigg\{\int_{\{\bar{t}=\text{\textgreek{t}}\}}r_{+}^{-\text{\textgreek{b}}}\Big\{\sum_{j=1}^{l}\sum_{i=0}^{l-j}r_{+}^{-2}|\nabla_{h_{\text{\textgreek{t}},N}}^{i}(N^{j}\text{\textgreek{f}})|_{h_{\text{\textgreek{t}},N}}^{2}+\\
 & +\sum_{0\le j_{1}+j_{2}\le l-2}\big(|\nabla_{h_{\text{\textgreek{t}},N}}^{j_{1}+1}(N^{j_{2}+1}\text{\textgreek{f}})|_{h_{\text{\textgreek{t}},N}}^{2}+|\nabla_{h_{\text{\textgreek{t}},N}}^{j_{1}+1}(N^{j_{2}}\text{\textgreek{f}})|_{h_{\text{\textgreek{t}},N}}^{2}+\big|\nabla_{h_{\text{\textgreek{t}},N}}^{j_{1}}(N^{j_{2}}\square\text{\textgreek{f}})\big|_{h_{\text{\textgreek{t}},N}}^{2}\big)\Big\}\, dh_{N}+\\
 & +\sum_{j=0}^{l-1}\Big|\int_{\partial\mathcal{M}_{\text{\textgreek{t}}}}h_{\partial\mathcal{M}_{\text{\textgreek{t}}}}\big(\nabla_{h_{\partial\mathcal{M}_{\text{\textgreek{t}}}}}^{j}(Yu),\nabla_{h_{\partial\mathcal{M}_{\text{\textgreek{t}}}}}^{j}u\big)\, dh_{\partial\mathcal{M}_{\text{\textgreek{t}}}}\Big|\Bigg\}.
\end{split}
\label{eq:EllipticEstimatesNonDegenerate-1}
\end{equation}
 
\end{lem*}
Our final assumtions on the geometry of $(\mathcal{M},g)$ in the
region $\{\bar{t}\ge0\}$ will be the following:

\begin{MyDescription}[leftmargin=3.0em]

\item [(G10)]{\label{enu:EllipticEstimates}The constants appearing
in the right hand side of (\ref{eq:EllipticEstimatesNonDegenerate-1})
do \underline{not} depend on $\text{\textgreek{t}}\ge0$.}

\item [(G11)]{\label{enu:EllipticEstimatesHorizon}The following elliptic
estimate holds for any $l\in\mathbb{N}$ and any smooth function $u$
on the submanifolds $\partial\mathcal{M}_{\text{\textgreek{t}}}$
with a constant $C_{l}$ not depending on $\text{\textgreek{t}}\ge0$:
\begin{equation}
\sum_{j=1}^{l}\int_{\partial\mathcal{M}_{\text{\textgreek{t}}}}\big|\nabla_{h_{\partial\mathcal{M}_{\text{\textgreek{t}}}}}^{j}u\big|_{h_{\partial\mathcal{M}_{\text{\textgreek{t}}}}}^{2}\, dh_{\partial\mathcal{M}_{\text{\textgreek{t}}}}\le C_{l}\cdot\int_{\partial\mathcal{M}_{\text{\textgreek{t}}}}\big|\nabla_{h_{\partial\mathcal{M}_{\text{\textgreek{t}}}}}^{l-2}(\text{\textgreek{D}}_{h_{\partial\mathcal{M}_{\text{\textgreek{t}}}}}u)\big|_{h_{\partial\mathcal{M}_{\text{\textgreek{t}}}}}^{2}\, dh_{\partial\mathcal{M}_{\text{\textgreek{t}}}}.\label{eq:EllipticEstimatesHorizon}
\end{equation}
}

\item [(G12)] {\label{enu:PoincareInequalities}For $C^{1}$ functions
$u$ on $\mathcal{M}$, the Poincare inequality 
\begin{equation}
\int_{\{\bar{t}=\text{\textgreek{t}}\}\cap\{r\le R\}}|u|^{2}\, dh_{N}\le C(R)\cdot\Big(\int_{\{\bar{t}=\text{\textgreek{t}}\}\cap\{r\le R\}}|\nabla_{h_{N}}u|_{h_{N}}^{2}\, dh_{N}+\int_{\{\bar{t}=\text{\textgreek{t}}\}\cap\{R\le r\le2R\}}|u|^{2}\, dh_{N}\Big)
\end{equation}
and the trace inequality 
\begin{equation}
\int_{\partial\mathcal{M}_{\text{\textgreek{t}}}}|(\nabla_{h_{\partial\mathcal{M}_{\text{\textgreek{t}}}}})^{1/2}u|^{2}\, dh_{\partial\mathcal{M}_{\text{\textgreek{t}}}}\le C(\text{\textgreek{e}})\cdot\Big(\int_{\partial\mathcal{M}\cap\{\text{\textgreek{t}}-\text{\textgreek{e}}\le\bar{t}\le\text{\textgreek{t}}\}}|\nabla_{h_{\partial\mathcal{M}}}u|_{h_{\partial\mathcal{M}}}^{2}\, dh_{\partial\mathcal{M}}+\int_{\partial\mathcal{M}\cap\{\text{\textgreek{t}}-\text{\textgreek{e}}\le\bar{t}\le\text{\textgreek{t}}\}}|u|^{2}\, dh_{\partial\mathcal{M}}\Big)
\end{equation}
(where $h_{\partial\mathcal{M}}$ is the Riemannian metric on $\partial\mathcal{M}$
induced by $h$) hold for constants $C(R)$ and $C(\text{\textgreek{e}})$
that do \underline{not} depend on $\text{\textgreek{t}}\ge0$. In
the above}

\item [(G13)] {\label{enu:SobolevInequality}The following Sobolev
inequality holds for smooth and compactly supported functions $u$
on the hypersurfaces $\{\bar{t}=\text{\textgreek{t}}\}_{\text{\textgreek{t}}\ge0}$
\[
\sup_{\{\bar{t}=\text{\textgreek{t}}\}}|u|^{2}\le C\cdot\sum_{j=0}^{\lceil\frac{d+1}{2}\rceil}\int_{\{\bar{t}=\text{\textgreek{t}}\}}|\nabla_{h_{N}}^{j}u|_{h_{N}}^{2}\, dh_{N}
\]
 (see \cite{Hebey1999}) with the constant $C$ in the right hand
side independent of $\text{\textgreek{t}}$.}

\end{MyDescription}
\begin{rem*}
Assumptions \ref{enu:EllipticEstimates}-\ref{enu:SobolevInequality}
are automatically satisfied in the case the spacetime $\mathcal{M}$
is near stationary or time periodic. Moreover, these assumptions are
only tied to the structure of the foliation in the region $\{r\lesssim1\}$.
\end{rem*}

\subsubsection{\label{sub:IntegratedLocalEnergyDecay}Integrated local energy decay
statement on $(\mathcal{M},g)$}

We assume that the following integrated local energy decay statement
holds on the spacetime $(\mathcal{M},g)$ under consideration:

\begin{MyDescription}[leftmargin=3.0em]

\item [(ILED1)]{\label{enu:IntegratedLocalEnergyDecay}\emph{Integrated
local energy decay with loss of derivatives: }We assume that there
exists a (non-empty) class $\mathcal{C}_{ILED}$ of boundary conditions
on $\partial_{tim}\mathcal{M}$, which is contained in the class of
admissible boundary conditions $\mathcal{C}_{adm}$, so that the following
statement holds: There exists an integer $k\ge0$ such that for any
$R,R_{f}>0$, any integer $m\ge0$, any $0<\text{\textgreek{h}}<a$,
any smooth $\text{\textgreek{f}}:\mathcal{M}\rightarrow\mathbb{C}$
solving $\square_{g}\text{\textgreek{f}}=F$ satisfying boundary conditions
on $\partial_{tim}\mathcal{M}$ belonging to the class $\mathcal{C}_{ILED}$
and any $0\le\text{\textgreek{t}}_{1}\le\text{\textgreek{t}}_{2}$
we can bound: 
\begin{equation}
\begin{split}\sum_{j=1}^{m} & \int_{\{\text{\textgreek{t}}_{1}\le\bar{t}\le\text{\textgreek{t}}_{2}\}\cap\{r\le R\}}\Big(|\nabla_{g}^{j}\text{\textgreek{f}}|_{h}^{2}+|\text{\textgreek{f}}|^{2}\Big)+\sum_{j=1}^{m}\int_{\{\text{\textgreek{t}}_{1}\le\bar{t}\le\text{\textgreek{t}}_{2}\}\cap\partial_{tim}\mathcal{M}}\Big(|\nabla_{g}^{j}\text{\textgreek{f}}|_{h}^{2}+|\text{\textgreek{f}}|^{2}\Big)\le\\
 & \hphantom{\int_{\{\text{\textgreek{t}}_{1}\le}}\le C_{m,\text{\textgreek{h}}}(R,R_{f})\sum_{j=0}^{m+k-1}\Big(\int_{\{\bar{t}=\text{\textgreek{t}}_{1}\}}J_{\text{\textgreek{m}}}^{N}(N^{j}\text{\textgreek{f}})\bar{n}^{\text{\textgreek{m}}}+\int_{\{\text{\textgreek{t}}_{1}\le\bar{t}\le\text{\textgreek{t}}_{2}\}}r_{+}^{1+\text{\textgreek{h}}}|\nabla_{g}^{j}F|_{h}^{2}\, dh\Big)+\\
 & \hphantom{\int_{\{\text{\textgreek{t}}_{1}\le}\le}+C_{m,\text{\textgreek{h}}}\sum_{j=0}^{m-1}\sum_{j_{1}+j_{2}=j}\int_{\{\text{\textgreek{t}}_{1}\le\bar{t}\le\text{\textgreek{t}}_{1}\}\cap\{r\ge R_{f}\}}r_{+}^{-1}\Big(\big|\nabla_{h_{\text{\textgreek{t}},N}}^{j_{1}+1}(N^{j_{2}}\text{\textgreek{f}})\big|^{2}+r_{+}^{-2}\big|\nabla_{h_{\text{\textgreek{t}},N}}^{j_{1}}(N^{j_{2}}\text{\textgreek{f}})\big|^{2}+r_{+}^{-2}\big|N^{j+1}\text{\textgreek{f}}\big|^{2}\Big)\, dg.
\end{split}
\label{eq:IntegratedLocalEnergyDecayImprovedDecay}
\end{equation}
 }

\end{MyDescription}

Using Lemma \ref{lem:MorawetzDrLemmaHyperboloids} and \ref{lem:Boundedness Hyperboloids},
as well as a trace theorem for $\int_{\partial_{tim}\mathcal{M}}|\text{\textgreek{f}}|^{2}$,
(\ref{eq:IntegratedLocalEnergyDecayImprovedDecay}) can be improved
into the statement that for any $0<\text{\textgreek{h}}<a$ and $R\gg1$:
\begin{equation}
\begin{split}\int_{\{\text{\textgreek{t}}_{1}\le\bar{t}\le\text{\textgreek{t}}_{2}\}} & r_{+}^{-1-\text{\textgreek{h}}}\big(\sum_{j=1}^{m}|\nabla_{g}^{j}\text{\textgreek{f}}|_{h}^{2}+r^{-2}|\text{\textgreek{f}}|^{2}\big)+\sum_{j=1}^{m}\int_{\{\text{\textgreek{t}}_{1}\le\bar{t}\le\text{\textgreek{t}}_{2}\}\cap\partial_{tim}\mathcal{M}}\Big(|\nabla_{g}^{j}\text{\textgreek{f}}|_{h}^{2}+|\text{\textgreek{f}}|^{2}\Big)\le\\
\le & C_{m,\text{\textgreek{h}}}(R)\sum_{j=0}^{m+k-1}\Big(\int_{\{\bar{t}=\text{\textgreek{t}}_{1}\}}J_{\text{\textgreek{m}}}^{N}(N^{j}\text{\textgreek{f}})\bar{n}^{\text{\textgreek{m}}}+\int_{\{\text{\textgreek{t}}_{1}\le\bar{t}\le\text{\textgreek{t}}_{2}\}}r_{+}^{1+\text{\textgreek{h}}}|\nabla_{g}^{j}F|_{h}^{2}\, dh\Big)+\\
 & +C_{m,\text{\textgreek{h}}}\cdot\sum_{j=0}^{m-1}\sum_{j_{1}+j_{2}=j}\int_{\{\text{\textgreek{t}}_{1}\le\bar{t}\le\text{\textgreek{t}}_{1}\}\cap\{r\ge R\}}r_{+}^{-1}\Big(\big|\nabla_{h_{\text{\textgreek{t}},N}}^{j_{1}+1}(N^{j_{2}}\text{\textgreek{f}})\big|^{2}+r_{+}^{-2}\big|\nabla_{h_{\text{\textgreek{t}},N}}^{j_{1}}(N^{j_{2}}\text{\textgreek{f}})\big|^{2}+r_{+}^{-2}\big|N^{j+1}\text{\textgreek{f}}\big|^{2}\Big)\, dg.
\end{split}
\label{eq:IntegratedLocalEnergyDecayImprovedDecayInfinity-1}
\end{equation}
 We should also notice that the results of this section can be readily
established in case one replaces Assumption \ref{enu:IntegratedLocalEnergyDecay}
by the following pair of integrated local energy decay statements: 
\begin{description}
\item [{\emph{Alternative~integrated~local~energy~decay~statement:}}] \noindent With
the notations as in Assumption \ref{enu:IntegratedLocalEnergyDecay},
we assume that there exists an integer $k\ge0$ and an $R_{c}>0$
such that for any $R,R_{f}>0$, any integer $m\ge0$, any $0<\text{\textgreek{h}}<a$,
any smooth $\text{\textgreek{f}}:\mathcal{M}\rightarrow\mathbb{C}$
solving $\square_{g}\text{\textgreek{f}}=F$ satisfying boundary conditions
on $\partial_{tim}\mathcal{M}$ belonging to the class $\mathcal{C}_{ILED}$
and any $0\le\text{\textgreek{t}}_{1}\le\text{\textgreek{t}}_{2}$
we can bound: 
\begin{equation}
\begin{split}\sum_{j=1}^{m} & \int_{\{\text{\textgreek{t}}_{1}\le\bar{t}\le\text{\textgreek{t}}_{2}\}\cap\{r\le R\}}\Big(|\nabla_{g}^{j}\text{\textgreek{f}}|_{h}^{2}+|\text{\textgreek{f}}|^{2}\Big)+\sum_{j=1}^{m}\int_{\{\text{\textgreek{t}}_{1}\le\bar{t}\le\text{\textgreek{t}}_{2}\}\cap\partial_{tim}\mathcal{M}}\Big(|\nabla_{g}^{j}\text{\textgreek{f}}|_{h}^{2}+|\text{\textgreek{f}}|^{2}\Big)\le\\
 & \hphantom{\int_{\{\text{\textgreek{t}}_{1}\le}}\le C_{m,\text{\textgreek{h}}}(R,R_{f})\sum_{j=0}^{m+k-1}\Big(\int_{\{\bar{t}=\text{\textgreek{t}}_{1}\}}J_{\text{\textgreek{m}}}^{N}(N^{j}\text{\textgreek{f}})\bar{n}^{\text{\textgreek{m}}}+\int_{\{\text{\textgreek{t}}_{1}\le\bar{t}\le\text{\textgreek{t}}_{2}\}}r_{+}^{1+\text{\textgreek{h}}}|\nabla_{g}^{j}F|_{h}^{2}\, dh\Big)+\\
 & \hphantom{\int_{\{\text{\textgreek{t}}_{1}\le}\le}+C_{m,\text{\textgreek{h}}}\sum_{j=0}^{m+k-1}\sum_{j_{1}+j_{2}=j}\int_{\{\text{\textgreek{t}}_{1}\le\bar{t}\le\text{\textgreek{t}}_{1}\}\cap\{r\ge R_{f}\}}r_{+}^{-1}\Big(\big|\nabla_{h_{\text{\textgreek{t}},N}}^{j_{1}+1}(N^{j_{2}}\text{\textgreek{f}})\big|^{2}+r_{+}^{-2}\big|\nabla_{h_{\text{\textgreek{t}},N}}^{j_{1}}(N^{j_{2}}\text{\textgreek{f}})\big|^{2}+r_{+}^{-2}\big|N^{j+1}\text{\textgreek{f}}\big|^{2}\Big)\, dg
\end{split}
\label{eq:IntegratedLocalEnergyDecayImprovedDecayAltLoss}
\end{equation}
and 
\begin{equation}
\begin{split}\sum_{j=1}^{m} & \int_{\{\text{\textgreek{t}}_{1}\le\bar{t}\le\text{\textgreek{t}}_{2}\}\cap\{R_{c}\le r\le R\}}\Big(|\nabla_{g}^{j}\text{\textgreek{f}}|_{h}^{2}+|\text{\textgreek{f}}|^{2}\Big)+\sum_{j=1}^{m}\int_{\{\text{\textgreek{t}}_{1}\le\bar{t}\le\text{\textgreek{t}}_{2}\}\cap\partial_{tim}\mathcal{M}}\Big(|\nabla_{g}^{j}\text{\textgreek{f}}|_{h}^{2}+|\text{\textgreek{f}}|^{2}\Big)\le\\
 & \hphantom{\int_{\{\text{\textgreek{t}}_{1}\le}}\le C_{m,\text{\textgreek{h}}}(R,R_{f})\sum_{j=0}^{m-1}\Big(\int_{\{\bar{t}=\text{\textgreek{t}}_{1}\}}J_{\text{\textgreek{m}}}^{N}(N^{j}\text{\textgreek{f}})\bar{n}^{\text{\textgreek{m}}}+\int_{\{\text{\textgreek{t}}_{1}\le\bar{t}\le\text{\textgreek{t}}_{2}\}}r_{+}^{1+\text{\textgreek{h}}}|\nabla_{g}^{j}F|_{h}^{2}\, dh\Big)+\\
 & \hphantom{\int_{\{\text{\textgreek{t}}_{1}\le}\le}+C_{m,\text{\textgreek{h}}}\sum_{j=0}^{m-1}\sum_{j_{1}+j_{2}=j}\int_{\{\text{\textgreek{t}}_{1}\le\bar{t}\le\text{\textgreek{t}}_{1}\}\cap\{r\ge R_{f}\}}r_{+}^{-1}\Big(\big|\nabla_{h_{\text{\textgreek{t}},N}}^{j_{1}+1}(N^{j_{2}}\text{\textgreek{f}})\big|^{2}+r_{+}^{-2}\big|\nabla_{h_{\text{\textgreek{t}},N}}^{j_{1}}(N^{j_{2}}\text{\textgreek{f}})\big|^{2}+r_{+}^{-2}\big|N^{j+1}\text{\textgreek{f}}\big|^{2}\Big)\, dg.
\end{split}
\label{eq:IntegratedLocalEnergyDecayImprovedDecayAlt}
\end{equation}

\end{description}
However, we will not pursue this issue again in the paper.

\subsection{\label{sub:lFirstPolynomialDecay}A first polynomial decay result}

On any spacetime $(\mathcal{M},g)$ satisfying the geometric assumptions
\ref{enu:AsymptoticFlatness}-\ref{enu:SobolevInequality} and Assumption
\ref{enu:IntegratedLocalEnergyDecay} on integrated local energy decay
with loss of derivatives, we will establish the following polynomial
decay estimates:
\begin{thm}
\label{thm:FirstPointwiseDecayNewMethod}Let $(\mathcal{M}^{d+1},g)$,
$d\ge3$, satisfy Assumptions \ref{enu:AsymptoticFlatness}-\ref{enu:SobolevInequality}
and \ref{enu:IntegratedLocalEnergyDecay}. For any smooth solution
$\text{\textgreek{f}}$ to $\square_{g}\text{\textgreek{f}}=F$ on
$J^{+}(\{\bar{t}=0\})$ with suitably decaying initial data on $\{\bar{t}=0\}$
(and satisfying boundary conditions on $\partial_{tim}\mathcal{M}$
belonging to the class $\mathcal{C}_{ILED}$), the following decay
estimates hold for any $\text{\textgreek{t}}\ge0$, any integer $m\ge0$
and any $\text{\textgreek{e}}>0$, $0<\text{\textgreek{h}}<a$, provided
(\ref{eq:NormFinitenessFriedlanderRadiation}) holds for all $0\le j\le m+1+d+3k$:
\begin{equation}
\sum_{j=0}^{m-1}\sum_{j_{1}+j_{2}=j}\int_{\{\bar{t}=\text{\textgreek{t}}\}}\Big(\big|\nabla_{h_{\text{\textgreek{t}},N}}^{j_{1}+1}(N^{j_{2}}\text{\textgreek{f}})\big|^{2}+r_{+}^{-2}\big|N^{j+1}\text{\textgreek{f}}\big|^{2}+r_{+}^{-2}|\text{\textgreek{f}}|^{2}\Big)\, dh_{N}\lesssim_{m,\text{\textgreek{e}},\text{\textgreek{h}}}\text{\textgreek{t}}^{-2+\text{\textgreek{e}}}\mathcal{E}_{bound}^{(2,m+3k)}[\text{\textgreek{f}}](0)+\mathcal{F}_{\text{\textgreek{h}},\text{\textgreek{e}}}^{(2,m,k)}[F](\text{\textgreek{t}}),\label{eq:FinalDecayFirstenergy}
\end{equation}
\begin{equation}
\sup_{\{\bar{t}=\text{\textgreek{t}}\}}\big|r_{+}^{\frac{d-2}{2}}\nabla_{g}^{m}\text{\textgreek{f}}\big|_{h}^{2}\lesssim_{m,\text{\textgreek{e}},\text{\textgreek{h}}}\text{\textgreek{t}}^{-2+\text{\textgreek{e}}}\mathcal{E}_{bound}^{(2,m+\lceil\frac{d+2}{2}\rceil+3k)}[\text{\textgreek{f}}](0)+\mathcal{F}_{\text{\textgreek{h}},\text{\textgreek{e}}}^{(2,m+\lceil\frac{d+2}{2}\rceil,k)}[F](\text{\textgreek{t}})\label{eq:PointwiseDecayFirst-1}
\end{equation}
and 
\begin{align}
\sup_{\{\bar{t}=\text{\textgreek{t}}\}}\big|r_{+}^{\frac{d-1}{2}}\nabla_{g}^{m}\text{\textgreek{f}}\big|_{h}^{2}\lesssim_{m,\text{\textgreek{h}}} & \text{\textgreek{t}}^{-1}\mathcal{E}_{bound}^{(2,m+\lceil\frac{d+2}{2}\rceil+2k)}[\text{\textgreek{f}}](0)+\mathcal{F}_{\text{\textgreek{h}}}^{(1,m+\lceil\frac{d+2}{2}\rceil,k)}[F](\text{\textgreek{t}}).\label{eq:PointwiseDecayRadiationField-1}
\end{align}
In the above, $k$ is the integer measuring the derivative loss in
the integrated local energy decay statement (\ref{eq:IntegratedLocalEnergyDecayImprovedDecay})
and for some fixed $R,C>1$: 
\begin{equation}
\mathcal{E}_{bound}^{(p,m)}[\text{\textgreek{f}}](\text{\textgreek{t}})=\sum_{j=0}^{m}\int_{\{\bar{t}=\text{\textgreek{t}}\}\cap\{r\le R+1\}}\big|\nabla_{h_{\text{\textgreek{t}},N}}^{j_{1}}(N^{j_{2}}\text{\textgreek{f}})\big|_{h_{\text{\textgreek{t}},N}}^{2}\, dh_{N}+\sum_{\substack{components\\
of\,\mathcal{N}_{af,\mathcal{M}}
}
}\sum_{j=0}^{m-1}\sum_{j_{1}+j_{2}+j_{3}=j}\mathcal{E}_{bound,R;\text{\textgreek{d}}}^{(p,1)}[r^{-j_{2}}\partial_{v}^{j_{1}}\partial_{\text{\textgreek{sv}}}^{j_{2}}\partial_{u}^{j_{3}}\text{\textgreek{f}}](\text{\textgreek{t}})
\end{equation}
(where the derivatives in the second term are with respect to the
$(u,v,\text{\textgreek{sv}})$ coordinate charts over each connected
component of $\mathcal{N}_{af,\mathcal{M}}$),
\begin{align}
\mathcal{F}_{\text{\textgreek{h}},\text{\textgreek{e}}}^{(2,m,k)}[F](\text{\textgreek{t}})= & \text{\textgreek{t}}^{-2+\text{\textgreek{e}}}\sum_{j=0}^{m+3k-1}\int_{\{0\le\bar{t}\le\text{\textgreek{t}}\}}r^{3}\big|\nabla_{g}^{j}F\big|_{h}^{2}\, dg+\text{\textgreek{t}}^{-1+\text{\textgreek{e}}}\sum_{j=0}^{m+2k-1}\int_{\{C^{-1}\text{\textgreek{t}}\le\bar{t}\le\text{\textgreek{t}}\}}r^{2}\big|\nabla_{g}^{j}F\big|_{h}^{2}\, dg+\\
 & +\sum_{j=0}^{m+2k-1}\int_{\{C^{-1}\text{\textgreek{t}}\le\bar{t}\le\text{\textgreek{t}}\}}(r^{1+\text{\textgreek{e}}}+r^{1+\text{\textgreek{h}}})\big|\nabla_{g}^{j}F\big|_{h}^{2}\, dg\nonumber 
\end{align}
and 
\begin{equation}
\mathcal{F}_{\text{\textgreek{h}}}^{(1,m,k)}[F](\text{\textgreek{t}})=\text{\textgreek{t}}^{-1}\sum_{j=0}^{m+2k-1}\int_{\{0\le\bar{t}\le\text{\textgreek{t}}\}}r^{3}\big|\nabla_{g}^{j}F\big|_{h}^{2}\, dg+\sum_{j=0}^{m+2k-1}\int_{\{C^{-1}\text{\textgreek{t}}\le\bar{t}\le\text{\textgreek{t}}\}}r^{2}\big|\nabla_{g}^{j}F\big|_{h}^{2}\, dg.
\end{equation}
\end{thm}
\begin{rem*}
In case there exists some small $\text{\textgreek{d}}_{0}>0$ such
that the deformation tensor of the vector field $T$ in the region
$\{r\gg1\}$ satisfies the following bound for any integer $m\ge1$:
\begin{align}
\mathcal{L}_{T}^{m}g=O(\bar{t}{}^{-\text{\textgreek{d}}_{0}})\Big\{ & O(r^{-1-a})dvdu+O(r)d\text{\textgreek{sv}}d\text{\textgreek{sv}}+O(1)dud\text{\textgreek{sv}}+\label{eq:DeformationTensorTAwaySlow}\\
 & +O(r^{-a})dvd\text{\textgreek{sv}}+O(r^{-1})du^{2}+O(r^{-2-a})dv^{2}\Big\}\nonumber 
\end{align}
and the last term in the right hand side of the integrated local energy
decay estimate (\ref{eq:IntegratedLocalEnergyDecayImprovedDecay})
is replaced by
\begin{equation}
C_{m,\text{\textgreek{h}}}\sum_{j=0}^{m-1}\sum_{j_{1}+j_{2}=j}\int_{\{\text{\textgreek{t}}_{1}\le\bar{t}\le\text{\textgreek{t}}_{1}\}\cap\{r\ge R_{f}\}}\bar{t}^{-\text{\textgreek{d}}_{0}}r_{+}^{-1}\Big(\big|\nabla_{h_{\text{\textgreek{t}},N}}^{j_{1}+1}(N^{j_{2}}\text{\textgreek{f}})\big|^{2}+r_{+}^{-2}\big|\nabla_{h_{\text{\textgreek{t}},N}}^{j_{1}}(N^{j_{2}}\text{\textgreek{f}})\big|^{2}+r_{+}^{-2}\big|N^{j+1}\text{\textgreek{f}}\big|^{2}\Big)\, dg
\end{equation}
(this is consistent with the remark below Lemma \ref{lem:MorawetzDrLemmaHyperboloids}),
then the $\text{\textgreek{e}}$ loss in the exponent of $\text{\textgreek{t}}$
in (\ref{eq:FinalDecayFirstenergy}) and (\ref{eq:PointwiseDecayFirst-1})
can be removed. This follows readily the fact that, in this case,
the second summand of the right hand side of (\ref{eq:BoundednessGeneralRadiative})
comes with a $O(\bar{t}{}^{-\text{\textgreek{d}}_{0}})$ factor, which
enables us to deduce (\ref{eq:FinalDecayFirstenergy}) from (\ref{eq:FinalStepDyadic})
and (\ref{eq:FinalStepDyadic-1}) for some $\text{\textgreek{e}}\ll\text{\textgreek{d}}_{0}$.
\end{rem*}
The proof of Theorem \ref{thm:FirstPointwiseDecayNewMethod} will
be presented in Section \ref{sub:ProofOfFirstPointwiseDecay}.

We will also establish the following generalisation of Theorem \ref{thm:FirstPointwiseDecayNewMethod}
with improved weights in $r$ associated to higher derivatives of
$\text{\textgreek{f}}$:
\begin{thm}
\label{thm:SlowPointwiseDecayHighDerivativesNewMethod}Let $(\mathcal{M}^{d+1},g)$,
$d\ge3$, satisfy Assumptions \ref{enu:AsymptoticFlatness}-\ref{enu:SobolevInequality}
and \ref{enu:IntegratedLocalEnergyDecay}. For any smooth solution
$\text{\textgreek{f}}$ to $\square_{g}\text{\textgreek{f}}=F$ on
$J^{+}(\{\bar{t}=0\})\subset(\mathcal{M},g)$ with suitably decaying
initial data on $\{\bar{t}=0\}$ (and satisfying boundary conditions
on $\partial_{tim}\mathcal{M}$ belonging to the class $\mathcal{C}_{ILED}$),
the following decay estimates hold for any $\text{\textgreek{t}}\ge0$,
any integers $q,m\ge1$ and any $\text{\textgreek{e}}>0$, $0<\text{\textgreek{h}}<a$
provided (\ref{eq:NormFinitenessFriedlanderRadiation}) holds for
all $0\le j\le q+m+d+3k$:
\begin{equation}
\begin{split}\sum_{0\le i_{1}+i_{2}\le m-1}\sum_{j=0}^{q-1}\sum_{j_{1}+j_{2}=j}\int_{\{\bar{t}=\text{\textgreek{t}}\}}\Big(r^{2j_{1}}\big|\nabla_{h_{\text{\textgreek{t}},N}}^{i_{1}+j_{1}+1}(N^{i_{2}+j_{2}}\text{\textgreek{f}})\big|^{2}+r_{+}^{-2} & \big|N^{j+i_{1}+i_{2}+1}\text{\textgreek{f}}\big|^{2}+r_{+}^{-2}|\text{\textgreek{f}}|^{2}\Big)\, dh_{N}\lesssim_{m,q,\text{\textgreek{e}},\text{\textgreek{h}}}\\
\lesssim_{m,q,\text{\textgreek{e}},\text{\textgreek{h}}}\, & \text{\textgreek{t}}^{-2+\text{\textgreek{e}}}\mathcal{E}_{bound}^{(2q,q,m-1+3k)}[\text{\textgreek{f}}](0)+\mathcal{F}_{\text{\textgreek{h}},\text{\textgreek{e}}}^{(2,q,m-1+k)}[F](\text{\textgreek{t}}),
\end{split}
\label{eq:FinalDecayFirstenergy-1}
\end{equation}
\begin{equation}
\sum_{0\le i_{1}+i_{2}\le m-1}\sum_{j_{1}+j_{2}=q}\sup_{\{\bar{t}=\text{\textgreek{t}}\}}\big|r_{+}^{\frac{d-2}{2}+j_{1}}\nabla_{h_{\text{\textgreek{t}},N}}^{j_{1}+i_{1}}(N^{j_{2}+i_{2}}\text{\textgreek{f}})\big|_{h}^{2}\lesssim_{m,q,\text{\textgreek{e}},\text{\textgreek{h}}}\text{\textgreek{t}}^{-2+\text{\textgreek{e}}}\mathcal{E}_{bound}^{(2q,q,m-1+\lceil\frac{d+2}{2}\rceil+3k)}[\text{\textgreek{f}}](0)+\mathcal{F}_{\text{\textgreek{h}},\text{\textgreek{e}}}^{(2,q,m-1+\lceil\frac{d+2}{2}\rceil,k)}[F](\text{\textgreek{t}})\label{eq:PointwiseDecayFirst-1-1}
\end{equation}
and 
\begin{align}
\sum_{0\le i_{1}+i_{2}\le m-1}\sum_{j_{1}+j_{2}=q}\sup_{\{\bar{t}=\text{\textgreek{t}}\}}\big|r_{+}^{\frac{d-1}{2}+j_{1}}\nabla_{h_{\text{\textgreek{t}},N}}^{j_{1}+i_{1}}(N^{j_{2}+i_{2}}\text{\textgreek{f}})\big|_{h}^{2}\lesssim_{m,q,\text{\textgreek{h}}} & \text{\textgreek{t}}^{-1}\mathcal{E}_{bound}^{(2q,q,m-1+\lceil\frac{d+2}{2}\rceil+2k)}[\text{\textgreek{f}}](0)+\mathcal{F}_{\text{\textgreek{h}}}^{(1,q,m-1+\lceil\frac{d+2}{2}\rceil,k)}[F](\text{\textgreek{t}}).\label{eq:PointwiseDecayRadiationField-1-1}
\end{align}
In the above, $k$ is the integer measuring the derivative loss in
the integrated local energy decay statement (\ref{eq:IntegratedLocalEnergyDecayImprovedDecay})
and for some fixed $R,C>1$: 
\begin{equation}
\mathcal{E}_{bound}^{(p,q,m)}[\text{\textgreek{f}}](\text{\textgreek{t}})=\sum_{j=0}^{q+m-1}\int_{\{\bar{t}=\text{\textgreek{t}}\}\cap\{r\le R+1\}}\big|\nabla_{h_{\text{\textgreek{t}},N}}^{j_{1}}(N^{j_{2}}\text{\textgreek{f}})\big|_{h_{\text{\textgreek{t}},N}}^{2}\, dh_{N}+\sum_{\substack{components\\
of\,\mathcal{N}_{af,\mathcal{M}}
}
}\sum_{j=0}^{m-1}\sum_{j_{1}+j_{2}+j_{3}=j}\mathcal{E}_{bound,R;\text{\textgreek{d}}}^{(p,q)}[r^{-j_{2}}\partial_{v}^{j_{1}}\partial_{\text{\textgreek{sv}}}^{j_{2}}\partial_{u}^{j_{3}}\text{\textgreek{f}}](\text{\textgreek{t}})
\end{equation}
(where the derivatives in the second term are with respect to the
$(u,v,\text{\textgreek{sv}})$ coordinate charts over each connected
component of $\mathcal{N}_{af,\mathcal{M}}$),
\begin{align}
\mathcal{F}_{\text{\textgreek{h}},\text{\textgreek{e}}}^{(2,q,m,k)}[F](\text{\textgreek{t}})= & \text{\textgreek{t}}^{-2+\text{\textgreek{e}}}\sum_{j=0}^{m+3k-1}\sum_{j_{1}+j_{2}=j}\sum_{i_{1}+i_{2}\le q-1}\int_{\{0\le\bar{t}\le\text{\textgreek{t}}\}}r^{3+2i_{1}}\big|\nabla_{h_{\text{\textgreek{t}},N}}^{j_{1}+i_{1}}(N^{j_{2}+i_{2}}F)\big|_{h_{\text{\textgreek{t}},N}}^{2}\, dg+\\
 & +\text{\textgreek{t}}^{-1+\text{\textgreek{e}}}\sum_{j=0}^{m+2k-1}\sum_{j_{1}+j_{2}=j}\sum_{i_{1}+i_{2}\le q-1}\int_{\{C^{-1}\text{\textgreek{t}}\le\bar{t}\le\text{\textgreek{t}}\}}r^{2+2i_{1}}\big|\nabla_{h_{\text{\textgreek{t}},N}}^{j_{1}+i_{1}}(N^{j_{2}+i_{2}}F)\big|_{h_{\text{\textgreek{t}},N}}^{2}\, dg+\nonumber \\
 & +\sum_{j=0}^{m+2k-1}\sum_{j_{1}+j_{2}=j}\sum_{i_{1}+i_{2}\le q-1}\int_{\{C^{-1}\text{\textgreek{t}}\le\bar{t}\le\text{\textgreek{t}}\}}(r^{1+\text{\textgreek{e}}+2i_{1}}+r^{1+\text{\textgreek{h}}})\big|\nabla_{h_{\text{\textgreek{t}},N}}^{j_{1}+i_{1}}(N^{j_{2}+i_{2}}F)\big|_{h_{\text{\textgreek{t}},N}}^{2}\, dg\nonumber 
\end{align}
and 
\begin{align}
\mathcal{F}_{\text{\textgreek{h}}}^{(1,q,m,k)}[F](\text{\textgreek{t}})= & \text{\textgreek{t}}^{-1}\sum_{j=0}^{m+2k-1}\sum_{j_{1}+j_{2}=j}\sum_{i_{1}+i_{2}\le q-1}\int_{\{0\le\bar{t}\le\text{\textgreek{t}}\}}r^{3+2i_{1}}\big|\nabla_{h_{\text{\textgreek{t}},N}}^{j_{1}+i_{1}}(N^{j_{2}+i_{2}}F)\big|_{h_{\text{\textgreek{t}},N}}^{2}\, dg+\\
 & +\sum_{j=0}^{m+2k-1}\sum_{j_{1}+j_{2}=j}\sum_{i_{1}+i_{2}\le q-1}\int_{\{C^{-1}\text{\textgreek{t}}\le\bar{t}\le\text{\textgreek{t}}\}}r^{2+2i_{1}}\big|\nabla_{h_{\text{\textgreek{t}},N}}^{j_{1}+i_{1}}(N^{j_{2}+i_{2}}F)\big|_{h_{\text{\textgreek{t}},N}}^{2}\, dg.\nonumber 
\end{align}
\end{thm}
\begin{rem*}
Notice that each derivative of $\text{\textgreek{f}}$ tangential
to $\{\bar{t}=\text{\textgreek{t}}\}$ carries an extra $r$-weight.
Again, as before, in case there exists some small $\text{\textgreek{d}}_{0}>0$
such that the deformation tensor of the vector field $T$ in the region
$\{r\gg1\}$ satisfies the bound (\ref{eq:DeformationTensorTAwaySlow})
for any $m\in\mathbb{N}$ and the last term in the right hand side
of the integrated local energy decay estimate (\ref{eq:IntegratedLocalEnergyDecayImprovedDecay})
is replaced by
\begin{equation}
C_{m,\text{\textgreek{h}}}\sum_{j=0}^{m-1}\sum_{j_{1}+j_{2}=j}\int_{\{\text{\textgreek{t}}_{1}\le\bar{t}\le\text{\textgreek{t}}_{1}\}\cap\{r\ge R_{f}\}}\bar{t}^{-\text{\textgreek{d}}_{0}}r_{+}^{-1}\Big(\big|\nabla_{h_{\text{\textgreek{t}},N}}^{j_{1}+1}(N^{j_{2}}\text{\textgreek{f}})\big|^{2}+r_{+}^{-2}\big|\nabla_{h_{\text{\textgreek{t}},N}}^{j_{1}}(N^{j_{2}}\text{\textgreek{f}})\big|^{2}+r_{+}^{-2}\big|N^{j+1}\text{\textgreek{f}}\big|^{2}\Big)\, dg,
\end{equation}
then the $\text{\textgreek{e}}$ loss in the exponent of $\text{\textgreek{t}}$
in (\ref{eq:FinalDecayFirstenergy-1}) and (\ref{eq:PointwiseDecayFirst-1-1})
can be removed. 
\end{rem*}
The proof of this theorem follows exactly as that of Theorem \ref{thm:FirstPointwiseDecayNewMethod},
the only difference being that Corollary \ref{cor:HigherOrderNewMethodAddition}
is used in place of \ref{cor:NewMethodHigherDerivativesNotImproved},
and Lemma \ref{lem:Boundedness Hyperboloids Higher Order} in place
of \ref{lem:Boundedness Hyperboloids}. Thus, the details of the proof
will be omitted.

\subsection{Energy boundedness with loss of derivatives}

The integrated local energy decay assumption (\ref{eq:IntegratedLocalEnergyDecayImprovedDecayInfinity-1})
on $(\mathcal{M},g)$ allows us to establish the following energy
boundedness statement with loss of derivatives on $(\mathcal{M},g)$: 
\begin{lem}
\label{lem:BoundednessWithLoss}For any $R>0$, any integer $m\ge1$
and any $0<\text{\textgreek{h}}<a$ we can bound for any $0\le\text{\textgreek{t}}_{1}\le\text{\textgreek{t}}_{2}$
and any smooth function $\text{\textgreek{f}}$ satisfying $\square\text{\textgreek{f}}=F$
on $(\mathcal{M},g)$ with finite energy norm on $\{\bar{t}=\text{\textgreek{t}}_{1}\}$
(and satisfying boundary conditions on $\partial_{tim}\mathcal{M}$
belonging to the class $\mathcal{C}_{ILED}$), we can bound provided
(\ref{eq:NormFinitenessFriedlanderRadiation}) holds for all $0\le j\le m+d+k$:
\begin{equation}
\begin{split}\sum_{j=0}^{m-1}\sum_{j_{1}+j_{2}=j}\int_{\mathcal{H}^{+}(\text{\textgreek{t}}_{1},\text{\textgreek{t}}_{2})}\Big(\big|\nabla_{h_{\mathcal{H}}}^{j_{1}+1}(N^{j_{2}}\text{\textgreek{f}})\big|_{h_{\mathcal{H}}}^{2} & +|\text{\textgreek{f}}|^{2}\Big)\, dh_{\mathcal{H}}+\sum_{j=0}^{m}\sum_{j_{1}+j_{2}=j}\int_{\partial_{tim}\mathcal{M}\cap\{\text{\textgreek{t}}_{1}\le\bar{t}\le\text{\textgreek{t}}_{2}\}}\big|\nabla_{g}^{j}\text{\textgreek{f}}\big|_{h}^{2}\, dh_{\partial\mathcal{M}_{tim}}+\\
+\sum_{j=0}^{m-1}\sum_{j_{1}+j_{2}=j}\int_{\{\bar{t}=\text{\textgreek{t}}_{2}\}}\Big(\big|\nabla_{h_{N}}^{j_{1}+1}(N^{j_{2}}\text{\textgreek{f}})\big|_{h_{N}}^{2} & +r_{+}^{-2}\big|N^{j+1}\text{\textgreek{f}}|^{2}\Big)\, dh_{N}\le\\
\le\, & C_{m,R}\sum_{j=0}^{m+k-1}\sum_{j_{1}+j_{2}=j}\int_{\{\bar{t}=\text{\textgreek{t}}_{1}\}}\Big(\big|\nabla_{h_{N}}^{j_{1}+1}(N^{j_{2}}\text{\textgreek{f}})\big|_{h_{N}}^{2}+r_{+}^{-2}\big|N^{j+1}\text{\textgreek{f}}|^{2}\Big)\, dh_{N}+\\
 & +\sum_{j=0}^{m-1}\sum_{j_{1}+j_{2}=j}\int_{\{\text{\textgreek{t}}_{1}\le\bar{t}\le\text{\textgreek{t}}_{2}\}\cap\{r\ge R\}}r_{+}^{-1}\Big(\big|\nabla_{h_{N}}^{j_{1}+1}(N^{j_{2}}\text{\textgreek{f}})\big|_{h_{N}}^{2}+r_{+}^{-2}\big|N^{j+1}\text{\textgreek{f}}|^{2}+r_{+}^{-2}|\text{\textgreek{f}}|^{2}\Big)\, dh+\\
 & +C_{m,R}\sum_{j=0}^{m+k-1}\int_{\{\text{\textgreek{t}}_{1}\le\bar{t}\le\text{\textgreek{t}}_{2}\}}r_{+}^{1+\text{\textgreek{h}}}|\nabla_{g}^{j}F|_{h}^{2}\, dh.
\end{split}
\label{eq:BoundednessWithLoss}
\end{equation}
 where $\mathcal{H}^{+}(\text{\textgreek{t}}_{1},\text{\textgreek{t}}_{2})=\mathcal{H}^{+}\cap\{\text{\textgreek{t}}_{1}\le\bar{t}\le\text{\textgreek{t}}_{2}\}$
and $k$ is the integer measuring the loss of derivatives in (\ref{eq:IntegratedLocalEnergyDecayImprovedDecay}).\end{lem}
\begin{proof}
Without loss of generality, we can assume that $R$ is large in terms
of the geometry of $(\mathcal{M},g)$. Let us fix a second smooth
vector field $N_{1}$ on $(\mathcal{M},g)$, such that:

\begin{itemize}

\item $N_{1}\equiv N$ on $\{r\ge1\}$

\item $|N_{1}|_{h}\le2$ on $(\mathcal{M},g)$

\item $g(N_{1},N_{1})\le-c<0$ everywhere on $(\mathcal{M},g)$

\item For any $l\in\mathbb{N}$: $\big|\mathcal{L}_{N_{1}}^{l}g\big|_{h}\le C_{l}$.%
\footnote{This is possible in view of Assumptions \ref{enu:DeformationTensor}
and \ref{enu:UniformityY}.%
}

\item $N_{1}$ and $N$ are linearly independent on $\mathcal{H}^{+}$
and their span intersects the tangent space of the surfaces $\mathcal{H}^{+}\cap\{\bar{t}=const\}$
transversally

\item $\Big|\big(g(N_{1},N)\big)^{2}-g(N_{1},N_{1})g(N,N)\Big|\ge c>0$
on $\mathcal{H}^{+}$.

\end{itemize}

\noindent Notice that in the case $\mathcal{H}^{+}=\emptyset$ we
can simply choose $N_{1}=N$.

Fix a smooth cut-off function $\text{\textgreek{q}}:[0,+\infty)\rightarrow[0,1]$
such that $\text{\textgreek{q}}(x)=1$ for $x\le1$ and $\text{\textgreek{q}}(x)=0$
for $x\ge2$, and define: 
\begin{equation}
\text{\textgreek{q}}_{R}(r)\doteq\text{\textgreek{q}}(\frac{r}{R}).
\end{equation}
 Define also the energy current 
\begin{equation}
J_{\text{\textgreek{m}}}^{m}[\text{\textgreek{f}}]=\sum_{j=0}^{m-1}\sum_{j_{1}+j_{2}=j}\big(J_{\text{\textgreek{m}}}^{N}(N_{1}^{j_{1}}N^{j_{2}}\text{\textgreek{f}})\big).\label{eq:EnergyBoundedness current}
\end{equation}
By integrating $\nabla_{g}^{\text{\textgreek{m}}}\big(\text{\textgreek{q}}_{R}J_{\text{\textgreek{m}}}^{m}[\text{\textgreek{f}}]\big)$
over the region $\{\text{\textgreek{t}}_{1}\le\bar{t}\le\text{\textgreek{t}}_{2}\}$,
we obtain (due to the boundedness of the derivatives of the deformation
tensors of $N_{1},N$): 
\begin{align}
\int_{\mathcal{H}^{+}(\text{\textgreek{t}}_{1},\text{\textgreek{t}}_{2})}J_{\text{\textgreek{m}}}^{m}[\text{\textgreek{f}}]n_{\mathcal{H}}^{\text{\textgreek{m}}}+\int_{\{\bar{t}=\text{\textgreek{t}}_{2}\}\cap\{r\le R\}}J_{\text{\textgreek{m}}}^{m}[\text{\textgreek{f}}]\bar{n}^{\text{\textgreek{m}}}\le & \int_{\{\bar{t}=\text{\textgreek{t}}_{1}\}\cap\{r\le2R\}}J_{\text{\textgreek{m}}}^{m}[\text{\textgreek{f}}]\bar{n}^{\text{\textgreek{m}}}+C_{m}\int_{J^{+}(\{\bar{t}=\text{\textgreek{t}}_{1}\})\cap\{r\le2R\}}\sum_{j=0}^{m}|\nabla_{g}^{j}\text{\textgreek{f}}|_{h}^{2}+\label{eq:IntegratingTheBoundednessCurrent}\\
 & +\Big|\int_{\partial_{tim}\mathcal{M}\cap\{\text{\textgreek{t}}_{1}\le\bar{t}\le\text{\textgreek{t}}_{2}\}}J_{\text{\textgreek{m}}}^{m}[\text{\textgreek{f}}]n_{\partial_{tim}\mathcal{M}}^{\text{\textgreek{m}}}\Big|.\nonumber 
\end{align}

Adding to (\ref{eq:IntegratingTheBoundednessCurrent}) the integrated
local energy decay statement (\ref{eq:IntegratedLocalEnergyDecayImprovedDecayInfinity-1}),
and using Lemma \ref{lem:Boundedness Hyperboloids}, we obtain:
\begin{equation}
\begin{split}\int_{\mathcal{H}^{+}(\text{\textgreek{t}}_{1},\text{\textgreek{t}}_{2})}J_{\text{\textgreek{m}}}^{m}[\text{\textgreek{f}}] & n_{\mathcal{H}}^{\text{\textgreek{m}}}+\int_{\{\bar{t}=\text{\textgreek{t}}_{2}\}}J_{\text{\textgreek{m}}}^{m}[\text{\textgreek{f}}]\bar{n}^{\text{\textgreek{m}}}+\int_{\{\text{\textgreek{t}}_{1}\le\bar{t}\le\text{\textgreek{t}}_{2}\}\cap\{r\le2R\}}\sum_{j=0}^{m}|\nabla_{g}^{j}\text{\textgreek{f}}|_{h}^{2}\, dh+\int_{\partial_{tim}\mathcal{M}\cap\{\text{\textgreek{t}}_{1}\le\bar{t}\le\text{\textgreek{t}}_{2}\}}\sum_{j=0}^{m}|\nabla_{g}^{j}\text{\textgreek{f}}|_{h}^{2}\, dh_{\partial_{tim}\mathcal{M}}\le\\
\le & C_{m,R}\int_{\{\bar{t}=\text{\textgreek{t}}_{1}\}}J_{\text{\textgreek{m}}}^{m+k}[\text{\textgreek{f}}]\bar{n}^{\text{\textgreek{m}}}+\sum_{j=0}^{m-1}\sum_{j_{1}+j_{2}=j}\int_{\{\text{\textgreek{t}}_{1}\le\bar{t}\le\text{\textgreek{t}}_{2}\}\cap\{r\ge2R\}}r_{+}^{-1}\Big(\big|\nabla_{h_{N}}^{j_{1}+1}(N^{j_{2}}\text{\textgreek{f}})\big|_{h_{N}}^{2}+r_{+}^{-2}\big|N^{j+1}\text{\textgreek{f}}|^{2}+r_{+}^{-2}|\text{\textgreek{f}}|^{2}\Big)\, dh+\\
 & +C_{m,R}\sum_{j=0}^{m+k-1}\int_{\{\text{\textgreek{t}}_{1}\le\bar{t}\le\text{\textgreek{t}}_{2}\}}r_{+}^{1+\text{\textgreek{h}}}|\nabla_{g}^{j}F|_{h}^{2}\, dh.
\end{split}
\label{eq:BeforeEllipticEstimates}
\end{equation}

In view of the assumptions on $N_{1},N$ (it is here that we make
use of the fact that their span is transversal to $\mathcal{H}^{+}\cap\{\bar{t}=const\}$),
as well Assumption \ref{enu:EllipticEstimatesHorizon} on the uniformity
of the elliptic estimates on sections of $\mathcal{H}^{+}$ (note
that Theorem \ref{thm:FriedlanderRadiation} applies to yield $\limsup_{r\rightarrow+\infty}\big|r^{\frac{d-1}{2}+j_{1}}\nabla_{h_{\text{\textgreek{t}},N}}^{j_{1}}(N^{j_{2}}\text{\textgreek{f}})\big|_{h_{\text{\textgreek{t}},N}}<+\infty$
for all $j_{1}+j_{2}\le m+k$), we can bound:
\begin{equation}
\int_{\mathcal{H}^{+}(\text{\textgreek{t}}_{1},\text{\textgreek{t}}_{2})}J_{\text{\textgreek{m}}}^{m}[\text{\textgreek{f}}]n_{\mathcal{H}}^{\text{\textgreek{m}}}+\sum_{j=0}^{m-2}\int_{\mathcal{H}^{+}(\text{\textgreek{t}}_{1},\text{\textgreek{t}}_{2})}|\nabla_{g}^{j}F|_{h}\, dh_{\mathcal{H}}\ge c_{m}\cdot\sum_{j=0}^{m-1}\sum_{j_{1}+j_{2}=j}\int_{\mathcal{H}^{+}(\text{\textgreek{t}}_{1},\text{\textgreek{t}}_{2})}\big|\nabla_{h_{\mathcal{H}}}^{j_{1}+1}(N^{j_{2}}\text{\textgreek{f}})\big|_{h_{\mathcal{H}}}^{2}\, dh_{\mathcal{H}}.\label{eq:TermsHorizon}
\end{equation}
 Thus, from (\ref{eq:BeforeEllipticEstimates}) and (\ref{eq:TermsHorizon}),
as well as a trace theorem on the horizon, we obtain:
\begin{equation}
\begin{split}\sum_{j=0}^{m-1}\sum_{j_{1}+j_{2}=j}\int_{\mathcal{H}^{+}(\text{\textgreek{t}}_{1},\text{\textgreek{t}}_{2})}\Big( & \big|\nabla_{h_{\mathcal{H}}}^{j_{1}+1}(N^{j_{2}}\text{\textgreek{f}})\big|_{h_{\mathcal{H}}}^{2}+|\text{\textgreek{f}}|^{2}\Big)\, dh_{\mathcal{H}}+\int_{\{\bar{t}=\text{\textgreek{t}}_{2}\}}J_{\text{\textgreek{m}}}^{m}[\text{\textgreek{f}}]\bar{n}^{\text{\textgreek{m}}}+\\
+\int_{\{\text{\textgreek{t}}_{1}\le\bar{t}\le\text{\textgreek{t}}_{2}\}\cap\{r\le2R\}} & \sum_{j=0}^{m}|\nabla_{g}^{j}\text{\textgreek{f}}|_{h}^{2}\, dh+\int_{\partial_{tim}\mathcal{M}\cap\{\text{\textgreek{t}}_{1}\le\bar{t}\le\text{\textgreek{t}}_{2}\}}\sum_{j=0}^{m}|\nabla_{g}^{j}\text{\textgreek{f}}|_{h}^{2}\, dh_{\partial_{tim}\mathcal{M}}\le\\
\le & C_{m,R}\int_{\{\bar{t}=\text{\textgreek{t}}_{1}\}}J_{\text{\textgreek{m}}}^{m+k}[\text{\textgreek{f}}]\bar{n}^{\text{\textgreek{m}}}+\sum_{j=0}^{m-1}\sum_{j_{1}+j_{2}=j}\int_{\{\text{\textgreek{t}}_{1}\le\bar{t}\le\text{\textgreek{t}}_{2}\}\cap\{r\ge R\}}r_{+}^{-1}\Big(\big|\nabla_{h_{N}}^{j_{1}+1}(N^{j_{2}}\text{\textgreek{f}})\big|_{h_{N}}^{2}+r_{+}^{-2}\big|N^{j+1}\text{\textgreek{f}}|^{2}+r_{+}^{-2}|\text{\textgreek{f}}|^{2}\Big)\, dh+\\
 & +C_{m,R}\sum_{j=0}^{m-2}\int_{\mathcal{H}^{+}(\text{\textgreek{t}}_{1},\text{\textgreek{t}}_{2})}|\nabla_{g}^{j}F|_{h}\, dh_{\mathcal{H}}.
\end{split}
\label{eq:BeforeEllipticEstimates-1}
\end{equation}

Moreover, using Assumptions \ref{enu:EllipticEstimates}, \ref{enu:EllipticEstimatesHorizon}
and \ref{enu:PoincareInequalities} on the uniformity of elliptic
estimates and trace inequalities, we can bound:
\begin{equation}
\begin{split}\sum_{j=0}^{m-1}\sum_{j_{1}+j_{2}=j}\int_{\{\bar{t}=\text{\textgreek{t}}_{2}\}}\Big(\big|\nabla_{h_{N}}^{j_{1}+1}(N^{j_{2}}\text{\textgreek{f}})\big|_{h_{N}}^{2} & +r_{+}^{-2}\big|N^{j+1}\text{\textgreek{f}}|^{2}\Big)\, dh_{N}\le C_{m}\Big(\int_{\{\bar{t}=\text{\textgreek{t}}_{2}\}}J_{\text{\textgreek{m}}}^{m}[\text{\textgreek{f}}]\bar{n}^{\text{\textgreek{m}}}+\\
 & +\sum_{j=0}^{m-1}\sum_{j_{1}+j_{2}=j}\int_{\mathcal{H}^{+}(\text{\textgreek{t}}_{1},\text{\textgreek{t}}_{2})}\Big(\big|\nabla_{h_{\mathcal{H}}}^{j_{1}+1}(N^{j_{2}}\text{\textgreek{f}})\big|_{h_{\mathcal{H}}}^{2}+|\text{\textgreek{f}}|^{2}\Big)\, dh_{\mathcal{H}}+\\
 & +\int_{\partial_{tim}\mathcal{M}\cap\{\text{\textgreek{t}}_{1}\le\bar{t}\le\text{\textgreek{t}}_{2}\}}\sum_{j=0}^{m}|\nabla_{g}^{j}\text{\textgreek{f}}|_{h}^{2}\, dh_{\partial_{tim}\mathcal{M}}+\sum_{j=0}^{m-2}\int_{\{\bar{t}=\text{\textgreek{t}}_{2}\}}|\nabla_{g}^{j}F|_{h}^{2}\, dh_{N}\Big).
\end{split}
\label{eq:BeforeSharpTrace}
\end{equation}
 Thus, using (\ref{eq:BeforeEllipticEstimates-1}) and (\ref{eq:BeforeSharpTrace})
and a trace theorem for the terms 
\[
\sum_{j=0}^{m-2}\Big(\int_{\{\bar{t}=\text{\textgreek{t}}_{2}\}}|\nabla_{g}^{j}F|_{h}^{2}\, dh_{N}+\int_{\mathcal{H}^{+}(\text{\textgreek{t}}_{1},\text{\textgreek{t}}_{2})}|\nabla_{g}^{j}F|_{h}\, dh_{\mathcal{H}}\Big),
\]
 we deduce the required energy boundedness estimate:
\begin{equation}
\begin{split}\sum_{j=0}^{m-1}\sum_{j_{1}+j_{2}=j}\int_{\mathcal{H}^{+}(\text{\textgreek{t}}_{1},\text{\textgreek{t}}_{2})}\Big(\big|\nabla_{h_{\mathcal{H}}}^{j_{1}+1}(N^{j_{2}}\text{\textgreek{f}})\big|_{h_{\mathcal{H}}}^{2}+ & |\text{\textgreek{f}}|^{2}\Big)\, dh_{\mathcal{H}}+\sum_{j=0}^{m-1}\sum_{j_{1}+j_{2}=j}\int_{\{\bar{t}=\text{\textgreek{t}}_{2}\}}\Big(\big|\nabla_{h_{N}}^{j_{1}+1}(N^{j_{2}}\text{\textgreek{f}})\big|_{h_{N}}^{2}+r_{+}^{-2}\big|N^{j+1}\text{\textgreek{f}}|^{2}\Big)\, dh_{N}+\\
+\int_{\partial_{tim}\mathcal{M}\cap\{\text{\textgreek{t}}_{1}\le\bar{t}\le\text{\textgreek{t}}_{2}\}}\sum_{j=0}^{m}|\nabla_{g}^{j}\text{\textgreek{f}}|_{h}^{2}\, dh_{\partial_{tim}\mathcal{M}}\le & C_{m,R}\sum_{j=0}^{m+k-1}\sum_{j_{1}+j_{2}=j}\int_{\{\bar{t}=\text{\textgreek{t}}_{1}\}}\Big(\big|\nabla_{h_{N}}^{j_{1}+1}(N^{j_{2}}\text{\textgreek{f}})\big|_{h_{N}}^{2}+r_{+}^{-2}\big|N^{j+1}\text{\textgreek{f}}|^{2}\Big)\, dh_{N}+\\
 & +\sum_{j=0}^{m-1}\sum_{j_{1}+j_{2}=j}\int_{\{\text{\textgreek{t}}_{1}\le\bar{t}\le\text{\textgreek{t}}_{2}\}\cap\{r\ge R\}}r^{-1}\Big(\big|\nabla_{h_{N}}^{j_{1}+1}(N^{j_{2}}\text{\textgreek{f}})\big|_{h_{N}}^{2}+r_{+}^{-2}\big|N^{j+1}\text{\textgreek{f}}|^{2}+r_{+}^{-2}|\text{\textgreek{f}}|^{2}\Big)\, dh+\\
 & +C_{m,R}\sum_{j=0}^{m+k-1}\int_{\{\text{\textgreek{t}}_{1}\le\bar{t}\le\text{\textgreek{t}}_{2}\}}r_{+}^{1+\text{\textgreek{h}}}|\nabla_{g}^{j}F|_{h}^{2}\, dh.
\end{split}
\label{eq:BoundednessWithLoss-1}
\end{equation}

\end{proof}

\subsection{\label{sub:ProofOfFirstPointwiseDecay}Proof of Theorem \ref{thm:FirstPointwiseDecayNewMethod}}

We will assume without loss of generality that $\text{\textgreek{f}}$
is real valued. We will also set $\text{\textgreek{F}}=\text{\textgreek{W}}\text{\textgreek{f}}.$

Fix an $R>0$ large enough in terms of the geometry of $(\mathcal{M},g)$.
Fix also a smooth cut-off $\text{\textgreek{q}}_{R}:\mathcal{M}\rightarrow[0,1]$
which is only a function of $r$, such that $\text{\textgreek{q}}_{R}\equiv0$
on $\{r\le R\}$ and $\text{\textgreek{q}}_{R}\equiv1$ on $\{r\ge R+1\}$.
Fix also a small number $0<\text{\textgreek{d}}\ll1$. 

We will use the following notations for the $r^{p}$-weighted energy
norms for any $\text{\textgreek{t}}\ge0$: 
\begin{equation}
\mathcal{E}_{bound}^{(p,m)}[\text{\textgreek{f}}](\text{\textgreek{t}})=\sum_{j=0}^{m}\int_{\{\bar{t}=\text{\textgreek{t}}\}\cap\{r\le R+1\}}\big|\nabla_{h_{\text{\textgreek{t}},N}}^{j_{1}}(N^{j_{2}}\text{\textgreek{f}})\big|_{h_{\text{\textgreek{t}},N}}^{2}\, dh_{N}+\sum_{\substack{components\\
of\,\mathcal{N}_{af,\mathcal{M}}
}
}\sum_{j=0}^{m-1}\sum_{j_{1}+j_{2}+j_{3}=j}\mathcal{E}_{bound,R;\text{\textgreek{d}}}^{(p,0)}[r^{-j_{2}}\partial_{v}^{j_{1}}\partial_{\text{\textgreek{sv}}}^{j_{2}}\partial_{u}^{j_{3}}\text{\textgreek{f}}](\text{\textgreek{t}}),
\end{equation}
\begin{align}
\mathcal{E}_{bulk,\text{\textgreek{h}}}^{(p,m)}[\text{\textgreek{f}}](\text{\textgreek{t}})= & \sum_{j=0}^{m}\int_{\{\bar{t}=\text{\textgreek{t}}\}\cap\{r\le R+1\}}\big|\nabla_{h_{\text{\textgreek{t}},N}}^{j_{1}}N^{j_{2}}\text{\textgreek{f}}\big|_{h_{\text{\textgreek{t}},N}}^{2}\, dh_{N}+\sum_{j=1}^{m}\sum_{j_{1}+j_{2}=j-1}\int_{\{\bar{t}=\text{\textgreek{t}}\}\cap\mathcal{H}}\big|\nabla_{h_{\mathcal{H}}}^{j_{1}+1}(N^{j_{2}}\text{\textgreek{f}})\big|_{h_{\mathcal{H}}}^{2}\, dh_{\mathcal{H}_{\text{\textgreek{t}}}}+\\
 & +\sum_{j=0}^{m}\int_{\{\bar{t}=\text{\textgreek{t}}\}\cap\partial_{tim}\mathcal{M}}\big|\nabla_{g}^{j}\text{\textgreek{f}}\big|_{h}^{2}\, dh_{\partial_{tim}\mathcal{M}_{,\text{\textgreek{t}}}}+\sum_{\substack{components\\
of\,\mathcal{N}_{af,\mathcal{M}}
}
}\sum_{j=0}^{m-1}\sum_{j_{1}+j_{2}+j_{3}=j}\mathcal{E}_{bulk,R,\text{\textgreek{h}};\text{\textgreek{d}}}^{(p,0)}[r^{-j_{2}}\partial_{v}^{j_{1}}\partial_{\text{\textgreek{sv}}}^{j_{2}}\partial_{u}^{j_{3}}\text{\textgreek{f}}](\text{\textgreek{t}})\nonumber 
\end{align}
and
\begin{equation}
\mathcal{E}_{en}^{(m)}[\text{\textgreek{f}}](\text{\textgreek{t}})=\sum_{j=1}^{m}\sum_{j_{1}+j_{2}=j-1}\int_{\{\bar{t}=\text{\textgreek{t}}\}}\Big(\big|\nabla_{h_{\text{\textgreek{t}},N}}^{j_{1}+1}(N^{j_{3}}\text{\textgreek{f}})\big|_{h_{\text{\textgreek{t}},N}}^{2}+r_{+}^{-2}\big(\big|N^{j}\text{\textgreek{f}}\big|^{2}+|\text{\textgreek{f}}|^{2}\big)\Big)\, dh_{N}.
\end{equation}
 We will also use the norm 
\begin{equation}
\mathcal{F}_{\text{\textgreek{h}}}^{(p,m)}[\text{\textgreek{f}}](\text{\textgreek{t}})=\sum_{j=0}^{m-1}\int_{\{\bar{t}=\text{\textgreek{t}}\}}(r^{p+1}+r^{1+\text{\textgreek{h}}})\big|\nabla_{g}^{j}(\square\text{\textgreek{f}})\big|_{h}^{2}\, dh_{N}.
\end{equation}

Using Corollary \ref{cor:NewMethodHigherDerivativesNotImproved} and
Lemma \ref{lem:BoundednessWithLoss}, we can bound for any integer
$m\ge0$, any $0\le\text{\textgreek{t}}_{1}\le\text{\textgreek{t}}_{2}$
and any $0<p\le2$:

\begin{equation}
\mathcal{E}_{bound}^{(p,m)}[\text{\textgreek{f}}](\text{\textgreek{t}}_{2})+\int_{\text{\textgreek{t}}_{1}}^{\text{\textgreek{t}}_{2}}\mathcal{E}_{bulk,\text{\textgreek{h}}}^{(p-1,m)}[\text{\textgreek{f}}](\text{\textgreek{t}})\, d\text{\textgreek{t}}\lesssim_{p,m,\text{\textgreek{h}}}\mathcal{E}_{bound}^{(p,m+k)}[\text{\textgreek{f}}](\text{\textgreek{t}}_{1})+\int_{\text{\textgreek{t}}_{1}}^{\text{\textgreek{t}}_{2}}\mathcal{F}_{\text{\textgreek{h}}}^{(p,m+k)}[\text{\textgreek{f}}](\text{\textgreek{t}})\, d\text{\textgreek{t}}.\label{eq:StartPointNewMethod}
\end{equation}

Starting from (\ref{eq:StartPointNewMethod}) for $\text{\textgreek{t}}_{1}=0$
and letting $\text{\textgreek{t}}_{2}\rightarrow+\infty$, we obtain
for any $m\in\mathbb{N}$ and any $T\ge0$: 
\begin{equation}
\int_{0}^{T}\mathcal{E}_{bulk,\text{\textgreek{h}}}^{(1,m)}[\text{\textgreek{f}}](\text{\textgreek{t}})\, d\text{\textgreek{t}}\lesssim_{m,\text{\textgreek{h}}}\mathcal{E}_{bound}^{(2,m+k)}[\text{\textgreek{f}}](0)+\int_{0}^{T}\mathcal{F}_{\text{\textgreek{h}}}^{(2,m+k)}[\text{\textgreek{f}}](\text{\textgreek{t}})\, d\text{\textgreek{t}}.\label{eq:FirstStepBeforePigeonhole}
\end{equation}
 An application of the pigeonhole principle on (\ref{eq:FirstStepBeforePigeonhole})
readily yields that there exists a sequence of positive numbers $\{\text{\textgreek{t}}_{n}\}_{n\in\mathbb{N}}$
with $\text{\textgreek{t}}_{0}\ge1$ and $2\text{\textgreek{t}}_{n}\le\text{\textgreek{t}}_{n+1}\le4\text{\textgreek{t}}_{n}$
such that 
\begin{equation}
\mathcal{E}_{bulk,\text{\textgreek{h}}}^{(1,m)}[\text{\textgreek{f}}](\text{\textgreek{t}}_{n})\lesssim_{m,\text{\textgreek{h}}}\text{\textgreek{t}}_{n}^{-1}\Big(\mathcal{E}_{bound}^{(2,m+k)}[\text{\textgreek{f}}](0)+\int_{0}^{\text{\textgreek{t}}_{n}}\mathcal{F}_{\text{\textgreek{h}}}^{(2,m+k)}[\text{\textgreek{f}}](\text{\textgreek{t}})\, d\text{\textgreek{t}}\Big).\label{eq:FirstStepDecay}
\end{equation}

Applying (\ref{eq:StartPointNewMethod}) in the intevals $\{\text{\textgreek{t}}_{n}\le\bar{t}\le\text{\textgreek{t}}_{n+1}\}$
for $p=1$ and using (\ref{eq:FirstStepDecay}) for $m+k$ in place
of $m$ to bound the first term of the right hand side,%
\footnote{Notice the trivial inequality $\mathcal{E}_{bound}^{(p,m)}[\text{\textgreek{f}}](\text{\textgreek{t}})\lesssim_{p,m,\text{\textgreek{h}}}\mathcal{E}_{bulk,\text{\textgreek{h}}}^{(p,m)}[\text{\textgreek{f}}](\text{\textgreek{t}})$.%
} we obtain for any $\bar{\text{\textgreek{t}}}\in[\text{\textgreek{t}}_{n},\text{\textgreek{t}}_{n+1}]$:
\begin{equation}
\int_{\text{\textgreek{t}}_{n}}^{\bar{\text{\textgreek{t}}}}\mathcal{E}_{bulk,\text{\textgreek{h}}}^{(0,m)}[\text{\textgreek{f}}](\text{\textgreek{t}})\, d\text{\textgreek{t}}\lesssim_{m,\text{\textgreek{h}}}\text{\textgreek{t}}_{n}^{-1}\Big(\mathcal{E}_{bound}^{(2,m+2k)}[\text{\textgreek{f}}](0)+\int_{0}^{\text{\textgreek{t}}_{n}}\mathcal{F}_{\text{\textgreek{h}}}^{(2,m+2k)}[\text{\textgreek{f}}](\text{\textgreek{t}})\, d\text{\textgreek{t}}\Big)+\int_{\text{\textgreek{t}}_{n}}^{\bar{\text{\textgreek{t}}}}\mathcal{F}_{\text{\textgreek{h}}}^{(1,m+k)}[\text{\textgreek{f}}](\text{\textgreek{t}})\, d\text{\textgreek{t}}.\label{eq:FirstStepBeforePigeonhole-1}
\end{equation}
Using the mean value theorem on the intervals $\{\text{\textgreek{t}}_{n}\le\bar{t}\le\text{\textgreek{t}}_{n+1}\}$,
we thus obtain on a possibly different sequence $\bar{\text{\textgreek{t}}}_{n}\in[\text{\textgreek{t}}_{n}+\frac{\text{\textgreek{t}}_{n+1}-\text{\textgreek{t}}_{n}}{4},\text{\textgreek{t}}_{n+1}-\frac{\text{\textgreek{t}}_{n+1}-\text{\textgreek{t}}_{n}}{4}]$:
\begin{equation}
\mathcal{E}_{bulk,\text{\textgreek{h}}}^{(0,m)}[\text{\textgreek{f}}](\bar{\text{\textgreek{t}}}_{n})\lesssim_{m,\text{\textgreek{h}}}\bar{\text{\textgreek{t}}}_{n}^{-2}\Big(\mathcal{E}_{bound}^{(2,m+2k)}[\text{\textgreek{f}}](0)+\int_{0}^{\bar{\text{\textgreek{t}}}_{n}}\mathcal{F}_{\text{\textgreek{h}}}^{(2,m+2k)}[\text{\textgreek{f}}](\text{\textgreek{t}})\, d\text{\textgreek{t}}\Big)+\bar{\text{\textgreek{t}}}_{n}^{-1}\int_{\frac{1}{4}\bar{\text{\textgreek{t}}}_{n}}^{\bar{\text{\textgreek{t}}}_{n}}\mathcal{F}_{\text{\textgreek{h}}}^{(1,m+k)}[\text{\textgreek{f}}](\text{\textgreek{t}})\, d\text{\textgreek{t}}.\label{eq:FinalStepDyadic}
\end{equation}

Notice that by interpolating between (\ref{eq:FinalStepDyadic}) for
$p=2$ on the intervals $\{\bar{\text{\textgreek{t}}}_{n}\le\bar{t}\le\bar{\text{\textgreek{t}}}_{n+1}\}$,
we can also bound: 
\begin{equation}
\mathcal{E}_{bound,\text{\textgreek{h}}}^{(\text{\textgreek{e}},m)}[\text{\textgreek{f}}](\bar{\text{\textgreek{t}}}_{n})\lesssim_{m,\text{\textgreek{e}}}\bar{\text{\textgreek{t}}}_{n}^{-2+\text{\textgreek{e}}}\Big(\mathcal{E}_{bound}^{(2,m+2k)}[\text{\textgreek{f}}](0)+\int_{0}^{\bar{\text{\textgreek{t}}}_{n}}\mathcal{F}_{\text{\textgreek{h}}}^{(2,m+2k)}[\text{\textgreek{f}}](\text{\textgreek{t}})\, d\text{\textgreek{t}}\Big)+\bar{\text{\textgreek{t}}}_{n}^{-1+\text{\textgreek{e}}}\int_{\frac{1}{4}\bar{\text{\textgreek{t}}}_{n}}^{\bar{\text{\textgreek{t}}}_{n}}\mathcal{F}_{\text{\textgreek{h}}}^{(1,m+k)}[\text{\textgreek{f}}](\text{\textgreek{t}})\, d\text{\textgreek{t}}.\label{eq:FinalStepDyadic-1}
\end{equation}
Applying (\ref{eq:StartPointNewMethod}) for $p=\text{\textgreek{e}}$
on the intervals $\{\bar{\text{\textgreek{t}}}_{n}\le\bar{t}\le\bar{\text{\textgreek{t}}}_{n+1}\}$
and using (\ref{eq:FinalStepDyadic-1}), we obtain for any $\bar{\text{\textgreek{t}}}\in[\bar{\text{\textgreek{t}}}_{n},\bar{\text{\textgreek{t}}}_{n+1}]$:
\begin{align}
\int_{\bar{\text{\textgreek{t}}}_{n}}^{\bar{\text{\textgreek{t}}}}\mathcal{E}_{bulk,\text{\textgreek{h}}}^{(-1+\text{\textgreek{e}},m)}[\text{\textgreek{f}}](\text{\textgreek{t}})\, d\text{\textgreek{t}}\lesssim_{m,\text{\textgreek{e}},\text{\textgreek{h}}} & \bar{\text{\textgreek{t}}}_{n}^{-2+\text{\textgreek{e}}}\Big(\mathcal{E}_{bound}^{(2,m+2k)}[\text{\textgreek{f}}](0)+\int_{0}^{\bar{\text{\textgreek{t}}}_{n}}\mathcal{F}_{\text{\textgreek{h}}}^{(2,m+2k)}[\text{\textgreek{f}}](\text{\textgreek{t}})\, d\text{\textgreek{t}}\Big)+\bar{\text{\textgreek{t}}}_{n}^{-1+\text{\textgreek{e}}}\int_{\frac{1}{4}\bar{\text{\textgreek{t}}}_{n}}^{\bar{\text{\textgreek{t}}}_{n}}\mathcal{F}_{\text{\textgreek{h}}}^{(1,m+k)}[\text{\textgreek{f}}](\text{\textgreek{t}})\, d\text{\textgreek{t}}+\label{eq:BoundForILEDErrors}\\
 & +\int_{\bar{\text{\textgreek{t}}}_{n}}^{\bar{\text{\textgreek{t}}}}\mathcal{F}_{\text{\textgreek{h}}}^{(\text{\textgreek{e}},m+k)}[\text{\textgreek{f}}](\text{\textgreek{t}})\, d\text{\textgreek{t}}.\nonumber 
\end{align}

Using the fact that 
\begin{equation}
\mathcal{E}_{en}^{(m)}[\text{\textgreek{f}}](\text{\textgreek{t}})\lesssim_{m,\text{\textgreek{h}}}\mathcal{E}_{bulk,\text{\textgreek{h}}}^{(0,m)}[\text{\textgreek{f}}](\text{\textgreek{t}}),
\end{equation}
 by applying Lemma \ref{lem:BoundednessWithLoss} on the intervals
$\{\bar{\text{\textgreek{t}}}_{n}\le\bar{t}\le\bar{\text{\textgreek{t}}}_{n+1}\}$
and using (\ref{eq:FinalStepDyadic}) and (\ref{eq:BoundForILEDErrors})
(for $m+k$ in place of $k$), we obtain (\ref{eq:FinalDecayFirstenergy})
for some $C>1$ and any $\text{\textgreek{t}}\ge1$ in view of the
fact that $\frac{5}{4}\bar{\text{\textgreek{t}}}_{n}\le\bar{\text{\textgreek{t}}}_{n+1}\le8\bar{\text{\textgreek{t}}}_{n}$:
\begin{align}
\mathcal{E}_{en}^{(m)}[\text{\textgreek{f}}](\text{\textgreek{t}})\lesssim_{m,\text{\textgreek{e}},\text{\textgreek{h}}} & \text{\textgreek{t}}^{-2+\text{\textgreek{e}}}\Big(\mathcal{E}_{bound}^{(2,m+3k)}[\text{\textgreek{f}}](0)+\int_{0}^{\text{\textgreek{t}}}\mathcal{F}_{\text{\textgreek{h}}}^{(2,m+3k)}[\text{\textgreek{f}}](\text{\textgreek{t}})\, d\text{\textgreek{t}}\Big)+\text{\textgreek{t}}^{-1+\text{\textgreek{e}}}\int_{C^{-1}\text{\textgreek{t}}}^{\text{\textgreek{t}}}\mathcal{F}_{\text{\textgreek{h}}}^{(1,m+2k)}[\text{\textgreek{f}}](\text{\textgreek{t}})\, d\text{\textgreek{t}}+\label{eq:FinalEnergyDecayNewMethod}\\
 & +\int_{C^{-1}\text{\textgreek{t}}}^{\text{\textgreek{t}}}\mathcal{F}_{\text{\textgreek{h}}}^{(\text{\textgreek{e}},m+2k)}[\text{\textgreek{f}}](\text{\textgreek{t}})\, d\text{\textgreek{t}}.\nonumber 
\end{align}

In view of Assumption \ref{enu:SobolevInequality}, we obtain from
(\ref{eq:FinalEnergyDecayNewMethod}) for any integer $m\ge0$ and
any $\text{\textgreek{t}}\ge0$ using the Sobolev embedding theorem:
\begin{align}
\sup_{\{\bar{t}=\text{\textgreek{t}}\}}\big|\nabla^{m}\text{\textgreek{f}}\big|_{h}^{2}\lesssim_{m,\text{\textgreek{e}},\text{\textgreek{h}}} & \text{\textgreek{t}}^{-2+\text{\textgreek{e}}}\Big(\mathcal{E}_{bound}^{(2,m+\lceil\frac{d+1}{2}\rceil+3k)}[\text{\textgreek{f}}](0)+\int_{0}^{\text{\textgreek{t}}}\mathcal{F}_{\text{\textgreek{h}}}^{(2,m+\lceil\frac{d+1}{2}\rceil+3k)}[\text{\textgreek{f}}](\text{\textgreek{t}})\, d\text{\textgreek{t}}\Big)+\text{\textgreek{t}}^{-1+\text{\textgreek{e}}}\int_{C^{-1}\text{\textgreek{t}}}^{\text{\textgreek{t}}}\mathcal{F}_{\text{\textgreek{h}}}^{(1,m+\lceil\frac{d+1}{2}\rceil+2k)}[\text{\textgreek{f}}](\text{\textgreek{t}})\, d\text{\textgreek{t}}+\label{eq:PointwiseDecayWeakInR}\\
 & +\int_{C^{-1}\text{\textgreek{t}}}^{\text{\textgreek{t}}}\mathcal{F}_{\text{\textgreek{h}}}^{(\text{\textgreek{e}},m+\lceil\frac{d+1}{2}\rceil+2k)}[\text{\textgreek{f}}](\text{\textgreek{t}})\, d\text{\textgreek{t}}.\nonumber 
\end{align}

Using the fundamental theorem of calculus and applying the product
rule for derivatives, we can bound for any function $\text{\textgreek{y}}$
on $\mathcal{M}$ and any $1\ll R_{1}<R_{2}$: 
\begin{equation}
\int_{\{\bar{t}=\text{\textgreek{t}}\}\cap\{r=R_{2}\}}|r^{\frac{d-2}{2}}\text{\textgreek{y}}|^{2}\, d\text{\textgreek{sv}}\lesssim\int_{\{\bar{t}=\text{\textgreek{t}}\}\cap\{r=R_{1}\}}|r^{\frac{d-2}{2}}\text{\textgreek{y}}|^{2}\, d\text{\textgreek{sv}}+\int_{\{\bar{t}=\text{\textgreek{t}}\}\cap\{R_{1}\le r\le R_{2}\}}|r^{\frac{d-3}{2}}\text{\textgreek{y}}|\big(|\partial_{v}(r^{\frac{d-1}{2}}\text{\textgreek{y}})|+r^{-2}|\partial_{u}(r^{\frac{d-1}{2}}\text{\textgreek{y}})|+|r^{\frac{d-3}{2}}\text{\textgreek{y}}|\big)\, dvd\text{\textgreek{sv}}.\label{eq:FromFundamentalTheoremOfCalculus}
\end{equation}
Hence, using a Sobolev inequality on $\mathbb{S}^{d-1}$, a trace
inequality for the first term of the right hand side of (\ref{eq:FromFundamentalTheoremOfCalculus})
and a Cauchy--Schwarz inequality for the second term, we infer from
(\ref{eq:FromFundamentalTheoremOfCalculus}) for any integer $m\ge0$:
\begin{equation}
\sup_{\{\bar{t}=\text{\textgreek{t}}\}}\big|r_{+}^{\frac{d-2}{2}}\nabla_{g}^{m}\text{\textgreek{f}}\big|_{h}^{2}\lesssim_{m}\mathcal{E}_{en}^{(m+\lceil\frac{d+2}{2}\rceil)}[\text{\textgreek{f}}](\text{\textgreek{t}}).\label{eq:ImprovedSobolev}
\end{equation}
 Thus, (\ref{eq:FinalEnergyDecayNewMethod}) yields: 
\begin{align}
\sup_{\{\bar{t}=\text{\textgreek{t}}\}}\big|r_{+}^{\frac{d-2}{2}}\nabla_{g}^{m}\text{\textgreek{f}}\big|_{h}^{2}\lesssim_{m,\text{\textgreek{e}},\text{\textgreek{h}}} & \text{\textgreek{t}}^{-2+\text{\textgreek{e}}}\Big(\mathcal{E}_{bound}^{(2,m+\lceil\frac{d+2}{2}\rceil+3k)}[\text{\textgreek{f}}](0)+\int_{0}^{\text{\textgreek{t}}}\mathcal{F}_{\text{\textgreek{h}}}^{(2,m+\lceil\frac{d+2}{2}\rceil+3k)}[\text{\textgreek{f}}](\text{\textgreek{t}})\, d\text{\textgreek{t}}\Big)+\label{eq:PointwiseDecayFirst}\\
 & +\text{\textgreek{t}}^{-1+\text{\textgreek{e}}}\int_{C^{-1}\text{\textgreek{t}}}^{\text{\textgreek{t}}}\mathcal{F}_{\text{\textgreek{h}}}^{(1,m+\lceil\frac{d+2}{2}\rceil+2k)}[\text{\textgreek{f}}](\text{\textgreek{t}})\, d\text{\textgreek{t}}+\int_{C^{-1}\text{\textgreek{t}}}^{\text{\textgreek{t}}}\mathcal{F}_{\text{\textgreek{h}}}^{(\text{\textgreek{e}},m+\lceil\frac{d+2}{2}\rceil+2k)}[\text{\textgreek{f}}](\text{\textgreek{t}})\, d\text{\textgreek{t}}.\nonumber 
\end{align}

Using (\ref{eq:FirstStepDecay}) and (\ref{eq:StartPointNewMethod})
for $p=1$ on the intervals $\{\text{\textgreek{t}}_{n}\le\bar{t}\le\text{\textgreek{t}}_{n+1}\}$,
we can bound for any integer $m\ge0$ and any $\text{\textgreek{t}}\ge0$
(for some fixed $C>1$): 
\begin{equation}
\mathcal{E}_{bound}^{(1,m)}[\text{\textgreek{f}}](\text{\textgreek{t}})\lesssim_{m,\text{\textgreek{h}}}\text{\textgreek{t}}^{-1}\Big(\mathcal{E}_{bound}^{(2,m+2k)}[\text{\textgreek{f}}](0)+\int_{0}^{\text{\textgreek{t}}}\mathcal{F}_{\text{\textgreek{h}}}^{(2,m+2k)}[\text{\textgreek{f}}](\text{\textgreek{t}})\, d\text{\textgreek{t}}\Big)+\int_{C^{-1}\text{\textgreek{t}}}^{\text{\textgreek{t}}}\mathcal{F}_{\text{\textgreek{h}}}^{(1,m+2k)}[\text{\textgreek{f}}](\text{\textgreek{t}})\, d\text{\textgreek{t}}.\label{eq:FirstStepDecay-1}
\end{equation}
 Again, using the fundamental theorem of calculus we can bound for
any any function $\text{\textgreek{y}}$ on $\mathcal{M}$ and any
$1\ll R_{1}<R_{2}$: 
\begin{equation}
\int_{\{\bar{t}=\text{\textgreek{t}}\}\cap\{r=R_{2}\}}|r^{\frac{d-1}{2}}\text{\textgreek{y}}|^{2}\, d\text{\textgreek{sv}}\lesssim\int_{\{\bar{t}=\text{\textgreek{t}}\}\cap\{r=R_{1}\}}|r^{\frac{d-1}{2}}\text{\textgreek{y}}|^{2}\, d\text{\textgreek{sv}}+\int_{\{\bar{t}=\text{\textgreek{t}}\}\cap\{R_{1}\le r\le R_{2}\}}|r^{\frac{d-1}{2}}\text{\textgreek{y}}|\big(|\partial_{v}(r^{\frac{d-1}{2}}\text{\textgreek{y}})|+r^{-2}|\partial_{u}(r^{\frac{d-1}{2}}\text{\textgreek{y}})|\big)\, dvd\text{\textgreek{sv}},\label{eq:FromFundamentalTheoremOfCalculus-1}
\end{equation}
 and thus from (\ref{eq:FirstStepDecay-1}), (\ref{eq:FromFundamentalTheoremOfCalculus-1})
and a Sobolev inequality on $\mathbb{S}^{d-1}$, we obtain the desired
decay rate for $r^{\frac{d-1}{2}}\text{\textgreek{f}}$: 
\begin{align}
\sup_{\{\bar{t}=\text{\textgreek{t}}\}}\big|r_{+}^{\frac{d-1}{2}}\nabla_{g}^{m}\text{\textgreek{f}}\big|_{h}^{2}\lesssim_{m,\text{\textgreek{h}}} & \text{\textgreek{t}}^{-1}\Big(\mathcal{E}_{bound}^{(2,m+\lceil\frac{d+2}{2}\rceil+2k)}[\text{\textgreek{f}}](0)+\int_{0}^{\text{\textgreek{t}}}\mathcal{F}_{\text{\textgreek{h}}}^{(2,m+\lceil\frac{d+2}{2}\rceil+2k)}[\text{\textgreek{f}}](\text{\textgreek{t}})\, d\text{\textgreek{t}}\Big)+\int_{C^{-1}\text{\textgreek{t}}}^{\text{\textgreek{t}}}\mathcal{F}_{\text{\textgreek{h}}}^{(1,m+\lceil\frac{d+2}{2}\rceil+2k)}[\text{\textgreek{f}}](\text{\textgreek{t}})\, d\text{\textgreek{t}}.\label{eq:PointwiseDecayRadiationField}
\end{align}

\qed

\section{\label{sec:Improved-polynomial-decay}Improved polynomial decay $\bar{t}^{-\frac{d}{2}}$
for solutions to $\square_{g}\text{\textgreek{f}}=0$ in dimensions
$d\ge3$}

In this Section, we will establish $\bar{t}{}^{-\frac{d}{2}}$ polynmial
decay estimates for solutions $\text{\textgreek{f}}$ to $\square\text{\textgreek{f}}=0$
on spacetimes $(\mathcal{M}^{d+1},g)$, $d\ge3$, satisfying Assumptions
\ref{enu:AsymptoticFlatness}-\ref{enu:SobolevInequality} and \ref{enu:IntegratedLocalEnergyDecay},
which in addition possess two ``almost Killing'' vector fields $T,K$
(not necessarily distinct) with timelike span on $\mathcal{M}\backslash\mathcal{H}^{+}$
and for which $\mathcal{H}^{+}$ becomes a non-degenerate ``almost
Killing'' horizon. These estimates extend the $\bar{t}{}^{-\frac{3}{2}+\text{\textgreek{e}}}$
decay rate established in the region $\{r\lesssim1\}$ of Schwarzschild
spacetime by Schlue in \cite{Schlue2013}.

\subsection{\label{sub:GeometricAssumtionsImprovedDecay}Further assumptions
and geometric constructions on $(\mathcal{M},g)$}

Let $(\mathcal{M}^{d+1},g)$, $d\ge3$, be a smooth Lorentzian manifold
with possibly non-empty piecewise smooth boundary $\partial\mathcal{M}$.
We assume that $(\mathcal{M},g)$ satisfies the geometric assumptions
\ref{enu:AsymptoticFlatness}-\ref{enu:SobolevInequality} stated
in Sections \ref{sec:FriedlanderRadiation} and \ref{sec:Firstdecay},
as well as Assumption \ref{enu:IntegratedLocalEnergyDecay} on integrated
local energy decay with loss of derivatives. We will use the same
notation as in Section \ref{sec:Firstdecay} for the subsets $\mathcal{N}_{af,\mathcal{M}}$,$\mathcal{H}$,
$\mathcal{H}^{+}$, $\partial_{tim}\mathcal{M}$, $\mathcal{H}_{\text{\textgreek{t}}}$,
$\partial_{tim}\mathcal{M}^{\text{\textgreek{t}}}$ of $\mathcal{M}$,
the functions $r$, $r_{+}$, $\bar{t}$ on $\mathcal{M}$ and the
vectorfields $N$, $T$, $Y$. 

Recall that in the coordinate chart $(u,v,\text{\textgreek{sv}})$
on each connected component of the region $\{r\gg1\}$, the metric
$g$ takes the form: 
\begin{align}
g= & -\Big(4+O(r^{-1-a})\Big)dvdu+r^{2}\cdot\Big(g_{\mathbb{S}^{d-1}}+h_{\mathbb{S}^{d-1}}\Big)+\Big(h^{as}(u,\text{\textgreek{sv}})+O(r^{-a})\Big)dud\text{\textgreek{sv}}+\label{eq:MetricUV-2-1}\\
 & +O(r^{-a})dvd\text{\textgreek{sv}}+4\Big(-\frac{2M(u,\text{\textgreek{sv}})}{r}+O(r^{-1-a})\Big)du^{2}+O(r^{-2-a})dv^{2},\nonumber 
\end{align}
with $h_{\mathbb{S}^{d-1}}=O(r^{-1})$.

\paragraph*{Assumptions on the vector fields $T$, $K$}

We will also assume that $(\mathcal{M},g)$ possesses two smooth vector
fields, $T$ and $K$, not necessarily distinct, such that:

\begin{MyDescription}[leftmargin=3.0em]

\item [(EG1)]{\label{enu:Transversality}The following equality holds:
$d\bar{t}(K)=d\bar{t}(T)=1$}

\item [(EG2)]{\label{enu:T,FAsymptotically}\label{enu:T,FAsymptotically}
In the region $\{r\gg1\}$, $T$ is as in Section \ref{sec:FriedlanderRadiation}
and $K=\text{\textgreek{F}}+T$, where $\text{\textgreek{F}}$ is
a generator of a rotation of $\mathbb{S}^{d-1}$ (possibly being identically
$0$).}

\item [(EG3)] {\label{enu:TimelikeSpan}The span of $T$ and $K$
everywhere on $\mathcal{M}\backslash\mathcal{H}$ contains a timelike
direction. }

\item [(EG4)] {\label{enu:NonDegenerateHorizon}The span of $\{T,K\}$
is tangential to $\mathcal{H}$. Moreover, $\mathcal{H}^{+}$ is non
degenerate with respect to $K$, in the sense that $g(K,K)=0$ and
$d(g(K,K))\neq0$ on $\mathcal{H}^{+}\cap\{\bar{t}\ge0\}$ and the
following red-shift type estimate holds for some $r_{1}>0$, any $0\le\text{\textgreek{t}}_{1}\le\text{\textgreek{t}}_{2}$,
any $l\in\mathbb{N}$ and any $\text{\textgreek{f}}\in C^{\infty}(\mathcal{M})$
(see also \cite{DafRod2}): 
\begin{align}
\sum_{j=1}^{l}\int_{\{\text{\textgreek{t}}_{1}\le\bar{t}\le\text{\textgreek{t}}_{2}\}\cap\{r\le r_{1}\}}\big|\nabla^{j}\text{\textgreek{f}}\big|_{h}^{2}\, dg\le & C_{l}\Big\{\sum_{j=1}^{l}\int_{\{\bar{t}=\text{\textgreek{t}}_{1}\}\cap\{r\le2r_{1}\}}\big|\nabla^{j}\text{\textgreek{f}}\big|_{h}^{2}\, dh_{N}+\sum_{j=1}^{l}\int_{\{\text{\textgreek{t}}_{1}\le\bar{t}\le\text{\textgreek{t}}_{2}\}\cap\{r_{1}\le r\le2r_{1}\}}\big|\nabla^{j}\text{\textgreek{f}}\big|_{h}^{2}\, dg+\label{eq:RedShiftAssumption}\\
 & +\sum_{j=0}^{l-1}\int_{\{\text{\textgreek{t}}_{1}\le\bar{t}\le\text{\textgreek{t}}_{2}\}\cap\{r\le2r_{1}\}}\big|\nabla^{j}(\square\text{\textgreek{f}})\big|_{h}^{2}\, dg\Big\}.\nonumber 
\end{align}
}

\end{MyDescription}

\emph{Convention. }Since we have not assumed that $T,K$ commute in
the near region $\{r\lesssim1\}$, iterated Lie derivatives in the
directions of $T,K$ will not necessarily commute. Hence, it will
be useful to introduce the following pointwise norms for smooth tensors
$\mathfrak{m}$ on $\mathcal{M}$ for any $l\in\mathbb{N}$ and any
two vector fields $X^{(0)},X^{(1)}$: 
\begin{equation}
\big|\mathfrak{m}\big|_{X^{(0)},X^{(1)}}^{(l)}\doteq\Big(\sum_{(e_{1},\ldots e_{d})\in\{0,1\}^{l}}\big|\mathcal{L}_{X^{(e_{1})}}\ldots\mathcal{L}_{X^{(e_{d})}}\mathfrak{m}\big|_{h}^{2}\Big)^{1/2}\label{eq:IteratedLieDerivatives}
\end{equation}
and for any $l,n\in\mathbb{N}$:
\begin{equation}
\big|\mathfrak{m}\big|_{X^{(0)},X^{(1)}}^{(l;n)}\doteq\Big(\sum_{j=0}^{n}\sum_{(e_{1},\ldots e_{d})\in\{0,1\}^{l}}\big|\nabla_{g}^{j}(\mathcal{L}_{X^{(e_{1})}}\ldots\mathcal{L}_{X^{(e_{d})}}\mathfrak{m})\big|_{h}^{2}\Big)^{1/2}\label{eq:IteratedLieDerivatives-1}
\end{equation}

\begin{MyDescription}[leftmargin=3.0em]

\item [(EG5)]{\label{enu:BoundsDeformationTensor}We assume that we
can bound for some small $\text{\textgreek{d}}_{0}>0$, any pair of
integers $j_{1},j_{2}\ge0$ and any $\text{\textgreek{t}}\ge0$: 
\begin{equation}
\sup_{\{\bar{t}=\text{\textgreek{t}}\}}\big|g\big|_{T,K}^{(j_{1};j_{2})}\lesssim_{j_{1},j_{2}}(1+\text{\textgreek{t}})^{-\max\{0,j_{1}+\text{\textgreek{d}}_{0}-1\}}.\label{eq:DeformationTensorTF}
\end{equation}
 Moreover, in the $(u,v,\text{\textgreek{sv}})$ coordinate chart
on each connected component of the region $\{r\gg1\}$ the following
precise bounds on the derivatives of $g$ are assumed to hold for
$j_{1}\ge1$ and any $(e_{1},\ldots,e_{j_{1}})\in\{0,1\}^{j_{1}}$:
\begin{align}
\mathcal{L}_{X^{(e_{1})}}\ldots\mathcal{L}_{X^{(e_{j_{1}})}}g=O(\bar{t}{}^{-(j_{1}+\text{\textgreek{d}}_{0}-1)})\Big\{ & O(r^{-1-a})dvdu+O(r)d\text{\textgreek{sv}}d\text{\textgreek{sv}}+O(1)dud\text{\textgreek{sv}}+\label{eq:DeformationTensorTAway}\\
 & +O(r^{-a})dvd\text{\textgreek{sv}}+O(r^{-1})du^{2}+O(r^{-2-a})dv^{2}\Big\}\nonumber 
\end{align}
 where $X^{(0)}=T$ and $X^{(1)}=K$.}

\end{MyDescription}
\begin{rem*}
For any $m_{0}\in\mathbb{N}$ and $\text{\textgreek{d}}_{0}\in(0,1)$,
inequalities (\ref{eq:DeformationTensorTF}) and (\ref{eq:DeformationTensorTAway})
can be relaxed to hold only for $j_{1},j_{2}$ less than some large
constant $M=M(m_{0},\text{\textgreek{d}}_{0})\in\mathbb{N}$ depending
on $m_{0}$, and then Theorem \ref{thm:ImprovedDecayEnergy} will
still hold provided $m$ is restricted to take values up to $m_{0}$.
This fact will also apply to all the assumptions regarding estimates
on the derivatives of the metric $g$ appearing in the text, and will
not be highlighted again. In all the assumptions that are stated in
this section, the number of derivatives of the metric $M=M(m_{0},\text{\textgreek{d}}_{0})$
that need to appear in the related estimates can be bounded from above
by 
\begin{equation}
m_{0}+12\lceil\text{\textgreek{d}}_{0}^{-1}\cdot d\rceil d\cdot k,
\end{equation}
 where $k$ is the number expressing the loss of derivatives in Assumption
\ref{enu:IntegratedLocalEnergyDecay}.
\end{rem*}
Let us define the vector field 
\begin{equation}
\text{\textgreek{F}}\doteq K-T\label{eq:PhiVectorField}
\end{equation}
 on $\mathcal{M}$. Recall that $\text{\textgreek{F}}$ is a rotation
vector field in the region $\{r\gg1\}$ in view of Assumption \ref{enu:T,FAsymptotically}.
Notice that (\ref{eq:DeformationTensorTF}) and (\ref{eq:DeformationTensorTAway})
hold for $\text{\textgreek{F}}$ in place of $K$. Moreover, due to
Assumption \ref{enu:Transversality} we have $d\bar{t}(\text{\textgreek{F}})=0$,
and thus $\text{\textgreek{F}}$ is tangent to the level sets of $\bar{t}$.

\paragraph*{Definition of the $K_{R_{c}}$ vector field and the $|\cdot|_{T,K,R_{c}}^{(l)}$
norm}

Let us introduce the (large) parameter $R_{c}>0$. This parameter
will be fixed after Lemma \ref{lem:EllipticEstimatesWithoutKillingError}
has been established. Fixing a smooth cut-off function $\text{\textgreek{q}}:[0,+\infty)\rightarrow[0,1]$
such that $\text{\textgreek{q}}\equiv1$ on $[0,1]$ and $\text{\textgreek{q}}\equiv0$
on $[2,+\infty)$, and defining $\text{\textgreek{q}}_{R_{c}}:\mathcal{M}\rightarrow[0,1]$
by the relation 
\begin{equation}
\text{\textgreek{q}}_{R_{c}}\doteq\text{\textgreek{q}}(\frac{r}{R_{c}}),
\end{equation}
 we introduce the following vector field: 
\begin{equation}
K_{R_{c}}\doteq\text{\textgreek{q}}_{R_{c}}K+(1-\text{\textgreek{q}}_{R_{c}})T.
\end{equation}
 Notice that for any $\text{\textgreek{t}}\ge0$, the set $\{\bar{t}=\text{\textgreek{t}}\}\cap\{g(K_{R_{c}},K_{R_{c}})>0\}$
is precompact and its closure does not intersect the horizon.

Provided $R_{c}\gg1$, due to (\ref{eq:DeformationTensorTF}) we can
bound for any pair of integers $j_{1},j_{2}\ge0$ (using the (\ref{eq:IteratedLieDerivatives-1})
norm) 
\begin{equation}
\sup_{\{\bar{t}=\text{\textgreek{t}}\}\cap\{r\notin[R_{c},2R_{c}]\}}\big|g\big|_{T,K_{R_{c}}}^{(j_{1};j_{2})}\lesssim_{j_{1},j_{2}}(1+\text{\textgreek{t}})^{-\max\{0,j_{1}+\text{\textgreek{d}}_{0}-1\}},\label{eq:DeformationTensorKtilde}
\end{equation}
while when restricted in the region $\{r\gg1\}$, the following refined
bounds hold (in the $(u,v,\text{\textgreek{sv}})$ coordinate chart
on each connected component of $\{r\gg1\}$) due to (\ref{eq:DeformationTensorTAway})
for $j_{1}\ge1$ and any $(e_{1},\ldots,e_{j_{1}})\in\{0,1\}^{j_{1}}$(with
$X^{(0)}=T$ and $X^{(1)}=K_{R_{c}}$):%
\footnote{For the calculation of the Lie derivatives in the direction of $K_{R_{c}}$
it is convenient to use the formula $\mathcal{L}_{f\cdot X}(\text{\textgreek{w}})=f\mathcal{L}_{X}(\text{\textgreek{w}})+df\cdot i_{X}(\text{\textgreek{w}})$
for any smooth function $f$, vector field $X$ and 1-form $\text{\textgreek{w}}$.%
}

\begin{itemize}

\item For $r\notin[R_{c},2R_{c}]$:
\begin{align}
\mathcal{L}_{X^{(e_{1})}}\ldots\mathcal{L}_{X^{(e_{j_{1}})}}g=O(\text{\textgreek{t}}^{-(j_{1}+\text{\textgreek{d}}_{0}-1)})\Big\{ & O(r^{-1-a})dvdu+O(r)d\text{\textgreek{sv}}d\text{\textgreek{sv}}+O(1)dud\text{\textgreek{sv}}+\label{eq:DeformationTensorTAway-1}\\
 & +O(r^{-a})dvd\text{\textgreek{sv}}+O(r^{-1})du^{2}+O(r^{-2-a})dv^{2}\Big\}.\nonumber 
\end{align}

\item For $\{R_{c}\le r\le2R_{c}\}$: 
\begin{align}
\mathcal{L}_{X^{(e_{1})}}\ldots\mathcal{L}_{X^{(e_{j_{1}})}}g=\Big(O(r^{-1})+O(\text{\textgreek{t}}^{-(j_{1}+\text{\textgreek{d}}_{0}-1)})\Big)\Big\{ & O(r)dvdu+O(r)d\text{\textgreek{sv}}d\text{\textgreek{sv}}+O(r^{2})dud\text{\textgreek{sv}}+\label{eq:DeformationTensorTAway-2}\\
 & +O(r^{2})dvd\text{\textgreek{sv}}+O(r)du^{2}+O(r)dv^{2}\Big\}.\nonumber 
\end{align}
\end{itemize}

Moreover, it will be convenient to introduce the following truncated
version of the pointwise norm (\ref{eq:IteratedLieDerivatives}) for
any smooth tensor $k$ on $\mathcal{M}$: 
\begin{equation}
\big|k\big|_{T,K,R_{c}}^{(l)}\doteq\big|\text{\textgreek{q}}_{R_{c}}k\big|_{T,K}^{(l)}+|\mathcal{L}_{T}^{l}k|_{h}.\label{eq:TruncatedNorm}
\end{equation}

\paragraph*{The Riemannian metric $h_{\text{\textgreek{t}},K_{R_{c}},\text{\textgreek{F}}}$
on the hyperboloids $\{\bar{t}=\text{\textgreek{t}}\}$}

According to Section \ref{sec:RiemannianMetric} of the Appendix,
for any $\text{\textgreek{t}}\ge0$ we can define the Riemannian metric
$h_{\text{\textgreek{t}},K_{R_{c}},\text{\textgreek{F}}}$ on the
hypersurfaces $\{\bar{t}=\text{\textgreek{t}}\}$ associated to $K_{R_{c}},\text{\textgreek{F}}$.
Due to the fact that $g(K,K)=0$ and $dg(K,K)\neq0$ on $\mathcal{H}^{+}$,
for any $\text{\textgreek{t}}\ge0$ there exists some $r_{0}(\text{\textgreek{t}})>0$
such that in a neighborhood of $\mathcal{H}_{\text{\textgreek{t}}}$
of the form $\{r\le r_{0}(\text{\textgreek{t}})\}\cap\{\bar{t}=\text{\textgreek{t}}\}$
the metric $h_{\text{\textgreek{t}},K_{R_{c}},\text{\textgreek{F}}}$
has the form: 
\begin{equation}
h_{\text{\textgreek{t}},K_{R_{c}},\text{\textgreek{F}}}=\big(r^{-1}+O(1)\big)dr^{2}+h_{\mathcal{H}_{\text{\textgreek{t}}},r},\label{eq:MetricNearHorizon-1}
\end{equation}
 where $h_{\mathcal{H}_{\text{\textgreek{t}}},r}$ is a symmetric
$(0,2)$-tensor on $\{r\le r_{0}(\text{\textgreek{t}})\}\cap\{\bar{t}=\text{\textgreek{t}}\}$,
smooth up to $\mathcal{H}_{\text{\textgreek{t}}}$, satisfying $h_{\mathcal{H}_{\text{\textgreek{t}}},r}(Y,\cdot)\equiv0$
(where the vector field $Y$ is transversal to $\mathcal{H}_{\text{\textgreek{t}}}$
and was defined in Section \ref{sub:AssumptionsFirstDecay}). Hence,
in the language of Section \ref{sec:Elliptic-estimates} of the Appendix,
$h_{\text{\textgreek{t}},K_{R_{c}},\text{\textgreek{F}}}$ corresponds
to the singular metric $h$. For notational simplicity, from now on
we will adopt the shorthand notation $h_{R_{c}}=h_{\text{\textgreek{t}},K_{R_{c}},\text{\textgreek{F}}}$.

While $h_{R_{c}}$ will be useful to define certain elliptic operators,
we will mostly measure the norms of tensors on the $\{\bar{t}=\text{\textgreek{t}}\}$
hypersurfaces with the previously defined non singular metric $h_{\text{\textgreek{t}},N}$
associated to $N$. 

In view of the aforementioned assumptions on the almost-Killing vector
fields $K,T,\text{\textgreek{F}}$, as well as the expression (\ref{eq:MetricUV-2-1})
for the asymptotics of $g$ and the expression (\ref{eq:SplittingWaveOperatorCoordinateInvariantTwo})
for $\square_{g}$, we can bound for any smooth $\text{\textgreek{f}}:\mathcal{M}\rightarrow\mathbb{C}$
and any integer $l\ge0$ and $R_{c}\gg1$: 
\begin{equation}
\big|\nabla_{h_{\text{\textgreek{t}},N}}^{l}(\text{\textgreek{D}}_{h_{R_{c}},mod}\text{\textgreek{f}})\big|_{h_{\text{\textgreek{t}},N}}^{2}\le C(\text{\textgreek{t}},R_{c})\big|\nabla_{h_{\text{\textgreek{t}},N}}^{l}(\square\text{\textgreek{f}})\big|_{h_{\text{\textgreek{t}},N}}^{2}+C(\text{\textgreek{t}},R_{c})\cdot\mathcal{T}_{T,K,R_{c}}^{(l+2)}[\text{\textgreek{f}}]+C(\text{\textgreek{t}})\text{\textgreek{q}}_{r\sim R_{c}}\sum_{j_{1}=0}^{l+1}\sum_{j_{2}=1}^{l+2-j_{1}}\Big(\big|\nabla_{h_{\text{\textgreek{t}},N}}^{j_{2}}\text{\textgreek{f}}\big|_{T,K,R_{c}}^{(j_{1})}\Big)^{2},\label{eq:UsingTheWaveEquationKilling}
\end{equation}
where 
\begin{align}
\mathcal{T}_{T,K,R_{c}}^{(m)}[\text{\textgreek{f}}]= & r_{+}^{-2}\Big(\big|\text{\textgreek{f}}\big|_{T,K,R_{c}}^{(m)}\Big)^{2}+\sum_{j_{1}+j_{2}=m-2}\Big(|\nabla_{h_{\text{\textgreek{t}},N}}^{j_{1}+1}\text{\textgreek{f}}|_{T,K,R_{c}}^{(j_{2}+1)}\Big)^{2}+\label{eq:EllipticEstimatesEnergy}\\
 & +\sum_{0\le j_{1}+j_{2}\le m-2}\big(r_{+}^{-2(m-j_{1}-j_{2}-1)}+r_{+}^{-2}\text{\textgreek{t}}^{-2(\text{\textgreek{d}}_{0}+m-j_{1}-j_{2}-2)}\big)\Big(|\nabla_{h_{\text{\textgreek{t}},N}}^{j_{1}}\text{\textgreek{f}}|_{T,K,R_{c}}^{(j_{2}+1)}\Big)^{2}+\nonumber \\
 & +\sum_{0\le j\le m-2}\big(r_{+}^{-2-2(m-2-j)}\text{\textgreek{t}}^{-2\text{\textgreek{d}}_{0}}+r_{+}^{-2}\text{\textgreek{t}}^{-2(\text{\textgreek{d}}_{0}+m-2-j)}\big)\big|\nabla_{h_{\text{\textgreek{t}},N}}^{j+1}\text{\textgreek{f}}\big|_{h_{\text{\textgreek{t}},N}}^{2},\nonumber 
\end{align}
 $\text{\textgreek{q}}_{r\sim R_{c}}$ is identically $1$ on $\{R_{c}\le r\le2R_{c}\}$
and $0$ elsewhere, and the elliptic operator $\text{\textgreek{D}}_{h_{R_{c}},mod}$
on $\{\bar{t}=\text{\textgreek{t}}\}\backslash\mathcal{H}_{\text{\textgreek{t}}}$
is defined as: 
\begin{equation}
\text{\textgreek{D}}_{h_{R_{c}},mod}\text{\textgreek{y}}=\mathfrak{w}_{\text{\textgreek{t}},K_{R_{c}},\text{\textgreek{F}}}^{-1}\cdot div_{h_{R_{c}}}\big(\mathfrak{w}_{\text{\textgreek{t}},K_{R_{c}},\text{\textgreek{F}}}\cdot d\text{\textgreek{y}}\big),\label{eq:PerturbedLaplacian-2-1}
\end{equation}
 with $\mathfrak{w}_{\text{\textgreek{t}},K_{R_{c}},\text{\textgreek{F}}}>0$
on $\{\bar{t}=\text{\textgreek{t}}\}\backslash\mathcal{H}_{\text{\textgreek{t}}}$
satisfying $\mathfrak{w}_{\text{\textgreek{t}},K_{R_{c}},\text{\textgreek{F}}}\sim_{\text{\textgreek{t}}}r^{\frac{1}{2}}$
near $\mathcal{H}^{+}$ and $\mathfrak{w}_{\text{\textgreek{t}},K_{R_{c}},\text{\textgreek{F}}}\sim_{R_{c}}1$
in the region $r\gg1$ (for its precise form, see Section \ref{sec:RiemannianMetric}
of the Appendix).%
\footnote{It is obivious that we do not need all the terms of $\mathcal{T}_{T,K,R_{c}}^{(m)}[\text{\textgreek{f}}]$
to bound the left hand side of (\ref{eq:UsingTheWaveEquationKilling}),
but since this expression will appear frequently in what follows,
we chose to introduce it here.%
}

It will also be convenient to introduce the higher order pointwise
norm 
\begin{equation}
\mathcal{T}_{T,K,R_{c}}^{(m,n)}[\text{\textgreek{f}}]\doteq\sum_{j_{1}+j_{2}\le n}\mathcal{T}_{T,K,R_{c}}^{(m)}[\nabla_{h_{\text{\textgreek{t}},N}}^{j_{1}}(T^{j_{2}}\text{\textgreek{f}})]\label{eq:EllipticEstimatesEnergyExtraDerivativesTerm}
\end{equation}
 for any pair of integers $m,n\ge0$ and $\text{\textgreek{f}}\in C^{\infty}(\mathcal{M})$.

Similarly, in the region $\{r\ge R_{0}\}$ for some $R_{0}\gg1$ we
can bound: 
\begin{align}
\big|\nabla_{h_{\text{\textgreek{t}},N}}^{l}(\text{\textgreek{D}}_{h_{\text{\textgreek{t}},T},T}\text{\textgreek{f}})\big|_{h_{\text{\textgreek{t}},N}}^{2}\le & C(\text{\textgreek{t}})\big|\nabla_{h_{\text{\textgreek{t}},N}}^{l}(\square\text{\textgreek{f}})\big|_{h_{\text{\textgreek{t}},N}}^{2}+C(\text{\textgreek{t}})\cdot\mathcal{T}_{T}^{(l+2)}[\text{\textgreek{f}}],\label{eq:UsingTheWaveEquationOnlyT}
\end{align}
where 
\begin{align}
\mathcal{T}_{T}^{(m)}[\text{\textgreek{f}}]= & r_{+}^{-2}\big|T^{m}\text{\textgreek{f}}\big|^{2}+\sum_{j_{1}+j_{2}=m-2}\big|\nabla_{h_{\text{\textgreek{t}},N}}^{j_{1}+1}(T^{j_{2}+1}\text{\textgreek{f}})\big|_{h_{\text{\textgreek{t}},N}}^{2}+\label{eq:EllipticEstimatesEnergyOnlyT}\\
 & +\sum_{0\le j_{1}+j_{2}\le m-2}\big(r_{+}^{-2(m-j_{1}-j_{2}-1)}+r_{+}^{-2}\text{\textgreek{t}}^{-2(\text{\textgreek{d}}_{0}+m-j_{1}-j_{2}-2)}\big)\big|\nabla_{h_{\text{\textgreek{t}},N}}^{j_{1}}(T^{j_{2}+1}\text{\textgreek{f}})\big|_{h_{\text{\textgreek{t}},N}}^{2}+\nonumber \\
 & +\sum_{0\le j\le m-2}\big(r_{+}^{-2-2(m-2-j)}\text{\textgreek{t}}^{-2\text{\textgreek{d}}_{0}}+r_{+}^{-2}\text{\textgreek{t}}^{-2(\text{\textgreek{d}}_{0}+m-2-j)}\big)\big|\nabla_{h_{\text{\textgreek{t}},N}}^{j+1}\text{\textgreek{f}}\big|_{h_{\text{\textgreek{t}},N}}^{2}\nonumber 
\end{align}
 and the operator $\text{\textgreek{D}}_{h_{\text{\textgreek{t}},T},T}$
on $\{\bar{t}=\text{\textgreek{t}}\}\cap\{r\ge R_{0}\}$ is defined
as: 
\begin{equation}
\text{\textgreek{D}}_{h_{\text{\textgreek{t}},T},T}\text{\textgreek{y}}=\big(\sqrt{-g(T,T)}\big)^{-1}\cdot div_{h_{\text{\textgreek{t}},T}}\big(\sqrt{-g(T,T)}\cdot d\text{\textgreek{y}}\big).\label{eq:PerturbedLaplacian-2-1-1}
\end{equation}

The following higher order pointwise norm will also be useful: 

\begin{equation}
\mathcal{T}_{T}^{(m,n)}[\text{\textgreek{f}}]\doteq\sum_{j_{1}+j_{2}\le n}\mathcal{T}_{T}^{(m)}[\nabla_{h_{\text{\textgreek{t}},N}}^{j_{1}}(T^{j_{2}}\text{\textgreek{f}})]\label{eq:EllipticEstimatesEnergyExtraDerivativesTerm-1}
\end{equation}
 for any pair of integers $m,n\ge0$ and $\text{\textgreek{f}}\in C^{\infty}(\mathcal{M})$.

\paragraph*{Assumptions on the uniformity of the degenerate elliptic estimates
on the hyperboloids $\{\bar{t}=\text{\textgreek{t}}\}$}

Let us introduce the functions $r_{hor},r_{tim}:\mathcal{M}\rightarrow[0,1]$
by the relations 
\begin{equation}
r_{hor}(x)=\frac{dist_{h}(x,\mathcal{H})}{1+dist_{h}(x,\mathcal{H})}
\end{equation}
 and 
\begin{equation}
r_{tim}(x)=\frac{dist_{h}(x,\partial_{tim}\mathcal{M})}{1+dist_{h}(x,\partial_{tim}\mathcal{M})}.
\end{equation}
 Notice that this definition does not guarantee that $r_{hor}$ and
$r_{tim}$ are smooth functions away from $\partial\mathcal{M}$.
For this reason, we will mollify $r_{hor}$ and $r_{tim}$ way from
$\partial\mathcal{M}$ so that they are smooth functions on $\mathcal{M}$,
and we will replace the original $r_{hor}$ and $r_{tim}$ by the
corresponding mollified functions.

Using Propositions \ref{prop:DegenerateEllipticEstimates} and \ref{Prop:ControlOfTheAngularTermsAway}
of the Appendix, together with (\ref{eq:MetricNearHorizon-1}), (\ref{eq:UsingTheWaveEquationKilling})
and (\ref{eq:UsingTheWaveEquationOnlyT}), we can deduce the following
elliptic estimate:
\begin{lem*}
\label{lem:DegenerateEllipticEstimates-1}For any $l\in\mathbb{N}$
with $l\le\lfloor\frac{d+1}{2}\rfloor$, any $n_{0}\in\mathbb{N}$,
any $\text{\textgreek{b}}\in(-\bar{\text{\textgreek{d}}}_{n_{0}},1)$
(for some $\bar{\text{\textgreek{d}}}_{n_{0}}>0$ depending on $n_{0}$
and $\text{\textgreek{t}}$) and any $R_{0}\gg1$ large in terms of
the geometry of $(\mathcal{M},g)$, we can bound for any $\text{\textgreek{f}}\in C^{\infty}(\mathcal{M})$
satisfying for all $j_{1}+j_{2}\le l+n_{0}$ the finite radiation
field condition $\limsup_{r\rightarrow+\infty}\big|r^{\frac{d-1}{2}+j_{1}}\nabla_{h_{\text{\textgreek{t}},N}}^{j_{1}}(N^{j_{2}}\text{\textgreek{f}})\big|_{h_{\text{\textgreek{t}},N}}<+\infty$:
\begin{equation}
\begin{split}\sum_{n=0}^{n_{0}}\sum_{j=0}^{l-1}\int_{\{\bar{t}=\text{\textgreek{t}}\}}r_{+}^{-\text{\textgreek{b}}-2j}| & \nabla_{h_{\text{\textgreek{t}},N}}^{n+l-j}\text{\textgreek{f}}|_{\big(1-\log(r_{tim})\big)\cdot h_{R_{c}}}^{2}\, dh_{N}\le C_{\text{\textgreek{b}},n_{0}}(\text{\textgreek{t}},R_{c})\int_{\{\bar{t}=\text{\textgreek{t}}\}}r_{+}^{-\text{\textgreek{b}}}\Big\{\sum_{n=0}^{n_{0}}\big|\nabla_{h_{\text{\textgreek{t}},N}}^{n+l-2}(\square\text{\textgreek{f}})\big|_{\big(1-\log(r_{tim})\big)\cdot h_{\text{\textgreek{t}},N}}^{2}+\mathcal{T}_{T,K,R_{c}}^{(l,n_{0})}[\text{\textgreek{f}}]\Big\}\, dh_{N}+\\
 & +C_{\text{\textgreek{b}},n_{0}}(\text{\textgreek{t}})\sum_{n=0}^{n_{0}}\sum_{j_{1}=0}^{l-1}\sum_{j_{2}=1}^{l-j_{1}}\int_{\{\bar{t}=\text{\textgreek{t}}\}\cap\{R_{c}\le r\le2R_{c}\}}r_{+}^{-\text{\textgreek{b}}}\Big(\big|\nabla_{h_{\text{\textgreek{t}},N}}^{j_{2}+n}\text{\textgreek{f}}\big|_{T,K,R_{c}}^{(j_{1})}\Big)^{2}\, dh_{N}+\\
 & +C_{\text{\textgreek{b}},n_{0}}(\text{\textgreek{t}})\sum_{j=0}^{1}\max\Big\{-Re\big\{\int_{\partial_{tim}\mathcal{M}^{\text{\textgreek{t}}}}h_{\partial_{tim}\mathcal{M}^{\text{\textgreek{t}}}}\Big(\nabla_{h_{\partial_{tim}\mathcal{M}^{\text{\textgreek{t}}}}}^{j}(Y\text{\textgreek{f}}),\nabla_{h_{\partial_{tim}\mathcal{M}^{\text{\textgreek{t}}}}}^{j}\bar{\text{\textgreek{f}}}\Big)\, dh_{\partial_{tim}\mathcal{M}^{\text{\textgreek{t}}}}\big\},0\Big\}
\end{split}
\label{eq:EllipticEstimatesDegenerate}
\end{equation}
 and for any $0<\text{\textgreek{e}}\ll1-\text{\textgreek{b}}$ (provided
that $R_{c}\gg R_{0}$):
\begin{equation}
\begin{split}\sum_{n=0}^{n_{0}}\sum_{j=1}^{l}\int_{\{\bar{t}=\text{\textgreek{t}}\}\cap\{R_{c}\le r\le2R_{c}\}} & r_{+}^{2(j-1)-\text{\textgreek{b}}+\frac{\text{\textgreek{e}}}{2}}\big|\nabla_{g}^{j+n}\text{\textgreek{f}}\big|_{h}^{2}\, dh_{N}\le\\
\le & C_{\text{\textgreek{b},\textgreek{e}},n_{0}}(\text{\textgreek{t}},R_{c})\int_{\{\bar{t}=\text{\textgreek{t}}\}\cap\{r\ge2R_{0}\}}r_{+}^{-\text{\textgreek{b}}}\Big\{\sum_{n=0}^{n_{0}}\big|\nabla_{h_{\text{\textgreek{t}},N}}^{n+l-2}(\square\text{\textgreek{f}})\big|_{h_{\text{\textgreek{t}},N}}^{2}+\mathcal{T}_{T}^{(l,n_{0})}[\text{\textgreek{f}}]\Big\}\, dh_{N}+\\
 & +C_{\text{\textgreek{b}},\text{\textgreek{e}},n_{0}}(\text{\textgreek{t}})R_{c}^{-\frac{\text{\textgreek{e}}}{2}}\sum_{j=1}^{n_{0}+l}\int_{\{\bar{t}=\text{\textgreek{t}}\}\cap\{2R_{0}\le r\le4R_{0}\}}|\nabla_{g}^{j}\text{\textgreek{f}}|_{h}^{2}\, dh_{N}.
\end{split}
\label{eq:EstimateFarAwayForAngularBulk-1-1}
\end{equation}
 
\end{lem*}
Our extra assumptions will be the following uniformity conditions
on the elliptic estimates (\ref{eq:EllipticEstimatesDegenerate})
and the Sobolev-type estimates of Lemmas \ref{lem:Gagliardo-NirenbergAway}
and \ref{lem:Gagliardo-NirenbergDegenerate}:

\begin{MyDescription}[leftmargin=3.0em]

\item [(EG6)]{\label{enu:UniformityOfDegenerateEllipticEstimates}The
constants $\bar{\text{\textgreek{d}}}_{n_{0}}$, $C_{\text{\textgreek{b}},n_{0}}(\text{\textgreek{t}})$,
$C_{\text{\textgreek{b}},\text{\textgreek{e}},n_{0}}(\text{\textgreek{t}})$
and $C_{\text{\textgreek{b}},\text{\textgreek{e}},n_{0}}(\text{\textgreek{t}},R_{c})$
in Lemma \ref{lem:DegenerateEllipticEstimates-1} can be chosen not
to depend on $\text{\textgreek{t}}$.}

\item [(EG7)]{\label{enu:UniformityOfGagliardoNirenberg}The constants
in the estimates of Lemmas \ref{lem:Gagliardo-NirenbergAway} and
\ref{lem:Gagliardo-NirenbergDegenerate} can be chosen not to depend
on $\text{\textgreek{t}}$.}

\end{MyDescription}

We will also need the following assumption regarding the volume of
the domains $\{\bar{t}=\text{\textgreek{t}}\}\cap\{r\le1\}$ as $\text{\textgreek{t}}\rightarrow+\infty$:

\begin{MyDescription}[leftmargin=3.0em]

\item [(EG8)]{\label{enu:BoundednessVolume}There exists some $C>0$
such that for any $\text{\textgreek{t}}\ge0$ the $h_{\text{\textgreek{t}},N}$-volume
of the $\{\bar{t}=\text{\textgreek{t}}\}\cap\{r\le1\}$ region is
uniformly bounded in $\text{\textgreek{t}}$: 
\begin{equation}
Vol_{h_{\text{\textgreek{t}},N}}\big(\{\bar{t}=\text{\textgreek{t}}\}\cap\{r\le1\}\big)\le C.\label{eq:UniformBoundVolume}
\end{equation}
}

\end{MyDescription}

\paragraph*{The higher order integrated local energy decay estimate assumption}

In Section \ref{sub:IntegratedLocalEnergyDecay} we defined the class
$\mathcal{C}_{ILED}$ of admissible boundary conditions on $\partial_{tim}\mathcal{M}$
for which the integrated local energy decay statement (\ref{eq:IntegratedLocalEnergyDecayImprovedDecay})
holds, and we assumed that $\mathcal{C}_{ILED}\neq\emptyset$. However,
since the class $C_{ILED}$ was not necessarily closed under differentiation
with respect to the vector fields $T$, $K$ introduced in this section,
we do not know a priori that (\ref{eq:IntegratedLocalEnergyDecayImprovedDecay})
holds for $T$ and $K$ derivatives of functions belonging to $C_{ILED}$.
the In this section, therefore, we will assume that a slightly stronger
form of the integrated local energy decay statement (\ref{eq:IntegratedLocalEnergyDecayImprovedDecay})
(including also higher order $T$ and $K$ derivatives of $\text{\textgreek{f}}$)
holds on $(\mathcal{M},g)$, which further restricts the class of
allowed boundary conditions imposed on $\partial_{tim}\mathcal{M}$:

\begin{MyDescription}[leftmargin=3.0em]

\item [(ILED2)]{\label{enu:IntegratedLocalEnergyDecayImproved}We
assume that there exists a non empty class $\mathcal{C}_{ILED}^{(T,K)}$
of boundary conditions on $\partial_{tim}\mathcal{M}$, which is contained
in the class of admissible boundary conditions $\mathcal{C}_{adm}$,
so that the following integrated local energy decay statement holds
on $(\mathcal{M},g)$: There exists an integer $k\ge0$ such that
for any $R,R_{f}>0$, any integers $m\ge0$ and $j_{1},j_{2}\ge0$,
any $0<\text{\textgreek{h}}<a$, any smooth $\text{\textgreek{f}}:\mathcal{M}\rightarrow\mathbb{C}$
solving $\square_{g}\text{\textgreek{f}}=F$ satisfying boundary conditions
on $\partial_{tim}\mathcal{M}$ belonging to the class $\mathcal{C}_{ILED}^{(T,K)}$
and any $0\le\text{\textgreek{t}}_{1}\le\text{\textgreek{t}}_{2}$,
we can bound
\begin{equation}
\begin{split}\sum_{j=0}^{m}\int_{\{\text{\textgreek{t}}_{1}\le\bar{t}\le\text{\textgreek{t}}_{2}\}\cap\{r\le R\}} & |\nabla_{g}^{j}(T^{j_{1}}K_{R_{c}}^{j_{2}}\text{\textgreek{f}})|_{h}^{2}+\sum_{j=1}^{m}\int_{\{\text{\textgreek{t}}_{1}\le\bar{t}\le\text{\textgreek{t}}_{2}\}\cap\partial_{tim}\mathcal{M}}|\nabla_{g}^{j}(T^{j_{1}}K_{R_{c}}^{j_{2}}\text{\textgreek{f}})|_{h}^{2}\le\\
\le & C_{m,\text{\textgreek{h}},i_{1},i_{2}}(R,R_{f})\cdot\sum_{j=0}^{m+k-1}\Big(\int_{\{\bar{t}=\text{\textgreek{t}}_{1}\}}J_{\text{\textgreek{m}}}^{N}(N^{j}(T^{j_{1}}K_{R_{c}}^{j_{2}}\text{\textgreek{f}}))\bar{n}^{\text{\textgreek{m}}}+\int_{\{\text{\textgreek{t}}_{1}\le\bar{t}\le\text{\textgreek{t}}_{2}\}}r_{+}^{1+\text{\textgreek{h}}}|\nabla_{g}^{j}\big(\square_{g}(T^{j_{1}}K_{R_{c}}^{j_{2}}\text{\textgreek{f}})\big)|_{h}^{2}\, dh\Big)+\\
 & +C_{m,\text{\textgreek{h}}i_{1},i_{2}}\sum_{i=0}^{m-1}\sum_{i_{1}+i_{2}=i}\int_{\{\text{\textgreek{t}}_{1}\le\bar{t}\le\text{\textgreek{t}}_{1}\}\cap\{r\ge R_{f}\}}\bar{t}^{-\text{\textgreek{d}}_{0}}r_{+}^{-1}\Big(\big|\nabla_{h_{\text{\textgreek{t}},N}}^{i_{1}+1}(T^{i_{2}+j_{1}}K_{R_{c}}^{j_{2}}\text{\textgreek{f}})\big|^{2}+r_{+}^{-2}\big|\nabla_{h_{\text{\textgreek{t}},N}}^{i_{1}}(T^{i_{2}+j_{1}}K_{R_{c}}^{j_{2}}\text{\textgreek{f}})\big|^{2}+\\
 & \hphantom{\hphantom{+C_{m,\text{\textgreek{h}}}\sum_{i=0}^{m-1}\sum_{i_{1}+i_{2}=i}\int_{\{\text{\textgreek{t}}_{1}\le\bar{t}\le\text{\textgreek{t}}_{1}\}\cap\{r\ge R_{f}\}}\bar{t}^{-\text{\textgreek{d}}_{0}}r_{+}^{-1}\Big(}}+r_{+}^{-2}\big|T^{i+j_{1}+1}K_{R_{c}}^{j_{2}}\text{\textgreek{f}}\big|^{2}\Big)\, dh.
\end{split}
\label{eq:IntegratedLocalEnergyDecayImprovedDecay-2}
\end{equation}
}

\end{MyDescription}
\begin{rem*}
Notice that (\ref{eq:IntegratedLocalEnergyDecayImprovedDecay-2})
follows trivially from (\ref{eq:IntegratedLocalEnergyDecayImprovedDecay})
(provided that the last term in the right hand side of (\ref{eq:IntegratedLocalEnergyDecayImprovedDecay})
is replaced by
\begin{equation}
C_{m,\text{\textgreek{h}}}\sum_{j=0}^{m-1}\sum_{j_{1}+j_{2}=j}\int_{\{\text{\textgreek{t}}_{1}\le\bar{t}\le\text{\textgreek{t}}_{1}\}\cap\{r\ge R_{f}\}}\bar{t}^{-\text{\textgreek{d}}_{0}}r_{+}^{-1}\Big(\big|\nabla_{h_{\text{\textgreek{t}},N}}^{j_{1}+1}(N^{j_{2}}\text{\textgreek{f}})\big|^{2}+r_{+}^{-2}\big|\nabla_{h_{\text{\textgreek{t}},N}}^{j_{1}}(N^{j_{2}}\text{\textgreek{f}})\big|^{2}+r_{+}^{-2}\big|N^{j+1}\text{\textgreek{f}}\big|^{2}\Big)\, dg,
\end{equation}
which is consistent with the remark below Lemma \ref{lem:MorawetzDrLemmaHyperboloids}
in view of the assumptions on the vector field $T$) in the case when
$\partial_{tim}\mathcal{M}=\emptyset$, or in the case when $T$ and
$K$ are tangential to $\partial_{tim}\mathcal{M}$ and $\mathcal{C}_{ILED}$
is restricted to include only the Dirichlet boundary conditions. In
these cases, $\mathcal{C}_{ILED}^{(T,K)}=\mathcal{C}_{ILED}$. Furthermore,
the results of this section also hold if one replaces Assumption \ref{enu:IntegratedLocalEnergyDecayImproved}
by a pair of alternative integrated local energy decay assumptions
similar to (\ref{eq:IntegratedLocalEnergyDecayImprovedDecayAltLoss})
and (\ref{eq:IntegratedLocalEnergyDecayImprovedDecayAlt}), but we
will not pursue this issue any further.
\end{rem*}

\subsection{\label{sub:ShorthandEnergyNorms}Shorthand notation for energy norms}

In order to state and prove Theorem \ref{thm:ImprovedDecayEnergy}
conveniently, it will be useful to introduce some shorthand notation
for a variety of energy norms on the hyperboloids $\{\bar{t}=\text{\textgreek{t}}\}$.
More precisely, we will make use of the following notations for $p\in\mathbb{R}$,
$\text{\textgreek{e}}\in(0,1)$ and $q,l,m,k\in\mathbb{N}$: 
\begin{enumerate}
\item The following energy norms will appear when using the $r^{p}$-weighted
estimates of Sections \ref{sec:The-new-method} and \ref{sec:The-improved--hierarchy}:
\begin{equation}
\mathcal{E}_{bulk,R}^{(p,q,m)}[\text{\textgreek{f}}](\text{\textgreek{t}})\doteq\sum_{j=0}^{m}\sum_{0\le j_{1}+j_{2}\le q+j}\int_{\{\bar{t}=\text{\textgreek{t}}\}\cap\{r\ge R\}}r^{p-2(q+j-j_{1})}\big|\nabla_{h_{\text{\textgreek{t}},N}}^{j_{1}}T^{j_{2}}(\text{\textgreek{W}}\text{\textgreek{f}})\big|_{h_{\text{\textgreek{t}},N}}^{2}\,\text{\textgreek{W}}^{-2}dh_{N}\label{eq:BulkNewMethod}
\end{equation}
 and
\begin{align}
\mathcal{E}_{bound,R}^{(p,q,m)}[\text{\textgreek{f}}](\text{\textgreek{t}})\doteq\sum_{j=0}^{m}\sum_{0\le j_{1}+j_{2}\le q+j-1} & \int_{\{\bar{t}=\text{\textgreek{t}}\}\cap\{r\ge R\}}r^{p-2(q+j-j_{1}-1)}\Big(\big|\nabla_{h_{\text{\textgreek{t}},N}}^{j_{1}}T^{j_{2}}\partial_{v}(\text{\textgreek{W}}\text{\textgreek{f}})\big|_{h_{\text{\textgreek{t}},N}}^{2}+\label{eq:BoundaryNewMethod}\\
 & +r^{-2}\big|\nabla_{h_{\text{\textgreek{t}},N}}^{j_{1}+1}T^{j_{2}}(\text{\textgreek{W}}\text{\textgreek{f}})\big|_{h_{\text{\textgreek{t}},N}}^{2}+r^{-4}\big|\nabla_{h_{\text{\textgreek{t}},N}}^{j_{1}}T^{j_{2}+1}(\text{\textgreek{W}}\text{\textgreek{f}})\big|_{h_{\text{\textgreek{t}},N}}^{2}+\nonumber \\
 & +r^{-4}\big|\nabla_{h_{\text{\textgreek{t}},N}}^{j_{1}}T^{j_{2}}(\text{\textgreek{W}}\text{\textgreek{f}})\big|_{h_{\text{\textgreek{t}},N}}^{2}\Big)\,\text{\textgreek{W}}^{-2}dh_{N}.\nonumber 
\end{align}

\item The following weighted non degenerate energy norms on $\{\bar{t}=\text{\textgreek{t}}\}$
will also appear frequently (fixing some $R_{1}$ large in terms of
the geometry of $(\mathcal{M},g)$): 
\begin{equation}
\mathcal{E}_{en}^{(p,q,m)}[\text{\textgreek{f}}](\text{\textgreek{t}})\doteq\sum_{j=0}^{m}\sum_{1\le j_{2}+j_{3}\le q+j-1}\int_{\{\bar{t}=\text{\textgreek{t}}\}}r_{+}^{p-2(q+j-1-j_{2}-j_{3})}\Big(\big|\nabla_{h_{\text{\textgreek{t}},N}}^{j_{2}+1}(N^{j_{3}}\text{\textgreek{f}})\big|_{h_{\text{\textgreek{t}},N}}^{2}+r_{+}^{-2}\big|N^{q+j}\text{\textgreek{f}}\big|^{2}\Big)\, dh_{N},\label{eq:EnergyHyperboloids}
\end{equation}
\begin{equation}
\mathcal{E}_{bulk}^{(p,q,m)}[\text{\textgreek{f}}](\text{\textgreek{t}})\doteq\mathcal{E}_{bulk,R_{1}}^{(p,q,m)}[\text{\textgreek{f}}](\text{\textgreek{t}})+\sum_{j=0}^{m}\sum_{0\le j_{1}+j_{2}\le q+j}\int_{\{\bar{t}=\text{\textgreek{t}}\}\cap\{r\le R_{1}\}}\big|\nabla_{h_{\text{\textgreek{t}},N}}^{j_{1}}N^{j_{2}}\text{\textgreek{f}}\big|_{h_{\text{\textgreek{t}},N}}^{2}\, dh_{N}\label{eq:BulkNormEverywhere}
\end{equation}
and 
\begin{equation}
\mathcal{E}_{bound}^{(p,q,m)}[\text{\textgreek{f}}](\text{\textgreek{t}})\doteq\mathcal{E}_{bound,R_{1}}^{(p,q,m)}[\text{\textgreek{f}}](\text{\textgreek{t}})+\sum_{j=0}^{m}\sum_{0\le j_{1}+j_{2}\le q+j}\int_{\{\bar{t}=\text{\textgreek{t}}\}\cap\{r\le R_{1}\}}\big|\nabla_{h_{\text{\textgreek{t}},N}}^{j_{1}}N^{j_{2}}\text{\textgreek{f}}\big|_{h_{\text{\textgreek{t}},N}}^{2}\, dh_{N}.\label{eq:BoundaryNormEverywhere}
\end{equation}

\item The following energy norm (associated to the vector fields $K_{R_{c}},\text{\textgreek{F}}$)
which degenerates at $\mathcal{H}$ and $\partial_{tim}\mathcal{M}$
will also appear: 
\begin{equation}
\mathcal{E}_{en,deg}^{(p,q,m)}[\text{\textgreek{f}}](\text{\textgreek{t}})\doteq\sum_{j=0}^{m}\sum_{1\le j_{2}+j_{3}\le l+j-1}\int_{\{\bar{t}=\text{\textgreek{t}}\}}r_{+}^{p-2(l+j-1-j_{2}-j_{3})}\Big(\big|\nabla_{h_{\text{\textgreek{t}},N}}^{j_{2}+1}(K_{R_{c}}^{j_{3}}\text{\textgreek{f}})\big|_{\big(1-\log(r_{tim})\big)\cdot h_{R_{c}}}^{2}+r_{+}^{-2}\big|K_{R_{c}}^{l+j}\text{\textgreek{f}}\big|^{2}\Big)\, dh_{N}.\label{eq:DegenerateEnergyNorm}
\end{equation}

\item Finally, we will make use of the following spacetime norms for the
source terms, for some $C=C_{l,q,m,k}>0$ that will be fixed in the
statement of Theorem \ref{thm:ImprovedDecayEnergy}: 
\begin{align}
\mathcal{F}_{\text{\textgreek{h}},\text{\textgreek{e}}}^{(2l,q,m,k)}[F](\text{\textgreek{t}})= & \text{\textgreek{t}}^{-2l+\text{\textgreek{e}}}\sum_{j=0}^{m+(2l+1)k-1}\sum_{j_{1}+j_{2}=j}\sum_{i_{1}+i_{2}\le q-1}\int_{\{0\le\bar{t}\le\text{\textgreek{t}}\}}r^{3+2i_{1}}\big|\nabla_{h_{\text{\textgreek{t}},N}}^{j_{1}+i_{1}}(N^{j_{2}+i_{2}}F)\big|_{h_{\text{\textgreek{t}},N}}^{2}\, dg+\\
 & +\sum_{s=1}^{l-1}\text{\textgreek{t}}^{-2s+\text{\textgreek{e}}}\sum_{j=0}^{m+(2s+1)k-1}\sum_{j_{1}+j_{2}=j}\sum_{i_{1}=0}^{q-(l-s)-1}\sum_{i_{2}=0}^{l-s}\int_{\{C^{-1}\text{\textgreek{t}}\le\bar{t}\le\text{\textgreek{t}}\}}r^{3+2i_{1}}\big|\nabla_{h_{\text{\textgreek{t}},N}}^{j_{1}+i_{1}}(N^{j_{2}+i_{2}}F)\big|_{h_{\text{\textgreek{t}},N}}^{2}\, dg+\nonumber \\
 & +\sum_{s=1}^{l}\text{\textgreek{t}}^{-2s+1+\text{\textgreek{e}}}\sum_{j=0}^{m+2sk-1}\sum_{j_{1}+j_{2}=j}\sum_{i_{1}=0}^{q-(l-s)-1}\sum_{i_{2}=0}^{l-s}\int_{\{C^{-1}\text{\textgreek{t}}\le\bar{t}\le\text{\textgreek{t}}\}}r^{2+2i_{1}}\big|\nabla_{h_{\text{\textgreek{t}},N}}^{j_{1}+i_{1}}(N^{j_{2}+i_{2}}F)\big|_{h_{\text{\textgreek{t}},N}}^{2}\, dg+\nonumber \\
 & +\sum_{j=0}^{m+2k-1}\sum_{j_{1}+j_{2}=j}\sum_{i_{1}=0}^{q-l}\sum_{i_{2}=0}^{l-1}\int_{\{C^{-1}\text{\textgreek{t}}\le\bar{t}\le\text{\textgreek{t}}\}}(r^{1+\text{\textgreek{e}}+2i_{1}}+r^{1+\text{\textgreek{h}}})\big|\nabla_{h_{\text{\textgreek{t}},N}}^{j_{1}+i_{1}}(N^{j_{2}+i_{2}}F)\big|_{h_{\text{\textgreek{t}},N}}^{2}\, dg\nonumber 
\end{align}
and 
\begin{align}
\mathcal{F}_{deg,\text{\textgreek{h}}}^{(2l,q,m,k)}[F](\text{\textgreek{t}})= & \text{\textgreek{t}}^{-2l}\sum_{j=0}^{m+(2l+1)k-1}\sum_{j_{1}+j_{2}=j}\sum_{i_{1}+i_{2}\le q-1}\int_{\{0\le\bar{t}\le\text{\textgreek{t}}\}}r^{3+2i_{1}}\big|\nabla_{h_{\text{\textgreek{t}},N}}^{j_{1}+i_{1}}(K_{R_{c}}^{j_{2}+i_{2}}F)\big|_{h_{\text{\textgreek{t}},N}}^{2}\, dg+\\
 & +\sum_{s=1}^{l-1}\text{\textgreek{t}}^{-2s}\sum_{j=0}^{m+(2s+1)k-1}\sum_{j_{1}+j_{2}=j}\sum_{i_{1}=0}^{q-(l-s)-1}\sum_{i_{2}=0}^{l-s}\int_{\{C^{-1}\text{\textgreek{t}}\le\bar{t}\le\text{\textgreek{t}}\}}r^{3+2i_{1}}\big|\nabla_{h_{\text{\textgreek{t}},N}}^{j_{1}+i_{1}}(K_{R_{c}}^{j_{2}+i_{2}}F)\big|_{h_{\text{\textgreek{t}},N}}^{2}\, dg+\nonumber \\
 & +\sum_{s=1}^{l}\text{\textgreek{t}}^{-2s+1}\sum_{j=0}^{m+2sk-1}\sum_{j_{1}+j_{2}=j}\sum_{i_{1}=0}^{q-(l-s)-1}\sum_{i_{2}=0}^{l-s}\int_{\{C^{-1}\text{\textgreek{t}}\le\bar{t}\le\text{\textgreek{t}}\}}r^{2+2i_{1}}\big|\nabla_{h_{\text{\textgreek{t}},N}}^{j_{1}+i_{1}}(K_{R_{c}}^{j_{2}+i_{2}}F)\big|_{h_{\text{\textgreek{t}},N}}^{2}\, dg+\nonumber \\
 & +\sum_{j=0}^{m+2k-1}\sum_{j_{1}+j_{2}=j}\sum_{i_{1}=0}^{q-l}\sum_{i_{2}=0}^{l-1}\int_{\{C^{-1}\text{\textgreek{t}}\le\bar{t}\le\text{\textgreek{t}}\}}(r^{1+2i_{1}}+r^{1+\text{\textgreek{h}}})\big|\nabla_{h_{\text{\textgreek{t}},N}}^{j_{1}+i_{1}}(K_{R_{c}}^{j_{2}+i_{2}}F)\big|_{h_{\text{\textgreek{t}},N}}^{2}\, dg.\nonumber 
\end{align}

\end{enumerate}

\subsection{Statement of the results on improved polynomial decay}

In this class of spacetimes $(\mathcal{M},g)$ we will establish the
following result:
\begin{thm}
\label{thm:ImprovedDecayEnergy}Assume that $(\mathcal{M}^{d+1},g)$,
$d\ge3$, satisfies Assumptions \ref{enu:AsymptoticFlatness}--\ref{enu:SobolevInequality},
\ref{enu:Transversality}--\ref{enu:BoundednessVolume} and \ref{enu:IntegratedLocalEnergyDecayImproved}.
Then for any smooth function $\text{\textgreek{f}}$ on $\mathcal{M}$
solving $\square_{g}\text{\textgreek{f}}=F$ on $J^{+}(\{\bar{t}=0\})$
with suitably decaying initial data on $\{\bar{t}=0\}$ and satisfying
boundary conditions on $\partial_{tim}\mathcal{M}$ belonging to the
class $\mathcal{C}_{ILED}^{(T,K)}$ and such that 
\begin{equation}
Re\Big\{\int_{\partial_{tim}\mathcal{M}^{\text{\textgreek{t}}}}Y\text{\textgreek{f}}\cdot\bar{\text{\textgreek{f}}}\, dh_{\partial_{tim}\mathcal{M}^{\text{\textgreek{t}}}}\Big\}\ge0\mbox{ and }Re\Big\{\int_{\partial_{tim}\mathcal{M}^{\text{\textgreek{t}}}}h_{\partial_{tim}\mathcal{M}^{\text{\textgreek{t}}}}\Big(\nabla_{h_{\partial_{tim}\mathcal{M}^{\text{\textgreek{t}}}}}(Y\text{\textgreek{f}}),\nabla_{h_{\partial_{tim}\mathcal{M}^{\text{\textgreek{t}}}}}\bar{\text{\textgreek{f}}}\Big)\, dh_{\partial_{tim}\mathcal{M}^{\text{\textgreek{t}}}}\Big\}\ge0,\label{eq:BoundaryConditionImprovedDecay}
\end{equation}
 the following bounds hold for any integer $1\le q\le\lfloor\frac{d+1}{2}\rfloor$,
any $0<\text{\textgreek{e}}\ll\text{\textgreek{d}}_{0}$, any integer
$m\ge1$ and any $\text{\textgreek{t}}\ge0$, provided (\ref{eq:NormFinitenessFriedlanderRadiation})
holds for all $0\le j\le q+\lceil\frac{d}{2}\rceil+m+\lceil\text{\textgreek{d}}_{0}^{-1}\cdot2(q-1)\rceil(3q+1)\cdot k$:
\begin{align}
\mathcal{E}_{en}^{(0,q,m)}[\text{\textgreek{f}}](\text{\textgreek{t}})+\mathcal{E}_{bound}^{(\text{\textgreek{e}},q,m)}[\text{\textgreek{f}}](\text{\textgreek{t}})+\int_{\text{\textgreek{t}}}^{+\infty}\mathcal{E}_{en}^{(-1+\text{\textgreek{e}},q,m)}[\text{\textgreek{f}}](s)\, ds\lesssim_{m,\text{\textgreek{e}}} & \text{\textgreek{t}}^{-2q+C_{m}\text{\textgreek{e}}}\mathcal{E}_{bound}^{(2q,q,m+\lceil\text{\textgreek{d}}_{0}^{-1}\cdot2(q-1)\rceil(3q+1)\cdot k)}[\text{\textgreek{f}}](0)+\label{eq:ImprovedNonDegenerateEnergyDecay}\\
 & +\mathcal{F}_{\text{\textgreek{e}},\text{\textgreek{e}}}^{(2q,q,m+\lceil\text{\textgreek{d}}_{0}^{-1}\cdot2(q-1)\rceil(3q+1)\cdot k,k)}[F](\text{\textgreek{t}})\nonumber 
\end{align}
and 
\begin{equation}
\mathcal{E}_{en,deg}^{(0,q,m)}[\text{\textgreek{f}}](\text{\textgreek{t}})\lesssim_{m,\text{\textgreek{e}}}\text{\textgreek{t}}^{-2q}\mathcal{E}_{bound}^{(2q,q,m+\lceil\text{\textgreek{d}}_{0}^{-1}\cdot2(q-1)\rceil(3q+1)\cdot k)}[\text{\textgreek{f}}](0)+\mathcal{F}_{deg,\text{\textgreek{e}}}^{(2q,q,m+\lceil\text{\textgreek{d}}_{0}^{-1}\cdot2(q-1)\rceil(3q+1)\cdot k,k)}[F](\text{\textgreek{t}}).\label{eq:ImprovedEnergyDecay}
\end{equation}
 See Section \ref{sub:ShorthandEnergyNorms} for the notations on
the energy norms used.\end{thm}
\begin{rem*}
Notice that the condition (\ref{eq:BoundaryConditionImprovedDecay})
is satisfied when $\text{\textgreek{f}}$ is subject to Dirichlet
or Neumann boundary conditions. Theorem \ref{thm:ImprovedDecayEnergy}
still holds if one replaces the condition (\ref{eq:BoundaryConditionImprovedDecay})
with any boundary condition for which Lemma \ref{lem:ILEDCommuted}
can still be established.

Let us also remark in the case the integrated local energy decay statement
in Assumption \ref{enu:IntegratedLocalEnergyDecayImproved} does not
lose derivatives (i.\,e.~$k=0$), we can replace Assumption \ref{enu:BoundsDeformationTensor}
on the $\bar{t}^{-\text{\textgreek{d}}_{0}}$ decay of the deformation
tensors of $T$ and $K$ with the following $O(\text{\textgreek{e}}_{0})$-smallness
assumption:

\medskip{}

\noindent \emph{Uniform~smallness~of~the~deformation~tensor:}
There exists some (small) $\text{\textgreek{e}}_{0}>0$ so that: 
\begin{equation}
\sup_{\{\bar{t}=\text{\textgreek{t}}\}}\Big(\big|\mathcal{L}_{T}g\big|+\big|\mathcal{L}_{K}g\big|\Big)=O(\text{\textgreek{e}}_{0})\label{eq:DeformationTensorUniformDecay-1-1}
\end{equation}
 and in the $(u,r,\text{\textgreek{sv}})$ coordinate chart on each
connected component of the region $\mathcal{N}_{af,\mathcal{M}}$:
\begin{align}
\mathcal{L}_{T}g=O(\text{\textgreek{e}}_{0})\Big\{ & O(r^{-1-a})drdu+O(r)d\text{\textgreek{sv}}d\text{\textgreek{sv}}+O(1)dud\text{\textgreek{sv}}+\label{eq:DeformationTensorTAwayDecay-1-1}\\
 & +O(r^{-a})drd\text{\textgreek{sv}}+O(r^{-1})du^{2}+O(r^{-2-a})dr^{2}\Big\}\nonumber 
\end{align}
and 
\begin{align}
\mathcal{L}_{K}g=O(\text{\textgreek{e}}_{0})\Big\{ & O(r^{-1-a})drdu+O(r)d\text{\textgreek{sv}}d\text{\textgreek{sv}}+O(1)dud\text{\textgreek{sv}}+\label{eq:DeformationTensorKAwayDecay-1-1}\\
 & +O(r^{-a})drd\text{\textgreek{sv}}+O(r^{-1})du^{2}+O(r^{-2-a})dr^{2}\Big\}.\nonumber 
\end{align}
Moreover, any further Lie differentiation of $g$ in the direction
of $T$ or $K$ should improve the above decay rates by a factor of
$\bar{t}^{-1}$. 

\smallskip{}

In this case, we also relax the integrated local energy decay assumption
\ref{enu:IntegratedLocalEnergyDecayImproved}, by replacing the $\bar{t}^{-\text{\textgreek{d}}_{0}}$
factor in the last term of the right hand side of (\ref{eq:IntegratedLocalEnergyDecayImprovedDecay-2})
with $\text{\textgreek{e}}_{0}$. Under these weaker assumptions,
we can still obtain (\ref{eq:ImprovedNonDegenerateEnergyDecay}) and
(\ref{eq:ImprovedEnergyDecay}), but with an $O(\text{\textgreek{e}}_{0})$
loss in the exponents of both inequalities. The proof in this case
is similar, and actually easier.
\end{rem*}
The proof of Theorem \ref{thm:ImprovedDecayEnergy} will be presented
in Section \ref{sub:ProofImprovedEnergyDecay}.

As a Corollary of Theorem \ref{thm:ImprovedDecayEnergy}, we will
establish improved pointwise decay rates for $\text{\textgreek{f}}$
and its derivatives:
\begin{cor}
\label{cor:ImprovedPointwiseDecay}Assume that $(\mathcal{M}^{d+1},g)$,
$d\ge3$, satisfies Assumptions \ref{enu:AsymptoticFlatness}-\ref{enu:SobolevInequality},
\ref{enu:Transversality}-\ref{enu:BoundednessVolume} and \ref{enu:IntegratedLocalEnergyDecayImproved}.
Then for any smooth function $\text{\textgreek{f}}$ on $\mathcal{M}$
solving $\square_{g}\text{\textgreek{f}}=F$ on $J^{+}(\{\bar{t}=0\})$
with suitably decaying initial data on $\{\bar{t}=0\}$ and satisfying
boundary conditions on $\partial_{tim}\mathcal{M}$ belonging to the
class $\mathcal{C}_{ILED}^{(T,K)}$ and such that 
\begin{equation}
Re\Big\{\int_{\partial_{tim}\mathcal{M}^{\text{\textgreek{t}}}}Y\text{\textgreek{f}}\cdot\bar{\text{\textgreek{f}}}\, dh_{\partial_{tim}\mathcal{M}^{\text{\textgreek{t}}}}\Big\}\ge0\mbox{ and }Re\Big\{\int_{\partial_{tim}\mathcal{M}^{\text{\textgreek{t}}}}h_{\partial_{tim}\mathcal{M}^{\text{\textgreek{t}}}}\Big(\nabla_{h_{\partial_{tim}\mathcal{M}^{\text{\textgreek{t}}}}}(Y\text{\textgreek{f}}),\nabla_{h_{\partial_{tim}\mathcal{M}^{\text{\textgreek{t}}}}}\bar{\text{\textgreek{f}}}\Big)\, dh_{\partial_{tim}\mathcal{M}^{\text{\textgreek{t}}}}\Big\}\ge0,\label{eq:BoundaryConditionImprovedDecay-1}
\end{equation}
 the following pointwise decay estimates hold for any $0<\text{\textgreek{e}}\ll\text{\textgreek{d}}_{0}$,
any integer $m\ge0$ and any $\text{\textgreek{t}}\ge0$, provided
(\ref{eq:NormFinitenessFriedlanderRadiation}) holds for all $0\le j\le d+1+m+\lceil\text{\textgreek{d}}_{0}^{-1}\cdot2(q-1)\rceil(3q+1)\cdot k$:

1. In case the dimesion $d$ is odd, we can bound: 
\begin{equation}
\sup_{\{\bar{t}=\text{\textgreek{t}}\}}\big|\text{\textgreek{f}}\big|^{2}\lesssim_{\text{\textgreek{e}}}\text{\textgreek{t}}^{-d}\cdot\mathcal{E}_{0,d}[\text{\textgreek{f}}](0)+\mathcal{F}_{deg,\text{\textgreek{e}},0,d}[F](\text{\textgreek{t}}),\label{eq:ImprovedPointwiseDecayAway}
\end{equation}

and if $m\ge1$: 
\begin{equation}
\sup_{\{\bar{t}=\text{\textgreek{t}}\}}\big|\nabla_{g}^{m}\text{\textgreek{f}}\big|_{h}^{2}\lesssim_{m,\text{\textgreek{e}}}\text{\textgreek{t}}^{-d-1}\cdot\mathcal{E}_{m+2,d}[\text{\textgreek{f}}](0)+\mathcal{F}_{deg,\text{\textgreek{e}},m+2,d}[F](\text{\textgreek{t}}).\label{eq:ImprovedDecayOddHigherOrder}
\end{equation}

2. In case the dimension $d$ is even, we can bound: 
\begin{equation}
\sup_{\{\bar{t}=\text{\textgreek{t}}\}}\big|\nabla_{g}^{m}\text{\textgreek{f}}\big|_{h}^{2}\lesssim_{m,\text{\textgreek{e}}}\text{\textgreek{t}}^{-d+C_{m}\text{\textgreek{e}}}\cdot\mathcal{E}_{m,d}[\text{\textgreek{f}}](0)+\mathcal{F}_{\text{\textgreek{e}},m,d}[F](\text{\textgreek{t}}).\label{eq:ImprovedPointwiseDecay}
\end{equation}
In the above, 
\begin{equation}
\mathcal{E}_{m,d}[\text{\textgreek{f}}](0)\doteq\mathcal{E}_{bound}^{(2\lceil\frac{d+1}{2}\rceil,\lceil\frac{d+1}{2}\rceil,m+\lceil\text{\textgreek{d}}_{0}^{-1}\cdot2(\lceil\frac{d+1}{2}\rceil-1)\rceil(3\lceil\frac{d+1}{2}\rceil+1)\cdot k)}[\text{\textgreek{f}}](0),
\end{equation}
\begin{equation}
\mathcal{F}_{deg,\text{\textgreek{h}},m,d}[F](\text{\textgreek{t}})\doteq\mathcal{F}_{deg,\text{\textgreek{h}}}^{(2\lceil\frac{d+1}{2}\rceil,\lceil\frac{d+1}{2}\rceil,m+\lceil\text{\textgreek{d}}_{0}^{-1}\cdot2(\lceil\frac{d+1}{2}\rceil-1)\rceil(3q+1)\cdot k,k)}[F](\text{\textgreek{t}})
\end{equation}
and 
\begin{equation}
\mathcal{F}_{\text{\textgreek{e}},m,d}[F(\text{\textgreek{t}})\doteq\mathcal{F}_{\text{\textgreek{e}},\text{\textgreek{e}}}^{(2\lceil\frac{d+1}{2}\rceil,\lceil\frac{d+1}{2}\rceil,m+\lceil\text{\textgreek{d}}_{0}^{-1}\cdot2(\lceil\frac{d+1}{2}\rceil-1)\rceil(3q+1)\cdot k,k)}[F](\text{\textgreek{t}}).
\end{equation}
See Section \ref{sub:ShorthandEnergyNorms} for the notations on the
energy norms used.\end{cor}
\begin{rem*}
Again, Corollary \ref{cor:ImprovedPointwiseDecay} still holds if
one replaces the condition (\ref{eq:BoundaryConditionImprovedDecay-1})
with any boundary condition for which Lemma \ref{lem:ILEDCommuted}
can still be established. Moreover, in case the integrated local energy
decay statement in Assumption \ref{enu:IntegratedLocalEnergyDecayImproved}
does not lose derivatives (i.\,e.~$k=0$), Assumption \ref{enu:BoundsDeformationTensor}
can be replaced by (\ref{eq:DeformationTensorUniformDecay-1-1})--(\ref{eq:DeformationTensorKAwayDecay-1-1})
and the integrated local energy decay assumption \ref{enu:IntegratedLocalEnergyDecayImproved}
can be relaxed by replacing the $\bar{t}^{-\text{\textgreek{d}}_{0}}$
factor in the last term of the right hand side of (\ref{eq:IntegratedLocalEnergyDecayImprovedDecay-2})
with $\text{\textgreek{e}}_{0}$. In this case, inequalities (\ref{eq:ImprovedPointwiseDecayAway}),
(\ref{eq:ImprovedDecayOddHigherOrder}) and (\ref{eq:ImprovedPointwiseDecay})
still hold with an $O(\text{\textgreek{e}}_{0})$ loss in the exponent
of $\text{\textgreek{t}}$.
\end{rem*}
The proof of Corollary \ref{cor:ImprovedPointwiseDecay} will be presented
in Section \ref{sub:ProofCorollaryImprovedDecay}.

\subsection{\label{sub:SketchOfProof}Sketch of the proof of Theorem \ref{thm:ImprovedDecayEnergy}
and Corollary \ref{cor:ImprovedPointwiseDecay}}

In this Section, we will first sketch the proof of Theorem \ref{thm:ImprovedDecayEnergy}
and Corollary \ref{cor:ImprovedPointwiseDecay} under some simplifying
assumptions on the structure of the spacetime $(\mathcal{M},g)$,
and then we will highlight the main difficulties arising in the general
case. 

Let us assume first that $d=3$, $\partial\mathcal{M}=\emptyset$,
$F=0$, $k=0$ (i.\,e.~there is no derivative loss in the integrated
local energy decay statement (\ref{eq:IntegratedLocalEnergyDecayImprovedDecay})),
$m=1$ and the vector field $T$ is globally timelike and Killing.
In this case, there is no condition on $\partial_{tim}\mathcal{M}$
that $\text{\textgreek{f}}$ is assumed to satisfy, and we can assume
without loss of generality that the vector field $K$ has been fixed
so that $T\equiv K$. Let us also note that in this case, the $\text{\textgreek{e}}$-loss
in (\ref{eq:ImprovedNonDegenerateEnergyDecay}) can be dropped, and
the estimates of Theorem \ref{thm:ImprovedDecayEnergy} and Corollary
\ref{cor:ImprovedPointwiseDecay} read as follows: 
\begin{equation}
\sum_{i_{1}+i_{2}=1}\int_{\{\bar{t}=\text{\textgreek{t}}\}}\big(\big|\nabla_{h_{\text{\textgreek{t}},T}}^{i_{1}+1}(T^{i_{2}}\text{\textgreek{f}})\big|_{h_{T}}^{2}+r^{-2}\big|T^{2}\text{\textgreek{f}}\big|^{2}\big)\, dh_{T}\lesssim\text{\textgreek{t}}^{-4}\mathcal{E}_{in}[\text{\textgreek{f}}](0),\label{eq:EnergyDecaySketch}
\end{equation}
\begin{equation}
\sup_{\{\bar{t}=\text{\textgreek{t}}\}}|\text{\textgreek{f}}|\lesssim\text{\textgreek{t}}^{-\frac{3}{2}}\label{eq:PointwiseDecaySketch}
\end{equation}
and: 
\begin{equation}
\sup_{\{\bar{t}=\text{\textgreek{t}}\}}|\nabla_{g}\text{\textgreek{f}}|_{h}\lesssim\text{\textgreek{t}}^{-2}.\label{eq:HigherPointwiseDecaySketch}
\end{equation}

The main idea for the proof of (\ref{eq:EnergyDecaySketch}) is the
following (assuming without loss of generality that $\text{\textgreek{f}}$
is real valued): From Theorem \ref{thm:SlowPointwiseDecayHighDerivativesNewMethod}
we deduce that:
\begin{equation}
\sum_{i_{1}+i_{2}=1}\int_{\{\bar{t}=\text{\textgreek{t}}\}}\big(r^{2}\big|\nabla_{h_{\text{\textgreek{t}},T}}^{2}\text{\textgreek{f}}\big|_{h_{T}}^{2}+\big|\nabla_{h_{\text{\textgreek{t}},T}}^{i_{1}}(T^{1+i_{2}}\text{\textgreek{f}})\big|_{h_{T}}^{2}\big)\, dh_{T}\lesssim\text{\textgreek{t}}^{-2}\mathcal{E}_{in}[\text{\textgreek{f}}](0).\label{eq:EnergyDecaySketch-1}
\end{equation}
Let us fix a vector field $L$ on $\mathcal{M}$ so that $[T,L]=0$
and $L=\partial_{v}$ in the $(\bar{t},v,\text{\textgreek{sv}})$
coordinate system on each connected component of the region $(r\gg1)$.
Using the expression for the equation $\square\text{\textgreek{f}}=0$,
from (\ref{eq:EnergyDecaySketch-1}) we deduce that 
\begin{equation}
\sum_{i_{1}+i_{2}=1}\int_{\{\bar{t}=\text{\textgreek{t}}\}}\big(r^{2}\big|L(T\text{\textgreek{f}})\big|^{2}+\big|\nabla_{h_{\text{\textgreek{t}},T}}^{i_{1}}(T^{1+i_{2}}\text{\textgreek{f}})\big|_{h_{T}}^{2}\big)\, dh_{T}\lesssim\text{\textgreek{t}}^{-2}\mathcal{E}_{in}[\text{\textgreek{f}}](0).\label{eq:IntermediateStepSketch}
\end{equation}

Fixing a dyadic sequence $\{\text{\textgreek{t}}_{n}\}_{n\in\mathbb{N}}$,
by repeating the proof Theorem \ref{thm:FirstPointwiseDecayNewMethod}
on the intervals $\{\text{\textgreek{t}}_{n}\le\bar{t}\le\text{\textgreek{t}}_{n+1}\}$
with $T\text{\textgreek{f}}$ in place of $\text{\textgreek{f}}$
and using the estimate (\ref{eq:IntermediateStepSketch}) (notice
that for the sketch of the proof we have assumed that $T$ is Killing,
and thus $[T,\square_{g}]=0$), we readily obtain: 
\begin{equation}
\int_{\{\bar{t}=\text{\textgreek{t}}\}}\big(\big|\nabla_{h_{\text{\textgreek{t}},T}}(T\text{\textgreek{f}})\big|_{h_{T}}^{2}+r^{-2}\big|T^{2}\text{\textgreek{f}}\big|^{2}\big)\, dh_{T}\lesssim\text{\textgreek{t}}^{-4}\mathcal{E}_{in}[\text{\textgreek{f}}](0).\label{eq:IntermediateStepSketch-1}
\end{equation}
Using again the expression for the wave equation $\square\text{\textgreek{f}}=0$,
from (\ref{eq:IntermediateStepSketch-1}) we deduce that 
\begin{equation}
\int_{\{\bar{t}=\text{\textgreek{t}}\}}(\text{\textgreek{D}}_{ell}\text{\textgreek{f}})^{2}\, dh_{T}\lesssim\text{\textgreek{t}}^{-4}\mathcal{E}_{in}[\text{\textgreek{f}}](0)\label{eq:EllipticEquation}
\end{equation}
for a suitable elliptic operator $\text{\textgreek{D}}_{ell}$ on
the hyperboloids $\{\bar{t}=\text{\textgreek{t}}\}$. The elliptic
estimates of Section \ref{sec:Elliptic-estimates} of the Appendix
then yield: 
\begin{equation}
\int_{\{\bar{t}=\text{\textgreek{t}}\}}\big|\nabla_{h_{\text{\textgreek{t}},T}}^{2}\text{\textgreek{f}}\big|_{h_{T}}\, dh_{T}\lesssim\text{\textgreek{t}}^{-4}\mathcal{E}_{in}[\text{\textgreek{f}}](0),
\end{equation}
 which combined with (\ref{eq:IntermediateStepSketch-1}) yields (\ref{eq:EnergyDecaySketch}).

The estimate (\ref{eq:HigherPointwiseDecaySketch}) for $\nabla_{g}\text{\textgreek{f}}$
follows readily from a Sobolev inequality applied on (\ref{eq:EnergyDecaySketch})
for $\text{\textgreek{f}}$ and $T\text{\textgreek{f}}$ (combined
with elliptic estimates). The zeroth order estimate (\ref{eq:PointwiseDecaySketch}),
on the other hand, follows from (\ref{eq:EnergyDecaySketch}), the
decay estimate from Theorem \ref{thm:FirstPointwiseDecayNewMethod}:
\[
\int_{\{\bar{t}=\text{\textgreek{t}}\}}\big(\big|\nabla_{h_{\text{\textgreek{t}},T}}\text{\textgreek{f}}\big|_{h_{T}}^{2}+r^{-2}\big|T\text{\textgreek{f}}\big|^{2}\big)\, dh_{T}\lesssim\text{\textgreek{t}}^{-2}\mathcal{E}_{in}[\text{\textgreek{f}}](0)
\]
and the following Gagliardo--Nirenberg type estimate on $\{\bar{t}=\text{\textgreek{t}}\}$
(see Section \ref{sub:GagliardoNirenberg}): 
\begin{align}
\sup_{\{\bar{t}=\text{\textgreek{t}}\}}|\text{\textgreek{f}}|^{2}\lesssim & \Big(\int_{\{\bar{t}=\text{\textgreek{t}}\}}\big(\big|\nabla_{h_{\text{\textgreek{t}},T}}\text{\textgreek{f}}\big|_{h_{T}}^{2}+r^{-2}\big|T\text{\textgreek{f}}\big|^{2}\big)\, dh_{T}\Big)^{1/2}\cdot\Big(\sum_{i_{1}+i_{2}=1}\int_{\{\bar{t}=\text{\textgreek{t}}\}}\big(\big|\nabla_{h_{\text{\textgreek{t}},T}}^{i_{1}+1}(T^{i_{2}}\text{\textgreek{f}})\big|_{h_{T}}^{2}+r^{-2}\big|T^{2}\text{\textgreek{f}}\big|^{2}\big)\, dh_{T}\Big)^{1/2}+\\
 & +\sum_{i_{1}+i_{2}=1}\int_{\{\bar{t}=\text{\textgreek{t}}\}}\big(\big|\nabla_{h_{\text{\textgreek{t}},T}}^{i_{1}+1}(T^{i_{2}}\text{\textgreek{f}})\big|_{h_{T}}^{2}+r^{-2}\big|T^{2}\text{\textgreek{f}}\big|^{2}\big)\, dh_{T}.\nonumber 
\end{align}

One important difficulty arising in the proof of Theorem \ref{thm:ImprovedDecayEnergy}
in the more general class of spacetimes $(\mathcal{M},g)$ under consideration
comes from the fact that $T$ is not in general a Killing vector field,
and in fact its deformation tensor decays only like $\bar{t}^{-\text{\textgreek{d}}_{0}}$
for some small $\text{\textgreek{d}}_{0}>0$. This results in a number
of error terms appearing each time $\square$ is commuted with $T$,
which can only be controlled in the final step of the estimates, using
also some refined elliptic estimates leading to the $\text{\textgreek{e}}$-loss
in (\ref{eq:ImprovedNonDegenerateEnergyDecay}) (however, we avoid
this loss in (\ref{eq:ImprovedEnergyDecay})). Furthermore, the slow
$O(\bar{t}^{-\text{\textgreek{d}}_{0}})$ decay of the deformation
tensor of $T$ combined with the loss of derivatives in the integrated
local energy decay statement (\ref{eq:IntegratedLocalEnergyDecayImprovedDecay})
require an iteration of the above procedure $\sim\text{\textgreek{d}}_{0}^{-1}\cdot k$
times, leading to the corresponding derivative loss in the estimates
of Theorem \ref{thm:ImprovedDecayEnergy}. Notice, however, that in
case where the integrated local energy decay statement (\ref{eq:IntegratedLocalEnergyDecayImprovedDecay})
does not lose derivatives, the same steps can be applied (without
the need of the extra iterations of the procedure) even when the deformation
tensor of $T$ does not decay at all, but is, instead, merely bounded
by some small constant $\text{\textgreek{e}}_{0}>0$. 

We should also notice that in the general case where $T$ is not everywhere
non-spacelike (and thus we necessarily have $K\neq T$), one extra
difficulty arises from the fact that, in order to avoid losing unnecessary
$r$-weights in the estimates of Theorem \ref{thm:ImprovedDecayEnergy},
we only commute $\square$ with $\text{\textgreek{q}}K$ instead of
$K$, for a cut-off function $\text{\textgreek{q}}$ supported in
the far away region $\{r\gg1\}$ (notice that $g(K,K)\sim r^{2}$
in case $K\neq T$ in the far away region). However, commutation with
$\text{\textgreek{q}}K$ leads to additional error terms which do
not decay in time. A key element in dealing with these terms are the
elliptic estimates of Lemma \ref{Prop:ControlOfTheAngularTermsAway}.

\subsection{Commutation with $T$, $K_{R_{c}}$ and control of the error terms}

The following Lemma will provide us with some useful estimates for
the commutator of $\square_{g}$ with the almost Killing vector field
$T$ and the vector field $K_{R_{c}}$ which fails to be almost Killing
in the region $r\sim R_{c}$ of $\mathcal{M}$.
\begin{lem}
\label{lem:CommutationsWithTandK}Provided that $1\ll R_{0}\ll R_{c}$,
let us fix $\text{\textgreek{q}}_{R_{0}},\text{\textgreek{q}}_{\sim R_{c}}:\mathcal{M}\rightarrow[0,1]$
so that 
\begin{equation}
\text{\textgreek{q}}_{R_{0}}=\text{\textgreek{q}}\circ(\frac{r}{R_{0}})
\end{equation}
 for some smooth $\text{\textgreek{q}}:[0,+\infty)\rightarrow[0,1]$
satisfying $\text{\textgreek{q}}\equiv1$ on $[0,1]$ and $\text{\textgreek{q}}\equiv0$
on $[2,+\infty)$ and
\begin{equation}
\text{\textgreek{q}}_{\sim R_{c}}=\begin{cases}
0, & \mbox{ on }\{r\notin[R_{c},2R_{c}]\}\\
1, & \mbox{ on }\{r\in[R_{c},2R_{c}]\}.
\end{cases}
\end{equation}
 Then the following commutation relations hold for any $\text{\textgreek{f}}\in C^{\infty}(\mathcal{M})$
and any integer $l\ge0$: 
\begin{align}
\square_{g}(T^{l}\text{\textgreek{f}})= & T^{l}(\square_{g}\text{\textgreek{f}})+\text{\textgreek{q}}_{R_{0}}\cdot\sum_{j=1}^{l}O(\text{\textgreek{t}}^{-(l-j)-\text{\textgreek{d}}_{0}})\big(|\nabla^{j+1}\text{\textgreek{f}}|_{h}+|\nabla^{j}\text{\textgreek{f}}|_{h}\big)+\label{eq:Tcommuted}\\
 & +(1-\text{\textgreek{q}}_{R_{0}})\mathscr{Err}_{(T,\ldots T)}^{(l)}[\text{\textgreek{f}}]\nonumber 
\end{align}
 and for any $(e_{1},\ldots e_{l})\in\{0,1\}^{l}$ setting $X^{(0)}=T$
and $X^{(1)}=K_{R_{c}}$: 
\begin{align}
\square_{g}(X^{(e_{1})}\cdots X^{(e_{l})}\text{\textgreek{f}})= & X^{(e_{1})}\cdots X^{(e_{l})}(\square_{g}\text{\textgreek{f}})+\text{\textgreek{q}}_{R_{0}}\cdot\sum_{j=1}^{l}O(\text{\textgreek{t}}^{-(l-j)-\text{\textgreek{d}}_{0}})\big(|\nabla^{j+1}\text{\textgreek{f}}|_{h}+|\nabla^{j}\text{\textgreek{f}}|_{h}\big)+\label{eq:Kcommuted}\\
 & +(1-\text{\textgreek{q}}_{R_{0}})(1-\text{\textgreek{q}}_{r\sim R_{c}})\mathscr{Err}_{(X^{(e_{1})},\ldots X^{(e_{l})})}^{(l)}[\text{\textgreek{f}}]+\nonumber \\
 & +\text{\textgreek{q}}_{r\sim R_{c}}\sum_{j=0}^{l}O(r^{j-1})\big|\nabla^{j+1}\text{\textgreek{f}}\big|_{h},\nonumber 
\end{align}
where (setting for simplicity $\text{\textgreek{f}}_{i_{1}\ldots i_{j}}=X^{(e_{i_{1}})}\cdots X^{(e_{i_{j}})}\text{\textgreek{f}}$)
\begin{align}
\mathscr{Err}_{(X^{(e_{1})},\ldots X^{(e_{l})})}^{(l)}[\text{\textgreek{f}}]=\sum_{\substack{components\\
of\,\mathcal{N}_{af,\mathcal{M}}
}
}\sum_{j=0}^{l-1}\sum_{\{i_{1},\ldots i_{j}\}\subset\{1,\ldots l\}}O(\text{\textgreek{t}}^{-(l-j-1)-\text{\textgreek{d}}_{0}})\cdot\Big\{ & O(r^{-1})\partial_{v}^{2}(\text{\textgreek{f}}_{i_{1}\ldots i_{j}})+O(r^{-2})\partial_{v}\partial_{\text{\textgreek{sv}}}(\text{\textgreek{f}}_{i_{1}\ldots i_{j}})+\label{eq:ErrorCommuteDt-1}\\
 & +O(r^{-3})\partial_{\text{\textgreek{sv}}}\partial_{\text{\textgreek{sv}}}(\text{\textgreek{f}}_{i_{1}\ldots i_{j}})+O(r^{-1-a})\partial_{u}\partial_{v}(\text{\textgreek{f}}_{i_{1}\ldots i_{j}})+\nonumber \\
 & +O(r^{-2-a})\partial_{u}\partial_{\text{\textgreek{sv}}}(\text{\textgreek{f}}_{i_{1}\ldots i_{j}})+O(r^{-2-a})\partial_{u}^{2}(\text{\textgreek{f}}_{i_{1}\ldots i_{j}})+\nonumber \\
 & +O(r^{-1-a})\partial_{v}(\text{\textgreek{f}}_{i_{1}\ldots i_{j}})+O(r^{-2-a})\partial_{\text{\textgreek{sv}}}(\text{\textgreek{f}}_{i_{1}\ldots i_{j}})+O(r^{-2-a})\partial_{u}(\text{\textgreek{f}}_{i_{1}\ldots i_{j}})\Big\}.\nonumber 
\end{align}
\end{lem}
\begin{proof}
The relations (\ref{eq:Tcommuted}) and (\ref{eq:Kcommuted}) follow
readily by differnetiating the expression for the wave equation 
\begin{equation}
\square\text{\textgreek{f}}=g^{\text{\textgreek{m}\textgreek{n}}}\partial_{\text{\textgreek{m}}}\partial_{\text{\textgreek{n}}}\text{\textgreek{f}}+\frac{1}{\sqrt{-\det(g)}}\partial_{\text{\textgreek{m}}}\big(\sqrt{-\det(g)}\cdot g^{\text{\textgreek{m}\textgreek{n}}}\big)\partial_{\text{\textgreek{n}}}\text{\textgreek{f}}
\end{equation}
 with respect to $T$ and $K_{R_{c}}$, using Assumption \ref{enu:BoundsDeformationTensor}
on the almost Killing vector fields $T,K$ and the relation $K_{R_{c}}=\text{\textgreek{q}}_{R_{c}}K+(1-\text{\textgreek{q}}_{R_{c}})T$.
\end{proof}
The error term obtained from the commutation of $\square$ with $K_{R_{c}}$
on the region $r\sim R_{c}$ (where $K_{R_{c}}$ fails to be almost
Killing) will be controlled with the use of suitable elliptic estimates.
In particular, we will establish the following Lemma:
\begin{lem}
\label{lem:EllipticEstimatesWithoutKillingError}For any $l\in\mathbb{N}$
with $l\le\lfloor\frac{d+1}{2}\rfloor$, any $n_{0}\in\mathbb{N}$,
any $\text{\textgreek{b}}\in(-\bar{\text{\textgreek{d}}}_{n_{0}},1)$
(for some $\bar{\text{\textgreek{d}}}_{n_{0}}>0$ depending on $n_{0}$)
and any $0<\text{\textgreek{e}}\ll1-\text{\textgreek{b}}$, if $R_{c}$
is large in terms of $\text{\textgreek{b}},\text{\textgreek{e}}$
and the geometry of $(\mathcal{M},g)$ we can bound for any $\text{\textgreek{f}}\in C^{\infty}(\mathcal{M})$
satisfying for all $j_{1}+j_{2}\le l+n_{0}$ the finite radiation
field condition $\limsup_{r\rightarrow+\infty}\big|r^{\frac{d-1}{2}+j_{1}}\nabla_{h_{\text{\textgreek{t}},N}}^{j_{1}}(N^{j_{1}}\text{\textgreek{f}})\big|_{h_{\text{\textgreek{t}},N}}<+\infty$:
\begin{equation}
\begin{split}\sum_{n=0}^{n_{0}}\sum_{j=0}^{l-3}\int_{\{\bar{t}=\text{\textgreek{t}}\}}r_{+}^{-\text{\textgreek{b}}-2j}| & \nabla_{h_{\text{\textgreek{t}},N}}^{n+l-j}\text{\textgreek{f}}|_{\big(1-\log(r_{tim})\big)\cdot h_{R_{c}}}^{2}\, dh_{N}+\sum_{n=0}^{n_{0}}\sum_{j=1}^{l}\int_{\{\bar{t}=\text{\textgreek{t}}\}\cap\{R_{c}\le r\le2R_{c}\}}r_{+}^{2(j-1)-\text{\textgreek{b}}+2\text{\textgreek{e}}}\big|\nabla_{g}^{j+n}\text{\textgreek{f}}\big|_{h}^{2}\, dh_{N}\le\\
\le & C_{\text{\textgreek{b}},n_{0}}(R_{c})\int_{\{\bar{t}=\text{\textgreek{t}}\}}r_{+}^{-\text{\textgreek{b}}}\Big\{\sum_{n=0}^{n_{0}}\big|\nabla_{h_{\text{\textgreek{t}},N}}^{n+l-2}(\square\text{\textgreek{f}})\big|_{\big(1-\log(r_{tim})\big)\cdot h_{\text{\textgreek{t}},N}}^{2}+\mathcal{T}_{T,K,R_{c}}^{(l,n_{0})}[\text{\textgreek{f}}]\Big\}\, dh_{N}+\\
 & +C_{\text{\textgreek{b}},n_{0}}(R_{c})\sum_{j=0}^{1}\max\Big\{-Re\big\{\int_{\partial_{tim}\mathcal{M}^{\text{\textgreek{t}}}}h_{\partial_{tim}\mathcal{M}^{\text{\textgreek{t}}}}\Big(\nabla_{h_{\partial_{tim}\mathcal{M}^{\text{\textgreek{t}}}}}^{j}(Y\text{\textgreek{f}}),\nabla_{h_{\partial_{tim}\mathcal{M}^{\text{\textgreek{t}}}}}^{j}\bar{\text{\textgreek{f}}}\Big)\, dh_{\partial_{tim}\mathcal{M}^{\text{\textgreek{t}}}}\big\},0\Big\},
\end{split}
\label{eq:EllipticEstimatesWithoutKillingError}
\end{equation}
where $\mathcal{T}_{T,K,R_{c}}^{(l,k)}[\text{\textgreek{f}}]$ is
given by (\ref{eq:EllipticEstimatesEnergy}).\end{lem}
\begin{rem*}
From now on, we will assume that $R_{c}$ has been fixed large in
terms of $\text{\textgreek{b}}$ and the geometry of $(\mathcal{M},g)$.
Hence, we will drop the dependence of constants on $R_{c}$ (replacing
it with dependence on the parameters on which $\text{\textgreek{b}}$
will depend on).\end{rem*}
\begin{proof}
Inequality (\ref{eq:EllipticEstimatesWithoutKillingError}) follows
readily after adding (\ref{eq:EllipticEstimatesDegenerate}) and (\ref{eq:EstimateFarAwayForAngularBulk-1-1})
(for $4\text{\textgreek{e}}$ in place of $\text{\textgreek{e}}$),
and absorbing the last term of the right hand side of (\ref{eq:EstimateFarAwayForAngularBulk-1-1})
into the left hand side of (\ref{eq:EllipticEstimatesDegenerate}),
which can be done provided $R_{c}$ has been fixed large in terms
of $\text{\textgreek{b}},R_{0}$ and $0<\text{\textgreek{e}}\ll1-\text{\textgreek{b}}$.
Recall that in view of Assumption \ref{enu:UniformityOfDegenerateEllipticEstimates},
the constants in the right hand sides of (\ref{eq:EllipticEstimatesDegenerate})
and (\ref{eq:EstimateFarAwayForAngularBulk-1-1}) do not depend on
$\text{\textgreek{t}}$. 
\end{proof}

\subsection{Integrated local energy decay after commuting with $T,K_{R_{c}}$}

Let us fix a vector field $L$ on $\mathcal{M}$ such that $[L,T]=0$
everywhere, $L\equiv0$ in the region $\{r\le2R_{c}\}$ and $L\equiv\partial_{v}$
in the coordinate chart $(\bar{t},v,\text{\textgreek{sv}})$ on each
connected component of the region $\{r\ge2R_{c}+1\}$.
\begin{lem}
\label{lem:ILEDCommuted}For any $l\in\mathbb{N}$ with $l\le\lfloor\frac{d+1}{2}\rfloor$,
and any $0<\text{\textgreek{e}}<1$, if $R_{c}$ is large in terms
of $\text{\textgreek{e}}$ and the geometry of $(\mathcal{M},g)$
we can bound for any integer $m\ge1$, any $0\le\text{\textgreek{t}}_{1}\le\text{\textgreek{t}}_{2}$
and any $\text{\textgreek{f}}\in C^{\infty}(\mathcal{M})$ solving
$\square\text{\textgreek{f}}=F$ satisfying boundary conditions on
$\partial_{tim}\mathcal{M}$ belonging to the class $\mathcal{C}_{ILED}^{(T,K)}$
and such that 
\begin{equation}
Re\Big\{\int_{\partial_{tim}\mathcal{M}^{\text{\textgreek{t}}}}Y\text{\textgreek{f}}\cdot\bar{\text{\textgreek{f}}}\, dh_{\partial_{tim}\mathcal{M}^{\text{\textgreek{t}}}}\Big\}\ge0\mbox{ and }Re\Big\{\int_{\partial_{tim}\mathcal{M}^{\text{\textgreek{t}}}}h_{\partial_{tim}\mathcal{M}^{\text{\textgreek{t}}}}\Big(\nabla_{h_{\partial_{tim}\mathcal{M}^{\text{\textgreek{t}}}}}(Y\text{\textgreek{f}}),\nabla_{h_{\partial_{tim}\mathcal{M}^{\text{\textgreek{t}}}}}\bar{\text{\textgreek{f}}}\Big)\, dh_{\partial_{tim}\mathcal{M}^{\text{\textgreek{t}}}}\Big\}\ge0\label{eq:BoundaryConditionImprovedDecay-1-1}
\end{equation}
(provided (\ref{eq:NormFinitenessFriedlanderRadiation}) holds for
all $0\le j\le m+l+k+\lceil\frac{d}{2}\rceil$): 
\begin{equation}
\begin{split}\sum_{j=0}^{m-1}\int_{\{\text{\textgreek{t}}_{1}\le\bar{t}\le\text{\textgreek{t}}_{2}\}}\Big(\, & \sum_{j_{2}+j_{3}=l+j-1}r_{+}^{-1+\text{\textgreek{e}}}\big|\nabla_{h_{\text{\textgreek{t}},N}}^{j_{2}+1}(N^{j_{3}}\text{\textgreek{f}})\big|_{h_{\text{\textgreek{t}},N}}^{2}+\sum_{1\le j_{2}+j_{3}\le l+j-1}r_{+}^{-1+\text{\textgreek{e}}-2(l+j-j_{2}-j_{3})}\big|\nabla_{h_{\text{\textgreek{t}},N}}^{j_{2}}(N^{j_{3}}\text{\textgreek{f}})\big|_{h_{\text{\textgreek{t}},N}}^{2}+r_{+}^{-1-\text{\textgreek{e}}}\big|N^{l+j}\text{\textgreek{f}}\big|^{2}\Big)\, dg+\\
+\sum_{j=0}^{m-1}\sum_{1\le j_{2}+j_{3}\le l+j-1} & \int_{\{\bar{t}=\text{\textgreek{t}}_{2}\}}r_{+}^{-2(l+j-1-j_{2}-j_{3})}\Big(r_{+}^{\text{\textgreek{e}}}|\mathcal{L}_{L}\nabla_{h_{\text{\textgreek{t}},N}}^{j_{2}}(N^{j_{3}}\text{\textgreek{f}})|^{2}+\big|\nabla_{h_{\text{\textgreek{t}},N}}^{j_{2}+1}(N^{j_{3}}\text{\textgreek{f}})\big|_{h_{\text{\textgreek{t}},N}}^{2}+r_{+}^{-2}\big|N^{l+j}\text{\textgreek{f}}\big|^{2}\Big)\, dh_{N}+\\
+\sum_{j=0}^{m-1}\sum_{1\le j_{2}+j_{3}\le l+j-1} & \int_{\mathcal{H}^{+}(\text{\textgreek{t}}_{1},\text{\textgreek{t}}_{2})}\big|\nabla_{h_{\mathcal{H}}}^{j_{2}+1}(K_{R_{c}}^{j_{3}}\text{\textgreek{f}})\big|_{h_{\mathcal{H}}}^{2}\, dh_{\mathcal{H}}\le\\
\le C_{\text{\textgreek{e}},m}\Big\{\, & \sum_{j=0}^{k+m-1}\sum_{1\le j_{2}+j_{3}\le l+j-1}\int_{\{\bar{t}=\text{\textgreek{t}}_{1}\}}r_{+}^{\text{\textgreek{e}}-2(l+j-1-j_{2}-j_{3})}\Big(\big|\nabla_{h_{\text{\textgreek{t}},N}}^{j_{2}+1}(N^{j_{3}}\text{\textgreek{f}})\big|_{h_{\text{\textgreek{t}},N}}^{2}+r_{+}^{-2}\big|N^{l+j}\text{\textgreek{f}}\big|^{2}\Big)\, dh_{N}+\\
 & +\text{\textgreek{t}}_{1}^{-\text{\textgreek{d}}_{0}}\int_{\{\text{\textgreek{t}}_{1}\le\bar{t}\le\text{\textgreek{t}}_{2}\}}r_{+}^{-1+\text{\textgreek{e}}}\mathcal{T}_{T,K,R_{c},sl}^{(l,m+k-1)}[\text{\textgreek{f}}]\, dg+\sum_{j=l-1}^{m+k+l-2}\int_{\{\text{\textgreek{t}}_{1}\le\bar{t}\le\text{\textgreek{t}}_{2}\}}r_{+}^{1+\text{\textgreek{e}}}\big|\nabla_{g}^{j}F\big|_{h}^{2}\, dg\Big\}
\end{split}
\label{eq:ILEDCommuted}
\end{equation}

(where $\mathcal{L}_{L}$ denotes the Lie derivative in the direction
of the outgoing vector field $L$) and 
\begin{equation}
\begin{split}\sum_{j=0}^{m-1}\sum_{1\le j_{2}+j_{3}\le l+j-1}\, & \int_{\{\bar{t}=\text{\textgreek{t}}_{2}\}}r_{+}^{-2(l+j-1-j_{2}-j_{3})}\Big(\big|\nabla_{h_{\text{\textgreek{t}},N}}^{j_{2}+1}(K_{R_{c}}^{j_{3}}\text{\textgreek{f}})\big|_{\big(1-\log(r_{tim})\big)\cdot h_{R_{c}}}^{2}+r_{+}^{-2}\big|K_{R_{c}}^{l+j}\text{\textgreek{f}}\big|^{2}\Big)\, dh_{N}\le\\
\le C_{\text{\textgreek{e}},m}\Big\{\, & \sum_{j=0}^{k+m-1}\sum_{1\le j_{2}+j_{3}\le l+j-1}\int_{\{\bar{t}=\text{\textgreek{t}}_{1}\}}r_{+}^{-2(l+j-1-j_{2}-j_{3})}\Big(\big|\nabla_{h_{\text{\textgreek{t}},N}}^{j_{2}+1}(N^{j_{3}}\text{\textgreek{f}})\big|_{h_{\text{\textgreek{t}},N}}^{2}+r_{+}^{-2}\big|N^{l+j}\text{\textgreek{f}}\big|^{2}\Big)\, dh_{N}+\\
 & +\text{\textgreek{t}}_{1}^{-\text{\textgreek{d}}_{0}}\sum_{j=0}^{m-1}\sum_{1\le j_{2}+j_{3}\le l+j-1}\int_{\{\bar{t}=\text{\textgreek{t}}_{1}\}}r_{+}^{\text{\textgreek{e}}-2(l+j-1-j_{2}-j_{3})}\Big(\big|\nabla_{h_{\text{\textgreek{t}},N}}^{j_{2}+1}(N^{j_{3}}\text{\textgreek{f}})\big|_{h_{\text{\textgreek{t}},N}}^{2}+r_{+}^{-2}\big|N^{l+j}\text{\textgreek{f}}\big|^{2}\Big)\, dh_{N}+\\
 & +\text{\textgreek{t}}_{1}^{-\text{\textgreek{d}}_{0}}\int_{\{\text{\textgreek{t}}_{1}\le\bar{t}\le\text{\textgreek{t}}_{2}\}}r_{+}^{-1+\text{\textgreek{e}}}\mathcal{T}_{T,K,R_{c},sl}^{(l,m+k-1)}[\text{\textgreek{f}}]\, dg+\sum_{j=l-1}^{m+k+l-2}\int_{\{\text{\textgreek{t}}_{1}\le\bar{t}\le\text{\textgreek{t}}_{2}\}}r_{+}^{1+\text{\textgreek{e}}}\big|\nabla_{g}^{j}F\big|_{h}^{2}\, dg\Big\}.
\end{split}
\label{eq:DegenerateBoundednessWithLoss}
\end{equation}
In the above, 
\begin{align}
\mathcal{T}_{T,K,R_{c},sl}^{(l,m+k-1)}[\text{\textgreek{f}}]\doteq & \mathcal{T}_{T,K,R_{c}}^{(l,m+k-1)}[\text{\textgreek{f}}]+\sum_{0\le i_{1}+i_{2}\le m}\sum_{0\le j\le l-2}\int_{\{\text{\textgreek{t}}_{1}\le\bar{t}\le\text{\textgreek{t}}_{2}\}}\big(r_{+}^{-2-2(l-2-j)}+r_{+}^{-2}\text{\textgreek{t}}^{-2(l-2-j)}\big)\big|\nabla_{h_{\text{\textgreek{t}},N}}^{j+i_{1}+1}(K_{R_{c}}^{i_{2}}\text{\textgreek{f}})\big|_{h_{\text{\textgreek{t}},N}}^{2}\, dg+\\
 & +\sum_{0\le i_{1}+i_{2}\le m}\sum_{0\le j_{1}+j_{2}\le l-2}\int_{\{\text{\textgreek{t}}_{1}\le\bar{t}\le\text{\textgreek{t}}_{2}\}}r_{+}^{-2}\text{\textgreek{t}}^{-2(l-j_{1}-j_{2}-2)}\Big(|\nabla_{h_{\text{\textgreek{t}},N}}^{j_{1}+i_{1}}\text{\textgreek{f}}|_{T,K,R_{c}}^{(j_{2}+i_{2}+1)}\Big)^{2}\, dg.\nonumber 
\end{align}
 \end{lem}
\begin{rem*}
Notice that in comparison to (\ref{eq:ILEDCommuted}), the left hand
side of (\ref{eq:DegenerateBoundednessWithLoss}) controls only a
degenerate energy norm of $\text{\textgreek{f}}$ on $\{\bar{t}=\text{\textgreek{t}}_{2}\}$
and no bulk terms. However, in the right hand side of (\ref{eq:DegenerateBoundednessWithLoss})
the dependence on the $r^{\text{\textgreek{e}}}$-weighted initial
energy comes with a factor decaying polynomially in time.\end{rem*}
\begin{proof}
Without loss of generality, we will assume that $\text{\textgreek{t}}_{1}$
is large in terms of $\text{\textgreek{e}},m$. Moreover, we will
assume that $\text{\textgreek{f}}$ is real valued. 

Recall that according to Assumption \ref{enu:DeformationTensor},
in each connected component of the region $\{r\gg1\}$ we have: 
\begin{align}
\mathcal{L}_{T}g=O(\text{\textgreek{t}}^{-\text{\textgreek{d}}_{0}}) & \Big\{ O(r^{-1-a})dvdu+O(r)d\text{\textgreek{sv}}d\text{\textgreek{sv}}+O(1)dud\text{\textgreek{sv}}+\label{eq:DeformationTensorTAway-3}\\
 & +O(r^{-a})dvd\text{\textgreek{sv}}+O(r^{-1})du^{2}+O(r^{-2-a})dv^{2}\Big\}.\nonumber 
\end{align}
Thus, in view of the integrated local energy decay assumption \ref{enu:IntegratedLocalEnergyDecayImproved}
(which is satisfied in view of our assumption that $\text{\textgreek{f}}$
has boundary conditions on $\partial_{tim}\mathcal{M}$ belonging
to the class $\mathcal{C}_{ILED}^{(T,K)}$) and Lemma \ref{lem:MorawetzDrLemmaHyperboloids}
in the region $\{r\gg1\}$ (see the remark below that Lemma for the
case when $T$ has deformation tensor with slow polynomial decay in
$\bar{t}$), we can bound for any $(e_{1},\ldots e_{l-1})\in\{0,1\}^{l-1}$
and any $R_{f}>0$ to be fixed later (setting $X^{(0)}=T$ and $X^{(1)}=K_{R_{c}}$):
\begin{equation}
\begin{split}\sum_{j=1}^{m}\int_{\{\text{\textgreek{t}}_{1}\le\bar{t}\le\text{\textgreek{t}}_{2}\}} & r_{+}^{-1-\text{\textgreek{e}}}\big(\big|\nabla_{g}^{j}(\mathcal{L}_{X^{(e_{1})}\ldots X^{(e_{l-1})}}\text{\textgreek{f}})\big|_{h}^{2}+r_{+}^{-2}\big|\mathcal{L}_{X^{(e_{1})}\ldots X^{(e_{l-1})}}\text{\textgreek{f}}\big|^{2}\big)\, dg+\\
+\sum_{j=0}^{m}\sum_{j_{1}+j_{2}=j} & \int_{\partial_{tim}\mathcal{M}\cap\{\text{\textgreek{t}}_{1}\le\bar{t}\le\text{\textgreek{t}}_{2}\}}\big|\nabla_{g}^{j}(\mathcal{L}_{X^{(e_{1})}\ldots X^{(e_{l-1})}}\text{\textgreek{f}})\big|_{h}^{2}\, dh_{\partial_{tim}\mathcal{M}}\le\\
\le\, & C_{\text{\textgreek{e}},m,R_{f}}\sum_{j=0}^{m+k-1}\sum_{j_{1}+j_{2}=j}\int_{\{\bar{t}=\text{\textgreek{t}}_{1}\}}\Big(\big|\nabla_{h_{N}}^{j_{1}+1}(N^{j_{2}}\mathcal{L}_{X^{(e_{1})}\ldots X^{(e_{l-1})}}\text{\textgreek{f}})\big|_{h_{N}}^{2}+r_{+}^{-2}\big|N^{j+1}\mathcal{L}_{X^{(e_{1})}\ldots X^{(e_{l-1})}}\text{\textgreek{f}}|^{2}\Big)\, dh_{N}+\\
 & +C_{\text{\textgreek{e}},m}\text{\textgreek{t}}_{1}^{-\text{\textgreek{d}}_{0}}\sum_{j=0}^{m-1}\sum_{j_{1}+j_{2}=j}\int_{\{\text{\textgreek{t}}_{1}\le\bar{t}\le\text{\textgreek{t}}_{2}\}\cap\{r\gtrsim1\}}r_{+}^{-1}\big(\big|\nabla_{h_{\text{\textgreek{t}},N}}^{1+j_{1}}\big(T^{j_{2}}(\mathcal{L}_{X^{(e_{1})}\ldots X^{(e_{l-1})}}\text{\textgreek{f}})\big)\big|_{h_{\text{\textgreek{t}},N}}^{2}+r_{+}^{-2}\Big(\big|T^{m}\big(\mathcal{L}_{X^{(e_{1})}\ldots X^{(e_{l-1})}}\text{\textgreek{f}}\big)\big|^{2}+\\
 & \hphantom{+C_{\text{\textgreek{e}},m}\text{\textgreek{t}}_{1}^{-\text{\textgreek{d}}_{0}}\sum_{j=0}^{m-1}\sum_{j_{1}+j_{2}=j}\int_{\{\text{\textgreek{t}}_{1}\le\bar{t}\le\text{\textgreek{t}}_{2}\}\cap\{r\gtrsim1\}}r_{+}^{-1}\big(}+\big|\mathcal{L}_{X^{(e_{1})}\ldots X^{(e_{l-1})}}\text{\textgreek{f}}\big|^{2}\Big)\big)\, dg+\\
 & +C_{\text{\textgreek{e}},m,R_{f}}\sum_{j=0}^{k+m-1}\int_{\{\text{\textgreek{t}}_{1}\le\bar{t}\le\text{\textgreek{t}}_{2}\}}r_{+}^{1+\text{\textgreek{e}}}\big|\nabla_{g}^{j}(\square\mathcal{L}_{X^{(e_{1})}\ldots X^{(e_{l-1})}}\text{\textgreek{f}})\big|_{h}^{2}\, dg\Big\}.
\end{split}
\label{eq:FirstIled}
\end{equation}

Using Lemma \ref{lem:BoundednessWithLoss} for $\mathcal{L}_{X^{(e_{1})}\ldots X^{(e_{l-1})}}\text{\textgreek{f}}$
in place of $\text{\textgreek{f}}$ (and adapting the proof a bit
so as to use (\ref{eq:FirstIled}) instead of simply (\ref{eq:IntegratedLocalEnergyDecayImprovedDecay})
and repeating the proof of Lemma \ref{lem:MorawetzDrLemmaHyperboloids}),
we obtain the following energy boundedness statement for any $\text{\textgreek{g}}_{0}>0$
to be fixed small later: 
\begin{equation}
\begin{split}\sum_{j=0}^{m-1}\sum_{j_{1}+j_{2}=j}\int_{\{\bar{t}=\text{\textgreek{t}}_{2}\}} & \Big(\big|\nabla_{h_{N}}^{j_{1}+1}(N^{j_{2}}\mathcal{L}_{X^{(e_{1})}\ldots X^{(e_{l-1})}}\text{\textgreek{f}})\big|_{h_{N}}^{2}+r_{+}^{-2}\big|N^{j+1}\mathcal{L}_{X^{(e_{1})}\ldots X^{(e_{l-1})}}\text{\textgreek{f}}|^{2}\Big)\, dh_{N}\le\\
\le C_{m,\text{\textgreek{e}}}\Bigg\{ & \sum_{j=0}^{m+k-1}\sum_{j_{1}+j_{2}=j}\int_{\{\bar{t}=\text{\textgreek{t}}_{1}\}}\Big(\big|\nabla_{h_{N}}^{j_{1}+1}(N^{j_{2}}\mathcal{L}_{X^{(e_{1})}\ldots X^{(e_{l-1})}}\text{\textgreek{f}})\big|_{h_{N}}^{2}+r_{+}^{-2}\big|N^{j+1}\mathcal{L}_{X^{(e_{1})}\ldots X^{(e_{l-1})}}\text{\textgreek{f}}|^{2}\Big)\, dh_{N}+\\
 & +\text{\textgreek{t}}_{1}^{-\text{\textgreek{d}}_{0}}\sum_{j=0}^{m-1}\sum_{j_{1}+j_{2}=j}\int_{\{\text{\textgreek{t}}_{1}\le\bar{t}\le\text{\textgreek{t}}_{2}\}}r_{+}^{-1+\text{\textgreek{e}}}\Big(\big|\nabla_{h_{N}}^{j_{1}+1}(N^{j_{2}}\mathcal{L}_{X^{(e_{1})}\ldots X^{(e_{l-1})}}\text{\textgreek{f}})\big|_{h_{N}}^{2}+r_{+}^{-2}\Big(\big|N^{j+1}\big(\mathcal{L}_{X^{(e_{1})}\ldots X^{(e_{l-1})}}\text{\textgreek{f}}\big)\big|^{2}+\\
 & \hphantom{+\text{\textgreek{t}}_{1}^{-\text{\textgreek{d}}_{0}}\sum_{j=0}^{m-1}\sum_{j_{1}+j_{2}=j}\int_{\{\text{\textgreek{t}}_{1}\le\bar{t}\le\text{\textgreek{t}}_{2}\}}r_{+}^{-1+\text{\textgreek{e}}}\Big(}+\big|\mathcal{L}_{X^{(e_{1})}\ldots X^{(e_{l-1})}}\text{\textgreek{f}}\big|^{2}\Big)\Big)\, dg+\\
 & +\sum_{j=0}^{m+k-1}\int_{\{\text{\textgreek{t}}_{1}\le\bar{t}\le\text{\textgreek{t}}_{2}\}}r_{+}^{1+\text{\textgreek{e}}}\big|\nabla_{g}^{j}(\square\mathcal{L}_{X^{(e_{1})}\ldots X^{(e_{l-1})}}\text{\textgreek{f}})\big|_{h}^{2}\, dg\Bigg\}.
\end{split}
\label{eq:BoundednessWithLoss-2}
\end{equation}

Fixing some $\text{\textgreek{d}},\text{\textgreek{h}}>0$ small in
terms of $\text{\textgreek{e}}$, we obtain from (\ref{eq:newMethodFinalStatementHyperboloids})
for $p=\text{\textgreek{e}}$ and for $T^{j}\text{\textgreek{f}}$
in place of $\text{\textgreek{f}}$ for $j=l-1,\ldots m+l-2$: 
\begin{equation}
\begin{split}\sum_{j=l-1}^{m+l-2}\int_{\{\text{\textgreek{t}}_{1}\le\bar{t}\le\text{\textgreek{t}}_{2}\}\cap\{r\ge2R_{c}\}} & \big(r_{+}^{-1+\text{\textgreek{e}}}\big|\nabla_{h_{\text{\textgreek{t}},N}}(T^{j}\text{\textgreek{f}})\big|_{h_{\text{\textgreek{t}},N}}^{2}+r_{+}^{-1-\text{\textgreek{e}}}|T^{j+1}\text{\textgreek{f}}|^{2}+r_{+}^{-3+\text{\textgreek{e}}}\big|T^{j}\text{\textgreek{f}}\big|^{2}\big)\, dg+\\
+\sum_{j=l-1}^{m+l-2}\int_{\{\bar{t}=\text{\textgreek{t}}_{2}\}\cap\{r\ge2R_{c}\}} & r_{+}^{\text{\textgreek{e}}}\big(\big|L(T^{j}\text{\textgreek{f}})\big|^{2}+r_{+}^{-2}\big|\nabla_{h_{\text{\textgreek{t}},N}}(T^{j}\text{\textgreek{f}})\big|_{h_{\text{\textgreek{t}},N}}^{2}+r_{+}^{-2}|T^{j}\text{\textgreek{f}}|^{2}\big)\, dh_{N}\le\\
\le C_{\text{\textgreek{e}},m}\Bigg\{ & \sum_{j=l-1}^{m+l-2}\int_{\{\bar{t}=\text{\textgreek{t}}_{1}\}\cap\{r\ge R_{c}\}}\big(r_{+}^{\text{\textgreek{e}}}\big|\nabla_{h_{\text{\textgreek{t}},N}}(T^{j}\text{\textgreek{f}})\big|_{h_{\text{\textgreek{t}},N}}^{2}+r_{+}^{-2}|T^{j+1}\text{\textgreek{f}}|^{2}+r_{+}^{-2+\text{\textgreek{e}}}\big|T^{j}\text{\textgreek{f}}\big|^{2}\big)\, dh_{N}+\\
 & +\sum_{j=l-1}^{m+l-2}\int_{\{\text{\textgreek{t}}_{1}\le\bar{t}\le\text{\textgreek{t}}_{2}\}\cap\{R_{c}\le r\le2R_{c}\}}r_{+}^{-1+\text{\textgreek{e}}}\big(\big|\nabla_{h_{\text{\textgreek{t}},N}}(T^{j}\text{\textgreek{f}})\big|_{h_{\text{\textgreek{t}},N}}^{2}+|T^{j+1}\text{\textgreek{f}}|^{2}+r_{+}^{-2}\big|T^{j}\text{\textgreek{f}}\big|^{2}\big)\, dg+\\
 & +\sum_{j=l-1}^{m+l-2}\int_{\{\text{\textgreek{t}}_{1}\le\bar{t}\le\text{\textgreek{t}}_{2}\}}r_{+}^{1+\text{\textgreek{e}}}\big|\square(T^{j}\text{\textgreek{f}})\big|^{2}\, dg\Bigg\}.
\end{split}
\label{eq:FromNewMethod}
\end{equation}

Using the elliptic estimate (\ref{eq:EllipticEstimatesWithoutKillingError})
repeatedly together with the expression (\ref{eq:Kcommuted}) for
the commutation of $\square$ with $\mathcal{L}_{X^{(e_{i})}}$ at
each step (notice that Theorem \ref{thm:FriedlanderRadiation} applies
to yield $\limsup_{r\rightarrow+\infty}\big|r^{\frac{d-1}{2}+j_{1}}\nabla_{h_{\text{\textgreek{t}},N}}^{j_{1}}(T^{j_{2}}\text{\textgreek{f}})\big|_{h_{\text{\textgreek{t}},N}}<+\infty$
for all $j_{1}+j_{2}\le l+m+k$), in view of the condition (\ref{eq:BoundaryConditionImprovedDecay-1-1})
for $\text{\textgreek{f}}$ at $\partial_{tim}\mathcal{M}$, we can
bound after summing over all possible combinations of $(e_{1},\ldots e_{l-1})\in\{0,1\}^{l-1}$:
\begin{equation}
\begin{split}\sum_{0\le j_{1}+j_{2}\le m-1}\, & \sum_{0\le i_{1}+i_{2}\le l-1}\int_{\{\text{\textgreek{t}}_{1}\le\bar{t}\le\text{\textgreek{t}}_{2}\}}\Big(r_{+}^{-1+\text{\textgreek{e}}-2(l-1-i_{1}-i_{2})}\big|\nabla_{h_{\text{\textgreek{t}},N}}^{i_{1}+j_{1}+1}(K_{R_{c}}^{i_{2}+j_{2}}\text{\textgreek{f}})\big|_{\big(1-\log(r_{tim})\big)\cdot h_{R_{c}}}^{2}+r_{+}^{-1-\text{\textgreek{e}}}|K_{R}^{l+j_{1}+j_{2}}\text{\textgreek{f}}|^{2}\Big)\, dg\le\\
\le\, & \sum_{(e_{1},\ldots e_{l-1})\in\{0,1\}^{l-1}}\Bigg\{\sum_{j=1}^{m}\int_{\{\text{\textgreek{t}}_{1}\le\bar{t}\le\text{\textgreek{t}}_{2}\}}r_{+}^{-1-\text{\textgreek{e}}}\big(\big|\nabla_{g}^{j}(\mathcal{L}_{X^{(e_{1})}\ldots X^{(e_{l-1})}}\text{\textgreek{f}})\big|_{h}^{2}+r_{+}^{-2}\big|\mathcal{L}_{X^{(e_{1})}\ldots X^{(e_{l-1})}}\text{\textgreek{f}}\big|^{2}\big)\, dg\Bigg\}+\\
 & +C_{\text{\textgreek{e}},m}\sum_{j=l-1}^{m+l-2}\int_{\{\text{\textgreek{t}}_{1}\le\bar{t}\le\text{\textgreek{t}}_{2}\}\cap\{r\ge2R_{c}\}}r_{+}^{-1+\text{\textgreek{e}}}\big(\big|\nabla_{h_{\text{\textgreek{t}},N}}(T^{j}\text{\textgreek{f}})\big|_{h_{\text{\textgreek{t}},N}}^{2}+r_{+}^{-2}\big|T^{j}\text{\textgreek{f}}\big|^{2}\big)\, dg+\\
 & +C_{\text{\textgreek{e}},m}\sum_{0\le i_{1}+i_{2}\le m-1}\sum_{0\le j\le l-2}\int_{\{\text{\textgreek{t}}_{1}\le\bar{t}\le\text{\textgreek{t}}_{2}\}}r_{+}^{-1+\text{\textgreek{e}}}\big(r_{+}^{-2-2(l-2-j)}\text{\textgreek{t}}^{-2\text{\textgreek{d}}_{0}}+r_{+}^{-2}\text{\textgreek{t}}^{-2(\text{\textgreek{d}}_{0}+l-2-j)}\big)\big|\nabla_{h_{\text{\textgreek{t}},N}}^{j+i_{1}+1}(K_{R_{c}}^{i_{2}}\text{\textgreek{f}})\big|_{h_{\text{\textgreek{t}},N}}^{2}\, dg+\\
 & +C_{\text{\textgreek{e}},m}\sum_{0\le i_{1}+i_{2}\le m-1}\sum_{0\le j_{1}+j_{2}\le l-2}\int_{\{\text{\textgreek{t}}_{1}\le\bar{t}\le\text{\textgreek{t}}_{2}\}}r_{+}^{-3+\text{\textgreek{e}}}\text{\textgreek{t}}^{-2(\text{\textgreek{d}}_{0}+l-j_{1}-j_{2}-2)}\Big(|\nabla_{h_{\text{\textgreek{t}},N}}^{j_{1}+i_{1}}\text{\textgreek{f}}|_{T,K,R_{c}}^{(j_{2}+i_{2}+1)}\Big)^{2}\, dg+\\
 & +C_{\text{\textgreek{e}},m}\sum_{j=l-2}^{m+l-3}\int_{\{\text{\textgreek{t}}_{1}\le\bar{t}\le\text{\textgreek{t}}_{2}\}}r_{+}^{-1+\text{\textgreek{e}}}\big|\nabla_{g}^{j}F\big|_{h}^{2}\, dg,
\end{split}
\label{eq:EllipticEstimatesForILEDBulk}
\end{equation}
 and 
\begin{equation}
\begin{split}\sum_{0\le j_{1}+j_{2}\le m-1}\, & \sum_{0\le i_{1}+i_{2}\le l-1}\int_{\{\bar{t}=\text{\textgreek{t}}_{2}\}}r_{+}^{-2(l-1-i_{1}-i_{2})}\Big(\big|\nabla_{h_{\text{\textgreek{t}},N}}^{i_{1}+j_{1}+1}(K_{R_{c}}^{i_{2}+j_{2}}\text{\textgreek{f}})\big|_{\big(1-\log(r_{tim})\big)\cdot h_{R_{c}}}^{2}+r_{+}^{-2}\big|K_{R_{c}}^{l+j}\text{\textgreek{f}}\big|^{2}\Big)\, dh_{N}\le\\
\le\, & C_{m}\sum_{(e_{1},\ldots e_{l-1})\in\{0,1\}^{l-1}}\Bigg\{\sum_{j=0}^{m-1}\sum_{j_{1}+j_{2}=j}\int_{\{\bar{t}=\text{\textgreek{t}}_{2}\}}\Big(\big|\nabla_{h_{N}}^{j_{1}+1}(N^{j_{2}}\mathcal{L}_{X^{(e_{1})}\ldots X^{(e_{l-1})}}\text{\textgreek{f}})\big|_{h_{N}}^{2}+r_{+}^{-2}\big|N^{j}\mathcal{L}_{X^{(e_{1})}\ldots X^{(e_{l-1})}}\text{\textgreek{f}}|^{2}\Big)\, dh_{N}\Bigg\}+\\
 & +C_{m}\sum_{0\le i_{1}+i_{2}\le m-1}\sum_{0\le j\le l-2}\int_{\{\bar{t}=\text{\textgreek{t}}_{2}\}}\big(r_{+}^{-2-2(l-2-j)}\text{\textgreek{t}}^{-2\text{\textgreek{d}}_{0}}+r_{+}^{-2}\text{\textgreek{t}}^{-2(\text{\textgreek{d}}_{0}+l-2-j)}\big)\big|\nabla_{h_{\text{\textgreek{t}},N}}^{j+i_{1}+1}(K_{R_{c}}^{i_{2}}\text{\textgreek{f}})\big|_{h_{\text{\textgreek{t}},N}}^{2}\, dh_{N}+\\
 & +C_{m}\sum_{0\le i_{1}+i_{2}\le m-1}\sum_{0\le j_{1}+j_{2}\le l-2}\int_{\{\bar{t}=\text{\textgreek{t}}_{2}\}}r_{+}^{-2}\text{\textgreek{t}}^{-2(\text{\textgreek{d}}_{0}+l-j_{1}-j_{2}-2)}\Big(|\nabla_{h_{\text{\textgreek{t}},N}}^{j_{1}+i_{1}}\text{\textgreek{f}}|_{T,K,R_{c}}^{(j_{2}+i_{2}+1)}\Big)^{2}\, dh_{N}+\\
 & +C_{m}\sum_{j=l-2}^{m+l-3}\int_{\{\bar{t}=\text{\textgreek{t}}_{2}\}}\big|\nabla_{g}^{j}F\big|_{h}^{2}\, dh_{N}.
\end{split}
\label{eq:EllipticEstimatesForBoundedness}
\end{equation}
 Moreover, using the expression (\ref{eq:ConformalWaveOperator})
for the wave equation in the region $\{r\gg1\}$ and Lemma \ref{lem:Commutator expressions}
to estimate the commutator of $\square$ with $L$ (and recalling
that $L$ is supported in the region $\{r\ge2R_{c}\}$), we can also
bound in view of Lemma \ref{lem:EllipticEstimatesWithoutKillingError}:
\begin{equation}
\begin{split}\sum_{j=0}^{m-1}\sum_{1\le j_{2}+j_{3}\le l+j-1}\, & \int_{\{\bar{t}=\text{\textgreek{t}}_{2}\}}r_{+}^{\text{\textgreek{e}}-2(l+j-1-j_{2}-j_{3})}\big|\mathcal{L}_{L}\nabla_{h_{\text{\textgreek{t}},N}}^{j_{2}}(K_{R_{c}}^{j_{3}}\text{\textgreek{f}})\big|_{\big(1-\log(r_{tim})\big)\cdot h_{R_{c}}}^{2}\, dh_{N}\le\\
\le C_{\text{\textgreek{e}},m}\Bigg\{ & \sum_{j=l-1}^{m+l-2}\int_{\{\bar{t}=\text{\textgreek{t}}_{2}\}\cap\{r\ge2R_{c}\}}r_{+}^{\text{\textgreek{e}}}\big(\big|L(T^{j}\text{\textgreek{f}})\big|^{2}+r_{+}^{-2}\big|\nabla_{h_{\text{\textgreek{t}},N}}(T^{j}\text{\textgreek{f}})\big|_{h_{\text{\textgreek{t}},N}}^{2}+r_{+}^{-2-\text{\textgreek{e}}}|T^{j+1}\text{\textgreek{f}}|^{2}\big)\, dh_{N}+\\
 & +\sum_{(e_{1},\ldots e_{l-1})\in\{0,1\}^{l-1}}\Bigg\{\sum_{j=0}^{m-1}\sum_{j_{1}+j_{2}=j}\int_{\{\bar{t}=\text{\textgreek{t}}_{2}\}}\Big(\big|\nabla_{h_{N}}^{j_{1}+1}(N^{j_{2}}\mathcal{L}_{X^{(e_{1})}\ldots X^{(e_{l-1})}}\text{\textgreek{f}})\big|_{h_{N}}^{2}+r_{+}^{-2}\big|N^{j}\mathcal{L}_{X^{(e_{1})}\ldots X^{(e_{l-1})}}\text{\textgreek{f}}|^{2}\Big)\, dh_{N}+\\
 & +\sum_{0\le i_{1}+i_{2}\le m-1}\sum_{0\le j\le l-2}\int_{\{\bar{t}=\text{\textgreek{t}}_{2}\}}\big(r_{+}^{\text{\textgreek{e}}-2-2(l-2-j)}\text{\textgreek{t}}^{-2\text{\textgreek{d}}_{0}}+r_{+}^{\text{\textgreek{e}}-2}\text{\textgreek{t}}^{-2(\text{\textgreek{d}}_{0}+l-2-j)}\big)\big|\nabla_{h_{\text{\textgreek{t}},N}}^{j+i_{1}+1}(K_{R_{c}}^{i_{2}}\text{\textgreek{f}})\big|_{h_{\text{\textgreek{t}},N}}^{2}\, dh_{N}+\\
 & +\sum_{0\le i_{1}+i_{2}\le m-1}\sum_{0\le j_{1}+j_{2}\le l-2}\int_{\{\bar{t}=\text{\textgreek{t}}_{2}\}}r_{+}^{\text{\textgreek{e}}-2}\text{\textgreek{t}}^{-2(\text{\textgreek{d}}_{0}+l-j_{1}-j_{2}-2)}\Big(|\nabla_{h_{\text{\textgreek{t}},N}}^{j_{1}+i_{1}}\text{\textgreek{f}}|_{T,K,R_{c}}^{(j_{2}+i_{2}+1)}\Big)^{2}\, dh_{N}+\\
 & +\sum_{j=l-2}^{m+l-3}\int_{\{\bar{t}=\text{\textgreek{t}}_{2}\}}r_{+}^{\text{\textgreek{e}}}\big|\nabla_{g}^{j}F\big|_{h}^{2}\, dh_{N}\Bigg\}.
\end{split}
\label{eq:EllipticEstimatesForImprovedBoundedness}
\end{equation}

Therefore, fixing $R_{f}$ large enough in terms of $\text{\textgreek{e}},m$,
after adding (for all possible combinations of $(e_{1},\ldots e_{l-1})\in\{0,1\}^{l-1}$)
the estimates (\ref{eq:FirstIled}), (\ref{eq:BoundednessWithLoss-2})
and a small multiple (in terms of $\text{\textgreek{e}},m$) of (\ref{eq:FromNewMethod}),
and using the expression (\ref{eq:Kcommuted}) for the commutation
of $\square$ with $\mathcal{L}_{X^{(e_{i})}}$, we obtain in view
of (\ref{eq:EllipticEstimatesForILEDBulk}), (\ref{eq:EllipticEstimatesForBoundedness})
and (\ref{eq:EllipticEstimatesForImprovedBoundedness}) as well as
a trace-type inequality for the terms $\int_{\{\bar{t}=\text{\textgreek{t}}_{2}\}}|\nabla^{j}F|^{2}$,
$\int_{\{\bar{t}=\text{\textgreek{t}}_{2}\}}\big|\nabla_{h_{\text{\textgreek{t}},N}}^{j+i_{1}+1}(K_{R_{c}}^{i_{2}}\text{\textgreek{f}})\big|_{h_{\text{\textgreek{t}},N}}^{2}$
and $\int_{\{\bar{t}=\text{\textgreek{t}}_{2}\}}\Big(|\nabla_{h_{\text{\textgreek{t}},N}}^{j_{1}+i_{1}}\text{\textgreek{f}}|_{T,K,R_{c}}^{(j_{2}+i_{2}+1)}\Big)^{2}$
(and recalling that $\text{\textgreek{t}}_{1}$was assumed to be large
in terms of $\text{\textgreek{e}},m$): 
\begin{equation}
\begin{split}\sum_{j=0}^{m-1}\int_{\{\text{\textgreek{t}}_{1}\le\bar{t}\le\text{\textgreek{t}}_{2}\}}\Big(\, & \sum_{j_{2}+j_{3}=l+j-1}r_{+}^{-1+\text{\textgreek{e}}}\big|\nabla_{h_{\text{\textgreek{t}},N}}^{j_{2}+1}(K_{R_{c}}^{j_{3}}\text{\textgreek{f}})\big|_{\big(1-\log(r_{tim})\big)\cdot h_{R_{c}}}^{2}+\\
 & +\sum_{1\le j_{2}+j_{3}\le l+j-1}r_{+}^{-1+\text{\textgreek{e}}-2(l+j-j_{2}-j_{3})}\big|\nabla_{h_{\text{\textgreek{t}},N}}^{j_{2}}(K_{R_{c}}^{j_{3}}\text{\textgreek{f}})\big|_{\big(1-\log(r_{tim})\big)\cdot h_{R_{c}}}^{2}+r_{+}^{-1-\text{\textgreek{e}}}\big|K_{R_{c}}^{l+j}\text{\textgreek{f}}\big|^{2}\Big)\, dg+\\
+\sum_{j=0}^{m-1}\sum_{1\le j_{2}+j_{3}\le l+j-1} & \int_{\{\bar{t}=\text{\textgreek{t}}_{2}\}}r_{+}^{-2(l+j-1-j_{2}-j_{3})}\Big(r^{\text{\textgreek{e}}}|\mathcal{L}_{L}\nabla_{h_{\text{\textgreek{t}},N}}^{j_{2}}(K_{R_{c}}^{j_{3}}\text{\textgreek{f}})|_{\big(1-\log(r_{tim})\big)\cdot h_{R_{c}}}^{2}+\\
 & \hphantom{\int_{\{\bar{t}=\text{\textgreek{t}}_{2}\}}r_{+}^{-2(l+j-1-j_{2}-j_{3})}\Big(}+\big|\nabla_{h_{R_{c}}}^{j_{2}+1}(K_{R_{c}}^{j_{3}}\text{\textgreek{f}})\big|_{\big(1-\log(r_{tim})\big)\cdot h_{R_{c}}}^{2}+r_{+}^{-2}\big|K_{R_{c}}^{l+j}\text{\textgreek{f}}\big|^{2}\Big)\, dh_{N}\le\\
\le C_{\text{\textgreek{e}},m}\Big\{ & \sum_{j=0}^{k+m-1}\sum_{1\le j_{2}+j_{3}\le l+j-1}\int_{\{\bar{t}=\text{\textgreek{t}}_{1}\}}r_{+}^{\text{\textgreek{e}}-2(l+j-1-j_{2}-j_{3})}\Big(\big|\nabla_{h_{\text{\textgreek{t}},N}}^{j_{2}+1}(N^{j_{3}}\text{\textgreek{f}})\big|_{h_{\text{\textgreek{t}},N}}^{2}+r_{+}^{-2}\big|N^{l+j}\text{\textgreek{f}}\big|^{2}\Big)\, dh_{N}+\\
 & +\text{\textgreek{t}}_{1}^{-\text{\textgreek{d}}_{0}}\int_{\{\text{\textgreek{t}}_{1}\le\bar{t}\le\text{\textgreek{t}}_{2}\}}r_{+}^{-1+\text{\textgreek{e}}}\mathcal{T}_{T,K,R_{c},sl}^{(l,m+k-1)}[\text{\textgreek{f}}]\, dg+\sum_{j=l-1}^{m+k+l-2}\int_{\{\text{\textgreek{t}}_{1}\le\bar{t}\le\text{\textgreek{t}}_{2}\}}r_{+}^{1+\text{\textgreek{e}}}\big|\nabla_{g}^{j}F\big|_{h}^{2}\, dg\Big\}.
\end{split}
\label{eq:FirstIled-1-1}
\end{equation}

Let $r_{1}>0$ be the small constant appearing in (\ref{eq:RedShiftAssumption}),
and fix a smooth cut-off function $\text{\textgreek{q}}_{r_{1}}:\mathcal{M}\rightarrow[0,1]$
such that $\text{\textgreek{q}}_{r_{1}}$ is a function of $r$ satisfying
$\text{\textgreek{q}}_{r_{1}}\equiv1$ on $\{r\le r_{1}\}$ and $\text{\textgreek{q}}_{r_{1}}\equiv0$
on $\{r\ge2r_{1}\}$. Integrating $\nabla^{\text{\textgreek{m}}}(\text{\textgreek{q}}_{r_{1}}J_{\text{\textgreek{m}}}^{m+l}[\text{\textgreek{f}}])$
over $\mathcal{M}$ (where $J_{\text{\textgreek{m}}}^{m+l}[\text{\textgreek{f}}]$
is given by \ref{eq:EnergyBoundedness current}), and using the red
shift estimate (\ref{eq:RedShiftAssumption}), we obtain: 
\begin{equation}
\begin{split}\int_{\mathcal{H}^{+}(\text{\textgreek{t}}_{1},\text{\textgreek{t}}_{2})}J_{\text{\textgreek{m}}}^{m+l}[\text{\textgreek{f}}]n_{\mathcal{H}}^{\text{\textgreek{m}}}+ & \int_{\{\bar{t}=\text{\textgreek{t}}_{2}\}\cap\{r\le r_{1}\}}J_{\text{\textgreek{m}}}^{m+l}[\text{\textgreek{f}}]\bar{n}^{\text{\textgreek{m}}}+\sum_{j=1}^{m+l-1}\int_{\{\text{\textgreek{t}}_{1}\le\bar{t}\le\text{\textgreek{t}}_{2}\}\cap\{r\le r_{1}\}}\big|\nabla_{g}^{j}\text{\textgreek{f}}\big|_{h}^{2}\, dg\le\\
\le C_{m}\Big\{ & \sum_{j=1}^{m+l-1}\int_{\{\bar{t}=\text{\textgreek{t}}_{1}\}\cap\{r\le2r_{1}\}}\big|\nabla_{g}^{j}\text{\textgreek{f}}\big|_{h}^{2}\, dh_{N}+\sum_{j=1}^{m+l-1}\int_{\{\text{\textgreek{t}}_{1}\le\bar{t}\le\text{\textgreek{t}}_{2}\}\cap\{r_{1}\le r\le2r_{1}\}}\big|\nabla_{g}^{j}\text{\textgreek{f}}\big|_{h}^{2}\, dg+\\
 & +\sum_{j=0}^{m+l-2}\int_{\{\text{\textgreek{t}}_{1}\le\bar{t}\le\text{\textgreek{t}}_{2}\}\cap\{r\le2r_{1}\}}\big|\nabla_{g}^{j}F\big|_{h}^{2}\, dg\Big\}.
\end{split}
\label{eq:RedShiftAssumption-1}
\end{equation}
 Adding to \ref{eq:FirstIled-1-1} a small multiple of the non degenerate
estimate \ref{eq:RedShiftAssumption-1}, and using Assumptions \ref{enu:EllipticEstimates}
and \ref{enu:EllipticEstimatesHorizon} (together with a Hardy type
inequality for the terms $\sum_{j=0}^{l-2}\int_{\{r\le2r_{1}\}}\big|\nabla_{g}^{j}F\big|_{h}^{2}\, dg$),
we obtain (\ref{eq:ILEDCommuted}): 
\begin{equation}
\begin{split}\sum_{j=0}^{m-1}\int_{\{\text{\textgreek{t}}_{1}\le\bar{t}\le\text{\textgreek{t}}_{2}\}}\Big(\, & \sum_{j_{2}+j_{3}=l+j-1}r_{+}^{-1+\text{\textgreek{e}}}\big|\nabla_{h_{\text{\textgreek{t}},N}}^{j_{2}+1}(N^{j_{3}}\text{\textgreek{f}})\big|_{h_{\text{\textgreek{t}},N}}^{2}+\sum_{1\le j_{2}+j_{3}\le l+j-1}r_{+}^{-1+\text{\textgreek{e}}-2(l+j-j_{2}-j_{3})}\big|\nabla_{h_{\text{\textgreek{t}},N}}^{j_{2}}(N^{j_{3}}\text{\textgreek{f}})\big|_{h_{\text{\textgreek{t}},N}}^{2}+r_{+}^{-1-\text{\textgreek{e}}}\big|N^{l+j}\text{\textgreek{f}}\big|^{2}\Big)\, dg+\\
+\sum_{j=0}^{m-1}\sum_{1\le j_{2}+j_{3}\le l+j-1} & \int_{\{\bar{t}=\text{\textgreek{t}}_{2}\}}r_{+}^{-2(l+j-1-j_{2}-j_{3})}\Big(r_{+}^{\text{\textgreek{e}}}|\mathcal{L}_{L}\nabla_{h_{\text{\textgreek{t}},N}}^{j_{2}}(N^{j_{3}}\text{\textgreek{f}})|^{2}+\big|\nabla_{h_{\text{\textgreek{t}},N}}^{j_{2}+1}(N^{j_{3}}\text{\textgreek{f}})\big|_{h_{\text{\textgreek{t}},N}}^{2}+r_{+}^{-2}\big|N^{l+j}\text{\textgreek{f}}\big|^{2}\Big)\, dh_{N}+\\
+\sum_{j=0}^{m-1}\sum_{1\le j_{2}+j_{3}\le l+j-1} & \int_{\mathcal{H}^{+}(\text{\textgreek{t}}_{1},\text{\textgreek{t}}_{2})}\big|\nabla_{h_{\mathcal{H}}}^{j_{2}+1}(K_{R_{c}}^{j_{3}}\text{\textgreek{f}})\big|_{h_{\mathcal{H}}}^{2}\, dh_{\mathcal{H}}\le\\
\le C_{\text{\textgreek{e}},m}\Big\{\, & \sum_{j=0}^{k+m-1}\sum_{1\le j_{2}+j_{3}\le l+j-1}\int_{\{\bar{t}=\text{\textgreek{t}}_{1}\}}r_{+}^{\text{\textgreek{e}}-2(l+j-1-j_{2}-j_{3})}\Big(\big|\nabla_{h_{\text{\textgreek{t}},N}}^{j_{2}+1}(N^{j_{3}}\text{\textgreek{f}})\big|_{h_{\text{\textgreek{t}},N}}^{2}+r_{+}^{-2}\big|N^{l+j}\text{\textgreek{f}}\big|^{2}\Big)\, dh_{N}+\\
 & +\text{\textgreek{t}}_{1}^{-\text{\textgreek{d}}_{0}}\int_{\{\text{\textgreek{t}}_{1}\le\bar{t}\le\text{\textgreek{t}}_{2}\}}r_{+}^{-1+\text{\textgreek{e}}}\mathcal{T}_{T,K,R_{c},sl}^{(l,m+k-1)}[\text{\textgreek{f}}]\, dg+\sum_{j=l-1}^{m+k+l-2}\int_{\{\text{\textgreek{t}}_{1}\le\bar{t}\le\text{\textgreek{t}}_{2}\}}r_{+}^{1+\text{\textgreek{e}}}\big|\nabla_{g}^{j}F\big|_{h}^{2}\, dg\Big\}.
\end{split}
\label{eq:FirstIled-1-1-1}
\end{equation}

Moreover, by adding (for all possible combinations of $(e_{1},\ldots e_{l-1})\in\{0,1\}^{l-1}$)
the estimates (\ref{eq:FirstIled}), (\ref{eq:BoundednessWithLoss-2})
and $\text{\textgreek{t}}_{1}^{-\text{\textgreek{d}}_{0}}$ times
a small multiple of (\ref{eq:FromNewMethod}), we also obtain the
degenerate energy boundedness statement (\ref{eq:DegenerateBoundednessWithLoss}):
\begin{equation}
\begin{split}\sum_{j=0}^{m-1}\sum_{1\le j_{2}+j_{3}\le l+j-1}\, & \int_{\{\bar{t}=\text{\textgreek{t}}_{2}\}}r_{+}^{-2(l+j-1-j_{2}-j_{3})}\Big(\big|\nabla_{h_{\text{\textgreek{t}},N}}^{j_{2}+1}(K_{R_{c}}^{j_{3}}\text{\textgreek{f}})\big|_{\big(1-\log(r_{tim})\big)\cdot h_{R_{c}}}^{2}+r_{+}^{-2}\big|K_{R_{c}}^{l+j}\text{\textgreek{f}}\big|^{2}\Big)\, dh_{N}\le\\
\le C_{\text{\textgreek{e}},m}\Big\{\, & \sum_{j=0}^{k+m-1}\sum_{1\le j_{2}+j_{3}\le l+j-1}\int_{\{\bar{t}=\text{\textgreek{t}}_{1}\}}r_{+}^{-2(l+j-1-j_{2}-j_{3})}\Big(\big|\nabla_{h_{\text{\textgreek{t}},N}}^{j_{2}+1}(N^{j_{3}}\text{\textgreek{f}})\big|_{h_{\text{\textgreek{t}},N}}^{2}+r_{+}^{-2}\big|N^{l+j}\text{\textgreek{f}}\big|^{2}\Big)\, dh_{N}+\\
 & +\text{\textgreek{t}}_{1}^{-\text{\textgreek{d}}_{0}}\sum_{j=0}^{m-1}\sum_{1\le j_{2}+j_{3}\le l+j-1}\int_{\{\bar{t}=\text{\textgreek{t}}_{1}\}}r_{+}^{\text{\textgreek{e}}-2(l+j-1-j_{2}-j_{3})}\Big(\big|\nabla_{h_{\text{\textgreek{t}},N}}^{j_{2}+1}(N^{j_{3}}\text{\textgreek{f}})\big|_{h_{\text{\textgreek{t}},N}}^{2}+r_{+}^{-2}\big|N^{l+j}\text{\textgreek{f}}\big|^{2}\Big)\, dh_{N}+\\
 & +\text{\textgreek{t}}_{1}^{-\text{\textgreek{d}}_{0}}\int_{\{\text{\textgreek{t}}_{1}\le\bar{t}\le\text{\textgreek{t}}_{2}\}}r_{+}^{-1+\text{\textgreek{e}}}\mathcal{T}_{T,K,R_{c},sl}^{(l,m+k-1)}[\text{\textgreek{f}}]\, dg+\sum_{j=l-1}^{m+k+l-2}\int_{\{\text{\textgreek{t}}_{1}\le\bar{t}\le\text{\textgreek{t}}_{2}\}}r_{+}^{1+\text{\textgreek{e}}}\big|\nabla_{g}^{j}F\big|_{h}^{2}\, dg\Big\}.
\end{split}
\label{eq:FirstIled-1-1-2}
\end{equation}
 
\end{proof}

\subsection{\label{sub:ProofImprovedEnergyDecay}Proof of Theorem \ref{thm:ImprovedDecayEnergy}
on improved polynomial decay}

Without loss of generality, we will assume that $\text{\textgreek{f}}$
is real valued. Moreover, in order to avoid confusion with unnecessarily
complicated notations, we will assume that $F\equiv0$, since the
proof of (\ref{eq:ImprovedNonDegenerateEnergyDecay}) and (\ref{eq:ImprovedEnergyDecay})
in the case $F\neq0$ follows by repeating exactly the same steps.

The proof of Theorem \ref{thm:ImprovedDecayEnergy} will proceed by
induction on $q$, from $q=1$ up to $q=\lfloor\frac{d+1}{2}\rfloor$.

More precisely, we will assume the following inductive hypothesis
for some integer $1\le q_{0}\le\lfloor\frac{d+1}{2}\rfloor$: 
\begin{description}
\item [{Inductive~hypothesis:}] For any integer $1\le q\le q_{0}-1$,
any $0<\text{\textgreek{e}}\ll\text{\textgreek{d}}_{0}$, any integer
$m\ge1$, any $0\le p<2q-1$ and any $0\le\text{\textgreek{t}}_{1}\le\text{\textgreek{t}}_{2}$
the following bounds hold: 
\begin{equation}
\mathcal{E}_{en}^{(p,q,m)}[\text{\textgreek{f}}](\text{\textgreek{t}}_{2})+\int_{\text{\textgreek{t}}_{2}}^{\text{\textgreek{t}}_{3}}\mathcal{E}_{en}^{(p-1,q,m)}[\text{\textgreek{f}}](\text{\textgreek{t}})\, d\text{\textgreek{t}}\lesssim_{m,\text{\textgreek{e}}}(\text{\textgreek{t}}_{2}-\text{\textgreek{t}}_{1})^{-2q+p+C_{m}\text{\textgreek{e}}}\mathcal{E}_{bound}^{(2q,q,m+\lceil\text{\textgreek{d}}_{0}^{-1}\cdot2(q-1)\rceil(3q+1)\cdot k)}[\text{\textgreek{f}}](\text{\textgreek{t}}_{1}),\label{eq:InductiveHypothesisHigherOrder}
\end{equation}
\begin{equation}
\Bigg\{\mathcal{E}_{en}^{(0,q,m)}[\text{\textgreek{f}}](\text{\textgreek{t}}_{2})+\mathcal{E}_{bound}^{(\text{\textgreek{e}},q,m)}[\text{\textgreek{f}}](\text{\textgreek{t}}_{2})+\int_{\text{\textgreek{t}}_{2}}^{+\infty}\mathcal{E}_{en}^{(-1+\text{\textgreek{e}},q,m)}[\text{\textgreek{f}}](s)\, ds\Bigg\}\lesssim_{m,\text{\textgreek{e}}}(\text{\textgreek{t}}_{2}-\text{\textgreek{t}}_{1})^{-2q+C_{m}\text{\textgreek{e}}}\mathcal{E}_{bound}^{(2q,q,m+\lceil\text{\textgreek{d}}_{0}^{-1}\cdot2(q-1)\rceil(3q+1)\cdot k)}[\text{\textgreek{f}}](\text{\textgreek{t}}_{1})\label{eq:InductiveHypothesisNonDegenerate}
\end{equation}
and 
\begin{equation}
\mathcal{E}_{en,deg}^{(0,q,m)}[\text{\textgreek{f}}](\text{\textgreek{t}}_{2})\lesssim_{m}(\text{\textgreek{t}}_{2}-\text{\textgreek{t}}_{1})^{-2q}\mathcal{E}_{bound}^{(2q,q,m+\lceil\text{\textgreek{d}}_{0}^{-1}\cdot2(q-1)\rceil(3q+1)\cdot k)}[\text{\textgreek{f}}](\text{\textgreek{t}}_{1}).\label{eq:InductiveHypothesisDegenerate}
\end{equation}

\end{description}
Granted this inductve hypothesis, the inductive step of our induction
scheme will be the following: 
\begin{description}
\item [{Inductive~step:}] For any $0<\text{\textgreek{e}}\ll\text{\textgreek{d}}_{0}$,
any integer $m\ge1$, any $0\le p<2q_{0}-1$ and any $0\le\text{\textgreek{t}}_{1}\le\text{\textgreek{t}}_{2}$
the following bounds hold: 
\begin{equation}
\mathcal{E}_{en}^{(p,q_{0},m)}[\text{\textgreek{f}}](\text{\textgreek{t}}_{2})+\int_{\text{\textgreek{t}}_{2}}^{\text{\textgreek{t}}_{3}}\mathcal{E}_{en}^{(p-1,q_{0},m)}[\text{\textgreek{f}}](\text{\textgreek{t}})\, d\text{\textgreek{t}}\lesssim_{m,\text{\textgreek{e}}}(\text{\textgreek{t}}_{2}-\text{\textgreek{t}}_{1})^{-2q_{0}+p+C_{m}\text{\textgreek{e}}}\mathcal{E}_{bound}^{(2q_{0},q_{0},m+\lceil\text{\textgreek{d}}_{0}^{-1}\cdot2(q_{0}-1)\rceil(3q_{0}+1)\cdot k)}[\text{\textgreek{f}}](\text{\textgreek{t}}_{1}),\label{eq:InductiveStepHigherOrder}
\end{equation}
\begin{equation}
\begin{split}\Bigg\{\mathcal{E}_{en}^{(0,q_{0},m)}[\text{\textgreek{f}}](\text{\textgreek{t}}_{2})+\mathcal{E}_{bound}^{(\text{\textgreek{e}},q_{0},m)}[\text{\textgreek{f}} & ](\text{\textgreek{t}}_{2})+\int_{\text{\textgreek{t}}_{2}}^{+\infty}\mathcal{E}_{en}^{(-1+\text{\textgreek{e}},q_{0},m)}[\text{\textgreek{f}}](s)\, ds\Bigg\}\lesssim_{m,\text{\textgreek{e}}}\\
\lesssim_{m,\text{\textgreek{e}}}\, & (\text{\textgreek{t}}_{2}-\text{\textgreek{t}}_{1})^{-2q_{0}+C_{m}\text{\textgreek{e}}}\mathcal{E}_{bound}^{(2q_{0},q_{0},m+\lceil\text{\textgreek{d}}_{0}^{-1}\cdot2(q_{0}-1)\rceil(3q_{0}+1)\cdot k)}[\text{\textgreek{f}}](\text{\textgreek{t}}_{1})
\end{split}
\label{eq:InductiveStepNonDegenerate}
\end{equation}
and 
\begin{equation}
\mathcal{E}_{en,deg}^{(0,q_{0},m)}[\text{\textgreek{f}}](\text{\textgreek{t}}_{2})\lesssim_{m}(\text{\textgreek{t}}_{2}-\text{\textgreek{t}}_{1})^{-2q_{0}}\mathcal{E}_{bound}^{(2q_{0},q_{0},m+\lceil\text{\textgreek{d}}_{0}^{-1}\cdot2(q_{0}-1)\rceil(3q_{0}+1)\cdot k)}[\text{\textgreek{f}}](\text{\textgreek{t}}_{1}).\label{eq:InductiveStepDegenerate}
\end{equation}

\end{description}
Note that for $q=1$ Theorem \ref{thm:ImprovedDecayEnergy} degenerates
to Theorem \ref{thm:SlowPointwiseDecayHighDerivativesNewMethod},
and thus the basis of the induction has already been established.

\subsubsection{Proof of the inductive step}

In order to establish inequalities (\ref{eq:InductiveStepHigherOrder}),
(\ref{eq:InductiveStepNonDegenerate}) and (\ref{eq:InductiveStepDegenerate})
of the inductive step, we will first need to prove a series of lemmas.
Without loss of generality, we will assume that $\text{\textgreek{t}}_{1}=0$
in (\ref{eq:InductiveStepHigherOrder}), (\ref{eq:InductiveStepNonDegenerate})
and (\ref{eq:InductiveStepDegenerate}).
\begin{lem}
\label{lem:InductionOnNumberOfKillingDerivatives}There exists a $C>1$
such that for any integer $1\le l\le q_{0}$, any integer $m\ge1$,
and any $\text{\textgreek{f}}\in C^{\infty}(\mathcal{M})$ solving
$\square\text{\textgreek{f}}=0$, there exists a sequence $\{\text{\textgreek{t}}_{n}\}_{n\in\mathbb{N}}$
of positive numbers satisfying $(1+C^{-1})\text{\textgreek{t}}_{n}\le\text{\textgreek{t}}_{n+1}\le(1+C)\text{\textgreek{t}}_{n}$
such that 
\begin{equation}
\begin{split}\sum_{(e_{1},\ldots e_{l-1})\in\{0,1\}^{l-1}}\mathcal{E}_{bulk}^{(2(q_{0}-l),q_{0}-l+1,m)} & [\mathcal{L}_{X^{(e_{1})}\ldots X^{(e_{l-1})}}\text{\textgreek{f}}](\text{\textgreek{t}}_{n})\lesssim_{m}\\
\lesssim_{m}\, & \text{\textgreek{t}}_{n}^{-2l}\mathcal{E}_{bound}^{(2q_{0},q_{0},m+3l\cdot k)}[\text{\textgreek{f}}](0)+\text{\textgreek{t}}_{n}^{-\text{\textgreek{d}}_{0}}\sum_{i=1}^{2l}\sum_{j=0}^{q_{0}-1}\text{\textgreek{t}}_{n}^{-2j-(2l-i)}\mathcal{E}_{en}^{(2q_{0}-i,q_{0}-j,m+3l\cdot k)}[\text{\textgreek{f}}](\text{\textgreek{t}}_{n})+\\
 & +\text{\textgreek{t}}_{n}^{-\text{\textgreek{d}}_{0}}\sum_{i=1}^{2l}\sum_{j=0}^{q_{0}-1}\text{\textgreek{t}}_{n}^{-2j-(2l-i)}\int_{\text{\textgreek{t}}_{n-2}}^{\text{\textgreek{t}}_{n+1}}\mathcal{E}_{en}^{(2q_{0}-i-1,q_{0}-j,m+3l\cdot k)}[\text{\textgreek{f}}](s)\, ds,
\end{split}
\label{eq:InductionOnNumberOfKillingDerivatives}
\end{equation}
 where $X^{(0)}=T$ and $X^{(1)}=K_{R_{c}}$.\end{lem}
\begin{proof}
The proof of (\ref{eq:InductionOnNumberOfKillingDerivatives}) will
follow by induction on $l$.

Inequality (\ref{eq:InductionOnNumberOfKillingDerivatives}) for $l=1$
follows directly from the proof of Lemma \ref{thm:SlowPointwiseDecayHighDerivativesNewMethod}
(using also Hardy-type inequalities of the form established in Lemma
\ref{lem:HardyForphiHyperboloids}) adapted to the case when $\big|\mathcal{L}_{T}^{j}g\big|_{h}\lesssim_{j}\text{\textgreek{t}}^{-\text{\textgreek{d}}_{0}}$
for all integers $j\ge0$ (see also the remark below Lemma \ref{thm:SlowPointwiseDecayHighDerivativesNewMethod}).
It thus remains to prove (\ref{eq:InductionOnNumberOfKillingDerivatives})
for any $1<l_{ind}\le q_{0}$, assuming the cases $l\le l_{ind}-1$
have been established. Without loss of generality, we will assume
that $l_{ind}=q_{0}$, since this is the hardest case. The case $l_{ind}<q_{0}$
will follow in the same way, and hence the relevant details will be
omitted.

Since we have assumed that (\ref{eq:InductionOnNumberOfKillingDerivatives})
holds for $l=q_{0}-1$ and any integer $m\ge1$, we can bound on a
sequence $\{\text{\textgreek{t}}_{n}\}_{n\in\mathbb{N}}$ satisfying
$(1+C^{-1})\text{\textgreek{t}}_{n}\le\text{\textgreek{t}}_{n+1}\le(1+C)\text{\textgreek{t}}_{n}$:
\begin{equation}
\begin{split}\sum_{(e_{1},\ldots e_{q_{0}-2})\in\{0,1\}^{q_{0}-2}} & \mathcal{E}_{bulk}^{(2,2,m)}[\mathcal{L}_{X^{(e_{1})}\ldots X^{(e_{q_{0}-2})}}\text{\textgreek{f}}](\text{\textgreek{t}}_{n})\lesssim_{m}\\
\lesssim_{m}\, & \text{\textgreek{t}}_{n}^{-2(q_{0}-1)}\mathcal{E}_{bound}^{(2q_{0},q_{0},m+3(q_{0}-1)\cdot k)}[\text{\textgreek{f}}](0)+\text{\textgreek{t}}_{n}^{-\text{\textgreek{d}}_{0}}\sum_{i=1}^{2(q_{0}-1)}\sum_{j=0}^{q_{0}-1}\text{\textgreek{t}}_{n}^{-2j-(2(q_{0}-1)-i)}\mathcal{E}_{en}^{(2q_{0}-i,q_{0}-j,m+3(q_{0}-1)\cdot k)}[\text{\textgreek{f}}](\text{\textgreek{t}}_{n})+\\
 & +\text{\textgreek{t}}_{n}^{-\text{\textgreek{d}}_{0}}\sum_{i=1}^{2(q_{0}-1)}\sum_{j=0}^{q_{0}-1}\text{\textgreek{t}}_{n}^{-2j-(2(q_{0}-1)-i)}\int_{\text{\textgreek{t}}_{n-2}}^{\text{\textgreek{t}}_{n+1}}\mathcal{E}_{en}^{(2q_{0}-i-1,q_{0}-j,m+3(q_{0}-1)\cdot k)}[\text{\textgreek{f}}](s)\, ds.
\end{split}
\label{eq:InductionHypothesisFirstStep}
\end{equation}
 Using the expression  (\ref{eq:SplittingWaveOperator}) for the wave
operator in the region $\{r\gg1\}$ (and hence morally interchanging
$\nabla_{h_{\text{\textgreek{t}},N}}^{2}+\square\rightarrow\partial_{v}T$),
we can bound (using also some Hardy type inequalities of the form
established in Lemma \ref{lem:HardyForphiHyperboloids}): 
\begin{multline}
\sum_{(e_{1},\ldots e_{q_{0}-1})\in\{0,1\}^{q_{0}-1}}\mathcal{E}_{bound}^{(2,1,m)}[\mathcal{L}_{X^{(e_{1})}\ldots X^{(e_{q_{0}-1})}}\text{\textgreek{f}}](\text{\textgreek{t}}_{n})\lesssim_{m}\sum_{(e_{1},\ldots e_{q_{0}-2})\in\{0,1\}^{q_{0}-2}}\mathcal{E}_{bulk}^{(2,2,m)}[\mathcal{L}_{X^{(e_{1})}\ldots X^{(e_{q_{0}-2})}}\text{\textgreek{f}}](\text{\textgreek{t}}_{n})+\\
+\sum_{(e_{1},\ldots e_{q_{0}-2})\in\{0,1\}^{q_{0}-2}}\sum_{j=0}^{m-1}\sum_{0\le j_{1}+j_{2}\le j}\int_{\{\bar{t}=\text{\textgreek{t}}_{n}\}\cap\{r\ge R_{c}\}}r^{2-2(q+j-2-j_{1})}\big|\nabla_{h_{\text{\textgreek{t}},N}}^{j_{1}}\big(T^{j_{2}}\square(\mathcal{L}_{X^{(e_{1})}\ldots X^{(e_{q_{0}-2})}}\text{\textgreek{f}})\big)\big|_{h_{\text{\textgreek{t}},N}}^{2}\, dh_{N}.\label{eq:AfterUsingTheWaveequation}
\end{multline}
 Thus, from (\ref{eq:InductionHypothesisFirstStep}) and (\ref{eq:AfterUsingTheWaveequation}),
we obtain in view of the expression (\ref{eq:Kcommuted}) for the
commutator of $\square$ with $\mathcal{L}_{X^{(e_{i})}}$: 
\begin{equation}
\begin{split}\sum_{(e_{1},\ldots e_{q_{0}-1})\in\{0,1\}^{q_{0}-1}} & \mathcal{E}_{bound}^{(2,1,m)}[\mathcal{L}_{X^{(e_{1})}\ldots X^{(e_{q_{0}-1})}}\text{\textgreek{f}}](\text{\textgreek{t}}_{n})\lesssim_{m}\\
\lesssim_{m}\, & \text{\textgreek{t}}_{n}^{-2(q_{0}-1)}\mathcal{E}_{bound}^{(2q_{0},q_{0},m+3(q_{0}-1)\cdot k)}[\text{\textgreek{f}}](0)+\text{\textgreek{t}}_{n}^{-\text{\textgreek{d}}_{0}}\sum_{i=1}^{2(q_{0}-1)}\sum_{j=0}^{q_{0}-1}\text{\textgreek{t}}_{n}^{-2j-(2(q_{0}-1)-i)}\mathcal{E}_{en}^{(2q_{0}-i,q_{0}-j,m+3(q_{0}-1)\cdot k)}[\text{\textgreek{f}}](\text{\textgreek{t}}_{n})+\\
 & +\text{\textgreek{t}}_{n}^{-\text{\textgreek{d}}_{0}}\sum_{i=1}^{2(q_{0}-1)}\sum_{j=0}^{q_{0}-1}\text{\textgreek{t}}_{n}^{-2j-(2(q_{0}-1)-i)}\int_{\text{\textgreek{t}}_{n-2}}^{\text{\textgreek{t}}_{n+1}}\mathcal{E}_{en}^{(2q_{0}-i-1,q_{0}-j,m+3(q_{0}-1)\cdot k)}[\text{\textgreek{f}}](s)\, ds.
\end{split}
\label{eq:InductionHypothesisTransition}
\end{equation}

Repeating the first steps of the proof of Theorem \ref{thm:FirstPointwiseDecayNewMethod}
for $\mathcal{L}_{X^{(e_{1})}\ldots X^{(e_{q_{0}-1})}}\text{\textgreek{f}}$
in place of $\text{\textgreek{f}}$, using Lemma \ref{lem:ILEDCommuted}
(and in particular the estimate (\ref{eq:ILEDCommuted})) in place
of the simple integrated local energy decay assumption \ref{eq:IntegratedLocalEnergyDecayImprovedDecay},
in view of (\ref{eq:Kcommuted}) we can bound on a sequence $\{\text{\textgreek{t}}_{n}\}_{n\in\mathbb{N}}$
with $(1+C^{-1})\text{\textgreek{t}}_{n}\le\text{\textgreek{t}}_{n+1}\le(1+C)\text{\textgreek{t}}_{n}$
(possibly different than the one appearing before): 
\begin{equation}
\begin{split}\sum_{(e_{1},\ldots e_{q_{0}-1})\in\{0,1\}^{q_{0}-1}} & \mathcal{E}_{bulk}^{(0,1,m)}[\mathcal{L}_{X^{(e_{1})}\ldots X^{(e_{q_{0}-2})}}\text{\textgreek{f}}](\text{\textgreek{t}}_{n})\lesssim_{m}\\
\lesssim_{m}\, & \text{\textgreek{t}}_{n}^{-2}\sum_{(e_{1},\ldots e_{q_{0}-1})\in\{0,1\}^{q_{0}-1}}\mathcal{E}_{bound}^{(2,1,m+2k)}[\mathcal{L}_{X^{(e_{1})}\ldots X^{(e_{q_{0}-1})}}\text{\textgreek{f}}](\text{\textgreek{t}}_{n})+\\
 & +\text{\textgreek{t}}_{n}^{-\text{\textgreek{d}}_{0}}\sum_{i=1}^{2q_{0}}\sum_{j=0}^{q_{0}-1}\text{\textgreek{t}}_{n}^{-2j-(2q_{0}-i)}\mathcal{E}_{en}^{(2q_{0}-i,q_{0}-j,m+(3q_{0}-1)\cdot k)}[\text{\textgreek{f}}](\text{\textgreek{t}}_{n})+\\
 & +\text{\textgreek{t}}_{n}^{-\text{\textgreek{d}}_{0}}\sum_{i=1}^{2q_{0}}\sum_{j=0}^{q_{0}-1}\text{\textgreek{t}}_{n}^{-2j-(2q_{0}-i)}\int_{\text{\textgreek{t}}_{n-2}}^{\text{\textgreek{t}}_{n+1}}\mathcal{E}_{en}^{(2q_{0}-i-1,q_{0}-j,m+(3q_{0}-1)\cdot k)}[\text{\textgreek{f}}](s)\, ds.
\end{split}
\label{eq:InductiveStepAlmostFinalDecay}
\end{equation}
 Using the bound (\ref{eq:InductionHypothesisTransition}) with $m+2k$
in place of $m$, we finally obtain the desired inequality: 
\begin{equation}
\begin{split}\sum_{(e_{1},\ldots e_{q_{0}-1})\in\{0,1\}^{q_{0}-1}} & \mathcal{E}_{bulk}^{(0,1,m)}[\mathcal{L}_{X^{(e_{1})}\ldots X^{(e_{q_{0}-1})}}\text{\textgreek{f}}](\text{\textgreek{t}}_{n})\lesssim_{m}\\
\lesssim_{m}\, & \text{\textgreek{t}}_{n}^{-2q_{0}}\mathcal{E}_{bound}^{(2q_{0},q_{0},m+(3q_{0}-1)\cdot k)}[\text{\textgreek{f}}](0)+\text{\textgreek{t}}_{n}^{-\text{\textgreek{d}}_{0}}\sum_{i=1}^{2q_{0}}\sum_{j=0}^{q_{0}-1}\text{\textgreek{t}}_{n}^{-2j-(2q_{0}-i)}\mathcal{E}_{en}^{(2q_{0}-i,q_{0}-j,m+(3q_{0}-1)\cdot k)}[\text{\textgreek{f}}](\text{\textgreek{t}}_{n})+\\
 & +\text{\textgreek{t}}_{n}^{-\text{\textgreek{d}}_{0}}\sum_{i=1}^{2q_{0}}\sum_{j=0}^{q_{0}-1}\text{\textgreek{t}}_{n}^{-2j-(2q_{0}-i)}\int_{\text{\textgreek{t}}_{n-2}}^{\text{\textgreek{t}}_{n+1}}\mathcal{E}_{en}^{(2q_{0}-i-1,q_{0}-j,m+(3q_{0}-1)\cdot k)}[\text{\textgreek{f}}](s)\, ds.
\end{split}
\label{eq:FinishedInductiveStep}
\end{equation}
 Thus, the proof of the Lemma is complete.\end{proof}
\begin{lem}
\label{lem:UniformDecayInductiveStep}For any integer $m\ge1$, any
$0<\text{\textgreek{e}}\ll\text{\textgreek{d}}_{0}$ small in terms
of $m$ and any smooth function $\text{\textgreek{f}}$ on $\mathcal{M}$
solving $\square\text{\textgreek{f}}=0$ we can bound: 
\begin{equation}
\begin{split}\Bigg\{\mathcal{E}_{en}^{(0,q_{0},m)}[\text{\textgreek{f}}](\text{\textgreek{t}}) & +\mathcal{E}_{bound}^{(\text{\textgreek{e}},q_{0},m)}[\text{\textgreek{f}}](\text{\textgreek{t}})+\int_{\text{\textgreek{t}}}^{+\infty}\mathcal{E}_{en}^{(-1+\text{\textgreek{e}},q_{0},m)}[\text{\textgreek{f}}](s)\, ds\Bigg\}\lesssim_{m,\text{\textgreek{e}}}\\
\lesssim_{m,\text{\textgreek{e}}}\, & \text{\textgreek{t}}^{-2q_{0}+C_{m}\cdot\text{\textgreek{e}}}\mathcal{E}_{bound}^{(2q_{0},q_{0},m+(3q_{0}+1)\cdot k)}[\text{\textgreek{f}}](0)+\text{\textgreek{t}}^{-\text{\textgreek{d}}_{0}}\sum_{i=1}^{2l}\sum_{j=0}^{q_{0}-1}\text{\textgreek{t}}^{-2j-(2q_{0}-i)}\mathcal{E}_{en}^{(2q_{0}-i+\text{\textgreek{e}},q_{0}-j,m+(3q_{0}+1)\cdot k)}[\text{\textgreek{f}}](\text{\textgreek{t}})+\\
 & +\text{\textgreek{t}}^{-\text{\textgreek{d}}_{0}}\sum_{i=1}^{2l}\sum_{j=0}^{q_{0}-1}\text{\textgreek{t}}^{-2j-(2q_{0}-i)}\int_{C_{\text{\textgreek{e}},m}^{-1}\text{\textgreek{t}}}^{+\infty}\mathcal{E}_{en}^{(2q_{0}-i-1+\text{\textgreek{e}},q_{0}-j,m+(3q_{0}+1)\cdot k)}[\text{\textgreek{f}}](s)\, ds
\end{split}
\label{eq:InductionOnNumberOfKillingDerivatives-1}
\end{equation}
 and 
\begin{align}
\mathcal{E}_{en,deg}^{(0,q_{0},m)}[\text{\textgreek{f}}](\text{\textgreek{t}})\lesssim_{m} & \text{\textgreek{t}}^{-2q_{0}}\mathcal{E}_{bound}^{(2q_{0},q_{0},m+(3q_{0}+1)\cdot k)}[\text{\textgreek{f}}](0)+\text{\textgreek{t}}^{-\text{\textgreek{d}}_{0}}\sum_{i=1}^{2l}\sum_{j=0}^{q_{0}-1}\text{\textgreek{t}}^{-2j-(2q_{0}-i)}\mathcal{E}_{en}^{(2q_{0}-i+\text{\textgreek{e}},q_{0}-j,m+(3q_{0}+1)\cdot k)}[\text{\textgreek{f}}](\text{\textgreek{t}})+\label{eq:DegenerateBoundedness-1}\\
 & +\text{\textgreek{t}}^{-\text{\textgreek{d}}_{0}}\sum_{i=1}^{2l}\sum_{j=0}^{q_{0}-1}\text{\textgreek{t}}^{-2j-(2q_{0}-i)}\int_{C_{\text{\textgreek{e}},m}^{-1}\text{\textgreek{t}}}^{+\infty}\mathcal{E}_{en}^{(2q_{0}-i-1+\text{\textgreek{e}},q_{0}-j,m+(3q_{0}+1)\cdot k)}[\text{\textgreek{f}}](s)\, ds.\nonumber 
\end{align}
\end{lem}
\begin{proof}
From Lemma (\ref{lem:InductionOnNumberOfKillingDerivatives}), we
can bound on a sequence $\{\text{\textgreek{t}}_{n}\}_{n\in\mathbb{N}}$
satisfying $(1+C^{-1})\text{\textgreek{t}}_{n}\le\text{\textgreek{t}}_{n+1}\le(1+C)\text{\textgreek{t}}_{n}$
for some large $C>0$: 
\begin{equation}
\begin{split}\sum_{(e_{1},\ldots e_{l-1})\in\{0,1\}^{q_{0}-1}}\mathcal{E}_{bulk}^{(0,1,m)} & [\mathcal{L}_{X^{(e_{1})}\ldots X^{(e_{q_{0}-1})}}\text{\textgreek{f}}](\text{\textgreek{t}}_{n})\lesssim_{m}\\
\lesssim_{m}\, & \text{\textgreek{t}}_{n}^{-2q_{0}}\mathcal{E}_{bound}^{(2q_{0},q_{0},m+3q_{0}\cdot k)}[\text{\textgreek{f}}](0)+\text{\textgreek{t}}_{n}^{-\text{\textgreek{d}}_{0}}\sum_{i=1}^{2q_{0}}\sum_{j=0}^{q_{0}-1}\text{\textgreek{t}}_{n}^{-2j-(2q_{0}-i)}\mathcal{E}_{en}^{(2q_{0}-i,q_{0}-j,m+3q_{0}\cdot k)}[\text{\textgreek{f}}](\text{\textgreek{t}}_{n})+\\
 & +\text{\textgreek{t}}_{n}^{-\text{\textgreek{d}}_{0}}\sum_{i=1}^{2q_{0}}\sum_{j=0}^{q_{0}-1}\text{\textgreek{t}}_{n}^{-2j-(2q_{0}-i)}\int_{\text{\textgreek{t}}_{n-2}}^{\text{\textgreek{t}}_{n+1}}\mathcal{E}_{en}^{(2q_{0}-i-1,q_{0}-j,m+3q_{0}\cdot k)}[\text{\textgreek{f}}](s)\, ds.
\end{split}
\label{eq:InductionOnNumberOfKillingDerivatives-2}
\end{equation}
 Using (\ref{eq:EllipticEstimatesWithoutKillingError}) successively
for $\mathcal{L}_{X^{(e_{2})}\ldots X^{(e_{q_{0}-1})}}\text{\textgreek{f}}$,
$\mathcal{L}_{X^{(e_{3})}\ldots X^{(e_{q_{0}-1})}}\text{\textgreek{f}},\ldots$,
$\text{\textgreek{f}}$ (in place of $\text{\textgreek{f}}$), making
also use of (\ref{eq:Kcommuted}) for the commutator of $\square$
with $\mathcal{L}_{X^{(e_{i})}}$, we can bound (assuming without
loss of generality that $\text{\textgreek{t}}_{1}$ is large enough
in terms of $m$): 
\begin{equation}
\mathcal{E}_{en}^{(0,q_{0},m)}[\text{\textgreek{f}}](\text{\textgreek{t}}_{n})\lesssim_{m}\sum_{(e_{1},\ldots e_{l-1})\in\{0,1\}^{l-1}}\mathcal{E}_{bulk}^{(0,1,m)}[\mathcal{L}_{X^{(e_{1})}\ldots X^{(e_{q_{0}-1})}}\text{\textgreek{f}}](\text{\textgreek{t}}_{n}).\label{eq:ControlOfTheNonDegenerateEnergy}
\end{equation}
 Therefore, from (\ref{eq:InductionOnNumberOfKillingDerivatives-2})
and (\ref{eq:ControlOfTheNonDegenerateEnergy}) we obtain: 
\begin{align}
\mathcal{E}_{en}^{(0,q_{0},m)}[\text{\textgreek{f}}](\text{\textgreek{t}}_{n})\lesssim_{m} & \text{\textgreek{t}}_{n}^{-2q_{0}}\mathcal{E}_{bound}^{(2q_{0},q_{0},m+3q_{0}\cdot k)}[\text{\textgreek{f}}](0)+\text{\textgreek{t}}_{n}^{-\text{\textgreek{d}}_{0}}\sum_{i=1}^{2q_{0}}\sum_{j=0}^{q_{0}-1}\text{\textgreek{t}}_{n}^{-2j-(2q_{0}-i)}\mathcal{E}_{en}^{(2q_{0}-i,q_{0}-j,m+3q_{0}\cdot k)}[\text{\textgreek{f}}](\text{\textgreek{t}}_{n})+\label{eq:DecayEnergyDyadicSubsequence}\\
 & +\text{\textgreek{t}}_{n}^{-\text{\textgreek{d}}_{0}}\sum_{i=1}^{2q_{0}}\sum_{j=0}^{q_{0}-1}\text{\textgreek{t}}_{n}^{-2j-(2q_{0}-i)}\int_{\text{\textgreek{t}}_{n-2}}^{\text{\textgreek{t}}_{n+1}}\mathcal{E}_{en}^{(2q_{0}-i-1,q_{0}-j,m+3q_{0}\cdot k)}[\text{\textgreek{f}}](s)\, ds.\nonumber 
\end{align}

From (\ref{eq:newMethodDu+DvPhi}) and (\ref{eq:IntegratedLocalEnergyDecayImprovedDecay}),
we also obtain after repeating the first steps of the proof of Theorem
\ref{thm:FirstPointwiseDecayNewMethod} on a sequence $\{\bar{\text{\textgreek{t}}}_{n}\}_{n\in\mathbb{N}}$
with $(1+C^{-1})\bar{\text{\textgreek{t}}}_{n}\le\bar{\text{\textgreek{t}}}_{n+1}\le(1+C)\bar{\text{\textgreek{t}}}_{n}$:
\begin{equation}
\mathcal{E}_{bulk}^{(2q_{0}-1-\text{\textgreek{e}},q_{0},m)}[\text{\textgreek{f}}](\bar{\text{\textgreek{t}}}_{n})\lesssim_{\text{\textgreek{e}},m}\bar{\text{\textgreek{t}}}_{n}^{-1}\mathcal{E}_{bound}^{(2q_{0},q_{0},m+k)}[\text{\textgreek{f}}](0).\label{eq:BeforeInterpolation}
\end{equation}
 Notice that, a priori, $\{\bar{\text{\textgreek{t}}}_{n}\}_{n\in\mathbb{N}}$
might be different than $\{\text{\textgreek{t}}_{n}\}_{n\in\mathbb{N}}$.
However, we can run the pigeonhole principle argument leading to the
choice of these sequences more carefully and arrange so that $\text{\textgreek{t}}_{n}=\bar{\text{\textgreek{t}}}_{n}$
provided $C>0$ had been fixed large enough in terms of $\text{\textgreek{e}}$,
$m$ and the geometry of $(\mathcal{M},g)$ (henceforth we will thus
assume without loss of generality that $\bar{\text{\textgreek{t}}}_{n}=\text{\textgreek{t}}_{n}$);
this follows from the following general fact: If $f_{1},f_{2}:(0,+\infty)\rightarrow(0,+\infty)$
are measurable functions satisfying 
\begin{equation}
\int_{0}^{+\infty}f_{i}(x)\, dx\le C_{i}
\end{equation}
for $i=1,2$, then there exists a sequence $\{x_{n}\}_{n\in\mathbb{N}}$
with $2x_{n}\le x_{n+1}\le4x_{n}$ such that for all $n\in\mathbb{N}$
and for $i=1,2$: 
\begin{equation}
f_{i}(x_{n})\le\log^{-1}(\frac{4}{3})\cdot\frac{C_{i}}{x_{n}}.\label{eq:DyadicCouple}
\end{equation}

\medskip{}

\noindent \emph{Proof of (\ref{eq:DyadicCouple}):} This is established
by contradiction: If there exists an interval $[a,2a]\subset(0,+\infty)$
such that the measurable sets 
\begin{equation}
I_{i}^{(a)}=\Big\{ x\in[a,2a]\,\big|\, f_{i}(x)\le\log^{-1}(\frac{4}{3})\cdot\frac{C_{i}}{x}\Big\}
\end{equation}
 for $i=1,2$ are disjoint, then there exists an $i_{0}\in\{1,2\}$
such that the complement of $I_{i_{0}}^{(a)}$ on $[a,2a]$ is at
least of half measure, i.\,e.
\begin{equation}
m\Big\{\big(I_{i_{0}}^{(a)}\big)^{c}\cap[a,2a]\Big\}\ge\frac{a}{2}.
\end{equation}
But then one obtains 
\begin{align}
C_{i_{0}} & \ge\int_{0}^{+\infty}f_{i}(x)\, dx\ge\int_{\big(I_{i_{0}}^{(a)}\big)^{c}\cap[a,2a]}f_{i_{0}}(x)\, dx>\\
 & >\int_{\big(I_{i_{0}}^{(a)}\big)^{c}\cap[a,2a]}\log^{-1}(\frac{4}{3})\cdot\frac{C_{i_{0}}}{x}\, dx\ge\log^{-1}(\frac{4}{3})\cdot C_{i_{0}}\int_{[\frac{3a}{2},2a]}\frac{1}{x}\, dx=C_{i_{0}},\nonumber 
\end{align}
(the second to last inequality following because the integral of $\frac{1}{x}$
over subsets of $[a,2a]$ of measure at least $\frac{a}{2}$ is minimized
over $[\frac{3a}{2},2a]$), which is a contradiction. Thus, there
exists an infinite sequence $\{x_{n}\}_{n\in\mathbb{N}}$ with $2x_{n}\le x_{n+1}\le4x_{n}$
on which (\ref{eq:DyadicCouple}) is satisfied.

\medskip{}

By interpolating between (\ref{eq:DecayEnergyDyadicSubsequence})
and (\ref{eq:BeforeInterpolation}) (note again that we have assumed
$\text{\textgreek{t}}_{n}=\bar{\text{\textgreek{t}}}_{n}$), we obtain:
\begin{align}
\mathcal{E}_{en}^{(\text{\textgreek{e}},q_{0},m)}[\text{\textgreek{f}}](\text{\textgreek{t}}_{n})\lesssim_{m} & \text{\textgreek{t}}_{n}^{-2q_{0}+C_{m}\text{\textgreek{e}}}\mathcal{E}_{bound}^{(2q_{0},q_{0},m+3q_{0}\cdot k)}[\text{\textgreek{f}}](0)+\text{\textgreek{t}}_{n}^{-\text{\textgreek{d}}_{0}}\sum_{i=1}^{2q_{0}}\sum_{j=0}^{q_{0}-1}\text{\textgreek{t}}_{n}^{-2j-(2q_{0}-i)}\mathcal{E}_{en}^{(2q_{0}-i,q_{0}-j,m+3q_{0}\cdot k)}[\text{\textgreek{f}}](\text{\textgreek{t}}_{n})+\label{eq:EpsilonDecayEnergyDyadicSubsequence}\\
 & +\text{\textgreek{t}}_{n}^{-\text{\textgreek{d}}_{0}}\sum_{i=1}^{2q_{0}}\sum_{j=0}^{q_{0}-1}\text{\textgreek{t}}_{n}^{-2j-(2q_{0}-i)}\int_{\text{\textgreek{t}}_{n-2}}^{\text{\textgreek{t}}_{n+1}}\mathcal{E}_{en}^{(2q_{0}-i-1,q_{0}-j,m+3q_{0}\cdot k)}[\text{\textgreek{f}}](s)\, ds.\nonumber 
\end{align}

By applying Lemma \ref{lem:ILEDCommuted} in the regions $\{\text{\textgreek{t}}_{n}\le\bar{t}\le\text{\textgreek{t}}_{n+1}\}$
and using (\ref{eq:DecayEnergyDyadicSubsequence}), (\ref{eq:EpsilonDecayEnergyDyadicSubsequence})
and the fact that $\text{\textgreek{t}}_{n}\sim_{\text{\textgreek{e}},m}\text{\textgreek{t}}_{n+1}$,
we thus obtain provided $\text{\textgreek{e}}$ has been chosen small
in terms of $m$ and $\text{\textgreek{d}}_{0}$: 
\begin{equation}
\begin{split}\Bigg\{\mathcal{E}_{en}^{(0,q_{0},m)}[\text{\textgreek{f}}](\text{\textgreek{t}}) & +\mathcal{E}_{bound}^{(\text{\textgreek{e}},q_{0},m)}[\text{\textgreek{f}}](\text{\textgreek{t}})+\int_{\text{\textgreek{t}}}^{+\infty}\mathcal{E}_{en}^{(-1+\text{\textgreek{e}},q_{0},m)}[\text{\textgreek{f}}](s)\, ds\Bigg\}\lesssim_{m,\text{\textgreek{e}}}\\
\lesssim_{m,\text{\textgreek{e}}}\, & \text{\textgreek{t}}^{-2q_{0}+C_{m}\text{\textgreek{e}}}\mathcal{E}_{bound}^{(2q_{0},q_{0},m+(3q_{0}+1)\cdot k)}[\text{\textgreek{f}}](0)+\text{\textgreek{t}}^{-\text{\textgreek{d}}_{0}}\sum_{i=1}^{2l}\sum_{j=0}^{q_{0}-1}\text{\textgreek{t}}^{-2j-(2q_{0}-i)}\mathcal{E}_{en}^{(2q_{0}-i+\text{\textgreek{e}},q_{0}-j,m+(3q_{0}+1)\cdot k)}[\text{\textgreek{f}}](\text{\textgreek{t}})+\\
 & +\text{\textgreek{t}}^{-\text{\textgreek{d}}_{0}}\sum_{i=1}^{2l}\sum_{j=0}^{q_{0}-1}\text{\textgreek{t}}^{-2j-(2q_{0}-i)}\int_{C_{\text{\textgreek{e}},m}^{-1}\text{\textgreek{t}}}^{+\infty}\mathcal{E}_{en}^{(2q_{0}-i-1+\text{\textgreek{e}},q_{0}-j,m+(3q_{0}+1)\cdot k)}[\text{\textgreek{f}}](s)\, ds
\end{split}
\label{eq:InductionOnNumberOfKillingDerivatives-1-1}
\end{equation}
 and 
\begin{align}
\mathcal{E}_{en,deg}^{(0,q_{0},m)}[\text{\textgreek{f}}](\text{\textgreek{t}})\lesssim_{m} & \text{\textgreek{t}}^{-2q_{0}}\mathcal{E}_{bound}^{(2q_{0},q_{0},m+(3q_{0}+1)\cdot k)}[\text{\textgreek{f}}](0)+\text{\textgreek{t}}^{-\text{\textgreek{d}}_{0}}\sum_{i=1}^{2l}\sum_{j=0}^{q_{0}-1}\text{\textgreek{t}}^{-2j-(2q_{0}-i)}\mathcal{E}_{en}^{(2q_{0}-i+\text{\textgreek{e}},q_{0}-j,m+(3q_{0}+1)\cdot k)}[\text{\textgreek{f}}](\text{\textgreek{t}})+\label{eq:DegenerateBoundedness}\\
 & +\text{\textgreek{t}}^{-\text{\textgreek{d}}_{0}}\sum_{i=1}^{2l}\sum_{j=0}^{q_{0}-1}\text{\textgreek{t}}^{-2j-(2q_{0}-i)}\int_{C_{\text{\textgreek{e}},m}^{-1}\text{\textgreek{t}}}^{+\infty}\mathcal{E}_{en}^{(2q_{0}-i-1+\text{\textgreek{e}},q_{0}-j,m+(3q_{0}+1)\cdot k)}[\text{\textgreek{f}}](s)\, ds.\nonumber 
\end{align}

\end{proof}
We are now ready to establish inequalities (\ref{eq:InductiveHypothesisNonDegenerate})
and (\ref{eq:InductiveHypothesisDegenerate}) of the inductive step:
\begin{lem}
\label{lem:IterationForFinalDecay} For any integer $m\ge1$, any
$0<\text{\textgreek{e}}\ll\text{\textgreek{d}}_{0}$ small in terms
of $m$ and any smooth function $\text{\textgreek{f}}$ on $\mathcal{M}$
solving $\square\text{\textgreek{f}}=0$ we can bound: 
\begin{equation}
\Bigg\{\mathcal{E}_{en}^{(0,q_{0},m)}[\text{\textgreek{f}}](\text{\textgreek{t}})+\mathcal{E}_{bound}^{(\text{\textgreek{e}},q_{0},m)}[\text{\textgreek{f}}](\text{\textgreek{t}})+\int_{\text{\textgreek{t}}}^{+\infty}\mathcal{E}_{en}^{(-1+\text{\textgreek{e}},q_{0},m)}[\text{\textgreek{f}}](s)\, ds\Bigg\}\lesssim_{m,\text{\textgreek{e}}}\text{\textgreek{t}}^{-2q_{0}+C_{m}\cdot\text{\textgreek{e}}}\mathcal{E}_{bound}^{(2q_{0},q_{0},m+\lceil\text{\textgreek{d}}_{0}^{-1}\cdot2(q_{0}-1)\rceil(3q_{0}+1)\cdot k)}[\text{\textgreek{f}}](0)\label{eq:NonDegenerateDecay}
\end{equation}
and 
\begin{equation}
\mathcal{E}_{en,deg}^{(0,q_{0},m)}[\text{\textgreek{f}}](\text{\textgreek{t}})\lesssim_{m}\text{\textgreek{t}}^{-2q_{0}}\mathcal{E}_{bound}^{(2q_{0},q_{0},m+\lceil\text{\textgreek{d}}_{0}^{-1}\cdot2(q_{0}-1)\rceil(3q_{0}+1)\cdot k)}[\text{\textgreek{f}}](0).\label{eq:DegenerateDecay}
\end{equation}
\end{lem}
\begin{proof}
From Lemma \ref{lem:InductionOnNumberOfKillingDerivatives} we can
bound: 
\begin{equation}
\begin{split}\Bigg\{\mathcal{E}_{en}^{(0,q_{0},m)}[\text{\textgreek{f}}](\text{\textgreek{t}}) & +\mathcal{E}_{bound}^{(\text{\textgreek{e}},q_{0},m)}[\text{\textgreek{f}}](\text{\textgreek{t}})+\int_{\text{\textgreek{t}}}^{+\infty}\mathcal{E}_{en}^{(-1+\text{\textgreek{e}},q_{0},m)}[\text{\textgreek{f}}](s)\, ds\Bigg\}\lesssim_{m,\text{\textgreek{e}}}\\
\lesssim_{m,\text{\textgreek{e}}}\, & \text{\textgreek{t}}^{-2q_{0}+C_{m}\cdot\text{\textgreek{e}}}\mathcal{E}_{bound}^{(2q_{0},q_{0},m+(3q_{0}+1)\cdot k)}[\text{\textgreek{f}}](0)+\text{\textgreek{t}}^{-\text{\textgreek{d}}_{0}}\sum_{i=1}^{2l}\sum_{j=0}^{q_{0}-1}\text{\textgreek{t}}^{-2j-(2l-i)}\mathcal{E}_{en}^{(2q_{0}-i+\text{\textgreek{e}},q_{0}-j,m+(3q_{0}+1)\cdot k)}[\text{\textgreek{f}}](\text{\textgreek{t}})+\\
 & +\text{\textgreek{t}}^{-\text{\textgreek{d}}_{0}}\sum_{i=1}^{2l}\sum_{j=0}^{q_{0}-1}\text{\textgreek{t}}^{-2j-(2q_{0}-i)}\int_{C_{\text{\textgreek{e}},m}^{-1}\text{\textgreek{t}}}^{+\infty}\mathcal{E}_{en}^{(2q_{0}-i-1+\text{\textgreek{e}},q_{0}-j,m+(3q_{0}+1)\cdot k)}[\text{\textgreek{f}}](s)\, ds
\end{split}
\label{eq:InductionOnNumberOfKillingDerivatives-1-2}
\end{equation}
 and 
\begin{align}
\mathcal{E}_{en,deg}^{(0,q_{0},m)}[\text{\textgreek{f}}](\text{\textgreek{t}})\lesssim_{m} & \text{\textgreek{t}}^{-2q_{0}}\mathcal{E}_{bound}^{(2q_{0},q_{0},m+(3q_{0}+1)\cdot k)}[\text{\textgreek{f}}](0)+\text{\textgreek{t}}^{-\text{\textgreek{d}}_{0}}\sum_{i=1}^{2l}\sum_{j=0}^{q_{0}-1}\text{\textgreek{t}}^{-2j-(2q_{0}-i)}\mathcal{E}_{en}^{(2q_{0}-i+\text{\textgreek{e}},q_{0}-j,m+(3q_{0}+1)\cdot k)}[\text{\textgreek{f}}](\text{\textgreek{t}})+\label{eq:DegenerateBoundedness-1-1}\\
 & +\text{\textgreek{t}}^{-\text{\textgreek{d}}_{0}}\sum_{i=1}^{2l}\sum_{j=0}^{q_{0}-1}\text{\textgreek{t}}^{-2j-(2q_{0}-i)}\int_{C_{\text{\textgreek{e}},m}^{-1}\text{\textgreek{t}}}^{+\infty}\mathcal{E}_{en}^{(2q_{0}-i-1+\text{\textgreek{e}},q_{0}-j,m+(3q_{0}+1)\cdot k)}[\text{\textgreek{f}}](s)\, ds.\nonumber 
\end{align}

Using the inductive hypothesis (i.\,e.~(\ref{eq:InductiveHypothesisHigherOrder}),
(\ref{eq:InductiveHypothesisNonDegenerate}) and (\ref{eq:InductiveHypothesisDegenerate}))
as well as Theorem \ref{thm:FirstPointwiseDecayNewMethod} and the
integrated local energy decay statement (\ref{eq:IntegratedLocalEnergyDecayImprovedDecay}),
we obtain from (\ref{eq:InductionOnNumberOfKillingDerivatives-1-2})
and (\ref{eq:DegenerateBoundedness-1-1}): 
\begin{equation}
\Bigg\{\mathcal{E}_{en}^{(0,q_{0},m)}[\text{\textgreek{f}}](\text{\textgreek{t}})+\mathcal{E}_{bound}^{(\text{\textgreek{e}},q_{0},m)}[\text{\textgreek{f}}](\text{\textgreek{t}})+\int_{\text{\textgreek{t}}}^{+\infty}\mathcal{E}_{en}^{(-1+\text{\textgreek{e}},q_{0},m)}[\text{\textgreek{f}}](s)\, ds\Bigg\}\lesssim_{m,\text{\textgreek{e}}}\text{\textgreek{t}}^{-\min\{2+\text{\textgreek{d}}_{0},2q_{0}-C_{m}\cdot\text{\textgreek{e}}\}}\mathcal{E}_{bound}^{(2q_{0},q_{0},m+(3q_{0}+1)\cdot k)}[\text{\textgreek{f}}](0)\label{eq:FirstNonDegenerateDecay}
\end{equation}
and 
\begin{equation}
\mathcal{E}_{en,deg}^{(0,q_{0},m)}[\text{\textgreek{f}}](\text{\textgreek{t}})\lesssim_{m}\text{\textgreek{t}}^{-\min\{2+\text{\textgreek{d}}_{0},2q_{0}\}}\mathcal{E}_{bound}^{(2q_{0},q_{0},m+(3q_{0}+1)\cdot k)}[\text{\textgreek{f}}](0).\label{eq:FirstDegenerateDecay}
\end{equation}
 Going back to (\ref{eq:InductionOnNumberOfKillingDerivatives-1-2})
and (\ref{eq:DegenerateBoundedness-1-1}) and using (\ref{eq:FirstNonDegenerateDecay})
(with $m+(3q_{0}+1)\cdot k)$ in place of $m$) for the error terms
in the right hand side, combined with the inductive hypothesis' inequalities
(\ref{eq:InductiveHypothesisHigherOrder}), (\ref{eq:InductiveHypothesisNonDegenerate})
for the lower order terms, we obtain the following improvement of
(\ref{eq:FirstNonDegenerateDecay}) and (\ref{eq:FirstDegenerateDecay}):
\begin{equation}
\Bigg\{\mathcal{E}_{en}^{(0,q_{0},m)}[\text{\textgreek{f}}](\text{\textgreek{t}})+\mathcal{E}_{bound}^{(\text{\textgreek{e}},q_{0},m)}[\text{\textgreek{f}}](\text{\textgreek{t}})+\int_{\text{\textgreek{t}}}^{+\infty}\mathcal{E}_{en}^{(-1+\text{\textgreek{e}},q_{0},m)}[\text{\textgreek{f}}](s)\, ds\Bigg\}\lesssim_{m,\text{\textgreek{e}}}\text{\textgreek{t}}^{-\min\{2+2\text{\textgreek{d}}_{0},2q_{0}-C_{m}\cdot\text{\textgreek{e}}\}}\mathcal{E}_{bound}^{(2q_{0},q_{0},m+2(3q_{0}+1)\cdot k)}[\text{\textgreek{f}}](0)\label{eq:FirstNonDegenerateDecay-1}
\end{equation}
and 
\begin{equation}
\mathcal{E}_{en,deg}^{(0,q_{0},m)}[\text{\textgreek{f}}](\text{\textgreek{t}})\lesssim_{m}\text{\textgreek{t}}^{-\min\{2+2\text{\textgreek{d}}_{0},2q_{0}\}}\mathcal{E}_{bound}^{(2q_{0},q_{0},m+2(3q_{0}+1)\cdot k)}[\text{\textgreek{f}}](0).\label{eq:FirstDegenerateDecay-1}
\end{equation}
 Repeating the same procedure $\lceil\text{\textgreek{d}}_{0}^{-1}\cdot2(q_{0}-1)\rceil$
times, we finally obtain the desired decay statement: 
\begin{equation}
\Bigg\{\mathcal{E}_{en}^{(0,q_{0},m)}[\text{\textgreek{f}}](\text{\textgreek{t}})+\mathcal{E}_{bound}^{(\text{\textgreek{e}},q_{0},m)}[\text{\textgreek{f}}](\text{\textgreek{t}})+\int_{\text{\textgreek{t}}}^{+\infty}\mathcal{E}_{en}^{(-1+\text{\textgreek{e}},q_{0},m)}[\text{\textgreek{f}}](s)\, ds\Bigg\}\lesssim_{m,\text{\textgreek{e}}}\text{\textgreek{t}}^{-2q_{0}+C_{m}\cdot\text{\textgreek{e}}}\mathcal{E}_{bound}^{(2q_{0},q_{0},m+\lceil\text{\textgreek{d}}_{0}^{-1}\cdot2(q_{0}-1)\rceil(3q_{0}+1)\cdot k)}[\text{\textgreek{f}}](0)\label{eq:FirstNonDegenerateDecay-1-1}
\end{equation}
and 
\begin{equation}
\mathcal{E}_{en,deg}^{(0,q_{0},m)}[\text{\textgreek{f}}](\text{\textgreek{t}})\lesssim_{m}\text{\textgreek{t}}^{-2q_{0}}\mathcal{E}_{bound}^{(2q_{0},q_{0},m+\lceil\text{\textgreek{d}}_{0}^{-1}\cdot2(q_{0}-1)\rceil(3q_{0}+1)\cdot k)}[\text{\textgreek{f}}](0).\label{eq:FirstDegenerateDecay-1-1}
\end{equation}

\end{proof}
Finally, we will establish inequality (\ref{eq:InductiveHypothesisHigherOrder})
of the inductive step:
\begin{lem}
\label{lem:ProofOfHigherRWeightedStep}For any $0\le p\le2q_{0}-2$,
any integer $m\ge1$, any $0<\text{\textgreek{e}}\ll\text{\textgreek{d}}_{0}$
small in terms of $m$ and any smooth function $\text{\textgreek{f}}$
on $\mathcal{M}$ solving $\square\text{\textgreek{f}}=0$ we can
bound: 
\begin{equation}
\mathcal{E}_{en}^{(p,q_{0},m)}[\text{\textgreek{f}}](\text{\textgreek{t}})+\int_{\text{\textgreek{t}}}^{+\infty}\mathcal{E}_{en}^{(p-1,q_{0},m)}[\text{\textgreek{f}}](s)\, ds\lesssim_{m,\text{\textgreek{e}}}\text{\textgreek{t}}^{-2q_{0}+p+C_{m}\text{\textgreek{e}}}\mathcal{E}_{bound}^{(2q_{0},q_{0},m+\lceil\text{\textgreek{d}}_{0}^{-1}\cdot2(q_{0}-1)\rceil(3q_{0}+1)\cdot k)}[\text{\textgreek{f}}](0).\label{eq:InductiveStepHigherOrder-1}
\end{equation}
\end{lem}
\begin{proof}
Inequality (\ref{eq:InductiveStepHigherOrder-1}) follows readily
after interpolating between (\ref{eq:NonDegenerateDecay}) and (\ref{eq:FinalDecayFirstenergy-1}).
\end{proof}

\subsection{\label{sub:GagliardoNirenberg}Gagliardo--Nirenberg type inequalities
on the hyperboloids $\{\bar{t}=const\}$}

In the proof of Corollary \ref{cor:ImprovedPointwiseDecay}, we need
to obtain refined pointwise control for functions $\text{\textgreek{y}}$
on the hyperboloids $\{\bar{t}=const\}$ by estimating the $L^{2}$
norms of higher order derivatives of them. To this end, we will make
use of the following Gagliardo--Nirenberg type estimates:
\begin{lem}
\label{lem:Gagliardo-NirenbergAway}For any $r_{0}\ge0$, any $\text{\textgreek{t}}>0$
and any smooth function $\text{\textgreek{y}}:\mathcal{S}_{\text{\textgreek{t}},r_{0}}\rightarrow\mathbb{C}$
(where $\mathcal{S}_{\text{\textgreek{t}},r_{0}}\doteq\{\bar{t}=\text{\textgreek{t}}\}\cap\{r\ge r_{0}\}$)
satisfying $r^{\frac{d-1}{2}}|\nabla_{h_{\text{\textgreek{t}},N}}^{l}\text{\textgreek{y}}|=O(1)$
for $0\le l\le\lceil\frac{d+1}{2}\rceil-1$ as $r\rightarrow+\infty$:

1. If $d$ is odd, we can bound:%
\footnote{Recall that $d$ is the dimension of the hypersurface $\{\bar{t}=\text{\textgreek{t}}\}$.%
} 
\begin{align}
\sup_{\mathcal{S}_{\text{\textgreek{t}},r_{0}}}|\text{\textgreek{y}}|^{2}\le C(\text{\textgreek{t}}) & \Bigg\{\Big(\int_{\mathcal{S}_{\text{\textgreek{t}},r_{0}}}|\nabla_{h_{\text{\textgreek{t}},N}}^{\frac{d-1}{2}}\text{\textgreek{y}}|_{h_{\text{\textgreek{t}},N}}^{2}\, dh_{N}\Big)^{\frac{1}{2}}\Big(\int_{\mathcal{S}_{\text{\textgreek{t}},r_{0}}}|\nabla_{h_{\text{\textgreek{t}},N}}^{\frac{d+1}{2}}\text{\textgreek{y}}|_{h_{\text{\textgreek{t}},N}}^{2}\, dh_{N}\Big)^{\frac{1}{2}}+\int_{\mathcal{S}_{\text{\textgreek{t}},r_{0}}}|\nabla_{h_{\text{\textgreek{t}},N}}^{\frac{d+1}{2}}\text{\textgreek{y}}|_{h_{\text{\textgreek{t}},N}}^{2}\, dh_{N}\Bigg\}.\label{eq:GagliardoNirenbergOddfinal-ForHigherDerivatives}
\end{align}

and
\begin{align}
\sup_{\mathcal{S}_{\text{\textgreek{t}},r_{0}}}|\text{\textgreek{y}}|^{2}\le C(\text{\textgreek{t}}) & \Bigg\{\Big(\int_{\mathcal{S}_{\text{\textgreek{t}},r_{0}}}|\nabla_{h_{\text{\textgreek{t}},N}}\text{\textgreek{y}}|_{h_{\text{\textgreek{t}},N}}^{2}\, dh_{N}\Big)^{\frac{1}{d-1}}\Big(\int_{\mathcal{S}_{\text{\textgreek{t}},r_{0}}}|\nabla_{h_{\text{\textgreek{t}},N}}^{\frac{d+1}{2}}\text{\textgreek{y}}|_{h_{\text{\textgreek{t}},N}}^{2}\, dh_{N}\Big)^{\frac{d-2}{d-1}}+\int_{\mathcal{S}_{\text{\textgreek{t}},r_{0}}}|\nabla_{h_{\text{\textgreek{t}},N}}^{\frac{d+1}{2}}\text{\textgreek{y}}|_{h_{\text{\textgreek{t}},N}}^{2}\, dh_{N}\Bigg\}\label{eq:GagliardoNirenbergOddfinal-1}
\end{align}

2. If $d$ is even, we can bound for any $\text{\textgreek{e}}>0$:
\begin{align}
\sup_{\mathcal{S}_{\text{\textgreek{t}},r_{0}}}|\text{\textgreek{y}}|^{2}\le C_{\text{\textgreek{e}}}(\text{\textgreek{t}}) & \Bigg\{\Big(\int_{\mathcal{S}_{\text{\textgreek{t}},r_{0}}}|\nabla_{h_{\text{\textgreek{t}},N}}\text{\textgreek{y}}|_{h_{\text{\textgreek{t}},N}}^{2}\, dh_{N}\Big)^{\frac{\text{\textgreek{e}}}{d-2}}\Big(\int_{\mathcal{S}_{\text{\textgreek{t}},r_{0}}}\big(|\nabla_{h_{\text{\textgreek{t}},N}}^{\frac{d}{2}+1}\text{\textgreek{y}}|_{h_{\text{\textgreek{t}},N}}^{2}+|\nabla_{h_{\text{\textgreek{t}},N}}^{\frac{d}{2}}\text{\textgreek{y}}|_{h_{\text{\textgreek{t}},N}}^{2}\big)\, dh_{N}\Big)^{1-\frac{\text{\textgreek{e}}\cdot(2d-2)}{2(d-2)}}+\label{eq:GagliardoNirenbergEvenFinal-1}\\
 & +\int_{\mathcal{S}_{\text{\textgreek{t}},r_{0}}}\big(|\nabla_{h_{\text{\textgreek{t}},N}}^{\frac{d}{2}+1}\text{\textgreek{y}}|_{h_{\text{\textgreek{t}},N}}^{2}+|\nabla_{h_{\text{\textgreek{t}},N}}^{\frac{d}{2}}\text{\textgreek{y}}|_{h_{\text{\textgreek{t}},N}}^{2}\big)\, dh_{N}+\int_{\mathcal{S}_{\text{\textgreek{t}},r_{0}}}r_{+}^{\text{\textgreek{e}}}|L^{\frac{d}{2}}\text{\textgreek{y}}|^{2}\, dh_{N}\Bigg\}.\nonumber 
\end{align}
 In the above, $L$ is a vector field which is identically $0$ in
the region $\{r\lesssim1\}$ and equals the coordinate vector field
$\partial_{v}$ of the coordinate system $(v,\text{\textgreek{sv}})$
on $\{\bar{t}=\text{\textgreek{t}}\}\cap\{r\gg1\}$.\end{lem}
\begin{proof}
Let us fix an $R_{0}>0$ large in terms of the geometry of $(\mathcal{M},g)$.

On $\mathbb{R}^{d}$, the following Gagliardo--Nirenberg type inequalities
hold for $f\in C_{0}^{\infty}(\mathbb{R}^{d})$ (see \cite{Nirenberg2011}
and \cite{Stein1971}):

1. If $d$ is odd: 
\begin{equation}
||f||_{L^{\infty}}\le C||f||_{\dot{H}{}^{\frac{d-1}{2}}(\mathbb{R}^{d})}^{\frac{1}{2}}\cdot||f||_{\dot{H}{}^{\frac{d+1}{2}}(\mathbb{R}^{d})}^{\frac{1}{2}}.\label{eq:GagliardoNirenbergFlatOdd-Higher}
\end{equation}

and 
\begin{equation}
||f||_{L^{\infty}}\le C||f||_{\dot{H}{}^{1}(\mathbb{R}^{d})}^{\frac{1}{d-1}}\cdot||f||_{\dot{H}{}^{\frac{d+1}{2}}(\mathbb{R}^{d})}^{\frac{d-2}{d-1}}\label{eq:GagliardoNirenbergFlatOdd}
\end{equation}

2. If $d$ is even: 
\begin{equation}
||f||_{L^{\infty}}\le C_{\text{\textgreek{e}}}||f||_{\dot{H}{}^{\frac{d}{2}-\text{\textgreek{e}}}(\mathbb{R}^{d})}^{\frac{1}{2}}\cdot||f||_{\dot{H}{}^{\frac{d}{2}+\text{\textgreek{e}}}(\mathbb{R}^{d})}^{\frac{1}{2}}.\label{eq:FirstGagliardoNirenbergFlatEven}
\end{equation}

\noindent In the above, we have used the homogeneous norms $\dot{H}{}^{a}(\mathbb{R}^{d})$
defined with the use of the Fourier transform as 
\begin{equation}
||f||_{\dot{H}{}^{a}(\mathbb{R}^{d})}^{2}\doteq\int_{\mathbb{R}^{d}}|\text{\textgreek{x}}|^{2a}\big|\hat{f}\big|^{2}\, d\text{\textgreek{x}},
\end{equation}
 where $\hat{f}$ is the Fourier transform of $f$. 

Using a simple interpolation argument, we can also bound 
\begin{equation}
||f||_{\dot{H}{}^{\frac{d}{2}-\text{\textgreek{e}}}(\mathbb{R}^{d})}\le C_{\text{\textgreek{e}}}||f||_{\dot{H}{}^{1}(\mathbb{R}^{d})}^{\frac{2\text{\textgreek{e}}}{d-2}}\cdot||f||_{\dot{H}{}^{\frac{d}{2}}(\mathbb{R}^{d})}^{\frac{d-2-2\text{\textgreek{e}}}{d-2}}
\end{equation}
 and 
\begin{equation}
||f||_{\dot{H}{}^{\frac{d}{2}+\text{\textgreek{e}}}(\mathbb{R}^{d})}\le C_{\text{\textgreek{e}}}||f||_{\dot{H}{}^{\frac{d}{2}}(\mathbb{R}^{d})}^{1-\text{\textgreek{e}}}\cdot||f||_{\dot{H}{}^{\frac{d}{2}+1}(\mathbb{R}^{d})}^{\text{\textgreek{e}}},
\end{equation}
 and thus we obtain from (\ref{eq:FirstGagliardoNirenbergFlatEven})
in the case $d$ is even: 
\begin{equation}
||f||_{L^{\infty}}\le C_{\text{\textgreek{e}}}||f||_{\dot{H}{}^{1}(\mathbb{R}_{+}^{d})}^{\frac{\text{\textgreek{e}}}{d-2}}\cdot||f||_{\dot{H}{}^{\frac{d}{2}}(\mathbb{R}_{+}^{d})}^{1-\frac{\text{\textgreek{e}}\cdot d}{2(d-2)}}\cdot||f||_{\dot{H}{}^{\frac{d}{2}+1}(\mathbb{R}_{+}^{d})}^{\frac{\text{\textgreek{e}}}{2}}.\label{eq:GagliardoNirenbergFlatEven}
\end{equation}

Let $\text{\textgreek{q}}_{R_{0}}:\{\bar{t}=\text{\textgreek{t}}\}\rightarrow[0,1]$
be a smooth cut-off function such that $\text{\textgreek{q}}_{R_{0}}$
is a function of $r$ satisfying $\text{\textgreek{q}}_{R_{0}}\equiv0$
on $\{r\le R_{0}\}$ and $\text{\textgreek{q}}_{R_{0}}\equiv1$ on
$\{r\ge2R_{0}\}$. Since $\text{\textgreek{q}}_{R_{0}}\text{\textgreek{y}}$
is supported on $\{r\ge R_{0}\}$, by pulling back through the diffeomorphism
$(r,\text{\textgreek{sv}})$ the operator $\nabla_{\mathbb{R}^{d}}$
on $\{\bar{t}=\text{\textgreek{t}}\}$ and using the bound 
\begin{equation}
\big|\nabla_{\mathbb{R}^{d}}^{l}\text{\textgreek{f}}\big|_{h_{\text{\textgreek{t}},N}}^{2}\le C_{l}\sum_{j=1}^{l}r^{-2(l-j)}\big|\nabla_{h_{\text{\textgreek{t}},N}}\text{\textgreek{f}}\big|_{h_{\text{\textgreek{t}},N}}^{2},
\end{equation}
 we obtain from (\ref{eq:GagliardoNirenbergFlatOdd-Higher}), (\ref{eq:GagliardoNirenbergFlatOdd})
and (\ref{eq:GagliardoNirenbergFlatEven}):

1. In case $d$ is odd: 
\begin{equation}
\sup_{\mathcal{S}_{\text{\textgreek{t}},r_{0}}}|\text{\textgreek{q}}_{R_{0}}\text{\textgreek{y}}|^{2}\le C\cdot\Big(\sum_{j=0}^{\frac{d-3}{2}}\int_{\mathcal{S}_{\text{\textgreek{t}},r_{0}}}|\nabla_{h_{\text{\textgreek{t}},N}}^{\frac{d-1}{2}-j}(\text{\textgreek{q}}_{R_{0}}\text{\textgreek{y}})|_{h_{\text{\textgreek{t}},N}}^{2}\, dh_{N}\Big)^{\frac{1}{2}}\Big(\sum_{j=0}^{\frac{d-1}{2}}\int_{\mathcal{S}_{\text{\textgreek{t}},r_{0}}}r^{-2j}|\nabla_{h_{\text{\textgreek{t}},N}}^{\frac{d+1}{2}-j}(\text{\textgreek{q}}_{R_{0}}\text{\textgreek{y}})|_{h_{\text{\textgreek{t}},N}}^{2}\, dh_{N}\Big)^{\frac{1}{2}}.\label{eq:GagliardoNirenbergAsymptoticOdd-2}
\end{equation}

and 
\begin{equation}
\sup_{\mathcal{S}_{\text{\textgreek{t}},r_{0}}}|\text{\textgreek{q}}_{R_{0}}\text{\textgreek{y}}|^{2}\le C\cdot\Big(\int_{\mathcal{S}_{\text{\textgreek{t}},r_{0}}}|\nabla_{h_{\text{\textgreek{t}},N}}(\text{\textgreek{q}}_{R_{0}}\text{\textgreek{y}})|_{h_{\text{\textgreek{t}},N}}^{2}\, dh_{N}\Big)^{\frac{1}{d-1}}\Big(\sum_{j=0}^{\frac{d-1}{2}}\int_{\mathcal{S}_{\text{\textgreek{t}},r_{0}}}r^{-2j}|\nabla_{h_{\text{\textgreek{t}},N}}^{\frac{d+1}{2}-j}(\text{\textgreek{q}}_{R_{0}}\text{\textgreek{y}})|_{h_{\text{\textgreek{t}},N}}^{2}\, dh_{N}\Big)^{\frac{d-2}{d-1}}\label{eq:GagliardoNirenbergAsymptoticOdd}
\end{equation}

2. In case $d$ is even: 
\begin{align}
\sup_{\mathcal{S}_{\text{\textgreek{t}},r_{0}}}|\text{\textgreek{q}}_{R_{0}}\text{\textgreek{y}}|^{2}\le & C_{\text{\textgreek{e}}}\cdot\Big(\int_{\mathcal{S}_{\text{\textgreek{t}},r_{0}}}|\nabla_{h_{\text{\textgreek{t}},N}}(\text{\textgreek{q}}_{R_{0}}\text{\textgreek{y}})|_{h_{\text{\textgreek{t}},N}}^{2}\, dh_{N}\Big)^{\frac{\text{\textgreek{e}}}{d-2}}\Big(\sum_{j=0}^{\frac{d-2}{2}}\int_{\mathcal{S}_{\text{\textgreek{t}},r_{0}}}r^{-2j}|\nabla_{h_{\text{\textgreek{t}},N}}^{\frac{d}{2}-j}(\text{\textgreek{q}}_{R_{0}}\text{\textgreek{y}})|_{h_{\text{\textgreek{t}},N}}^{2}\, dh_{N}\Big)^{1-\frac{\text{\textgreek{e}}\cdot d}{2(d-2)}}\cdot\label{eq:GagliardoNirenbergAsymptoticEven}\\
 & \hphantom{C_{\text{\textgreek{e}}}}\cdot\Big(\sum_{j=0}^{\frac{d}{2}}\int_{\mathcal{S}_{\text{\textgreek{t}},r_{0}}}r^{-2j}|\nabla_{h_{\text{\textgreek{t}},N}}^{\frac{d}{2}+1-j}(\text{\textgreek{q}}_{R_{0}}\text{\textgreek{y}})|_{h_{\text{\textgreek{t}},N}}^{2}\, dh_{N}\Big)^{\frac{\text{\textgreek{e}}}{2}}.\nonumber 
\end{align}
 By applying the Leibnitz rule and using Hardy type inequalities of
the form established in Lemma \ref{lem:HardyForphiHyperboloids},
we obtain from (\ref{eq:GagliardoNirenbergAsymptoticOdd-2}), (\ref{eq:GagliardoNirenbergAsymptoticOdd})
and (\ref{eq:GagliardoNirenbergAsymptoticEven}):

1. In case $d$ is odd: 
\begin{equation}
\sup_{\mathcal{S}_{\text{\textgreek{t}},r_{0}}}|\text{\textgreek{q}}_{R_{0}}\text{\textgreek{y}}|^{2}\le C\cdot\Bigg\{\Big(\int_{\mathcal{S}_{\text{\textgreek{t}},r_{0}}}|\nabla_{h_{\text{\textgreek{t}},N}}^{\frac{d-1}{2}}\text{\textgreek{y}}|_{h_{\text{\textgreek{t}},N}}^{2}\, dh_{N}\Big)^{\frac{1}{2}}\Big(\int_{\mathcal{S}_{\text{\textgreek{t}},r_{0}}}|\nabla_{h_{\text{\textgreek{t}},N}}^{\frac{d+1}{2}}\text{\textgreek{y}}|_{h_{\text{\textgreek{t}},N}}^{2}\, dh_{N}\Big)^{\frac{1}{2}}+\int_{\mathcal{S}_{\text{\textgreek{t}},r_{0}}\cap\{R_{0}\le r\le2R_{0}\}}|\text{\textgreek{y}}|^{2}\, dh_{N}\Bigg\}.\label{eq:GagliardoNirenbergAsymptoticOdd-1-2}
\end{equation}

and 
\begin{equation}
\sup_{\mathcal{S}_{\text{\textgreek{t}},r_{0}}}|\text{\textgreek{q}}_{R_{0}}\text{\textgreek{y}}|^{2}\le C\cdot\Bigg\{\Big(\int_{\mathcal{S}_{\text{\textgreek{t}},r_{0}}}|\nabla_{h_{\text{\textgreek{t}},N}}\text{\textgreek{y}}|_{h_{\text{\textgreek{t}},N}}^{2}\, dh_{N}\Big)^{\frac{1}{d-1}}\Big(\int_{\mathcal{S}_{\text{\textgreek{t}},r_{0}}}|\nabla_{h_{\text{\textgreek{t}},N}}^{\frac{d+1}{2}}\text{\textgreek{y}}|_{h_{\text{\textgreek{t}},N}}^{2}\, dh_{N}\Big)^{\frac{d-2}{d-1}}+\int_{\mathcal{S}_{\text{\textgreek{t}},r_{0}}\cap\{R_{0}\le r\le2R_{0}\}}|\text{\textgreek{y}}|^{2}\, dh_{N}\Bigg\}\label{eq:GagliardoNirenbergAsymptoticOdd-1}
\end{equation}

2. In case $d$ is even: 
\begin{equation}
\sup_{\mathcal{S}_{\text{\textgreek{t}},r_{0}}}|\text{\textgreek{q}}_{R_{0}}\text{\textgreek{y}}|^{2}\le C_{\text{\textgreek{e}}}\cdot\Bigg\{\Big(\int_{\mathcal{S}_{\text{\textgreek{t}},r_{0}}}|\nabla_{h_{\text{\textgreek{t}},N}}\text{\textgreek{y}}|_{h_{\text{\textgreek{t}},N}}^{2}\, dh_{N}\Big)^{\frac{\text{\textgreek{e}}}{d-2}}\Big(\int_{\mathcal{S}_{\text{\textgreek{t}},r_{0}}}\big(|\nabla_{h_{\text{\textgreek{t}},N}}^{\frac{d}{2}+1}\text{\textgreek{y}}|_{h_{\text{\textgreek{t}},N}}^{2}+|\nabla_{h_{\text{\textgreek{t}},N}}^{\frac{d}{2}}\text{\textgreek{y}}|_{h_{\text{\textgreek{t}},N}}^{2}\big)\, dh_{N}\Big)^{1-\frac{\text{\textgreek{e}}\cdot(2d-2)}{2(d-2)}}+\int_{\mathcal{S}_{\text{\textgreek{t}},r_{0}}\cap\{R_{0}\le r\le2R_{0}\}}|\text{\textgreek{y}}|^{2}\, dh_{N}\Bigg\}.\label{eq:GagliardoNirenbergAsymptoticEven-1}
\end{equation}
 Since $(1-\text{\textgreek{q}}_{R_{0}})$ is supported in the compact
region $\{\bar{t}=\text{\textgreek{t}}\}\cap\{r\le2R_{0}\}$ and $R_{0}$
is fixed in terms of the geometry of $\{\bar{t}=\text{\textgreek{t}}\}$,
using the Sobolev inequality (see e.\,g.~\cite{Hebey1999}) 
\begin{equation}
\sup_{\mathcal{S}_{\text{\textgreek{t}},r_{0}}}|(1-\text{\textgreek{q}}_{R_{0}})\text{\textgreek{y}}|^{2}\le C(\text{\textgreek{t}})\sum_{j=0}^{\lceil\frac{d+1}{2}\rceil}\int_{\mathcal{S}_{\text{\textgreek{t}},r_{0}}\cap\{r\le2R_{0}\}}\big|\nabla_{h_{\text{\textgreek{t}},N}}^{j}\text{\textgreek{y}}\big|_{h_{\text{\textgreek{t}},N}}^{2}\, dh_{N}
\end{equation}
 and the Poincare inequality 
\begin{equation}
\int_{\mathcal{S}_{\text{\textgreek{t}},r_{0}}\cap\{r\le2R_{0}\}}|\text{\textgreek{y}}|^{2}\, dh_{N}\le C(\text{\textgreek{t}})\int_{\mathcal{S}_{\text{\textgreek{t}},r_{0}}\cap\{r\le2R_{0}\}}\big|\nabla_{h_{\text{\textgreek{t}},N}}\text{\textgreek{y}}\big|_{h_{\text{\textgreek{t}},N}}^{2}\, dh_{N}+\int_{\mathcal{S}_{\text{\textgreek{t}},r_{0}}\cap\{R_{0}\le r\le2R_{0}\}}|\text{\textgreek{y}}|^{2}\, dh_{N},
\end{equation}
 we obtain from (\ref{eq:GagliardoNirenbergAsymptoticOdd-1-2}), (\ref{eq:GagliardoNirenbergAsymptoticOdd-1})
and (\ref{eq:GagliardoNirenbergAsymptoticEven-1}): 

1. In case $d$ is odd: 
\begin{align}
\sup_{\mathcal{S}_{\text{\textgreek{t}},r_{0}}}|\text{\textgreek{y}}|^{2}\le C(\text{\textgreek{t}}) & \Bigg\{\Big(\int_{\mathcal{S}_{\text{\textgreek{t}},r_{0}}}|\nabla_{h_{\text{\textgreek{t}},N}}^{\frac{d-1}{2}}\text{\textgreek{y}}|_{h_{\text{\textgreek{t}},N}}^{2}\, dh_{N}\Big)^{\frac{1}{2}}\Big(\int_{\mathcal{S}_{\text{\textgreek{t}},r_{0}}}|\nabla_{h_{\text{\textgreek{t}},N}}^{\frac{d+1}{2}}\text{\textgreek{y}}|_{h_{\text{\textgreek{t}},N}}^{2}\, dh_{N}\Big)^{\frac{1}{2}}+\label{eq:GagliardoNirenbergAsymptoticOdd-1-1-1}\\
 & +\int_{\mathcal{S}_{\text{\textgreek{t}},r_{0}}}|\nabla_{h_{\text{\textgreek{t}},N}}^{\frac{d+1}{2}}\text{\textgreek{y}}|_{h_{\text{\textgreek{t}},N}}^{2}\, dh_{N}+\int_{\mathcal{S}_{\text{\textgreek{t}},r_{0}}\cap\{R_{0}\le r\le2R_{0}\}}|\text{\textgreek{y}|}^{2}\, dh_{N}\Bigg\}.\nonumber 
\end{align}

and 
\begin{align}
\sup_{\mathcal{S}_{\text{\textgreek{t}},r_{0}}}|\text{\textgreek{y}}|^{2}\le C(\text{\textgreek{t}}) & \Bigg\{\Big(\int_{\mathcal{S}_{\text{\textgreek{t}},r_{0}}}|\nabla_{h_{\text{\textgreek{t}},N}}\text{\textgreek{y}}|_{h_{\text{\textgreek{t}},N}}^{2}\, dh_{N}\Big)^{\frac{1}{d-1}}\Big(\int_{\mathcal{S}_{\text{\textgreek{t}},r_{0}}}|\nabla_{h_{\text{\textgreek{t}},N}}^{\frac{d+1}{2}}\text{\textgreek{y}}|_{h_{\text{\textgreek{t}},N}}^{2}\, dh_{N}\Big)^{\frac{d-2}{d-1}}+\label{eq:GagliardoNirenbergAsymptoticOdd-1-1}\\
 & +\int_{\mathcal{S}_{\text{\textgreek{t}},r_{0}}}|\nabla_{h_{\text{\textgreek{t}},N}}^{\frac{d+1}{2}}\text{\textgreek{y}}|_{h_{\text{\textgreek{t}},N}}^{2}\, dh_{N}+\int_{\mathcal{S}_{\text{\textgreek{t}},r_{0}}\cap\{R_{0}\le r\le2R_{0}\}}|\text{\textgreek{y}|}^{2}\, dh_{N}\Bigg\}\nonumber 
\end{align}

2. In case $d$ is even: 
\begin{align}
\sup_{\mathcal{S}_{\text{\textgreek{t}},r_{0}}}|\text{\textgreek{y}}|^{2}\le C_{\text{\textgreek{e}}}(\text{\textgreek{t}}) & \Bigg\{\Big(\int_{\mathcal{S}_{\text{\textgreek{t}},r_{0}}}|\nabla_{h_{\text{\textgreek{t}},N}}\text{\textgreek{y}}|_{h_{\text{\textgreek{t}},N}}^{2}\, dh_{N}\Big)^{\frac{\text{\textgreek{e}}}{d-2}}\Big(\int_{\mathcal{S}_{\text{\textgreek{t}},r_{0}}}\big(|\nabla_{h_{\text{\textgreek{t}},N}}^{\frac{d}{2}+1}\text{\textgreek{y}}|_{h_{\text{\textgreek{t}},N}}^{2}+|\nabla_{h_{\text{\textgreek{t}},N}}^{\frac{d}{2}}\text{\textgreek{y}}|_{h_{\text{\textgreek{t}},N}}^{2}\big)\, dh_{N}\Big)^{1-\frac{\text{\textgreek{e}}\cdot(2d-2)}{2(d-2)}}+\label{eq:GagliardoNirenbergAsymptoticEven-1-1}\\
 & +\int_{\mathcal{S}_{\text{\textgreek{t}},r_{0}}}\big(|\nabla_{h_{\text{\textgreek{t}},N}}^{\frac{d}{2}+1}\text{\textgreek{y}}|_{h_{\text{\textgreek{t}},N}}^{2}+|\nabla_{h_{\text{\textgreek{t}},N}}^{\frac{d}{2}}\text{\textgreek{y}}|_{h_{\text{\textgreek{t}},N}}^{2}\big)\, dh_{N}+\int_{\mathcal{S}_{\text{\textgreek{t}},r_{0}}\cap\{R_{0}\le r\le2R_{0}\}}|\text{\textgreek{y}}|^{2}\, dh_{N}\Bigg\}.\nonumber 
\end{align}

Finally, in case $d$ is odd, using Lemma \ref{lem:HardyCritical}
(as well as Lemma \ref{lem:HardyForphiHyperboloids}) we can bound
\begin{equation}
\int_{\mathcal{S}_{\text{\textgreek{t}},r_{0}}\cap\{R_{0}\le r\le2R_{0}\}}|\text{\textgreek{y}|}^{2}\, dh_{N}\le C\cdot\Big(\int_{\mathcal{S}_{\text{\textgreek{t}},r_{0}}}|\nabla_{h_{\text{\textgreek{t}},N}}\text{\textgreek{y}}|_{h_{\text{\textgreek{t}},N}}^{2}\, dh_{N}\Big)^{\frac{1}{2}}\Big(\int_{\mathcal{S}_{\text{\textgreek{t}},r_{0}}}|\nabla_{h_{\text{\textgreek{t}},N}}^{\frac{d+1}{2}}\text{\textgreek{y}}|_{h_{\text{\textgreek{t}},N}}^{2}\, dh_{N}\Big)^{\frac{1}{2}},
\end{equation}
 and 
\begin{equation}
\int_{\mathcal{S}_{\text{\textgreek{t}},r_{0}}\cap\{R_{0}\le r\le2R_{0}\}}|\text{\textgreek{y}|}^{2}\, dh_{N}\le C\cdot\Big(\int_{\mathcal{S}_{\text{\textgreek{t}},r_{0}}}|\nabla_{h_{\text{\textgreek{t}},N}}\text{\textgreek{y}}|_{h_{\text{\textgreek{t}},N}}^{2}\, dh_{N}\Big)^{\frac{1}{d-1}}\Big(\int_{\mathcal{S}_{\text{\textgreek{t}},r_{0}}}|\nabla_{h_{\text{\textgreek{t}},N}}^{\frac{d+1}{2}}\text{\textgreek{y}}|_{h_{\text{\textgreek{t}},N}}^{2}\, dh_{N}\Big)^{\frac{d-2}{d-1}}
\end{equation}
and thus (\ref{eq:GagliardoNirenbergAsymptoticOdd-1-1}) and (\ref{eq:GagliardoNirenbergAsymptoticOdd-1-1-1})
become 
\begin{align}
\sup_{\mathcal{S}_{\text{\textgreek{t}},r_{0}}}|\text{\textgreek{y}}|^{2}\le C(\text{\textgreek{t}}) & \Bigg\{\Big(\int_{\mathcal{S}_{\text{\textgreek{t}},r_{0}}}|\nabla_{h_{\text{\textgreek{t}},N}}^{\frac{d-1}{2}}\text{\textgreek{y}}|_{h_{\text{\textgreek{t}},N}}^{2}\, dh_{N}\Big)^{\frac{1}{2}}\Big(\int_{\mathcal{S}_{\text{\textgreek{t}},r_{0}}}|\nabla_{h_{\text{\textgreek{t}},N}}^{\frac{d+1}{2}}\text{\textgreek{y}}|_{h_{\text{\textgreek{t}},N}}^{2}\, dh_{N}\Big)^{\frac{1}{2}}+\int_{\mathcal{S}_{\text{\textgreek{t}},r_{0}}}|\nabla_{h_{\text{\textgreek{t}},N}}^{\frac{d+1}{2}}\text{\textgreek{y}}|_{h_{\text{\textgreek{t}},N}}^{2}\, dh_{N}\Bigg\}.\label{eq:GagliardoNirenbergOddfinal-2}
\end{align}
and 
\begin{align}
\sup_{\mathcal{S}_{\text{\textgreek{t}},r_{0}}}|\text{\textgreek{y}}|^{2}\le C(\text{\textgreek{t}}) & \Bigg\{\Big(\int_{\mathcal{S}_{\text{\textgreek{t}},r_{0}}}|\nabla_{h_{\text{\textgreek{t}},N}}\text{\textgreek{y}}|_{h_{\text{\textgreek{t}},N}}^{2}\, dh_{N}\Big)^{\frac{1}{d-1}}\Big(\int_{\mathcal{S}_{\text{\textgreek{t}},r_{0}}}|\nabla_{h_{\text{\textgreek{t}},N}}^{\frac{d+1}{2}}\text{\textgreek{y}}|_{h_{\text{\textgreek{t}},N}}^{2}\, dh_{N}\Big)^{\frac{d-2}{d-1}}+\int_{\mathcal{S}_{\text{\textgreek{t}},r_{0}}}|\nabla_{h_{\text{\textgreek{t}},N}}^{\frac{d+1}{2}}\text{\textgreek{y}}|_{h_{\text{\textgreek{t}},N}}^{2}\, dh_{N}\Bigg\}\label{eq:GagliardoNirenbergOddfinal}
\end{align}

In case $d$ is even, on the other hand, from Lemma \ref{lem:HardyForphiHyperboloids}
we can readily bound 
\begin{equation}
\int_{\mathcal{S}_{\text{\textgreek{t}},r_{0}}\cap\{R_{0}\le r\le2R_{0}\}}|\text{\textgreek{y}|}^{2}\, dh_{N}\le C_{\text{\textgreek{e}}}\int_{\mathcal{S}_{\text{\textgreek{t}},r_{0}}}r_{+}^{\text{\textgreek{e}}}|L^{\frac{d}{2}}\text{\textgreek{y}}|^{2}\, dh_{N},
\end{equation}
 and hence (\ref{eq:GagliardoNirenbergAsymptoticEven-1-1}) becomes:
\begin{align}
\sup_{\mathcal{S}_{\text{\textgreek{t}},r_{0}}}|\text{\textgreek{y}}|^{2}\le C_{\text{\textgreek{e}}}(\text{\textgreek{t}}) & \Bigg\{\Big(\int_{\mathcal{S}_{\text{\textgreek{t}},r_{0}}}|\nabla_{h_{\text{\textgreek{t}},N}}\text{\textgreek{y}}|_{h_{\text{\textgreek{t}},N}}^{2}\, dh_{N}\Big)^{\frac{\text{\textgreek{e}}}{d-2}}\Big(\int_{\mathcal{S}_{\text{\textgreek{t}},r_{0}}}\big(|\nabla_{h_{\text{\textgreek{t}},N}}^{\frac{d}{2}+1}\text{\textgreek{y}}|_{h_{\text{\textgreek{t}},N}}^{2}+|\nabla_{h_{\text{\textgreek{t}},N}}^{\frac{d}{2}}\text{\textgreek{y}}|_{h_{\text{\textgreek{t}},N}}^{2}\big)\, dh_{N}\Big)^{1-\frac{\text{\textgreek{e}}\cdot(2d-2)}{2(d-2)}}+\label{eq:GagliardoNirenbergAsymptoticEven-1-1-1}\\
 & +\int_{\mathcal{S}_{\text{\textgreek{t}},r_{0}}}\big(|\nabla_{h_{\text{\textgreek{t}},N}}^{\frac{d}{2}+1}\text{\textgreek{y}}|_{h_{\text{\textgreek{t}},N}}^{2}+|\nabla_{h_{\text{\textgreek{t}},N}}^{\frac{d}{2}}\text{\textgreek{y}}|_{h_{\text{\textgreek{t}},N}}^{2}\big)\, dh_{N}+\int_{\mathcal{S}_{\text{\textgreek{t}},r_{0}}}r_{+}^{\text{\textgreek{e}}}|L^{\frac{d}{2}}\text{\textgreek{y}}|^{2}\, dh_{N}\Bigg\}.\nonumber 
\end{align}
\end{proof}
\begin{lem}
\label{lem:Gagliardo-NirenbergDegenerate}For any $\text{\textgreek{t}}>0$
and any smooth function $\text{\textgreek{y}}:\{\bar{t}=\text{\textgreek{t}}\}\cap\{r\le1\}\rightarrow\mathbb{C}$
we can bound 
\begin{align}
\sup_{\{\bar{t}=\text{\textgreek{t}}\}\cap\{r\le\frac{1}{2}\}}\Big(\big(-\log(r)+1\big)^{-2\lceil\frac{d-1}{2}\rceil}\cdot|\text{\textgreek{y}}|^{2}\Big)\le C(\text{\textgreek{t}})\Bigg\{ & \int_{\{\bar{t}=\text{\textgreek{t}}\}\cap\{r\le1\}}\big|\nabla_{h_{\text{\textgreek{t}},N}}^{\lceil\frac{d+1}{2}\rceil}\text{\textgreek{y}}\big|_{\big(1-\log(r_{tim})\big)\cdot h_{R_{c}}}^{2}\, dh_{N}+\label{eq:SobolevNearHorizon}\\
 & +\sum_{j=0}^{\lceil\frac{d-1}{2}\rceil}\int_{\{\bar{t}=\text{\textgreek{t}}\}\cap\{\frac{1}{2}\le r\le1\}}\big|\nabla_{h_{\text{\textgreek{t}},N}}^{j}\text{\textgreek{y}}\big|_{h_{\text{\textgreek{t}},N}}^{2}\, dh_{N}\Bigg\}.\nonumber 
\end{align}
\end{lem}
\begin{rem*}
Notice that the energy norm in the right hand side of (\ref{eq:SobolevNearHorizon})
degenerates polynomially at $\mathcal{H}^{+}$ and logarithmically
at $\partial_{tim}\mathcal{M}$.\end{rem*}
\begin{proof}
Using the fundamental theorem of calculus and a Cauchy--Schwarz inequality,
we can bound for any smooth $\text{\textgreek{f}}:[0,1]\rightarrow\mathbb{C}$
and any $x_{0}\in(0,1]$: 
\begin{align}
|\text{\textgreek{f}}(x_{0})| & \le\int_{x_{0}}^{1}\big|\frac{d\text{\textgreek{f}}}{dx}\big|\, dx+|\text{\textgreek{f}}(1)|\le\label{eq:FundamentalTheoremOfCalculusDegenerate}\\
 & \le\Big(\int_{x_{0}}^{1}x^{-1}\, dx\Big)^{\frac{1}{2}}\Big(\int_{x_{0}}^{1}x\big|\frac{d\text{\textgreek{f}}}{dx}\big|^{2}\, dx\Big)^{\frac{1}{2}}+|\text{\textgreek{f}}(1)|\le\nonumber \\
 & \le\big(-\log(x_{0})\big)\cdot\Big(\int_{0}^{1}x\big|\frac{d\text{\textgreek{f}}}{dx}\big|^{2}\, dx\Big)^{\frac{1}{2}}+|\text{\textgreek{f}}(1)|.\nonumber 
\end{align}

Let us fix a small $r_{0}=r_{0}(\text{\textgreek{t}})$, so that $\{\bar{t}=\text{\textgreek{t}}\}\cap\{r\le2r_{0}\}$
is diffeomorphic to $[0,2r_{0}]\times\mathcal{H}_{\text{\textgreek{t}}}$
and $h_{R_{c}}$ has the form (\ref{eq:MetricNearHorizon-1}) there.
It then readily follows in view of (\ref{eq:FundamentalTheoremOfCalculusDegenerate})
and a Sobolev inequality on the surfaces $\{r=const\}$ that 
\begin{align}
\sup_{\{\bar{t}=\text{\textgreek{t}}\}\cap\{r\le r_{0}\}}\Big(\big(-\log(r)\big)^{-2\lceil\frac{d-1}{2}\rceil}\cdot|\text{\textgreek{y}}|^{2}\Big)\le C(\text{\textgreek{t}})\Bigg\{ & \sum_{j=0}^{\lceil\frac{d-1}{2}\rceil}\int_{\{\bar{t}=\text{\textgreek{t}}\}\cap\{r\le2r_{0}\}}r_{hor}\big(-\log(r_{tim})\big)^{-2\lceil\frac{d-1}{2}\rceil}\big|\nabla_{\mathcal{H}_{\text{\textgreek{t}}}}^{j}(\partial_{r}\text{\textgreek{y}})\big|^{2}\, dh_{N}+\label{eq:AlmostDoneDegenerateSobolev}\\
 & +\sum_{j=0}^{\lceil\frac{d-1}{2}\rceil}\int_{\{\bar{t}=\text{\textgreek{t}}\}\cap\{r_{0}\le r\le2r_{0}\}}\big|\nabla_{h_{\text{\textgreek{t}},N}}^{j}\text{\textgreek{y}}\big|_{h_{\text{\textgreek{t}},N}}^{2}\, dh_{N}\Bigg\},\nonumber 
\end{align}
 where $\nabla_{\mathcal{H}_{\text{\textgreek{t}}}}$ is the covariant
derivative on the surfaces $\{r=const\}$ associated to the metric
$h_{\mathcal{H}_{\text{\textgreek{t}}}}$ in (\ref{eq:MetricNearHorizon-1}).
Inequality (\ref{eq:SobolevNearHorizon}) now readily follows from
(\ref{eq:AlmostDoneDegenerateSobolev}) in view of the Sobolev-inequality
in the region $\{r_{0}\le r\le\frac{1}{2}\}$: 
\begin{equation}
\sup_{\{\bar{t}=\text{\textgreek{t}}\}\cap\{r_{0}\le r\le\frac{1}{2}\}}|\text{\textgreek{y}}|^{2}\le C(\text{\textgreek{t}})\sum_{j=0}^{\lceil\frac{d+1}{2}\rceil}\int_{\{\bar{t}=\text{\textgreek{t}}\}\cap\{r_{0}\le r\le\frac{1}{2}\}}\big|\nabla_{h_{\text{\textgreek{t}},N}}^{j}\text{\textgreek{y}}\big|_{h_{\text{\textgreek{t}},N}}^{2}\, dh_{N},
\end{equation}
the Hardy-type inequality: 
\begin{align}
\sum_{j=0}^{\lceil\frac{d-1}{2}\rceil}\int_{\{\bar{t}=\text{\textgreek{t}}\}\cap\{r\le2r_{0}\}}r_{hor}\big(-\log(r_{tim})\big)^{-2\lceil\frac{d-1}{2}\rceil}\big|\nabla_{\mathcal{H}_{\text{\textgreek{t}}}}^{j}(\partial_{r}\text{\textgreek{y}})\big|^{2}\, dh_{N}\le C(\text{\textgreek{t}})\Bigg\{ & \int_{\{\bar{t}=\text{\textgreek{t}}\}\cap\{r\le1\}}\big|\nabla_{h_{\text{\textgreek{t}},N}}^{\lceil\frac{d+1}{2}\rceil}\text{\textgreek{y}}\big|_{\big(1-\log(r_{tim})\big)\cdot h_{R_{c}}}^{2}\, dh_{N}+\\
 & +\sum_{j=0}^{\lceil\frac{d-1}{2}\rceil}\int_{\{\bar{t}=\text{\textgreek{t}}\}\cap\{\frac{1}{2}\le r\le1\}}\big|\nabla_{h_{\text{\textgreek{t}},N}}^{j}\text{\textgreek{y}}\big|_{h_{\text{\textgreek{t}},N}}^{2}\, dh_{N}\Bigg\},\nonumber 
\end{align}
 and the Poincare-type inequality 
\begin{equation}
\sum_{j=0}^{\lceil\frac{d+1}{2}\rceil}\int_{\{\bar{t}=\text{\textgreek{t}}\}\cap\{r_{0}\le r\le\frac{1}{2}\}}\big|\nabla_{h_{\text{\textgreek{t}},N}}^{j}\text{\textgreek{y}}\big|_{h_{\text{\textgreek{t}},N}}^{2}\, dh_{N}\le C(\text{\textgreek{t}})\Bigg\{\int_{\{\bar{t}=\text{\textgreek{t}}\}\cap\{r\le1\}}\big|\nabla_{h_{\text{\textgreek{t}},N}}^{\lceil\frac{d+1}{2}\rceil}\text{\textgreek{y}}\big|_{h_{R_{c}}}^{2}\, dh_{N}+\sum_{j=0}^{\lceil\frac{d-1}{2}\rceil}\int_{\{\bar{t}=\text{\textgreek{t}}\}\cap\{\frac{1}{2}\le r\le1\}}\big|\nabla_{h_{\text{\textgreek{t}},N}}^{j}\text{\textgreek{y}}\big|_{h_{\text{\textgreek{t}},N}}^{2}\, dh_{N}\Bigg\}.
\end{equation}

\end{proof}

\subsection{\label{sub:ProofCorollaryImprovedDecay}Proof of Corollary \ref{cor:ImprovedPointwiseDecay}}

Inequality (\ref{eq:ImprovedPointwiseDecay}) for any dimension $d$
follows readily from the Gagliardo--Nirenberg type inequalities (\ref{eq:GagliardoNirenbergOddfinal-1})
and (\ref{eq:GagliardoNirenbergEvenFinal-1}) for $r_{0}=0$, Assumption
\ref{enu:UniformityOfGagliardoNirenberg} and the decay estimate (\ref{eq:ImprovedNonDegenerateEnergyDecay})
of Theorem \ref{thm:ImprovedDecayEnergy} for $q=1,\lfloor\frac{d+1}{2}\rfloor$
and $m+1$ in place of $m$.

In case the dimension $d$ is odd, after adding the Gagliardo--Nirenberg
type inequalities (\ref{eq:GagliardoNirenbergOddfinal-1}) (for some
small $r_{0}>0$) and (\ref{eq:SobolevNearHorizon}), and using a
Poincare type inequality and the ``critical'' Hardy-type estimate
(\ref{eq:CriticalHardy}) (together with Lemma \ref{lem:HardyForphiHyperboloids})
to bound 
\begin{equation}
\sum_{j=0}^{\lceil\frac{d-1}{2}\rceil}\int_{\{\bar{t}=\text{\textgreek{t}}\}\cap\{\frac{1}{2}\le r\le1\}}\big|\nabla_{h_{\text{\textgreek{t}},N}}^{j}\text{\textgreek{f}}\big|_{h_{\text{\textgreek{t}},N}}^{2}\, dh_{N}\le C\cdot\Big(\mathcal{E}_{en,deg}^{(0,\frac{d+1}{2},0)}[\text{\textgreek{f}}](\text{\textgreek{t}})+\Big(\mathcal{E}_{en,deg}^{(0,1,0)}[\text{\textgreek{f}}](\text{\textgreek{t}})\Big)^{\frac{1}{d-1}}\Big(\mathcal{E}_{en,deg}^{(0,\frac{d+1}{2},0)}[\text{\textgreek{f}}](\text{\textgreek{t}})\Big)^{\frac{d-2}{d-1}}\Big),
\end{equation}
 we obtain in view of Assumption \ref{enu:UniformityOfGagliardoNirenberg}:
\begin{equation}
\sup_{\{\bar{t}=\text{\textgreek{t}}\}\cap\{r\le\frac{1}{2}\}}\Big(\big(-\log(r)+1\big)^{-2\lceil\frac{d-1}{2}\rceil}\cdot|\text{\textgreek{f}}|^{2}\Big)\le C\cdot\Bigg\{\Big(\mathcal{E}_{en,deg}^{(0,1,0)}[\text{\textgreek{f}}](\text{\textgreek{t}})\Big)^{\frac{1}{d-1}}\Big(\mathcal{E}_{en,deg}^{(0,\frac{d+1}{2},0)}[\text{\textgreek{f}}](\text{\textgreek{t}})\Big)^{\frac{d-2}{d-1}}+\mathcal{E}_{en,deg}^{(0,\frac{d+1}{2},0)}[\text{\textgreek{f}}](\text{\textgreek{t}})\Bigg\}.\label{eq:AlmostDonePerfectDecay}
\end{equation}
 Therefore, from (\ref{eq:AlmostDonePerfectDecay}) and (\ref{eq:ImprovedEnergyDecay})
we deduce that 
\begin{equation}
\sup_{\{\bar{t}=\text{\textgreek{t}}\}}\Big(\big(-\log(r)+1\big)^{-2\lceil\frac{d-1}{2}\rceil}\cdot|\text{\textgreek{f}}|^{2}\Big)\lesssim_{\text{\textgreek{e}}}\text{\textgreek{t}}^{-d}\cdot\mathcal{E}_{bound}^{(d+1,\frac{d+1}{2},\lceil\text{\textgreek{d}}_{0}^{-1}\cdot2(d-1)\rceil(3\frac{d+1}{2}+1)\cdot k)}[\text{\textgreek{f}}](0)+\mathcal{F}_{deg,\text{\textgreek{e}}}^{(d+1,\frac{d+1}{2},m+\lceil\text{\textgreek{d}}_{0}^{-1}\cdot(d-1)\rceil(3\frac{d+1}{2}+1)\cdot k,k)}[F](\text{\textgreek{t}}).\label{eq:DegenerateDecayImproved}
\end{equation}

Moreover, using a standard Sobolev estimate we can also bound: 
\begin{equation}
\sup_{\{\bar{t}=\text{\textgreek{t}}\}\cap\{r\le1\}}|\text{\textgreek{f}}|^{2}\le C\Big\{\mathcal{E}_{en}^{(0,\frac{d+1}{2},0)}[\text{\textgreek{f}}](\text{\textgreek{t}})+\int_{\{\bar{t}=\text{\textgreek{t}}\}\cap\{r\le1\}}|\text{\textgreek{f}}|^{2}\, dh_{N}\Big\},
\end{equation}
 which, in view of a Poincare inequality and the Hardy inequality
(\ref{eq:CriticalHardy}) (together with Lemma \ref{lem:HardyForphiHyperboloids}),
yields: 
\begin{equation}
\sup_{\{\bar{t}=\text{\textgreek{t}}\}\cap\{r\le1\}}|\text{\textgreek{f}}|^{2}\le C\cdot\Bigg\{\mathcal{E}_{en}^{(0,\frac{d+1}{2},0)}[\text{\textgreek{f}}](\text{\textgreek{t}})+\Big(\mathcal{E}_{en,deg}^{(0,\frac{d-1}{2},0)}[\text{\textgreek{f}}](\text{\textgreek{t}})\Big)^{\frac{1}{2}}\Big(\mathcal{E}_{en,deg}^{(0,\frac{d+1}{2},0)}[\text{\textgreek{f}}](\text{\textgreek{t}})\Big)^{\frac{1}{2}}\Bigg\}.\label{eq:AlmostImprovedNonDegenerateDecay}
\end{equation}
 Therefore, from Theorem \ref{thm:ImprovedDecayEnergy} and (\ref{eq:AlmostImprovedNonDegenerateDecay})
we obtain 
\begin{equation}
\sup_{\{\bar{t}=\text{\textgreek{t}}\}\cap\{r\le1\}}|\text{\textgreek{f}}|^{2}\lesssim_{\text{\textgreek{e}}}\text{\textgreek{t}}^{-\frac{d}{2}}\mathcal{E}_{bound}^{(d+1,\frac{d+1}{2},\lceil\text{\textgreek{d}}_{0}^{-1}\cdot2(d-1)\rceil(3\frac{d+1}{2}+1)\cdot k)}[\text{\textgreek{f}}](0)+\mathcal{F}_{deg,\text{\textgreek{e}},0,d}[F](\text{\textgreek{t}}),
\end{equation}
 which combined with (\ref{eq:DegenerateDecayImproved}) yields the
required estimate (\ref{eq:ImprovedPointwiseDecayAway}) for the pointwise
decay of $\text{\textgreek{f}}$.

In case $m\ge1$, we can bound in view of the Gagliardo-Nirenberg
inequality (\ref{eq:GagliardoNirenbergOddfinal-ForHigherDerivatives})
for some fixed small $r_{0}>0$: 
\begin{equation}
\sum_{i=1}^{m}\sup_{\{\bar{t}=\text{\textgreek{t}}\}\cap\{r\ge r_{0}\}}|\nabla_{g}^{i}\text{\textgreek{f}}|_{h}^{2}\le C_{m,r_{0}}\Bigg\{\Big(\mathcal{E}_{en,deg}^{(0,\frac{d+1}{2},m)}[\text{\textgreek{f}}](\text{\textgreek{t}})\Big)^{\frac{1}{2}}\Big(\mathcal{E}_{en,deg}^{(0,\frac{d+1}{2},m+1)}[\text{\textgreek{f}}](\text{\textgreek{t}})\Big)^{\frac{1}{2}}+\mathcal{E}_{en,deg}^{(0,\frac{d+1}{2},m+1)}[\text{\textgreek{f}}](\text{\textgreek{t}})\Bigg\}.\label{eq:BoundHigherOrderOdd-1}
\end{equation}
Therefore, from Theorem \ref{thm:ImprovedDecayEnergy} and (\ref{eq:BoundHigherOrderOdd-1})
we obtain: 
\begin{equation}
\sum_{i=1}^{m}\sup_{\{\bar{t}=\text{\textgreek{t}}\}\cap\{r\ge r_{0}\}}|\nabla_{g}^{i}\text{\textgreek{f}}|_{h}^{2}\lesssim_{m,r_{0},\text{\textgreek{e}}}\text{\textgreek{t}}^{-d-1}\mathcal{E}_{bound}^{(d+1,\frac{d+1}{2},m+1+\lceil\text{\textgreek{d}}_{0}^{-1}\cdot2(d-1)\rceil(3\frac{d+1}{2}+1)\cdot k)}[\text{\textgreek{f}}](0)+\mathcal{F}_{deg,\text{\textgreek{e}},m+2,d}[F](\text{\textgreek{t}}).\label{eq:BoundDerivativeAway}
\end{equation}

Moreover, using Lemma \ref{lem:Gagliardo-NirenbergDegenerate} for
$\partial^{i}\text{\textgreek{f}}$ in place of $\text{\textgreek{y}}$,
we can bound (in view also of Theorem \ref{thm:ImprovedDecayEnergy}):
\begin{align}
\sum_{i=1}^{m+1}\sup_{\{\bar{t}=\text{\textgreek{t}}\}\cap\{r\le\frac{1}{2}\}}\Big(\big(1-\log(r)\big)^{-2\lceil\frac{d-1}{2}\rceil}|\nabla_{g}^{i}\text{\textgreek{f}}|_{h}^{2}\Big) & \lesssim_{m}\mathcal{E}_{en,deg}^{(0,\frac{d+1}{2},m+2)}[\text{\textgreek{f}}](\text{\textgreek{t}})\lesssim_{m,\text{\textgreek{e}}}\label{eq:DegenerateSobolev}\\
 & \lesssim_{m,\text{\textgreek{e}}}\text{\textgreek{t}}^{-d-1}\mathcal{E}_{bound}^{(d+1,\frac{d+1}{2},m+2+\lceil\text{\textgreek{d}}_{0}^{-1}\cdot2(d-1)\rceil(3\frac{d+1}{2}+1)\cdot k)}[\text{\textgreek{f}}](0)+\mathcal{F}_{deg,\text{\textgreek{e}},m+1,d}[F](\text{\textgreek{t}}).\nonumber 
\end{align}
 Since the function $\big(-\log(r)\big)^{2\lceil\frac{d-1}{2}\rceil}$
is integrable near $r=0$, using the fundamental theorem of calculus
and a Cauchy--Schwarz inequality we can bound for any smooth function
$\text{\textgreek{y}}$ on $\{\bar{t}=\text{\textgreek{t}}\}\cap\{r\le\frac{1}{2}\}$:%
\footnote{Here we have also used Assumption \ref{enu:BoundednessVolume} on
the boundedness of the volume of the region $\{\bar{t}=\text{\textgreek{t}}\}\cap\{r\lesssim1\}$.%
} 
\begin{equation}
\sup_{\{\bar{t}=\text{\textgreek{t}}\}\cap\{r\le\frac{1}{2}\}}|\text{\textgreek{y}}|^{2}\le C\Bigg\{\sup_{\{\bar{t}=\text{\textgreek{t}}\}\cap\{r\le\frac{1}{2}\}}\Big(\big(1-\log(r)\big)^{-2\lceil\frac{d-1}{2}\rceil}|\nabla_{h_{\text{\textgreek{t}},N}}\text{\textgreek{y}}|_{h_{\text{\textgreek{t}},N}}^{2}\Big)+\sup_{\{\bar{t}=\text{\textgreek{t}}\}\cap\{\frac{1}{4}\le r\le\frac{1}{2}\}}|\text{\textgreek{y}}|^{2}\Bigg\}.\label{eq:BoundWithMoreLossOfDerivatives}
\end{equation}
 Thus, from (\ref{eq:BoundDerivativeAway}), (\ref{eq:DegenerateSobolev})
and (\ref{eq:BoundWithMoreLossOfDerivatives}) (for $\partial^{i}\text{\textgreek{f}}$
in place of $\text{\textgreek{y}}$) we deduce the desired bound (\ref{eq:ImprovedDecayOddHigherOrder}):
\begin{equation}
\sum_{i=1}^{m}\sup_{\{\bar{t}=\text{\textgreek{t}}\}}\big|\nabla^{m}\text{\textgreek{f}}\big|_{h}^{2}\lesssim_{m,\text{\textgreek{e}}}\text{\textgreek{t}}^{-d-1}\cdot\mathcal{E}_{bound}^{(d+1,\frac{d+1}{2},m+2+\lceil\text{\textgreek{d}}_{0}^{-1}\cdot2(d-1)\rceil(3\frac{d+1}{2}+1)\cdot k)}[\text{\textgreek{f}}](0)+\mathcal{F}_{deg,\text{\textgreek{e}},m+2,d}[F](\text{\textgreek{t}}).\label{eq:ImprovedDecayOddHigherOrder-1}
\end{equation}
 \qed

\appendix

\section{\label{sec:RiemannianMetric}Construction of the natural Riemannian
metrics $h_{\text{\textgreek{t}},N}$ and $h_{\text{\textgreek{t}},K,\text{\textgreek{F}}}$}

Let $(\mathcal{M}^{d+1},g)$, $d\ge1$, be a Lorentzian manifold and
$\bar{t}:\mathcal{M}\rightarrow\mathbb{R}$ a smooth function with
acausal level sets.%
\footnote{This manifold $\mathcal{M}$ will correspond to the manifold $\mathcal{M}\backslash\partial\mathcal{M}$
in the language of Sections \ref{sec:Firstdecay}-\ref{sec:Improved-polynomial-decay}.%
} For any timelike vector field $N$ on $\mathcal{M}$ such that $d\bar{t}(N)=1$,
there exists a special Riemannian metric $h_{\text{\textgreek{t}},N}$
defined on the $S_{\text{\textgreek{t}}}=\{\bar{t}=\text{\textgreek{t}}\}$
hypersurfaces naturally associated to the vector field $N$. This
metric does not necessarily coincide with the induced metric on $S_{\text{\textgreek{t}}}$,
but its usefulness lies in the fact that the Laplace operator associated
to $h_{\text{\textgreek{t}},N}$ appears naturally in a ``useful''
decomposition of the wave operator $\square_{g}$. A similar construction
of a Riemannian metric $h_{\text{\textgreek{t}},K,\text{\textgreek{F}}}$
(not natural in this case) can be constructed in case one has two
vector fields $K,\text{\textgreek{F}}$ with merely timelike span,
such that $K$ satisfies $d\bar{t}(K)=1$ and becomes non timelike
only on a set with compact intersection with each $S_{\text{\textgreek{t}}}$,
and $\text{\textgreek{F}}$ is tangent to the level sets of $\bar{t}$.
We will now proceed with the details of the construction of these
metrics.

\subsection{Construction of $h_{\text{\textgreek{t}},N}$}

Let $N$ be a timelike vector field on $\mathcal{M}$ with $d\bar{t}(N)=1$.
Let $\mathcal{V}_{\mathcal{M}}$ denote vector bundle on $S_{\text{\textgreek{t}}}$
defined as the pullback of $T\mathcal{M}$ through the inclusion $i:S_{\text{\textgreek{t}}}\rightarrow\mathcal{M}$,
the latter giving rise to the natural inclusion $TS_{\text{\textgreek{t}}}\hookrightarrow\mathcal{V}_{\mathcal{M}}$.
The vector bundle $\mathcal{V}_{\mathcal{M}}$ inherits from $T\mathcal{M}$
the Lorentzian metric $i^{*}g$ and the timelike section $i^{*}N$,
which will be denoted as $g$ and $N$ respectively for notational
simplicity. Similarly, the one form $d\bar{t}$ on $\mathcal{M}$
can also be viewed as a section of the dual bundle $\mathcal{V}_{\mathcal{M}}^{*}$.
Let us denote with $g^{-1}$ the Lorentzian metric on $\mathcal{V}_{\mathcal{M}}^{*}$
associated to $g$ on $\mathcal{V}_{\mathcal{M}}$.

Let $\mathcal{K}_{N}\hookrightarrow\mathcal{V}_{\mathcal{M}}^{*}$
be the vector subbundle of $\mathcal{V}_{\mathcal{M}}^{*}$ defined
as the set of $v\in\text{\textgreek{G}}(\mathcal{V}_{\mathcal{M}}^{*})$
such that $v(N)=0$. Then the Lorentzian metric $g^{-1}$ on $\mathcal{V}_{\mathcal{M}}^{*}$
induces a metric $h_{inv,\text{\textgreek{t}},N}$ on $\mathcal{K}_{N}$.
It is easy to verify that $h_{int,\text{\textgreek{t}},N}$ is positive
definite, owing to the fact that $g(N,N)<0$.

We can naturally identify $T^{*}S_{\text{\textgreek{t}}}$ and $\mathcal{K}_{N}$
in the following way: Since $\mathcal{V}_{\mathcal{M}}$ can be split
as $TS_{\text{\textgreek{t}}}\oplus(N)$ (where $(N)$ is the line
bundle spanned by $N$), any $\text{\textgreek{w}}\in\text{\textgreek{G}}(T^{*}S_{\text{\textgreek{t}}})$
can be extended to a section $\text{\textgreek{w}}_{N}\in\text{\textgreek{G}}(\mathcal{V}_{\mathcal{M}}^{*})$
by demanding that $\text{\textgreek{w}}_{N}|_{TS_{\text{\textgreek{t}}}}=\text{\textgreek{w}}$
and $\text{\textgreek{w}}_{N}(N)=0$. But then, since $\text{\textgreek{w}}_{N}(N)=0$,
$\text{\textgreek{w}}_{N}$ is a section of $\mathcal{K}$. It is
easy to verify that the mapping $\text{\textgreek{w}}\rightarrow\text{\textgreek{w}}_{N}$
is a vector bundle isomorphism. Thereofore, $T^{*}S_{\text{\textgreek{t}}}$
inherits the Lorentzian metric $h_{inv,\text{\textgreek{t}},N}$ of
$\mathcal{K}_{N}$.

Finally, we define the (positive definite) metric $h_{\text{\textgreek{t}},N}$
on $TS_{\text{\textgreek{t}}}$ as the dual metric of $h_{inv,\text{\textgreek{t}},N}$
on $T^{*}S_{\text{\textgreek{t}}}$. Thus, $(S_{\text{\textgreek{t}}},h_{\text{\textgreek{t}},N})$
becomes a Riemannian manifold. Notice that the following relation
holds:
\begin{equation}
dvol_{g}=\sqrt{-g(N,N)}\cdot d\bar{t}\wedge dvol_{h_{\text{\textgreek{t}},N}}\label{eq:RadonNikodym-1}
\end{equation}
 where $dvol_{g}$ is the natural volume form on $\mathcal{M}$ associated
with $g$, while $dvol_{h_{\text{\textgreek{t}},N}}$ is the natural
volume form on $(S_{\text{\textgreek{t}}},h_{\text{\textgreek{t}},N})$
extended to a $d$-form on $\mathcal{M}$ by the requirement that
$i_{N}dvol_{h_{\text{\textgreek{t}},N}}=0$. 

The connection of the metric $h_{\text{\textgreek{t}},N}$ with the
wave operator $\square_{g}$ on $\mathcal{M}$ is the following: In
any local coordinate system $(x^{1},\ldots,x^{d})$ on $S_{\text{\textgreek{t}}}$,
extended to a local coordinate system $(\bar{t},x^{1},\ldots,x^{d})$
on $\mathcal{M}$ by the requirement that $N(x^{i})=0$, the wave
operator $\square_{g}$ on $\mathcal{M}$ around $S_{\text{\textgreek{t}}}$
satisfies: 
\begin{align}
\square_{g}= & (\sqrt{-g})^{-1}N\big(\sqrt{-g}g^{\bar{t}\bar{t}}N\big)+(\sqrt{-g})^{-1}\partial_{x^{i}}\big(\sqrt{-g}g^{x^{i}\bar{t}}N\big)+(\sqrt{-g})^{-1}N\big(\sqrt{-g}g^{x^{i}\bar{t}}\partial_{x^{i}}\big)+\label{eq:SplittingWaveOperator}\\
 & +\text{\textgreek{D}}_{h_{\text{\textgreek{t}},N},N},\nonumber 
\end{align}
 where the operator $\text{\textgreek{D}}_{h_{\text{\textgreek{t}},N},N}$
on $S_{\text{\textgreek{t}}}$ is defined as: 
\begin{equation}
\text{\textgreek{D}}_{h_{\text{\textgreek{t}},N},N}=\frac{1}{\sqrt{-g(N,N)}}div_{h_{\text{\textgreek{t}},N}}\big(\sqrt{-g(N,N)}d\big)\label{eq:PerturbedLaplacian}
\end{equation}
and $div_{h_{\text{\textgreek{t}},N}}$ acting on the one form $\text{\textgreek{w}}$
on $S_{\text{\textgreek{t}}}$ is defined as the divergence (with
respect to $h_{\text{\textgreek{t}},N})$ of the dual vector field
$\text{\textgreek{w}}^{\sharp}$.%
\footnote{Equivalently, it is the dual of the gradient operator on functions
with respect to the inner product $\int_{S_{\text{\textgreek{t}}}}\langle\cdot,\cdot\rangle_{h_{\text{\textgreek{t}},N}}\, dvol_{h_{\text{\textgreek{t}},N}}$.%
} Equivalently, $\square_{g}\text{\textgreek{f}}$ takes the following
covariant form for any $\text{\textgreek{y}}\in C^{\infty}(\mathcal{M})$
(assuming that $\mathcal{M}$ is orientable): 
\begin{equation}
\square_{g}\text{\textgreek{y}}=div_{g}\big((N\text{\textgreek{y}})\cdot d\bar{t}\big)+\star\mathcal{L}_{N}\big((N_{tan}\text{\textgreek{y}})\cdot dvol_{g}\big)+\text{\textgreek{D}}_{h_{\text{\textgreek{t}},N},N}\text{\textgreek{y}},\label{eq:SplittingWaveOperatorCoordinateInvariant}
\end{equation}
where $\star$ is the Hodge star operator on $(\mathcal{M},g)$ and
$N_{tan}$ is the projection of $\nabla\bar{t}$ on $S_{\text{\textgreek{t}}}$
along $N$.

\subsection{Construction of $h_{\text{\textgreek{t}},K,\text{\textgreek{F}}}$}

Let $K,\text{\textgreek{F}}$ be two vector fields on $\mathcal{M}$
with timelike span such that $d\bar{t}(K)=1$ and $\text{\textgreek{F}}$
is tangent to the leaves of the foliation $\{S_{\text{\textgreek{t}}}\}_{\text{\textgreek{t}}\in\mathbb{R}}$.
Assume also that for any $\text{\textgreek{t}}\in\mathbb{R}$, the
set $\mathfrak{A}_{\text{\textgreek{t}}}=\{p\in S_{\text{\textgreek{t}}}:g(K,K)\ge0\}$
is compact. In this case, we will also define a Riemannian metric
$h_{\text{\textgreek{t}},K,\text{\textgreek{F}}}$ on $S_{\text{\textgreek{t}}}$
which will prove helpful in decomposing the wave operator $\square_{g}$,
but contrary to the metric $h_{\text{\textgreek{t}},N}$ constructed
in the previous Section, the construction of $h_{\text{\textgreek{t}},K,\text{\textgreek{F}}}$
will not be natural.

Proceeding as in the previous Section, we can define the natural metric
$h_{\text{\textgreek{t}},K}$ on $S_{\text{\textgreek{t}}}$, but
$h_{\text{\textgreek{t}},K}$ will now not be Riemannian on $\mathfrak{A}_{\text{\textgreek{t}}}$.
Moreover, $h_{\text{\textgreek{t}},K}$ will be singular at the points
where $g(K,K)=0$, although its inverse (i.\,e.~the associated metric
on the dual bndle $T^{*}S_{\text{\textgreek{t}}}$) will be smooth
everywhere on $S_{\text{\textgreek{t}}}$.

Since the span of $K,\text{\textgreek{F}}$ is everywhere timelike
on $\mathcal{M}$ and $\mathfrak{A}_{\text{\textgreek{t}}}$ is compact,
for any $\text{\textgreek{t}}\in\mathbb{R}$ there exists some $C_{\text{\textgreek{t}}}>0$
such that everywhere on $S_{\text{\textgreek{t}}}$:
\begin{equation}
\frac{g(K,K)}{g(\text{\textgreek{F}},\text{\textgreek{F}})g(K,K)-\big(g(\text{\textgreek{F}},K)\big)^{2}}>-C_{\text{\textgreek{t}}}.\label{eq:BoundPlaneNorm}
\end{equation}
Notice that the bound (\ref{eq:BoundPlaneNorm}) holds trivially outside
$\mathfrak{A}_{\text{\textgreek{t}}}$. Thus, if $\text{\textgreek{F}}_{\bot}$
denotes the projecton of $\text{\textgreek{F}}$ on the orthogonal
complement of $K$, we have
\begin{equation}
\frac{1}{g(\text{\textgreek{F}}_{\bot},\text{\textgreek{F}}_{\bot})}=\frac{g(K,K)}{g(\text{\textgreek{F}},\text{\textgreek{F}})g(K,K)-\big(g(\text{\textgreek{F}},K)\big)^{2}}>-C_{\text{\textgreek{t}}},\label{eq:OrthogonalComplementF}
\end{equation}
and 
\begin{equation}
g(\text{\textgreek{F}}_{\bot},\text{\textgreek{F}}_{\bot})\le0\label{eq:TimelikeF}
\end{equation}
 on $\cup_{\text{\textgreek{t}}}\mathfrak{A}_{\text{\textgreek{t}}}$.
Recall also that for any vector fields $X,Y$ with timelike (or null)
span the following inverted Cauchy inequality holds pointwise: 
\begin{equation}
\big(g(X,Y)\big)^{2}\ge|g(X,X)|\cdot|g(Y,Y)|.\label{eq:InvertedCauchy}
\end{equation}
Thus, for any vector field $X$ on $\mathcal{M}$ with $g(X,K)=0$
and $g(X,X)\le0$ (notice that such a vector field must be identically
$0$ outside $\cup_{\text{\textgreek{t}}}\mathfrak{A}_{\text{\textgreek{t}}}$)
we can bound due to (\ref{eq:OrthogonalComplementF}), (\ref{eq:TimelikeF})
and (\ref{eq:InvertedCauchy}): 
\begin{equation}
g(X,X)\ge\frac{\big(g(X,\text{\textgreek{F}})\big)^{2}}{-|g(\text{\textgreek{F}}_{\bot},\text{\textgreek{F}}_{\bot})|}=\frac{\big(g(X,\text{\textgreek{F}})\big)^{2}}{g(\text{\textgreek{F}}_{\bot},\text{\textgreek{F}}_{\bot})}>-C_{\text{\textgreek{t}}}\big(g(X,\text{\textgreek{F}})\big)^{2}.\label{eq:BoundForRiemannian}
\end{equation}

Recall that for any $\text{\textgreek{w}}\in\text{\textgreek{G}}(T^{*}S_{\text{\textgreek{t}}})$,
extended to an element of $\text{\textgreek{G}}(\mathcal{V}_{M}^{*})$
by the condition $\text{\textgreek{w}}(K)=0$, we have everywhere
on $S_{\text{\textgreek{t}}}$: 
\begin{equation}
h_{\text{\textgreek{t}},K}^{-1}(\text{\textgreek{w}},\text{\textgreek{w}})=g^{-1}(\text{\textgreek{w}},\text{\textgreek{w}}).
\end{equation}
The bound (\ref{eq:BoundForRiemannian}) then readily implies that
for any $\text{\textgreek{w}}\in\text{\textgreek{G}}(T^{*}S_{\text{\textgreek{t}}})$
we have 
\begin{equation}
h_{\text{\textgreek{t}},K}^{-1}(\text{\textgreek{w}},\text{\textgreek{w}})+C_{\text{\textgreek{t}}}\big(\text{\textgreek{w}}(\text{\textgreek{F}})\big)^{2}>0.\label{eq:RiemannianDual}
\end{equation}
Therefore, the symmetric $(2,0)$-tensor 
\begin{equation}
h_{inv,\text{\textgreek{t}},K,\text{\textgreek{F}}}\doteq h_{\text{\textgreek{t}},K}^{-1}+C_{\text{\textgreek{t}}}\cdot\text{\textgreek{F}}\otimes\text{\textgreek{F}}\label{eq:FinalRiemannianDual}
\end{equation}
is a positive definite metric on the dual bundle $T^{*}S_{\text{\textgreek{t}}}$,
and its inverse $h_{\text{\textgreek{t}},K,\text{\textgreek{F}}}$
is a Riemannian metric on $S_{\text{\textgreek{t}}}$. 

Notice that in this case, $S_{\text{\textgreek{t}}}$ carries two
volume forms, $dvol_{h_{\text{\textgreek{t}},K,\text{\textgreek{F}}}}$
and $i_{K}dvol_{g}$, and their Radon-Nikodym derivative, i.\,e.
\begin{equation}
\mathfrak{m}_{\text{\textgreek{t}},K,\text{\textgreek{F}}}\doteq\frac{\text{\textgreek{d}}(i_{K}dvol_{g})}{\text{\textgreek{d}}(dvol_{h_{\text{\textgreek{t}},K,\text{\textgreek{F}}}})},\label{eq:RadonNikodym}
\end{equation}
is a smooth function on $S_{\text{\textgreek{t}}}$ satisfying 
\begin{equation}
\mathfrak{m}_{\text{\textgreek{t}},K,\text{\textgreek{F}}}\sim\sqrt{-\Big(g(K,K)-\frac{\big(g(K,\text{\textgreek{F}})\big)^{2}}{g(\text{\textgreek{F}},\text{\textgreek{F}})}\Big)}\label{eq:Degeneration}
\end{equation}
(compare with (\ref{eq:RadonNikodym-1})).

Using the Riemannian metric $h_{\text{\textgreek{t}},K,\text{\textgreek{F}}}$,
we can decompose the wave operator as: 
\begin{equation}
\square_{g}\text{\textgreek{y}}=div_{g}\big((K\text{\textgreek{y}})\cdot d\bar{t}\big)+\star\mathcal{L}_{K}\big((K_{tan}\text{\textgreek{y}})\cdot dvol_{g}\big)-C_{\text{\textgreek{t}}}\big(\star\mathcal{L}_{\text{\textgreek{F}}}\big((\text{\textgreek{F}}\text{\textgreek{y}})\cdot dvol_{g}\big)\big)+\text{\textgreek{D}}_{h_{\text{\textgreek{t}},K,\text{\textgreek{F}}},mod}\text{\textgreek{y}},\label{eq:SplittingWaveOperatorCoordinateInvariantTwo}
\end{equation}
 where $\star$ is the Hodge star operator on $(\mathcal{M},g)$,
$K_{tan}$ is the projection of $\nabla\bar{t}$ on $S_{\text{\textgreek{t}}}$
along $K$ and the elliptic operator $\text{\textgreek{D}}_{h_{\text{\textgreek{t}},K,\text{\textgreek{F}}},mod}$
on $S_{\text{\textgreek{t}}}$ is defined as
\begin{equation}
\text{\textgreek{D}}_{h_{\text{\textgreek{t}},K,\text{\textgreek{F}}},mod}\text{\textgreek{y}}\doteq\mathfrak{w}_{\text{\textgreek{t}},K,\text{\textgreek{F}}}^{-1}\cdot div_{h_{\text{\textgreek{t}},K,\text{\textgreek{F}}}}\big(\mathfrak{w}_{\text{\textgreek{t}},K,\text{\textgreek{F}}}\cdot d\text{\textgreek{y}}\big).\label{eq:EllipticOperatorTwoVectorFields}
\end{equation}

\section{\label{sec:Elliptic-estimates}Elliptic estimates on asymptotically
Euclidean Riemannian manifolds with boundary }

In this section, we will establish some general elliptic estimates
for asymptotically Euclidean Riemannian manifolds. This class of manifolds
will include, in particular, the slices $\{\bar{t}=const\}$ of the
hyperboloidal foliation of the spacetimes $(\mathcal{M},g)$ appearing
in Sections \ref{sec:Firstdecay} and \ref{sec:Improved-polynomial-decay},
equipped with a Riemannian metric of the form that was introduced
in Section \ref{sec:RiemannianMetric} of the Appendix. However, the
results of the current section might also be of independent interest. 

Let $\mathcal{S}^{d}$ be a smooth manifold with boundary of dimension
$d\ge3$, with smooth compact boundary $\partial\mathcal{S}^{d-1}$
(allowed to be empty). We will assume that $\partial\mathcal{S}$
splits into two (not necessarily connected) components: 
\begin{equation}
\partial\mathcal{S}=\partial_{tim}\mathcal{S}\cup\partial_{hor}\mathcal{S}.
\end{equation}
 The reason for assuming such a splitting for $\partial\mathcal{S}$
is that the hypersurfaces $\{\bar{t}=\text{\textgreek{t}}\}$ of the
spacetimes $(\mathcal{M},g)$ in Sections \ref{sec:Firstdecay} and
\ref{sec:Improved-polynomial-decay} (i.\,e.~the hypersurfaces on
which the elliptic estimates of the current Section will be applied)
have boundary $\{\bar{t}=\text{\textgreek{t}}\}\cap\partial\mathcal{M}$,
which is split as the disjoint union of $\{\bar{t}=\text{\textgreek{t}}\}\cap\partial_{tim}\mathcal{M}$
and $\{\bar{t}=\text{\textgreek{t}}\}\cap\mathcal{H}$.

Let $h$ be a Riemannian metric on $\mathcal{S}\backslash\partial\mathcal{S}$.
We would like $h$ to model the Riemannian metric $h_{\text{\textgreek{t}},K_{R_{c}},\text{\textgreek{F}}}$
on the hypersurfaces $\{\bar{t}=\text{\textgreek{t}}\}$ of Section
\ref{sec:Improved-polynomial-decay} (constructed as in the previous
Section of the Appendix). To this end, we assume that $(\mathcal{S}\backslash\partial\mathcal{S},h)$
is asymptotically flat in the sense that there exists a compact subset
$\mathcal{K}\subset\mathcal{S}$ containing $\partial\mathcal{S}$
such that $\mathcal{S}\backslash\mathcal{K}$ has a finite number
of connected components, each mapped diffeomorphically onto $(R_{0},+\infty)\times\mathbb{S}^{d-1}$
through a coordinate chart $(r,\text{\textgreek{sv}})$, and in this
chart $h$ has the expression 
\begin{equation}
h=dr^{2}+r^{2}g_{\mathbb{S}^{d-1}}+h_{as},\label{eq:EuclideanAsymptotics}
\end{equation}
 with 
\begin{equation}
|\nabla^{m}h_{as}|_{h}=O(r^{-m-1})\label{eq:EuclideanAsymptoticsRates}
\end{equation}
 for all integers $m\ge0$. In what follows, $\nabla$ will denote
the covariant connection with respect to $h$. 

We extend $r$ smoothly on the whole of $\mathcal{S}$, so that it
is strictly positive on $\mathcal{S}\backslash\partial\mathcal{S}$
and satisfies $r=0$ and $dr\neq0$ on $\partial\mathcal{S}$, if
$\partial\mathcal{S}\neq\emptyset$. If $\partial\mathcal{S}=\emptyset$,
we simply require that $r\ge1$ everywhere on $\mathcal{S}$. Notice
that the assumtion $dr\neq0$ on $\partial\mathcal{S}$ together with
the compactness of $\partial\mathcal{S}$ imply that for $\{r\ll1\}$
the level sets of the function $r$ are smooth hypersurfaces of $\mathcal{S}$,
and $r$ can be used as a coordinate function.

As for the behaviour of $h$ near the boundary $\partial\mathcal{S}$,
we impose the following assumptions (in accordance with the behaviour
of metric $h_{\text{\textgreek{t}},K_{R_{c}},\text{\textgreek{F}}}$
on the hypersurfaces $\{\bar{t}=\text{\textgreek{t}}\}$ of Section
\ref{sec:Improved-polynomial-decay}): Let us denote by $h_{tan}$
the induced metric on the $\{r=const\}$ hypersurfaces for $\{0<r\ll1\}$,
and extend it to a symmetric $(0,2)$-form on $\mathcal{S}\backslash\partial\mathcal{S}$
by the requirement that $h_{tan}(\nabla r,\cdot)\equiv0$. Then we
\underline{assume} that the metric $h$ in the region $\{r\ll1\}$
takes the following form:
\begin{enumerate}
\item Near the $\partial_{hor}\mathcal{S}$ component of the boundary:
\begin{equation}
h=\big(r^{-1}+O(1)\big)dr^{2}+h_{tan},\label{eq:MetricNearHorizon}
\end{equation}
 and $h_{tan}$ extends smoothly on $\partial_{hor}\mathcal{S}$,
with $h_{tan}|_{\partial_{hor}\mathcal{S}}$ being positive definite. 
\item Near the $\partial_{tim}\mathcal{S}$ component of the boundary: 
\begin{equation}
h=\big(1+O(r)\big)dr^{2}+h_{tan},\label{eq:MetricNearTimelikeBoundary}
\end{equation}
and $h_{tan}$ extends smoothly on $\partial_{tim}\mathcal{S}$, with
$h_{tan}|_{\partial_{tim}\mathcal{S}}$ being positive definite.
\end{enumerate}
We will also assume that we are given a continuous function $\text{\textgreek{w}}:\mathcal{S}\rightarrow[0,+\infty)$
which is smooth on $\mathcal{S}\backslash\partial_{hor}\mathcal{S}$,
such that 
\begin{itemize}
\item $\text{\textgreek{w}}=0$ only on $\partial_{hor}\mathcal{S}$
\item $\text{\textgreek{w}}=cr^{\frac{1}{2}}\big(1+O(r)\big)$ for some
$c>0$ near $\partial_{hor}\mathcal{S}$
\item $\text{\textgreek{w}}>0$ near $\partial_{tim}\mathcal{S}$ and
\item $\text{\textgreek{w}}=1+O(r^{-1})$ in the region $\{r\gg1\}$.
\end{itemize}
For this class of Riemannian manifolds we will establish a series
of elliptic estimates in the following sections. But before that,
we have to carry out the construction of some geometric objects that
will appear in our estimates.

\subsection{Geometric constructions on $(\mathcal{S},h)$}

\noindent We will define the perturbed Laplacian $\text{\textgreek{D}}_{h,\text{\textgreek{w}}}$
associated to $h,\text{\textgreek{w}}$ by the relation: 
\begin{equation}
\text{\textgreek{D}}_{h,\text{\textgreek{w}}}=\text{\textgreek{w}}^{-1}div_{h}\big(\text{\textgreek{w}}\cdot d\big).\label{eq:PerturbedLaplacian-1}
\end{equation}
 This operator models the operator (\ref{eq:EllipticOperatorTwoVectorFields})
associated to the metric $h_{\text{\textgreek{t}},K_{R_{c}},\text{\textgreek{F}}}$
on the hypersurfaces $\{\bar{t}=\text{\textgreek{t}}\}$ of the spacetimes
$(\mathcal{M},g)$ of Section \ref{sec:Improved-polynomial-decay}.

On each connected component of the region $\{r\gg1\}$ in the $(r,\text{\textgreek{sv}})$
coordinate chart, (\ref{eq:PerturbedLaplacian-1}) takes the form:
\begin{equation}
\text{\textgreek{D}}_{h,\text{\textgreek{w}}}=\text{\textgreek{D}}_{h}+O(r^{-2})\partial_{r}+O(r^{-3})\partial_{\text{\textgreek{sv}}}.\label{eq:PerturbedLaplacianAway}
\end{equation}
On the other hand, near the boundary $\partial\mathcal{S}$ we have
the following relations:
\begin{enumerate}
\item According to (\ref{eq:MetricNearHorizon}), if $Y=|\nabla r|_{h}^{-2}\cdot\nabla r$,%
\footnote{Notice that $Y$ can be extended smoothly on $\partial\mathcal{S}$.%
} near $\partial_{hor}\mathcal{S}$ we have: 
\begin{equation}
\text{\textgreek{D}}_{h,\text{\textgreek{w}}}=a^{-1}Y\big(ra(1+O(r))\cdot Y\big)+\text{\textgreek{D}}_{h_{tan}}+O(r)\cdot X,\label{eq:LaplacianSingularMetric}
\end{equation}
 where $a$ is a positive function smooth up to $\partial_{hor}\mathcal{S}$,
$\text{\textgreek{D}}_{h_{tan}}$ is the Laplacian of the induced
metric on the $\{r=const\}$ hypersuerfaces and $X$ is a vector field
in a neighborhood of $\partial\mathcal{S}$ smooth up to $\partial\mathcal{S}$.
Hence, $\text{\textgreek{D}}_{h,\text{\textgreek{w}}}$ is degenerate
elliptic near $\partial_{hor}\mathcal{S}$.
\item According to (\ref{eq:MetricNearTimelikeBoundary}), near $\partial_{tim}\mathcal{S}$
we have 
\begin{equation}
\text{\textgreek{D}}_{h,\text{\textgreek{w}}}=Y\big((1+O(r))\cdot Y\big)+\text{\textgreek{D}}_{h_{tan}}+X.\label{eq:LaplacianSingularMetric-2}
\end{equation}
Hence, $\text{\textgreek{D}}_{h,\text{\textgreek{w}}}$ is uniformly
elliptic near $\partial_{tim}\mathcal{S}$.
\end{enumerate}
It will be convenient to have a canonical coordinate ``chart'' near
the boundary of $\partial\mathcal{S}$. If $r_{1}$ is small enough,
we can define the diffeomorphism 
\begin{equation}
\mathcal{J}:\{r\le r_{1}\}\subset\mathcal{S}\rightarrow[0,r_{1}]\times\partial\mathcal{S},
\end{equation}
 so that for any point $p\in\{r\le r_{1}\}$:
\begin{itemize}
\item $\text{\textgreek{p}}_{1}(\mathcal{J}(p))=r(p)$ 
\item $\text{\textgreek{p}}_{2}(\mathcal{J}(p))$ is the unique point on
$\partial\mathcal{S}$ connected with $p$ through an integral line
of the vector field $Y=|\nabla r|_{h}^{-2}\cdot\nabla r$. 
\end{itemize}
\noindent In the above, $\text{\textgreek{p}}_{1}:[0,r_{1}]\times\partial\mathcal{S}\rightarrow[0,r_{1}]$
and $\text{\textgreek{p}}_{2}:[0,r_{1}]\times\partial\mathcal{S}\rightarrow\partial\mathcal{S}$
are the projections onto the first and second factor of $[0,r_{1}]\times\partial\mathcal{S}$
respectively. 

Notice that through this map, the vector field $Y$ is mapped to the
vector field $\partial_{r}$, i.\,e.~the coordinate vector field
on $[0,r_{1}]\times\partial\mathcal{S}$ which is tangent to the lines
$\{\text{\textgreek{p}}_{2}=const\}$ and is mapped to $\frac{d}{dr}$
by $\text{\textgreek{p}}_{1}$. Moreover, $h_{tan}$ is a smooth non-degenerate
Riemannian metric on the hypersurfaces $\{\text{\textgreek{p}}_{1}=\text{\textgreek{r}}\}$,
varying smoothly with $\text{\textgreek{r}}$. From now on, we will
assume that $\{r\le r_{1}\}$ and $[0,r_{1}]\times\partial\mathcal{S}$
have been identified through $\mathcal{J}$.

\paragraph*{Construction of the auxilliary metrics $\tilde{h}$, $\tilde{h}_{tim}$}

Since $h$ is singular on $\partial_{hor}\mathcal{S}$, it will be
useful to have a second Riemannian metric $\tilde{h}$ on $\mathcal{S}$
that is smooth up to $\partial_{hor}\mathcal{S}$,%
\footnote{$h$ is already smooth up to $\partial_{tim}\mathcal{S}$%
} so that we can measure the norms of tensors with the use of $\tilde{h}$.
The metric $\tilde{h}$ will also be used to define covariant derivative
operators and geometric volume forms which are regular up to $\partial_{hor}\mathcal{S}$
(where the associated constructions with $h$ will either be singular
or degenerate). This metric will model the metric $h_{\text{\textgreek{t}},N}$
on the hypersurfaces $\{\bar{t}=\text{\textgreek{t}}\}$ of the spacetimes
$(\mathcal{M},g)$ of Sections \ref{sec:Firstdecay} and \ref{sec:Improved-polynomial-decay}.
To this end, we define a Riemannian metric $\tilde{h}$ on $\mathcal{S}$
such that $\tilde{h}\equiv h$ on $\{r\ge1\}$, and 
\begin{equation}
\tilde{h}=dr^{2}+h_{tan}\label{eq:DistortedMetricNearHorizon}
\end{equation}
 in the region $r\ll1$. 

It will be convenient to define the smooth functions 
\begin{equation}
r_{+}=(1+r^{2})^{1/2}.
\end{equation}
 and 
\begin{equation}
r_{-}=\big(1+\frac{1}{r}\big)^{-1}.
\end{equation}
 Notice that $r_{+}\sim r$ for $r\gg1$ and $r_{+}\sim1$ for $r\lesssim1$,
while $r_{-}\sim r$ near $\partial\mathcal{S}$ and $r_{-}\sim1$
away from $\partial\mathcal{S}$. 

We also define the smooth funcions $r_{hor},r_{tim}:\mathcal{S}\rightarrow[0,1]$
by the following requirements 
\begin{itemize}
\item $r_{hor}=r$ in the region $\{dist_{\tilde{h}}(\cdot,\partial_{hor}\mathcal{S})\ll1\}$
(where $dist_{\tilde{h}}(\cdot,\partial_{hor}\mathcal{S})$ is smooth),
$r_{hor}>0$ on $\mathcal{S}\backslash\partial_{hor}\mathcal{S}$
and $r_{hor}\equiv1$ on $\{dist_{\tilde{h}}(\cdot,\partial_{hor}\mathcal{S})\ge1\}$. 
\item $r_{tim}=r$ in the region $\{dist_{\tilde{h}}(\cdot,\partial_{tim}\mathcal{S})\ll1\}$,
$r_{tim}>0$ on $\mathcal{S}\backslash\partial_{tim}\mathcal{S}$
and $r_{tim}\equiv1$ on $\{dist_{\tilde{h}}(\cdot,\partial_{tim}\mathcal{S})\ge1\}$. 
\end{itemize}
Finally, define the Riemannian metric $\tilde{h}_{tim}$ on $\mathcal{S}\backslash\partial_{tim}\mathcal{S}$
so that $\tilde{h}_{tim}\equiv\tilde{h}$ on $\{r_{tim}\ge\frac{1}{2}\}$
and 
\begin{equation}
\tilde{h}_{tim}=r^{-2}dr^{2}+\big(-\log(r_{tim})\big)\cdot h_{tan}\label{eq:TimelikeMetric}
\end{equation}
 in the region $\{r_{tim}\ll1\}$. This metric will only be used to
handle difficulties appearing near $\partial_{tim}\mathcal{S}$ in
the derivation of the elliptic estimates in this Section, and is not
associated with any geometric construction performed in Sections \ref{sec:Firstdecay}
and \ref{sec:Improved-polynomial-decay}.
\begin{rem*}
We will raise and lower indices only with the use of the singular
metric $h$. The non-singular metric $\tilde{h}$ will only be used
to measure norms of tensors on $\mathcal{S}$. Covariant derivatives
with respect to $h$ will be simply denoted by $\nabla$, while the
ones associated to $\tilde{h}$ and $\tilde{h}_{tim}$ will be denoted
by $\nabla^{(\tilde{h})}$ and $\nabla^{(\tilde{h}_{tim})}$ respectively.
\end{rem*}

\subsection{Elliptic estimates on $(\mathcal{S},h)$}

In this Section, we will establish elliptic estimates on the Riemannian
manifolds $(\mathcal{S},h)$ associated to the elliptic operators
(\ref{eq:PerturbedLaplacian-1}) and (\ref{eq:NonDegeneratePerturbedLaplacian}).
To this end, we will need some shorthand notation for weighted norms
of derivatives of functions $u$ on $\mathcal{S}$, with weights which
are either regular or degenerate on parts of the boundary $\partial\mathcal{S}$.
These weights are tied naturally to the use of differential operators
associated to the different Riemannian metrics $h$, $\tilde{h}$
and $\tilde{h}_{tim}$ that we have already defined on $\mathcal{S}\backslash\partial\mathcal{S}$.
With this motivation, after fixing some smooth cut-off functions $\text{\textgreek{q}}_{\ge r_{1}},\text{\textgreek{q}}_{\le r_{1}}:\mathcal{S}\rightarrow[0,1]$
so that $\text{\textgreek{q}}_{\ge r_{1}}\equiv1$ on $\{r\ge r_{1}\}$
and $\text{\textgreek{q}}_{\ge r_{1}}\equiv0$ on $\{r\le r_{1}/2\}$
and $\text{\textgreek{q}}_{\le r_{1}}=1-\text{\textgreek{q}}_{\ge r_{1}}$,
we introduce the following definition:
\begin{defn*}
We will introduce the following pointwise norm on $\mathcal{S}\backslash\partial\mathcal{S}$
for any pair of Riemannian metrics $h_{1}$, $h_{2}$ on $\mathcal{S}\backslash\partial\mathcal{S}$
and any integer $m\ge2$: 
\begin{align}
|u|_{h_{1},h_{2};m}^{2}\doteq & \text{\textgreek{q}}_{\ge r_{1}}\big|\big(\nabla^{(h_{1})}\big)^{m}u\big|_{h_{2}}^{2}+\text{\textgreek{q}}_{\le r_{1}}\cdot r_{hor}^{2}\big|\big(\nabla^{(h_{1})}\big)^{m-2}(Y^{2}u)\big|_{h_{2}}^{2}+\text{\textgreek{q}}_{\le r_{1}}\cdot r_{hor}\big|\big(\nabla^{(h_{1})}\big)^{m-2}\big(i_{*}\nabla^{(h_{tan})}(Yu)\big)\big|_{h_{2}}^{2}+\label{eq:DefinitionNormPointwise}\\
 & +\text{\textgreek{q}}_{\le r_{1}}\cdot\big|\big(\nabla^{(h_{1})}\big)^{m-2}\big((i_{*}\nabla^{(h_{tan})})^{2}u\big)\big|_{h_{2}}^{2}.\nonumber 
\end{align}
In the above the $i_{*}$ notation is used as follows: For any $(0,k)$-tensor
$\mathfrak{m}$ on the surfaces $\{r=const\}$ on $\{r\le r_{1}\}$
varying smoothly with $r$ we denote with $i_{*}\mathfrak{m}$ the
unique tensor on $\{r\le r_{1}\}$ such that for any local frame $\{E_{0},E_{1},\ldots E_{d-1}\}$
on $\{r\le r_{1}\}$ with $E_{0}=Y$ and $E_{1}r=\ldots=E_{d-1}r=0$
and any $i_{1},\ldots i_{k}\in\{0,\ldots d-1\}$, the following relation
holds: 
\begin{equation}
i_{*}\mathfrak{m}(E_{i_{1}},\ldots E_{i_{k}})=\begin{cases}
0, & \mbox{ if some of the }i_{k}\mbox{'s is }0\\
\mathfrak{m}(E_{i_{1}},\ldots E_{i_{k}}) & \mbox{ if none of the }i_{k}\mbox{'s is }0.
\end{cases}
\end{equation}
 \end{defn*}
\begin{rem*}
In most instances where the notation (\ref{eq:DefinitionNormPointwise})
will appear, $h_{1}$ will be the everywhere regular metric $\tilde{h}$,
while $h_{2}$ will be a metric which is singular on $\partial_{tim}\mathcal{S}$. 
\end{rem*}
We will establish the following lemma: 
\begin{prop}
\label{prop:DegenerateEllipticEstimates}\textbf{(Degenerate elliptic estimates).}
For any $l\in\mathbb{N}$ with $2\le l\le\lfloor\frac{d+1}{2}\rfloor$,
any $k_{0}\in\mathbb{N}$ and any $\text{\textgreek{b}}\in(-\bar{\text{\textgreek{d}}}_{k_{0}},1)$
for some $\bar{\text{\textgreek{d}}}_{k_{0}}>0$ depending on $k_{0}$,
there exists a (small) $r_{0}>0$ so that we can bound for any $u\in C^{\infty}(\mathcal{S})$
with $\limsup_{r\rightarrow+\infty}\big|r^{\frac{d-1}{2}+j}\nabla^{j}u\big|_{h}<+\infty$
for any $j\le l+k_{0}$: 
\begin{equation}
\begin{split}\sum_{k=0}^{k_{0}}\Big\{\int_{\mathcal{S}}r_{+}^{-\text{\textgreek{b}}}|u|_{\tilde{h},\tilde{h}_{tim};k+l}^{2} & \, dvol_{\tilde{h}}+\sum_{j=1}^{l-1}\int_{\{r\ge r_{0}\}}r_{+}^{-2j-\text{\textgreek{b}}}\cdot|\big(\nabla^{(\tilde{h})}\big)^{k+l-j}u|_{\tilde{h}}^{2}\, dvol_{\tilde{h}}+\\
\hphantom{\sum_{k=0}^{k_{0}}\Big\{}\sum_{j=1}^{l-1}\int_{\{r\le r_{0}\}}\Big(|\big(\nabla^{(\tilde{h})} & \big)^{k+l-j-1}(Yu)|_{\tilde{h}_{tim}}^{2}+r_{-}^{-1}\log^{-2}(\frac{r_{-}}{2})\cdot|\big(\nabla^{(\tilde{h})}\big)^{k+l-j-1}\big(i_{*}(\nabla^{(h_{tan})}u)\big)|_{\tilde{h}_{tim}}^{2}\Big)\, dvol_{\tilde{h}}\Big\}\le\\
\le & C_{\text{\textgreek{b}},k_{0}}\sum_{k=0}^{k_{0}}\int_{\mathcal{S}}r_{+}^{-\text{\textgreek{b}}}\big|\big(\nabla^{(\tilde{h})}\big)^{k+l-2}(\text{\textgreek{D}}_{h,\text{\textgreek{w}}}u)\big|_{\tilde{h}_{tim}}^{2}\, dvol_{\tilde{h}}+\\
 & +C_{\text{\textgreek{b}},k_{0}}\sum_{j=0}^{1}\max\Big\{-Re\big\{\int_{\partial_{tim}\mathcal{S}}h_{tan}\Big(\big(\nabla^{(h_{tan})}\big)^{j}(Yu),\big(\nabla^{(h_{tan})}\big)^{j}\bar{u}\Big)\, dh_{tan}\big\},0\Big\}
\end{split}
\label{eq:EllipticEstimatesDegenerateAtTimelike}
\end{equation}
 and 
\begin{equation}
\begin{split}\sum_{k=0}^{k_{0}}\Big\{\int_{\mathcal{S}}r_{+}^{-\text{\textgreek{b}}}|u|_{\tilde{h},\big(1-\log(r_{tim})\big)\cdot\tilde{h};k+l}^{2} & \, dvol_{\tilde{h}}+\sum_{j=1}^{l-1}\int_{\{r\ge r_{0}\}}r_{+}^{-2j-\text{\textgreek{b}}}\cdot|\big(\nabla^{(\tilde{h})}\big)^{k+l-j}u|_{\tilde{h}}^{2}\, dvol_{\tilde{h}}+\\
\hphantom{\sum_{k=0}^{k_{0}}}+\sum_{j=1}^{l-1}\int_{\{r\le r_{0}\}}\Big(|\big(\nabla^{(\tilde{h})}\big)^{k+l-j-1} & (Yu)|_{\big(1-\log(r_{tim})\big)\tilde{h}}^{2}+r_{-}^{-1}\log^{-2}(\frac{r_{-}}{2})\cdot|\big(\nabla^{(\tilde{h})}\big)^{k+l-j-1}\big(i_{*}(\nabla^{(h_{tan})}u)\big)|_{\big(1-\log(r_{tim})\big)\tilde{h}}^{2}\Big)\, dvol_{\tilde{h}}\Big\}\le\\
\le & C_{\text{\textgreek{b}},k_{0}}\sum_{k=0}^{k_{0}}\int_{\mathcal{S}}r_{+}^{-\text{\textgreek{b}}}\big|\big(\nabla^{(\tilde{h})}\big)^{k+l-2}(\text{\textgreek{D}}_{h,\text{\textgreek{w}}}u)\big|_{\big(1-\log(r_{tim})\big)\cdot\tilde{h}}^{2}\, dvol_{\tilde{h}}+\\
 & +C_{\text{\textgreek{b}},k_{0}}\sum_{j=0}^{1}\max\Big\{-Re\big\{\int_{\partial_{tim}\mathcal{S}}h_{tan}\Big(\big(\nabla^{(h_{tan})}\big)^{j}(Yu),\big(\nabla^{(h_{tan})}\big)^{j}\bar{u}\Big)\, dh_{tan}\big\},0\Big\},
\end{split}
\label{eq:EllipticEstimatesDegenerateRiemannian}
\end{equation}
 where the constant $C_{\text{\textgreek{b}},k_{0}}$ of the right
hand side depends only on $\text{\textgreek{b}}$, $k_{0}$ and on
the geometry of $(\mathcal{S},h)$, $\tilde{h}$ and $\text{\textgreek{w}}$. \end{prop}
\begin{rem*}
Notice that the boundary terms on the right hand side of (\ref{eq:EllipticEstimatesDegenerateAtTimelike})
and (\ref{eq:EllipticEstimatesDegenerateRiemannian}) contain only
terms on the $\partial_{tim}\mathcal{S}$ part of the boundary. These
terms vanish when $u$ satisfies the Dirichlet or Neumann conditions
on $\partial_{tim}\mathcal{S}$.\end{rem*}
\begin{proof}
Withot loss of generality, we can assume that $u$ is real valued.
We will treat $r_{0}$ and $\bar{\text{\textgreek{d}}}_{k_{0}}$ as
small parameters (which will be fixed later in the proof). 

The region $\{r\le2r_{0}\}$ will be identified with $[0,2r_{0}]\times\partial\mathcal{S}$
through the diffeomorphism $\mathcal{J}$. We will assume that $r_{0}$
is small enough so that $[0,2r_{0}]\times\partial_{hor}\mathcal{S}$
and $[0,2r_{0}]\times\partial_{tim}\mathcal{S}$ are disjoint. Let
us define two smooth cut-off functions $\text{\textgreek{q}}_{\le r_{0}}^{hor},\text{\textgreek{q}}_{\le r_{0}}^{tim}:\mathcal{S}\rightarrow[0,1]$
such that

1. $supp(\text{\textgreek{q}}_{\le r_{0}}^{hor})\subseteq[0,2r_{0}]\times\partial_{hor}\mathcal{S}$
and $\text{\textgreek{q}}_{\le r_{0}}^{hor}\equiv1$ on $[0,r_{0}]\times\partial_{hor}\mathcal{S}$
and

2. $supp(\text{\textgreek{q}}_{\le r_{0}}^{tim})\subseteq[0,2r_{0}]\times\partial_{tim}\mathcal{S}$
and $\text{\textgreek{q}}_{\le r_{0}}^{tim}\equiv1$ on $[0,r_{0}]\times\partial_{tim}\mathcal{S}$. 

\noindent We will also set $\text{\textgreek{q}}_{\ge r_{0}}=1-\text{\textgreek{q}}_{\le r_{0}}^{hor}-\text{\textgreek{q}}_{\le r_{0}}^{tim}$.
Let us also define $\text{\textgreek{q}}_{\ge R_{0}}=1-\text{\textgreek{q}}(\frac{r}{R_{0}})$
for some $R_{0}$ large and fixed in terms of $\text{\textgreek{b}}$
and the geometry of $(\mathcal{S},h)$.

Without loss of generality, we can assume that $k_{0}=0$ (since the
proof in the case $k_{0}\ge1$ follows in exactly the same way). In
order to present the necessary ideas in a simpler form, we will first
establish the case when $l=2$. 

\smallskip{}

\noindent \emph{Remark. }Notice that for $l=2$, it is most difficult
to treat the case $d=3$, since in that case one is not able to obtain
control over weighted $L^{2}$ norms of $u$ using only the finiteness
of $\int_{\mathcal{S}}r_{+}^{-\text{\textgreek{b}}}|\nabla^{2}u|^{2}\, dvol_{\tilde{h}}$
together with Hardy and Poincare inequalities. The same difficulties
occur generally in dimension $d$ for $l=\lfloor\frac{d+1}{2}\rfloor$.

\smallskip{}

We will start by establishing the following elliptic estimate on each
connected component of the far away region $\{r\gg1\}$: 
\begin{align}
\int_{\mathcal{S}}\text{\textgreek{q}}_{\ge R_{0}}r^{-\text{\textgreek{b}}}\Big(r^{-(d-1)}\partial_{r}\big(r^{d-1}\partial_{r}u\big)+r^{-2}\text{\textgreek{D}}_{\mathbb{S}^{d-1}}u\Big)^{2}\, dvol_{h}\ge & c_{\text{\textgreek{b}}}\int_{\mathcal{S}}\text{\textgreek{q}}_{\ge R_{0}}r^{-\text{\textgreek{b}}}\Big(|\nabla^{2}u|_{h}^{2}+r^{-2}|\nabla u|_{h}^{2}\Big)\, dvol_{h}+\label{eq:FinalEllipticEstimateAway-1}\\
 & +\int_{\mathcal{S}}O(|\nabla\text{\textgreek{q}}_{\ge R_{0}}|+|\nabla^{2}\text{\textgreek{q}}_{\ge R_{0}}|)|\nabla u|_{h}^{2}\, dvol_{h}.\nonumber 
\end{align}
 Without loss of generality, it suffices to establish (\ref{eq:FinalEllipticEstimateAway-1})
in the case $0<1-\text{\textgreek{b}}\ll1$, since the case $\text{\textgreek{b}}\in(-\bar{\text{\textgreek{d}}}_{k_{0}},0]$
(provided $\bar{\text{\textgreek{d}}}_{k_{0}}$ is small enough) follows
by a straightforward integration by parts scheme (and thus (\ref{eq:FinalEllipticEstimateAway-1})
for intermediate values of $\text{\textgreek{b}}$ will follow by
an easy interpolation argument).

On each connected component of the region $\{r\ge R_{0}\}$ in the
coordinate chart $(r,\text{\textgreek{sv}})$ we calculate 
\begin{equation}
\text{\textgreek{D}}_{h,\text{\textgreek{w}}}u=r^{-(d-1)}\partial_{r}\big(r^{d-1}\partial_{r}u\big)+r^{-2}\text{\textgreek{D}}_{\mathbb{S}^{d-1}}u+O_{\text{\textgreek{m}\textgreek{n}}}(r^{-1})\nabla^{\text{\textgreek{m}}}\nabla^{\text{\textgreek{n}}}u+O_{\text{\textgreek{m}}}(r^{-2})\nabla^{\text{\textgreek{m}}}u.\label{eq:AsymptoticLaplacian}
\end{equation}
 Assuming (\ref{eq:FinalEllipticEstimateAway-1}) has been established,
from (\ref{eq:FinalEllipticEstimateAway-1}) and (\ref{eq:AsymptoticLaplacian})
the following estimate is readily deduced (provided that $R_{0}$
has been fixed large in terms of $\text{\textgreek{b}}$ and the geometry
of $(\mathcal{S},h)$): 
\begin{align}
\int_{\mathcal{S}}\text{\textgreek{q}}_{\ge R_{0}}r^{-\text{\textgreek{b}}}\big(\text{\textgreek{D}}_{h,\text{\textgreek{w}}}u\big)^{2}\, dvol_{h}\ge & c_{\text{\textgreek{b}}}\int_{\mathcal{S}}\text{\textgreek{q}}_{\ge R_{0}}r^{-\text{\textgreek{b}}}\Big(|\nabla^{2}u|_{h}^{2}+r^{-2}|\nabla u|_{h}^{2}\Big)\, dvol_{h}+\label{eq:FinalEllipticEstimateAway-1-1}\\
 & +\int_{\mathcal{S}}O(|\nabla\text{\textgreek{q}}_{\ge R_{0}}|+|\nabla^{2}\text{\textgreek{q}}_{\ge R_{0}}|)|\nabla u|_{h}^{2}\, dvol_{h}.\nonumber 
\end{align}

We will now establish (\ref{eq:FinalEllipticEstimateAway-1}). By
expanding the square we can trivially infer for any $a_{s}\in[0,1]$
(which will be fixed later): 
\begin{equation}
\begin{split}\int_{\mathcal{S}}\text{\textgreek{q}}_{\ge R_{0}}r^{-\text{\textgreek{b}}} & \Big(r^{-(d-1)}\partial_{r}\big(r^{d-1}\partial_{r}u\big)+r^{-2}\text{\textgreek{D}}_{\mathbb{S}^{d-1}}u\Big)^{2}\, r^{d-1}drd\text{\textgreek{sv}}=\\
= & a_{s}\int_{\mathcal{S}}\text{\textgreek{q}}_{\ge R_{0}}r^{-\text{\textgreek{b}}}\Big\{\big(r^{-(d-1)}\partial_{r}\big(r^{d-1}\partial_{r}u\big)\big)^{2}+2r^{-(d-1)}\partial_{r}\big(r^{d-1}\partial_{r}u\big)\cdot r^{-2}\text{\textgreek{D}}_{\mathbb{S}^{d-1}}u+r^{-4}\big(\text{\textgreek{D}}_{\mathbb{S}^{d-1}}u\big)^{2}\Big\}\, r^{d-1}drd\text{\textgreek{sv}}+\\
 & +(1-a_{s})\int_{\mathcal{S}}\text{\textgreek{q}}_{\ge R_{0}}r^{-\text{\textgreek{b}}}\Big\{\big(r^{-(d-1)}\partial_{r}\big(r^{d-1}\partial_{r}u\big)\big)^{2}+2r^{-(d-1)}\partial_{r}\big(r^{d-1}\partial_{r}u\big)\cdot r^{-2}\text{\textgreek{D}}_{\mathbb{S}^{d-1}}u+r^{-4}\big(\text{\textgreek{D}}_{\mathbb{S}^{d-1}}u\big)^{2}\Big\}\, r^{d-1}drd\text{\textgreek{sv}}.
\end{split}
\label{eq:CriticalEstimateFirstStep-2}
\end{equation}
 Using a Cauchy--Schwarz inequality 
\begin{equation}
2a_{s}r^{-(d-1)}\partial_{r}\big(r^{d-1}\partial_{r}u\big)\cdot r^{-2}\text{\textgreek{D}}_{\mathbb{S}^{d-1}}u\ge-\frac{a_{s}^{2}}{b_{s}}\big(r^{-(d-1)}\partial_{r}\big(r^{d-1}\partial_{r}u\big)\big)^{2}-b_{s}r^{-4}\big(\text{\textgreek{D}}_{\mathbb{S}^{d-1}}u\big)^{2}
\end{equation}
 (for a parameter $b_{s}>0$ to be fixed later) in the first summand
of the right hand side of (\ref{eq:CriticalEstimateFirstStep-2}),
and then applying the product rule and an integration by parts in
the $\partial_{r}$ and the angular directions for both summands,
we obtain from (\ref{eq:CriticalEstimateFirstStep-2}): 
\begin{equation}
\begin{split}\int_{\mathcal{S}}\text{\textgreek{q}}_{\ge R_{0}}r^{-\text{\textgreek{b}}}\Big(r^{-(d-1)} & \partial_{r}\big(r^{d-1}\partial_{r}u\big)+r^{-2}\text{\textgreek{D}}_{\mathbb{S}^{d-1}}u\Big)^{2}\, r^{d-1}drd\text{\textgreek{sv}}\ge\\
\ge & \int_{\mathcal{S}}\text{\textgreek{q}}_{\ge R_{0}}r^{-\text{\textgreek{b}}}\Big\{\big(1-\frac{a_{s}^{2}}{b_{s}}\big)\Big(\big(\partial_{r}^{2}u\big)^{2}+(d-1)(1+\text{\textgreek{b}})r^{-2}\big(\partial_{r}u\big)^{2}\Big)+\\
 & \hphantom{\int_{\mathcal{N}}\text{\textgreek{q}}_{\ge R_{0}}r^{-\text{\textgreek{b}}}\Big\{}+2(1-a_{s})r^{-2}\big|\nabla_{\mathbb{S}^{d-1}}\partial_{r}u\big|_{g_{\mathbb{S}^{d-1}}}^{2}+\big(1-b_{s}\big)r^{-4}\big(\text{\textgreek{D}}_{\mathbb{S}^{d-1}}u\big)^{2}\Big\}\, r^{d-1}drd\text{\textgreek{sv}}-\\
 & -\int_{\mathcal{S}}\text{\textgreek{q}}_{\ge R_{0}}2(1-a_{s})(2+\text{\textgreek{b}})r^{-3-\text{\textgreek{b}}}\cdot\partial_{r}u\cdot\text{\textgreek{D}}_{\mathbb{S}^{d-1}}u\, r^{d-1}drd\text{\textgreek{sv}}+\\
 & +\int_{\mathcal{S}}O(|\nabla\text{\textgreek{q}}_{\ge R_{0}}|)\Big((\partial_{r}u)^{2}+\partial_{r}u\cdot\text{\textgreek{D}}_{\mathbb{S}^{d-1}}u\Big)\, r^{d-1}drd\text{\textgreek{sv}}.
\end{split}
\label{eq:CriticalEstimateFirstStep}
\end{equation}
 Using the Hardy type inequality 
\begin{equation}
\int_{\mathcal{S}}\text{\textgreek{q}}_{\ge R_{0}}r^{d-3-\text{\textgreek{b}}}\big(\partial_{r}u\big)^{2}\, drd\text{\textgreek{sv}}\le\frac{4}{(d-2-\text{\textgreek{b}})^{2}}\int_{\mathcal{S}}\text{\textgreek{q}}_{\ge R_{0}}r^{d-1-\text{\textgreek{b}}}\big(\partial_{r}^{2}u\big)^{2}\, drd\text{\textgreek{sv}}+C_{\text{\textgreek{b}}}\int_{\mathcal{S}}O(|\nabla\text{\textgreek{q}}_{\ge R_{0}}|)(\partial_{r}u)^{2}\, r^{d-1}drd\text{\textgreek{sv}},
\end{equation}
we obtain from (\ref{eq:CriticalEstimateFirstStep}) for any $0<\text{\textgreek{d}}_{H}<1$
(also to be fixed later): 
\begin{equation}
\begin{split}\int_{\mathcal{\mathcal{S}}}\text{\textgreek{q}}_{\ge R_{0}}r^{-\text{\textgreek{b}}}\Big(r^{-(d-1)} & \partial_{r}\big(r^{d-1}\partial_{r}u\big)+r^{-2}\text{\textgreek{D}}_{\mathbb{S}^{d-1}}u\Big)^{2}\, r^{d-1}drd\text{\textgreek{sv}}\ge\\
\ge & \int_{\mathcal{S}}\text{\textgreek{q}}_{\ge R_{0}}r^{-\text{\textgreek{b}}}\Big\{\big(1-\frac{a_{s}^{2}}{b_{s}}\big)\Big(\text{\textgreek{d}}_{H}\big(\partial_{r}^{2}u\big)^{2}+\big\{\big(d-1-\frac{d-2-\text{\textgreek{b}}}{2}\big)^{2}+O(\text{\textgreek{d}}_{H})\Big\} r^{-2}\big(\partial_{r}u\big)^{2}\Big)+\\
 & \hphantom{\int_{\mathcal{N}}\text{\textgreek{q}}_{\ge R_{0}}r^{-\text{\textgreek{b}}}\Big\{}+2(1-a_{s})r^{-2}\big|\nabla_{\mathbb{S}^{d-1}}\partial_{r}u\big|_{g_{\mathbb{S}^{d-1}}}^{2}+\big(1-b_{s}\big)r^{-4}\big(\text{\textgreek{D}}_{\mathbb{S}^{d-1}}u\big)^{2}\Big\}\, r^{d-1}drd\text{\textgreek{sv}}-\\
 & -\int_{\mathcal{S}}\text{\textgreek{q}}_{\ge R_{0}}2(1-a_{s})(2+\text{\textgreek{b}})r^{-3-\text{\textgreek{b}}}\cdot\partial_{r}u\cdot\text{\textgreek{D}}_{\mathbb{S}^{d-1}}u\, r^{d-1}drd\text{\textgreek{sv}}+\\
 & +\int_{\mathcal{S}}O(|\nabla\text{\textgreek{q}}_{\ge R_{0}}|)\Big((\partial_{r}u)^{2}+\partial_{r}u\cdot\text{\textgreek{D}}_{\mathbb{S}^{d-1}}u\Big)\, r^{d-1}drd\text{\textgreek{sv}}.
\end{split}
\label{eq:CriticalEstimateFirstStep-3}
\end{equation}

Setting 
\begin{equation}
u_{0}\doteq\frac{1}{Vol(\mathbb{S}^{d-1})}\int_{\mathbb{S}^{d-1}}u\, d\text{\textgreek{sv}}\label{eq:SphericallySymmetric}
\end{equation}
 and 
\begin{equation}
u_{\ge1}=u-u_{0},\label{eq:NonSphericallySymmetric}
\end{equation}
 and noting that $u_{0}$ and $u_{\ge1}$ are orthogonal with respect
to the $L^{2}(dg_{\mathbb{S}^{d-1}})$ inner product, we obtain from
(\ref{eq:CriticalEstimateFirstStep-3}): 
\begin{equation}
\begin{split}\int_{\mathcal{S}}\text{\textgreek{q}}_{\ge R_{0}}r^{-\text{\textgreek{b}}}\Big(r^{-(d-1)}\partial_{r}\big(r^{d-1} & \partial_{r}u\big)+r^{-2}\text{\textgreek{D}}_{\mathbb{S}^{d-1}}u\Big)^{2}\, r^{d-1}drd\text{\textgreek{sv}}\ge\big(1-\frac{a_{s}^{2}}{b_{s}}\big)\text{\textgreek{d}}_{H}\cdot\int_{\mathcal{S}}\text{\textgreek{q}}_{\ge R_{0}}r^{-\text{\textgreek{b}}}\big(\partial_{r}^{2}u\big)^{2}\, r^{d-1}drd\text{\textgreek{sv}}+\\
 & +\int_{\mathcal{S}}\text{\textgreek{q}}_{\ge R_{0}}r^{-\text{\textgreek{b}}}\Big\{\big(1-\frac{a_{s}^{2}}{b_{s}}\big)\Big(\big\{\big(d-1-\frac{d-2-\text{\textgreek{b}}}{2}\big)^{2}+O(\text{\textgreek{d}}_{H})\Big\} r^{-2}\big(\partial_{r}u_{\ge1}\big)^{2}\Big)+\\
 & \hphantom{+\int_{\mathcal{S}}\text{\textgreek{q}}_{\ge R_{0}}r^{-\text{\textgreek{b}}}\Big\{}+2(1-a_{s})r^{-2}\big|\nabla_{\mathbb{S}^{d-1}}\partial_{r}u_{\ge1}\big|_{g_{\mathbb{S}^{d-1}}}^{2}+\big(1-b_{s}\big)r^{-4}\big(\text{\textgreek{D}}_{\mathbb{S}^{d-1}}u_{\ge1}\big)^{2}\Big\}\, r^{d-1}drd\text{\textgreek{sv}}-\\
 & -\int_{\mathcal{S}}\text{\textgreek{q}}_{\ge R_{0}}2(1-a_{s})(2+\text{\textgreek{b}})r^{-3-\text{\textgreek{b}}}\cdot\partial_{r}u_{\ge1}\cdot\text{\textgreek{D}}_{\mathbb{S}^{d-1}}u_{\ge1}\, r^{d-1}drd\text{\textgreek{sv}}+\\
 & +\int_{\mathcal{S}}O(|\nabla\text{\textgreek{q}}_{\ge R_{0}}|)\Big((\partial_{r}u)^{2}+\partial_{r}u\cdot\text{\textgreek{D}}_{\mathbb{S}^{d-1}}u\Big)\, r^{d-1}drd\text{\textgreek{sv}}.
\end{split}
\label{eq:CriticalEstimateFirstStep-3-1}
\end{equation}
 Since the first non zero eigenvalue of $\text{\textgreek{D}}_{\mathbb{S}^{d-1}}$
equals $d-1$, we can bound (in view of (\ref{eq:SphericallySymmetric})
and (\ref{eq:NonSphericallySymmetric})): 
\begin{equation}
\int_{\mathbb{S}^{d-1}}\big|\nabla_{\mathbb{S}^{d-1}}\partial_{r}u_{\ge1}\big|_{g_{\mathbb{S}^{d-1}}}^{2}\, ds\ge(d-1)\int_{\mathbb{S}^{d-1}}(\partial_{r}u_{\ge1})^{2}\, ds
\end{equation}
and thus, since $d\ge3$, from (\ref{eq:CriticalEstimateFirstStep-3-1})
we obtain after setting $\text{\textgreek{e}}_{\text{\textgreek{b}}}\doteq1-\text{\textgreek{b}}>0$:
\begin{equation}
\begin{split}\int_{\mathcal{S}}\text{\textgreek{q}}_{\ge R_{0}}r^{-\text{\textgreek{b}}}\Big(r^{-(d-1)} & \partial_{r}\big(r^{d-1}\partial_{r}u\big)+r^{-2}\text{\textgreek{D}}_{\mathbb{S}^{d-1}}u\Big)^{2}\, r^{d-1}drd\text{\textgreek{sv}}\ge c(a_{s},b_{s},\text{\textgreek{d}}_{H})\int_{\mathcal{S}}\text{\textgreek{q}}_{\ge R_{0}}r^{-\text{\textgreek{b}}}\big(\partial_{r}^{2}u\big)^{2}\, r^{d-1}drd\text{\textgreek{sv}}+\\
 & +\int_{\mathcal{S}}\text{\textgreek{q}}_{\ge R_{0}}r^{-\text{\textgreek{b}}}\Big\{\big(1-\frac{a_{s}^{2}}{b_{s}}\big)\big(2-\frac{\text{\textgreek{e}}_{\text{\textgreek{b}}}}{2}\big)^{2}+4(1-a_{s})+O(\text{\textgreek{d}}_{H})\Big\} r^{-2}\big(\partial_{r}u_{\ge1}\big)^{2}+\big(1-b_{s}\big)r^{-4}\big(\text{\textgreek{D}}_{\mathbb{S}^{d-1}}u_{\ge1}\big)^{2}\Big\}\, r^{d-1}drd\text{\textgreek{sv}}-\\
 & -\int_{\mathcal{S}}\text{\textgreek{q}}_{\ge R_{0}}2(1-a_{s})(3-\text{\textgreek{e}}_{\text{\textgreek{b}}})r^{-3-\text{\textgreek{b}}}\cdot\partial_{r}u_{\ge1}\cdot\text{\textgreek{D}}_{\mathbb{S}^{d-1}}u_{\ge1}\, r^{d-1}drd\text{\textgreek{sv}}+\\
 & +c(a_{s},b_{s})\int_{\mathcal{S}}\text{\textgreek{q}}_{\ge R_{0}}r^{-2-\text{\textgreek{b}}}\big(\partial_{r}u_{0}\big)^{2}\Big)\, r^{d-1}drd\text{\textgreek{sv}}+\int_{\mathcal{S}}O(|\nabla\text{\textgreek{q}}_{\ge R_{0}}|)\Big((\partial_{r}u)^{2}+\partial_{r}u\cdot\text{\textgreek{D}}_{\mathbb{S}^{d-1}}u\Big)\, r^{d-1}drd\text{\textgreek{sv}}.
\end{split}
\label{eq:CriticalEstimateFirstStep-3-1-1}
\end{equation}

Using a Cauchy--Schwarz inequality for the third line of the right
hand side of (\ref{eq:CriticalEstimateFirstStep-3-1-1}): 
\begin{align}
\Big|\int_{\mathcal{S}}\text{\textgreek{q}}_{\ge R_{0}}2(1-a_{s})(3-\text{\textgreek{e}}_{\text{\textgreek{b}}})r^{-3-\text{\textgreek{b}}}\cdot\partial_{r}u_{\ge1}\cdot\text{\textgreek{D}}_{\mathbb{S}^{d-1}}u_{\ge1}\, r^{d-1}drd\text{\textgreek{sv}}\Big|\le & \int_{\mathcal{S}}\text{\textgreek{q}}_{\ge R_{0}}(1-a_{s})^{2}\frac{(3-\text{\textgreek{e}}_{\text{\textgreek{b}}})^{2}}{(1-b_{s}-\text{\textgreek{d}}_{H})}r^{-2-\text{\textgreek{b}}}\big(\partial_{r}u_{\ge1}\big)^{2}\, r^{d-1}drd\text{\textgreek{sv}}+\\
 & +\int_{\mathcal{S}}\text{\textgreek{q}}_{\ge R_{0}}(1-b_{s}-\text{\textgreek{d}}_{H})\cdot r^{-4-\text{\textgreek{b}}}\big(\text{\textgreek{D}}_{\mathbb{S}^{d-1}}u_{\ge1}\big)^{2}\, r^{d-1}drd\text{\textgreek{sv}},\nonumber 
\end{align}
we infer from (\ref{eq:CriticalEstimateFirstStep-3-1-1}): 
\begin{equation}
\begin{split}\int_{\mathcal{S}}\text{\textgreek{q}}_{\ge R_{0}}r^{-\text{\textgreek{b}}}\Big(r^{-(d-1)}\partial_{r}\big(r^{d-1}\partial_{r}u\big) & +r^{-2}\text{\textgreek{D}}_{\mathbb{S}^{d-1}}u\Big)^{2}\, r^{d-1}drd\text{\textgreek{sv}}\ge c(a_{s},b_{s},\text{\textgreek{d}}_{H})\int_{\mathcal{S}}\text{\textgreek{q}}_{\ge R_{0}}r^{-\text{\textgreek{b}}}\big(\partial_{r}^{2}u\big)^{2}\, r^{d-1}drd\text{\textgreek{sv}}+\\
 & +\int_{\mathcal{S}}\text{\textgreek{q}}_{\ge R_{0}}r^{-\text{\textgreek{b}}}\Big\{ A_{co}\cdot r^{-2}\big(\partial_{r}u_{\ge1}\big)^{2}+\text{\textgreek{d}}_{H}r^{-4}\big(\text{\textgreek{D}}_{\mathbb{S}^{d-1}}u_{\ge1}\big)^{2}\Big\}\, r^{d-1}drd\text{\textgreek{sv}}+\\
 & +c(a_{s},b_{s})\int_{\mathcal{S}}\text{\textgreek{q}}_{\ge R_{0}}r^{-2-\text{\textgreek{b}}}\big(\partial_{r}u_{0}\big)^{2}\Big)\, r^{d-1}drd\text{\textgreek{sv}}+\int_{\mathcal{S}}O(|\nabla\text{\textgreek{q}}_{\ge R_{0}}|)\Big((\partial_{r}u)^{2}+\partial_{r}u\cdot\text{\textgreek{D}}_{\mathbb{S}^{d-1}}u\Big)\, r^{d-1}drd\text{\textgreek{sv}},
\end{split}
\label{eq:CriticalEstimateFirstStep-3-1-1-1}
\end{equation}
 where 
\begin{equation}
A_{co}\doteq\big(1-\frac{a_{s}^{2}}{b_{s}}\big)\big(2-\frac{\text{\textgreek{e}}_{\text{\textgreek{b}}}}{2}\big)^{2}+4(1-a_{s})+O(\text{\textgreek{d}}_{H})-(1-a_{s})^{2}\frac{(3-\text{\textgreek{e}}_{\text{\textgreek{b}}})^{2}}{(1-b_{s}-\text{\textgreek{d}}_{H})}.\label{eq:FinalCoefficient}
\end{equation}

It thus remains to show that the parameters $a_{s},b_{s},\text{\textgreek{d}}_{H}$
can be suitably chosen in terms of $\text{\textgreek{e}}_{\text{\textgreek{b}}}$
(provided that $\text{\textgreek{e}}_{\text{\textgreek{b}}}\ll1$,
which we have assumed without loss of generality for the proof of
(\ref{eq:FinalEllipticEstimateAway-1})) so that $A_{coeff}>0$, and
then (\ref{eq:FinalEllipticEstimateAway-1}) will follow. Setting
$a_{s}=1-\text{\textgreek{d}}_{1}$ and $b_{s}=1-\frac{3}{2}\text{\textgreek{d}}_{1}$,
we can directly calculate from (\ref{eq:FinalCoefficient}): 
\begin{align}
A_{coeff} & =\big(\frac{1}{2}\text{\textgreek{d}}_{1}+O(\text{\textgreek{d}}_{1}^{2})\big)\big(4-2\text{\textgreek{e}}_{\text{\textgreek{b}}}+O(\text{\textgreek{e}}_{\text{\textgreek{b}}}^{2})\big)+4\text{\textgreek{d}}_{1}+O(\text{\textgreek{d}}_{H})-\frac{2}{3}\text{\textgreek{d}}_{1}\big(9-6\text{\textgreek{e}}_{\text{\textgreek{b}}}+O(\text{\textgreek{e}}_{\text{\textgreek{b}}}^{2})\big)\big(1+O(\text{\textgreek{d}}_{1}^{-1}\text{\textgreek{d}}_{H})\big)=\\
 & =+3\text{\textgreek{e}}_{\text{\textgreek{b}}}\text{\textgreek{d}}_{1}\Big(1+O(\text{\textgreek{e}}_{\text{\textgreek{b}}})+O(\text{\textgreek{e}}_{\text{\textgreek{b}}}^{-1}\text{\textgreek{d}}_{1})+O(\text{\textgreek{e}}_{\text{\textgreek{b}}}^{-1}\text{\textgreek{d}}_{1}^{-1}\text{\textgreek{d}}_{H})\Big),\nonumber 
\end{align}
 and thus it follows that $A_{coeff}>0$ provided that $\text{\textgreek{d}}_{H}\ll\text{\textgreek{d}}_{1}\ll\text{\textgreek{e}}_{\text{\textgreek{b}}}\ll1$
($\text{\textgreek{d}}_{H},\text{\textgreek{d}}_{1}$ can be fixed
in terms of $\text{\textgreek{e}}_{\text{\textgreek{b}}}$). Therefore,
from (\ref{eq:CriticalEstimateFirstStep-3-1-1-1}) (using also an
integration by parts for the last term of the right hand side of (\ref{eq:CriticalEstimateFirstStep-3-1-1-1}),
as well as the fact that $r^{d-1}drd\text{\textgreek{sv}}\sim dvol_{h}$
in that region) we finally obtain the desired estimate (\ref{eq:FinalEllipticEstimateAway-1})
on each connected component of the region $\{r\gg1\}$: 
\begin{align}
\int_{\mathcal{\mathcal{S}}}\text{\textgreek{q}}_{\ge R_{0}}r^{-\text{\textgreek{b}}}\Big(r^{-(d-1)}\partial_{r}\big(r^{d-1}\partial_{r}u\big)+r^{-2}\text{\textgreek{D}}_{\mathbb{S}^{d-1}}u\Big)^{2}\, r^{d-1}drd\text{\textgreek{sv}}\ge & c_{\text{\textgreek{b}}}\int_{\mathcal{S}}\text{\textgreek{q}}_{\ge R_{0}}r^{-\text{\textgreek{b}}}\Big(|\nabla^{2}u|_{h}^{2}+r^{-2}|\nabla u|_{h}^{2}\Big)\, dvol_{h}+\label{eq:FinalEllipticEstimateAway}\\
 & +\int_{\mathcal{S}}O(|\nabla\text{\textgreek{q}}_{\ge R_{0}}|+|\nabla^{2}\text{\textgreek{q}}_{\ge R_{0}}|)|\nabla u|_{h}^{2}\, dvol_{h}.\nonumber 
\end{align}

We will now proceed to establish estimates in the region $\{r\lesssim R_{0}\}$.
Proceeding through integrations by parts using the formula 
\begin{equation}
\nabla_{\text{\textgreek{m}}}\nabla_{\text{\textgreek{n}}}X_{\text{\textgreek{a}}_{1}\ldots\text{\textgreek{a}}_{k}}-\nabla_{\text{\textgreek{n}}}\nabla_{\text{\textgreek{m}}}X_{\text{\textgreek{a}}_{1}\ldots\text{\textgreek{a}}_{k}}=h^{\text{\textgreek{b}}_{1}\text{\textgreek{g}}_{1}}\mathcal{R}_{\text{\textgreek{m}\textgreek{n}}\text{\textgreek{b}}_{1}\text{\textgreek{a}}_{1}}X_{\text{\textgreek{g}}_{1}\text{\textgreek{a}}_{2}\ldots\text{\textgreek{a}}_{k}}+\ldots+h^{\text{\textgreek{b}}_{k}\text{\textgreek{g}}_{k}}\mathcal{R}_{\text{\textgreek{m}\textgreek{n}}\text{\textgreek{b}}_{k}\text{\textgreek{a}}_{k}}X_{\text{\textgreek{a}}_{1}\text{\textgreek{a}}_{2}\ldots\text{\textgreek{g}}_{k}},
\end{equation}
 we readily obtain that: 
\begin{align}
\int_{\mathcal{S}}\text{\textgreek{q}}_{\ge r_{0}}(1-\text{\textgreek{q}}_{\ge R_{0}})\cdot(\text{\textgreek{D}}_{h,\text{\textgreek{sv}}}u)^{2}\, dvol_{h} & =\int_{\mathcal{S}}\text{\textgreek{q}}_{\ge r_{0}}(1-\text{\textgreek{q}}_{\ge R_{0}})\cdot\Big(\nabla^{\text{\textgreek{m}}}\nabla_{\text{\textgreek{m}}}u+O_{\text{\textgreek{m}}}(r^{-2})\nabla^{\text{\textgreek{m}}}u\Big)\cdot\Big(\nabla^{\text{\textgreek{n}}}\nabla_{\text{\textgreek{n}}}u+O_{\text{\textgreek{n}}}(r^{-2})\nabla^{\text{\textgreek{n}}}u\Big)\, dvol_{h}\ge\label{eq:FirstIntegrationByPartsElliptic}\\
 & \ge\frac{1}{2}\int_{\mathcal{S}}\text{\textgreek{q}}_{\ge r_{0}}(1-\text{\textgreek{q}}_{\ge R_{0}})\cdot|\nabla^{2}u|_{h}^{2}\, dvol_{h}-C\int_{\{r_{0}\le r\le2R_{0}\}}|\nabla u|_{\tilde{h}}^{2}\, dvol_{h}.\nonumber 
\end{align}
 Notice that here we have used the volume form associated with $h$
(in place of $\tilde{h}$).

Therefore, from (\ref{eq:FirstIntegrationByPartsElliptic}) and (\ref{eq:FinalEllipticEstimateAway-1-1}),
and recalling that in the region $\{r\ge r_{0}\}$ we have $|\cdot|_{h}\sim_{r_{0}}|\cdot|_{\tilde{h}}$
and $dvol_{h}\sim_{r_{0}}dvol_{\tilde{h}}$, we deduce that there
exists some large $R_{1}>0$ depending only on on $\text{\textgreek{b}}$
and on the geometry of $(\mathcal{S},h)$ in the region $\{r\ge1\}$
such that: 
\begin{align}
\int_{\mathcal{S}}\text{\textgreek{q}}_{\ge r_{0}}r_{+}^{-\text{\textgreek{b}}}(\text{\textgreek{D}}_{h,\text{\textgreek{w}}}u)^{2}\, dvol_{h}\ge & c_{\text{\textgreek{b}}}\int_{\mathcal{S}}\text{\textgreek{q}}_{\ge r_{0}}r_{+}^{-\text{\textgreek{b}}}|\nabla^{2}u|_{\tilde{h}}^{2}\, dvol_{\tilde{h}}+c_{\text{\textgreek{b}}}\cdot\int_{\{r\ge R_{1}\}}r^{-2-\text{\textgreek{b}}}|\nabla u|_{h}^{2}\, dvol_{\tilde{h}}-\label{eq:EllipticEstimatesFarAway}\\
 & -C_{r_{0},\text{\textgreek{b}}}\int_{\{r_{0}\le r\le R_{1}\}}|\nabla u|_{h}^{2}\, dvol_{h}.\nonumber 
\end{align}

In the region $\{r\le2r_{0}\}$ we will perform the same integration
by parts procedure, but here we will use the explicit forms (\ref{eq:LaplacianSingularMetric})
and (\ref{eq:LaplacianSingularMetric-2}) for the Laplacian of $h$
near $\partial_{hor}\mathcal{S}$ and $\partial_{tim}\mathcal{S}$
respectively, as well as the form (\ref{eq:DistortedMetricNearHorizon})
for the non-singular metric $\tilde{h}$. 

On $[0,2r_{0}]\times\partial_{hor}\mathcal{S}$, the perturbed Laplacian
$\text{\textgreek{D}}_{h,\text{\textgreek{w}}}$ takes the form (\ref{eq:LaplacianSingularMetric}):
\begin{equation}
\text{\textgreek{D}}_{h,\text{\textgreek{w}}}=a^{-1}\partial_{r}\big(r\tilde{a}\cdot\partial_{r}\big)+\text{\textgreek{D}}_{h_{tan}}+O(r)\cdot X,\label{eq:LaplacianSingularMetric-1}
\end{equation}
 where $\tilde{a}=a\cdot(1+O(r))$. 

Let us define the weight function $w:[0,2r_{0}]\times\partial_{hor}\mathcal{S}\rightarrow(0,+\infty)$
by the relation: 
\begin{equation}
w(r,\text{\textgreek{j}})=a(r,\text{\textgreek{j}})\cdot\big(1+\int_{0}^{r}\tilde{a}^{-1}(\text{\textgreek{r}},\text{\textgreek{j}})\, d\text{\textgreek{r}}\big).\label{eq:WeightFunction}
\end{equation}
Since $\text{\textgreek{q}}_{\le r_{0}}^{hor}$ is supported in $[0,2r_{0}]\times\partial_{hor}\mathcal{S}$,
we calculate: 
\begin{align}
\int_{\mathcal{S}}\text{\textgreek{q}}_{\le r_{0}}^{hor}w\big(\text{\textgreek{D}}_{h,\text{\textgreek{w}}}u\big)^{2}\, dvol_{\tilde{h}} & =\int_{0}^{2r_{0}}\int_{\partial_{hor}\mathcal{S}}\text{\textgreek{q}}_{\le r_{0}}^{hor}w\Big(a^{-1}\partial_{r}\big(r\tilde{a}\partial_{r}u\big)+\text{\textgreek{D}}_{h_{tan}}u+O(r)\cdot Xu\Big)^{2}\, dh_{tan}dr=\label{eq:BeforeEllipticEstimatesSingular}\\
 & =\int_{0}^{2r_{0}}\int_{\partial_{hor}\mathcal{S}}\text{\textgreek{q}}_{\le r_{0}}^{hor}w\Big(a^{-2}\big(\partial_{r}\big(r\tilde{a}\partial_{r}u\big)\big)^{2}+2a^{-1}\big(\partial_{r}\big(r\tilde{a}\partial_{r}u\big)\big)\big(\text{\textgreek{D}}_{h_{tan}}u\big)+\big(\text{\textgreek{D}}_{h_{tan}}u\big)^{2}\Big)\, dh_{tan}dr+\nonumber \\
 & \hphantom{=}+\int_{0}^{2r_{0}}\int_{\partial_{hor}\mathcal{S}}\text{\textgreek{q}}_{\le r_{0}}^{hor}O(r)\Big|a^{-1}\partial_{r}\big(r\tilde{a}\partial_{r}u\big)+\text{\textgreek{D}}_{h_{tan}}u\Big|\cdot\big|Xu\big|+O(r^{2})\big|Xu\big|^{2}\Big)\, dh_{tan}dr.\nonumber 
\end{align}
Integrating by parts three times in the mixed $\partial_{r}(r\partial_{r}u)\cdot\text{\textgreek{D}}_{h_{tan}}u$
term, and using the fact that 
\[
\partial_{r}(wa^{-1})=\tilde{a}^{-1},
\]
 we estimate (notice that the resulting boundary terms at $r=0$ vanish
because $r\partial_{r}u$ and $r\nabla^{(h_{tan})}u$ vanish there):
\begin{align}
\int_{0}^{2r_{0}}\int_{\partial_{hor}\mathcal{S}}\text{\textgreek{q}}_{\le r_{0}}^{hor}wa^{-1}\big(\partial_{r}\big(r\tilde{a}\partial_{r}u\big)\big)\big(\text{\textgreek{D}}_{h_{tan}}u\big)\, dh_{tan}dr\ge & \int_{0}^{2r_{0}}\int_{\partial_{hor}\mathcal{S}}\text{\textgreek{q}}_{\le r_{0}}^{hor}\tilde{a}\cdot r(1+O(r))\cdot\big|\nabla^{(h_{tan})}\partial_{r}u\big|_{h_{s}}^{2}\, dh_{tan}dr-\label{eq:IntegrationByPartsMixedTerms}\\
 & -C\int_{0}^{r_{0}}\int_{\partial_{hor}\mathcal{S}}\big(r^{2}|\partial_{r}u|^{2}+|\nabla^{(h_{tan})}u|_{h_{tan}}^{2}\big)\, dh_{tan}dr-\nonumber \\
 & -C_{r_{0}}\int_{r_{0}}^{2r_{0}}\int_{\partial_{hor}\mathcal{S}}|\nabla u|_{\tilde{h}}^{2}\, dh_{tan}dr.\nonumber 
\end{align}
 Thus, from (\ref{eq:BeforeEllipticEstimatesSingular}) and (\ref{eq:IntegrationByPartsMixedTerms})
we obtain: 
\begin{align}
\int_{\mathcal{S}}\text{\textgreek{q}}_{\le r_{0}}^{hor}\big(\text{\textgreek{D}}_{h,\text{\textgreek{w}}}u\big)^{2}\, dvol_{\tilde{h}}\ge & c\int_{0}^{2r_{0}}\int_{\partial_{hor}\mathcal{S}}\text{\textgreek{q}}_{\le r_{0}}^{hor}\Big(\big(\partial_{r}\big(r\tilde{a}\partial_{r}u\big)\big)^{2}+r\big|\nabla^{(h_{tan})}\partial_{r}u\big|_{h_{tan}}^{2}+\big(\text{\textgreek{D}}_{h_{tan}}u\big)^{2}\Big)\, dh_{tan}dr-\label{eq:EllipticEstimatesInTheNearRegion}\\
 & -C\int_{0}^{r_{0}}\int_{\partial_{hor}\mathcal{S}}\big(r^{2}|\partial_{r}u|^{2}+|\nabla^{(h_{tan})}u|_{h_{tan}}^{2}\big)\, dh_{tan}dr-C_{r_{0}}\int_{r_{0}}^{2r_{0}}\int_{\partial_{hor}\mathcal{S}}|\nabla u|_{\tilde{h}}^{2}\, dh_{tan}dr+\nonumber \\
 & +\int_{0}^{2r_{0}}\int_{\partial_{hor}\mathcal{S}}\text{\textgreek{q}}_{\le r_{0}}^{hor}\Big(O(r)\Big|a^{-1}\partial_{r}\big(r\tilde{a}\partial_{r}u\big)+\text{\textgreek{D}}_{h_{tan}}u\Big|\cdot\big|Xu\big|+O(r^{2})\big|Xu\big|^{2}\Big)\, dh_{tan}dr.\nonumber 
\end{align}

By expanding 
\begin{equation}
\big(\partial_{r}\big(r\tilde{a}\partial_{r}u\big)\big)^{2}=r^{2}\tilde{a}^{2}(\partial_{r}^{2}u)^{2}+2r\tilde{a}^{2}(1+O(r))\partial_{r}^{2}u\cdot\partial_{r}u+\tilde{a}^{2}\big(1+O(r)\big)(\partial_{r}u)^{2}\label{eq:TrickBeforehardy}
\end{equation}
 and integrating by parts in the resulting $\partial_{r}^{2}u\partial_{r}u$
term, we obtain from (\ref{eq:EllipticEstimatesInTheNearRegion}):
\begin{align}
\int_{\mathcal{S}}\text{\textgreek{q}}_{\le r_{0}}^{hor}\big(\text{\textgreek{D}}_{h,\text{\textgreek{w}}}u\big)^{2}\, dvol_{\tilde{h}}\ge & c\int_{0}^{2r_{0}}\int_{\partial_{hor}\mathcal{S}}\text{\textgreek{q}}_{\le r_{0}}^{hor}\Big(r^{2}\big(\partial_{r}^{2}u\big)^{2}-Cr\cdot\big(\partial_{r}u\big)^{2}+r\big|\nabla^{(h_{tan})}\partial_{r}u\big|_{h_{tan}}^{2}+\big(\text{\textgreek{D}}_{h_{tan}}u\big)^{2}\Big)\, dh_{tan}dr-\label{eq:EllipticEstimatesInTheNearRegion-2}\\
 & -C\int_{0}^{r_{0}}\int_{\partial_{hor}\mathcal{S}}\big(r^{2}|\partial_{r}u|^{2}+|\nabla^{(h_{tan})}u|_{h_{tan}}^{2}\big)\, dh_{tan}dr-C_{r_{0}}\int_{r_{0}}^{2r_{0}}\int_{\partial_{hor}\mathcal{S}}|\nabla u|_{\tilde{h}}^{2}\, dh_{tan}dr+\nonumber \\
 & +\int_{0}^{2r_{0}}\int_{\partial_{hor}\mathcal{S}}\text{\textgreek{q}}_{\le r_{0}}^{hor}\Big(O(r)\Big|a^{-1}\partial_{r}\big(r\tilde{a}\partial_{r}u\big)+\text{\textgreek{D}}_{h_{tan}}u\Big|\cdot\big|Xu\big|+O(r^{2})\big|Xu\big|^{2}\Big)\, dh_{tan}dr.\nonumber 
\end{align}
Using, now a Hardy type inequality (of the form established in Lemma
\ref{lem:HardyGeneral}) for the$\partial_{r}^{2}u$ and $\nabla^{(h_{tan})}\partial_{r}u$
terms in the right hand side of (\ref{eq:EllipticEstimatesInTheNearRegion}),
as well as elliptic estimates for the $\text{\textgreek{D}}_{h_{tan}}u$
term (using here the compactness of the level sets of $r$) we obtain
from (\ref{eq:EllipticEstimatesInTheNearRegion-2}): 

\begin{align}
\int_{\mathcal{S}}\text{\textgreek{q}}_{\le r_{0}}^{hor}\big(\text{\textgreek{D}}_{h,\text{\textgreek{w}}}u\big)^{2}\, dvol_{\tilde{h}}\ge & c\int_{0}^{2r_{0}}\int_{\partial_{hor}\mathcal{S}}\text{\textgreek{q}}_{\le r_{0}}^{hor}\Big(r^{2}\big(\partial_{r}^{2}u\big)^{2}+r\big|\nabla^{(h_{tan})}\partial_{r}u\big|_{h_{tan}}^{2}+\big|\big(\nabla^{(h_{tan})}\big)^{2}u\big|_{h_{tan}}^{2}\Big)\, dh_{tan}dr+\label{eq:EllipticEstimatesInTheNearRegion-1}\\
 & +c\int_{0}^{2r_{0}}\int_{\partial_{hor}\mathcal{S}}\text{\textgreek{q}}_{\le r_{0}}^{hor}\Big(\big(\partial_{r}u\big){}^{2}+r^{-1}\log^{-2}(r)\cdot\big|\big(\nabla^{(h_{tan})}\big)u\big|_{h_{tan}}^{2}\Big)\, dh_{tan}dr+\nonumber \\
 & -C\int_{0}^{r_{0}}\int_{\partial_{hor}\mathcal{S}}\big(r^{2}|\partial_{r}u|^{2}+|\nabla^{(h_{tan})}u|_{h_{tan}}^{2}\big)\, dh_{tan}dr-C_{r_{0}}\int_{r_{0}}^{2r_{0}}\int_{\partial_{hor}\mathcal{S}}|\nabla u|_{\tilde{h}}^{2}\, dh_{tan}dr+\nonumber \\
 & +\int_{0}^{2r_{0}}\int_{\partial_{hor}\mathcal{S}}\text{\textgreek{q}}_{\le r_{0}}^{hor}\Big(O(r)\Big|a^{-1}\partial_{r}\big(r\tilde{a}\partial_{r}u\big)+\text{\textgreek{D}}_{h_{tan}}u\Big|\cdot\big|Xu\big|+O(r^{2})\big|Xu\big|^{2}\Big)\, dh_{tan}dr.\nonumber 
\end{align}
 Applying a Cauchy--Schwarz inequality on the last term of the rght
hand side of (\ref{eq:EllipticEstimatesInTheNearRegion-1}), and absorbing
all lower order (with respect to decaying powers in $r$) terms into
their top order counterparts, we obtain provided thar $r_{0}$ is
small enough:

\begin{align}
\int_{\mathcal{S}}\text{\textgreek{q}}_{\le r_{0}}^{hor}\big(\text{\textgreek{D}}_{h,\text{\textgreek{w}}}u\big)^{2}\, dvol_{\tilde{h}}\ge & c\int_{0}^{2r_{0}}\int_{\partial_{hor}\mathcal{S}}\text{\textgreek{q}}_{\le r_{0}}^{hor}\Big(r^{2}\big(\partial_{r}^{2}u\big)^{2}+r\big|\nabla^{(h_{tan})}\partial_{r}u\big|_{h_{tan}}^{2}+\big|\big(\nabla^{(h_{tan})}\big)^{2}u\big|_{h_{tan}}^{2}\Big)\, dh_{tan}dr+\label{eq:EllipticEstimatesInTheNearHorizonRegion}\\
 & +c\int_{0}^{2r_{0}}\int_{\partial_{hor}\mathcal{S}}\text{\textgreek{q}}_{\le r_{0}}^{hor}\Big(\big(\partial_{r}u\big){}^{2}+r^{-1}\log^{-2}(r)\cdot\big|\big(\nabla^{(h_{tan})}\big)u\big|_{h_{tan}}^{2}\Big)\, dh_{tan}dr-\nonumber \\
 & -C_{r_{0}}\int_{r_{0}}^{2r_{0}}\int_{\partial_{hor}\mathcal{S}}|\nabla u|_{\tilde{h}}^{2}\, dh_{tan}dr.\nonumber 
\end{align}

In the region $[0,2r_{0}]\times\partial_{tim}\mathcal{S}$, on the
other hand, we have 
\begin{equation}
\text{\textgreek{D}}_{h,\text{\textgreek{w}}}=\partial_{r}\big((1+O(r))\partial_{r}\big)+\text{\textgreek{D}}_{h_{tan}}+X.\label{eq:PurturbedLaplacianNearTimelikeBound}
\end{equation}
 Hence, we calculate after expanding the square (and applying the
product rule for derivatives): 
\begin{align}
\int_{\mathcal{S}}\text{\textgreek{q}}_{\le r_{0}}^{tim}\big(\text{\textgreek{D}}_{h,\text{\textgreek{w}}}u\big)^{2}\, dvol_{\tilde{h}} & =\int_{0}^{2r_{0}}\int_{\partial_{tim}\mathcal{S}}\text{\textgreek{q}}_{\le r_{0}}^{tim}\Big(\partial_{r}\big((1+O(r))\partial_{r}u\big)+\text{\textgreek{D}}_{h_{tan}}u+Xu\Big)^{2}\, dh_{tan}dr=\label{eq:BeforeEllipticEstimatesTimelike}\\
 & \ge\int_{0}^{2r_{0}}\int_{\partial_{tim}\mathcal{S}}\text{\textgreek{q}}_{\le r_{0}}^{tim}\Big((1+O(r))\big(\partial_{r}^{2}u\big)^{2}+2(1+O(r))\partial_{r}^{2}u\cdot\text{\textgreek{D}}_{h_{tan}}u+\big(\text{\textgreek{D}}_{h_{tan}}u\big)^{2}\Big)\, dh_{tan}dr-\nonumber \\
 & \hphantom{\ge}-C\int_{0}^{2r_{0}}\int_{\partial_{tim}\mathcal{S}}\text{\textgreek{q}}_{\le r_{0}}^{tim}\Big(\big(|\partial_{r}^{2}u|+|\text{\textgreek{D}}_{h_{tan}}u|\big)\cdot\big(|\partial_{r}u|+\big|\nabla^{(h_{tan})}u\big|_{h_{tan}}\big)+\big(|\partial_{r}u|^{2}+\big|\nabla^{(h_{tan})}u\big|_{h_{tan}}^{2}\big)\Big)\, dh_{tan}dr.\nonumber 
\end{align}
 After integrating by parts in the $\partial_{r}^{2}u\cdot\text{\textgreek{D}}_{h_{S}}u$
term and using elliptic estimates for $\text{\textgreek{D}}_{h_{S}}$
on the surfaces $\{r=const\}$, we obtain from (\ref{eq:BeforeEllipticEstimatesTimelike}):
\begin{align}
\int_{\mathcal{S}}\text{\textgreek{q}}_{\le r_{0}}^{tim}\big(\text{\textgreek{D}}_{h,\text{\textgreek{w}}}u\big)^{2}\, dvol_{\tilde{h}}\ge & c\int_{0}^{2r_{0}}\int_{\partial_{tim}\mathcal{S}}\text{\textgreek{q}}_{\le r_{0}}^{tim}\Big(\big(\partial_{r}^{2}u\big)^{2}+\big|\nabla^{(h_{tan})}\partial_{r}u\big|_{h_{tan}}^{2}+\big|\big(\nabla^{(h_{tan})}\big)^{2}u\big|_{h_{tan}}^{2}\Big)\, dh_{tan}dr+\label{eq:BeforeEllipticEstimatesTimelike-1}\\
 & +2\int_{\partial_{tim}\mathcal{S}}h_{tan}\big(\nabla^{(h_{tan})}\partial_{r}u,\nabla^{(h_{tan})}u\big)\, dh_{tan}-\nonumber \\
 & -C\int_{0}^{2r_{0}}\int_{\partial_{tim}\mathcal{S}}\text{\textgreek{q}}_{\le r_{0}}^{tim}\Big(\big(|\partial_{r}^{2}u|+|\text{\textgreek{D}}_{h_{tan}}u|\big)\cdot\big(|\partial_{r}u|+\big|\nabla^{(h_{tan})}u\big|_{h_{tan}}\big)+\big(|\partial_{r}u|^{2}+\big|\nabla^{(h_{tan})}u\big|_{h_{tan}}^{2}\big)\Big)\, dh_{tan}dr-\nonumber \\
 & -C\int_{0}^{2r_{0}}\int_{\partial_{tim}\mathcal{S}}\text{\textgreek{q}}_{\le r_{0}}^{tim}|\nabla^{2}u|_{\tilde{h}}\cdot|\nabla u|_{\tilde{h}}\, dh_{tan}dr-\nonumber \\
 & -C\int_{0}^{2r_{0}}\int_{\partial_{tim}\mathcal{S}}\big(|\nabla^{2}\text{\textgreek{q}}_{\le r_{0}}^{tim}|+|\nabla\text{\textgreek{q}}_{\le r_{0}}^{tim}|+1\big)|\nabla u|_{\tilde{h}}^{2}\, dh_{tan}dr.\nonumber 
\end{align}
 Applying a Cauchy--Schwarz inequality for the second and fourth terms
of the right hand side and using the Hardy type inequality
\begin{align}
\int_{0}^{2r_{0}}\int_{\partial_{tim}\mathcal{S}}\text{\textgreek{q}}_{\le r_{0}}^{tim}\Big(\big(\partial_{r}^{2}u\big)^{2}+\big|\nabla^{(h_{tan})}\partial_{r}u\big|_{h_{tan}}^{2}\Big)\, dh_{tan}dr\ge & c\int_{0}^{r_{0}}\int_{\partial_{tim}\mathcal{S}}r_{-}^{-1}\log^{-2}(r_{-})\Big(\big(\partial_{r}u\big)^{2}+\big|\nabla^{(h_{tan})}u\big|_{h_{tan}}^{2}\Big)\, dh_{tan}dr-\\
 & -C_{r_{0}}\int_{r_{0}}^{2r_{0}}\int_{\partial_{tim}\mathcal{S}}\Big(\big(\partial_{r}u\big)^{2}+\big|\nabla^{(h_{tan})}u\big|_{h_{tan}}^{2}\Big)\, dh_{tan}dr,\nonumber 
\end{align}
 we deduce from (\ref{eq:BeforeEllipticEstimatesTimelike-1}) provided
$r_{0}$ is small enough: 
\begin{align}
\int_{\mathcal{S}}\text{\textgreek{q}}_{\le r_{0}}^{tim}\big(\text{\textgreek{D}}_{h,\text{\textgreek{w}}}u\big)^{2}\, dvol_{\tilde{h}}\ge & c\int_{0}^{2r_{0}}\int_{\partial_{tim}\mathcal{S}}\text{\textgreek{q}}_{\le r_{0}}^{tim}\Big(\big(\partial_{r}^{2}u\big)^{2}+\big|\nabla^{(h_{tan})}\partial_{r}u\big|_{h_{tan}}^{2}+\big|\big(\nabla^{(h_{tan})}\big)^{2}u\big|_{h_{tan}}^{2}\Big)\, dh_{tan}dr+\label{eq:EllipticEstimatesTimelike}\\
 & +c\int_{0}^{r_{0}}\int_{\partial_{tim}\mathcal{S}}r_{-}^{-1}\log^{-2}(r_{-})\Big(\big(\partial_{r}u\big)^{2}+\big|\nabla^{(h_{tan})}u\big|_{h_{tan}}^{2}\Big)\, dh_{tan}dr-\nonumber \\
 & -C_{r_{0}}\int_{r_{0}}^{2r_{0}}\int_{\partial_{tim}\mathcal{S}}\Big(\big(\partial_{r}u\big)^{2}+\big|\nabla^{(h_{tan})}u\big|_{h_{tan}}^{2}\Big)\, dh_{tan}dr-\nonumber \\
 & +2\int_{\partial_{tim}\mathcal{S}}h_{tan}\big(\nabla^{(h_{tan})}\partial_{r}u,\nabla^{(h_{tan})}u\big)\, dh_{tan}.\nonumber 
\end{align}

By adding (\ref{eq:EllipticEstimatesFarAway}), (\ref{eq:EllipticEstimatesInTheNearHorizonRegion})
and (\ref{eq:EllipticEstimatesTimelike}) we infer that: 
\begin{align}
\int_{\mathcal{S}}r_{+}^{-\text{\textgreek{b}}}\big(\text{\textgreek{D}}_{h,\text{\textgreek{w}}}u\big)^{2}\, dvol_{\tilde{h}}\ge & c\int_{\mathcal{S}}r_{+}^{-\text{\textgreek{b}}}\big|\big(\nabla^{(\tilde{h})}\big)^{2}u\big|_{h}^{2}\, dvol_{\tilde{h}}+c\int_{\mathcal{S}}r_{+}^{-2-\text{\textgreek{b}}}\Big(|(\text{\textgreek{q}}_{\le r_{0}}^{hor}+\text{\textgreek{q}}_{\le r_{0}}^{tim})Yu|^{2}+r_{-}^{-1}\log^{-2}(\frac{r_{-}}{2})\big|\nabla^{(\tilde{h})}u\big|_{h}^{2}\Big)\, dvol_{\tilde{h}}-\label{eq:FullEllipticEstimatesBeforeCompactness}\\
 & -C_{r_{0}}\int_{\{r_{0}\le r\le R_{1}\}}|\nabla u|_{\tilde{h}}^{2}\, dvol_{\tilde{h}}+2\int_{\partial_{tim}\mathcal{S}}h_{tan}\big(\nabla^{(h_{tan})}Yu,\nabla^{(h_{tan})}u\big)\, dh_{tan}.\nonumber 
\end{align}
From now on, we will assume that $r_{0}$ has been fixed, and we will
drop the dependence of constants on it.

Let us denote with $\mathcal{H}^{2}(\mathcal{S},h,\tilde{h})$ the
semi-norm space consisting of the functions $\text{\textgreek{y}}$
of $H_{loc}^{2}(\mathcal{S})$ with $||\text{\textgreek{y}}||_{\mathcal{H}^{2}(\mathcal{S},h,\tilde{h})}<+\infty$,
where 
\begin{equation}
||\text{\textgreek{y}}||_{\mathcal{H}^{2}(\mathcal{S},h,\tilde{h})}^{2}\doteq\int_{\mathcal{S}}r_{+}^{-\text{\textgreek{b}}}\big|\big(\nabla^{(\tilde{h})}\big)^{2}\text{\textgreek{y}}\big|_{h}^{2}\, dvol_{\tilde{h}}+\int_{\mathcal{S}}r_{+}^{-2-\text{\textgreek{b}}}\Big(|(\text{\textgreek{q}}_{\le r_{0}}^{hor}+\text{\textgreek{q}}_{\le r_{0}}^{tim})Y\text{\textgreek{y}}|^{2}+r_{-}^{-1}\log^{-2}(\frac{r_{-}}{2})\big|\nabla\text{\textgreek{y}}\big|_{h}^{2}\Big)\, dvol_{\tilde{h}}
\end{equation}
(this is the semi-norm appearing in the right hand side of (\ref{eq:FullEllipticEstimatesBeforeCompactness})).
Notice that $||\cdot||_{\mathcal{H}^{2}(\mathcal{S},h,\tilde{h})}$
becomes an actual norm if we mod out the constant functions (the resulting
normed space being a Hilbert space). It will be also convenient to
introduce the semi-norm space $\mathcal{H}^{1}(\mathcal{S},h,\tilde{h})$
defined by the norm: 
\begin{equation}
||\text{\textgreek{y}}||_{\mathcal{H}^{1}(\mathcal{S},h,\tilde{h})}^{2}\doteq\int_{\mathcal{S}}r_{+}^{-\text{\textgreek{b}}}\big|\nabla\text{\textgreek{y}}\big|_{h}^{2}\, dvol_{\tilde{h}}.
\end{equation}

Using the Rellich--Kondrachov theorem for smooth manifolds with boundary
(see e.\,g.~\cite{Hebey1999}), we infer that for any $\text{\textgreek{e}}_{0}>0$
(which will be fixed small with respect to all the constants, and
their inverses, appearing in (\ref{eq:FullEllipticEstimatesBeforeCompactness}),
as well as the restriction of the weights in the integrals of (\ref{eq:FullEllipticEstimatesBeforeCompactness})
over $\{r_{0}\le r\le R_{1}\}$), the set $\mathcal{D}_{\text{\textgreek{e}}_{0}}$
of functions $\text{\textgreek{y}}\in\mathcal{H}^{2}(\mathcal{S},h,\tilde{h})$
satisfying 
\begin{equation}
||\text{\textgreek{y}}||_{\mathcal{H}^{1}(\mathcal{S},h,\tilde{h})}^{2}=1\label{eq:UnitBallLowSpace}
\end{equation}
and 
\begin{equation}
\int_{\{r_{0}\le r\le R_{1}\}}|\nabla\text{\textgreek{y}}|_{\tilde{h}}^{2}\, dvol_{\tilde{h}}\ge\text{\textgreek{e}}_{0}||\text{\textgreek{y}}||_{\mathcal{H}^{2}(\mathcal{S},h,\tilde{h})}^{2}\label{eq:RellichKondrachov}
\end{equation}
 is a precompact subset of the semi-norm space $\mathcal{H}^{1}(\mathcal{S},h,\tilde{h})$. 

From Lemma \ref{lem:AbscenceOfHarmonicFunctions} we deduce that any
non-constant function $\text{\textgreek{y}}\in\mathcal{H}^{2}(\mathcal{S},h,\tilde{h})$
satisfies for any $\text{\textgreek{x}}\ge\text{\textgreek{b}}$:
\begin{equation}
||\text{\textgreek{D}}_{h,\text{\textgreek{w}}}\text{\textgreek{y}}||_{L_{\text{\textgreek{x}}}^{2}(\mathcal{S},\tilde{h})}+\max\big\{-\int_{\partial_{tim}\mathcal{S}}Y\text{\textgreek{y}}\cdot\text{\textgreek{y}}\, dh_{tan},0\}>0,\label{eq:PositivityLaplacian}
\end{equation}
 where 
\begin{equation}
||f||_{L_{\text{\textgreek{x}}}^{2}(\mathcal{S},\tilde{h})}^{2}\doteq\int_{\mathcal{S}}r_{+}^{-\text{\textgreek{x}}}f^{2}\, dvol_{\tilde{h}}.
\end{equation}
 Therefore, since $\mathcal{D}_{\text{\textgreek{e}}_{0}}$ is precompact
and no constant function lies in its closure in the seminorm space
$\mathcal{H}^{1}(\mathcal{S},h,\tilde{h})$ (due to (\ref{eq:UnitBallLowSpace})),
we infer that we can bound for any $\text{\textgreek{y}}\in\mathcal{D}$:
\begin{equation}
\int_{\{r_{0}\le r\le R_{1}\}}|\nabla\text{\textgreek{y}}|_{\tilde{h}}^{2}\, dvol_{\tilde{h}}\le C_{\text{\textgreek{e}}_{0}}\Big(||\text{\textgreek{D}}_{h,\text{\textgreek{w}}}\text{\textgreek{y}}||_{L_{\text{\textgreek{b}}}^{2}(\mathcal{S},\tilde{h})}^{2}+\max\big\{-\int_{\partial_{tim}\mathcal{S}}Y\text{\textgreek{y}}\cdot\text{\textgreek{y}}\, dh_{tan},0\}\Big).\label{eq:BoundLowFrequency}
\end{equation}

Thus, fixing $\text{\textgreek{e}}_{0}$ small enough in terms of
$r_{0}$, $\text{\textgreek{b}}$ and the geometry of $(\mathcal{S},h)$,
returning to our original function $u$ we distinguish between two
cases:

\begin{enumerate}

\item In case $u$ is constant or $||u||_{\mathcal{H}^{1}(\mathcal{S},h,\tilde{h})}^{-1}\cdot u\in\mathcal{D}_{\text{\textgreek{e}}_{0}}$,
from (\ref{eq:FullEllipticEstimatesBeforeCompactness}) and (\ref{eq:BoundLowFrequency})
we can bound: 
\begin{equation}
\begin{split}\int_{\mathcal{S}}r_{+}^{-\text{\textgreek{b}}}\big|\big(\nabla^{(\tilde{h})}\big)^{2}u\big|_{h}^{2}\, dvol_{\tilde{h}}+ & \int_{\mathcal{S}}r_{+}^{-2-\text{\textgreek{b}}}\Big(|(\text{\textgreek{q}}_{\le r_{0}}^{hor}+\text{\textgreek{q}}_{\le r_{0}}^{tim})Yu|^{2}+r_{-}^{-1}\log^{-2}(\frac{r_{-}}{2})\big|\nabla u\big|_{h}^{2}\Big)\, dvol_{\tilde{h}}\le\\
\le & C_{\text{\textgreek{b}}}\int_{\mathcal{S}}r_{+}^{-\text{\textgreek{b}}}\big(\text{\textgreek{D}}_{h,\text{\textgreek{w}}}u\big)^{2}\, dvol_{\tilde{h}}+C_{\text{\textgreek{b}}}\max\big\{-\int_{\partial_{tim}\mathcal{S}}h_{tan}\big(\nabla^{(h_{tan})}Yu,\nabla^{(h_{tan})}u\big)\, dh_{tan},0\big\}+\\
 & +C_{\text{\textgreek{b}}}\max\big\{-\int_{\partial_{tim}\mathcal{S}}Y\text{\textgreek{y}}\cdot\text{\textgreek{y}}\, dh_{tan},0\}.
\end{split}
\end{equation}

\item In case $u$ is not constant and $||u||_{\mathcal{H}^{1}(\mathcal{S},h,\tilde{h})}^{-1}\cdot u\notin\mathcal{D}_{\text{\textgreek{e}}_{0}}$,
from the definition of $\mathcal{D}_{\text{\textgreek{e}}_{0}}$ (i.\,e.~(\ref{eq:RellichKondrachov}))
we can bound: 
\begin{equation}
\int_{\{r_{0}\le r\le R_{1}\}}|\nabla u|_{\tilde{h}}^{2}\, dvol_{\tilde{h}}\le\text{\textgreek{e}}_{0}||u||_{\mathcal{H}^{2}(\mathcal{S},h,\tilde{h})}^{2}.\label{eq:RellichKondrachov-1}
\end{equation}

\end{enumerate}

Thus, if $\text{\textgreek{e}}_{0}$ has been fixed small enough in
terms of $r_{0}$, $\text{\textgreek{b}}$ and the geometry of $(\mathcal{S},h)$,
from (\ref{eq:FullEllipticEstimatesBeforeCompactness}) and (\ref{eq:RellichKondrachov-1})
we deduce that: 
\begin{equation}
\begin{split}\int_{\mathcal{S}}r_{+}^{-\text{\textgreek{b}}}\big|\big(\nabla^{(\tilde{h})}\big)^{2}u\big|_{h}^{2}\, dvol_{\tilde{h}}+ & \int_{\mathcal{S}}r_{+}^{-2-\text{\textgreek{b}}}\Big(|\text{\textgreek{q}}_{\le r_{0}}Yu|^{2}+r_{-}^{-1}\log^{-2}(\frac{r_{-}}{2})\big|\nabla u\big|_{h}^{2}\Big)\, dvol_{\tilde{h}}\le\\
\le\, & C_{\text{\textgreek{b}}}\int_{\mathcal{S}}r_{+}^{-\text{\textgreek{b}}}\big(\text{\textgreek{D}}_{h,\text{\textgreek{w}}}u\big)^{2}\, dvol_{\tilde{h}}+C_{\text{\textgreek{b}}}\max\big\{-\int_{\partial_{tim}\mathcal{S}}h_{tan}\big(\nabla^{(h_{tan})}Yu,\nabla^{(h_{tan})}u\big)\, dh_{tan},0\big\}+\\
 & +C_{\text{\textgreek{b}}}\max\big\{-\int_{\partial_{tim}\mathcal{S}}Y\text{\textgreek{y}}\cdot\text{\textgreek{y}}\, dh_{tan},0\}.
\end{split}
\label{eq:FinalStepEllipticEstimates}
\end{equation}
 Therefore, the elliptic estimates (\ref{eq:EllipticEstimatesDegenerateAtTimelike})
and (\ref{eq:EllipticEstimatesDegenerateRiemannian}) in the case
$l=2$ (and $k_{0}=0$) have been established.

The case when $2<l\le\lfloor\frac{d+1}{2}\rfloor$ follows in an analogous
way: In order to derive the analogue of (\ref{eq:FullEllipticEstimatesBeforeCompactness}),
one needs to commute $l-2$ times with $\nabla^{(\tilde{h}_{tim})}$
(the curvature terms appearing in this way are treated exactly as
we did for the simple curvature terms in the $l=2$ case using the
flat asymptotics of $(\mathcal{S},h)$). By applying a Hardy type
inequality near $\partial\mathcal{S}$ (using the form of the metric
$\tilde{h}_{tim}$ there) in order to obtain an estimate of the form
\begin{equation}
\sum_{j=0}^{l-3}\int_{\{r\le r_{0}\}}\big|\big(\nabla^{(\tilde{h})}\big)^{j}(\text{\textgreek{D}}_{h,\text{\textgreek{w}}}u)\big|_{\tilde{h}_{tim}}^{2}\, dvol_{\tilde{h}}\lesssim\int_{\{r\le2r_{0}\}}\big|\big(\nabla^{(\tilde{h})}\big)^{l-2}(\text{\textgreek{D}}_{h,\text{\textgreek{w}}}u)\big|_{\tilde{h}_{tim}}^{2}\, dvol_{\tilde{h}}+\sum_{j=0}^{l-3}\int_{\{r_{0}\le r\le2r_{0}\}}\big|\big(\nabla^{(\tilde{h})}\big)^{j}(\text{\textgreek{D}}_{h,\text{\textgreek{w}}}u)\big|_{\tilde{h}_{tim}}^{2}\, dvol_{\tilde{h}},
\end{equation}
and then integrating by parts as before, one readily obtains the following
estimate: 
\begin{equation}
\begin{split}\int_{\mathcal{S}}r_{+}^{-\text{\textgreek{b}}}\big|\big(\nabla^{(\tilde{h})}\big)^{l-2} & (\text{\textgreek{D}}_{h,\text{\textgreek{w}}}u)\big|_{\tilde{h}_{tim}}^{2}\, dvol_{\tilde{h}}\ge c\int_{\mathcal{S}}r_{+}^{-\text{\textgreek{b}}}|u|_{\tilde{h},\tilde{h}_{tim};l}^{2}\, dvol_{\tilde{h}}+c\sum_{j=1}^{l-1}\int_{\{r\ge r_{0}\}}r_{+}^{-2j-\text{\textgreek{b}}}\cdot|\big(\nabla^{(\tilde{h})}\big)^{l-j}u|_{\tilde{h}}^{2}\, dvol_{\tilde{h}}+\\
 & +c\sum_{j=1}^{l-1}\int_{\{r\le r_{0}\}}\Big(|\big(\nabla^{(\tilde{h})}\big)^{l-j-1}(Yu)|_{\tilde{h}_{tim}}^{2}+r_{-}^{-1}\log^{-2}(\frac{r_{-}}{2})\cdot|\big(\nabla^{(\tilde{h})}\big)^{l-j-1}\big(i_{*}(\nabla^{(h_{tan})}u)\big)|_{\tilde{h}_{tim}}^{2}\Big)\, dvol_{\tilde{h}}\Big\}-\\
 & -C_{r_{0}}\sum_{j=1}^{l-1}\int_{\{r_{0}\le r\le R_{1}\}}|\nabla^{l-j-1}u|_{\tilde{h}}^{2}\, dvol_{\tilde{h}}+c\cdot\min\Big\{\int_{\partial_{tim}\mathcal{S}}h_{tan}\Big(\big(\nabla^{(h_{tan})}\big)(Yu),\big(\nabla^{(h_{tan})}\big)u\Big)\, dh_{tan}\big\},0\Big\}.
\end{split}
\label{eq:FullEllipticEstimatesBeforeCompactness-1}
\end{equation}
 Using the relation 
\begin{equation}
\text{\textgreek{D}}_{h,\text{\textgreek{w}}}u=Y\big((1+O(r))\cdot Y\big)u+\text{\textgreek{D}}_{h_{tan}}u+Xu
\end{equation}
near $\partial_{tim}\mathcal{S}$, together with a Hardy-type inequality,
we can immediately estimate: 
\begin{align}
\int_{\{r\le r_{0}\}}|u|_{\tilde{h},\big(1-\log(r_{tim})\big)\tilde{h};l}^{2}\, dvol_{\tilde{h}}\lesssim & \int_{\{r\le2r_{0}\}}|u|_{\tilde{h},\tilde{h}_{tim};l}^{2}\, dvol_{\tilde{h}}+\int_{\{r\le2r_{0}\}}\big|\big(\nabla^{(\tilde{h})}\big)^{l-2}(\text{\textgreek{D}}_{h,\text{\textgreek{w}}}u)\big|_{\big(1-\log(r_{tim})\big)\tilde{h}}^{2}\, dvol_{\tilde{h}}+\\
 & +\sum_{j=1}^{l-1}\int_{\{r_{0}\le r\le2r_{0}\}}|\nabla^{l-j-1}u|_{\tilde{h}}^{2}\, dvol_{\tilde{h}}\nonumber 
\end{align}
and thus from (\ref{eq:FullEllipticEstimatesBeforeCompactness-1})
we also obtain: 
\begin{equation}
\begin{split}\int_{\mathcal{S}} & r_{+}^{-\text{\textgreek{b}}}\big|\big(\nabla^{(\tilde{h})}\big)^{l-2}(\text{\textgreek{D}}_{h,\text{\textgreek{w}}}u)\big|_{\big(1-\log(r_{tim})\big)\tilde{h}}^{2}\, dvol_{\tilde{h}}\ge c\int_{\mathcal{S}}r_{+}^{-\text{\textgreek{b}}}|u|_{\tilde{h},\big(1-\log(r_{tim})\big)\tilde{h};l}^{2}\, dvol_{\tilde{h}}+c\sum_{j=1}^{l-1}\int_{\{r\ge r_{0}\}}r_{+}^{-2j-\text{\textgreek{b}}}\cdot|\big(\nabla^{(\tilde{h})}\big)^{l-j}u|_{\tilde{h}}^{2}\, dvol_{\tilde{h}}+\\
 & +c\sum_{j=1}^{l-1}\int_{\{r\le r_{0}\}}\Big(|\big(\nabla^{(\tilde{h})}\big)^{l-j-1}(Yu)|_{\big(1-\log(r_{tim})\big)\tilde{h}}^{2}+r_{-}^{-1}\log^{-2}(\frac{r_{-}}{2})\cdot|\big(\nabla^{(\tilde{h})}\big)^{l-j-1}\big(i_{*}(\nabla^{(h_{tan})}u)\big)|_{\big(1-\log(r_{tim})\big)\tilde{h}}^{2}\Big)\, dvol_{\tilde{h}}\Big\}-\\
 & -C_{r_{0}}\sum_{j=1}^{l-1}\int_{\{r_{0}\le r\le R_{1}\}}|\nabla^{l-j-1}u|_{\tilde{h}}^{2}\, dvol_{\tilde{h}}+c\cdot\min\Big\{\int_{\partial_{tim}\mathcal{S}}h_{tan}\Big(\big(\nabla^{(h_{tan})}\big)(Yu),\big(\nabla^{(h_{tan})}\big)u\Big)\, dh_{tan}\big\},0\Big\}.
\end{split}
\label{eq:FullEllipticEstimatesBeforeCompactness-1-1}
\end{equation}
 Using the same Fredholm-type technique as before, we can absorb the
$\int_{\{r_{0}\le r\le R_{1}\}}|\nabla^{l-j-1}u|_{\tilde{h}}^{2}$
terms in the right hand side of (\ref{eq:FullEllipticEstimatesBeforeCompactness-1})
and (\ref{eq:FullEllipticEstimatesBeforeCompactness-1-1}) after adding
to the left hand side of each of these estimates a large multiple
of 
\begin{equation}
\int_{\mathcal{S}}r_{+}^{-\text{\textgreek{b}}}\big|\big(\nabla^{(\tilde{h})}\big)^{l-2}(\text{\textgreek{D}}_{h,\text{\textgreek{w}}}u)\big|_{\tilde{h}_{tim}}^{2}\, dvol_{\tilde{h}}+\max\Big\{-\int_{\partial_{tim}\mathcal{S}}Yu\cdot u\, dh_{tan}\big\},0\Big\}
\end{equation}
and 
\begin{equation}
\int_{\mathcal{S}}r_{+}^{-\text{\textgreek{b}}}\big|\big(\nabla^{(\tilde{h})}\big)^{l-2}(\text{\textgreek{D}}_{h,\text{\textgreek{w}}}u)\big|_{\big(1-\log(r_{tim})\big)\tilde{h}}^{2}\, dvol_{\tilde{h}}+\max\Big\{-\int_{\partial_{tim}\mathcal{S}}Yu\cdot u\, dh_{tan}\big\},0\Big\}
\end{equation}
respectively, thus obtaining inequalities (\ref{eq:EllipticEstimatesDegenerateAtTimelike})
and (\ref{eq:EllipticEstimatesDegenerateRiemannian}). We will omit
the details. 

The case $k_{0}\ge1$ follows in exactly the same way. 
\end{proof}
Let us assume that we are given a smooth function $\text{\textgreek{w}}_{nd}:\mathcal{S}\rightarrow(0,+\infty)$
with $\text{\textgreek{w}}_{nd}=1+O(r^{-1})$ in $\{r\gg1\}$ (notice
that we necessarily have $\text{\textgreek{w}}_{nd}\neq0$ on $\partial\mathcal{S}$),
and let us define the non-degenerate elliptic operator 
\begin{equation}
\text{\textgreek{D}}_{\tilde{h},\text{\textgreek{w}}_{nd}}\doteq\text{\textgreek{w}}_{nd}^{-1}div_{\tilde{h}}\big(\text{\textgreek{w}}_{nd}\cdot d\big).\label{eq:NonDegeneratePerturbedLaplacian}
\end{equation}
This operator will model the operator (\ref{eq:PerturbedLaplacian})
associated to the metric $h_{\text{\textgreek{t}},N}$ on the hypersurfaces
$\{\bar{t}=\text{\textgreek{t}}\}$ of the spacetimes $(\mathcal{M},g)$
of Section \ref{sec:Firstdecay}. The following non-degenerate variant
of Proposition \ref{prop:DegenerateEllipticEstimates} holds for (\ref{eq:NonDegeneratePerturbedLaplacian}):
\begin{prop}
\label{Prop:NonDegenerateEllipticEstimates}\textbf{(Non-degenerate elliptic estimates).}
For any $l\in\mathbb{N}$ with $2\le l\le\lfloor\frac{d+1}{2}\rfloor$,
any $k_{0}\in\mathbb{N}$ and any $\text{\textgreek{b}}\in(-\bar{\text{\textgreek{d}}}_{k_{0}},1)$
for some $\bar{\text{\textgreek{d}}}_{k_{0}}>0$ depending on $k_{0}$,
we can bound for any $u\in C^{\infty}(\mathcal{S})$ satisfying $\limsup_{r\rightarrow+\infty}\big|r^{\frac{d-1}{2}+j}\nabla^{j}u\big|_{h}<+\infty$
for any $j\le l+k_{0}$: 
\begin{equation}
\begin{split}\sum_{k=0}^{k_{0}}\Big\{\int_{\mathcal{S}}r_{+}^{-\text{\textgreek{b}}}|\big(\nabla^{(\tilde{h})}\big)^{k+l}u|_{\tilde{h}}^{2}\, dvol_{\tilde{h}}+ & \sum_{j=1}^{l-1}\int_{\mathcal{S}}r_{+}^{-2j-\text{\textgreek{b}}}r_{-}^{-1}\log^{-2}(\frac{r_{-}}{2})\cdot|\big(\nabla^{(\tilde{h})}\big)^{k+l-j}u|_{\tilde{h}}^{2}\, dvol_{\tilde{h}}\Big\}\le\\
\le\, & C_{\text{\textgreek{b}},k_{0}}\sum_{k=0}^{k_{0}}\int_{\mathcal{S}}r_{+}^{-\text{\textgreek{b}}}\big|\big(\nabla^{(\tilde{h})}\big)^{k+l-2}(\text{\textgreek{D}}_{\tilde{h},\text{\textgreek{w}}_{nd}}u)\big|_{\tilde{h}}^{2}\, dvol_{\tilde{h}}+\\
 & +C_{\text{\textgreek{b}},k_{0}}\sum_{j=0}^{k_{0}+l-1}\max\Big\{-Re\big\{\int_{\partial\mathcal{S}}h_{tan}\Big(\big(\nabla^{(h_{tan})}\big)^{j}(Yu),\big(\nabla^{(h_{tan})}\big)^{j}\bar{u}\Big)\, dh_{tan}\big\},0\Big\},
\end{split}
\label{eq:EllipticEstimatesNonDegenerate}
\end{equation}
 where $\nabla^{(\tilde{h})}$ denotes the covariant derivative with
respect to the metric $\tilde{h}$ and the constant $C_{\text{\textgreek{b}}}$
of the right hand side depends only on $\text{\textgreek{b}}$ and
on the geometry of $(\mathcal{S},h)$, $\tilde{h}$ and $\text{\textgreek{w}}_{nd}$.
 \end{prop}
\begin{rem*}
Notice that in contrast to (\ref{eq:EllipticEstimatesDegenerateRiemannian}),
the right hand side of (\ref{eq:EllipticEstimatesNonDegenerate})
contains terms on the whole of $\partial\mathcal{S}$.
\end{rem*}
The proof of Proposition \ref{Prop:NonDegenerateEllipticEstimates}
follows in the same lines as that for Lemma \ref{prop:DegenerateEllipticEstimates}
(using everywhere the metric $\tilde{h}$ in place of $h$ and $\tilde{h}_{tim}$),
and hence it will be omitted. 

It will also be useful to establish the following estimate in the
region $\{r\gg1\}$ (this estimate will be used in Section \ref{sec:Improved-polynomial-decay}
to control the error terms arising when commuting $\square_{g}$ with
the vector field $K_{R_{c}}$ which is not settling to a Killing field
in the region $\{r\sim R_{c}\}$):
\begin{prop}
\label{Prop:ControlOfTheAngularTermsAway}\textbf{(Improved control of derivatives in the far away region).}There
exists some $R_{0}\ge1$ large in terms of the geometry of $(\mathcal{S},h)$,
such that for any $R>R_{0}$ the following bound holds for any $l\in\mathbb{N}$
with $l\le\lfloor\frac{d+1}{2}\rfloor$, any $k_{0}\in\mathbb{N}$,
any $\text{\textgreek{b}}\in(-\bar{\text{\textgreek{d}}}_{k_{0}},1)$
for some $\bar{\text{\textgreek{d}}}_{k_{0}}>0$ depending on $k_{0}$,
any $0<\text{\textgreek{e}}<1$ and any $u\in C^{\infty}(\mathcal{N})$
with $\sum_{k=0}^{k_{0}}\sum_{j=1}^{l}\int_{\{r\ge R_{0}\}}r^{-2(l-j)}|\nabla^{j+k}u|_{h}^{2}\, dvol_{h}<\infty$:
\begin{multline}
\sum_{k=0}^{k_{0}}\Big\{\sum_{j=1}^{l}\int_{\{R\le r\le2R\}}r^{2(j-1)-\text{\textgreek{b}}+\text{\textgreek{e}}}|\nabla^{j+k}u|_{h}^{2}\, dvol_{h}+\sum_{j=1}^{l}\int_{\{r\ge R_{0}\}}r^{-2(l-j)-\text{\textgreek{b}}}|\nabla^{j+k}u|_{h}^{2}\, dvol_{h}\Big\}\le\\
\le C_{R,\text{\textgreek{b}},\text{\textgreek{e}},k_{0}}\sum_{k=0}^{k_{0}}\int_{\{r\ge R_{0}\}}r^{-\text{\textgreek{b}}}\big|\nabla^{k+l-2}(\text{\textgreek{D}}_{h,\text{\textgreek{w}}}u)\big|_{h}^{2}\, dvol_{h}+C_{\text{\textgreek{b}},\text{\textgreek{e}},k_{0}}\sum_{j=1}^{l+k_{0}}\int_{\{R_{0}\le r\le2R_{0}\}}|\nabla^{j}u|_{h}^{2}\, dvol_{h}.\label{eq:EstimateFarAwayForAngularBulk}
\end{multline}
\end{prop}
\begin{rem*}
Notice that the constant in front of the last term of the right hand
side of (\ref{eq:EstimateFarAwayForAngularBulk}) does not depend
on $R$.\end{rem*}
\begin{proof}
Without loss of generality, we can assume that $k_{0}=0$, since the
proof in the case $k_{0}\ge1$ follows in exactly the same way.

It suffices to establish the following estimate on $\mathbb{R}^{d}$
for $R_{0}$, $R$, $l$, $\text{\textgreek{b}}$, $\text{\textgreek{e}}$
as above and $u\in C^{\infty}(\mathbb{R}^{d})$ with $\limsup_{r\rightarrow+\infty}\big|r^{\frac{d-1}{2}+j}\nabla^{j}u\big|_{h}<+\infty$
for $j\le l$:

\begin{multline}
\sum_{j=1}^{l}\int_{\{R\le r\le2R\}}r^{2(j-1)-\text{\textgreek{b}}+\text{\textgreek{e}}}\big|\nabla^{j}u\big|_{e}^{2}\, dvol_{e}+\sum_{j=1}^{l}\int_{\{r\ge R_{0}\}}r^{-2(l-j)-\text{\textgreek{b}}}|\nabla^{j}u|_{e}^{2}\, dvol_{e}\le\\
\le C_{R,R_{0},\text{\textgreek{b}},\text{\textgreek{e}}}\int_{\{r\ge R_{0}\}}r^{-\text{\textgreek{b}}}\big|\nabla^{l-2}(\text{\textgreek{D}}_{\mathbb{R}^{d}}u)\big|_{e}^{2}\, dvol_{e}+C_{R_{0},\text{\textgreek{b}},\text{\textgreek{e}}}\sum_{j=1}^{l}\int_{\{R_{0}\le r\le2R_{0}\}}|\nabla^{j}u|_{e}^{2}\, dvol_{e}.\label{eq:EstimateFarAwayFlat}
\end{multline}
 Assuming that (\ref{eq:EstimateFarAwayFlat}) holds, by substituting
\begin{equation}
\text{\textgreek{D}}_{h,\text{\textgreek{w}}}u=\text{\textgreek{D}}_{\mathbb{R}^{d}}u+O(r^{-1})\nabla^{2}u+O(r^{-2})\nabla u
\end{equation}
 one obtains (\ref{eq:EstimateFarAwayForAngularBulk}) by absorbing
the resulting error terms into the left hand side, provided that $R_{0}$
has been fixed sufficiently large in terms of the geometry of $(\mathcal{S},h)$
(in view also of the flat asymptotics of $(\mathcal{S},h)$).

As we did in the proof of Proposition \ref{prop:DegenerateEllipticEstimates},
we will prove (\ref{eq:EstimateFarAwayFlat}) in detail in the case
$l=2$, and omit the details for the case $l>2$ (which follows in
a similar, albeit notationally more complicated, way). 

Fix a function $w:[0,+\infty)\rightarrow[0,+\infty)$ such that

\begin{itemize}

\item $w(x)=x^{2+\text{\textgreek{e}}}$ for $x\le1$,%
\footnote{Let us note that for more general $l$, one should choose $w=x^{2(l-1)+\text{\textgreek{e}}}$
for $x\le1$.%
}

\item $w(x)=1$ for $x\ge2$ and

\item $\frac{dw}{dx}\ge0$ on $[0,+\infty)$.

\end{itemize}

\noindent We then define the function $w_{R}:\mathbb{R}^{d}\rightarrow(0,+\infty)$
by the relation 
\begin{equation}
w_{R}=R^{2+\text{\textgreek{e}}}w(\frac{r}{R}).
\end{equation}

Fixing also a smooth cut-off $\text{\textgreek{q}}_{R_{0}}:\mathbb{R}^{d}\rightarrow[0,1]$
such that $\text{\textgreek{q}}_{R_{0}}\equiv0$ on $\{r\le R_{0}\}$
and $\text{\textgreek{q}}_{R_{0}}\equiv1$ on $\{r\ge R_{0}\}$, we
obtain after integrating by parts: 
\begin{equation}
\begin{split}\int_{\mathbb{R}^{d}}\text{\textgreek{q}}_{R_{0}}w_{R}\cdot & r^{-\text{\textgreek{b}}}(\text{\textgreek{D}}_{\mathbb{R}^{d}}u)^{2}\, dvol_{e}=\int_{\mathbb{R}^{d}}\text{\textgreek{q}}_{R_{0}}w_{R}\cdot r^{-\text{\textgreek{b}}}|\nabla^{2}u|_{e}^{2}\, dvol_{e}+\int_{\mathbb{R}^{d}}\text{\textgreek{q}}_{R_{0}}w_{R}\cdot\big(\nabla_{\text{\textgreek{m}}}\nabla_{\text{\textgreek{n}}}r_{+}^{-\text{\textgreek{b}}}-\frac{1}{2}(\text{\textgreek{D}}_{h}r_{+}^{-\text{\textgreek{b}}})h_{\text{\textgreek{m}\textgreek{n}}}\big)\cdot\nabla^{\text{\textgreek{m}}}u\cdot\nabla^{\text{\textgreek{n}}}u\, dvol_{e}+\\
 & +\int_{\mathbb{R}^{d}}O_{R_{0}}(|\nabla\text{\textgreek{q}}_{R_{0}}|_{e}+|\nabla^{2}\text{\textgreek{q}}_{R_{0}}|_{e})|\nabla u|_{e}^{2}\, dvol_{e}+\int_{\mathbb{R}^{d}}\text{\textgreek{q}}_{R_{0}}O(|\nabla w_{R}|_{e}\cdot|\nabla r^{-\text{\textgreek{b}}}|_{e}+|\nabla^{2}w_{R}|_{e}r^{-\text{\textgreek{b}}})|\nabla u|_{e}^{2}\, dvol_{e}.
\end{split}
\label{eq:IntegrationByPartsEasy}
\end{equation}
 Notice that the boudary terms at infinity obtained through this integration
by parts procedure vanish. This follows from the fact that, because
$\sum_{j=1}^{l}\int_{\{r\ge R_{0}\}}r^{-2(l-j)-\text{\textgreek{b}}}|\nabla^{j}u|_{e}^{2}\, dvol_{e}<\infty$,
exactly as in the proof of Lemma \ref{lem:AbscenceOfHarmonicFunctions},
we can find a sequence of positive numbers $\{R_{n}\}_{n\in\mathbb{N}}$
tending to $+\infty$ so that 
\begin{equation}
\lim_{n\rightarrow\infty}\Big(R_{n}\sum_{j=1}^{l}\int_{\{r=R_{n}\}}r^{-2(l-j)-\text{\textgreek{b}}}|\nabla^{j}u|_{e}^{2}\, dvol_{\{r=R_{n}\}}\Big)=0.
\end{equation}
Thus, in view of the relation 
\begin{equation}
\nabla_{\text{\textgreek{m}}}\nabla_{\text{\textgreek{n}}}r_{+}^{-\text{\textgreek{b}}}-\frac{1}{2}(\text{\textgreek{D}}_{h}r_{+}^{-\text{\textgreek{b}}})h_{\text{\textgreek{m}\textgreek{n}}}=\frac{1}{2}\text{\textgreek{b}}\Big(\big(d-\text{\textgreek{b}}-4\big)h_{\text{\textgreek{m}\textgreek{n}}}+2(\text{\textgreek{b}}+2)dr_{\text{\textgreek{m}}}\cdot dr_{\text{\textgreek{n}}}\Big)r^{-2-\text{\textgreek{b}}}+o(r^{-2-\text{\textgreek{b}}}),\label{eq:LowerOrdeByWeights}
\end{equation}
 using a Hardy type inequality we readily obtain from (\ref{eq:IntegrationByPartsEasy})
after adding to both sides of (\ref{eq:IntegrationByPartsEasy}) the
quantity $\int_{\{R_{0}\le r\le2R_{0}\}}\big(|\nabla^{2}u|_{e}^{2}+|\nabla u|_{e}^{2}\big)$
(notice that we have assumed $\bar{\text{\textgreek{d}}}_{k_{0}}>0$
to be small enough): 
\begin{multline}
\int_{\{R\le r\le2R\}}r^{2+\text{\textgreek{e}}}\big(r^{-\text{\textgreek{b}}}|\nabla^{2}u|_{e}^{2}+r^{-2-\text{\textgreek{b}}}|\nabla u|_{e}^{2}\big)\, dvol_{e}+\int_{\{r\ge R_{0}\}}\big(r^{-\text{\textgreek{b}}}|\nabla^{2}u|_{e}^{2}+r^{-2-\text{\textgreek{b}}}|\nabla u|_{e}^{2}\big)\, dvol_{e}\le C_{R,R_{0},\text{\textgreek{b}},\text{\textgreek{e}}}\int_{\{r\ge R_{0}\}}r^{-\text{\textgreek{b}}}(\text{\textgreek{D}}_{\mathbb{R}^{d}}u)^{2}\, dvol_{e}+\\
+C_{R_{0},\text{\textgreek{b}},\text{\textgreek{e}}}\int_{\{R_{0}\le r\le2R\}}r^{\text{\textgreek{e}}-\text{\textgreek{b}}}|\nabla u|_{e}^{2}\, dvol_{e}+C_{R_{0},\text{\textgreek{b}},\text{\textgreek{e}}}\int_{\{R_{0}\le r\le2R_{0}\}}\big(|\nabla^{2}u|_{e}^{2}+|\nabla u|_{e}^{2}\big)\, dvol_{e}.\label{eq:StraightForwardEllipticEstimates}
\end{multline}

Let us consider the seminorm space $\mathcal{H}_{R_{0},\text{\textgreek{b}}}^{2}$
defined as the completion of the space $C^{\infty}(\{r\ge R_{0}\})$
with the seminorm 
\begin{equation}
||\text{\textgreek{y}}||_{\mathcal{H}_{R_{0},\text{\textgreek{b}}}^{2}}^{2}\doteq\int_{\{r\ge R_{0}\}}\big(r^{-\text{\textgreek{b}}}|\nabla^{2}\text{\textgreek{y}}|_{e}^{2}+r^{-2-\text{\textgreek{b}}}|\nabla\text{\textgreek{y}}|_{e}^{2}\big)\, dvol_{e}.\label{eq:SemiDefiniteNorm}
\end{equation}
Notice that $\mathcal{H}_{R_{0},\text{\textgreek{b}}}^{2}$ modulo
the constant functions becomes a Hilbert space.

The subspace of harmonic functions 
\begin{equation}
\mathcal{V}_{hrm}\doteq\Big\{\text{\textgreek{y}}\in\mathcal{H}_{R_{0},\text{\textgreek{b}}}^{2}\big|\text{\textgreek{D}}_{\mathbb{R}^{d}}\text{\textgreek{y}}=0\Big\}
\end{equation}
 is a closed subspace of $\mathcal{H}_{R_{0},\text{\textgreek{b}}}^{2}$.
If we introduce the following semi-definite inner product on $\mathcal{H}_{R_{0},\text{\textgreek{b}}}^{2}$:
\begin{equation}
\langle\text{\textgreek{y}}_{1},\text{\textgreek{y}}_{2}\rangle_{R_{0}}\doteq\int_{\{R_{0}\le r\le2R_{0}\}}\Big(\nabla^{\text{\textgreek{m}}}\nabla^{\text{\textgreek{n}}}\text{\textgreek{y}}_{1}\cdot\nabla_{\text{\textgreek{m}}}\nabla_{\text{\textgreek{n}}}\text{\textgreek{y}}_{2}+\nabla^{\text{\textgreek{m}}}\text{\textgreek{y}}_{1}\cdot\nabla_{\text{\textgreek{m}}}\text{\textgreek{y}}_{2}\Big)\, dvol_{e},\label{eq:BoundaryInnerProduct}
\end{equation}
then $\langle\cdot,\cdot\rangle_{R_{0}}$ is continuous with respect
to $||\cdot||_{\mathcal{H}_{R_{0},\text{\textgreek{b}}}^{2}}$, and
for any $\text{\textgreek{y}}\in\mathcal{V}_{hrm}$ which is \underline{not}
a constant we can bound 
\begin{equation}
\langle\text{\textgreek{y}},\text{\textgreek{y}}\rangle_{R_{0}}>0.
\end{equation}
 This follows from the fact that if $\langle\text{\textgreek{y}},\text{\textgreek{y}}\rangle_{R_{0}}=0$,
i.\,e.~if $\text{\textgreek{y}}$ is constant on $\{R_{0}<r<2R_{0}\}$,
and $\text{\textgreek{y}}$ is harmonic (i.\,e.~belongs to $\mathcal{V}_{hrm}$),
then $\text{\textgreek{y}}$ must be identically constant. Therefore,
the orthogonal complement of $\mathcal{V}_{hrm}$ with respect to
$\langle\cdot,\cdot\rangle_{R_{0}}$, that is the subspace 
\begin{equation}
\mathcal{V}_{orth}\doteq\Big\{\text{\textgreek{y}}\in\mathcal{H}_{R_{0},\text{\textgreek{b}}}^{2}\,\big|\,\forall\text{\textgreek{f}}\in\mathcal{V}_{hrm}:\,\mbox{\ensuremath{\langle\text{\textgreek{y}}},\textgreek{f}\ensuremath{\rangle_{R_{0}}}}=0\Big\},
\end{equation}
 is a closed subspace of $\mathcal{H}_{R_{0},\text{\textgreek{b}}}^{2}$,
satisfying $\mathcal{V}_{hrm}\cap\mathcal{V}_{orth}=<1>$. Moreover,
we can decompose any $\text{\textgreek{y}}\in\mathcal{H}_{R_{0},\text{\textgreek{b}}}^{2}$
as 
\begin{equation}
\text{\textgreek{y}}=\text{\textgreek{y}}_{hrm}+\text{\textgreek{y}}_{orth},
\end{equation}
 uniquely modulo addition of some constant function, where $\text{\textgreek{y}}_{hrm}\in\mathcal{V}_{hrm}$
and $\text{\textgreek{y}}_{orth}\in\mathcal{V}_{orth}$. 

\noindent \emph{Remark. }We should emphasize that we will \underline{not}
need to establish that the resulting projection of $\mathcal{H}_{R_{0},\text{\textgreek{b}}}^{2}/<1>$
onto $\mathcal{V}_{orth}/<1>$ along $\mathcal{V}_{hrm}/<1>$ is continuous
with respect to the topology of $\mathcal{H}_{R_{0},\text{\textgreek{b}}}^{2}/<1>$.

Let us introduce the semi-norm 
\begin{equation}
||\text{\textgreek{y}}||_{\mathcal{H}_{R_{0},\text{\textgreek{b}}}^{1}}^{2}\doteq\int_{\{r\ge R_{0}\}}r^{-\text{\textgreek{b}}}|\nabla\text{\textgreek{y}}|_{e}^{2}\, dvol_{e}.\label{eq:SemiDefiniteNorm-1-1}
\end{equation}
Moreover, fixing a sufficiently small $\text{\textgreek{e}}_{0}>0$,
using the Rellich--Kondrachov theorem (see e.\,g.~\cite{Hebey1999}),
we can establish that the subset $\mathcal{D}$ of functions $\text{\textgreek{y}}$
in $\mathcal{V}_{orth}$ satisfying 
\begin{equation}
||\text{\textgreek{y}}||_{\mathcal{H}_{R_{0},\text{\textgreek{b}}}^{1}}=1\label{eq:UnitBall}
\end{equation}
 and 
\begin{equation}
\int_{\{R_{0}\le r\le2R\}}|\nabla\text{\textgreek{y}}|_{e}^{2}\, dvol_{e}\ge\text{\textgreek{e}}_{0}R^{-2-\text{\textgreek{e}}}||\text{\textgreek{y}}||_{\mathcal{H}_{R_{0},\text{\textgreek{b}}}^{2}}^{2}
\end{equation}
 is a pre-compact subset of the semi-norm space $\mathcal{H}_{R_{0},\text{\textgreek{b}}}^{1}$.
Therefore, since no constant function lies in the closure of $\mathcal{D}$
with the semi norm $||\cdot||_{\mathcal{H}_{R_{0},\text{\textgreek{b}}}^{1}}$
(due to (\ref{eq:UnitBall})) and for any non constant $\text{\textgreek{y}}\in\mathcal{V}_{orth}$
we have $||\text{\textgreek{D}}_{\mathbb{R}^{d}}\text{\textgreek{y}}||_{L_{R_{0},\text{\textgreek{b}}}^{2}}>0$,
where 
\begin{equation}
||\text{\textgreek{f}}||_{L_{R_{0},\text{\textgreek{b}}}^{2}}\doteq\int_{\{r\ge R_{0}\}}r^{-\text{\textgreek{b}}}\text{\textgreek{f}}^{2}\, dvol_{e},
\end{equation}
 there exists some large $C_{R}>0$ so that we can bound for any $\text{\textgreek{y}}\in\mathcal{D}$:
\begin{equation}
\int_{\{R_{0}\le r\le2R\}}r^{\text{\textgreek{e}}-\text{\textgreek{b}}}|\nabla\text{\textgreek{y}}|_{e}^{2}\, dvol_{e}\le C_{R,R_{0},\text{\textgreek{b}},\text{\textgreek{e}}}||\text{\textgreek{D}}_{\mathbb{R}^{d}}\text{\textgreek{y}}||_{L_{R_{0},\text{\textgreek{b}}}^{2}}^{2}.\label{eq:LowFrequencyEstimate}
\end{equation}

Therefore, for any $\text{\textgreek{y}}\in\mathcal{H}_{R_{0},\text{\textgreek{b}}}^{2}$,
using (\ref{eq:StraightForwardEllipticEstimates}) in case $\text{\textgreek{y}}_{orth}\notin\{\text{\textgreek{l}}\mathcal{D}\,|\,\text{\textgreek{l}}\ge0\}$
and (\ref{eq:LowFrequencyEstimate}) in case $\text{\textgreek{y}}_{orth}\in\{\text{\textgreek{l}}\mathcal{D}\,|\,\text{\textgreek{l}}\ge0\}$
(and recalling the definition of $\mathcal{D}$), we obtain (provided
$\text{\textgreek{e}}_{0}$ was chosen sufficiently small): 
\begin{multline}
\int_{\{R\le r\le2R\}}r^{2+\text{\textgreek{e}}}\big(r^{-\text{\textgreek{b}}}|\nabla^{2}\text{\textgreek{y}}_{orth}|_{e}^{2}+r^{-2-\text{\textgreek{b}}}|\nabla\text{\textgreek{y}}_{orth}|_{e}^{2}\big)\, dvol_{e}+\int_{\{r\ge R_{0}\}}\big(r^{-\text{\textgreek{b}}}|\nabla^{2}\text{\textgreek{y}}_{orth}|_{e}^{2}+r^{-2-\text{\textgreek{b}}}|\nabla\text{\textgreek{y}}_{orth}|_{e}^{2}\big)\, dvol_{e}\le\\
\le C_{R,R_{0},\text{\textgreek{b}},\text{\textgreek{e}}}\int_{\{r\ge R_{0}\}}r^{-\text{\textgreek{b}}}(\text{\textgreek{D}}_{\mathbb{R}^{d}}\text{\textgreek{y}}_{orth})^{2}\, dvol_{e}+C_{R_{0},\text{\textgreek{b}},\text{\textgreek{e}}}<\text{\textgreek{y}}_{orth},\text{\textgreek{y}}_{orth}>.\label{eq:EllipticOrthogonal}
\end{multline}

We will now establish the necessary estimates for functions belonging
to $\mathcal{V}_{hrm}$. For any function $\text{\textgreek{y}}$
solving $\text{\textgreek{D}}_{\mathbb{R}^{d}}\text{\textgreek{y}}=0$
on $\{r\ge R_{0}\}$ and having finite $||\cdot||_{\mathcal{H}_{R_{0},\text{\textgreek{b}}}^{2}}$
norm, we readily deduce after decomposing it into spherical harmonics
and using (\ref{eq:PlainODE}) and (\ref{eq:ExplicitSlution}) (as
well as the fact that $2+\text{\textgreek{e}}-\text{\textgreek{b}}<d$)
that 
\begin{equation}
\int_{\{R\le r\le2R\}}r^{2-\text{\textgreek{b}}+\text{\textgreek{e}}}\big(|\nabla^{2}\text{\textgreek{y}}|_{e}^{2}+r^{-2}|\nabla\text{\textgreek{y}}|_{e}^{2}\big)\, dvol_{e}\le C_{R_{0},\text{\textgreek{b}},\text{\textgreek{e}}}\int_{\{R_{0}\le r\le2R_{0}\}}\big(|\nabla^{2}\text{\textgreek{y}}|_{e}^{2}+|\nabla\text{\textgreek{y}}|_{e}^{2}\big)\, dvol_{e},
\end{equation}
while the estimate 
\begin{equation}
\int_{\{r\ge R_{0}\}}\big(r^{-\text{\textgreek{b}}}|\nabla^{2}\text{\textgreek{y}}|_{e}^{2}+r^{-2-\text{\textgreek{b}}}|\nabla\text{\textgreek{y}}|_{e}^{2}\big)\, dvol_{e}\le C_{R_{0},\text{\textgreek{b}}}\int_{\{R_{0}\le r\le2R_{0}\}}\big(|\nabla^{2}\text{\textgreek{y}}|_{e}^{2}+|\nabla\text{\textgreek{y}}|_{e}^{2}\big)\, dvol_{e}
\end{equation}
follows readily after integrating by parts in the expression 
\begin{equation}
\int\text{\textgreek{q}}_{R_{0}}r^{-\text{\textgreek{b}}}(\text{\textgreek{D}}_{\mathbb{R}^{d}}\text{\textgreek{y}})^{2}\, dvol_{e}=0.
\end{equation}
 Thus, for any $\text{\textgreek{y}}\in\mathcal{H}_{R_{0},\text{\textgreek{b}}}^{2}$
we can bound 
\begin{equation}
\int_{\{R\le r\le2R\}}r^{2-\text{\textgreek{b}}+\text{\textgreek{e}}}\big(|\nabla^{2}\text{\textgreek{y}}_{hrm}|_{e}^{2}+r^{-2}|\nabla\text{\textgreek{y}}_{hrm}|_{e}^{2}\big)\, dvol_{e}+\int_{\{r\ge R_{0}\}}\big(r^{-\text{\textgreek{b}}}|\nabla^{2}\text{\textgreek{y}}_{hrm}|_{e}^{2}+r^{-2-\text{\textgreek{b}}}|\nabla\text{\textgreek{y}}_{hrm}|_{e}^{2}\big)\, dvol_{e}\le C_{R_{0},\text{\textgreek{b}},\text{\textgreek{e}}}\langle\text{\textgreek{y}}_{hrm},\text{\textgreek{y}}_{hrm}\rangle_{R_{0}}.\label{eq:HarmonicEstimate}
\end{equation}

Therefore, adding (\ref{eq:EllipticOrthogonal}) and (\ref{eq:HarmonicEstimate})
for $u$ in place of $\text{\textgreek{y}}$ and using a triangle
inequality and the fact that 
\begin{equation}
\langle u,u\rangle_{R_{0}}=\langle u_{orth},u_{orth}\rangle_{R_{0}}+\langle u_{hrm},u_{hrm}\rangle_{R_{0}},
\end{equation}
 we readily deduce the desired bound (\ref{eq:EstimateFarAwayFlat})
for $l=2$: 
\begin{multline*}
\int_{\{R\le r\le2R\}}r^{2+\text{\textgreek{e}}}\big(r^{-\text{\textgreek{b}}}|\nabla^{2}u|_{e}^{2}+r^{-2-\text{\textgreek{b}}}|\nabla u|_{e}^{2}\big)\, dvol_{e}+\int_{\{r\ge R_{0}\}}\big(r^{-\text{\textgreek{b}}}|\nabla^{2}u|_{e}^{2}+r^{-2-\text{\textgreek{b}}}|\nabla u|_{e}^{2}\big)\, dvol_{e}\le\\
\le C_{R,R_{0},\text{\textgreek{b}},\text{\textgreek{e}}}\int_{\{r\ge R_{0}\}}r^{-\text{\textgreek{b}}}(\text{\textgreek{D}}_{\mathbb{R}^{d}}u)^{2}\, dvol_{e}+C_{R_{0},\text{\textgreek{b}},\text{\textgreek{e}}}\int_{\{R_{0}\le r\le2R_{0}\}}\big(|\nabla^{2}u|_{e}^{2}+|\nabla u|_{e}^{2}\big)\, dvol_{e}.
\end{multline*}

Inequality (\ref{eq:EstimateFarAwayFlat}) for $2<l\le\lfloor\frac{d+1}{2}\rfloor$
follows in a similar way, and hence the details will be omitted.
\end{proof}

\subsection{A lemma on harmonic functions on $(\mathcal{S},h)$}
\begin{lem}
\label{lem:AbscenceOfHarmonicFunctions}Let $\bar{\text{\textgreek{d}}}>0$
be small in terms of the geometry of $(\mathcal{S},h)$. If a function
$u:\mathcal{S}\rightarrow\mathbb{R}$ with 
\begin{equation}
\int_{\mathcal{S}}r_{+}^{-\text{\textgreek{b}}}\big|\nabla^{2}u\big|_{h}^{2}\, dvol_{\tilde{h}}+\int_{\mathcal{S}}\big(-\log(r_{-})+1\big)^{-2}r_{-}^{-1}r_{+}^{-2-\text{\textgreek{b}}}\big|\nabla u\big|_{h}^{2}\, dvol_{\tilde{h}}+\int_{\{r\le r_{0}\}}|Yu|^{2}\, dvol_{\tilde{h}}<+\infty\label{eq:FinitenessHighNorm}
\end{equation}
for some $\text{\textgreek{b}}\in(-\bar{\text{\textgreek{d}}},1)$
and $r_{0}>0$ solves 
\begin{equation}
\text{\textgreek{D}}_{h,\text{\textgreek{w}}}u=0
\end{equation}
 satisfying the following boundary condition on $\partial_{tim}\mathcal{S}$
: 
\begin{equation}
\int_{\partial_{tim}\mathcal{S}}u\cdot Yu\, dh_{tan}\ge0,\label{eq:BoundaryConditionEllipticWithConstant}
\end{equation}
then $u$ is necessarily a constant function.\end{lem}
\begin{proof}
By standard elliptic regularity results (see i.\,e.~\cite{Gilbarg2001}),
$u\in C^{\infty}(\mathcal{S}\backslash\partial\mathcal{S})$. Let
us fix some $\bar{\text{\textgreek{d}}}<\text{\textgreek{e}}\ll1$
small enough in terms of $1-\text{\textgreek{b}}$ and the geometry
of $(\mathcal{S},h)$ (this is possible since $\bar{\text{\textgreek{d}}}$
was considered small in terms of $(\mathcal{S},h)$). In this way,
$0<\text{\textgreek{b}}+\text{\textgreek{e}}<1$. 

Suppose first that 
\begin{equation}
\int_{\mathcal{S}}\Big(r_{+}^{-\text{\textgreek{b}}-\text{\textgreek{e}}}|\nabla u|_{h}^{2}+r_{+}^{-2-\text{\textgreek{b}}-\text{\textgreek{e}}}u^{2}\Big)\, dvol_{h}<+\infty.\label{eq:EnergyBoundednessEpsilon}
\end{equation}
 In this case, we will show that $u\equiv0$. 

We can define for $\text{\textgreek{r}}\ge R_{0}$ (large in terms
of the geometry of $(\mathcal{S},h)$) the function: 
\begin{equation}
f(\text{\textgreek{r}})\doteq\int_{\{r=\text{\textgreek{r}}\}}r^{-\text{\textgreek{b}}-\text{\textgreek{e}}}\Big(|\nabla u|_{h}^{2}+r^{-2}u^{2}\Big)\, dvol_{h,\{r=\text{\textgreek{r}}\}}.
\end{equation}
 Due to the flat asymptotics of $(\mathcal{S},h)$, from (\ref{eq:EnergyBoundednessEpsilon})
we deduce that: 
\begin{equation}
\int_{R_{0}}^{+\infty}f(\text{\textgreek{r}})\, d\text{\textgreek{r}}<+\infty.\label{eq:BoundednessIntegral}
\end{equation}
 Therefore, since the function $\frac{1}{\text{\textgreek{r}}}$ is
not integrable in $\{\text{\textgreek{r}}\gg1\}$ , using the pigeonhole
principle we deduce that there exists a sequence $\{R_{n}\}_{n\in\mathbb{N}}$
tending to $+\infty$ such that 
\begin{equation}
R_{n}\cdot f(R_{n})\rightarrow0.\label{eq:DecayOnNearDyadic}
\end{equation}

Without loss of generality, we can assume that $r_{0}$ small in terms
of the geometry of $(\mathcal{S},h)$. In the $[0,r_{0}]\times\partial_{hor}\mathcal{S}$
coordinate chart near $\partial_{hor}\mathcal{S}$, $h$ takes the
form 
\begin{equation}
h=\big(r^{-1}+O(1)\big)dr^{2}+h_{tan},
\end{equation}
and thus our main assumption 
\begin{equation}
\int_{\mathcal{S}}\big(-\log(r_{-})+1\big)^{-2}r_{-}^{-1}r_{+}^{-2-\text{\textgreek{b}}}\big|\nabla u\big|_{h}^{2}\, dvol_{\tilde{h}}+\int_{\{r\le r_{0}\}}\big|Yu\big|^{2}\, dvol_{\tilde{h}}<+\infty\label{eq:AssumptionBoundednessEnergy}
\end{equation}
 implies (through a Hardy inequality for the zeroth order term) that
\begin{equation}
\int_{0}^{r_{0}}\int_{\partial_{hor}\mathcal{S}}r^{-1}\big(-\log(r)\big)^{-2}\Big(r(\partial_{r}u)^{2}+|\nabla^{(h_{tan})}u|_{h_{tan}}^{2}+|u|^{2}\Big)\, dh_{tan}dr<+\infty.\label{eq:BoundednessInCoordinatesHorizon}
\end{equation}

Similarly, in the $[0,r_{0}]\times\partial_{tim}\mathcal{S}$ region,
$h$ takes the form
\begin{equation}
h=\big(1+O(r)\big)dr^{2}+h_{tan},
\end{equation}
 and thus (\ref{eq:AssumptionBoundednessEnergy}) implies that 
\begin{equation}
\int_{0}^{r_{0}}\int_{\partial_{tim}\mathcal{S}}r^{-1}\big(-\log(r)\big)^{-2}\Big((\partial_{r}u)^{2}+|\nabla^{(h_{tan})}u|_{h_{tan}}^{2}+|u|^{2}\Big)\, dh_{tan}dr<+\infty.\label{eq:BoundednessInCoordinatesTimelike}
\end{equation}

Setting for $\text{\textgreek{r}}\ll1$: 
\begin{equation}
g(\text{\textgreek{r}})\doteq\int_{\{r=\text{\textgreek{r}}\}}\big(-\log(r)\big)^{-1}\Big(r_{hor}(\partial_{r}u)^{2}+|\nabla^{(h_{tan})}u|_{h_{tan}}^{2}+|u|^{2}\Big)\, dh_{tan},\label{eq:BoundaryFunctionToBeZero}
\end{equation}
 since $f(r)=r^{-1}\log^{-1}(r)$ is not integrable around $0$,  from
$\eqref{eq:BoundednessInCoordinatesHorizon}$, (\ref{eq:BoundednessInCoordinatesTimelike})
and a trivial pigeonhole principle argument we infer the existence
of a sequence $\{r_{n}\}_{n\in\mathbb{N}}$ tending to $0$ so that:
\begin{equation}
g(r_{n})\rightarrow0.\label{eq:SequenceNearBoundary}
\end{equation}

Let us fix a $C^{1}$ and piecewise $C^{2}$ function $w_{R_{1}}:\mathcal{S}\rightarrow(0,1]$
which satisfies the following properties for some $c>0$ depending
only on the geometry of $\tilde{h}$:

\begin{itemize}

\item $w_{R_{1}}$ is a function of $r$,

\item $w_{R_{1}}=1$ for $r<R_{1}$,

\item $\partial_{r}w_{R_{1}}\le0$, $\partial_{r}^{2}w_{R_{1}}<-c\cdot R_{1}^{-2}$
for $R_{1}<r<2R_{1}$ and 

\item $w_{R_{1}}=\big(\frac{r}{R_{1}}\big)^{-\text{\textgreek{b}}-\text{\textgreek{e}}}$
for $r>2R_{1}$. 

\end{itemize}

\noindent Since $(\mathcal{S},h)$ has a finite number of asymptotically
flat regions, we compute that in the $(r,\text{\textgreek{sv}})$
coordinate system on each connected component of the region $\{r\gg1\}$:
\begin{equation}
\text{\textgreek{D}}_{h,\text{\textgreek{w}}}w_{R_{1}}=\big(1+O(r^{-1})\big)\cdot\partial_{r}^{2}w_{R_{1}}+\big(d-1+O(r^{-1})\big)r^{-1}\partial_{r}w_{R_{1}}.\label{eq:LaplacianOnw}
\end{equation}

\noindent On the other hand, in the region $\{r<R_{1}\}$ we have
$\text{\textgreek{D}}_{h,\text{\textgreek{w}}}w_{R_{1}}\equiv0$ since
$w_{R_{1}}$ is a constant there.

Therefore, if $R_{1}$ is fixed large and $r_{1}$ is fixed small
in terms the geometry of $(\mathcal{S},h)$ and $\text{\textgreek{w}}$,
from (\ref{eq:LaplacianOnw}) and the properties of $w_{R_{1}}$ (and
the fact that $0<\text{\textgreek{b}}+\text{\textgreek{e}}<1$) we
infer that 
\begin{equation}
\text{\textgreek{D}}_{h,\text{\textgreek{w}}}w_{R_{1}}\le0\label{eq:NegativeLaplacian}
\end{equation}
 almost everywhere on $\mathcal{S}$ (recall that $w_{R_{1}}$ is
$C^{1}$ but only piecewise $C^{2}$).

Let us fix the vector fields $\text{\textgreek{n}}_{1}=-|\nabla r|_{h}^{-1}\nabla r$
in the region $\{r\ll1\}$ and $\text{\textgreek{n}}_{2}=|\nabla r|_{h}^{-1}\nabla r$
in the region $\{r\gg1\}$. Since $\text{\textgreek{D}}_{h,\text{\textgreek{w}}}u=0$
by assumption, after integrating by parts (and using the fact that
$w_{R_{1}}$ is $C^{1}$ and piecewise $C^{2}$) we obtain for any
integer $n$ sufficiently large: 
\begin{align}
0=-\int_{\{r_{n}\le r\le R_{n}\}}w_{R_{1}}\text{\textgreek{w}}\cdot\text{\textgreek{D}}_{h,\text{\textgreek{w}}}u\cdot u\, dvol_{h}= & \int_{\{r_{n}\le r\le R_{n}\}}\Big(\text{\textgreek{w}}\cdot w_{R_{1}}|\nabla u|_{h}^{2}-\frac{1}{2}\text{\textgreek{w}}\cdot\text{\textgreek{D}}_{h,\text{\textgreek{w}}}w_{R_{1}}\cdot u^{2}\Big)\, dvol_{h}-\label{eq:BeforeTakingLimitsGreen}\\
 & -\int_{\{r=r_{n}\}}\Big(\text{\textgreek{w}}\cdot w_{R_{1}}\text{\textgreek{n}}_{1}u\cdot u-\frac{1}{2}\text{\textgreek{w}}(\text{\textgreek{n}}_{1}w_{R_{1}})\cdot u^{2}\Big)\, dvol_{h,\{r=r_{n}\}}-\nonumber \\
 & -\int_{\{r=R_{n}\}}\Big(\text{\textgreek{w}}\cdot w_{R_{1}}\text{\textgreek{n}}_{2}u\cdot u-\frac{1}{2}\text{\textgreek{w}}(\text{\textgreek{n}}_{2}w_{R_{1}})\cdot u^{2}\Big)dvol_{h,\{r=R_{n}\}}.\nonumber 
\end{align}

In view if the form of $h$ in the $[0,r_{0}]\times\partial\mathcal{S}$
coordinate chart (see (\ref{eq:MetricNearHorizon}) and (\ref{eq:MetricNearTimelikeBoundary})),
we compute that 
\begin{equation}
dvol_{h,\{r=r_{n}\}}=dh_{tan}|_{r=r_{n}},
\end{equation}
\begin{equation}
\text{\textgreek{n}}_{1}=-r_{hor}^{\frac{1}{2}}(1+O(r))\partial_{r},
\end{equation}
and 
\begin{equation}
|\nabla u|_{h}^{2}=r_{hor}(\partial_{r}u)^{2}+|\nabla^{(h_{tan})}u|_{h_{tan}}^{2}.
\end{equation}
 Since $\int_{0}^{r_{0}}\int_{\partial_{tim}\mathcal{S}}|\big(\nabla^{(\tilde{h})}\big)^{2}u|_{h}^{2}\, dh_{tan}dr<+\infty$
and $h$ is regular up to $\partial_{tim}\mathcal{S}$, by applying
a trace theorem we infer that $Yu$ has a well defined limit on $L^{2}(\partial_{tim}\mathcal{S},h_{S})$.
Thus, since $w_{R_{1}}=1$ and $\text{\textgreek{w}}\sim r_{hor}^{\frac{1}{2}}$
near $\partial\mathcal{S}$, we obtain in view of (\ref{eq:BoundaryFunctionToBeZero})
and (\ref{eq:SequenceNearBoundary}): 
\begin{equation}
\limsup_{n\rightarrow+\infty}\int_{\{r=r_{n}\}}\text{\textgreek{w}}\Big(w_{R_{1}}\text{\textgreek{n}}_{1}u\cdot u-\frac{1}{2}(\text{\textgreek{n}}_{1}w_{R_{1}})\cdot u^{2}\Big)\, dvol_{h,\{r=r_{n}\}}\le-\int_{\partial_{tim}\mathcal{S}}Yu\cdot u\, dh_{tan}+0\le0,\label{eq:FirstBoundaryTrivial}
\end{equation}
because of (\ref{eq:BoundaryConditionEllipticWithConstant}).

Similarly, by a Cauchy--Schwarz inequality we have because of (\ref{eq:DecayOnNearDyadic})
as $n\rightarrow+\infty$: 
\begin{align}
\Big|\int_{\{r=R_{n}\}}\text{\textgreek{w}}\Big(w_{r_{1},R_{1}}\text{\textgreek{n}}_{2}u\cdot u-\frac{1}{2}(\text{\textgreek{n}}_{2}w_{r_{1},R_{1}})\cdot u^{2}\Big)dvol_{h,\{r=R_{n}\}}\Big| & \le C_{R_{1},\text{\textgreek{e}},\text{\textgreek{b}}}\cdot\int_{\{r=\text{\textgreek{r}}\}}r^{-\text{\textgreek{e}}-\text{\textgreek{b}}}\Big(|\nabla u|_{h}\cdot u+r^{-1}u^{2}\Big)\, dvol_{h,\{r=\text{\textgreek{r}}\}}\label{eq:SecondBoundaryTrivial}\\
 & \le C_{R_{1},\text{\textgreek{e}},\text{\textgreek{b}}}R_{n}\cdot f(R_{n})\rightarrow0\nonumber 
\end{align}
 Thus, by letting $n\rightarrow+\infty$ in (\ref{eq:BeforeTakingLimitsGreen}),
from (\ref{eq:FirstBoundaryTrivial}), (\ref{eq:SecondBoundaryTrivial})
and (\ref{eq:NegativeLaplacian}) we deduce that: 
\begin{equation}
\int_{\mathcal{S}}\text{\textgreek{w}}\Big(w_{R_{1}}|\nabla u|_{h}^{2}+w_{2}\cdot u^{2}\Big)\, dvol_{h}\le0
\end{equation}
for some suitably decaying $w_{2}\ge0$ which is not identically $0$.
Therefore, $u\equiv0$.

In order to establish Lemma \ref{lem:AbscenceOfHarmonicFunctions},
therefore, it suffices to establish that for any $u\in C^{\infty}(\mathcal{S})$
satisfying $\text{\textgreek{D}}_{h,\text{\textgreek{w}}}u=0$ and
(\ref{eq:FinitenessHighNorm}), there exists some constant $c_{u}$
so that $u-c_{u}$ satisfies (\ref{eq:EnergyBoundednessEpsilon}).
In view of (\ref{eq:FinitenessHighNorm}), it suffices to show that
\begin{equation}
\int_{\{r\ge2R_{u}\}}r^{-\text{\textgreek{e}}-\text{\textgreek{b}}}\big(|\nabla u|_{h}^{2}+r^{-2}(u-c_{u})^{2}\big)\, dvol_{h}<+\infty\label{eq:WeakerEpsilonBoundedness}
\end{equation}
 for some $R_{u}$ large depending on $u$ itself. 

We will work in the $(r,\text{\textgreek{sv}})$ coordinate system
on a single connected component $\mathcal{N}_{1}$ of the region $\{r\gg1\}$
(since the proof for each component is identical, this is not actually
a restriction). In this coordinate system on $\mathcal{N}_{1}$, we
will define the coordinate flat metric 
\begin{equation}
e=dr^{2}+r^{2}g_{\mathbb{S}^{d-1}},\label{eq:FlatMetric}
\end{equation}
 and the associated flat Laplacian: 
\begin{equation}
\text{\textgreek{D}}_{\mathbb{R}^{d}}\doteq r^{-(d-1)}\partial_{r}\big(r^{d-1}\partial_{r}\big)+r^{-2}\text{\textgreek{D}}_{\mathbb{S}^{d-1}}.\label{eq:FlatLaplacian}
\end{equation}
 Since $\text{\textgreek{D}}_{h,\text{\textgreek{w}}}u=0$ and $\text{\textgreek{D}}_{h,\text{\textgreek{w}}}$
has the asymptotics (\ref{eq:PerturbedLaplacianAway}), while $h$
is asymptotically of the form (\ref{eq:EuclideanAsymptotics}), $u$
satisfies the equation 
\begin{equation}
\text{\textgreek{D}}_{\mathbb{R}^{d}}u=F_{u},\label{eq:FlatEquation}
\end{equation}
 where 
\begin{equation}
F_{u}=O(r^{-1})\nabla^{2}u+O(r^{-2})\nabla u.\label{eq:NonHomogeneity}
\end{equation}

We will first show that there exists a smooth solution $u_{1}$ to
the boundary value problem 
\begin{equation}
\begin{cases}
\text{\textgreek{D}}_{\mathbb{R}^{d}}u_{1}=F_{u}, & \{r\ge R_{u}\}\cap\mathcal{N}_{1}\\
u_{1}=0 & \mbox{ on }\{r=R_{u}\}\cap\mathcal{N}_{1}
\end{cases}\label{eq:BoundaryValueProblem}
\end{equation}
 with: 
\begin{equation}
\int_{\{r\ge R_{u}\}\cap\mathcal{N}_{1}}r^{-\text{\textgreek{b}}}\big(|\nabla u_{1}|_{e}^{2}+r^{-2}u_{1}^{2}\big)\, dvol_{e}<+\infty.\label{eq:BoundednessEnergyU1}
\end{equation}

To this end, considering the sequence of cut-off functions $\text{\textgreek{q}}_{R_{n}}=\text{\textgreek{q}}(\frac{r}{R_{n}})$
for some dyadic sequence $\{R_{n}\}_{n\in\mathbb{N}}$ and some fixed
smooth cut-off function $\text{\textgreek{q}}:\mathbb{R}\rightarrow[0,1]$
satisfying $\text{\textgreek{q}}(x)=1$ for $x\le1$ and $\text{\textgreek{q}}(x)=0$
for $x\ge2$, we first solve the boundary value problems: 
\begin{equation}
\begin{cases}
\text{\textgreek{D}}_{\mathbb{R}^{d}}u_{1,n}=\text{\textgreek{q}}_{R_{n}}\cdot F_{u}, & \{r>R_{u}\}\cap\mathcal{N}_{1}\\
u_{1,n}=0 & \mbox{ on }\{r=R_{u}\}\cap\mathcal{N}_{1}.
\end{cases}\label{eq:SequenceBoundaryValueProblems}
\end{equation}

The existence of solutions $u_{1,n}$ to (\ref{eq:SequenceBoundaryValueProblems})
readily follows using the variational approach: Let us define the
Hilbert space $\mathcal{H}_{R_{u}}^{1}$ as the completion of the
vector space $C_{0,R_{u}}^{\infty}$, 
\begin{equation}
C_{0,R_{u}}^{\infty}=\Big\{\text{\textgreek{y}}\in C_{0}^{\infty}(\{r\ge R_{u}\}\cap\mathcal{N}_{1})\big|\text{\textgreek{y}}|_{r=R_{u}}=0\Big\}
\end{equation}
 with the norm 
\begin{equation}
||\text{\textgreek{y}}||_{\mathcal{H}_{R_{u}}^{1}}\doteq\Big(\int_{\{r>R_{u}\}\cap\mathcal{N}_{1}}\big(|\nabla\text{\textgreek{y}}|_{e}^{2}+r^{-2}\text{\textgreek{y}}^{2}\big)\, dvol_{e}\Big)^{1/2}.
\end{equation}
Then the function $\mathcal{F}_{n}:C_{0,R_{u}}^{\infty}\rightarrow\mathbb{R}$
given by the relation 
\begin{equation}
\mathcal{F}_{n}(\text{\textgreek{y}})\doteq\int_{\{r>R_{u}\}\cap\mathcal{N}_{1}}\big(\frac{1}{2}|\nabla\text{\textgreek{y}}|_{e}^{2}-\text{\textgreek{q}}_{R_{n}}F_{u}\cdot\text{\textgreek{y}}\big)\, dvol_{e}
\end{equation}
extends to a continuous function on $\mathcal{H}_{R_{u}}^{1}$, since
$\text{\textgreek{q}}_{R_{n}}F_{u}$ is compactly supported. Thus,
we can bound through a Cauchy--Schwarz inequality for any $v\in\mathcal{H}_{R_{u}}^{1}$
\begin{equation}
\int_{\{r>R_{u}\}\cap\mathcal{N}_{1}}\big|\text{\textgreek{q}}_{R_{n}}F_{u}\cdot v\big|\, dvol_{e}\lesssim_{R_{n}}||v||_{\mathcal{H}_{R_{u}}^{1}}.
\end{equation}
 Moreover, $\mathcal{F}_{n}$ is convex and satisfies $\mathcal{F}_{n}(\text{\textgreek{y}})\rightarrow+\infty$
as $||\text{\textgreek{y}}||_{\mathcal{H}_{R_{u}}^{1}}\rightarrow+\infty$.
By minimizing $\mathcal{F}_{n}$ over $\mathcal{H}_{R_{u}}^{1}$ using
the usual techniques (and using elliptic regularity results), we arrive
at a strong solution $u_{1,n}$ of (\ref{eq:SequenceBoundaryValueProblems})
belonging to $\mathcal{H}_{R_{u}}^{1}$.

Since $||u_{1,n}||_{\mathcal{H}_{R_{u}}^{1}}<+\infty$, using the
approach leading to (\ref{eq:DecayOnNearDyadic}) we deduce that there
exists a dyadic sequence $\{\text{\textgreek{r}}_{m}^{(n)}\}_{m\in\mathbb{N}}$
such that 
\begin{equation}
\text{\textgreek{r}}_{m}^{(n)}\cdot\int_{\{r=\text{\textgreek{r}}_{m}^{(n)}\}\cap\mathcal{N}_{1}}\Big(|\nabla u|_{e}^{2}+r^{-2}u^{2}\Big)\, dvol_{e,\{r=\text{\textgreek{r}}\}}\le1.
\end{equation}
 Hence, by multiplying (\ref{eq:SequenceBoundaryValueProblems}) with
$r^{-\text{\textgreek{b}}}u_{1,n}$ and integrating by parts over
$\{R_{u}<r<\text{\textgreek{r}}_{m}^{(n)}\}\cap\mathcal{N}_{1}$,
we obtain after taking the limit $m\rightarrow+\infty$: 
\begin{equation}
\int_{\{r>R_{u}\}\cap\mathcal{N}_{1}}\big(r^{-\text{\textgreek{b}}}|\nabla u_{1,n}|_{e}^{2}-\frac{1}{2}(\text{\textgreek{D}}_{\mathbb{R}^{d}}r^{-\text{\textgreek{b}}})u_{1,n}^{2}\big)\, dvol_{e}\le\int_{\{r>R_{u}\}\cap\mathcal{N}_{1}}\text{\textgreek{q}}_{R_{n}}r^{-\text{\textgreek{b}}}F_{u}\cdot u_{1,n}\, dvol_{e}.\label{eq:EnergyDegenerateAtInfinity}
\end{equation}

Notice that 
\begin{equation}
\text{\textgreek{D}}_{\mathbb{R}^{d}}r^{-\text{\textgreek{b}}}=-\text{\textgreek{b}}(d-2-\text{\textgreek{b}})\cdot r^{-\text{\textgreek{b}}-2}\le\begin{cases}
0, & \text{\textgreek{b}}\in[0,1)\\
(d-1)\cdot\bar{\text{\textgreek{d}}}r^{-\text{\textgreek{b}}-2} & \text{\textgreek{b}}\in(-\bar{\text{\textgreek{d}}},0]
\end{cases}\label{eq:SmallPositiveBound}
\end{equation}
 since $d\ge3$. Using also a Hardy-type inequality of the form established
in Lemma (\ref{lem:HardyForphi}), we deduce from (\ref{eq:EnergyDegenerateAtInfinity})
and (\ref{eq:SmallPositiveBound}) (and the assumption that $\bar{\text{\textgreek{d}}}$
is small in terms of the geometry of $(\mathcal{S},h)$):
\begin{equation}
\int_{\{r>R_{u}\}\cap\mathcal{N}_{1}}\big(r^{-\text{\textgreek{b}}}|\nabla u_{1,n}|_{e}^{2}+r^{-2-\text{\textgreek{b}}}u_{1,n}^{2}\big)\, dvol_{e}\le C_{\text{\textgreek{b}}}\int_{\{r>R_{u}\}\cap\mathcal{N}_{1}}\text{\textgreek{q}}_{R_{n}}r^{-\text{\textgreek{b}}}F_{u}\cdot u_{1,n}\, dvol_{e}.\label{eq:EnergyDegenerateAtInfinity-1}
\end{equation}
 Using a Cauchy--Schwarz inequality, from (\ref{eq:EnergyDegenerateAtInfinity-1})
we infer that 
\begin{equation}
\int_{\{r>R_{u}\}\cap\mathcal{N}_{1}}\big(r^{-\text{\textgreek{b}}}|\nabla u_{1,n}|_{e}^{2}+r^{-2-\text{\textgreek{b}}}u_{1,n}^{2}\big)\, dvol_{e}\le C_{\text{\textgreek{b}}}\int_{\{r>R_{u}\}\cap\mathcal{N}_{1}}\text{\textgreek{q}}_{R_{n}}r^{2-\text{\textgreek{b}}}|F_{u}|^{2}\, dvol_{e}\le C_{\text{\textgreek{b}},u}<+\infty\label{eq:EnergyDegenerateAtInfinityFinal}
\end{equation}
 with the constants not depending on $n$. Similarly, we can also
bound: 
\begin{equation}
\int_{\{r>R_{u}\}\cap\mathcal{N}_{1}}\big(r^{-\text{\textgreek{b}}}|\nabla(u_{1,n+1}-u_{1,n})|_{e}^{2}+r^{-2-\text{\textgreek{b}}}(u_{1,n+1}-u_{1,n})^{2}\big)\, dvol_{e}\le C_{\text{\textgreek{b}}}\int_{\{r>R_{u}\}\cap\mathcal{N}_{1}}(\text{\textgreek{q}}_{R_{n+1}}-\text{\textgreek{q}}_{R_{n}})r^{2-\text{\textgreek{b}}}|F_{u}|^{2}\, dvol_{e}.\label{eq:EnergyBooundDifference}
\end{equation}

From (\ref{eq:NonHomogeneity}) and (\ref{eq:FinitenessHighNorm})
we can bound 
\begin{equation}
\int_{\{r>R_{u}\}\cap\mathcal{N}_{1}}r^{2-\text{\textgreek{b}}}|F_{u}|^{2}\, dvol_{e}\le C_{u}<\infty.
\end{equation}
 Thus, from (\ref{eq:EnergyDegenerateAtInfinityFinal}) and (\ref{eq:EnergyBooundDifference})
we infer that the sequence of functions $u_{1,n}$ converges to a
function $u_{1}$ which is a weak (and hence strong, due to elliptic
regularity) solution of (\ref{eq:BoundaryValueProblem}) satisfying
(\ref{eq:BoundednessEnergyU1}).

Moreoever, after integrating by parts in the expression 
\[
\int_{\{r\ge R_{u}\}\cap\mathcal{N}_{1}}r^{-\text{\textgreek{b}}}(1-\text{\textgreek{q}}_{2R_{u}})(\text{\textgreek{D}}_{\mathbb{R}^{d}}u_{1})^{2}\, dvol_{e}=\int_{\{r\ge R_{u}\}\cap\mathcal{N}_{1}}r^{-\text{\textgreek{b}}}(1-\text{\textgreek{q}}_{2R_{u}})F_{u}^{2}\, dvol_{e}
\]
 and using (\ref{eq:BoundednessEnergyU1}), (\ref{eq:FinitenessHighNorm})
and the fact that $(\{r\ge R_{u}\}\cap\mathcal{N}_{1},e)$ is flat,
we obtain that: 
\begin{equation}
\int_{\{r\ge2R_{u}\}\cap\mathcal{N}_{1}}r^{-\text{\textgreek{b}}}|\nabla^{2}u_{1}|_{e}^{2}\, dvol_{e}<+\infty.\label{eq:BoundsSecondDerivative}
\end{equation}

We will now return to our function $u\in C^{\infty}(\mathcal{S})$
solving $\text{\textgreek{D}}_{h,\text{\textgreek{w}}}u=0$ and satisfying
(\ref{eq:FinitenessHighNorm}). Setting on $\{r\ge2R_{u}\}\cap\mathcal{N}_{1}$
\begin{equation}
u_{2}=u-u_{1},
\end{equation}
the new function $u_{2}$ will satisfy on $\{r>2R_{u}\}\cap\mathcal{N}_{1}$:
\begin{equation}
\text{\textgreek{D}}_{\mathbb{R}^{d}}u_{2}=0.\label{eq:HarmonicEquation}
\end{equation}
 (being smooth by elliptic regularity) and 
\begin{equation}
\int_{\{r\ge2R_{u}\}\cap\mathcal{N}_{1}}r^{-\text{\textgreek{b}}}\big(|\nabla^{2}u_{2}|_{e}^{2}+r^{-2}|\nabla u_{2}|_{e}^{2}\big)\, dvol_{e}<+\infty\label{eq:HighBoundU2}
\end{equation}
 because of (\ref{eq:FinitenessHighNorm}), (\ref{eq:BoundsSecondDerivative})
and (\ref{eq:BoundednessEnergyU1}).

For any $\text{\textgreek{r}}\ge2R_{u}$, $u_{2}|_{\{r=\text{\textgreek{r}}\}\cap\mathcal{N}_{1}}$
is a smooth function on $\mathbb{S}^{d-1}$ (by elliptic regularity)
and hence we can decompose $u{}_{2}|_{\{r=\text{\textgreek{r}}\}\cap\mathcal{N}_{1}}$
in spherical harmonics. This turns (\ref{eq:HarmonicEquation}) into
the following system of ODE's: 
\begin{equation}
r^{-(d-1)}\frac{d}{dr}\big(r^{d-1}\frac{d}{dr}(u_{2})_{m}\big)-\text{\textgreek{L}}_{m}r^{-2}(u_{2})_{m}=0,\label{eq:PlainODE}
\end{equation}
where the integer $m$ corresponds to an enumeration of the spherical
haromics, $\text{\textgreek{y}}_{m}$ denotes the projection of a
function $\text{\textgreek{y}}\in L^{2}(\mathbb{S}^{d-1})$ on $m$-th
spherical harmonic, and $-\text{\textgreek{L}}_{m}$ is the eigenvalue
of $\text{\textgreek{D}}_{\mathbb{S}^{d-1}}$ corresponding to the
$m$-th spherical harmonic (by convention $\text{\textgreek{L}}_{m}$
is increasing in $m$ and $\text{\textgreek{L}}_{0}=0$). Therefore,
we can explicitly solve: 
\begin{equation}
(u_{2})_{m}=c_{m}r^{a_{m}}+d_{m}r^{-(d-2)-a_{m}}\label{eq:ExplicitSlution}
\end{equation}
 where:

\begin{itemize}

\item $c_{m},d_{m}\in\mathbb{R}$

\item $a_{m}\in\mathbb{N}$ is an increasing sequence with $a_{0}=0$
and $a_{1}=1$.

\end{itemize}

By (\ref{eq:HighBoundU2}) we deduce that 
\begin{equation}
\sum_{m\in\mathbb{N}}\int_{r=2R_{u}}^{+\infty}r^{-\text{\textgreek{b}}}\Big\{\big(\frac{d^{2}(u_{2})_{m}}{dr^{2}}\big)^{2}+(\text{\textgreek{L}}_{m}+1)r^{-2}\big(\frac{d(u_{2})_{m}}{dr}\big)^{2}+\big(\text{\textgreek{L}}_{m}^{2}+1\big)r^{-4}(u_{2})_{m}^{2}\Big\}\, r^{d-1}dr<+\infty,\label{eq:BoundInFrequency}
\end{equation}
 which, in view of the fact that $\text{\textgreek{b}}\in(-\bar{\text{\textgreek{d}}},1)$,
forces $c_{m}=0$ for $m\neq0$ in (\ref{eq:ExplicitSlution}). Therefore,
from (\ref{eq:ExplicitSlution}), (\ref{eq:BoundInFrequency}) and
the fact that $c_{m}=0$ for $m\neq0$, we infer that for all $m\in\mathbb{N}$:
\begin{equation}
\sum_{m\in\mathbb{N}}\int_{r=2R_{u}}^{+\infty}r^{-\text{\textgreek{b}}-\text{\textgreek{e}}}\Big\{\big(\frac{d(u_{2})_{m}}{dr}\big)^{2}+\text{\textgreek{L}}_{m}r^{-2}(u_{2})_{m}^{2}\Big\}\, r^{d-1}dr<+\infty\label{eq:BoundFirstEnergy}
\end{equation}
(recall that $\text{\textgreek{L}}_{0}=0$). Therefore 
\begin{equation}
\int_{\{r\ge2R_{u}\}\cap\mathcal{N}_{1}}r^{-\text{\textgreek{b}}-\text{\textgreek{e}}}|\nabla u_{2}|_{e}^{2}\, dvol_{e}<+\infty.\label{eq:BoundFirstEnergyU@}
\end{equation}

Furthermore, using standard ode theory (and elliptic regularity),
we can establish in this case that the series $\sum_{m}(u_{2})_{m}e_{m}$
converges to $u_{2}$ pointwise (together with all its derivatives).
Moreover, because of (\ref{eq:ExplicitSlution}) and the fact that
$c_{m}=0$ for $m\neq0$, we deduce that: 
\begin{equation}
\lim_{r\rightarrow+\infty}u_{2}=c_{0}.\label{eq:LimitAtInfinity}
\end{equation}
 Thus, (\ref{eq:BoundFirstEnergyU@}) implies that (through a Hardy-type
inequality) that: 
\begin{equation}
\int_{\{r\ge2R_{u}\}\cap\mathcal{N}_{1}}r^{-\text{\textgreek{b}}-\text{\textgreek{e}}}\Big(|\nabla u_{2}|_{e}^{2}+r^{-2}(u_{2}-c_{0})^{2}\Big)\, dvol_{e}<+\infty.\label{eq:U2BoundFirstEnergy}
\end{equation}

All in all, from (\ref{eq:BoundednessEnergyU1}) and (\ref{eq:U2BoundFirstEnergy})
and the fact that $u=u_{1}+u_{2}$, we finally infer (\ref{eq:WeakerEpsilonBoundedness})
for $c_{u}=c_{0}$: 
\begin{equation}
\int_{\{r\ge2R_{u}\}\cap\mathcal{N}_{1}}r^{-\text{\textgreek{b}}-\text{\textgreek{e}}}\Big(|\nabla u|_{e}^{2}+r^{-2}(u-c_{0})^{2}\Big)\, dvol_{e}<+\infty.
\end{equation}
 Hence, the proof of the Lemma is complete.
\end{proof}

\section{Hardy type inequalities}

In this Section, we will state the main Hardy type inequalities that
are used throughout this paper. We will start with the following lemma:
\begin{lem}
\label{lem:HardyGeneral}For any strictly increasing function $g:(0,+\infty)\rightarrow\mathbb{R}$,
the following inequality holds for any $u\in C^{1}\big((0,+\infty)\big)$
and $0<a<b$: 
\begin{equation}
\int_{a}^{b}\frac{dg}{dx}\cdot|u|^{2}\, dx+g(a)|u(a)|^{2}\le C_{g}\cdot\Big(\int_{a}^{b}g^{2}\Big(\frac{dg}{dx}\Big)^{-1}\cdot\Big|\frac{du}{dx}\Big|^{2}\, dx+g(b)|u(b)|^{2}\Big)\label{eq:GeneralHardyIncreasing}
\end{equation}
 while for $g$ strictly decreasing we have: 
\begin{equation}
\int_{a}^{b}\big(-\frac{dg}{dx}\big)|u|^{2}\, dx+g(b)|u(b)|^{2}\le C_{g}\cdot\Big(\int_{a}^{b}\big(-g^{2}\Big(\frac{dg}{dx}\Big)^{-1}\big)\cdot\Big|\frac{du}{dx}\Big|^{2}\, dx+g(a)|u(a)|^{2}\Big).\label{eq:GeneralHardyDecreasing}
\end{equation}
\end{lem}
\begin{proof}
The proof follows readily after performing an integration by parts
in the terms $\int_{a}^{b}\frac{dg}{dx}|u|^{2}$ and $\int_{a}^{b}\big(-\frac{dg}{dx}\big)|u|^{2}$
respectively, and then using a Cauchy--Schwarz inequality.\end{proof}
\begin{lem}
\label{lem:HardyForphi}For any $k\in\mathbb{N}$ and $a\in\mathbb{R}$,
there exists a constant $C_{k,a}>0$ such that for any function $\text{\textgreek{y}}\in C^{\infty}(\mathbb{R}^{d})$
and any $0<R_{1}<R_{2}$ we can bound 
\begin{multline}
\sum_{j=1}^{\min\{\lfloor\frac{d-1+a}{2}\rfloor,k\}}\Big\{\int_{\{R_{1}\le r\le R_{2}\}}r^{a-2j}|\partial_{r}^{k-j}\text{\textgreek{y}}|^{2}\, dvol_{e}+\int_{\{r=R_{1}\}}r^{a+1-2j}|\partial_{r}^{k-j}\text{\textgreek{y}}|^{2}\, dvol_{\{r=R_{1}\}}\Big\}\le\\
\le C_{k,a}\cdot\Big\{\int_{\{R_{1}\le r\le R_{2}\}}r^{a}|\partial_{r}^{k}\text{\textgreek{y}}|^{2}\, dvol_{e}+\sum_{j=1}^{\lfloor\frac{d-1+a}{2}\rfloor}\int_{\{r=R_{2}\}}r^{a+1-2j}|\partial_{r}^{k-j}\text{\textgreek{y}}|^{2}\, dvol_{\{r=R_{2}\}}\Big\},\label{eq:HardyPrototypeHyperboloids}
\end{multline}
 where $\partial_{r}$ denotes the gradient of the polar distance
$r$ on $\mathbb{R}^{d}$, $dvol_{e}$ is the natural volume form
of the flat metric $e$ on $\mathbb{R}^{d}$ and $dvol_{\{r=R_{i}\}}$
is the natural volume form of the induced metric on the sphere $\{r=R_{i}\}$.\end{lem}
\begin{proof}
From (\ref{eq:GeneralHardyIncreasing}) for $g(x)=x^{b+1}$ for any
$b>-1$, it follows that for any function $h\in C^{\infty}(\mathbb{R})$:
\begin{equation}
\int_{R_{1}}^{R_{2}}x^{b}|h|^{2}\, dx+x^{b+1}|h(R_{1})|^{2}\le C_{b}\cdot\Big(\int_{R_{1}}^{R_{2}}x^{b+2}|\frac{dh}{dx}|^{2}\, dx+x^{b+1}|h(R_{2})|^{2}\Big)\label{eq:Hardy1}
\end{equation}
(with the constant $C_{b}$ depending only on $b$). Therefore, recalling
that in polar coordinates $(r,\text{\textgreek{sv}})$ on $\mathbb{R}^{d}$
we can write for any $b>0$ and $h\in C^{\infty}(\mathbb{R}^{d})$:
\[
\int_{\{R_{1}\le r\le R_{2}\}}r^{b}|h|^{2}=\int_{R_{1}}^{R_{2}}r^{b+d-1}\Big(\int_{\mathbb{S}^{d-1}}h(r,\text{\textgreek{sv}})\, d\text{\textgreek{sv}}\Big)\, dr,
\]
 inequality (\ref{eq:HardyPrototypeHyperboloids}) follows readily
after repeated applications of (\ref{eq:Hardy1}).
\end{proof}
Let $(\mathcal{M},g)$ be as in Section \ref{sec:Firstdecay}. Recall
that we denote with $\nabla_{h_{\text{\textgreek{t}},N}}$ the covariant
derivative associated to the Riemannian metric $h_{\text{\textgreek{t}},N}$
on the hyperboloid $\{\bar{t}=\text{\textgreek{t}}\}$ (see Section
\ref{sec:Firstdecay}). Let us also denote with $L$ the $\partial_{v}$
coordinate vector field in the $(v,\text{\textgreek{sv}})$ coordinate
chart on each connected component of $\{\bar{t}=\text{\textgreek{t}}\}\cap\{r\gg1\}$
(notice that this vector field does \underline{not} coincide with
the $\partial_{v}$ vector field on $\mathcal{M}$ in the $(u,v,\text{\textgreek{sv}})$
coordinate chart). Using Lemma \ref{lem:HardyForphi}, we can establish
the following Hardy type inequality for functions $\text{\textgreek{y}}$
on the hyperboloids $\{\bar{t}=const\}$ with tame behaviour near
$\mathcal{I}^{+}$:
\begin{lem}
\label{lem:HardyForphiHyperboloids} There exists an $R_{0}>0$ (large)
such that for any $\text{\textgreek{t}}>0$, any $k\in\mathbb{N}$,
any integer $0\le m\le k$, any $a\in\mathbb{R}$, any smooth tensor
field $\text{\textgreek{Y}}$ on $\{\bar{t}=\text{\textgreek{t}}\}$
satisfying $\lim_{r\rightarrow+\infty}r^{\frac{d+a-2}{2}}|\mathcal{L}_{L}^{l}\text{\textgreek{Y}}|_{h_{\text{\textgreek{t}},N}}=0$
for $0\le l\le k$, and any $R_{1}>R_{0}$ we can bound 
\begin{multline}
\sum_{j=1}^{\min\{\lfloor\frac{d-1+a}{2}\rfloor,k\}}\Big\{\int_{\{r\ge R_{1}\}\cap\{\bar{t}=\text{\textgreek{t}}\}}r^{a-2j}|\mathcal{L}_{L}^{k-j}\text{\textgreek{Y}}|^{2}\,\text{\textgreek{W}}^{2}dvd\text{\textgreek{sv}}+\int_{\{r=R_{1}\}\cap\{\bar{t}=\text{\textgreek{t}}\}}r^{a+1-2j}|\mathcal{L}_{L}^{k-j}\text{\textgreek{Y}}|^{2}\,\text{\textgreek{W}}^{2}d\text{\textgreek{sv}}\Big\}\le\\
\le C_{k,a}\cdot\int_{\{r\ge R_{1}\}\cap\{\bar{t}=\text{\textgreek{t}}\}}r^{a}|\mathcal{L}_{L}^{k}\text{\textgreek{Y}}|^{2}\,\text{\textgreek{W}}^{2}dvd\text{\textgreek{sv}},\label{eq:HardyHyperboloids}
\end{multline}
 where $\mathcal{L}_{L}$ denotes the Lie derivative in the direction
of $L$.\end{lem}
\begin{proof}
Fix an $R_{2}>R_{1}$. An application of Lemma \ref{lem:HardyForphi}
on each connected component of $\{R_{1}\le r\le R_{2}\}\cap\{\bar{t}=\text{\textgreek{t}}\}$(using
the coordinate chart $(v,\text{\textgreek{sv}})$ in the region $\{r\ge R_{1}\}\cap\{\bar{t}=\text{\textgreek{t}}\}$,
so that $\mathcal{L}_{L}$ becomes differentiation with respect to
$\partial_{v}$) readily yields that:

\begin{multline}
\sum_{j=1}^{\min\{\lfloor\frac{d-1+a}{2}\rfloor,k\}}\Big\{\int_{\{R_{1}\le r\le R_{2}\}\cap\{\bar{t}=\text{\textgreek{t}}\}}r^{a-2j}|\mathcal{L}_{L}^{k-j}\text{\textgreek{Y}}|^{2}\,\text{\textgreek{W}}^{2}dvd\text{\textgreek{sv}}+\int_{\{r=R_{1}\}\cap\{\bar{t}=\text{\textgreek{t}}\}}r^{a+1-2j}|\mathcal{L}_{L}^{k-j}\text{\textgreek{Y}}|^{2}\,\text{\textgreek{W}}^{2}d\text{\textgreek{sv}}\Big\}\le\\
\le C_{k,a}\cdot\Big(\int_{\{R_{1}\le r\le R_{2}\}\cap\{\bar{t}=\text{\textgreek{t}}\}}r^{a}|\mathcal{L}_{L}^{k}\text{\textgreek{Y}}|^{2}\,\text{\textgreek{W}}^{2}dvd\text{\textgreek{sv}}+\sum_{j=1}^{\lfloor\frac{d-1+a}{2}\rfloor}\int_{\{r=R_{2}\}\cap\{\bar{t}=\text{\textgreek{t}}\}}r^{a+1-2j}|\mathcal{L}_{L}^{k-j}\text{\textgreek{Y}}|^{2}\,\text{\textgreek{W}}^{2}d\text{\textgreek{sv}}\Big).\label{eq:HardyHyperboloids-1}
\end{multline}
Since $\lim_{r\rightarrow+\infty}r^{\frac{d+a-2}{2}}|\mathcal{L}_{L}^{l}\text{\textgreek{Y}}|_{h_{\text{\textgreek{t}},N}}=0$
for $0\le l\le k$, we infer: 
\begin{equation}
\lim_{R_{2}\rightarrow+\infty}\sum_{j=1}^{\lfloor\frac{d-1+a}{2}\rfloor}\int_{\{r=R_{2}\}\cap\{\bar{t}=\text{\textgreek{t}}\}}r^{a+1-2j}|\mathcal{L}_{L}^{k-j}\text{\textgreek{Y}}|^{2}\,\text{\textgreek{W}}^{2}d\text{\textgreek{sv}}=0.
\end{equation}
 Thus, (\ref{eq:HardyHyperboloids}) follows from (\ref{eq:HardyHyperboloids-1})
after taking the limit $R_{2}\rightarrow+\infty$.
\end{proof}
We will also need the following ``critical'' Hardy type inequality:
\begin{lem}
\label{lem:HardyCritical}For any $R>0$ large in terms of the geometry
of $\{\bar{t}=\text{\textgreek{t}}\}$ and any smooth function $\text{\textgreek{y}}$
on $\{\bar{t}=\text{\textgreek{t}}\}$ satisfying $\lim_{r\rightarrow+\infty}|\text{\textgreek{y}}|\rightarrow0$
we can bound: 
\begin{equation}
\int_{\{R\le r\le2R\}\cap\{\bar{t}=\text{\textgreek{t}}\}}r^{-d}|\text{\textgreek{y}}|^{2}\, dh_{N}\le C\cdot\Big(\int_{\{r\ge R\}\cap\{\bar{t}=\text{\textgreek{t}}\}}r^{-d}|\text{\textgreek{y}}|^{2}\, dh_{N}\Big)^{\frac{1}{2}}\Big(\int_{\{r\ge R\}\cap\{\bar{t}=\text{\textgreek{t}}\}}r^{-d+2}|L\text{\textgreek{y}}|^{2}\, dh_{N}\Big)^{\frac{1}{2}}.\label{eq:CriticalHardy}
\end{equation}
\end{lem}
\begin{rem*}
The exponent of the second term of the right hand side of (\ref{eq:CriticalHardy})
can not be increased, since then the inequality would not be satisfied
by the function $\text{\textgreek{y}}_{R_{0}}(r)=\frac{\log(R_{0})}{\log(r)+\log(R_{0})}$
for some large enough $R_{0}>0$.\end{rem*}
\begin{proof}
It suffices to establish the following inequality on $\mathbb{R}^{d}$
for any $R>0$ and any real $\text{\textgreek{y}}\in C^{\infty}(\mathbb{R}^{d})$
sarisfying $\lim_{r\rightarrow+\infty}|\text{\textgreek{y}}|\rightarrow0$:
\begin{equation}
\int_{\{R\le r\le2R\}}r^{-d}\text{\textgreek{y}}^{2}\, dvol_{e}\le C\cdot\Big(\int_{\{r\ge R\}}r^{-d}\text{\textgreek{y}}^{2}\, dvol_{e}\Big)^{\frac{1}{2}}\Big(\int_{\{r\ge R\}}r^{-d+2}(\partial_{r}\text{\textgreek{y}})^{2}\, dvol_{e}\Big)^{\frac{1}{2}}.\label{eq:CriticalHardyFlatSpace}
\end{equation}

In turn, (\ref{eq:CriticalHardyFlatSpace}) will follow (using polar
coordinates) by the following ``critical'' Hardy-type inequality
on $\mathbb{R}$ for any function $\text{\textgreek{y}}\in C^{\infty}(\mathbb{R})$
with $\lim_{r\rightarrow+\infty}|\text{\textgreek{y}}|\rightarrow0$:
\begin{equation}
\int_{R}^{2R}r^{-1}\text{\textgreek{y}}^{2}\, dr\le C\cdot\Big(\int_{R}^{+\infty}r^{-1}\text{\textgreek{y}}^{2}\, dr\Big)^{\frac{1}{2}}\Big(\int_{R}^{+\infty}r\cdot(\partial_{r}\text{\textgreek{y}})^{2}\, dr\Big)^{\frac{1}{2}}.\label{eq:CriticalHardyR}
\end{equation}

In order to establish (\ref{eq:CriticalHardyR}), let us fix a continuous
and piecewise $C^{1}$ function $\text{\textgreek{q}}_{1}:[0,+\infty)\rightarrow[0,1]$
by the relation:
\begin{equation}
\text{\textgreek{q}}_{1}(x)=\begin{cases}
0, & x\le1\\
x-1, & 1\le x\le2\\
1, & x\ge2,
\end{cases}
\end{equation}
 and define the function $\text{\textgreek{q}}_{1,R}:[0,+\infty)\rightarrow[0,1]$
by the relation: 
\begin{equation}
\text{\textgreek{q}}_{1,R}(r)\doteq\text{\textgreek{q}}_{1}(\frac{r}{R}).\label{eq:Cut-off}
\end{equation}

Using the fact that $\text{\textgreek{q}}_{1}$ is contnuous and piecewise
$C^{1}$, we obtain after integrating by parts (in view of the fact
that $\lim_{r\rightarrow+\infty}\text{\textgreek{y}}=0$): 
\begin{align}
\int_{0}^{+\infty}\partial_{r}\text{\textgreek{q}}_{1,R}\cdot\text{\textgreek{y}}^{2}\, dr & =-2\int_{0}^{+\infty}\text{\textgreek{q}}_{1,R}\cdot\text{\textgreek{y}}\partial_{r}\text{\textgreek{y}}\, dr\le\label{eq:AlmostCriticalHardy-1}\\
 & \le\Big(\int_{0}^{+\infty}\text{\textgreek{q}}_{1,R}r^{-1}\text{\textgreek{y}}^{2}\, dr\Big)^{\frac{1}{2}}\Big(\int_{0}^{+\infty}\text{\textgreek{q}}_{1,R}r\cdot(\partial_{r}\text{\textgreek{y}})^{2}\, dr\Big)^{\frac{1}{2}}.\nonumber 
\end{align}
 Thus, (\ref{eq:CriticalHardyR}) follows from (\ref{eq:AlmostCriticalHardy-1})
in view of the fact that $\text{\textgreek{q}}_{1,R}$ is supported
on $[R,+\infty)$ and $\partial_{r}\text{\textgreek{q}}_{1,R}$ is
identically $1$ on $[R,2R]$ and $0$ elsewhere.
\end{proof}
\bibliographystyle{plain}
\bibliography{DatabaseExample}

\end{document}